\documentclass[notitlepage, abstract=on]{scrreprt}
\usepackage[english]{babel}
\usepackage[utf8]{inputenc}
\usepackage{amsmath, amssymb}
\usepackage[T1]{fontenc}
\usepackage{amsfonts}
\usepackage{graphicx}
\usepackage{extarrows}
\usepackage{mathtools}
\usepackage[all]{xy}
\usepackage[center]{caption}
\usepackage{subfig}
\usepackage{mathrsfs}
\usepackage{mathtools}
\usepackage{slashed}
\usepackage{amsthm}
\usepackage{pifont}
\usepackage{fancyhdr}
\usepackage{hyperref}
\usepackage{titling}
\usepackage[backend=biber,
style=alphabetic,
sorting=nyt,
firstinits=true]{biblatex}
\newtheorem{thm}{Theorem}[chapter]
\newtheorem{lemma}[thm]{Lemma}
\newtheorem{prop}[thm]{Proposition}
\newtheorem{cor}[thm]{Corollary}
\newtheorem*{conjecture}{Conjecture}
\newtheorem*{problem}{Problem}
\newtheorem*{convention}{Convention}

\theoremstyle{definition}
\newtheorem{definition}[thm]{Definition}

\theoremstyle{remark}
\newtheorem{remark}[thm]{Remark}
\newtheorem{example}[thm]{Example}

\theoremstyle{bonuslemma}
\newtheorem{bonuslemma}{Lemma}

\newcommand{\fK}{\mathfrak{K}}
\newcommand{\fB}{\mathfrak{B}}
\newcommand{\fQ}{\mathfrak{Q}}
\newcommand{\fL}{\mathfrak{L}}
\newcommand{\fA}{\mathfrak{A}}
\newcommand{\fD}{\mathfrak{D}}

\newcommand{\fC}{\mathfrak{C}}

\newcommand{\fBH}{\mathfrak{B}(H)}
\newcommand{\fKH}{\mathfrak{K}(H)}

\newcommand{\fP}{\mathfrak{P}}
\newcommand{\cU}{\mathcal{U}}
\newcommand{\cL}{\mathcal{L}}
\newcommand{\cF}{\mathcal{F}}

\newcommand{\cH}{\mathcal{H}}
\newcommand{\cO}{\mathcal{O}}

\newcommand{\RR}{\mathbb{R}}
\newcommand{\CC}{\mathbb{C}}
\newcommand{\NN}{\mathbb{N}}

\newcommand{\TT}{\mathbb{T}}
\newcommand{\HH}{\mathbb{H}}
\newcommand{\ZZ}{\mathbb{Z}}
\newcommand{\QQ}{\mathbb{Q}}

\newcommand{\olM}{\overline{M}}
\newcommand{\olO}{\overline{\Omega}}
\newcommand{\dO}{\partial\Omega}
\newcommand{\dom}{\mathrm{dom}}
\newcommand{\im}{\mathrm{im}}
\newcommand{\id}{\mathrm{id}}
\newcommand{\Ad}{\mathrm{Ad}}
\newcommand{\Spin}{\mathrm{Spin}}

\newcommand{\LLip}{L\text{-}\operatorname{Lip}}
\newcommand{\diam}{\mathrm{diam}}
\newcommand{\relInd}{\mathrm{Ind}^\mathrm{rel}_u}
\newcommand{\Ind}{\mathrm{Ind}_u}
\newcommand{\supp}{\mathrm{supp}}
\newcommand{\ind}{\mathrm{index}}
\newcommand{\colim}{\mathrm{colim}}
\newcommand{\exc}{\mathrm{exc}}
\newcommand{\Hom}{\mathrm{Hom}}
\newcommand{\Dirac}{\slashed{D}}

\newcommand{\lsim}{\lesssim}
\newcommand{\gsim}{\gtrsim}
\newcommand{\sua}{\sim_\mathrm{ua}}
\newcommand{\cmark}{\ding{51}}
\newcommand{\xmark}{\ding{55}}

\title{Uniformly elliptic boundary value problems\thanks{The present document is a slightly modified version of the author's dissertation, which was submitted to and accepted by the University of Greifswald, and will be referred to as \emph{thesis} throughout. The modifications consist of minor corrections and elaborations. The original version of the dissertation is available as the first version of this document, or at the University of Greifswald's library server at \href{https://nbn-resolving.org/urn:nbn:de:gbv:9-opus-143203}{urn:nbn:de:gbv:9-opus-143203}. The indexing of theorems, lemmas, and related statements was kept consistent between versions. Statements present in this version but not the original one are indexed with letters instead of numbers.}}
\author{Matti Lyko\thanks{Universität Greifswald, Walther-Rathenau-Straße 47, 17489 Greifswald, Germany. \href{mailto:matti.lyko@uni-greifswald.de}{matti.lyko@uni-greifswald.de}. I acknowledge financial support by the Deutsche Forschungsgemeinschaft DFG (project number 510244757).}}
\date{}

\begin{document}

\begin{titlingpage}
\maketitle

\begin{abstract}
We study boundary conditions for elliptic operators on non-compact manifolds with boundary via uniform K-homology, a version of K-homology sensitive to the large-scale geometry of the manifold. To that end, we develop the theory of relative uniform K-homology. We show that boundary conditions for uniformly elliptic differential operators define classes in the relative and non-relative uniform K-homology groups of the manifold, depending on the assumed regularity of the boundary condition. Moreover, we define and study a relative index map on relative uniform K-homology that combines uniform coarse information on the interior with secondary information on the boundary. As an application, we compute that on a spin manifold with product structure and uniformly positive scalar curvature on the boundary the image of the relative uniform K-homology class of the Dirac operator under this relative index map is closely connected to a uniform version of the higher $\rho$-invariant of the boundary. In particular, a delocalized APS-index theorem of Piazza and Schick is proved in the uniform setting.
\end{abstract}
\end{titlingpage}

\newpage

\addtocontents{toc}{\protect\thispagestyle{empty}}

\tableofcontents

\thispagestyle{empty}

\newpage

\pagenumbering{arabic}

\chapter{Introduction}
The behavior of differential equations is influenced by the properties of the underlying space. Knowledge of a space can be translated to knowledge of differential equations on it, and vice versa. The study of this interplay blossomed into a field of mathematics which connects analysis, geometry, topology, and more. In the introduction to this thesis I want to describe some of the history of the field, introduce important objects and themes, and sketch the humble contribution this thesis aims to make.
\subsection*{The Atiyah-Singer index theorem}
Suppose we are given a linear differential equation $Du=v$. The two fundamental questions regarding this equation are existence and uniqueness of solutions. Uniqueness is captured by the kernel of the differential operator $D$. Its dimension measures the number of independent solutions of the homogeneous equation $Du=0$. Existence is captured by the cokernel of $D$. Its dimension measures the number of independent $v$ for which the equation $Du=v$ does \emph{not} admit a solution. 

The dimension of both the kernel and the cokernel are sensitive quantities; they generally change upon continuous perturbations of the operator $D$. It turns out that their difference, the so-called \emph{index} 
$$
\ind(D):=\dim(\ker(D))-\dim(\mathrm{coker}(D)) \; ,
$$
is unchanged upon continuous perturbations of $D$. While existence and uniqueness of solutions can be changed by perturbing $D$, this trade-off will always be a zero-sum game. Thus the index provides a more robust measure of these fundamental aspects of the behavior of $D$.

Generally, neither the kernel nor the cokernel of $D$ need to be finite-dimensional, in which case the index is not a well-defined integer. One class of differential operators that do admit a well-behaved index theory are so-called \emph{elliptic} operators, defined by the requirement that their leading-order coefficients be invertible in a suitable way. It is a consequence of the many favorable properties of elliptic operators that they have finite-dimensional kernel and cokernel, at least over a closed manifold. The prototypical example of an elliptic operator is the Laplace operator on a Riemannian manifold. Certain geometric first-order operators closely connected to the Laplacian are elliptic as well, for example the Euler characteristic operator and signature operator on differential forms, Dirac operator associated to spin (or spin$^c$ structures, or the Cauchy-Riemann operator on a complex manifold.

The indices of these geometric first-order operators turn out to compute characteristic numbers of the underlying manifold. For example the index of the Euler operator is the Euler characteristic, that of the signature operator the signature of the manifold. The index of the Cauchy-Riemann operator computes the Euler characteristic in sheaf cohomology. In turn each of these characteristic numbers can be expressed as an integral over certain characteristic classes. Therefore the indices of these classical geometric
operators, which are analytic objects, equal integrals of characteristic classes, which are topological objects. The Atiyah-Singer index is a vast generalization of these index-equals-topology statements. 
\begin{thm}[Atiyah-Singer index theorem, \cite{AS1968}]
Let $M$ be a closed manifold, $E\to M$ a vector bundle, and $D$ an elliptic differential operator over $M$. Denote by $D_E$ the operator $D$ twisted by $E$. Then,
\begin{equation}
\ind(D_E) = \int_M \mathrm{AS}(D,M) \wedge \mathrm{ch}(E) \; . \label{AtiyahSinger}
\end{equation}
Here $\mathrm{AS}(D,M)$ is a characteristic class universally constructed from $D$ and $M$, and $\mathrm{ch}(E)$ is the Chern character of $E$. 
\end{thm}
The class $\mathrm{AS}(D,M)$ depends on $D$ only through the K-theory class of its symbol, a topological quantity derived from $D$. For that reason the right-hand side of \eqref{AtiyahSinger} is called the \emph{topological index} of $D_E$. 

There are many ways to prove the Atiyah-Singer index theorem. Historically the first proof given by Atiyah-Singer was the so-called bordism proof \cite{Palais1965}: One shows that analytic and topological index is invariant under bordisms of the underlying manifold, thus it suffices to prove the index theorem on the generators of the bordism groups. This reduces the proof to explicit calculations. The second proof, the one first published by Atiyah and Singer, is the K-theoreticl proof \cite{AS1968}: As mentioned earlier the topological index depends only on the K-theory class of the symbol of the elliptic operator. The topological index therefore defines an integer-valued map on the K-theory of $M$. One shows that the same holds for the analytic index, essentially expressing the fact that the index of an elliptic operator depends only on its symbol. The topological index map on K-theory is uniquely characterized by its invariance under push-forward maps on K-theory, plus its value at a point. Then, equality of analytic and topological index is established by proving that the analytic index map has the same properties.

A conceptual explanation of why these proofs work might go something like this: The analytic quantity of interest turn out to have enough stability (bordism invariance, dependence only on the K-theory class of the symbol) to define an element of a suitable algebraic-topological object. Then the machinery of algebraic topology can be used to obtain the desired result.

There is a third method of proof of the index theorem, the heat kernel proof. For this one proves that the trace of the heat kernel $e^{-tD^2}$ associated to the operator $D$ is constant as a function of $t$. Its limit for $t\to\infty$ is the analytic index of $D$. An asymptotic expansion for small $t$ shows that the limit for $t\to 0$ is the topological index. The heat kernel proof uses practically no algebraic topology, making it somewhat mysterious that it is able to prove the same result as the other two proofs. However, we see the K-theoretical significance of the heat kernel below.
\subsection*{From the index theorem to K-homology}
Reconsider the index formula \eqref{AtiyahSinger}, and shift the focus onto the vector bundle $E$. It enters the topological index through its Chern character, which only depends on the K-theory class of $E$. Therefore the index extends to a homomorphism
$$
K^0(M) \longrightarrow \ZZ \; , \; [E] \longmapsto \ind(D_E) \; .
$$
For each elliptic operator over $M$ the index provides an element of $\Hom(K^0(M),\ZZ)$. K-theory, viewed as a (generalized) cohomology theory, has a dual homology theory called K-homology. This duality manifests for example through a bilinear pairing $K_0(M)\times K^0(M) \to \ZZ$ or, equivalently, a homomorphism $K_0(M) \to \Hom(K^0(M),\ZZ)$. This leads one to ask: Can we understand elliptic operators as elements of $K_0(M)$? If so, is every element of $K_0(M)$ defined by some elliptic operator?

The answer is affirmative. In fact there is a model for K-homology called analytic K-homology, whose elements are equivalence classes of abstract elliptic operators. The definition of such an object goes as follows: A \emph{Fredholm module} over the compact space $X$ is a triple $(H,\rho,F)$ consisting of a Hilbert space $H$, a representation $\rho: \, C(X) \to \fBH$, and a bounded operator $F\in \fBH$ such that
$$
[\rho(f),F] \quad , \quad \left( I-F^2\right) \rho(f) \quad , \quad \left(F-F^* \right)\rho(f) 
$$
are compact operators for all $f\in C(X)$. The first property is referred to as \emph{pseudo-locality} for $F$. Operators $T$ such that $T\rho(f)$ and $\rho(f)T$ are compact will be referred to as \emph{locally compact}. Thus the operators $I-F^2$ and $F-F^*$ are required to be locally compact. Note that on compact spaces local compactness is equivalent to compactness, but this is no longer true on non-compact spaces.

To explain why elliptic operators give rise to Fredholm modules we need to go through some analytic background. Suppose that $D$ is a symmetric elliptic operator on a closed manifold $M$. There is a scale of  Hilbert spaces $H^k(M)$, $k\in\ZZ_{\geq 0}$, on $M$ called Sobolev spaces. The space $H^k(M)$ is a Hilbert space of $k$-times (weakly) differentiable functions, and $H^0(M)$ is simply the $L^2$-space $L^2(M)$. Moreover, $H^{k+1}(M)$ is continuously included in $H^k(M)$, and the important Rellich-Kondrachov theorem states that this inclusion map is compact. If the differential operator $D$ has order $m$, then it boundedly maps $H^{k+m}(M)$ to $H^k(M)$ for all $k$. It is a core property of elliptic operators that they have quasi-inverses of negative order, meaning there exists an operator $Q$ that boundedly maps $H^k(M)$ to $H^{k+m}$, such that $QD-I$ and $DQ-I$ are smoothing operators, meaning they boundedly map $H^k(M)$ to $H^{k+l}(M)$ for all $l\geq 0$. Note that by the Rellich-Kondrachov theorem this implies that elliptic operators are invertible modulo compact operators, hence are Fredholm operators. This is the aforementioned fact that any elliptic operator over a closed manifold has a well-defined index.

To build a Fredholm module out of the symmetric elliptic operator $D$ let $H$ denote the Hilbert space $L^2(M)$, and $\rho$ the representation of $C(X)$ on $H$ by multiplication operators. Then define $F:=D(D^2+I)^{-\frac{1}{2}}$. First of all we have to argue that $F$ is a bounded operator. Since $D$ is elliptic, so is $D^2+I$. The latter elliptic operator is strictly positive and in particular invertible, and its inverse and quasi-inverse agree up to smoothing operators. In particular $(D^2+I)^{-1}$ boundedly maps $H^k(M)$ to $H^{k+2m}(M)$. Then, the square root $(D^2+I)^{-\frac{1}{2}}$ boundedly maps $H^k(M)$ to $H^{k+m}(M)$.\footnote{In the case $k=0$ -which suffices for our purposes here- we can view $D$ as an unbounded self-adjoint operator on $L^2(M)$ with domain $H^m(M)$. A routine calculation for unbounded operators shows that the square-root $(D^2+I)^{-\frac{1}{2}}$ boundedly maps $L^2(M)$ to the domain of $D$, i.e. to $H^k(M)$. See also Appendix \ref{AppendixUnbdOps}. For general $k$ one has to appeal to the much deeper result that the square root of a positive elliptic operator -in this case $(D^2+I)^{-1}$- is an elliptic operator of half the order \cite{Seeley1967}.} It follows that $F$ is the composition of the bounded maps $(D^2+I)^{-\frac{1}{2}}: \, L^2(M)\to H^m(M)$ and $D: \, H^m(M)\to L^2(M)$. Thus $F$ is bounded. 

To check the required compactness statements note first that $F$ is self-adjoint, so that $(F-F^*)\rho(f)$ is trivially compact. Furthermore, $I-F^2$ is calculated to coincide with $(D^2+I)^{-1}$. We have already noted above that this operator boundedly maps $L^2(M)$ to $H^{2m}(M)$, hence is compact by Rellich-Kondrachov. Lastly, for the compactness of the commutator $[\rho(f),F]$ we appeal to the general fact that the commutator of a pseudo-differential operator (which both $D$ and $(D^2+I)^{-\frac{1}{2}}$ are) has order one lower than the original operator. In our situation this means that 
$$
[\rho(f),D]: \; H^m(M)\to H^1(M) \quad , \quad [\rho(f),(D^2+I)^{-\frac{1}{2}}]: \; L^2(M) \to H^{m+1}(M) \; .
$$
Then, another appeal to the Rellich-Kondrachov theorem provides the compactness of 
$$
[\rho(f),F] = [\rho(f),D] (D^2+I)^{-\frac{1}{2}} + D [\rho(f),(D^2+I)^{-\frac{1}{2}}] \; .
$$
This concludes the proof that $(H,\rho,F)$ is a Fredholm modules. Of course the concept of a Fredholm module was set up precisely so that it would be an abstraction of this fundamental example. 

Analytic K-homology of $X$ is defined as the free group generated by all Fredholm modules under direct sum, modulo relations given by unitary equivalences and homotopies of operators. It provides a model for the homology theory dual to topological K-theory. In particular there is a pairing of analytic K-homology with K-theory defined precisely so that the pairing of the class of the elliptic operator $D$ with the class of a vector bundle $E$ results in the index of the twisted operator $D_E$. The index of an elliptic operator can also be recovered from its K-homology class in a different way. Given a closed connected manifold $M$ the push-forward map $M\to \mathrm{pt}$ to a point induces a homomorphism $K_0(M) \to K_0(\mathrm{pt})$. The group $K_0(\mathrm{pt})$ coincides with $\ZZ$, and the resulting map $K_0(M) \to K_0(\mathrm{pt})=\ZZ$ computes the index. In this way understanding the push-forward to a point furnishes index computations, and it is possible to give a proof of the Atiyah-Singer index theorem in this framework \cite{Baum2018}, \cite{Baum2016}.

Note that the compact space $X$ enters the the definition of analytic K-homology only through its function algebra $C(X)$. Thus analytic K-homology can be defined equally well for an arbitrary $C^*$-algebra. In this way analytic K-homology acquires an operator-algebraic flavor. Indeed, analytic K-homology is a special case of Kasparov's bivariant KK-theory, which unites K-theory and K-homology in a common operator-algebraic framework, and provides a potent tool for applications in geometry, topology, representation theory and others \cite{Kasparov1981}, \cite{Blackadar1998}. 

The operator-algebraic nature allows for a very direct generalization to non-compact spaces: If $X$ is non-compact, we can instead use the algebra $C_0(X)$ of functions vanishing at infinity to define the analytic K-homology $K_*(X)$. If $M$ is a non-compact manifold, then elliptic operators over $M$ still define Fredholm modules and hence K-homology classes. This is in stark difference to the index itself, which is always defined over a closed manifold, but generally not over a non-compact one.

Throughout this introduction we will emphasize the K-homological approach to index theory, and ultimately employ it in our studies in the main part of this thesis. This means that wherever possible index theorems should take the form of an identity for the K-homology classes of elliptic operators, and actual index formulas obtained as a corollary by plugging the K-homological identity into a suitable index map.
\subsection*{An application: Positive scalar curvature metrics}
Let us turn to an important application of the index theorem to a geometric problem. This problem, the question of the existence of metrics of positive scalar curvature, will motivate much of the following discussion. It goes like this: Let $M$ be a manifold, compact and without boundary. Does $M$ admit a Riemannian metric which has positive scalar curvature (psc for short)? 

In the case that $M$ is an even-dimensional spin manifold, the index theorem provides a necessary condition. Consider the spin Dirac operator $\Dirac$ associated to a Riemannian metric on $M$. The Lichnerowicz formula asserts that 
$$
\Dirac^2 = \nabla^* \nabla + \frac{1}{4}\cdot \mathrm{scal}\; ,
$$
where $\nabla$ is the Levi-Civita connection and $\mathrm{scal}_g$ the scalar curvature. If the latter is strictly positive, $\mathrm{scal}>0$, then $\Dirac^2$ is strictly positive, so that $\Dirac$ has trivial kernel. This cannot happen if the index of $\Dirac$ (strictly speaking the index of one of its graded components) is not equal to zero. Thus, the index of the Dirac operator is an obstruction to the existence of a psc metric on $M$. The Atiyah-Singer index theorem allows the computation of this obstruction in terms of purely topological data on $M$; the index of the Dirac operator equals the so-called $\hat{A}$-genus $\hat{A}(M)$. In particular $\ind(\Dirac)=\hat{A}(M)$ is independent of the metric $g$. Combining the Lichnerowicz-formula with the Atiyah-Singer index theorem thus provides a purely topological obstruction to the existence of a psc-metric:
\begin{thm}[Lichnerowicz, \cite{Lichnerowicz1963}]
Let $M$ be an even-dimensional closed spin manifold. If $M$ admits a psc metric, then 
$$
\hat{A}(M) = 0\; .
$$
\end{thm}
\subsection*{Higher indices}
While a useful obstruction, the $\hat{A}$-genus is not powerful enough to answer the psc-question on all manifolds of interest. Consider for example the torus $\TT^n$. It has trivial tangent bundle, so that all characteristic classes of the torus vanish, including the $\hat{A}$-genus. It was suspected that $\TT^n$ does \emph{not} admit a psc metric. Thus one would wish to find an obstruction that is fine enough to detect this. This is achieved via coarse (or higher) index theory.

Let $M$ be a compact connected Riemannian spin manifold, and $\tilde{M}\to M$ its universal cover. The fundamental group $\Gamma$ of $M$ acts on $\tilde{M}$ by deck transformations. The metric and spin structure on $M$ lift to corresponding $\Gamma$-invariant structure on $\tilde{M}$. Moreover, the Dirac operator $\Dirac$ on $M$ lifts to a $\Gamma$-invariant operator $\tilde{\Dirac}$ acting on the spinor bundle on $\tilde{M}$.

It is not immediately clear what we get by lifting to $\tilde{M}$. Indeed, since $\tilde{\Dirac}$ is $\Gamma$-invariant, it is completely determined by its action over a fundamental domain. It turns out that $\tilde{\Dirac}$ has a richer index theory nonetheless. Assume for the moment that $M$ is even-dimensional, so that the spinor bundle on $\tilde{M}$ is graded and $\tilde{\Dirac}$ is an odd operator. Consider the heat operator $\exp(-t\tilde{\Dirac}^2)$. This is a smoothing operator, and the limit of operator of finite propagation\footnote{Intuitively speaking, while heat propagates arbitrarily quickly, the greater the distance, the fewer heat is propagated across that distance. In fact the decay is exponential. Discarding the heat propagated across distances above large cut-offs provides arbitrarily accurate approximations of the heat operator by finite propagation operators. Those operators can be constructed to still be smoothing operators.}. By the Rellich-Kondrachov theorem smoothing operators are locally compact. Thus the heat operator is an element in the $C^*$-algebra $C^*(\tilde{M})$ generated by locally compact finite-propagation operators over $\tilde{M}$. The algebra $C^*(\tilde{M})$ is called the \emph{Roe algebra} of $\tilde{M}$. In a suitable sense the heat operator is actually close to being a projection in the Roe algebra. Thus we obtain an element
$$
\mathrm{Ind}(\tilde{\Dirac}):=\left[ \exp(-t\tilde{\Dirac}^2) \right] \in K_0\left( C^*(\tilde{M})\right) \; .
$$
It is independent of the metric on $M$, and of the parameter $t$. Compare this with what we would have gotten had we not lifted to $\tilde{M}$. We still have the heat operator $\exp(-t\Dirac^2)$, a smoothing operator over $M$. A smoothing operator over a compact manifold is compact, and a compact projection is essentially determined by its (graded) trace. The (graded) trace of the heat operator is nothing but the index of $\Dirac$. Thus, on the level of K-theory the heat operator on $M$ is the same as the index of $\Dirac$.\footnote{To be a bit more precise the heat operator over $M$ defines an element in the Roe algebra of $M$. However, every operator over the compact space $M$ has finite propagation, and it is locally compact if and only if it is compact. Hence $C^*(M)$ is simply the algebra of compact operators, whose zeroth K-theory is isomorphic to $\ZZ$ via the trace.} We think of $\mathrm{Ind}(\tilde{\Dirac})$ as an analogue to the index of $\Dirac$ that lives in a richer K-theory group, motivating the notation. We call it the \emph{coarse index}; the terminology will be motivated below. It can also be defined on odd-dimensional manifolds, then it is an element of $K_1(C^*(\tilde{M}))$. 

How does the coarse index help with the psc-question? Suppose $M$ admits a psc-metric. That metric lifts to a metric of uniformly positive scalar curvature (upsc) on $\tilde{M}$. Thus $\tilde{\Dirac}$ has a spectral gap around zero. Taking $t$ sufficiently large the function $\exp(-t\lambda^2)$ essentially vanishes on the spectrum of $\tilde{\Dirac}$, so that the class of $\exp(-t\tilde{\Dirac}^2)$ vanishes. (Of course this argument is not at all rigorous.) Thus $\mathrm{Ind}(\tilde{\Dirac})=0$ if $M$ admits a psc-metric. Consider the torus, so that $\tilde{M}=\RR^n$. One can show that $K_n(C^*(\RR^n))\cong\ZZ$, with $\mathrm{Ind}(\tilde{\Dirac})$ being a generator. In particular the coarse index does not vanish, so that the torus cannot admit a psc-metric. Thus the coarse index is able to detect information the index misses.  

In the construction of the coarse index we discarded the information that $\tilde{\Dirac}$ is actually a $\Gamma$-invariant operator, so that the heat operator becomes a $\Gamma$-invariant smoothing operators. As such it defines an element in the $\Gamma$-invariant Roe algebra $C^*(\tilde{M})^\Gamma$. This algebra turns out to be Morita-equivalent to the reduced group $C^*$-algebra $C^*_r\Gamma$. Thus there is also a \emph{$\Gamma$-equivariant index} (also called \emph{higher index})
$$
\mathrm{Ind}^\Gamma(\tilde{\Dirac}) \; \in \; K_*\left(C^*_r\Gamma \right) \; .
$$
Taking gradings and real structure into account we can produce from the heat operator an element of the real K-theory $KO_*(C^*_r\Gamma)$. This higher index class is called the \emph{Rosenberg index} \cite{Rosenberg1983}.\footnote{Rosenberg's approach to this index class was different than the one presented here. Instead of lifting to the universal cover, Rosenberg twists the Dirac operator on $M$ with a line bundle of Hilbert $C^*_r\Gamma$-modules called the Mi\v{s}\v{c}enko-Fomenko line bundle, see \cite{Mishchenko1980}. The theory of elliptic operators with coefficients in Hilbert $C^*$-module bundles is largely analogous to the usual elliptic theory. In particular, kernel and cokernel of this twisted Dirac operator are finite-dimensional Hilbert $C^*_r\Gamma$-modules, and their formal difference defines the Rosenberg index class in $KO_*(C^*_r\Gamma)$. This trades in the coarse geometry in the form of Roe algebras for Hilbert $C^*$-modules. Of course, for us the coarse geometry is a feature not a bug.} It is a powerful obstruction to the existence of a psc-metric. Indeed, the so-called Gromov-Lawson-Rosenberg conjecture posits that a closed spin manifold of dimension at least 5 admits a psc-metric if and only if its Rosenberg index vanishes. It is known to be false in general \cite{SchickGLR}, though has been proven for a number of groups $\Gamma$, in particular for $\Gamma=1$, i.e. simply-connected manifolds \cite{Stolz1992}. There is also a stable version of the Gromov-Lawson-Rosenberg conjecture stating the vanishing of the Rosenberg index of $M$ to be equivalent to the existence of a psc-metric on $M\times B^k$ for some $k$, where $B$ is the Bott manifold. It is currently open, and in fact implied by the Baum-Connes conjecture for $\Gamma$ \cite{Stolz2001}, which will be discussed below. See \cite{Stolz2023} for a review of these results and other things psc-related. See also \cite{Willett2019} for a book-length treatment of higher indices.
\subsection*{Coarse index theory}
It is the K-homological approach to index theory that index statements should be derived by applying index maps to the K-homology classes of elliptic operators. Therefore we ask: Is there a K-homological interpretation of the coarse index? The answer is affirmative. 

Let $M$ be a possibly non-compact manifold. The K-homology of $M$ can be identified with the K-theory of a certain quotient $C^*$-algebra: It holds that $K_*(M)$ is isomorphic to $K_{*+1}(D^*(M)/C^*(M))$. Here $C^*(M)$ is the Roe algebra already introduced above. It is the smallest $C^*$-algebra containing all locally compact finite-propagation operators. The algebra $D^*(M)$ is the smallest $C^*$-algebra containing all pesudo-local finite-propagation operators. It is called the \emph{structure algebra} of $M$. This isomorphism is a combination of a purely algebraic fact called \emph{Paschke duality} and the observation that every K-homology class can be represented by an operator with finite propagation. The boundary map in the long exact sequence in the K-theory of the pair $(D^*(M),C^*(M))$ provides us with a mapping
$$
\mathrm{Ind}: \; K_*(M) \cong K_{*+1}(D^*(M)/C^*(M)) \xlongrightarrow{\partial} K_*(C^*(M)) \; .
$$
We call it the \emph{coarse index map}. It maps the K-homology class of an elliptic operator over $M$ exactly to its coarse index as defined above.

What does the coarse index measure, and how should we understand its target, the K-theory of $C^*(M)$? Let us give Roe's initial intuition. We noted above that the coarse index of an elliptic operator $D$ defined via the heat operator is independent of the time parameter $t$, so that we are free to let $t$ be large. As time passes the distribution of heat on $M$ tends towards an equidistribution in the kernel of $D$, no matter the initial distribution. This trend towards equidistribution smears out local features of $M$. All that remains are large-scale tendencies of how heat propagates according to the dynamics of $D$.

The Roe algebra reflects this property of the coarse index to only remember large-scale information. Indeed $K_*(C^*(M))$ is an invariant of the large-scale geometry of $M$. One also speaks of the \emph{coarse} geometry of $M$, hence the term coarse index. Intuitively speaking the coarse geometry of a space is what remains of it when looked at from far away with squinted eyes. For example the coarse geometry of Euclidean space $\RR^n$ is the same as that of the lattice $\ZZ^n$; for the large-scale behavior the void between lattice points is irrelevant. The large-scale nature of $K_*(C^*(M))$ nicely matches with our intuition that heat flow smears out local features.

Return to the compact case for a moment. In this case $\exp(-tD^2)$ actually converges to the projection onto $\ker(D)$ in norm. As discussed earlier, this means that all K-theoretic information about the heat operator is already contained in the (graded) trace of that projection, i.e. the index of $D$. This is reflected by the Roe algebra: Over a compact space every operator has finite propagation, and an operator is locally compact if and only if it is compact. Thus $C^*(M)$ is isomorphic to the algebra of compact operators, and $K_0(C^*(M))$ to the integers. In this case the coarse index map $K_0(M)\to K_0(C^*(M))\cong\ZZ$ simply computes the index.\footnote{This is because the boundary map of the 6-term exact sequence of the pair $(D^*(M),C^*(M)=\fKH)$ factors through the boundary map $K_1(\fBH/\fKH)\to K_0(\fKH)=\ZZ$. It is a standard example in operator K-theory that this latter boundary map is an isomorphism given by computing the index.} In terms of heat flow we might say that the equidistribution towards which the heat flow converges is determined by its value at a point, and the (signed) number of different values is counted by the (graded) dimension of $\ker(D)$, i.e. its index.

Coarse index theory is the study of the coarse index map. Many properties of the index of elliptic operators on compact spaces generalize to the coarse setting. As we have seen above the coarse index of the Dirac operator obstructs the existence of a uniformly positive scalar curvature metric on $M$. There is a subtlety involved in this. Namely, the Roe algebra depends on the coarse structure of $M$, and this structure is metric-dependent. Different metrics can lead to inequivalent coarse structures. Thus, we should technically not speak of \emph{the} Roe algebra of $M$, but instead view the Roe algebra (and thus its K-theory) as dependent on a reference metric. Then, the coarse index of the Dirac operator obstructs the existence of a upsc-metric in the coarse equivalence class of a reference metric. We didn't encounter this problem when discussing coarse indices on a universal cover, because we only considered metrics arising as lifts of metrics on the compact base. Any two such metrics are bi-Lipschitz equivalent, and in particular coarsely equivalent.
\subsection*{The Baum-Connes conjecture}
Coarse index maps are closely related to assembly maps and the Baum-Connes conjecture. Though not immediately relevant to the contents of this thesis the Baum-Connes conjecture and the surrounding circle of ideas connects many different areas of mathematics and presented a guiding question for the last decades of mathematical research in analysis, geometry and group theory. Hence it is an important part of the larger context of our investigations. 

We begin with the coarse version, as it is closest to what we have discussed so far. Observe that the domain of the coarse index map
$$
\mathrm{Ind}: \; K_*(M) \longrightarrow K_*(C^*(M))
$$
depends on the topology of $M$ while the target depends only on the coarse geometry. To get a map between purely coarse-geometric objects one can replace the domain by a universal measure of the K-homologies within the coarse equivalence class of $M$. This object is called the \emph{coarse K-homology} of $M$. It can be realized for example as a colimit of K-homology over the coarse equivalence class of $M$. This yields a universal index map from the coarse K-homology of $M$ to the K-theory of its Roe algebra. This map is called the \emph{coarse assembly map}. The \emph{coarse Baum-Connes conjecture} asserts that the coarse assembly map is an isomorphism for every Riemannian manifold $M$, indeed for every coarse space. It is proved for a large class of spaces \cite{Yu2000}, but known to be false in general \cite{Higson2002}, though counterexamples to date are still somewhat convoluted, and finding more natural ones is a matter of ongoing research, see for example \cite{Schick2025}.

This discussion can also be carried out equivariantly. Suppose for simplicity that $M$ is a compact aspherical manifold, meaning that its universal cover is contractible. Then $M$ is a classifying space $B\Gamma$ for its fundamental group $\Gamma$, and the universal cover $\tilde{M}$ is the classifying space $E\Gamma$ for proper free $\Gamma$-actions. As discussed earlier a $\Gamma$-invariant elliptic operator over $\tilde{M}$ has an equivariant higher index in $K_*(C^*(\tilde{M})^\Gamma)=K_*(C^*_r\Gamma)$. This assignment gives rise to an equivariant index map from the $\Gamma$-equivariant K-homology $K_*^\Gamma(\tilde{M})$ to $K_*(C^*_r\Gamma)$. This equivariant K-homology turns out to be the same as the regular K-homology $K_*(M)=K_*(B\Gamma)$. Thus we obtain a map
$$
K_*(B\Gamma) \longrightarrow K_*(C^*_r\Gamma) \; .
$$ 
It can be defined for general discrete groups, not only those arising as the fundamental groups of closed aspherical manifolds. The map is \emph{almost} the Baum-Connes assembly map, but a slight modification is needed. Consider the classifying space $\underline{E}\Gamma$ for proper $\Gamma$-actions. The equivariant index map refines to a map
$$
\mu: \; K^\Gamma_*(\underline{E}\Gamma) \longrightarrow K_*\left(C^*_r\Gamma\right) \; .
$$
This is the \emph{Baum-Connes assembly map}. Note that by forgetting about freeness one obtains a canonical map
$$
K_*(B\Gamma)=K^\Gamma_*(E\Gamma) \longrightarrow K^\Gamma_*(\underline{E}\Gamma) \; ,
$$
whose composition with the Baum-Connes assembly map is the equivariant index map defined above. If $\Gamma$ is torsionfree, then this map is an isomorphism, so that the assembly map is a map $\mu: \, K_*(B\Gamma)\to K_*(C^*_r\Gamma)$.

The \emph{Baum-Connes conjecture} asserts that the Baum-Connes assembly map is an isomorphism for all groups $\Gamma$. It has been established for a large number of groups \cite{Higson2001},\cite{Kasparov2003}, \cite{Lafforgue2012}, \cite{Yu2000}, and no counterexamples are known. It contains multiple different important conjectures as corollaries: If the Baum-Connes conjecture holds for $\Gamma$, then so does the stable Gromov-Lawson-Rosenberg conjecture mentioned above, the the Novikov conjecture on the homotopy invariance of higher signaturs, and the Kadison-Kaplanski conjecture on the non-existence of non-trivial idempotents in $C^*_r\Gamma$. We refer to \cite{Aparicio2019} for a survey of the Baum-Connes conjecture, its consequences, and different methods of proving it.

Let us also mention that there is a version of the Baum-Connes assembly map with coefficients in a $\Gamma$-$C^*$-algebra $A$. Its domain is the the K-homology $K^\Gamma_*(\underline{E}\Gamma;A)$ with coefficients in $A$, and its target the K-theory $K_*(A\rtimes_r \Gamma)$ of the reduced crossed product of $\Gamma$ with $A$. Choosing $A=\CC$ with the trivial action returns the Baum-Connes assembly map without coefficients. There is a multitude of relations between the Baum-Connes conjecture and its variants. For example there is a method called the descent principle connecting the coarse assembly map to the Baum-Connes assembly map (see for example \cite[Chapter 12]{HigsonRoe2000}). Moreover, it is a result of Yu \cite{Yu1995} that a group $\Gamma$ satisfies the coarse Baum-Connes conjecture if and only if it satisfies the Baum-Connes conjecture with $\ell^\infty(\Gamma;\fK)$-coefficients, where $\fK$ denotes the algebra of compact operators. This last result was actually a motivation for the introduction of uniform K-homology, more on that below.
\subsection*{Index theory in the presence of boundaries}
So far we have only considered manifolds without boundary. In the presence of a boundary index theory acquires a somewhat different flavor. To appreciate why that is we first reconsider the boundaryless case. Let $M$ be a closed manifold and $D$ an elliptic operator over $M$. The operator $D$ may be viewed as an unbounded operator on $L^2(M)$ with domain given by the space of smooth sections $\Gamma(M)\subseteq L^2(M)$. This operator is not closed, but it has a canonical closed extension, namely to the domain $H^m(M)\subset L^2(M)$. It is a consequence of the Rellich-Kondrachov theorem that this canonical extension has a finite-dimensional kernel. Moreover, the Hilbert space adjoint $D^*$ is the canonical extension of the formal adjoint of $D$, hence its kernel is also finite-dimensional. We conclude that $\mathrm{index}(D)=\ker(D)-\ker(D^*)$ is a well-defined integer.

Now assume that $\olO$ is a compact manifold with boundary. Upon trying to replicate the above argument we run into two problems. The elliptic operator $D$ still defines an unbounded operator on $L^2(\olO)$, with domain given by the space $\Gamma_{cc}(\Omega)$ of smooth sections supported in the interior $\Omega$ of $\olO$. Again, this is not a closed operator. Here, the first problem appears: There is no canonical closed extension anymore, but instead many different ones. Indeed, each closed extension $D_e$ corresponds to a choice of boundary condition for $D$. Upon choosing a boundary condition we run into the second problem: Not all boundary conditions lead to domains that are (a subspace of) some Sobolev space, in fact the kernel need not be finite-dimensional. Even worse: Even if $\ker(D_e)$ is finite-dimensional, $\ker(D_e^*)$ need not be. Thus, we get a well-defined index only for so-called \emph{regular} boundary conditions for which $\ker(D_e)$ and $\ker(D_e^*)$ are finite-dimensional, which in practice means $\dom(D_e),\dom(D_e^*)\subseteq H^m(\olO)$. Importantly, the index of $D_e$ generally depends on the choice of boundary condition.\footnote{As an example the Euler characteristic operator over $\olO$ admits two local regular boundary conditions called the absolute and relative boundary conditions, whose indices are the absolute Euler characteristic $\chi(\olO)$ and the relative Euler characteristic $\chi(\olO,\dO)$ respectively. The two generally do not agree. See for example \cite{Gilkey1996}.}

Classically one has considered differential boundary conditions such as Dirichlet boundary conditions (the function vanishes on the boundary) or Neumann boundary conditions (the normal derivative vanishes on the boundary). More generally a \emph{local} boundary condition is one given by $Bu|_{\dO}=0$, where $B$ is a differential operator of order $<m$ over $\dO$. Indeed, one usually considers \emph{local elliptic} boundary conditions, where the symbol of $B$ satisfies a certain isomorphism assumption. This type of boundary condition is always regular, so that the index of the associated extension $D_B$ is well-defined. For such boundary conditions Atiyah and Bott proved the following index theorem:
\begin{thm}[Atiyah-Bott, \cite{Palais1965}, \cite{AtiyahBott1964}]
Let $\olO$ be a compact manifold with boundary. Let $D$ be an elliptic operator over $\olO$, and let $B$ define a local elliptic boundary condition for $D$. Then,
$$
\mathrm{index}(D_B) = \int_{\olO} \mathrm{AS}(D,B,\olO) \; .
$$ 
Here $\mathrm{AS}(D,B,\olO)$ is a universal characteristic class over $M$ that is built from the symbols of $D$ and $B$, as well as characteristic classes of $\olO$.
\end{thm}
The index theorem of Atiyah-Bott is very similar to the Atiyah-Singer theorem for closed manifolds. It expresses the index in terms of topological data; indeed the $\mathrm{D_B}$ depends on the boundary condition $B$ only through the topology of its symbol. Thus the Atiyah-Bott theorem provides a satisfying extension of the Atiyah-Singer theorem to local elliptic boundary conditions. There is a problem, however: Not all $D$ admit local elliptic boundary conditions. There is a topological obstruction whose non-vanishing forbids the existence of such boundary conditions. For the spin Dirac operator this obstruction usually does not vanish. Thus we cannot expect to be able to apply the Atiyah-Bott theorem to get obstructions to positive scalar curvature on manifolds with boundary.
\subsection*{The APS index theorem and secondary invariants}
It was an important insight of Atiyah, Patodi and Singer that first-order elliptic operators admit interesting \emph{non-local} regular boundary conditions. Assume that $\olO$ carries a metric that has product structure near the boundary. That means that around $\dO$ the space looks like a cylinder. Assume we are given a first-order elliptic operator $D$ on $M$ that also has product form near $\partial M$, meaning that near the boundary it has the form $D=\partial_t + D_{\dO}$, where $t$ denotes the height-variable of the cylinder and $D_{\dO}$ is a symmetric elliptic operator over $\dO$. The operator $D_{\dO}$ has real spectrum, and we can consider the projection $P$ onto the subspace of $L^2(\dO)$ corresponding to the non-negative part of its spectrum. Then, a function $u$ on $\olO$ is said to satisfy the \emph{APS-boundary condition} if $Pu|_{\dO}=0$. It is a regular boundary condition, and therefore has a well-defined index $\ind(D_{APS})$. This index is computed by the APS-theorem.
\begin{thm}[Atiyah-Patodi-Singer, \cite{APS1}]
Let $\olO$ be a compact manifold with boundary, and assume that $\olO$ has product structure near the boundary. Let $D$ be a first-order elliptic operator over $\olO$ that has product structure $D=\partial_t + D_{\dO}$ near the boundary. Then,
$$
\ind(D_{APS}) = \int_{\olO} \mathrm{AS}(D,\olO) - \frac{\eta(D_{\dO})+\dim(\ker(D_{\dO}))}{2} \; .
$$
Here $\mathrm{AS}(D,\olO)$ is the same characteristic class that appears in the Atiyah-Singer theorem, and $\eta(\dO)$ is the so-called $\eta$-invariant. 
\end{thm} 
The APS-index formula contains the same topological term as the Atiyah-Singer formula on a closed manifold. It also contains the dimension of the kernel of $D_{\dO}$; this is not a topological quantity. The newcomer is the $\eta$-invariant $\eta(D_{\dO})$. It can be defined for any self-adjoint elliptic operator $A$ over a closed manifold $\Sigma$. Such an $A$ has discrete spectrum of finite-multiplicity eigenvalues. In fact the series
\begin{equation}\label{EtaInvariant}
\eta(s) := \sum_{\lambda} \mathrm{sign}(\lambda)\cdot \mathrm{mult}(\lambda) \cdot |\lambda|^{-s}
\end{equation}
converges for all $s\in\CC$ with $\mathrm{Re}(s)$ sufficiently large. The function $\eta(s)$ admits a meromorphic continuation to the entire complex plane that turns out to have $s=0$ as a regular value. The $\eta$-invariant is defined to be $\eta(A):=\eta(0)$. As might be guessed by formally taking $s=0$ in \eqref{EtaInvariant} it is a measure of the spectral asymmetry of $A$ about $0$. 

The $\eta$-invariant is a rather sensitive measure of the spectral properties of $A$. For example, if $t\mapsto A_t$ is a homotopy of elliptic operators, then $\eta(A_1)-\eta(A_0)$ is equal to the so-called spectral flow, the signed number of eigenvalues that cross $0$ during the homotopy. Due to the presence of both $\dim(\ker(D_{\dO}))$ and $\eta(D_{\dO})$ in the APS-index formula, the APS-index is \emph{not} a topological quantity. This is in contrast to both the Atiyah-Singer and the Atiyah-Bott formula.

It is possible to produce more stable spectral invariants from the $\eta$-invariant. As mentioned above $\eta(A)$ changes along a homotopy of operators according to the spectral flow, an integer quantity. Thus $\eta(A)\in \RR/\ZZ$ is invariant under homotopies. There are also relative $\eta$-invariants constructed as follows: Let $f: \, \Sigma\to B\Gamma$ be a continuous map classifying a $\Gamma$-covering $\Sigma_f\to \Sigma$. Here $\Gamma$ is the fundamental group of $\Sigma$. Moreover, let $\alpha_1,\alpha_2: \, \Gamma \to U(k)$ be unitary representations of the same dimension. To the data $(f,\alpha_1,\alpha_2)$ there are associated bundles $E_j:= \Sigma_f \times_\Gamma \CC^k$, where $\Gamma$ acts on $\CC^k$ according to $\alpha_j$, $j=1,2$. Twisting $A$ by $E_j$ leads to the twisted $\eta$-invariants $\eta_{\alpha_j}(A):=\eta(A_{E_j})$. Define the $\rho$-invariant as\footnote{Different sign conventions are possible here. We stick to the conventions of the original work \cite{APS1}, but point out that a different convention is used in \cite{HigsonRoe2010}, which will be discussed later. Beyond this footnote we will not worry about signs.}
$$
\rho(A,f,\alpha_1,\alpha_2) := \frac{1}{2} \left( \eta_{\alpha_1}(A) - \eta_{\alpha_2}(A) \right) - \frac{1}{2} \left( \dim(\ker(A_{E_1})) - \dim(\ker(A_{E_2})) \right)
$$
Its reduction to $\RR/\ZZ$ is a homotopy invariant of $A$. Furthermore, the APS-index theorem implies\footnote{The difference between the integrals appearing in the APS-formula for $D_{E_j}$ vanishes. This is best seen via the heat kernel approach, which realizes these integrals as integrals over local data coming from the asymptotic expansion of the heat kernel of $D_{E_j}$. Locally $D_{E_1}$ and $D_{E_2}$ behave the same, because both arise from twisting $D$ with a flat vector bundle of the same dimension. Thus the integrands coming from the heat kernels agree. See \cite{APS2}.} that $\rho(A,f,\alpha_1,\alpha_2)$ is an integer if $\Sigma$ is the boundary of a manifold $\olO$, $A=D_{\dO}$ the adapted boundary operator of an elliptic operator $D$ over $M$, and the map $f$ extends to $M$. Thus $\rho(A,f,\alpha_1,\alpha_2) \in \RR/\ZZ$ turns out to be an invariant of only the $\Gamma$-bordism class of $\Sigma$. Atiyah-Patodi-Singer provide an $\RR/\ZZ$-valued index theorem computing $\rho(A,f,\alpha_1,\alpha_2)$ in terms of topological data associated to the $\alpha_j$ and the symbol of $A$.

There are also applications to the psc-question. For example, under the assumption that the \emph{maximal} Baum-Connes assembly map $K_*(B\Gamma)\to K_*(C^*_{\max}\Gamma)$ (more on this below) is an isomorphism Mathai \cite{Mathai1992} proved that if $\Sigma$ admits a psc-metric, then $\rho(\Dirac_\Sigma,f,\alpha_1,\alpha_2)=0$ for all $\alpha_1$, $\alpha_2$, where $f$ is taken to be the map classifying the universal cover of $\Sigma$. We stress that this result holds \emph{without} reducing to $\RR/\ZZ$.
\subsection*{Higher secondary invariants and mapping geometry to analysis}
Just as there are higher indices living in the K-theory of Roe algebras, there are also higher secondary invariants. Let $M$ be a closed spin manifold, and $\tilde{M}$ its universal cover. We again let $\tilde{\Dirac}$ denote the lift of the spin Dirac operator to $\tilde{M}$. Recall that there is a coarse index map $\mathrm{Ind}: \, K_*(\tilde{M})\to K_*(C^*(\tilde{M}))$ that arises as the boundary map in operator K-theory obtained from the isomorphism $K_*(\tilde{M})\cong K_{*+1}(D^*(\tilde{M})/C^*(\tilde{M}))$. The index map therefore fits into an exact sequence
$$
\cdots \longrightarrow K_{*+1}(D^*(\tilde{M})) \longrightarrow K_*(\tilde{M}) \xlongrightarrow{\mathrm{Ind}} K_*(C^*(\tilde{M})) \longrightarrow \cdots \; .
$$
Now suppose that $M$ carries a psc-metric, which we lift to a upsc-metric on $\tilde{M}$. Then $\mathrm{Ind}(\tilde{\Dirac})=0\in K_*(C^*(\tilde{M}))$. By exactness of the above sequence the class $[\tilde{\Dirac}]\in K_*(\tilde{M})$ has a pre-image in $K_{*+1}(D^*(\tilde{M}))$. In fact there exists a canonical pre-image obtained by manually lifting $[\tilde{\Dirac}] \in K_{*+1}(D^*(\tilde{M})/C^*(\tilde{M}))$ to an element of $\rho(\tilde{\Dirac})\in K_{*+1}(D^*(\tilde{M}))$, using the spectral gap provided by the upsc-assumption. This canonical lift $\rho(\tilde{\Dirac})$ is called the \emph{higher $\rho$-invariant}. Generally it depends on the choice of metric.

We defer an explanation of the connecion of the higher $\rho$-invariant to the $\rho$-invariant discussed before to discuss the geometric significance of the former. The higher $\rho$-invariants are measures of upsc-metrics up to upsc-bordism. Indeed, there is a more general framework of mapping geometry to analysis, of which higher indices and $\rho$-invariant are a part. Broadly speaking, mapping geometry to analysis means mapping an exact sequence of geometric content to one of analytic meaning, specifically the exact sequence involving K-homology, the K-theory of the Roe algebra and the index map, by means of higher invariants. This framework was originally devised by Higson and Roe in the context of surgery and higher signatures \cite{HigsonRoe2004a}, \cite{HigsonRoe2004b}, \cite{HigsonRoe2004c}. Piazza and Schick then carried out the same in the context of positive scalar curvature \cite{Piazza2014}. Let us describe that latter instance in more detail. The geometric part is \emph{positive-scalar-curvature sequence}
$$
\cdots \longrightarrow \mathrm{Pos}^{\mathrm{spin}}_n(B\Gamma) \longrightarrow \Omega^\mathrm{spin}_n(B\Gamma) \longrightarrow R^{\mathrm{spin}}_n(B\Gamma) \longrightarrow \mathrm{Pos}^{\mathrm{spin}}_{n-1}(B\Gamma) \longrightarrow \cdots
$$
introduced by Stolz \cite{Stolz1998} as an analogue to the surgery exact sequence in the psc-setting. Here, $\Gamma$ is a fixed discrete group, and:
\begin{itemize}
\item $\mathrm{Pos}^\mathrm{spin}_n(B\Gamma)$ consists closed $n$-dimensional spin-manifolds equipped with a psc-metric and continuous maps to $B\Gamma$, modulo spin-bordisms that are also equipped with a psc-metric and a continuous map to $B\Gamma$,
\item $\Omega^\mathrm{spin}_n(B\Gamma)$ is the spin-bordism group of closed $n$-dimensional spin manifolds with maps to $B\Gamma$, modulo spin bordisms with maps to $B\Gamma$,
\item $R^{\mathrm{spin}}_n(B\Gamma)$ consists of compact spin manifolds with possibly non-empty boundary equipped with a psc-metric on the boundary if existent and a map to $B\Gamma$, modulo a restricted notion of psc-bordism of the boundaries.
\end{itemize}
Note that in each of these groups the map $M\to B\Gamma$ gives rise to a $\Gamma$-covering of the manifold $M$ in question. The analytic counterpart is the \emph{analytic surgery exact sequence}
$$
\cdots \longrightarrow K_{n+1}\left(C^*_r\Gamma\right) \longrightarrow K_{n+1}\left(D^*(E\Gamma)^\Gamma\right) \longrightarrow K_n(B\Gamma) \xlongrightarrow{\mathrm{Ind}^\Gamma} K_n\left(C^*_r\Gamma\right) \longrightarrow \cdots \; ,
$$
which is the equivariant version of the exact sequence containing the coarse index map considered above. Piazza and Schick construct a commutative diagram \\
\centerline{\xymatrix{
\cdots \ar[r] & \mathrm{Pos}^\mathrm{spin}_{n}(B\Gamma) \ar[d] \ar[r] & \Omega_n^\mathrm{spin}(B\Gamma) \ar[d] \ar[r] & R^\mathrm{spin}_n(B\Gamma) \ar[r] \ar[d] & \cdots \\
\cdots \ar[r] & K_{n+1}\left(D^*(E\Gamma)^\Gamma \right) \ar[r] & K_n(B\Gamma) \ar[r] & K_n\left(C^*_r\Gamma \right) \ar[r] & \cdots
}}
from the positive-scalar-curvature sequence to the analytic surgery exact sequence. Let us describe the vertical maps. Recall that an element in $\mathrm{Pos}^\mathrm{spin}_{n+1}(B\Gamma)$ is represented by a spin-manifold $M$ with psc-metric and a map $M\to B\Gamma$. As mentioned above this map gives rise to a $\Gamma$-covering $\tilde{M}\to M$. Pulling back the spin-structure psc-metric to $\tilde{M}$ makes that manifold into spin-manifold with upsc-metric, where every part of the structure is $\Gamma$-invariant. Thus, from the $\Gamma$-invariant Dirac operator on $\tilde{M}$ we get a $\Gamma$-invariant higher $\rho$-invariant on $\tilde{M}$. Pushing it forward to $E\Gamma$ results in an element $\rho^\Gamma(\tilde{M})\in K_{n+1}(D^*(E\Gamma)^\Gamma)$. It is unchanged by psc-bordisms. Thus, the higher $\rho$-invariant gives rise to the map $\mathrm{Pos}^\mathrm{spin}_{n+1}(B\Gamma) \to K_{n+1}(D^*(E\Gamma)^\Gamma)$. We note that it is only interesting if $\Gamma$ has torsion. Otherwise the map $K_*(B\Gamma)\to K_*(C^*_r\Gamma)$ is the assembly map, which the Baum-Connes conjecture claims to be an isomorphism. By exactness it follows that $K_*(D^*(E\Gamma)^\Gamma)=0$.

The maps $\Omega_n^\mathrm{spin}(B\Gamma)\to K_n(B\Gamma)$ and $R^\mathrm{spin}_n(B\Gamma)\to K_n(C^*_r\Gamma)$ are defined similarly. For the former one lifts the Dirac operator to the cover coming from the map to $B\Gamma$, then pushes its K-homology class forward to $K^\Gamma_n(E\Gamma)=K_n(B\Gamma)$. Here, well-definedness under the bordism-relation follows from boundary-of-Dirac-is-Dirac formula, which we will describe below. The latter map is most easily described in the case of a cycle for $R^\mathrm{spin}_n(B\Gamma)$ that has no boundary. Then the associated element of $K_n(C^*_r\Gamma)$ is simply its $\Gamma$-invariant higher index. In this case bordism invariance is essentially bordism invariance of the index. In the presence of a psc-boundary one attaches cylindrical ends to obtain an APS-like index class. We will give the precise description in another context below. Bordism invariance is less straight-forward here.

Mapping geometry to analysis provides an intellectually satisfying framework connecting the exact sequences on both sides, and encoding various bordism invariance statements as well as index theorems. Apart from that it also has concrete applications. For example, the size of the image of the map $\mathrm{Pos}^\mathrm{spin}_*(B\Gamma)\to K_{*+1}(D^*(E\Gamma)^\Gamma)$ provides a lower bound for the size of $\mathrm{Pos}^\mathrm{spin}_*(B\Gamma)$. In this way one can estimate the size of the moduli space of psc-metrics, see for example \cite{Xie2021}. 

Let us now explain the connection of the higher $\rho$-invariant to secondary invariants. In \cite{HigsonRoe2010} Higson and Roe construct from a representation $\alpha_1: \, \Gamma \to U(k)$ a trace map $\mathrm{tr}_{\alpha_1}: \, K_0(C^*_{\max}\Gamma)\to \ZZ$. Given another representation $\alpha_2: \, \Gamma \to U(k)$ of the same dimension, it turns out that the difference $\mathrm{tr}_{\alpha_1}-\mathrm{tr}_{\alpha_2}$ factors through the K-theory of a certain variant of the structure algebra $D^*(E\Gamma)^\Gamma$. Indeed, there is an analogue 
$$
\cdots \longrightarrow K_{*+1}\left(C^*_{\max}\Gamma\right) \longrightarrow K_{*+1}\left(D^*_{\max}(B\Gamma)\right) \longrightarrow K_*(B\Gamma) \longrightarrow K_{*}\left(C^*_{\max}\Gamma\right) \longrightarrow\cdots
$$
to the analytic surgery exact sequence, which involves the maximal Baum-Connes assembly map $K_*(B\Gamma)\to K_*(C^*_{\max}\Gamma)$ instead. The algebra $D^*_{\max}(B\Gamma)$ is an adaptation of the usual structure algebra $D^*(E\Gamma)^\Gamma$. Higson and Roe also construct a relative trace map $\mathrm{tr}_{\alpha_1,\alpha_2}: \, K_0(D^*_{\max}(B\Gamma)) \to \RR$, as well as a relative $\RR/\ZZ$-valued index map $\mathrm{Ind}_{\alpha_1,\alpha_2}: \, K_1(B\Gamma)\to \RR/\ZZ$. These maps fit into a commutative diagram \\
\centerline{\xymatrix{
K_0(B\Gamma) \ar[r]\ar[d] & K_0(C^*_{\max}\Gamma) \ar[d]_{\mathrm{tr}_{\alpha_1}-\mathrm{tr}_{\alpha_2}}\ar[r] & K_0(D^*_{\max}(B\Gamma))\ar[d]_{\mathrm{tr}_{\alpha_1,\alpha_2}} \ar[r] & K_1(B\Gamma)\ar[d]_{\mathrm{Ind}_{\alpha_1,\alpha_2}} \ar[r] & K_1(C^*_{\max}\Gamma)\ar[d] \\
0 \ar[r] & \ZZ \ar[r] & \RR \ar[r] & \RR/\ZZ \ar[r] & 0
}}
If $M$ is a closed odd-dimensional spin manifold with fundamental group $\Gamma$ and positive scalar curvature, then the higher $\rho$-invariant $\rho(\tilde{\Dirac}_M) \in K_0(D^*(\tilde{M})^\Gamma)$ can be pushed forward along the classifying map $f$ of the $\Gamma$-bundle $\tilde{M}\to M$ to an element of $K_0(D^*_{\max}(B\Gamma))$. Application of the relative trace $\mathrm{tr}_{\alpha_1,\alpha_2}$ produces precisely the $\rho$-invariant $\rho(\Dirac_M,f,\alpha_1,\alpha_2)$ introduced earlier. In particular, if the maximal assembly map $K_*(B\Gamma)\to K_*(C^*_{\max}\Gamma)$ is an isomorphism, then $K_0(D^*_{\max}(B\Gamma))$ vanishes, and hence $\rho(\Dirac_M,f,\alpha_1,\alpha_2)=0$. This vanishing result was proved originally by Mathai \cite{Mathai1992}. 

This explains the connection between higher $\rho$-invariants and secondary invariants. It also illustrates a general principal by which one extracts numerical invariants from K-theoretical ones, namely by means of traces. In fact, traces are special cases of cocycles in cyclic cohomology. K-theory admits a noncommutative Chern character to cyclic homology, and the latter can be paired with a cyclic cocycles. In this way one obtains numerical invariants from K-theoretic objects -say an index class- by pairing with cyclic cocycles.
\subsection*{Boundary conditions and K-homology}
Let us return to index theory on manifolds with boundary. Our commitment to the K-homological perspective leads us look for an interpretation of the index theory in the presence of boundaries as well as secondary invariants via K-homology. The basis of this interpretation would be the establishment of K-homology classes associated to boundary conditions for an elliptic operator $D$ over a (compact) manifold $\olO$ with boundary. This leads to two related questions: In what K-homology groups should these classes live, and what boundary conditions should be considered?

The most direct guess is this: Since we wish to do index theory over $\olO$ we should look for classes in $K_*(\olO)$, so that the push-forward to a point can compute the index. First of all this means that we should restrict to regular boundary conditions, since only these have an index. Moreover, it is necessary to assume that the boundary condition is local in a suitable way, otherwise the resulting operator $D_e(D_e^*D_e+I)^{-1/2}$ will generally not be pseudo-local. Intuitively this is because solutions to non-local boundary value problems can propagate arbitrarily quickly.\footnote{See Section 5.5 below for an elementary example.} Thus we get K-homology classes from local elliptic boundary conditions as appear in the Atiyah-Bott theorem, but not for the APS-boundary condition. It seems like a feasible goal to give a proof of the Atiyah-Bott theorem using K-homology, but the author knows of no such proof in the literature. Let us mention here that there is a generalization due to Boutet de Monvel of Atiyah-Bott to suitably elliptic boundary conditions for elliptic pseudo-differential operators in a certain algebra of such boundary value problems \cite{Monvel1971}. While there are proofs of Boutet de Monvel's index theorem using operator-algebraic techniques \cite{Melo2006}, the author again knows of no proof appealing to K-homology directly.

As discussed earlier not all elliptic operators admit local regular boundary conditions. For generic elliptic operators we are therefore forced to retreat to a weaker setting. Over the interior $\Omega$ all extensions of $D$ behave in the same way. Indeed, an extension corresponds to a boundary condition, which is not seen over the interior alone. In fact, over $\Omega$ elliptic operators behave just as regularly as they do over closed manifolds. This is known as \emph{interior regularity}. As a consequence any closed extension defines a Fredholm module over $\Omega$. The resulting K-homology class turns out to be independent of the choice of boundary condition. 

There is a third possibility. Consider a local boundary condition $D_e$ that is \emph{semi-regular}, meaning that $\dom(D_e)\subseteq H^m(\olO)$, but $\dom(D_e^*)\subseteq H^m(\olO)$ is not required. Then, it can still be shown that $D_e(D_e^*D_e+I)^{-1/2}$ is pseudo-local over all of $\olO$, while the remaining assumptions for a Fredholm module generally only hold on the interior $\Omega$. Such an operator is called a \emph{relative Fredholm module} over $(\olO,\dO)$. Relative Fredholm modules define classes in relative K-homology $K_*(\olO,\dO)$. There is an excision isomorphism $K_*(\olO,\dO)\cong K_*(\Omega)$, which identifies the relative class of a local semi-regular boundary condition with its absolute class over $\Omega$. Surprisingly, the relative K-homology class is therefore also independent of the choice of local semi-regular boundary condition, if such boundary conditions exist. As opposed to local regular boundary conditions they always do, an example is the minimal boundary coundition $u|_{\dO}=0$.

From an index-theoretic viewpoint it may seem unsettling that we can obtain K-homology classes from non-regular boundary conditions, and that these are independent of the choice of boundary conditions. This is resolved by the fact that the K-homology groups $K_*(\Omega)$ and $K_*(\olO,\dO)$ actually carry no index maps. Indeed, on K-homology index maps come from pushing forward to a point, but the constant map $\Omega\to \mathrm{pt}$ is not proper, and therefore induces no push-forward map on K-homology. Indeed, Fredholm modules over non-compact spaces need not be Fredholm operators.

One may then rightfully ask if there is any index-theoretic significance to these classes in $K_*(\olO,\dO)\cong K_*(\Omega)$. There is, and it comes from the boundary map on K-homology. The long exact sequence in the K-homology of the pair $(\olO,\dO)$ reads as
$$
\cdots \longrightarrow K_*(\olO) \longrightarrow K_*(\olO,\dO) \xlongrightarrow{\partial} K_{*-1}(\dO) \longrightarrow K_{*-1}(\olO) \longrightarrow \cdots \; .
$$
Assuming for simplicity that $D=\partial_t+ D_{\dO}$ has product form, with $D_{\dO}$ a differential operator over $\dO$, the famed \emph{boundary-of-Dirac-is-Dirac formula} states that the K-homological boundary map takes the relative class $[D]\in K_*(\olO,\dO)$ of an elliptic operator to the absolute class $[D_{\dO}]\in K_{*-1}(\dO)$, i.e.
$$
\partial [D]= [D_{\dO}] \; .
$$
Among the consequences of this formula is the bordism invariance of the index. Indeed, pushing forward to a point rewards us with the following commutative diagram: \\
\centerline{ \xymatrix{
K_*(\olO,\dO) \ar[r]^\partial & K_{*-1}(\dO) \ar[r] \ar[d]_{\mathrm{index}} & K_{*-1}(\olO) \ar[d]^{\mathrm{index}} \\
& \ZZ \ar@{=}[r] & \ZZ
}}
From exactness of the top row we conclude that $\ind(D_{\dO})=\ind(\partial[D])=0$. This is the bordism invariance of the index: The index of a boundary vanishes. Bordism invariance is an important result, and the key to the bordism proof of the Atiyah-Singer index theorem.

We also note that by exactness the kernel of the boundary map is the image of the restriction map $K_*(\olO) \to K_*(\olO,\dO)$. Thus, if $D$ admits a local regular boundary condition, so that it actually defines a class in $K_*(\olO)$, then the associated relative class lies in the kernel of the boundary map. We conclude that $\partial [D]=[D_{\dO}]\in K_{*-1}(\partial \Omega)$ obstructs the existence of local regular boundary conditions for $D$. In the case that $D$ is a Dirac operator, $\partial [D]$ is the class of the Dirac operator on the boundary, which actually generates $K_{*-1}(\dO)$. This proves the earlier claim that Dirac operators do not admit local elliptic (or more generally regular) boundary conditions.

It is a feature of the K-homological approach to index theory that non-compact spaces can usually be treated in the same way as compact spaces. The existence of relative K-homology classes proved by Baum, Douglas and Taylor presents an exception to this. Their proof in \cite{BDT1989} uses compactness (as opposed to local compactness) of resolvents in a very direct way via an appeal to a discrete spectrum of finite multiplicity eigenvalues.\footnote{\cite[Exercises 10.9.8-10.9.13]{HigsonRoe2000}, sketch a somewhat different argument than \cite{BDT1989}, which also does not generalize to the non-compact case straight-forwardly. There one would have to investigate the mapping properties of operators of the form $\mathrm{sign}(D_e)$. If $D_e$ has a discrete spectrum of finite-multiplicity eigenvalues, it is a unitary up to a compact perturbation. Otherwise one needs some machinery to deal with it. This is precisely what the formal operator calculus we employ in this thesis is designed to do.} Hence one cannot get away by simply replacing compactness with local compactness everywhere. To the author's knowledge the only work treating the existence of relative classes in the non-compact case is a very recent article by Fries \cite{Fries2025}, where it is proved using a localization principle that the minimal boundary condition for certain elliptic differential operators define relative K-homology classes over open subset of possibly non-compact manifolds.\footnote{The statement in \cite{Fries2025} is given for differential operators that are elliptic in one of a number of elliptic calculi. One of those admissible calculi is that of properly supported pseudo-differential operators with bounded symbols on manifolds of bounded geometry, a situation that is rather close to the one considered in this thesis. It should be said that even in the bounded-geometry case the results of \cite{Fries2025} are most likely not restricted to this particular calculus, they should apply equally well to the calculi used in \cite{Engel2018} or \cite{Taylor2008}. However, a systematic treatment of the non-compact setting was not the main focus of \cite{Fries2025}, but instead the extension of the work of Baum, Douglas and Taylor to higher-order operators.} Of course, earlier statements of the existence of relative classes in the non-compact setting in the literature could simply have eluded the author's search, or it could have previously existed only as a folk theorem.

Let us also mention a relative analogue of equivariant (or higher) indices in the presence of boundaries. In \cite{Chang2015} Chang, Weinberger and Yu construct a relative index map 
$$
\mathrm{Ind}^\mathrm{rel}: \; K_*(\olO,\partial \Omega) \longrightarrow K_* \left( C^*_{\max}(\pi_1(\olO),\pi_1(\dO))\right) 
$$
from the relative K-homology of a compact manifold $\olO$ with boundary to the K-theory of the relative maximal group $C^*$-algebra. They prove that if $\olO$ admits a psc-metric with product structure on the boundary, then $\mathrm{Ind}^\mathrm{rel}(\Dirac_\Omega)=0$. They propose that their relative index (or rather its refinement to Real K-homology/K-theory) should be thought of as a relative version of the Rosenberg index, and the relative index map on the relative K-homology of the classifying spaces as a relative version of the Baum-Connes assembly map. Chang, Weinberger and Yu use localization algebra techniques for their construction, but other approaches have since been considered. Deeley and Goffeng use geometric K-homology to define and study the relative index \cite{Deeley2017}, see also \cite{Deeley2013}, and Kubota uses relative Mi\v{s}\v{c}enko-Fomenko line bundles as well as techniques from KK-theory \cite{Kubota2020}. Lastly, we mention also a work by Schick and Seyedhosseini \cite{Schick2021} working largely in the localization algebra framework, providing a comparatively more elementary approach and addressing a gap present in \cite{Chang2015}.
\subsection*{Uniform K-homology and uniform Roe-algebras}
In the context of the Baum-Connes conjecture with coefficients we mentioned a result by Yu \cite{Yu1995} stating the coarse Baum-Connes conjecture for a group $\Gamma$ to be equivalent to the Baum-Connes conjecture for $\Gamma$ with $\ell^\infty(\Gamma,\fK)$-coefficients. Part of that equivalence is an isomorphism between the Roe algebra $C^*(\Gamma)$ and the crossed product $\ell^\infty(\Gamma,\fK)\rtimes_r \Gamma$. Associated to certain metric spaces like $\Gamma$ there is also another $C^*$-algebra called the \emph{uniform Roe algebra} $C^*_u(\Gamma)$. It too is a measure of the coarse equivalence class of $\Gamma$. It is more rigid than the Roe algebra, but more flexible than the equivariant Roe algebra. Indeed, one always has comparison maps 
$$
K_*\left(C^*_r\Gamma\right) = K_*\left( C^*(\Gamma)^\Gamma\right) \longrightarrow K_*\left(C^*_u(\Gamma)\right) \longrightarrow K_*\left( C^*(\Gamma) \right)
$$ 
It is a feature of the uniform Roe algebra that its K-theory tends to be much richer than those of both the group $C^*$-algebra and the Roe algebra. Intuitively speaking, it is more flexible than the former and thus allows for a richer K-theory, whereas the latter is so flexible that its K-theory is simple again.\footnote{For a simple example to illustrate this point is $\Gamma=\ZZ$. Then $C^*_r\Gamma= C(S^1)$ via Fourier transform, so that $K_*(C^*_r\Gamma)=\ZZ$ in both even and odd degrees. On the other hand we have already mentioned above that $K_*(C^*(\ZZ))=K_*(C^*(\RR))$ vanishes in even degrees, and is given by $\ZZ$ in odd degrees (generated by the coarse index of the Dirac operator). In contrast examination of the Pimsner-Voiculescu sequence (see for example \cite{Blackadar1998}) shows the K-theory of $C^*_u(\ZZ)\cong \ell^\infty(\ZZ)\rtimes_r\ZZ$ to be an infinitely generated quotient of $\ell^\infty_\ZZ(\ZZ)$ in even degrees, though it is still given by $\ZZ$ in odd degrees. See Example \ref{ExUniRoeIntegers}. Analogous remarks apply to $\Gamma=\ZZ^n$.}

For a group $\Gamma$ there is an isomorphism between $C^*_u(\Gamma)$ and the crossed product $\ell^\infty(\Gamma)\rtimes_r \Gamma$. Therefore it seems natural in light of Yu's result to ask if there is a uniform analogue to the coarse assembly map, such that the Baum-Connes conjecture with $\ell^\infty(\Gamma)$-coefficients is equivalent to this map being an isomorphism. This question was investigated by \v{S}pakula, who introduced a uniform version of analytic K-homology. This uniform K-homology theory admits an index map to the K-theory of the uniform Roe algebra. Analogously to the construction of the coarse assembly map from the coarse index map one thus obtains a uniform coarse assembly map. For a group $\Gamma$ this map is an isomorphism if and only if the Baum-Connes conjecture with $\ell^\infty(\Gamma)$-coefficients holds for $\Gamma$. This was proved for torsion-free groups by \v{S}pakula \cite{Spakula2009}, with the general case being provided by Engel \cite{Engel2019}.

Uniform K-homology $K^{u}_*(X)$ is defined for metric spaces $X$. It retains information about the metric geometry of the space, as opposed to analytic K-homology, which is a purely topological invariant. It is based on a refinement of the concept of Fredholm modules requiring the compactness of the relevant operators to be uniform in a sense dictated by the metric structure. While originally introduced in the context of uniform Roe algebras and the Baum-Connes conjecture it is also relevant to index theory: If $M$ is a manifold of bounded geometry, then elliptic operators over $M$ define classes in uniform K-homology, provided their ellipticity is uniform in a suitable sense. This was proved for Dirac-type operators by \v{S}pakula \cite{Spakula2009}, and for general uniformly elliptic pseudo-differential operators by Engel \cite{Engel2018}. Therefore, such a uniformly elliptic operator $D$ admits a \emph{uniform coarse index} 
$$
\Ind(D) \in K_*(C^*_u(M))
$$
refining the coarse index, in the sense that the inclusion map $K_*(C^*_u(M))\to K_*(C^*(M))$ takes $\Ind(D)$ to the coarse index $\mathrm{Ind}(D)$. Because the K-theory of $C^*_u(M)$ is generally richer than that of $C^*(M)$ we expect the uniform coarse index to contain more information than its coarse counterpart. On the other hand it is only defined under additional uniformness conditions on both the elliptic operator and the underlying space.

If $M$ is a spin manifold of bounded geometry that has uniformly positive scalar curvature, then the uniform coarse index $\Ind(\Dirac)$ of its Dirac index vanishes. Thus the uniform coarse index provides a finer obstruction to a more restricted version of uniformly positive scalar curvature. We must note however that if $M$ is the universal cover of a closed spin manifold with fundamental group $\Gamma$, then the uniform coarse index is the image of the equivariant index $\mathrm{Ind}^\Gamma(\Dirac)$ under the inclusion map $K_*(C^*(M)^\Gamma)\to K_*(C^*_u(M))$. Thus, in this case $\Ind(\Dirac)$ cannot detect non-existence of uniformly positive scalar curvature that is not already detected by the equivariant index. It is however possible to twist the equivariant Dirac operator in a non-equivariant way that still retains uniformness, and the resulting uniform coarse index is an obstruction to positive scalar curvature outside the reach of equivariant index theory. However, beyond some discussion in the outlook to this thesis, this direction will not be investigated here.

In any case uniform coarse indices could provide helpful tools in the study of geometry questions about spaces of bounded geometry. The K-homological approach to index theory thus dictates that one should try to understand the uniform K-homology classes of uniformly elliptic operators. This presupposes understanding of uniform K-homology itself. As an adaptation of analytic K-homology it shares the basic structure of its analytic relative. \v{S}pakula proved Paschke duality and the existence of a Mayer-Vietoris sequence \cite{Spakula2009}, \cite{Spakula2010}, and Engel constructed the Kasparov product on uniform K-homology and used it to derive homotopy invariance \cite{Engel2019}. Moreover, Engel proved Poincaré duality for uniform K-homology \cite{Engel2019}, and constructed Chern characters on uniform K-homology and uniform K-theory and used them to derive index formulas \cite{Engel2015}. But we have seen above that in the K-homological setting certain properties of the index like bordism invariance are best deduced from properties of K-homology classes of elliptic operators on manifolds with boundaries. Moreover, while the uniform coarse index as an obstruction to uniformly positive scalar curvature has been discussed in \cite{Engel2014}, higher secondary invariants have not yet been considered in this setting. 
\subsection*{Results of this thesis}
Let us give an overview of the content and results of this thesis. In an attempt to achieve a clearer and less cluttered overview we do not state all results at the level the generality in which they will be proved in the main part. Moreover, we do not always list the precise technical assumptions of our results, though we indicate that they are present. The precise statements are of course given in the main part of the text.
\subsubsection*{Uniform K-homology classes in the presence of boundaries}
In broad strokes this thesis aims to study uniformly elliptic operators on manifolds with boundary and bounded geometry, as well as their indices and secondary invariants, through the lens of uniform K-homology. The conceptual starting point is the generalization of the results of Baum-Douglas-Taylor regarding K-homology classes associated to an elliptic operator on a manifold with boundary. Concretely, we prove the following.
\begin{thm}\label{IntroThmClasses}
Let $\olO$ be a manifold with boundary and bounded geometry. Let $D$ be a uniformly elliptic differential operator over $D$, and let $D_e$ be a boundary condition for $D$. 
\begin{itemize}
\item[(i)] $D_e$ defines a class in the uniform K-homology $K_*^{u}(\Omega)$ of the interior $\Omega=\olO\setminus \partial \Omega$.
\item[(ii)] If $D_e$ is local and semi-regular, then it defines a class in the relative uniform K-homology $K_*^{u}(\olO,\partial\Omega)$.
\item[(iii)] If $D_e$ is local, self-adjoint and regular, then it defines a class in the uniform K-homology $K^{u}_*(\olO)$ of $\olO$.
\end{itemize}
The classes in (i) and (ii) are independent of the choice of boundary condition, and depend only on the principal symbol of $D$.
\end{thm}
We will say more about the regularity conditions below. The relative uniform K-homology groups appearing in item (ii) have not previously been considered, though their definition is a straight-forward adaption of the relative analytic K-homology groups. There is a natural excision isomorphism $K^{u}_*(\olO,\partial\Omega) \cong K^{u}_*(\Omega)$ that identifies the classes of items (i) and (ii) in Theorem \ref{IntroThmClasses}. Moreover, there is a long exact sequence relating the relative uniform K-homology of $(\olO,\partial \Omega)$ to the uniform K-homology of $\olO$ and $\partial\Omega$. In the compact setting it was a central result of Baum-Douglas-Taylor that the boundary map on analytic K-homology maps the relative class of the Dirac operator to the class of the Dirac operator on the boundary. We generalize this boundary-of-Dirac-is-Dirac formula to our setting.
\begin{thm}
Let $\olO$ be an $(n+1)$-dimensional spin manifold with boundary and bounded geometry which has product structure near the boundary. Then,
$$
\partial [\Dirac_\Omega] = (-1)^{n+1} [\Dirac_{\partial\Omega}] \; \in \; K^{u}_n(\partial \Omega) \; ,
$$
where $\partial : \, K^{u}_{*+1}(\olO,\partial \Omega) \to K^{u}_*(\partial\Omega)$ denotes the boundary map on uniform K-homology.\footnote{The sign in our boundary-of-Dirac-is-Dirac formula is an artifact of the construction of the boundary map used here. See Remark \ref{SignBdryDirac}.}
\end{thm}
In the compact case boundary-of-Dirac-is-Dirac implied the bordism invariance of the index. This is also true in uniform K-homology. However, some elaboration is required. First of all, if $\Sigma$ is a spin manifold of bounded geometry\footnote{Whenever we speak of a \emph{manifold of bounded geometry}, this means that this manifold has no boundary. Bounded-geometry manifolds with boundary will be referred to explicitly as \emph{manifolds with boundary and bounded geometry}.}, then its Dirac operator has a uniform coarse index $\Ind(\Dirac_M) \in K_*(C^*_u(\Sigma))$. If $\Sigma'$ is another such manifold, then a priori it makes no sense to ask if $\Ind(\Dirac_\Sigma)$ and $\Ind(\Dirac_{\Sigma'})$ are equal, since they live in different K-theory groups. This is true in principle even if $\Sigma$ and $\Sigma'$ are connected by a bordism $\olO$. If, however, $\olO$ has finite width, meaning that each point has uniformly bounded distance to both $\Sigma$ and $\sigma'$, then the inclusions $\Sigma\hookrightarrow \olO$ and $\Sigma'\hookrightarrow \olO$ are uniform coarse equivalences. The K-theory of the uniform Re algebras is invariant under such equivalences, meaning that there is a canonical identification $K_*(C^*_u(\Sigma))\cong K_*(C^*_u(\olO))\cong K_*(C^*_u(\Sigma'))$. Then it makes sense to ask whether $\Ind(\Dirac_\Sigma)$ and $\Ind(\Dirac_{\Sigma'})$ agree under this identification. An adaptation of the argument given in the compact case above then shows that they are.
\begin{cor}
Let $\Sigma$ and $\Sigma'$ be spin-manifolds of bounded geometry that are connected by a spin-bordism of bounded geometry and finite width. Then,
$$
\Ind(\Dirac_\Sigma) = \Ind(\Dirac_{\Sigma'}) \; .
$$
\end{cor}
Note also that some sort of finite-width requirement for bordisms is necessary for an interesting bordism theory, because otherwise any $\Sigma$ is null-bordant via the infinite-width bordism $\olO=\RR_{\geq 0}\times\Sigma$.

The uniform coarse index $\Ind(\Dirac_\Sigma)$ obstructs the existence of a metric of bounded geometry with uniformly positive scalar curvature \cite{Engel2014}. The bordism invariance of the uniform coarse index then implies that a spin manifold that is bordant to one with non-vanishing index cannot carry such a metric.
\begin{cor}
Let $\Sigma$ be a spin manifold of bounded geometry. If $\Sigma$ is connected to a spin manifold $\Sigma'$ via a spin-bordism of bounded geometry and finite width such that
$$
\Ind(\Dirac_{\Sigma'}) \neq 0 \; ,
$$ 
then $\Sigma$ does not admit a metric of bounded geometry and uniformly positive scalar curvature in its uniform coarse equivalence class.
\end{cor}
\subsubsection*{Regularity and rough domains}
We have deferred further discussion of (semi-)regularity to this point. In the compact case this was defined via the inclusion of $\dom(D_e)$ in the Sobolev space $H^m(\olO)$ ($m$ being the order of $D$). The same is done here. Sobolev spaces on manifolds with boundary and bounded geometry were introduced by Schick \cite{Schick1998}. They behave mostly analogously to Sobolev spaces on compact manifolds, so that the nice behavior of Sobolev spaces leads to nice behavior of regular boundary condition.

We actually prove Theorem \ref{IntroThmClasses} in greater generality than is stated here. Instead of taking $\Omega$ to be a manifold with boundary and bounded geometry we allow it to be an arbitrary open subset of a manifold $M$ with bounded geometry (without boundary). To discuss regularity of boundary conditions in this more general setting we first need a sufficiently strong theory of Sobolev spaces on such domains. This is developed in Chapter 4. Again, the behavior of these Sobolev spaces is analogous to that of domains in Euclidean space or compact manifolds. Since we want to cover domains with arbitrarily rough boundaries, the theory gets somewhat technical in places, though it is subtle even for smooth boundaries. 

On sufficiently rough domains mere regularity will actually not suffice to prove the existence of uniform K-homology classes. Instead we need to assume a stronger version that -roughly speaking- guarantees regularity in intermediate degrees as well. The reason for this will become clear in a second. On sufficiently smooth domains (in particular manifolds with boundary and bounded geometry) regularity in the above sense will automatically imply this stronger form of regularity, but this is not generally the case. This actually poses the question of existence: Does there always exist a boundary condition satisfying the necessary conditions for Theorem \ref{IntroThmClasses} (ii)? The answer is affirmative; the minimal boundary condition always does the job. This recovers the analysis of Fries \cite{Fries2025}, who proved that the minimal boundary condition for a uniformly elliptic operator over any subset of a manifold with bounded geometry defines a class in (non-uniform) relative K-homology.

To prove Theorem \ref{IntroThmClasses} the uniform compactness of various operators constructed from $D_e$ needs to be established. We will actually prove the stronger statement that these operators are uniformly finitely summable. Uniform summability comes from a quantitative strengthening of the Rellich-Kondrachov compactness theorem. Establishing sufficient degrees of summability requires the use of Sobolev spaces of fractional degrees. This is due to the fact that we treat more general functions of $D_e$ by reduction to resolvents via integral representations. Making these integrals converge in operator norm often requires us to spend an infinitessimal amount of Sobolev degree. Thus integer-degree Sobolev spaces are not sufficient for our purposes. Fractional degrees are more subtle in many respects, but the complex interpolation method introduced by Calderón \cite{Calderon1964} provides a framework for the transfer from integer degrees to intermediate ones, hence we will employ it often. Our main objectives in our investigation of Sobolev spaces on domains are therefore uniform compactness and interpolation behavior. 

Once these facts on Sobolev spaces are established the actual calculations leading to uniform K-homology classes are organized around a formal operator calculus intended to mimic the pseudo-differential operator calculus in the presence of boundaries. The use of such an operator calculus is well-established in the analysis of boundary conditions, and its use in index theory goes back at least to Connes-Moscovici \cite{Connes1995}. It was also used by Fries in \cite{Fries2025} to extend the Baum-Douglas-Taylor result to higher-order differential operators.
\subsubsection*{Relative uniform K-homology}
The fundamental properties of relative uniform K-homology -excision, long exact sequence, suspension- are established by arguments that are lengthy but largely analogous to those for analytic K-homology. Certain boundedness assumptions on the geometry of the underlying spaces are required, however. These fundamental properties of relative uniform K-homology lay the abstract groundwork on which much of the calculus of the uniform K-homology classes of uniformly elliptic operators is built.

An important prerequisite to the establishment of these properties is a result called \emph{Paschke duality}. It relates K-homology to the K-theory of certain quotient $C^*$-algebra. It is well-established in the case of analytic K-homology (see for example \cite{HigsonRoe2000}), and was proved for non-relative uniform K-homology by \v{S}pakula \cite{Spakula2009}. For relative uniform K-homology it reads as follows.
\begin{thm}[Relative uniform Paschke duality]
Let $X$ be a metric space and $Z\subseteq X$ a closed subspace. Then, there is a natural isomorphism 
$$
K^{u}_*(X,Z) \cong K_{*+1}\left(\fD^{u}(X)/\fD^{u}(X,Z)\right) \; . 
$$
\end{thm}
Here $\fD^{u}(X)$ is the $C^*$-algebra generated by all uniformly pseudolocal operators over $X$, and $\fD^{u}(X,Z)$ the ideal in $\fD^{u}(X)$ containing all operators that are additionally uniformly locally compact over $X\setminus Z$. The algebras $\fD^{u}(X)$ and $\fD^{u}(X,Z)$ are called the \emph{uniform dual algebra} and \emph{relative uniform dual algebra} respectively.

Paschke duality allows one to apply the machinery of operator K-theory to (relative uniform) K-homology. Indeed, many of the fundamental results on relative uniform K-homology will be derived from a statement about (relative) uniform dual algebras. We prove the following structural properties of relative uniform K-homology.
\begin{thm}
Let $X$ be a metric space, and $Z \subseteq X$ a closed subspace. Under certain boundedness conditions on the geometry of $(X,Z)$ the following holds:
\begin{itemize}
\item[(i)] There is a natural excision isomorphism
$$
K_*^{u}(X,Z) \cong K_*^{u}(X\setminus Z) \; .
$$
\item[(ii)] There is a natural long exact sequence
$$
\cdots \longrightarrow K^{u}_*(Z) \longrightarrow K^{u}_*(X) \longrightarrow K^{u}_*(X,Z) \xlongrightarrow{\partial} K^{u}_{*-1}(Z) \longrightarrow \cdots \; .
$$
\item[(iii)] The boundary map of the pair $([0,1)\times X,X)$ provides a natural suspension isomorphism
$$
K^{u}_{*+1}\left( (0,1)\times X \right)\cong K^{u}_{*+1}([0,1)\times X,X) \longrightarrow K_*^{u}(X) \; ,
$$
whose inverse is the Kasparov product $K^{u}_*(X) \to K_{*+1}^{u}((0,1)\times X), \, \alpha\mapsto d\times \alpha$ with the generator $d\in K^{u}_{-1}((0,1))\cong\ZZ$.
\end{itemize}
\end{thm}
We have already mentioned above that the excision isomorphism connects the relative class of a boundary condition for a uniformly elliptic operator to the non-relative class over the interior. For the latter we will be able to show that the uniform K-homology class is independent of the choice of boundary condition, and depends only on the principal symbol of the operator. The excision isomorphism thus implies the same for the relative class. Moreover, the long exact sequence is a prerequisite to even formulate the boundary-of-Dirac-is-Dirac statement, and the suspension isomorphism is a key ingredient to its proof. Thus, these formal properties of relative uniform K-homology will have concrete consequences for the classes of uniformly elliptic boundary conditions later on.
\subsubsection*{Index maps on relative uniform K-homology}
We explore two possible approaches to endow relative uniform K-homology with an index map. The first is formally satisfying, but unfortunately useless for many practical purposes. The second carries interesting secondary information, and turns out to be a refinement of the first. As a trade-off it is harder to grasp from a structural standpoint. Due to its fascinating secondary nature this latter index will be the focus of most geometric applications discussed below.

The target of the first index map is the K-theory of the \emph{relative uniform Roe algebra}. Under mild geometric assumptions this algebra fits into a long exact sequence
$$
\cdots \longrightarrow K_*\left( C^*_u(Z) \right) \longrightarrow K_*\left( C^*_u(X) \right) \longrightarrow K_*\left( C^*_u(X,Z) \right) \xlongrightarrow{\partial} K_{*-1}\left( C^*_u(Z) \right) \longrightarrow \cdots \; .
$$
We construct a uniform coarse index map $\Ind: \, K_*^{u}(X,Z)\to K_*(C^*_u(X,Z))$, which then fits into a commutative diagram\\
\centerline{\xymatrix{
K_*^{u}(Z) \ar[r] \ar[d]_{\Ind} & K^{u}_*(X) \ar[r] \ar[d]_{\Ind} & K^{u}_*(X,Z) \ar[r]^\partial \ar[d]_{\Ind} & K_{*-1}^{u}(Z) \ar[d]_{\Ind} \\
K_*\left(C^*_u(Z)\right) \ar[r] & K_*\left(C^*_u(X)\right) \ar[r] & K_*\left(C^*_u(X,Z)\right) \ar[r]_\partial & K_{*-1}\left(C^*_u(Z)\right)
}}
This is satisfying from a formal standpoint, as it exhibits the uniform coarse index map as a transformation between homology theories, (relative) uniform K-homology on one side, K-theory of (relative) uniform Roe algebras on the other, that respects long exact sequences. Unfortunately $K_*\left(C^*_u(X,Z)\right)$ and hence the index map into this group turn out to vanish whenever $Z\hookrightarrow X$ is a uniform coarse equivalence. This happens in particular when $X$ is a finite-width bordism and $Z$ its boundary, making this index map insufficient for our purposes. 

As an alternative we introduce an index map that is in some sense closer to the uniform coarse index map on non-relative K-homology. We introduce the \emph{(relative) uniform structure algebras} $D^*_u(X)$ and $D^*_u(X,Z)$, which are finite-propagation counterparts to the (relative) uniform dual algebras. Under mild assumptions the quotient $D^*_u(X)/D^*_u(X,Z)$ turns out to be isomorphic to $\fD^{u}(X)/\fD^{u}(X,Z)$. Thus, using Paschke duality we may define the \emph{relative uniform index map} as 
$$
\relInd: \; K^{u}_*(X,Z) \cong K_{*+1}\left( D^*_u(X)/D^*_u(X,Z) \right) \xlongrightarrow{\partial} K_*\left(D^*_u(X,Z)\right) \; .
$$
The relative uniform structure algebra and hence the relative uniform index turn out to be hybrid objects combining the uniform coarse geometry of $X$ with the metric geometry of $Z$. Indeed, we will derive an algebraic Mayer-Vietoris sequence
$$
\cdots \to K_*\left(C^*_u(Z)\right) \to K_*\left(C^*_u(X)\right) \oplus K_*\left(D^*_u(Z)\right) \to K_*\left(D^*_u(X,Z)\right) \xrightarrow{\delta} K_{*-1}\left(C^*_u(Z)\right) \to \cdots
$$
expressing the K-theory of $D^*_u(X,Z)$ in terms of that of $C^*_u(Z)$, $C^*_u(X)$ and $D^*_u(Z)$.  While the uniform Roe algebras measure only the uniform coarse geometry, the uniform structure algebra $D^*_u(Z)$ depends on the metric geometry of $Z$. Moreover, we have seen in our discussion of higher $\rho$-invariants that the K-theory of structure algebras is the natural habitat of secondary invariants. Thus we see the relative uniform index map as retaining secondary information over $Z$. This intuition will be made explicit when we compute the relative uniform index of Dirac operators below.

The index map to $K_*(C^*_u(X,Z))$ factors through that to $K_*(D^*_u(X,Z))$. Indeed, there is a forgetful map $K_*(D^*_u(X,Z))\to K_*(C^*_u(X,Z))$ discarding the metric information on $Z$, so that the index map $\Ind: \, K^{u}_*(X,Z)\to K_*(C^*_u(X,Z))$ factors as 
$$
K^{u}_*(X,Z) \xlongrightarrow{\relInd} K_*(D^*_u(X,Z)) \longrightarrow K_*(C^*_u(X,Z)) \; .
$$
In fact the transformation from the long exact sequence in uniform K-homology to that in K-theory of the uniform Roe algebras factors through the algebraic Mayer-Vietoris, though in a non-commuting way. We summarize the situation in the following diagram: \\
\centerline{\xymatrix{
K_*^{u}(Z) \ar[r] \ar[d]_{\Ind} & K^{u}_*(X) \ar[r] \ar@{.>}[d]_{\Ind\oplus 0} & K^{u}_*(X,Z) \ar[r]^\partial \ar[d]_{\relInd} & K_{*-1}^{u}(Z) \ar[d]_{\Ind} \\ 
K_*\left(C^*_u(Z)\right) \ar[r] \ar@{=}[d] & K_*\left(C^*_u(X)\right)\oplus K_*\left(D^*_u(Z)\right) \ar@{.>}[d]_{\id\oplus 0} \ar[r] & K_*\left(D^*_u(X,Z)\right) \ar[d] \ar[r]^\delta & K_{*-1}\left(C^*_u(Z)\right)\ar@{=}[d] \\
K_*\left(C^*_u(Z)\right) \ar[r] & K_*\left(C^*_u(X)\right) \ar[r] & K_*\left(C^*_u(X,Z)\right) \ar[r]_\partial & K_{*-1}\left(C^*_u(Z)\right)
}}
The dotted vertical arrows indicate that the small squares containing them do \emph{not} commute due to the presence of the summand $K_*(D^*_u(Z))$. Were one to omit that summand the entire diagram would commute, but the second row would no longer be exact. The large squares from the first to the third row do however commute. The relative uniform index map to the relative uniform structure algebra thus presents an interesting but structurally unusual refinement of that to the uniform Roe algebra. In our applications we will mostly be interested in the former.
\subsubsection*{Partitioning hypersurfaces}
As a first application of our machinery we derive partitioned-manifold index theorems for uniform indices. The study of partitioned-manifold indices goes back to Roe \cite{Roe1989}. The goal is the reduction of index-theoretical information on a manifold to an index on a hypersurface partitioning the manifold into two halves; with the upshot that understanding the index on the hypersurface allows conclusion about the original manifold. In Roe's original version one considers compact partitioning hypersurface, so that that the index of the Dirac operator on it is simply given by the $\hat{A}$-genus. If this number is non-zero, then the entire manifold cannot admit a metric of uniformly positive scalar curvature. The partitioned-manifold index has since also been extended to non-compact partitioning hypersurfaces \cite{Hochs2025}.

We approach partitioned-manifold indices through Mayer-Vietoris sequences. We prove a Mayer-Vietoris sequence both for relative uniform K-homology (a corresponding sequence in non-relative uniform K-homology has already been constructed by \v{S}pakula \cite{Spakula2009}) and for the K-theory of relative uniform structure algebras. The key abstract prerequisite is that the relative uniform index coarse map takes one to the other.
\begin{thm}
Let $(X,Z)=(A_1\cup A_2,B_1\cup B_2)$ be a Mayer-Vietoris decomposition with $B_j=Z\cap A_j$ satisfying certain boundedness conditions. Then, the relative uniform index maps is a transformation from the Mayer-Vietoris sequence in relative uniform K-homology to the one in K-theory of the relative uniform structure algebras. In particular, there is a commutative diagram \\
\centerline{\xymatrix{
K^{u}_{*+1}(X,Z) \ar[r]^-\delta \ar[d]_{\relInd} & K^{u}_*(A_1\cap A_2,B_1\cap B_2) \ar[d]^{\relInd} \\
K_{*+1}\left( D^*_u(X,Z) \right) \ar[r]_-\delta & K_*\left( D^*_u(A_1\cap A_2,B_1\cap B_2) \right) 
}}
where $\delta$ are the respective Mayer-Vietoris boundary maps.
\end{thm}
Suppose that $\olO$ is an $(n+1)$-dimensional spin manifold with boundary and bounded geometry, and suppose that $N\subseteq \olO$ is a suitably bounded hypersurface partitioning $\olO$ into two pieces such that $\olO$ has product structure near $N$. The hypersurface potentially also has a boundary given by $\partial N=N\cap \dO$. Then $N$ gives rise to a Mayer-Vietoris decomposition. In particular there is a commutative diagram \\
\centerline{\xymatrix{
K^{u}_{n+1}(\olO,\dO) \ar[r]^\delta \ar[d]_{\relInd} & K^{u}_n(N,\partial N) \ar[d]^{\relInd} \\
K_{n+1}\left( D^*_u(\olO,\dO) \right) \ar[r]_\delta & K_n\left( D^*_u(N,\partial N) \right) 
}}
We define the \emph{relative uniform partitioned-manifold index} $\Ind^\mathrm{rel,pm}(\Dirac_\Omega)$ as the image of $[\Dirac_\Omega]\in K^{u}_{n+1}(\olO,\dO)$ under the composition
$$
K^{u}_{n+1}(\olO,\dO) \xlongrightarrow{\relInd} K_{n+1}\left( D^*_u(\olO,\dO) \right) \xlongrightarrow{\delta} K_n\left( D^*_u(N,\partial N) \right) \; .
$$
The above diagram tells us that we can compute the relative uniform partitioned-manifold index via the relative uniform index of $\delta[\Dirac_\Omega] \in K^{u}_n(N,\partial N)$. An adaptation of boundary-of-Dirac-is-Dirac shows this class to coincide with $(-1)^{n+1}[\Dirac_N] \in K^{u}_n(N,\partial N)$. We conclude the relative uniform partitioned-manifold index theorem.
\begin{thm}[Relative uniform partitioned-manifold index theorem]
Let $\olO$ be an $(n+1)$-dimensional spin manifold with boundary and bounded geometry, and $N$ a bounded-geometry hypersurface that partitions $\olO$ into two connected components. Assume that $\olO$ has product structure near $N$. Then,
$$
\Ind^\mathrm{rel,pm}(\Dirac_\Omega) = (-1)^{n+1} \relInd(\Dirac_N) \; \in \; K_n\left( D^*_u(N,\partial N) \right) \; .
$$
\end{thm}
As a corollary we also obtain a non-relative version in the case of an empty boundary. In this case there is a Mayer-Vietoris sequence in K-theory of the uniform Roe algebras, and the uniform coarse index map provides a transformation from the non-relative Mayer-Vietoris sequence in uniform K-homology to this one. Suppose $M$ is a spin manifold of bounded geometry, and $N\subseteq M$ a bounded-geometry partitioning hypersurface such that $M$ has product structure near $N$. In this case we define the \emph{uniform partitoned-manifold index} $\Ind^\mathrm{pm}(\Dirac_M)$ as the image of $[\Dirac_M]\in K^{u}_{n+1}(M)$ under the composition
$$
K^{u}_{n+1}(M) \xlongrightarrow{\Ind} K_{n+1}\left(C^*_u(M)\right) \xlongrightarrow{\delta} K_n\left(C^*_u(N)\right) \; ,
$$
where $\delta$ is the boundary map of the Mayer-Vietoris sequence in K-theory of the uniform Roe algebras. Then, specializing the relative version to the case of empty boundary yields the uniform partioned-manifold index theorem.
\begin{cor}[Uniform partioned-manifold index theorem]
Let $M$ be an $(n+1)$-dimensional spin manifold of bounded geometry, and $N$ a bounded-geometry hypersurface that partitions $M$ into two connected components. Assume that $M$ has product structure near $N$. Then,
$$
\Ind^\mathrm{pm}(\Dirac_M) = (-1)^{n+1} \Ind(\Dirac_N) \; \in \; K_n\left( C^*_u(\Sigma) \right) \; .
$$
\end{cor}
The uniform partitoned-manifold index $\Ind^{pm}(\Dirac_M)$ is a function of the uniform coarse index $\Ind(\Dirac_M)$. In particular, if the latter vanishes, then so does the former. Thus, analogously to the original non-uniform version the uniform partitioned-manifold index theorem provides obstructions to the existence of a metric of uniformly positive scalar curvature on $M$.
\begin{cor}
Let $M$ be an $(n+1)$-dimensional spin manifold of bounded geometry. Suppose there exists a bounded-geometry partitioning hypersurface $N\subseteq M$ such that $M$ has product structure near $N$, and
$$
\Ind(\Dirac_N) \neq 0 \; .
$$ 
Then $M$ does not admit a metric of bounded geometry and uniformly positive scalar curvature in its uniform coarse equivalence class.
\end{cor}
\subsubsection*{Relative uniform indices and uniform higher secondary invariants}
As a more concrete geometric application we compute the relative uniform index of the Dirac operator on a spin manifold with boundary and bounded geometry under different assumptions on the scalar curvature. A reasonable first guess would be that $\relInd(\Dirac_\Omega)$ vanishes if $\olO$ has uniformly positive scalar curvature. This is not the case. Instead $\relInd(\Dirac_\Omega)$ turns out to be a somewhat subtle object, relating to both certain APS-indices and secondary invariants on the boundary. This discussion is close in spirit to the work of Piazza and Schick on the delocalized APS-index theorem \cite{Piazza2014}.

Our starting point is a localization principle. If $M$ is a complete spin manifold and $Z\subseteq M$ a closed subspace such that $M$ has uniformly positive scalar curvature outside $Z$, then the Dirac operator on $M$ admits a localized index $\mathrm{Ind}^Z(\Dirac_M) \in K_*(C^*(Z))$ that gets mapped to the coarse index $\mathrm{Ind}(\Dirac_M)\in K_*(C^*(M))$ under the inclusion map $K_*(C^*(Z)))\to K_*(C^*(M))$ (see for example \cite{Roe2016}). The same is true uniformly, and we prove it in the presence of a boundary as well. 
\begin{thm}
Let $\olO$ be a spin manifold with boundary and bounded geometry. Let $Z\subseteq \olO$ be a closed subset such that $\olO$ has uniformly positive scalar curvature outside $Z$. Then, the Dirac operator on $\olO$ admits a localized relative uniform index
$$
\Ind^Z(\Dirac_\Omega) \; \in \; K_*\left( D^*_u(Z\cup \partial \Omega, \partial \Omega)\right)
$$
that gets mapped to $\relInd(\Dirac_\Omega)$ under the inclusion map $K_*( D^*_u(Z\cup \partial \Omega, \partial \Omega))\to K_*(D^*_u(\olO;\partial\Omega))$.

In particular, if $\olO$ has uniformly positive scalar curvature, then the Dirac operator on $\olO$ admits a localized relative uniform index
$$
\Ind^{\partial\Omega}(\Dirac_\Omega) \; \in \; K_*\left( D^*_u(\partial \Omega)\right)
$$
that gets mapped to $\relInd(\Dirac_\Omega)$ under the inclusion map $K_*( D^*_u(\partial \Omega))\to K_*(D^*_u(\olO;\partial\Omega))$.
\end{thm}
Note that even if the subset $Z\subseteq \olO$ is disjoint from the boundary, or $Z$ is empty for that matter, the index is still localized at the boundary as well. This is because $\Dirac_\Omega$ is not invertible near the boundary even if $\olO$ has uniformly positive scalar curvature there, due to the presence of boundary terms in the Lichnerowicz formula.

The localization principle has a consequence for the uniform coarse index in $K_*(C^*_u(\olO,\dO))$. Indeed, it is a consequence of our description of that index in terms of the relative uniform index in $K_*(D^*_u(\olO,\dO))$ that if the latter is in the image of $K_*(D^*_u(\dO))$, then the former vanishes. We conclude that the uniform coarse index in $K_*(C^*_u(\olO,\dO))$ obstructs uniformly positive scalar curvature metric.
\begin{cor}
Let $\olO$ be a spin manifold with boundary and bounded geometry. If $\olO$ has uniformly positive scalar curvature, then
$$
\Ind(\Dirac_\Omega) = 0 \; \in \; K_*\left(C^*_u(\olO,\dO)\right) \; .
$$
\end{cor}
Recall however that $K_*(C^*_u(\olO,\dO))$ vanishes whenever $\olO$ has finite width. In this case the upsc-obstruction provided by the uniform coarse index vanishes trivially, hence is useless. Even supposing infinite width we do not expect that this obstruction is all that interesting, because all information near the boundary is discarded. Thus we would be interested to know what the relative uniform index in $K_*(D^*_u(\olO,\dO))$ looks like in this case.

Recall that if the spin manifold $\Sigma$ carries a metric of uniformly positive scalar curvature, then the class $[\Dirac_\Sigma] \in K_*(\Sigma)$ has a canonical lift to $K_{*+1}(D^*(\Sigma))$, namely the higher $\rho$-invariant. The same is true for uniform K-homology as well, i.e. the uniform K-homology of the Dirac operator on a spin manifold of bounded geometry and uniformly positive scalar curvature admits a canonical lift to a class $\rho^{u}(\Dirac_\Sigma) \in K_{*+1}(D^*_u(\Sigma))$, which we call the \emph{uniform higher $\rho$-invariant} of $\Sigma$. Thus, if $\olO$ is a spin manifold with boundary and bounded geometry that has uniformly positive scalar curvature and product structure near the boundary, then $\partial\Omega$ is a spin manifold of bounded geometry with uniformly positive scalar curvature and thus admits a uniform $\rho$-invariant $\rho^{u}(\Dirac_{\partial \Omega})$. We compute that in this scenario the localized index of $\Dirac_\Omega$ coincides with the uniform $\rho$-invariant of the boundary - at least if $\olO$ is even-dimensional.
\begin{thm}
Let $\olO$ be an even-dimensional spin manifold with boundary and bounded geometry that has uniformly positive scalar curvature and product structure near the boundary. Then, 
$$
\Ind^{\partial \Omega} (\Dirac_\Omega) = \rho^{u}(\Dirac_{\partial \Omega}) \; \in \; K_0\left( D^*_u(\partial\Omega) \right) \; .
$$
\end{thm}
The author is convinced that this result is equally true in odd dimensions. The issue is that in these dimensions the K-theory boundary map, via which the localized index is defined, is much harder to compute. Thus the problem is a lack of computational techniques, not a conceptual one. This is discussed further in Chapter 7.\footnote{Piazza and Schick also run into this issue that computations using the Fredholm module picture of higher indices and secondary invariants are sometimes feasible only in even dimensions. They prove the delocalized APS-index theorem in these dimensions only \cite{Piazza2014}. The proof for odd dimensions was given by Xie and Yu using localization algebras \cite{Xie2014}. Again, the issues presented by the Fredholm module picture and the advantages of localization algebras in these instances are discussed in Chapter 7.}

The equality $\Ind^{\partial \Omega} (\Dirac_\Omega) = \rho^{u}(\Dirac_{\partial \Omega})$ is a refinement of the K-homological boundary-of-Dirac-is-Dirac formula to the level of uniform structure algebras. Note that it in particular implies that $\Ind^{\partial \Omega}(\Dirac_\Omega)$ never vanishes even on a manifold with uniformly positive scalar curvature. This is because $[\Dirac_{\partial\Omega}]\in K^{u}_*(\partial\Omega)$ is never zero,\footnote{The comparison map $K^{u}_*(\partial\Omega)\to K_*(\partial\Omega)$ to analytic K-homology sends the class of the Dirac operator in uniform K-homology to that in analytic K-homology. \cite[Lemma 12.2.4]{HigsonRoe2000} states that the latter never vanishes, hence neither does the former.} hence neither is its pre-image $\rho^{u}(\Dirac_{\partial\Omega})$.

Nonetheless we can derive a vanishing statement. Recall from the definition of the relative uniform index map that the composition
$$
K_*^{u}(\olO,\partial\Omega) \xlongrightarrow{\relInd} K_*\left( D^*_u(\olO,\partial\Omega) \right) \longrightarrow K_*\left( D^*_u(\olO) \right)
$$
is zero by exactness. Thus, pushing the above index formula forward to the K-theory of the uniform structure algebra $D^*_u(\olO)$ results in the vanishing of the uniform $\rho$-invariant.
\begin{cor}
Let $\olO$ be an even-dimensional spin manifold with boundary and bounded geometry that has uniformly positive scalar curvature and product structure near the boundary. Then, 
$$
\rho^{u}(\Dirac_{\partial \Omega}) = 0 \; \in \; K_0\left( D^*_u(\Omega) \right) \; .
$$
\end{cor}
This is a form of bordism invariance of the uniform higher $\rho$-invariants: If a spin manifold of bounded geometry with uniformly positive scalar curvature is null-bordant via a spin bordism with bounded geometry and uniformly positive scalar curvature, then its uniform higher $\rho$-invariant vanishes upon push-forward to the uniform structure algebra of the bordism. Having a vanishing statement in the K-theory of the structure algebra of the null-bordism is somewhat artificial. The statement becomes more natural if we consider our spaces to be endowed with maps to some reference space $X$, and the indices and secondary invariants to also be pushed forward to $X$. Then, if we require the null-bordism to extend the reference map to $X$, the push-forward to $X$ factors through the null-bordism. Thus we conclude that the push-forward of the uniform higher $\rho$-invariant to $X$ vanishes. Piazza and Schick consider manifolds and bordisms with actions of a fixed group $\Gamma$, and take equivariant reference maps to $X=E\Gamma$. In the context of mapping positive scalar curvature to analysis it is precisely this form of bordism invariance that yields the well-definedness of the map $\mathrm{Pos}^\mathrm{spin}_*(B\Gamma) \to K_{*+1}(D^*(E\Gamma)^\Gamma)$. We will not treat group actions or maps to classifying spaces in this thesis though.

In the preceeding discussion we considered manifolds $\olO$ that have uniformly positive curvature. However, if we only require $\olO$ to have uniformly positive scalar curvature on the boundary, then the uniform $\rho$-invariant $\rho^{u}(\Dirac_{\partial\Omega})$ can still be defined. However, it will no longer be the only contribution to the relative uniform index of $\Dirac_\Omega$, but there is also an index contribution coming from the interior. To define this additional term let $M:= \olO \cup_{\partial\Omega} (\RR_{\leq 0}\times \partial \Omega)$ be the spin manifold of bounded geometry obtained by attaching cylinders to the boundary of $\olO$. Note that we have no issues with smoothness because we assumed $\olO$ to have product structure near the boundary. Assuming uniformly positive scalar curvature on $\partial\Omega$, $M$ has uniformly positive scalar curvature outside the closed subset $\olO\subseteq M$. Thus the Dirac operator on $M$ admits a localized index in $K_*(C^*_u(\olO))$. We denote this localized index by
$$
\Ind^{APS}(\Dirac_\Omega) \; \in \; K_*\left( C^*_u(\olO) \right) \; ,
$$
and call it the \emph{uniform APS-index}. This terminology derives from the fact that the APS-boundary condition geometrically corresponds to attaching infinite cylinders to the boundary. We also remark that the non-uniform equivariant version of this APS-index gives rise to the map $R_*^\mathrm{spin}(B\Gamma)\to K_*(C^*_r\Gamma)$ in the framework of mapping positive scalar curvature to analysis.
\begin{thm}\label{RelIndexFormula}
Let $\olO$ be an even-dimensional spin manifold with boundary and bounded geometry. Assume that $\olO$ has product structure near the boundary, and $\partial\Omega$ has uniformly positive scalar curvature. Then,
$$
\relInd(\Dirac_\Omega) = \rho^{u}(\Dirac_{\partial\Omega}) + \Ind^{APS}(\Dirac_\Omega)  \; \in \; K_{0}\left( D^*_u(\olO,\partial\Omega) \right) \; .
$$
\end{thm}
Again, we can push forward this index formula to the K-theory of $D^*_u(\olO)$, resulting in the vanishing of the relative uniform index. Thus we conclude the following.
\begin{cor}[Uniform delocalized APS-index theorem]
Let $\olO$ be an even-dimensional spin manifold with boundary and bounded geometry. Assume that $\olO$ has product structure near the boundary, and $\partial\Omega$ has uniformly positive scalar curvature. Then,
$$
\Ind^{APS}(\Dirac_\Omega) = -\rho^{u}(\Dirac_{\partial\Omega}) \; \in \; K_{0}\left( D^*_u(\olO) \right) \; .
$$
\end{cor}
This is the uniform version of what Piazza and Schick called the \emph{delocalized APS-index theorem}.\footnote{Actually, the index theorem in \cite{Piazza2014} differs from ours by as sign. This is due to a different choice of orientation. In \cite{Piazza2014} the manifold $M$ is oriented such that a frame $(\nu,e_1,\cdots,e_n)$ consisting of an oriented frame $(e_1,\cdots,e_n)$ of $T\Sigma$ and a normal vector pointing in the cylinder direction is positively oriented, we choose the orientation in such a way that it is negatively oriented. Our orientation has the effect that it coincides with the given orientation of $\olO$, so that $\Dirac_M$ and $\Dirac_\Omega$ agree over $\olO$. This is the natural choice for our purposes, since first and foremost we wish to compute $\relInd(\Dirac_\Omega)$. On the other hand, the orientation used by Piazza and Schick coincides with the natural orientation on the cylinder, which is in turn natural for their purposes of reduction to a computation on the cylinder.} We refer to \cite[Section 2.3]{Piazza2014} for an explanation of why it should be viewed as a \emph{delocalized} index theorem. Again, the author is convinced that these results hold in the odd-dimensional case. The reason why we can prove the uniform delocalized APS-index theorem only in even dimensions is that contributions from the attached cylinder are treated via the formula $\Ind^{\partial \Omega}(\Dirac_{\RR_{\geq 0}\times\partial\Omega})=\rho^{u}(\Dirac_{\partial\Omega})$, which we can only prove in even dimensions. We discuss the problems involved in the odd case in Chapter 7.
\subsubsection*{Non-uniform counterparts}
The uniform theory is a finer version of the non-uniform version under additional boundedness assumptions, and forgetting uniformness allows passage from a uniform result to its non-uniform counterpart in analytic K-homology and coarse index theory. For the results of this thesis these non-uniform versions are in most cases well-established, and our contribution is simply to note that they hold uniformly as well. This is the case for example for boundary-of-Dirac-is-Dirac and its consequences, and of course the entire theory of relative uniform K-homology. There are however results which the author could not find stated even non-uniformly in the literature. This includes our Theorem \ref{IntroThmClasses}. As already discussed above to the author's knowledge the only work dealing with relative K-homology classes of elliptic operators on non-compact manifolds with boundary is a recent article by Fries \cite{Fries2025}. Concerning manifolds of bounded geometry it is shown that the minimal boundary condition for elliptic operators in a certain elliptic calculus on a manifold of bounded geometry define a relative K-homology class over any open subset. This is a special case of the non-uniform version of Theorem \ref{IntroThmClasses} (ii), since the minimal boundary condition is local and semi-regular. Thus our Theorem \ref{IntroThmClasses} extends the result of \cite{Fries2025} to more general boundary conditions, though it must be pointed out that \cite{Fries2025} also covers geometric settings besides the bounded-geometry one that are not covered by our discussion.

Another aspect of which the author has not encountered a non-uniform counterpart is our approach to relative index maps into the K-theory of relative structure algebras. The closest analogue of our $\relInd$ is the relative index of Chang, Weinberger and Yu. However, the latter cannot be an equivariant refinement of our relative index, since the Chang-Weinberger-Yu index vanishes in the presence of positive scalar curvature, whereas our relative index reproduces the $\rho$-invariant, which never vanishes. Thus, the index of Chang, Weinberger and Yu would likely be closer to an equivariant version of the index in $K_*(C^*(\olO,\dO))$, since the latter at least vanishes on upsc-manifolds. However, the rather forgetful nature of the index in $K_*(C^*(\olO,\dO))$ makes the author unsure if this is really to be expected. Additional speculation on the potential relation of our indices to that of Chang, Weinberger and Yu is to be found in the outlook to this thesis.

Let us also mention that there is a sharpening of the delocalized APS-index theorem for finite-width manifolds due to Zeidler \cite{Zeidler2016} that is formally analogous to our Theorem \ref{RelIndexFormula}. There the role of the relative index is played by a $\rho$-invariant of the bordism $\olO$, which is constructed from a unique lift of $[\Dirac_M] \in K_*(M=(\RR_{\geq 0}\times\dO)\cup_{\dO} \olO)$ to $K_*(D^*(M))$ that exists due to the assumption of finite width. Zeidler's $\rho$-invariant of the bordism must coincide with our relative index, simply because both equate (modulo sign conventions) to the right-hand side of the index formula of non-uniform version of Theorem \ref{RelIndexFormula}. However, seeing this equality directly seems not entirely trivial, given that Zeidler works with localization algebras and defines the $\rho$-invariant of the bordism less constructively. Also, Zeidler's definition only works for finite-width $\olO$, a restriction not present in our set-up. 

In any case, if our approach is indeed novel, it is most certainly not because it is somehow more non-trivial. Quite the opposite: It is a straight-forward adaptation of the usual construction of the coarse index map. It does, however, rely crucially on the fact that, say, Dirac operators on manifolds with boundary define relative K-homology cycles, which -as discussed again just now- might not have been shown in the literature for non-compact manifolds until very recently. Also, a lot of recent work on coarse index theory uses localization algebras instead of Fredholm modules to represent K-homology, and the translation between the two versions is not always obvious. It is entirely possible that an equivalent version of our relative index map using localization algebras already exists, but was either not found or not recognized as such by the author. However, in this scenario the knowledge of an equivalent description using Fredholm modules instead of localization algebras would still be interesting.
\subsubsection*{Outline}
Let us give a brief outline of the content of each subsequent chapter. More detailed descriptions can be found in the introduction to each chapter.
\begin{itemize}
\item In Chapter 2 relative uniform K-homology is introduced and its properties are established. Via Paschke duality relative uniform K-homology is related to the K-theory of quotients of uniform dual algebras. This picture is used to prove the excision isomorphism, and to derive the long exact sequence as well as the Mayer-Vietoris sequence for relative uniform K-homology. Moreover the suspension isomorphism is proved. This theory is applied to compute the uniform K-homology of graphs of bounded geometry.
\item Chapter 3 treats index maps on relative uniform K-homology. After a review of uniform Roe algebras we introduce the relative uniform Roe algebra and the uniform coarse index map to its K-theory. Then the relative uniform structure algebras are defined and the functoriality of their K-theory is discussed. Then Mayer-Vietoris sequences for the relative uniform structure algebras are derived. The relative uniform index maps are defined, and their naturality and compatibility with Mayer-Vietoris boundary maps is proved.
\item Chapter 4 revolves around Sobolev spaces on domains in manifolds of bounded geometry. After establishing the geometric basics different Sobolev spaces on domains are defined and compared. Then, extensions and trace maps are discussed, and their behavior under interpolation is investigated. Lastly, the uniformly summable Rellich-Kondrachov theorem is proved.
\item Chapter 5 discusses uniformly elliptic operators on domains. After their definition boundary conditions are discussed, and notions of locality and regularity are studied. These properties are investigated for geometrically relevant examples. Finite-propagation results for the case of local self-adjoint boundary conditions for uniformly elliptic first-order differential operators are given.
\item In Chapter 6 the formal operator calculus is set up, with the help of which the existence of the different uniform K-homology classes associated to a uniformly elliptic operator on a domain is proved. Finite-propagation representatives in the case of first-order differential operators are discussed. Behavior of these classes under products and restriction to subsets is discussed, culminating in a proof of boundary-of-Dirac-is-Dirac and the bordism invariance of the uniform coarse index, as well as the uniform versions of the partitioned-manifold index theorem.
\item In Chapter 7 the relative uniform index of the Dirac operator on a spin manifold with boundary and bounded geometry is investigated under different assumptions on the scalar curvature. A localization principle is proved. Assuming uniformly positive scalar curvature everywhere the index localized at the boundary is computed to agree in even dimensions with the uniform $\rho$-invariant of the boundary. Assuming only uniformly positive scalar curvature on the boundary the uniform delocalized APS-index theorem is proved in even dimensions.
\item The appendices contain a combination of reviews of prerequisite material, and of technical discussions that were deemed too distracting for the flow of the main part of this thesis. Appendix \ref{MultigradingAppendix} contains a review of (multi)gradings, as well as spinor bundles and Dirac operators. Appendix \ref{AppendixUA} contains a collection of basic facts about uniformly approximable families and finitely summable operators that are used throughout this thesis. Appendix \ref{AppendixUnbdOps} contains a review of unbounded operators. Basic definitions and facts about operator K-theory are collected in Appendix \ref{AppendixKTheory}. Lastly, Appendix \ref{AppendixUnifPairs} contains sufficient conditions for metric pairs to fulfill the boundedness required for excision and the long exact sequence in uniform K-homology.
\end{itemize}
\chapter{Relative uniform K-homology}\label{ChapterKHom}
In this chapter we define and investigate relative uniform K-homology. This relative version of uniform K-homology provides the receptacle for the relative classes built from boundary conditions for uniformly elliptic operators in Chapter 6. The properties of relative uniform K-homology we prove here are mainly the usual suspects one would expect from any relative homology theory, i.e. excision, long exact sequence, suspension isomorphism, and so on. However, we point to our discussion of the work of Baum, Douglas and Taylor on the K-homology classes of elliptic operators on manifolds with boundary in the introduction. There these abstract properties had concrete consequences for the classes of elliptic operators. This will be no different in the uniform case. Thus the results of this chapter lay the abstract groundwork for later calculations with the uniform K-homology classes of uniformly elliptic operators over manifolds with boundary.

Many tools from analytic K-homology have already been adapted to analytic K-homology, e.g. Paschke duality and the Mayer-Vietoris sequence \cite{Spakula2009}, as well as the Kasparov product \cite{Engel2019}. We review some background on (non-relative) uniform K-homology in Section 2.1. Then, we define relative uniform K-homology and work out its basic properties in Section 2.2. Section 2.3 is devoted to Paschke duality for relative uniform K-homology. This is an algebraic statement connecting relative uniform K-homology to the K-theory of uniform dual algebras. This allows for the transfer of results in operator K-theory to (relative uniform) K-homology. Indeed, in Section 2.4 we will deduce excision in uniform K-homology from a property of these dual algebras, and the long exact sequence and relative Mayer-Vietoris sequence in uniform K-homology from the corresponding sequence in K-theory. In Section 2.5 we discuss the suspension isomorphism and its connection to the Kasparov product. The results of this section will be central to our proof of boundary-of-Dirac-is-Dirac in uniform K-homology. Lastly, supplementing this abstract chapter with a more concrete result we compute the uniform K-homology of graphs of bounded geometry in Section 2.6.

Formally uniform K-homology is very much analogous to analytic K-homology, which is of course by design. At the risk of over-simplification we may divide results in analytic K-homology into two categories. Those in the first category are purely algebraic, requiring only formal properties of the algebra of compact operators and some matrix algebra. Those in the second require some non-trivial analytic or geometric tool, say the existence of a suitable partition of unity, Voiculescu's theorem, or Kasparov's technical theorem. The proofs of results in the first category carry over verbatim to uniform K-homology, with compact operators replaced by uniformly compact ones. Given that uniformly compact operators have the same formal properties as compact ones, this causes no issues. Proofs of results in the second category require checking that the analytic/geometric tool can be strengthened to apply to the uniform setting as well. Usually this forces one to assume some version of bounded geometry for the underlying metric space. 

Our discussion of relative uniform K-homology and its basic properties (Section \ref{SecRelUKHom}), as well as of relative uniform Paschke duality (Section \ref{SecPaschke}) fall entirely into the first category. We loosely follow Chapters 5 and 8 of \cite{HigsonRoe2000}, as well as  \cite{Higson1995}. Given that proofs are essentially the same algebraic calculations as in analytic K-homology, we will move more quickly here. The proofs of excision and the long exact sequence belong to the second category; relying on the existence of an adapted partition of unity and of a uniform completely positive section respectively. The existence of these objects for a given pair is not guaranteed, but requires certain conditions on the geometry (Definition \ref{ConditionGeomSubspace}). Since it uses the Kasparov product the suspension isomorphism would also fall into the second category. We are in the lucky position that the Kasparov product on uniform K-homology has already been constructed by Engel \cite{Engel2019}, so that we are free to simply use its properties, so that our discussion of the suspension isomorphism will be rather analogous to the analytic case, see \cite[Chapter 9]{HigsonRoe2000}.
\section{Uniform K-homology}
Let us begin with a review of uniform K-homology. First of all, to save ink throughout we introduce the following convention.
\begin{convention}
Throughout this thesis all metric spaces $X$ are assumed to be locally compact and second-countable, so that $C_0(X)$ is a separable $C^*$-algebra. Moreover, all Hilbert spaces are assumed to be separable.
\end{convention}
Uniform K-homology is defined via strengthening the relation of equality modulo compact operators to a uniform version. Thus we first introduce this uniform equality modulo compacts. Let $X$ be a metric space. For fixed $L,R\geq 0$ define 
\begin{equation}\label{LLip}
\LLip_R(X):= \left\{ \, f\in C_c(X) \, |  \, ||f||\leq 1, \, f \, \mathrm{is} \, L\text{-}\mathrm{Lipschitz}, \, \diam(\supp(f))\leq R \, \right\} \; .
\end{equation}
\begin{definition}
Let $H$ be a Hilbert space. Let $\sigma: \, C_0(X)\to\fBH$ be a map. We write $\sigma\sua 0$ if for all $L,R \geq 0$ the set $\sigma(\LLip_R(X))$ is a \emph{uniformly approximable} family of compact operators, meaning that for all $\varepsilon>0$ there exists $N\in\NN$ such that for all $f\in\LLip_R(X)$ there exists a rank-$N$ operator $K\in\fBH$ such that $\| \sigma(f)- K\|<\varepsilon$.

If $\sigma': \, C_0(X)\to \fBH$ is another map, then we write $\sigma\sua\sigma'$ if $\sigma-\sigma'\sua 0$.
\end{definition}
The equivalence relation $\sua$ is the analogue of equality modulo compact operators used in analytic K-homology. It enjoys the same formal properties, in that it is closed under sums, product, limits, adjoints, and functional calculus. Most of these properties can be seen immediately via $\frac{\varepsilon}{2}$-arguments. Nonetheless Appendix \ref{AppendixUA} collects basic facts concerning uniformly approximable families and the relation $\sua$. 

With the uniform equality modulo compacts defined we can turn to the definition of uniform Fredholm modules. It uses the language of multigradings, whose definition along with some basic facts and constructions are collected in Appendix \ref{MultigradingAppendix}.
\begin{definition}
Let $X$ be a metric space, and $p\in \ZZ_{\geq -1}$. A \emph{$p$-multigraded uniform Fredholm module} over $X$ is a triple $(H,\rho,F)$ consisting of
\begin{itemize}
\item a $p$-multigraded Hilbert space $H$, 
\item a representation $\rho: \, C_0(X)\to \fBH$ by even $p$-multigraded operators, 
\item and an odd $p$-multigraded operator $F\in \fBH$,
\end{itemize} such that
$$
[F,\rho] \sua 0 \; , \; (I-F^2)\cdot \rho \sua 0 \; , \; (F-F^*)\cdot \rho \sua 0 \; .
$$
If $p=0$, we speak of a \emph{graded uniform Fredholm module}. If $p=-1$, we speak of an \emph{ungraded uniform Fredholm module}. In these cases the condition that $\rho$ and $F$ be $p$-multigraded are void. In the ungraded case the condition that $\rho$ and $F$ be even and odd, respectively, are also void.
\end{definition} 
We note that since $\sua$ is stronger than equality modulo compacts a uniform Fredholm module is automatically a Fredholm module in the usual sense.
\begin{example}[{\cite[Proposition 2.8]{Engel2019}}]
If $X$ is totally bounded, then every Fredholm module over $X$ is automatically a uniform Fredholm module.
\end{example}
\begin{example}
\v{S}pakula proved in \cite{Spakula2009} that if $D$ is a Dirac type operator acting on a Clifford bundle $S$ over a (boundaryless) manifold of bounded geometry $M$, then $\left(L^2(M;S),\rho,D(D^2+I)^{-1/2}\right)$ is a uniform Fredholm module over $M$. Here $\rho: \, C_0(M)\to \fB(L^2(M;S))$ is the representation by multiplication operators. This result was generalized by Engel to $D$ being any symmetric and uniformly elliptic pseudo-differential operator of positive order acting on a vector bundle of bounded geometry over a manifold of bounded geometry \cite{Engel2018}.
\end{example}
\begin{example}\label{ExGrahKhom}
Let $X$ be a graph of bounded geometry, meaning that there is a uniform bound on how many edges meet a given vertex, equipped with the path metric. Let $V$ denote the vertex set of $X$. Set $H:=\ell^2(V)\otimes \ell^2(\NN)$. There is a representation $\rho$ of $C_0(X)$ on $\ell^2(V)$ by multiplication operators. Set $\tilde{\rho}:=\rho \otimes I: \, C_0(X)\to\fBH$. Let $S\in\fB(\ell^2(\NN))$ denote the left-shift operator, and set $\tilde{S}:= I\otimes S$. Then, 
$$
\left( H\oplus H, \tilde{\rho}\oplus\tilde{\rho}, \begin{pmatrix}
0 & \tilde{S}^* \\ \tilde{S} & 0
\end{pmatrix} \right)
$$
is a graded uniform Fredholm module over $X$. \v{S}pakula shows in \cite[Section 11]{Spakula2009} that its uniform K-homology class vanishes if and only if $X$ is amenable. Compare Proposition \ref{AmenableGraphs}.
\end{example}
The relations imposed on uniform Fredholm modules to obtain uniform K-homology are the same as for analytic K-homology. Suppose that $(H,\rho,F)$ and $(H',\rho',F')$ are $p$-multigraded uniform Fredholm modules. Their direct sum $(H\oplus H', \rho\oplus\rho',F\oplus F')$ is again a $p$-multigraded uniform Fredholm module. Moreover, they are called \emph{unitarily equivalent} if there is a unitary $U: \, H\to H'$, preserving (multi)gradings if present, such that $\rho'=U\rho U^*$ and $F'=UFU^*$. Lastly, if $H=H'$ and $\rho=\rho'$, we call the two modules \emph{operator-homotopic} if there is a continuous path $t\mapsto F_t$ of operators with $F_0=F$ and $F_1=F'$ such that $(H,\rho,F_t)$ is a $p$-multigraded uniform Fredholm module for every $t$.
\begin{definition}
Let $X$ be a metric space and $p\in\ZZ_{\geq -1}$. The \emph{$(-p)$th uniform K-homology group} $K^{u}_{-p}(X)$ is defined as the free group on the set of unitary equivalence classes of $p$-multigraded uniform Fredholm modules over $X$, modulo the relations:
\begin{itemize}
\item $[H,\rho,F]+[H',\rho',F']=[H\oplus H',\rho\oplus\rho',F\oplus F']$,
\item $[H,\rho,F]=[H',\rho',F']$ if $(H,\rho,F)$ and $(H',\rho',F')$ are operator-homotopic.
\end{itemize}
\end{definition}
The uniform K-homology groups enjoy the same basic properties as analytic K-homology. These will be discussed in the more general context of relative uniform K-homology below. Let us only mention here that there is a formal two-fold periodicity $K_{-(p+2)}^{u}(X)\cong K_{-p}^{u}(X)$ coming from the two-fold periodicity of mutigradings.
\begin{example}
If $X$ is totally bounded, then $K^{u}_*(X)=K_*(X)$.
\end{example}
\begin{example}[{\cite[Lemma 2.18]{Engel2019}}]\label{uniKHomDiscrete}
Suppose $X$ is a uniformly discrete metric space of bounded geometry, meaning that $\inf_{x\neq y} d(x,y)>0$ and $\sup_x|B_r(x)|<\infty$ for all $r\geq 0$. Then, $K_0^{u}(X)\cong \ell^\infty_\ZZ(X)$ and $K^{u}_0(X)=0$.
\end{example}
Let us turn to the functoriality of uniform K-homology.
\begin{definition}
Let $X$ and $Y$ be metric spaces. A $^*$-homomorphism $\phi: \, C_0(Y)\to C_0(X)$ will be called a \emph{filtered map} from $X$ to $Y$ if for all $L,R\geq 0$ there exist $L',R'$ such that
$$
\phi\left( \LLip_R(Y) \right) \subseteq L'\text{-}\mathrm{Lip}_{R'}(X) \; .
$$ 
\end{definition}
Note the reversal of direction; a filtered map $X\to Y$ is a  homomorphism $C_0(Y)\to C_0(X)$. This is because every continuous map $X\to Y$ induces a homomorphism $C_0(Y)\to C_0(X)$. But there are relevant homomorphisms $C_0(Y)\to C_0(X)$ that do not come from continuous maps $X\to Y$, for example the inclusion map $C_0(Y)\hookrightarrow C_0(X)$ if $Y$ is an open subset of $X$. Of course, it is a consequence of Gelfand duality that homomorphisms $C_0(Y)\to C_0(X)$ correspond uniquely to continuous maps $X^+\to Y^+$ between the respective one-point compactifications that preserve the point at infinity. Under this correspondence the inclusion $C_0(Y)\hookrightarrow C_0(X)$ corresponds to the map $X^+\to Y^+$ collapsing the complement of $Y$ to a point. While this equivalent topological description exists, we find it easier to work with algebra homomorphisms instead.
\begin{example}[{\cite[Discussion after Definition 2.15]{Spakula2009}}] \label{Lipfilter}
Let $g: \, X\to Y$ be a proper Lipschitz map such that $\sup_{y\in Y} \diam(g^{-1}B_R(y))<\infty$ for all $R\geq 0$. A map with this property is said to be \emph{uniformly cobounded} or \emph{effectively proper}. Then $g^*: \, C_0(Y)\to C_0(X)$ is a filtered map from $X$ to $Y$ ().
\end{example}
So far uniform K-homology has only been considered as functorial under uniformly cobounded Lipschitz maps \cite{Spakula2009}, \cite{Engel2019}. However, the additional flexibility provided by the notion of filtered maps can be useful. We will see an example in the discussion surrounding Lemma \ref{DifferentMetrics}. There we consider two metrics that are bi-Lipschitz equivalent on uniformly small scales but not globally. In this context the identity map from the space with one metric to the other is not generally uniformly cobounded and Lipschitz, but it is a filtered map.

It is immediate that if $\phi: \, C_0(Y)\to C_0(X)$ is a filtered map, and $\sigma: \, C_0(X)\to \fBH$ satisfies $\sigma\sua 0$, then $\sigma\circ \phi\sua 0$ as well. Thus, if $(H,\rho,F)$ is a uniform Fredholm module over $X$, then $\phi_*(H,\rho,F):=(H,\rho\circ\phi,F)$ is a uniform Fredholm module over $Y$. The operation $\phi_*$ is compatible with unitary equivalence, direct sums, and operator-homotopies. We conclude:
\begin{prop}
If $\phi: \, X\to Y$ is a filtered map, then 
$$
\phi_*: \; K^{u}_*(X)\longrightarrow K^{u}_*(Y) , \; [H,\rho,F]\longmapsto [H,\rho\circ\phi,F]
$$ 
is a well-defined group homomorphism. Thus uniform K-homology is a functor from the category of metric spaces\footnote{Recall that by our convention we really mean the category of all locally compact and second-countable metric spaces.} and filtered maps to the category of abelian groups.
\end{prop}
\begin{remark}\label{ApproxFilter}
To get the implication $\sigma\sua 0 \; \Rightarrow \sigma\circ\phi\sua 0$ it is actually not necessary that $\phi$ be a filtered map. It suffices to assume that $\phi: \, C_0(Y)\to C_0(X)$ is \emph{approximately filtered}, in the sense that for $L,R\geq 0$ and $\varepsilon>0$ there exist $L',R'$ such that $\phi\left( \LLip_R(Y) \right)$ is contained in the $\varepsilon$-neighborhood of $L'\text{-}\mathrm{Lip}_{R'}(X)$. This additional bit of flexibility will be useful when showing that any pair $(X,Z)$ with $X$ having jointly bounded geometry has uniformly controlled geometry. We will also show in Lemma \ref{UCapprox} below that under geometric boundedness assumptions Example \ref{Lipfilter} can be relaxed from Lipschitz continuity to uniform continuity, only that the resulting map will generally be approximately filtered. Let us also note that being approximately filtered behaves nicer with respect to limits. For example, if $g_n: \, X\to Y$ is a uniformly convergent sequence of maps, then the limit map $g$ is an approximately filtered map if every $g_n$ is approximately filtered. On the other hand $g$ need not be filtered even if every $g_n$ is filtered.
\end{remark}
Next, we recall the definition of jointly bounded geometry for metric spaces.
\begin{definition}[{\cite[Definition 2.20]{Engel2019}}]\label{DefJBB}
A metric space $X$ is said to have \emph{jointly bounded geometry} if it admits a countable Borel decomposition $X=\bigcup_i X_i$ and a discrete subset $\Gamma\subseteq X$ such that:
\begin{itemize}
\item each $X_i$ is totally bounded and has non-empty interior,
\item for all $\varepsilon>0$ there is $N\in\NN$ such that each $X_i$ admits an $\varepsilon$-net of cardinality at most $N$,
\item $\Gamma$ is a quasi-lattice, meaning that there is $c\geq 0$ such that $N_c(\Gamma)=X$, and for all $r\geq 0$ it holds that $\sup_{x\in X} |\Gamma\cap B_r(x)|<\infty$,
\item for all $r\geq 0$ it holds that $\sup_{y\in\Gamma}|\{ i \, | \, X_i\cap B_r(y)\neq\emptyset\}|<\infty$.
\end{itemize}
\end{definition}
This condition guarantees a number of favorable properties for uniform K-homology, chief among them the existence of the Kasparov product. However, we will first demonstrate how jointly bounded geometry allows us to extend the functoriality of uniform K-homology from Lipschitz-continuous to uniformly continuous maps.
\begin{bonuslemma}\hspace{-5pt}\footnote{The alphabetic index means that this lemma was not present in the original version of this thesis. Alphabetic indexing was chosen so that the numerical indexing is consistent with the original version.}\label{UCapprox}
Let $g: \, X\to Y$ be a map between metric spaces. Suppose that $X$ has jointly bounded geometry, and that $g$ is uniformly continuous, proper, and uniformly cobounded. Then, $g: \, X\to Y$ is approximately filtered in the sense of Remark \ref{ApproxFilter}, and hence defines a homomorphism $g_*: \, K^{u}_*(X)\to K^{u}_*(Y)$. 
\end{bonuslemma}
\begin{proof}
Fix $L,R \geq 0$ and $\varepsilon>0$. We need to produce $L',R'\geq 0$ such that $f\circ g$ is within $\varepsilon$ of a function in $L'\text{-Lip}_{R'}(X)$ for all $f\in\LLip_R(Y)$. First of all, note that since $g$ is uniformly cobounded, the diameter of the support of $f\circ g$ is uniformly bounded by a constant $R'$ depending only on $R$. 

Next, fix some $\delta>0$. Because $X$ has jointly bounded geometry, we can find a uniformly discrete subset $\Gamma_\delta \subseteq X$ such that the open $\delta$-balls around $\Gamma_\delta$ form a uniformly locally finite open cover of $X$ with positive Lebesgue number. Indeed, let $X=\bigcup_i X_i$ be a Borel decomposition as in Definition \ref{DefJBB}, and define $\Gamma_\delta$ to be the union of uniformly finite $\frac{\delta}{2}$-nets for the $X_i$ as per the second bullet point of that definition. Then, the open $\delta$-balls around $\Gamma_\delta$ cover $X$, and the third and fourth bullet points combine to show that this cover is uniformly locally finite. Because $\Gamma_\delta$ is made of $\frac{\delta}{2}$-nets, each $x\in X$ is at most $\frac{\delta}{2}$ away from some element of $\Gamma_\delta$. Then, the $\delta$-ball around that element contains the $\frac{\delta}{2}$-ball around $x$. We conclude that the cover has Lebesgue number at least $\frac{\delta}{2}$.

Next, we fix a partition of unity $\{\psi_\gamma\}_{\gamma\in\Gamma_\delta}$ subordinate to the cover by $\delta$-balls around $\Gamma_\delta$. Because the cover has positive Lebesgue number, we can choose this partition of unity to be uniformly Lipschitz (see the proof of Lemma \ref{LemmaLipPU} in the appendix). Set 
$$
f_\delta(x) := \sum_{\gamma\in\Gamma_\delta} (f\circ g)(\gamma) \cdot \psi_\gamma(x) \; .
$$  
Due to the uniform local finiteness of $\Gamma_\delta$ and the uniform Lipschitzness of the $\psi_\gamma$, the function $f_\delta$ has Lipschitz constant depending only on $\delta$. Now, compute that
\begin{align*}
\left| (f\circ g)(x) - f_\delta(x) \right| &= \left| \sum_{\gamma\in\Gamma_\delta} \left( (f\circ g)(x) - (f\circ g)(\gamma) \right) \cdot \psi_\gamma(x) \right| \\
&\leq \sum_{\gamma\in\Gamma_\delta} \left| (f\circ g)(x) - (f\circ g)(\gamma) \right| \cdot \psi_\gamma(x) \\
&\leq L \cdot \sum_{\gamma\in\Gamma_\delta} d\left( g(x),g(\gamma)\right) \cdot \psi_\gamma(x) \; .
\end{align*}
For any given $x$ this sum has only uniformly finitely many non-zero summands, namely those $\gamma$ that contain $x$ in their $\delta$-neighborhood. By uniform continuity of $g$ we may choose $\delta>0$ so that for such $\gamma$ we have $d\left( g(x),g(\gamma)\right)<\frac{\varepsilon}{L}$. In this case we conclude that
$$
\left| (f\circ g)(x) - f_\delta(x) \right| \leq \frac{L\varepsilon}{L}\cdot \sum_{\gamma\in\Gamma_\delta} \psi_\gamma(x) = \varepsilon \; .
$$
Hence, for adequately chosen $\delta$ the function $f_\delta$ is $L''(\delta)$-Lipschitz with support of diameter bounded by $R'+2\delta$ such that $f\circ g$ is within $\varepsilon$ of $f_\delta$.

In the last step we need to modify $f_\delta$ to ensure it has compact support. Fix $\eta>0$, and a function $\psi_\eta:\, X \to [0,1]$ with $\psi_\eta=0$ outside $g^{-1}\supp(f)$, and $\psi_\eta=1$ on the complement of the $\eta$-neighborhood of $X\setminus g^{-1}\supp(f)$, such that $\psi_\eta$ has Lipschitz constant $L(\eta)$ dependent only on $\eta$. Because $g$ is proper, the pre-image of $\supp(f)$ under $g$ is compact, hence so is the support of $\psi_\eta$. Therefore, $\psi_\eta f_\delta$ is a compactly supported function with Lipschitz constant depending only on $\eta$ and $\delta$, and support of diameter less than $R''$. We wish to estimate the norm $\| \psi_\eta f_\delta - f\circ g \|$. If $d(x,X\setminus g^{-1}\supp(f)) \leq\eta$, then $|f_\delta(x)|$ is bounded by the supremum of $|(f\circ g)(y)|$ with $y$ in the $\delta$-neighborhood of $x$. In thise case, $y$ can have distance at most $\eta+\delta$ from the complement of the support of $f\circ g$. Making $\eta$ sufficiently small and potentially also decreasing $\delta$ further, the Lipschitzness of $f$ and the uniform continuity of $g$ combine to yield that $|(f\circ g)(y)|$ can be at most $\frac{\varepsilon}{2}$. The same reasoning shows that $|f(x)| \leq\frac{\varepsilon}{2}$ as well. It follows that 
$$
\left| \psi_\eta(x)f_\delta(x) - (f\circ g)(x) \right| \leq |f_\delta(x)| + |f(x)| \leq \varepsilon \; .
$$
On the other hand, if $d(x,X\setminus g^{-1}\supp(f)) > \eta$, then
$$
\left| \psi_\eta(x)f_\delta(x) - (f\circ g)(x) \right| = \left| f_\delta(x) - (f\circ g)(x) \right| \leq \varepsilon \; .
$$
We conclude that $f\circ g$ is within $\varepsilon$ of the function $\psi_\eta f_\delta$, which is in $L'\text{-Lip}_{R'}(X)$ for $L':=L(\eta)\cdot L''(\delta)$. This concludes the proof.
\end{proof}
It is important for technical purposes to remark that uniform K-homology classes over spaces of jointly bounded geometry can be represented by Fredholm modules with operators that have finite propagation\footnote{To get finite-propagation representatives it actually suffices to assume that $X$ admits a quasilattice. This is all that is needed for the construction of the partition of unity used below. Moreover, if for any $c>0$ there exists a quasilattice $\Gamma$ with $N_c(\Gamma)=X$, then the supports of the partition of unity can be made uniformly small. As a consequence any K-homology class can be represented by a Fredholm module with arbitrarily small propagation. This is the case for manifolds of bounded geometry for example, see Lemma \ref{CoverBdGeom} below.}. Indeed, using the bounded geometry of $X$ one can produce a partition of unity $\{\psi_i\}_i$ such that $\sqrt{\psi_i}$ is uniformly Lipschitz. Now, if $(H,\rho,F)$ is a uniform Fredholm module over $X$, then 
\begin{equation}\label{Truncation}
\left( H,\rho, \mathrm{trunc}(F):=\sum_i \rho(\sqrt{\psi_i})F\rho(\sqrt{\psi_i}) \right)
\end{equation}
is also a uniform Fredholm module over $X$, and it represents the same uniform K-homology class as $(H,\rho,F)$ (\cite[Section 7]{Spakula2010}). Moreover, the truncated Fredholm module \eqref{Truncation} can also be made non-degenerate\footnote{A Fredholm module $(H,\rho,F)$ is called \emph{non-degenerate} if $\rho(C_0(X))H=H$.} without loosing the finite-propagation property. Thus, if $X$ has jointly bounded geometry, then any uniform K-homology class over $X$ can be represented by a non-degenerate finite-propagation uniform Fredholm module.

We come to the definition of the Kasparov product. This requires usage of tensor products of (multi)graded spaces. The reader is directed to Appendix \ref{MultigradingAppendix} for the relevant definitions. As mentioned above, it also requires the assumption of jointly bounded geometry (the Kasparov product falls into the second category of uniform K-homology results discussed in the introduction to this chapter).
\begin{definition}\label{AlignedModules}
Let $X$, $X'$ be metric spaces. Let $(H,\rho,F)$ and $(H',\rho',F')$ be uniform Fredholm modules over $X$ and $X'$, respectively, the former being $p$-multigraded and the latter $p'$-multigraded for $p,p'\geq 0$. Assume both Fredholm modules to be non-degenerate and have finite propagation. Set $\tilde{H}:=H\hat{\otimes}H'$ and $\tilde{\rho}:=\rho\hat{\otimes}\rho'$. This is a $(p+p')$-multigraded representation of $X\times X'$. A $(p+p')$-multigraded uniform Fredholm module $(\tilde{H},\tilde{\rho},\tilde{F})$ over $X\times X'$ is said to be \emph{aligned} with $(H,\rho,F)$ and $(H',\rho',F')$ if:
\begin{itemize}
\item[(a)] $\tilde{F}$ has finite propagation,
\item[(b)] for all $f\in C_0(X\times X')$ the operators 
$$
\tilde{\rho}(f)\left( \tilde{F}(F\hat{\otimes}I)+(F\hat{\otimes}I)\tilde{F}\right)\tilde{\rho}(\bar{f})\; , \; \tilde{\rho}(f)\left( \tilde{F}(I\hat{\otimes}F')+(I\hat{\otimes}F')\tilde{F}\right)\tilde{\rho}(\bar{f})
$$
are positive modulo compact operators,
\item[(c)] for all $f\in C_0(X\times X')$ the operator $\tilde{\rho}(f)\tilde{F}$ derives $\fK(H)\hat{\otimes} \fB(H')$, i.e.
$$
\left[ \tilde{\rho}(f)\tilde{F} , \fK(H)\hat{\otimes} \fB(H') \right] \subseteq \fK(H)\hat{\otimes} \fB(H') \; .
$$
\end{itemize}
\end{definition}
We observe that apart from the finite-propagation requirement the definition is exactly the same as for the definition of the Kasparov product in analytic K-homology.
Under the assumption of jointly bounded geometry aligned uniform Fredholm modules can be shown to always exist, and be unique up to operator-homotopy. The construction uses a uniform version of Kasparov's technical theorem adapted to this context.
\begin{prop}[{\cite[Proposition 2.24]{Engel2019}}]
Suppose $X$ and $X'$ have jointly bounded geometry. Then, for any non-degenerate finite-propagation uniform Fredholm modules $(H,\rho,F)$ and $(H',\rho',F')$ over $X$ and $X'$ there exists uniform Fredholm module $(\tilde{H},\tilde{\rho},\tilde{F})$ aligned with $(H,\rho,F)$ and $(H',\rho',F')$. The aligned Fredholm module is unique up to operator-homotopy, and its operator-homotopy class depends only on the operator-homotopy classes of $(H,\rho,F)$ and $(H',\rho',F')$.
\end{prop}
The above proposition enables the definition of the Kasparov product on uniform K-homology. If $\alpha\in K_{-p}^{u}(X)$ and $\alpha'\in K^{u}_{-p'}(X')$ are represented by non-degenerate finite-propagation Fredholm modules $x$ and $x'$, respectively, with $p,p'\geq 0$, then $\alpha\times \alpha' \in K^{u}_{-(p+p')}(X\times X')$ is defined as the uniform K-homology class of any uniform Fredholm module aligned with $x$ and $x'$. For $p$ or $p'$ equal to $-1$ we can appeal to formal periodicity to define the Kasparov product.
\begin{thm}[{\cite[Theorem 2.26]{Engel2019}}]\label{KasparovProperties}
Let $X$, $X'$ be metric spaces of jointly bounded geometry. The Kasparov product
$$
K^{u}_{-p}(X) \times K^{u}_{-p'}(X') \longrightarrow K^{u}_{-(p+p')}(X\times X')
$$
is bilinear and associative. Moreover, it enjoys the following properties:
\item[(i)] The flip map $\tau: \, X\times X'\to X'\times X$ satisfies 
$$
\tau_*\left( \alpha \times \alpha' \right) = (-1)^{pp'} \alpha'\times \alpha \quad \forall \; \alpha\in K^{u}_{-p}(X), \, \alpha'\in K_{-p'}^{u}(X') \; .
$$ 
\item[(ii)] If $\phi: \, X'\to X''$ is a filtered map between jointly-bounded-geometry spaces, then
$$
(\id_X\times \phi)_*\left( \alpha \times \alpha' \right) =  \alpha\times \phi_*\alpha \quad \forall \; \alpha\in K^{u}_{-p}(X), \, \alpha'\in K_{-p'}^{u}(X') \; .
$$ 
\item[(iii)] Let $1\in K^{u}_0(\mathrm{pt})\cong\ZZ$ denote the generator defined by an index-$1$ operator. Then,
$$
1 \times \alpha = \alpha \quad \forall \; \alpha \in K^{u}_*(X)=K^{u}_*(\mathrm{pt}\times X)
$$
\end{thm}
\begin{remark}\label{PropagationKasparov}
Let us discuss the finite-propagation assumptions involved in the construction of the Kasparov product. We observe first that the truncation procedure of \eqref{Truncation} behaves well in regards to alignment of Fredholm modules. Indeed, if in the setting of Definition \ref{AlignedModules} $(\tilde{H},\tilde{\rho},\tilde{F})$ satisfies points (b) and (c), then it is not hard to see that $(\tilde{H},\tilde{\rho},\mathrm{trunc}(\tilde{F}))$ also satisfies these points, and it has finite propagation. Since truncation does not change the uniform K-homology class we conclude that
$$
\left[ H,\rho,F \right] \times\left[ H',\rho',F' \right] = \left[ \tilde{H},\tilde{\rho},\mathrm{trunc}(\tilde{F}) \right] = \left[ \tilde{H},\tilde{\rho},\tilde{F} \right] \; .
$$
In fact we can even drop the requirement that $(H,\rho,F)$ and $(H',\rho',F')$ have finite propagation. Again, it is not hard to see that $(\tilde{H},\tilde{\rho},\tilde{F})$ satisfies points (b) and (c) with respect to $F$ and $F'$ if and only if it satisfies them with respect to $\mathrm{trunc}(F)$ and $\mathrm{trunc}(F')$. Since truncation does not change the uniform K-homology class we conclude that in this situation it still holds that 
$$
\left[ H,\rho,F \right] \times\left[ H',\rho',F' \right] = \left[ \tilde{H},\tilde{\rho},\tilde{F} \right] \; .
$$
This relieves us from the burden of passing to finite-propagation representatives when working with the Kasparov products.
\end{remark}
Using the Kasparov product one can derive homotopy invariance of uniform K-homology. Two filtered maps $\phi_0,\phi_1: \, C_0(Y)\to C_0(X)$ from $X$ to $Y$ are called \emph{homotopic} if there exists a filtered map $\Phi: \, C_0(Y)\to C_0([0,1]\times X)$ such that $\phi_j=\mathrm{eval}_j\circ \Phi$, $j=0,1$. 
\begin{thm}[{\cite[Theorem 2.27]{Engel2019}}]\label{HomotopyInvariance}
If $\phi_0,\phi_1: \, X\to X'$ are homotopic filtered maps between metric spaces of jointly bounded geometry, then
$$
(\phi_0)_*=(\phi_1)_*: \; K^{u}_*(X) \longrightarrow K^{u}_*(X') \; .
$$
\end{thm} 
We note that Engel only proved this statement for filtered maps coming from uniformly cobounded Lipschitz maps $X\to X'$ that are homotopic via a uniformly cobounded Lipschitz map $[0,1]\times X\to X'$. However, all that is needed for the proof is that the maps and the homotopy have the property of being filtered.
\section{Relative uniform K-homology}\label{SecRelUKHom}
In this section we define relative uniform K-homology and derive some basic facts. Let $X$ be a metric space\footnote{Recall our convention that all metric spaces are assumed locally compact and second-countable.}, and $Z\subseteq X$ be a closed subspace.
\begin{definition}
Let $p\in\mathbb{Z}_{\geq -1}$. A \textit{$p$-multigraded relative uniform Fredholm module} over $(X,Z)$ is a triple $(H,\rho,F)$ consisting of
\begin{itemize}
\item a $p$-multigraded Hilbert space $H$, 
\item a representation $\rho: \, C_0(X)\to \fBH$ by even $p$-multigraded operators, 
\item and an odd $p$-multigraded operator $F\in \fBH$,
\end{itemize}
such that 
\begin{align*}
& \left\{ [F,\rho(f)] \, | \, f\in \LLip_R(X) \right\} \\
& \left\{ (F^2-I)\cdot\rho(g),\, \rho(g)\cdot (F^2-I) \, | \, g\in \LLip_R(X\setminus Z) \right\} \\
& \left\{ (F-F^*)\cdot\rho(g),\, \rho(g)\cdot (F-F^*) \, | \, g\in \LLip_R(X\setminus Z) \right\}
\end{align*}
are uniformly approximable families of compact operators for all $L,R\geq 0$, or, written succinctly,
\begin{align*}
[F,\rho] &\sua 0\quad \mathrm{over} \; X , \\
(I-F^2)\rho &\sua 0 \quad  \mathrm{over} \; X\setminus Z  ,\\
 (F-F^*)\rho &\sua 0 \quad \mathrm{over} \; X\setminus Z .
\end{align*} 
If $p=0$ or $p=-1$, we speak of \emph{graded} and \emph{ungraded} Fredholm modules respectively, and the conditions on the multigrading and in the latter case on the grading are void.
\end{definition}
Direct sums, unitary equivalences and operator-homotopies of relative uniform Fredholm modules are defined analogously to the non-relative case, and the definition of relative uniform K-homology is also analogous.
\begin{definition}
Let $p\in\mathbb{Z}_{\geq -1}$. The $p^\mathrm{th}$ \textit{relative uniform K-homology group} $K_p^{u}(X,Z)$ of the pair $(X,Z)$ is defined as the free abelian group generated by the unitary equivalence classes of $p$-multigraded relative uniform Fredholm modules over $(X,Z)$, modulo the relations
\begin{itemize}
\item $[x]=[y]$ if $x$ and $y$ are operator-homotopic,
\item $[x\oplus y]=[x]+[y]$.
\end{itemize}
\end{definition}
Let us turn to the functoriality of relative uniform K-homology.
\begin{definition}\label{DefFilteredMapsPairs}
A \emph{filtered map of pairs} $(X,Z)\to(Y,W)$ is a $^*$-homomorphism $\phi: \, C_0(Y) \to C_0(X)$ such that 
\begin{itemize}
\item $\phi\left( C_0(Y\setminus W) \right) \subseteq C_0(X\setminus Z)$.
\item For all $L,R\geq 0$ there exist $L',R'\geq 0$  such that $\phi(\LLip_R(Y))\subseteq L'\text{-}\mathrm{Lip}_{R'}(X)$.
\end{itemize}
\end{definition}
If $g: \, X\to Y$ is a uniformly cobounded Lipschitz map such that $g(Z)\subseteq W$, then $g$ induces a filtered map $(X,Z)\to (Y,W)$ in the above sense. 

Just as in the non-relative case one sees that if $\phi: \, (X,Z)\to (Y,W)$ is a filtered map of pairs, and $(H,\rho,F)$ is a relative uniform Fredholm module over $(X,Z)$, then $\phi_*(H,\rho,F):=(H,\rho\circ\phi,F)$ is a relative uniform Fredholm module over $(Y,W)$. This operation is again compatible with the relations imposed by relative uniform K-homology, and thus gives rise to functoriality for pairs.
\begin{prop}
Relative uniform K-homology is a functor from the category of metric pairs\footnote{By our convention we really mean the category of locally compact second-countable metric pairs.} and filtered maps of pairs to abelian groups.
\end{prop}
We present of series of elementary facts about relative uniform K-homology. Analogous statements are standard fare for analytic K-homology, and the proofs carry over essentially verbatim, see for example Chapter 8 of \cite{HigsonRoe2000}. We will therefore not reproduce these proofs here. 
\begin{lemma}\label{CptPerturb}
Let $(H,\rho,F)$ be a $p$-multigraded relative uniform Fredholm module for $(X,Z)$, and let $K\in \fBH$ be an odd $p$-multigraded operator\footnote{As usual the multigrading condition is void if $p=0$ or $p=-1$, with the condition of being odd also void in the latter case.} such that 
$$
[K,\rho] \sua 0 \quad \mathrm{over} \; X \quad , \quad K\rho \sua 0 \; \mathrm{over} \; X\setminus Z .
$$
Then, $(H,\rho,F+K)$ is also a $p$-multigraded relative uniform Fredholm module over $(X,Z)$, and
$$
\left[ (H,\rho,F) \right] \; = \; \left[ (H,\rho,F+K) \right] \quad \mathrm{in} \; K^{u}_{-p}(X,Z) \; .
$$
\end{lemma}
A relative uniform Fredholm module $(H,\rho,F)$ is called \emph{degenerate} if $[\rho,F]=(I-F^2)\rho=(F-F^*)\rho=0$. The usual Eilenberg swindle shows that the class of a degenerate Fredholm module in $K^{u}_*$ vanishes. Also recall that if $H$ is a graded space, we denote by $H^\mathrm{op}$ the same space endowed with the opposite grading. Then, if $T$ is a graded operator on $H$, we denote the same operator on $H^\mathrm{op}$ by $T^\mathrm{op}$.
\begin{lemma}\label{Inverses}
Let $(H,\rho,F)$ be a $p$-multigraded relative uniform Fredholm module. Then, $(H^\mathrm{op},\rho^\mathrm{op},-F^\mathrm{op})$ represents the inverse of $[(H,\rho,F)]$ in $K^{u}_{-p}(X,Z)$.
\end{lemma}
As a consequence of Lemma \ref{Inverses} each element in $K^{u}_*(X,Z)$ may be represented by a relative uniform Fredholm module.

Just as for non-relative uniform K-homology we shall find it useful to be able to represent relative uniform K-homology classes by finite-propagation operators. Recall from the discussion surrounding \eqref{Truncation} that if $X$ has jointly bounded geometry\footnote{Again, coarsely bounded geometry, i.e. the existence of a quasi-lattice in $X$ would actually suffice here. However, jointly bounded geometry, is essentially a standing assumption for the more non-trivial applications of (relative) uniform K-homology, and it certainly holds for manifolds of bounded geometry, so we will not worry too much about only making this weaker assumption here.}, then it admits a uniformly locally finite partition of unity $\{\psi_i\}$ with supports of uniformly bounded diameters such that $\sqrt{\psi_i}$ is uniformly Lipschitz. Just as in the non-relative case truncation using the $\psi_i$ produces finite-propagation representatives for relative uniform K-homology classes.
\begin{lemma}\label{RelTrunc}
Let $X$ be a metric space of jointly bounded geometry, and $Z$ a closed subspace. If $(H,\rho,F)$ is a relative uniform Fredholm module over $(X,Z)$, then so is 
$$
\left( H,\rho,\mathrm{trunc}(F)=\sum_i \rho(\sqrt{\psi_i})F\rho(\sqrt{\psi_i}) \right) \; ,
$$
and the two differ only by a uniformly compact perturbation. In particular, any class in $K^{u}_*(X,Z)$ has a finite-propagation representative.
\end{lemma}
\begin{proof}
The argument is completely analogous to that in the non-relative case, see Proposition 7.4 of \cite{Spakula2009}. We first show that $\mathrm{trunc}(F)$ is uniformly pseudo-local. Fix $L,R\geq 0$, and let $f\in\LLip_R(X)$. Then,
$$
[\rho(f),\mathrm{trunc}(F)] = \sum_i \rho(\sqrt{\psi_i})[\rho(f),F]\rho(\sqrt{\psi_i}) \; .
$$
Due to the jointly bounded geometry of $X$, the support of $f$ -which is bounded in diameter by $R$- intersects the supports of a uniformly finite number of $\sqrt{\psi_i}$. Hence, the right-hand side is a uniformly finite sum, and each summand is uniformly approximable, since $[\rho,F]\sua 0$.

Now we show that $(F-\mathrm{trunc}(F))\rho \sua 0$ over $X\setminus Z$. To that end let $g\in \LLip_R(X\setminus Z)$, with $L,R\geq 0$ again fixed. Using that $F$ is uniformly pseudo-local, calculate that
$$
\mathrm{trunc}(F)\rho(g) = \sum_i \rho(\sqrt{\psi_i}) F \rho(\sqrt{\psi_i}g) \sua \sum_i \rho(\psi_ig) F = \rho(g)F \sua F\rho(g) \; .
$$
Here we again used that the support of $g$ intersects the supports of a uniformly finite number of $\psi_i$. An appeal to Lemma \ref{CptPerturb} concludes the proof.
\end{proof}
A graded relative uniform Fredholm module $(H,\rho,F)$ is called \emph{balanced} if $H$ is the direct sum $H'\oplus H'$ of an ungraded Hilbert space $H'$, $H$ is graded by this direct sum decomposition, and $\rho=\rho'\oplus\rho'$ for a representation $\rho'$ on $H'$. A graded unitary between balanced Fredholm modules is called \emph{balanced} if it is a direct sum of a single unitary. \\
We define $K^{u,b}_0(X,Z)$ as the free group on balanced unitary equivalence classes of balanced relative uniform Fredholm modules over $(X,Z)$, modulo direct sums and operator-homotopies through balanced Fredholm modules. There is a natural map $K^{u,b}_0(X,Z) \to K^{u}_0(X,Z)$ induced by the identity on the level of Fredholm modules. The following proposition states that this map is an isomorphism. We say that the zeroth relative uniform K-homology can be \emph{normalized} by the requirement that all Fredholm modules be balanced. For Paschke duality it will be more convenient to work with this balanced normalization.
\begin{prop}\label{NormBal}
The identity on Fredholm modules induces a natural isomorphism
$$
K^{u,b}_0(X,Z) \xlongrightarrow{\sim} K^{u}_0(X,Z)  \; .
$$
\end{prop}
The statement of the next lemma sounds very natural, but is actually not without subtlety. Indeed, item (ii) can fail for non-balanced graded Fredholm modules. Issues arise from a potential lack of degenerate Fredholm modules with a given representation. We refer to the discussion around Proposition 8.3.14, as well as Exercise 8.8.17, in \cite{HigsonRoe2000}.
\begin{lemma}\label{UArep}
Let $\rho,\rho': \, C_0(X) \to\fBH$ be representations such that $\rho \sim_{ua} \rho'$. 
\begin{itemize}
\item[(i)] Let $(H,\rho,F)$ be an ungraded relative uniform Fredholm module. Then, $(H,\rho',F)$ is also a relative uniform Fredholm module, and 
$$
\left[ (H,\rho',F)\right] \; = \; \left[ (H,\rho,F)\right] \quad \in \; K^{u}_1(X,Z) \; . 
$$
\item[(ii)] Let $(H\oplus H,\rho\oplus\rho,F)$ be a balanced graded relative uniform Fredholm module. Then, $(H\oplus H,\rho'\oplus \rho',F)$ is also a relative uniform Fredholm module, and 
$$
\left[H\oplus H,\rho'\oplus \rho',F\right] \; = \; \left[ H\oplus H,\rho\oplus \rho,F\right] \quad \in \; K^{u,b}_0(X,Z) \; . 
$$
\end{itemize}
\end{lemma}
Lastly the formal twofold periodicity of multigraded spaces and operators provides a formal twofold periodicity on relative uniform K-homology.
\begin{prop}
Let $p\geq \ZZ_{\geq -1}$. There is a natural formal periodicity isomorphism
$$
K^{u}_{-(p+2)}(X,Z) \cong K^{u}_{-p}(X,Z) \; .
$$
\end{prop}
\section{Paschke duality}\label{SecPaschke}
In this section we prove Paschke duality for relative uniform K-homology. Paschke duality relates the K-homology groups $K^{u}_*(X,Z)$ to the operator K-theory of quotients of the so-called uniform dual algebra $\fD^{u}(X)$. We refer to Appendix \ref{AppendixKTheory} for some basics on K-theory for $C^*$-algebras.

Consider for the moment an absolute cycle $(H,\rho,F)$ over $X$ and assume that $F^2=I$ and $F^*=F$. First consider the case that the cycle is ungraded. The operator $F$ has spectrum contained in $\{\pm 1\}$, and we may consider the projection $P$ onto the $(+1)$-eigenspace, it is given by $P=\frac{1}{2}(F+1)$. Defining $\fD^{u}_\rho(X)$ as the algebra of those operators $T$ satisfying $[T,\rho]\sua 0$, we get that $P$ is a projection in the algebra $\fD^{u}_\rho(X)$ and thus defines an element in its zeroth K-theory.

If the cycle $(H,\rho,F)$ is graded and balanced instead, and we write $F=\begin{pmatrix}
0 & U^* \\ U &0
\end{pmatrix}$, then $U$ is a unitary, which also lies in the algebra $\fD^{u}_\rho(X)$. Thus $U$ defines an element of the first K-theory of $\fD^{u}_\rho(X)$.

Now consider $(H,\rho,F)$ to be a relative cycle. In the ungraded case $P=\frac{1}{2}(F+I)$ need not be a projection, but it satisfies $(P^2-P)\rho\sua 0$ and $(P-P^*)\rho\sua 0$ over $X\setminus Z$. Thus $P$ still defines a projection in the quotient $\fD^{u}_\rho(X)/\fD^{u}_\rho(X,Z)$, where $\fD^{u}_\rho(X,Z)$ is the subalgebra of $\fD^{u}_\rho(X)$ consisting of those operators $T$ for which $T\rho\sua 0$ over $X\setminus Z$. Thus $P$ defines an element of the zeroth K-theory of that quotient. Analogously in the graded case $U$ is a unitary in the quotient $\fD^{u}_\rho(X)/\fD^{u}_\rho(X,Z)$ and hence defines an element in its first K-theory. It is this assignment that will give us Paschke duality.

There is one technical obstacle to overcome: Cycles with different representations generally give rise to projections/unitaries in different dual algebras. Two different approaches exist to handle this issue in non-uniform K-homology. There exists an (essentially unique) representation of $C_0(X)$ which contain every other representation as sub-representations up to unitary equivalence and compact perturbation. This is a consequence of a theorem of Voiculescu\footnote{In more detail: A representation $\rho: \, A\to\fB(H^\rho)$ of a $C^*$-algebra $\fA$ is called \emph{ample} if the image of $\rho$ contains no compact operator besides $0$. Voiculescu's theorem states that if $\rho$ is ample and $\sigma: \, A\to\fB(H^\sigma)$ is any other representation, then there exists an isometry $V: \, H^\sigma\to H^\rho$ such that $V^*\rho V - \sigma$ is compact. It then follows easily that any two ample representations are unitarily equivalent up to a compact perturbation.}. Then any Fredholm module over $X$ is stably isomorphic to one over this universal representation. A uniform version of Voiculescu's theorem exists, but is valid only under the assumption of jointly bounded geometry.

The other approach taking the colimit of the K-theories of the dual algebra quotients over all representations instead. While a bit more involved this works in full generality also in the uniform setting, making this path preferable for our purposes. The colimit-approach to Paschke duality can be found for example in \cite{Higson1995}, and we will loosely follow this source.

The proof of Paschke duality falls into what we have called the first category of proofs in K-homology, in that it is entirely algebraic, relying only on matrix algebra and the formal property of equality up to (uniformly) compact error. Paschke duality for uniform K-homology is therefore proved in exactly the same way as it is for analytic K-homology. Again, we thus present a few proofs as demonstrations that things in uniform K-homology indeed work exactly as they do in analytic K-homology, and skip the rest. 
\subsection{Uniform covering isometries}
We begin by introducing a notion for comparing different representations in the uniform setting. This will allow us to take colimits over all representations later on. This is a modification of the usual concept of covering isometry, and was introduced in \cite{Spakula2009}. First, let us quickly introduce a small bit of notation: For a representation $\rho: \, C_0(X)\to \fBH$ we shall also write $(H,\rho)$,\footnote{Recall our convention that all Hilbert spaces are assumed to be separable.} and for the sake of brevity we shall also say that $(H,\rho)$ is a representation of $X$ (instead of $C_0(X)$).
\begin{definition}\label{DefUCI}
Let $X$ and $Y$ be metric spaces, and $\phi: \, C_0(Y) \to C_0(X)$ a $^*$-homomorphism. Let $(H^\rho,\rho)$ and $(H^\sigma,\sigma)$ be representations of $X$ and $Y$ respectively. We say that an isometry $V: \, H^\rho \to H^\sigma$ \emph{uniformly covers} $\phi$ if
$$
V^* \sigma(f) V \sua \rho(\phi(f))  \; .
$$
In this case we shall also write $V: \, (H^\rho,\rho)\to (H^\sigma,\sigma)$ to stress the representations involved.
\end{definition}
The first lemma states that if $V: \, (H^\rho,\rho)\to (H^\sigma,\sigma)$ is a uniform covering isometry, then up to a uniformly compact error $\sigma$ splits as the direct sum of $\rho$ and an action on the complement of $V(H^\rho)$.
\begin{lemma}[Compare {\cite[Lemma 3.1.6]{HigsonRoe2000}}]\label{UEDiag}
Let $X$ and $Y$ be metric spaces, $\phi: \, C_0(Y) \to C_0(X)$ a $^*$-homomorphism, and $(H^\rho,\rho)$ and $(H^\sigma,\sigma)$ representations of $X$ and $Y$ respectively. Let $V: \, (H^\rho,\rho) \to (H^\sigma,\sigma)$ be an isometry uniformly covering $\phi$. Let $P:= VV^*$ denote the orthogonal projection onto the image of $V$. Then,
$$
P\sigma(I-P), \; (I-P)\sigma P \; \sua \; 0 \; .
$$ 
\end{lemma}
\begin{proof}
Decompose $H^\sigma=PH^\sigma \oplus (I-P)H^\sigma$, and write 
$$
\sigma \; = \; \begin{pmatrix}
 \sigma_{11} & \sigma_{12} \\ \sigma_{21} & \sigma_{22}
\end{pmatrix}
$$
with respect to that decomposition. This means that $\sigma_{11}=P\sigma P$, $\sigma_{12}= P \sigma (I-P)$ and so on. Let $f\in C_0(X)$. Because $\sigma$ is a $^*$-homomorphism, it holds that
\begin{align}
&\begin{pmatrix}
\sigma_{11}(\bar{f}f) & \sigma_{12}(\bar{f}f) \\ \sigma_{21}(\bar{f}f) & \sigma_{22}(\bar{f}f)
\end{pmatrix} \; =  \nonumber \\
& \begin{pmatrix}
\sigma_{11}(f)^*\sigma_{11}(f)+\sigma_{12}(f)^*\sigma_{12}(f) & \sigma_{11}(f)^*\sigma_{21}(f)+\sigma_{12}(f)^*\sigma_{22}(f) \\ \sigma_{21}(f)^*\sigma_{11}(f)+\sigma_{22}(f)^*\sigma_{12}(f) & \sigma_{22}(f)^*\sigma_{22}(f)+\sigma_{21}(f)^*\sigma_{21}(f)
\end{pmatrix} \; . \label{EqHomom}
\end{align}
Observe that $\sigma_{11}=P\sigma P \sua V(\rho\circ\phi) V^*$. It follows that\footnote{If $f\in\LLip_R(X)$, then $\bar{f}f\in (2L)\text{-}\mathrm{Lip}_R(X)$. The first relation therefore follows from $\sigma_{11}\sua V(\rho\circ\phi) V^*$. The last relation holds because if $\varphi_1 \sua \psi_1$ and $\varphi_2\sua\psi_2$, then $\varphi_1\varphi_2 \sua \psi_1\psi_2$, as can be shown with the help of Lemma \ref{UAIdeal} (iv).}
$$
\sigma_{11}(\bar{f}f) \; \sua \; V \rho(\phi(\bar{f})\phi(f)) V^* \; = \; (V\rho(\phi(f))V^*)^*(V\rho(\phi(f))V^*) \; \sua \; \sigma_{11}(f)^* \sigma_{11}(f) \; .
$$
The upper-left corner of \eqref{EqHomom} therefore reduces to $\sigma_{12}(f)^*\sigma_{12}(f) \sua 0$. It follows that $\sigma_{12} \sua0$ as well\footnote{We note in Appendix \ref{AppendixUA} that the set of $\sigma$ with $\sigma\sua 0$ forms a closed two-sided $^*$-ideal in the $C^*$-algebra of bounded maps from the unit ball in $C_0(X)$ to $\fBH$. It is a general fact that if $J$ is such an ideal in a $C^*$-algebra $A$, and $a\in A$ is such that $a^*a\in  J$, then $a\in J$.}, hence also $\sigma_{21}(f)=\sigma_{12}(\bar{f})^*\sua0$. 
\end{proof}
We can use the uniform covering isometry $V$ to view $H^\rho$ as a subspace of $H^\sigma$. The previous lemma states that $\sigma \sua (\rho\circ\phi) \oplus \sigma_{22}$, where $\sigma_{22}$ acts on the orthogonal complement of $H^\rho$. Assume that $\phi=\id_X$, and let $(H^\rho,\rho,F)$ be a uniform Fredholm module. Then add a degenerate Fredholm module on the orthogonal complement of $H^\rho$, for example $\pm I$. Since $\sigma$ agrees with $\rho\oplus\sigma_{22}$ up to a uniformly compact perturbation, this yields a uniform Fredholm module on $(H^\sigma,\sigma)$, and since it differs from $(H^\rho,\rho,F)$ only by a degenerate module, it will have the same uniform K-homology class. The following lemma makes this argument precise.
\begin{lemma}\label{UCI_Khom}
Let $(H^\rho,\rho)$ and $(H^\sigma,\sigma)$ be representations of $X$, and suppose that $V: \,( H^\rho,\rho) \to (H^{\sigma},\sigma)$ is an isometry uniformly covering the identity of $X$. Let $\kappa\in\{\pm 1\}$.
\begin{itemize}
\item[(i)] Let $x=(H^\rho,\rho,F)$ be an ungraded relative uniform Fredholm module over $(X,Z)$. Then, 
$$
\Ad_\kappa(V)(x):=(H^\sigma,\sigma,VFV^*+\kappa(I-VV^*))
$$ 
is also an ungraded relative uniform Fredholm module over $(X,Z)$, and 
$$
[\Ad_\kappa(V)(x)]=[x] \in K^{u}_1(X,Z) \; .
$$
\item[(ii)] Let $x=(H^\rho\oplus H^\rho,\rho\oplus\rho,\begin{pmatrix}
0 & F_- \\ F_+ & 0
\end{pmatrix})$ be a balanced graded relative uniform Fredholm module over $(X,Z)$. Then, 
$$
\Ad_\kappa(V)(x):=\left(H^\sigma\oplus H^\sigma,\sigma\oplus\sigma,\begin{pmatrix}
0 & VF_-V^*+\kappa(I-VV^*) \\ VF_+V^*+\kappa(I-VV^*) & 0
\end{pmatrix}\right)
$$ 
is also a balanced graded relative uniform Fredholm module over $(X,Z)$, and 
$$
[\Ad_\kappa(V)(x)]=[x] \in K^{u}_0(X,Z) \; .
$$
\end{itemize}
\end{lemma}
\begin{proof}
We start with (i). Set $P:=VV^*$, the orthogonal projection in $H^\sigma$ onto the image of $V$. Lemma \ref{UEDiag} and the fact that $V$ is a uniform covering isometry imply that 
$$
\sigma \; \sua \; P \sigma P + (I-P)\sigma (I-P) \; \sua \; V \rho V^* + (I-P)\sigma (I-P) \; .
$$
Because $(I-P)V=0$ and $V^*(I-P)=0$, it follows for $f\in C_0(X)$ and $g\in C_0(X\setminus Z)$ that
\begin{align*}
[\sigma(f),VFV^*+\kappa(I-P)] \; &\sua \; V [\rho(f),F] V^* \\
&\sua \; 0 \\
\left( (VFV^*+\kappa(I-P))^2-I \right)\sigma(g) \; &= \; \left( V F^2 V^* -P \right) \sigma(g) \\
&\sua V F^2 \rho(g) V^* - V \rho(g) V^* \\
&\sua 0 \\
\left( (VFV^*+\kappa(I-P))^*-VFV^*-\kappa(I-P) \right)\sigma(g) \; &\sua V(F^*-F)\rho(g) V^* \\
&\sua 0 \; .
\end{align*}
Hence, $(H^\sigma,\sigma,VFV^*+\kappa(I-P))$ is a relative uniform Fredholm module. To see that it defines the same uniform K-homology class as $(H^\rho,\rho,F)$, consider the self-adjoint unitary
$$
U \; := \; \begin{pmatrix}
0 & V^* \\ V & I-P
\end{pmatrix}
$$
acting on $H^\rho\oplus H^\sigma$. It holds that 
$$
U\begin{pmatrix}
F & 0 \\ 0 & \kappa I
\end{pmatrix} U \; = \; \begin{pmatrix}
\kappa I&0 \\ 0 &VFV^*+\kappa(I-P)
\end{pmatrix} \; .
$$
Thus, in $K^{u}_1(X,Z)$ we have the equalities
\begin{align*}
\left[ H^\rho,\rho,F \right] \; &= \; \left[H^\rho\oplus H^\sigma,\rho \oplus \sigma,F\oplus\kappa I \right] \\
&= \; \left[ H^\rho\oplus H^\sigma,U(\rho \oplus \sigma)U,\kappa I\oplus \left(VFV^*+(I-P)\right)\right] \; .
\end{align*}
Now, compute that
$$
U^*(\rho\oplus\sigma)U \; = \; \begin{pmatrix}
V^*\sigma V & V^*\sigma (I-P) \\
(I-P)\sigma V & V\rho V^* + (I-P)\sigma (I-P)
\end{pmatrix} \; .
$$
It thus follows from Lemma \ref{UEDiag} and the fact that $V$ is a uniform covering isometry that $U^*(\rho\oplus\sigma)U \sua \rho\oplus\sigma$. With the help of Lemma \ref{UArep} we conclude that 
\begin{align*}
\left[ H^\rho,\rho,F \right] \; &= \;  \left[ H^\rho\oplus H^\sigma,U^*(\rho \oplus \sigma)U,\kappa I\oplus (VFV^*+(I-P))\right] \\
&= \; \left[ H^\rho\oplus H^\sigma,\rho \oplus \sigma,\kappa I\oplus (VFV^*+(I-P))\right] \\
&= \; \left[ H^\sigma,\sigma,VFV^*+(I-P)\right] 
\end{align*} 
in $K^{u}_1(X,Z)$. \\
The proof of (ii) is analogous. The balanced unitary $U\oplus U$ on $(H^\rho\oplus H^\sigma)\oplus(H^\rho \oplus H^\sigma)$, with $U$ as above, provides a unitary equivalence between
$$
\begin{pmatrix}
 & & F_- & 0 \\ & & 0 & \kappa I \\ F_+ & 0 && \\ 0 & \kappa I & & 
\end{pmatrix} \; \mathrm{and} \; \begin{pmatrix}
 & & \kappa I & 0 \\ & & 0 & VF_-V^*+\kappa(I-P) \\ \kappa I & 0 && \\ 0 & VF_+V^*+\kappa(I-P) & & 
\end{pmatrix} \; ,
$$
and $(U\oplus U)((\rho \oplus \sigma)\oplus(\rho\oplus\sigma))(U\oplus U)\sua (\rho \oplus \sigma)\oplus(\rho\oplus\sigma)$. Equality of the uniform K-homology classes again follows from Lemma \ref{UArep}.
\end{proof}
\subsection{Uniform dual algebras}\label{SecDualAlgebra}
Next we recall the definition of (relative) uniform dual algebras, which were introduced in \cite{Spakula2009}. Their properties and compatibility with uniform covering isometries were partially covered there already, and the parts that were not are straight-forward adaptations of the non-uniform counterparts. Throughout, let $X$ be a metric space, and $Z$ a closed subspace.
\begin{definition} 
Let $(H,\rho)$ be a representation of $X$. Define $\mathfrak{D}^{u}_\rho(X)$ as the subset of $\fB(H)$ consisting of all operators $T$ such that $[T,\rho]\sua 0$ over $X$. Furthermore, let $\fD^{u}_\rho(X,Z)$ denote the subset of $\fD^{u}_\rho(X)$ of all operators $T$ such that $T\rho, \, \rho T \sua 0$ over $X\setminus Z$.
\end{definition}
We call $\fD^{u}_\rho(X)$ the \emph{uniform dual algebra} and $\fD^{u}_\rho(X,Z)$ the \emph{relative uniform dual algebra} (associated to $\rho$). If $Z=\emptyset$, we will write $\fC^{u}_\rho(X)$ for $\fD^{u}_\rho(X,\emptyset)$, the algebra of uniformly locally compact operators. Both the uniform dual algebra and its relative version are $C^*$-algebra, with the latter being an ideal in the former.
\begin{prop}[{\cite[Lemma 4.2]{Spakula2009}}]
$\fD^{u}_\rho(X)$ is a $C^*$-subalgebra of $\fB(H)$, and $\fD^{u}_\rho(X,Z)$ is an ideal\footnote{By an ideal of a $C^*$-algebra we shall always mean a closed two-sided $^*$-ideal} of $\fD^{u}_\rho(X)$. 
\end{prop}
We can therefore form the quotient $\fQ_\rho^{u}(X,Z):=\fD^{u}_\rho(X)/\fD^{u}_\rho(X,Z)$, this is a well-defined $C^*$-algebra. Paschke duality states that the K-theory of this quotient algebra represents relative uniform K-homology, at least after taking the colimit over $\rho$. Again, we refer to Appendix \ref{AppendixKTheory} for the basics of operator K-theory. This colimit can be implemented using uniform covering isometries, as is stated by the following proposition.
\begin{prop}[{\cite[Lemmas 4.8, 5.4]{Spakula2009}}] \label{AdDual}
Let $(X,Z)\xrightarrow{\phi} (Y,W)$ be a filtered map of pairs, let $(H^\rho,\rho)$ and $(H^\sigma,\sigma)$ be representations of $X$ and $Y$ respectively, and let $V: \, (H^\rho,\rho) \to (H^\sigma,\sigma)$ be an isometry uniformly covering $\phi$. Then, there is a well-defined $^*$-homomorphism
$$
\mathrm{Ad}(V): \; \fD^{u}_\rho(X) \longrightarrow \fD^{u}_\sigma(Y) \; , \; T \longmapsto V T V^*
$$ 
such that $\mathrm{Ad}(V)\fD^{u}_\rho(X,Z) \subseteq \fD^{u}_\sigma(Y,W)$.

The induced homomorphisms 
\begin{align*}
\Ad(V)_*: \; &K_*\left( \fD^{u}_\rho(X) \right) \longrightarrow K_*\left( \fD^{u}_\sigma(Y) \right) \\
\Ad(V)_*: \; &K_*\left( \fD^{u}_\rho(X,Z) \right) \longrightarrow K_*\left( \fD^{u}_\sigma(Y,W) \right) \\
\Ad(V)_*: \; &K_*\left( \fQ^{u}_\rho(X,Z) \right) \longrightarrow K_*\left( \fQ^{u}_\sigma(Y,W) \right)
\end{align*}
are independent of $V$, in that any other isometry $(H^\rho,\rho)\to (H^\sigma,\sigma)$ uniformly covering $\phi$ induces the same homomorphisms.
\end{prop}
Define $\fP(X)$ to be the category whose objects are representations $(H^\rho,\rho)$ of $X$,\footnote{All Hilbert spaces are assumed separable.} and whose morphisms are isometries $(H^\rho,\rho)\to (H^\sigma,\sigma)$ uniformly covering the identity. This category is directed, because for any two representations $(H^\rho,\rho)$ and $(H^\sigma,\sigma)$, the inclusion into $(H^\rho\oplus H^\sigma,\rho\oplus\sigma)$ as the first and second direct summand respectively is a morphism in $\fP(X)$.

By Proposition \ref{AdDual} the assignments mapping $(H^\rho,\rho)$ to either of the algebras $\fD^{u}_\rho(X)$, $\fD^{u}_\rho(X,Z)$ or $\fQ^{u}_\rho(X,Z)$ are well-defined functors from $\fP(X)$ to the category of $C^*$-algebras and $^*$-homomorphisms. Arbitrary colimits exist in that category, so that we can define $\fD^{u}(X)$, $\fD^{u}(X,Z)$ and $\fQ^{u}(X,Z)$ as the colimits of the respective functors.

We will be interested in the K-theory of these colimits. Operator K-theory commutes with colimits, so that
\begin{align*}
K_*(\fD^{u}(X)) &\cong \colim \; K_*(\fD_\rho^{u}(X)) \; , \\ 
K_*(\fD^{u}(X,Z)) &\cong \colim \; K_*(\fD_\rho^{u}(X,Z)) \; , \\
K_*(\fQ^{u}(X,Z)) &\cong \colim \; K_*(\fQ_\rho^{u}(X,Z)) \; ,
\end{align*}
where all colimits are taken over $\fP(X)$. If one wishes to avoid colimits of $C^*$-algebras, one could also take the respective right-hand sides as definitions for $K_*(\fD^{u}(X))$, $K_*(\fD^{u}(X,Z))$ and $K_*(\fQ^{u}(X,Z))$. These colimits of K-theory groups have a somewhat simpler structure, because all uniform covering isometries between fixed representations induce the same homomorphism on K-theory. The colimit over $\fP(X)$ of the K-theory groups therefore factors as the colimit over the poset whose elements are representations $(H^\rho,\rho)$, and where $\rho\leq\sigma$ holds if there is a uniform covering isometry $(H^\rho,\rho)\to(H^\sigma,\sigma)$.

Suppose that $X$ has jointly bounded geometry, and that $(H^\sigma,\sigma)$ is an ample representation.\footnote{Recall that $\sigma$ is called ample if $\sigma(C_0(X))H^\sigma=H^\sigma$, i.e. $\sigma$ is non-degenerate, and $\sigma(f)$ is compact only if $f=0$.} Then, the uniform version of Voiculescu's theorem proved by \v{S}pakula implies that for any other representation $(H^\rho,\rho)$ there exists a uniform covering isometry $(H^\rho,\rho)\to (H^\sigma,\sigma)$ (\cite[Corollary 3.3]{Spakula2010}).\footnote{The uniform Voiculescu theorem stated in \cite{Spakula2010} requires slightly different properties of the space $X$. See the discussion in \cite[Section 4.5]{Engel2014} for the observation that jointly bounded geometry is required instead of locally and coarsely bounded geometry separately, and that properness of the underlying space, which was required in \cite{Spakula2010}, is not essential.} While this does not quite make $\fD^{u}_\sigma(X)$ a cofinal object of the diagram $(H,\rho) \mapsto \fD^{u}_\rho(X)$, since there may be many such uniform covering isometries, this is no issue on the level of K-theory. Thus the K-theory of the uniform dual algebra associated to an ample representation is cofinal for the K-theories of uniform dual algebras associated to any representation, same for the relative uniform dual algebras and their quotients. Thus we conclude:
\begin{prop}
Let $X$ have jointly bounded geometry, and let $(H^\sigma,\sigma)$ be an ample representation of $X$. Then, 
$$
K_*(\fD^{u}(X)) \cong  K_*(\fD_\sigma^{u}(X)) \; , \; K_*(\fD^{u}(X,Z)) \cong  K_*(\fD_\sigma^{u}(X,Z)) \; , \; K_*(\fQ^{u}(X,Z)) \cong  K_*(\fQ_\sigma^{u}(X,Z)) \; .
$$
\end{prop}
The K-theory of the dual algebra of a space relative to the empty set vanishes, and the same holds in the uniform setting. Viewed through the lense of Paschke duality this expresses that the uniform K-homology of a space $X$ and that of $X$ relative to $\emptyset$ agree. Recall that we denote $\fD^{u}(X,\emptyset)$ by $\fC^{u}(X)$.
\begin{prop}[Compare {\cite[Lemma 1.6]{Higson1995}}]\label{RelEmptyset}
It holds that $K_*(\fC^{u}(X))=0$. Therefore the quotient map induces a natural isomorphism
$$
K_* \left( \fD^{u}(X) \right)\cong K_*\left( \fQ^{u}(X,\emptyset) \right) \; .
$$
\end{prop}
Lastly we turn to the functoriality in the underlying pair. This is essentially covered by Proposition \ref{AdDual}. In fact, if $\phi: \, (X,Z)\to (Y,W)$ is a filtered map of pairs and $(H^\rho,\rho)$ a representation of $X$, then $(H^\rho,\rho\circ\phi)$ is a representation of $Y$, and $\id_{H^\rho}: \, (H^\rho,\rho)\to (H^\rho,\rho\circ\phi)$ uniformly covers $\phi$. Thus we have inclusions $\fD^{u}_\rho(X)\subseteq \fD^{u}_{\rho\circ\phi}(Y)$ and $\fD^{u}_\rho(X,Z)\subseteq \fD^{u}_{\rho\circ\phi}(Y,W)$, and therefore a canonical map $\fQ^{u}_\rho(X,Z)\to \fQ^{u}_{\rho\circ\phi}(Y,W)$. Moreover, if $V: \, (H^\rho,\rho)\to (H^\sigma,\sigma)$ uniformly covers $\id_{X}$, then $V: \, (H^\rho,\rho\circ\phi)\to (H^\sigma,\sigma\circ\phi)$ uniformly covers $\id_Y$. These data therefore assemble into well-defined maps $\fD^{u}(X)\to \fD^{u}(Y)$, $\fD^{u}(X,Z)\to \fD^{u}(Y,W)$ and $\fQ^{u}(X,Z)\to \fQ^{u}(Y,W)$ in a functorial way. Passing to K-theory lets us conclude:
\begin{prop}
The assignments
$$
X \mapsto K_*\left( \fD^{u}(X) \right) \; , \; (X,Z) \mapsto K_*\left( \fD^{u}(X,Z) \right) \; , \; (X,Z) \mapsto K_*\left( \fQ^{u}(X,Z) \right)
$$
are functors from the category of metric spaces\footnote{recall that we really mean the category of second-countable locally compact metric spaces.} and filtered maps to abelian groups in the first case, and from the category of pairs and filtered maps of pairs to abelian groups in the latter cases.
\end{prop}
\subsection{Paschke duality for relative uniform K-homology}
The proof of Paschke duality for relative uniform K-homology is the same calculation as for analytic K-homology. We will therefore only indicate the steps involved instead of reproducing the entire proof.

We begin with the definition of the map $K_{*+1}(\fQ^{u}(X,Z))\to K_*^{u}(X,Z)$. To that end fix a representation $(H^\rho,\rho)$ of $X$. Let $p \in \fQ^{u}_\rho(X,Z)=\fD^{u}_\rho(X)/\fD^{u}_\rho(X,Z)$ be a projection, and choose a representative $P\in \fD^{u}_\rho(X)$ of $p$. Then, $P$ being a projection modulo $\fD^{u}_\rho(X,Z)$ translates into $(H^\rho,\rho,2P-I)$ being an ungraded relative uniform Fredholm module over $(X,Z)$. Denote its relative uniform K-homology class by $\varphi_1^\rho(p)$. The class $\varphi_1^\rho$ can be defined equally well for projections in a matrix algebra over $\fQ^{u}_\rho(X,Z)$.

Had we chosen a different representative $P'$ of $p$, then $P-P'$ would lie in $\fD^{u}_\rho(X,Z)$, meaning that the Fredholm modules $(H^\rho,\rho,2P-I)$ and $(H^\rho,\rho,2P'-I)$ differ only by a uniformly compact perturbation. By Lemma \ref{CptPerturb} the class $\varphi_1^\rho(p)$ is therefore independent of the choice of representative of $p$. Moreover, homotopies of $p$ can be lifted to homotopies of representatives, yielding operator-homotopies of the corresponding Fredholm modules. Moreover, the assignment $\varphi_1^\rho$ is compatible with direct sums, and sends the zero-projection to the class of a degenerate module, i.e. to zero. We conclude that $\varphi_1^\rho$ descends to a well-defined homomorphism
$$
\varphi_1^\rho: \, K_0\left( \fQ^{u}_\rho(X,Z) \right) \longrightarrow K_1^{u}\left( X,Z\right) \; .
$$
We can do the same with degree-1 K-theory: Let $u\in \fQ^{u}_\rho(X,Z)$ be a unitary, and $U\in \fD^{u}_\rho(X)$ a representative. Then, $\left( H^\rho\oplus H^\rho,\rho\oplus\rho, \begin{pmatrix}
0 & U^* \\ U & 0
\end{pmatrix} \right)$ is a balanced graded relative uniform Fredholm module over $(X,Z)$. Denote its class by $\varphi_0^\rho(u) \in K_1(\fQ^{u}_\rho(X,Z))$. It can be defined analogously for unitaries in a matrix algebra over $\fQ^{u}_\rho(X,Z)$. Different choices of representative again yield uniformly compact perturbations, homotopies yield operator-homotopies, the identity is mapped to zero, and the mapping respects direct sums. We therefore conclude that $\varphi_0^\rho$ descends to a well-defined homomorphism
$$
\varphi_0^\rho: \, K_1\left( \fQ^{u}_\rho(X,Z) \right) \longrightarrow K_0^{u}\left( X,Z\right) \; .
$$
We wish to argue that the maps $\varphi_*^\rho$ give rise to maps on $K_*(\fQ^{u}(X,Z))$. To that end suppose that $V: \, (H^\rho,\rho)\to (H^\sigma,\sigma)$ is a uniform covering isometry between representations of $X$. Using the definitions it is straight-forwardly calculated that
$$
\Ad_{(-1)^p}(V)\circ \varphi^\rho_p = \varphi_p^\sigma \circ \Ad(V)
$$
holds on the level of Fredholm modules. We recall that $\Ad_{(-1)^p}(V)$ is a suitable conjugation by $V$ on Fredholm modules, as defined in Lemma \ref{UCI_Khom}. There it is also shown that $[\Ad_{(-1)^p}(V)x]=[x]$ for any $[x]\in K_*^{u}(X,Z)$. We conclude that
$$
\varphi_*^\sigma\left( \Ad(V)a\right) = \varphi^\rho_* (a) \; \in K_*^{u}(X,Z) \; ,
$$ 
where $a\in \fQ^{u}_\rho(X,Z)$ is either a projection or a unitary, depending on the degree. This lets us conclude that the maps $\varphi_p^\rho$ assemble into a map $\varphi_p: \, K_{p+1}(\fQ^{u}(X,Z))\to K^{u}_p(X,Z)$ on the colimit. We also note that this map is natural with respect to filtered maps of pairs, given that functoriality for both (relative) uniform dual algebras and uniform K-homology is obtained by pre-composition. Thus we have the following:
\begin{lemma}
The maps $\varphi_p$, $p=0,1$, constructed above give rise to natural homomorphisms
$$
K_{1}\left( \fQ^{u}(X,Z) \right) \xlongrightarrow{\varphi_0} K_0^{u}\left( X,Z\right) \quad , \quad K_{0}\left( \fQ^{u}(X,Z) \right) \xlongrightarrow{\varphi_1} K_1^{u}\left( X,Z\right) \; .
$$
\end{lemma} 
The homomorphisms $\varphi_p$ implement Paschke duality. We describe the construction of their inverses. Let $x=(H^\rho,\rho,F)$ be an ungraded relative uniform Fredholm module over $(X,Z)$. Then, $P:=\frac{1}{2}(F+I)$ defines a projection in $\fQ^{u}_\rho(X,Z)$. Denote its K-theory class in the colimit by $\beta_1(x)\in K_0(\fQ^{u}(X,Z))$. The mapping $\beta_1$ takes operator-homotopies to homotopies of projections, and it respects unitary equivalences and direct sums. Therefore it gives rise to a well-defined homomorphism
$$
\beta_1: \; K_1^{u}(X,Z) \longrightarrow K_{0}\left( \fQ^{u}(X,Z) \right) \; .
$$
Analogously, if $x=(H^\rho\oplus H^\rho,\rho\oplus\rho, \begin{pmatrix}
0 & F_- \\ F_+ & 0
\end{pmatrix})$ is a balanced graded relative uniform Fredholm module over $(X,Z)$, then $F_+$ defines a unitary in $\fQ^{u}_\rho(X,Z)$. Denote its K-theory class in the colimit by $\beta_0(x)\in K_1(\fQ^{u}(X,Z))$. We note that this construction really requires the Fredholm module to be balanced. Again, the mapping $\beta_0$ descends to a well-defined homomorphism
$$
\beta_0: \; K_0^{u,b}(X,Z) \longrightarrow K_1\left( \fQ^{u}(X,Z) \right) \; .
$$
The maps $\beta_p$ and $\varphi_p$ are mutual inverses. This is essentially by definition, though for $p=0$ ever so slight attention has to be paid to the formal difference between $K_0^{u}(X,Z)$ and $K^{u,b}_0(X,Z)$. This is covered by Lemma \ref{NormBal} though. Thus we have proved Paschke duality for relative uniform K-homology: 
\begin{thm}
The natural homomorphisms 
$$
K^{u}_0(X,Z) \stackrel{\xrightharpoonup{\quad \beta_0 \quad}}{\xleftharpoondown[\quad\varphi_0\quad]{}} K_{1}\left( \fQ^{u}(X,Z) \right) \quad , \quad K^{u}_1(X,Z) \stackrel{\xrightharpoonup{\quad \beta_1 \quad}}{\xleftharpoondown[\quad\varphi_1\quad]{}} K_{0}\left( \fQ^{u}(X,Z) \right)
$$
are mutually inverse isomorphisms. In particular,
$$
K^{u}_0(X,Z) \cong K_{1}\left( \fQ^{u}(X,Z) \right) \quad \; \quad K^{u}_1(X,Z) \cong K_{0}\left( \fQ^{u}(X,Z) \right) \; .
$$
\end{thm}
From twofold periodicity we conclude the following.
\begin{cor}
Let $p\in \ZZ_{\geq -1}$. There is a natural isomorphism
$$
K^{u}_{-p}(X,Z) \cong K_{(p+1) \, \mathrm{mod}\, 2}\left(\fQ^{u}(X,Z) \right) \; .
$$
\end{cor}
Using the K-theory in degree $(p+1) \, \mathrm{mod}\, 2$ on the right-hand side of this isomorphism is technically correct, since that is what we have proved. It is however cumbersome from a notational standpoint. Using Bott periodicity for K-theory we also get an isomorhism
$$
K^{u}_{-p}(X) \cong K_{-p+1}\left( \fQ^{u}(X,Z) \right)  \; .
$$
In the rest of this thesis we will use this degree convention for all things related to Paschke duality. 
\section{Excision and the long exact sequence}
In this section we prove the core theorems relating relative and absolute uniform K-homology: excision and the long exact sequence. As eluded to in the introduction to this chapter these theorems fall into the second category of K-homology results, meaning that they require non-trivial geometric/analytic input. In order to make the proofs work in (relative) uniform K-homology we must make geometric assumptions on the underlying pair. We discuss these first.
\begin{definition}\label{ConditionGeomSubspace}
Let $X$ be a metric space, and $Z\subseteq X$ a closed subspace. We say that $X$ is \emph{uniformly resolvable} near $Z$ if there is a partition of unity $\{\psi_i\}_i$ on $X\setminus Z$ with the following properties:
\begin{itemize}
\item[(1)] For all $\delta>0$ there exists $L(\delta)\geq 0$ such that $\sqrt{\psi_i}$ is $L(\delta)$-Lipschitz if $d(\supp(\psi_i),Z)\geq\delta$.
\item[(2)] For all $R\geq 0$ and $\delta>0$, there exists $M(R,\delta)\in\NN$ such that for all compact $K\subseteq X\setminus Z$ with $\diam(K)\leq R$ and $d(K,Z)\geq\delta$ it holds that
$$
 \left| \left\{ i \, | \, \psi_i|_K \neq 0 \right\}\right| \leq M(R,\delta) \; .
$$
\item[(3)] For all $\varepsilon>0$ there exists $\delta>0$ such that for all $i$ it holds that
$$
\diam(\supp(\psi_i)) \geq \varepsilon \; \Rightarrow \; d(\supp(\psi_i),Z) \geq \delta \; .
$$
\end{itemize}
Furthermore we say that $X$ is \emph{transversally controlled} near $Z$ if the restriction $C_0^+(X)\to C_0^+(Z)$ admits an approximately filtered completely positive section\footnote{Any positive map $\sigma: \, C_0^+(X) \to \fBH$ is already completely positive. Thus we will also simply speak of positive section.}, meaning that there exists a completely positive map $\sigma: \, C_0^+(Z)\to C_0^+(X)$ such that $\sigma(f)|_Z=f$ for all $f\in C_0^+(Z)$, and for all $L,R\geq 0$ and $\varepsilon>0$ there exist $L',R'\geq 0$ such that $\sigma\big(\LLip_R(Z)\big)$ is contained in the $\varepsilon$-neighborhood of $L'\text{-}\mathrm{Lip}_{R'}(X)$. Here $C^+_0(Z)$ and $C^+_0(X)$ denote the algebras of functions that are constant at infinity.
\end{definition}
While the author was tempted to call a pair $(X,Z)$ with these properties a pair of bounded geometry, he ultimately opted for a term that is less over-loaded and hopefully more descriptive. Indeed, one could think of a partition of unity as in the definition as providing a resolution of the geometry of $X$ in terms of bounded pieces that shrink uniformly as one approaches $Z$, hence the name \emph{uniformly resolvable}. Moreover, we may think of a completely positive section $C_0^+(Z)\to C_0^+(X)$ as providing extensions of functions on $Z$ to the direction \glqq transversal to $Z$\grqq\, to obtain functions on $X$. This extension map being approximately filtered expresses that we have (approximate) control over Lipschitzness and the size of supports of the extensions, hence the name \emph{transversally controlled}. Of course this is only a vague intuition, since we cannot speak of transversal directions in general metric spaces.

These definitions are technical wishlists, they exactly provide the tools needed to adapt the proofs of excision and the long exact sequence from analytic K-homology to the uniform setting. Fortunately all pairs relevant for the purposes of this thesis are covered by the definitions.
\begin{lemma}
Let $(X,Z)$ be one of the following:
\begin{itemize}
\item $X$ is a manifold of bounded geometry, and $Z$ a bounded-geometry hypersurface (see Definition \ref{BGHypersurface}). This includes the case that $Z$ is the boundary of a manifold $X$ with boundary and bounded geometry.
\item $X$ is a simplicial complex of bounded geometry, and $Z$ is a subcomplex.
\item $Z$ has jointly bounded geometry, and $X$ is the cone $CZ=[0,1)\times Z$, which contains $Z$ as $Z\times\{0\} \subset CZ$.
\end{itemize}
Then $X$ is uniformly resolvable and transversally controlled near $Z$.

Moreover, if $X$ has jointly bounded geometry, then it is uniformly resolvable near any closed subspace $Z\subseteq X$.
\end{lemma}
The proof along with more general sufficient conditions can be found in Appendix \ref{AppendixUnifPairs}. There it is also shown that these really are conditions on the geometry near $Z$, in that if $U$ is a uniformly thick neighborhood of $Z$ in $X$ such that $U$ is uniformly resolvable/transversally controlled near $Z$, then so is $X$.
\subsection{Excision}
We begin with our proof of the excision isomorphism. We adapt the arguments given in Chapter 5 of \cite{HigsonRoe2000}. The proof relies on a uniform version of the Kasparov lemma, which gives an equivalent condition for pseudo-locality. In the uniform setting it was already discussed by Engel for proper metric spaces \cite[Lemma 4.45]{Engel2014}. However, the spaces we are interested in here, namely complements of closed subspaces, are rarely proper, requiring us to prove a version also applicable to our case.

Let $X$ be a metric space. Denote by $C_0^+(X)$ the $C^*$-algebra of continuous functions on $X$ which are constant at infinity. It equals the unitization of $C_0(X)$, every function $f\in C_0^+(X)$ can be written uniquely as $f=\lambda \cdot \mathbf{1}_X + g$, where $\lambda:=\lim_{x\to\infty}f(x)$ and $g\in C_0(X)$. Moreover, denote the $C^*$-algebra of bounded measurable functions on $X$ by $\mathcal{L}^\infty(X)$. Any representation $(H,\rho)$ of $C_0(X)$ extends to a representation of $\mathcal{L}^\infty(X)$. We will denote this extension by $\rho$ as well.
\begin{prop}\label{KasparovLemma}
Let $(H,\rho)$ be a representation of $X$, and let $T\in \fBH$. Then, the following are equivalent: 
\begin{itemize}
\item[(i)] $[\rho(f),T]\sua 0$ for $f\in C_0(X)$.
\item[(ii)] For all $f\in C_0^+(X)$ it holds that $\rho(f)T\rho(g)\sua 0$ and $\rho(g)T\rho(f)\sua0$ for $g\in C_0(X\setminus \supp(f))$.
\item[(iii)] For all $f\in C_b(X)$ it holds that $\rho(f)T\rho(g)\sua 0$ and $\rho(g)T\rho(f)\sua0$ for $g\in C_0(X\setminus \supp(f))$.
\item[(iv)] For all $f\in \mathcal{L}^\infty(X)$ it holds that $\rho(f)T\rho(g)\sua 0$ and $\rho(g)T\rho(f)\sua0$ for $g\in C_0(X\setminus \supp(f))$.
\end{itemize}
\end{prop}
\begin{proof}
We start by proving (i)$\Rightarrow$(iv). Suppose $[T,\rho(g)]\sua 0$ for $g\in C_0(X)$. Then for any $f\in \mathcal{L}^\infty(X)$, we have $\rho(f)T\rho(g) \sua \rho(f)\rho(g)T=0$ for $g\in C_0(X\setminus\supp(f))$, proving this implication. The implications (iv)$\Rightarrow$(iii)$\Rightarrow$(ii) are obvious.

It remains to prove (ii)$\Rightarrow$(i). Suppose that for all $f\in C_0^+(X)$ it holds that $\rho(f)T\rho(g)\sua 0$ for $g\in C_0(X\setminus \supp(f))$. It suffices to prove that $[T,\rho(g)]\sua0$ for self-adjoint $g\in C_0(X)$. Fix $\varepsilon>0$, $L,R\geq 0$, and let $g\in \LLip_R(X)$ be self-adjoint. Construct a partition of the image $g(X)\subseteq [-1,1]$ into Borel sets $U_1,\cdots,U_n$ as follows: Let $\tilde{U}_1$ be the half-open interval $(-\frac{\varepsilon}{4},\frac{\varepsilon}{4}]$. Let $\tilde{U}_2, \cdots, \tilde{U}_m$ be a partition of $[-1,1]\setminus\tilde{U}_1$ into half-open intervals $\tilde{U}_j$ with $\diam(\tilde{U}_j)<\varepsilon$. The $\tilde{U}_j$ can be chosen such that if the closures of $\tilde{U}_j$ and $\tilde{U}_i$ don't intersect, then the distance between $\tilde{U}_j$ and $\tilde{U}_i$ is $\geq \varepsilon/2$. Set $U_j:=\tilde{U}_j\cap g(X)$. By discarding empty intersections and relabeling the index $j$ in order of ascending supremum of $U_j$, we obtain a partition of $g(X)$ into non-empty Borel sets $U_1,\cdots,U_n$ with the properties that $\diam(U_j)<\varepsilon$, and $d(U_j,U_i)\geq\varepsilon/2$ if $|j-i|> 1$. Moreover, the cardinality $n$ of the partition is bounded by $C/\varepsilon$ for some universal constant $C>0$. 

We denote by $g_j\in \cL^\infty(X)$ the characteristic function of $g^{-1}(U_j)$. Fix $x_j\in g^{-1}(U_j)$, and set $\tilde{g}:= g(x_1)\cdot g_1+\cdots + g(x_n)\cdot g_n$, $\tilde{g}\in\cL^\infty(X)$. By construction, we have $||g-\tilde{g}||_\infty<\varepsilon$, from which it follows that $||[T,\rho(g)]-[T,\rho(\tilde{g})]||<2||T||\cdot \varepsilon$. Suppose $|i-j|>1$. From our assumptions on $g$ and $U_j,U_i$ it can be deduced that the distance of $g^{-1}(U_j)$ and $g^{-1}(U_i)$ is at least $\varepsilon/2L$. Set $L':=8L/\varepsilon$. Then, for all $j$, define the function
$$
\tilde{g}_j(x) \; := \; \min \left\{ 1, \; L'\cdot d\left(x, X\setminus U_{1/L'}(g^{-1}(U_j)) \right) \right\} \; .
$$
Here $U_c(A)$ denotes the open $c$-neighborhood of a subset $A\subseteq X$. These functions have the following properties:
\begin{itemize}
\item[(a)] $\tilde{g}_j\equiv 1$ on $\supp(g_j)$, and $||\tilde{g}_j||_\infty\leq 1$,
\item[(b)] $\tilde{g}_j$ is $L'$-Lipschitz,
\item[(c)] If $0\notin U_j$, then $\tilde{g}_j$ is compactly supported and $\diam(\supp(\tilde{g}_j))\leq R':= R+2/L'$. If $0\in U_j$, then $\tilde{g}_j$ is constant at infinity.
\item[(d)] the supports of $\tilde{g}_j$ and $\tilde{g}_i$ are disjoint if $|j-i|>1$. 
\end{itemize}
Properties (c) and (d) deserve justification. We start with the observation that if $x\in \supp(\tilde{g}_j)$, then $d(x,g^{-1}(U_j)) \leq 1/L'$, whence it follows from the Lipschitz property of $g$ that $d(g(x),U_j) \leq L/L'=\varepsilon/8$. Now, the function $g\in\LLip_R(X)$ is compactly supported. For $\delta>0$, the set $\{ x\, | \, |g(x)|\geq\delta \}\subseteq X$ is therefore also compact, and, because of the above observation, contains the support of $\tilde{g}_j$ for $\delta$ sufficiently small and $0\notin U_j$. This proves the compactness of the support of these $\tilde{g}_j$, and the diameter bound is clear from the definition of $\tilde{g}_j$. On the other hand, if $0\in U_j$, then $\tilde{g}_j \equiv 1$ outside the compact support of $g$, which immediately implies $\tilde{g}_j\in C_0^+(X)$. This proves (c). Property (d) follows from the observation above, together with $d(U_j,U_i)\geq \varepsilon/2$. \\
The properties (a)-(d) in particular imply that $\tilde{g}_j\in L'\text{-}\mathrm{Lip}_{R+2/L'}(X)$ for all except at most one $j$, and that for that exception, it still holds that $\tilde{g}_j \in C_0^+(X)$. Using the assumption on $T$, we deduce from (d) that there exists $N=N(L,R,\varepsilon)$ such that for $|i-j|>1$ the operator\footnote{Note that since $|i-j|>1$, at least one of the sets $U_i$ and $U_j$ does not contain $0$. Thus, $\tilde{g}_i$ or $\tilde{g}_i$ is in $L'\text{-}\mathrm{Lip}_{R+2/L'}(X)$, and in particular $\rho(\tilde{g}_j)=\rho(\tilde{g}_j)$ or $\rho(\tilde{g}_i)=\rho(\tilde{g}_i)$.}
$$
\rho(g_j)T \rho(g_i) \; = \; \rho(g_j)\big(\rho(\tilde{g}_j) T \rho(\tilde{g}_i) \big) \rho(g_i) 
$$
is within $\varepsilon^3$ of a rank-$N$ operator. Now, using the fact that the $g_i$ sum to $\mathbf{1}_X$, we write
\begin{align*}
[T,\rho(\tilde{g})] \; &= \; \left( \sum_i \rho(g_i)\cdot T \cdot \rho(\tilde{g})\right) - \left( \sum_j \rho(\tilde{g})\cdot T \cdot \rho(g_j)\right) \\
&= \; \sum_{i,j} \left( g(x_j)- g(x_i)\right) \rho(g_i) \cdot T \cdot \rho(g_j) \\
&= \; \sum_{|i-j|=1} \left( g(x_j)- g(x_i)\right) \rho(g_i) \cdot T \cdot \rho(g_j) \\
& \quad +\sum_{|i-j|>1} \left( g(x_j)- g(x_i)\right) \rho(g_i) \cdot T \cdot \rho(g_j) \; .
\end{align*}
Every summand in the second sum is within $\varepsilon^3$ of a rank-$N$ operator, and the sum has less than $n^2$ summands.  Recall that $n$ itself is bounded by $C/\varepsilon$, hence the second sum is within $C^2\varepsilon$ of a rank-$N$ operator. We split the first sum into a sum with $j=i+1$ and a sum with $j=i-1$. Moreover, observe that $\rho(g_i)$ is a projection, and $\rho(g_i)$ and $\rho(g_j)$ are orthogonal for $i\neq j$. Thus the sum with $j=i\pm1$ is a direct sum of operators from $\rho(g_{i\pm 1})H$ to $\rho(g_i)H$. The norm of this direct sum is the supremum of the norms of the summands. It follows that the norm of the sum over $|i-j|=1$ is bounded by
$$
 \big(\sup_i |g(x_{i+1})-g(x_i)| + \sup_i |g(x_{i-1})-g(x_i)| \big)\cdot ||T|| \; \leq \; 2\varepsilon\cdot ||T|| \; . 
$$
Combining these estimates, we conclude that $[T,\rho(\tilde{g})]$ is within $(2||T||+C^2)\varepsilon$ of a rank-$N$ operator, which in turn implies that $[T,\rho(g)]$ is within $(4||T||+C^2)\varepsilon$ of a rank-$N$ operator. Since $C$ is an independent constant, and $N$ only depends on $L,R$ and $\varepsilon$, this concludes the proof.
\end{proof}
With this uniform version of the Kasparov lemma in hand we can proceed to prove the excision isomorphism. The crucial technical calculation is an isomorphism between quotients of the (relative) uniform dual algebras. This isomorphism requires the assumption of $X$ being uniformly resolvable near $Z$.
\begin{prop}
Let $(X,Z)$ be a pair such that $X$ uniformly resolvable near $Z$. Then, for any representation $(H,\rho)$ of $X$ the inclusion $\fD^{u}_\rho(X)\hookrightarrow \fD^{u}_\rho(X\setminus Z)$ induces an isomorphism
$$
\fD^{u}_\rho(X)/\fD^{u}_\rho(X,Z) \; \longrightarrow \fD^{u}_\rho(X\setminus Z)/\fD^{u}_\rho(X\setminus Z,\emptyset)\; .
$$
\end{prop}
\begin{proof}
It holds by definition that $\fD^{u}_\rho(X,Z)=\fD^{u}_\rho(X)\cap \fD^{u}_\rho(X\setminus Z,\emptyset)$. This shows that the inclusion $\fD^{u}_\rho(X)\hookrightarrow \fD^{u}_\rho(X\setminus Z)$ induces a well-defined map on the quotients, and that this map is injective. To show that it is also surjective, we need to prove that $\fD^{u}_\rho(X\setminus Z)=\fD^{u}_\rho(X) + \fD^{u}_\rho(X\setminus Z,\emptyset)$. By our assumption of uniform resolvability there exists a partition of unity $(\psi_i)_{i\in I}$ satisfying the properties (1), (2), and (3) of Definition \ref{ConditionGeomSubspace}. Set $T':=\sum_{i\in I} \rho(\sqrt{\psi_i})T\rho(\sqrt{\psi_i})$. We will show that $T'\in \fD^{u}_\rho(X)$ and $T-T'\in\fD^{u}_\rho(X\setminus Z;\emptyset)$. 

To prove that $T'\in \fD^{u}_\rho(X)$, by Proposition \ref{KasparovLemma} it suffices to show that if $f\in C_b(X)$ is a bounded continuous function, then $\rho(f)T'\rho(g)\sua0$ for $g\in C_0(X\setminus\supp(f))$. Thus, fix $f\in C_b(X)$, let $L,R\geq 0$, $\varepsilon>0$, and let $g\in \LLip_R(X)$ with support disjoint from $\supp(f)$. Fix $\delta>0$ with $L\delta<\varepsilon$, and a function $\psi: \, X \to [0,1]$ with satisfies $\psi= 0$ on the closed $\delta$-neighborhood of $\supp(f)$, $\psi = 1$ on the complement in $X$ of the closed $2\delta$-neighborhood of $\supp(f)$, and such that $\psi$ is $L(\delta)$-Lipschitz\footnote{For example, one could take $\psi(x)$ to be the minimum of $1$ and $d(x,W)/\delta$, where $W$ denotes the closed $\delta$-neighborhood of $\supp(f)$. In this case $L(\delta)=\delta$.} for some $L(\delta)$. Set $\tilde{g}:= \psi g$, a function in $(L(\delta)L)\text{-}\mathrm{Lip}_R(X)$. It holds that $||g-\tilde{g}||\leq \sup_{d(x,\supp(f))\leq 2\delta} |g(x)| \leq 2\delta\cdot L <2\varepsilon$. Consequently, if $\rho(f)T\rho(\tilde{g})$ is within $C\varepsilon$ of some rank-$N$ operator, then $\rho(f)T\rho(g)$ is within $(C+2||f||\cdot||T||)\varepsilon$ of that rank-$N$ operator. Thus, to establish $T'\in\fD^{u}_\rho(X)$, it suffices to show that there exists $N=N(L,R,\varepsilon)$ such that $\rho(f)T\rho(\tilde{g})$ is within $C\varepsilon$ of a rank-$N$ operator for some $C=C(L,R)$. 

We have the equality
\begin{equation} \label{TKasp}
\rho(f) T' \rho(\tilde{g}) \; = \; \sum_i \rho(\sqrt{\psi_i}) \left( \rho(f|_{X\setminus Z}) T \rho(\sqrt{\psi_i }\tilde{g})\right) \; .
\end{equation}
Observe that the $i$-th summand in \eqref{TKasp} can only be non-zero if $\sqrt{\psi_i}f\neq 0$ and $\sqrt{\psi_i}\tilde{g}\neq 0$. In particular the support of $\psi_i$ has to intersect the supports of both $f$ and $\tilde{g}$. By construction, these have a distance $\geq \delta$, whence $\supp(\psi_i)$ has to have a diameter $\geq \delta$. From property (3) of the partition of unity it follows that there exists $\delta'=\delta'(\delta)$ such that $d(\supp(\psi_i),Z)\geq \delta'$. Therefore, the compactum $K:= \supp(\tilde{g})\cap \overline{\bigcup_{\psi_if\neq 0} \supp(\psi_i)}$ has distance $\geq\delta'$ from $Z$, and diameter $\leq R$. Property (2) provides $M(R,\delta')\in\NN$ such that at most $M(R,\delta')$ many $\psi_i$ are non-zero on $K$. In conjunction with property (3), the above discussion yields that at most $M(R,\delta')$ many summands in \eqref{TKasp} are non-zero, and $\sqrt{\psi_i}$ is $L(\delta')$-Lipschitz for those summands by property (1).

The operator $T$ was assumed to lie in $\fD^{u}_\rho(X\setminus Z)$, which by Proposition \ref{KasparovLemma} is equivalent\footnote{Via extension by zero, $C_b(X\setminus Z)$ becomes a subalgebra of $\mathcal{L}^\infty(X)$, thus it makes sense to insert a function from $C_b(X\setminus Z)$ into the representation $\rho$.} to $\rho(f')T\rho(h)\sua0$ for all fixed $f'\in C_b(X\setminus Z)$ and $h\in C_0(X\setminus Z)$ with support disjoint from $\supp(f')$. Thus, there exists a natural number $N'=N'(L,R,\varepsilon)$ such that $\rho(f|_{X\setminus Z})T\rho(h) $ is within $\varepsilon$ of a rank-$N'$ operator for any $h\in(L(\delta')L)\text{-}\mathrm{Lip}_R(X\setminus Z)$ with support disjoint from $\supp(f)$. Since $\sqrt{\psi_i}\tilde{g}$ is in $(L(\delta')L)\text{-}\mathrm{Lip}_R(X\setminus Z)$ for those $i$ contributing to a non-zero summand in \eqref{TKasp}, it follows from \eqref{TKasp} that $\rho(f)T\rho(\tilde{g})$ is within $M(R,\delta')\varepsilon$ of a rank-$N$ operator, where $N=N(L,R,\varepsilon)$ is given by $M(R,\delta')\cdot N'$. Note that since $\delta'$ depends only on $\delta$, which in turn depends only on $L$ and $\varepsilon$, the number $N$ indeed depends only on $L$, $R$ and $\varepsilon$. Thus, the proof of $T'\in\fD^{u}_\rho(X)$ is complete.

To prove that $T-T' \in \fD^{u}_\rho(X\setminus Z,\emptyset)$ we again fix $L,R\geq 0$, $\varepsilon>0$, and take $g\in \LLip_R(X\setminus Z)$. Fix $\delta>0$ with $L\delta<\varepsilon$, and a function $\varphi: \, X\setminus Z \to [0,1]$ with satisfies $\varphi = 0$ on the closed $\delta$-neighborhood of $Z$, $\varphi = 1$ on the complement in $X$ of the closed $2\delta$-neighborhood of $Z$, and such that $\varphi$ is $L(\delta)$-Lipschitz\footnote{For example, one could take $\varphi(x)$ to be the minimum of $1$ and $d(x,Z_\delta)/\delta$, where $Z_\delta$ denotes the closed $\delta$-neighborhood of $Z$. In this case $L(\delta)=\delta$.} for some $L(\delta)$. Set $\tilde{g}:= \varphi g$, a function in $(L(\delta)L)\text{-}\mathrm{Lip}_R(X\setminus Z)$. It holds that $||g-\tilde{g}||\leq \sup_{d(x,Z)\leq 2\delta} |g(x)| \leq 2\delta\cdot L <2\varepsilon$. Consequently, if $(T-T')\rho(\tilde{g})$ is within $C\varepsilon$ of some rank-$N$ operator, then $(T-T')\rho(g)$ is within $(C+2||T-T'||)\varepsilon$ of that rank-$N$ operator. Thus, to establish $T-T'\in\fD^{u}_\rho(X\setminus Z,\emptyset)$, it suffices to show that there exists $N=N(L,R,\varepsilon)$ such that $(T-T')\rho(\tilde{g})$ is within $C\varepsilon$ of a rank-$N$ operator for some $C=C(L,R)$. 

To that end, note that because $\supp(\tilde{g})$ has diameter $\leq R$, and distance $\geq\delta$ from $Z$, there exists $M(R,\delta)\in \NN$ such that $\sqrt{\psi_i}\tilde{g}$ is non-zero for at most $M(R,\delta)$ many $i$. Moreover, for those $\sqrt{\psi_i}g$ is in $(L(\delta)L)\text{-}\mathrm{Lip}_R(X\setminus Z)$. These facts follow from properties (1) and (3). Calculate that
\begin{align*}
T'\rho(\tilde{g}) \; &= \; \sum_i \rho(\sqrt{\psi_i}) T \rho (\sqrt{\psi_i}\tilde{g}) \\
&= \; \sum_{i} \left( \rho(\psi_i \tilde{g})T + \rho(\sqrt{\psi_i})\cdot [T,\rho(\sqrt{\psi_i}\tilde{g})] \right) \\
&= T\rho(\tilde{g}) - [T,\rho(\tilde{g})]+ \sum_i \rho(\sqrt{\psi_i})\cdot [T,\rho(\sqrt{\psi_i}\tilde{g})] \; ,
\end{align*}
which gives
\begin{equation}
 \label{TminusTPrime}
(T-T')\rho(\tilde{g}) \; = \; [T,\rho(\tilde{g})] - \sum_i \rho(\sqrt{\psi_i})\cdot [T,\rho(\sqrt{\psi_i}\tilde{g})] \; .
\end{equation}
Note that the sum on the right-hand side has at most $M(R,\delta)$ non-zero summands. Now, let $N'=N'(L,R,\varepsilon)$ be such that $[T,\rho(h)]$ is within $\varepsilon$ of a rank-$N'$ operator for all $h\in (L')\text{-}\mathrm{Lip}_R(X\setminus Z)$, where $L'=\max\{ L,L(\delta)L\}$, and set $N:= M(R,\delta)\cdot N'$. Since $\delta$ depends only on $L$ and $\varepsilon$, the constant $N$ depends only on $L$, $R$ and $\varepsilon$. Since $\tilde{g}$ and $\sqrt{\psi_i}g$ are in $(L')\text{-}\mathrm{Lip}_R(X\setminus Z)$, equation \eqref{TminusTPrime} implies that $(T-T')\rho(\tilde{g})$ is within $(1+M(R,\delta))\varepsilon$ of a rank-$N$ operator. This concludes the proof.
\end{proof}
We combine this isomorphism at the level of uniform dual algebras with Paschke duality to prove the excision isomorphism. The excision map on relative uniform K-homology is easily defined: If $(H,\rho,F)$ is a relative uniform Fredholm module over $(X,Z)$, then $(H,\rho|_{C_0(X\setminus Z)},F)$ is a uniform Fredholm module over $X\setminus Z$. This follows immediately from the definition of relative uniform Fredhol modules. This restriction process preserves unitary equivalences, operator-homotopies and direct sums, and therefore induces a well-defined homomorphism
$$
\exc: \; K_*^{u}(X,Z) \longrightarrow K_*^{u}(X\setminus Z) \; .
$$
The proof that it is an isomorphism is a straight-forward reduction to the corresponding isomorphism statement for the dual algebras.
\begin{thm}
Let $(X,Z)$ be pair such that $X$ is uniformly resolvable near $Z$. Then, the excision map $\exc: \, K_*^{u}(X,Z) \to K^{u}_*(X\setminus Z)$ is a natural isomorphism.
\end{thm}
\begin{proof}
Let $(H,\rho)$ be a representation of $X$. The previous proposition provides an isomorphism 
$$
\fD^{u}_\rho(X)/\fD^{u}_\rho(X,Z) \xlongrightarrow{\sim} \fD^{u}_\rho(X\setminus Z)/\fD^{u}_\rho(X\setminus Z,\emptyset) \; .
$$
If $V: \, (H^\rho,\rho)\to (H^\sigma,\sigma)$ is an isometry uniformly covering the identitiy on $X$, then the diagram\\
\centerline{
\xymatrix{ \fD^{u}_\rho(X)/\fD^{u}_\rho(X,Z)\ar[r] \ar[d]_{\Ad(V)} & \fD^{u}_\rho(X\setminus Z)/\fD^{u}_\rho(X\setminus Z,\emptyset)\ar[d]^{\Ad(V)} \\
\fD^{u}_\sigma(X)/\fD^{u}_\sigma(X,Z) \ar[r] & \fD^{u}_\sigma(X\setminus Z)/\fD^{u}_\sigma(X\setminus Z,\emptyset)
}}
commutes. Thus, the isomorphism passes to the colimits, and upon taking the K-theory, we obtain an isomorphism\footnote{There is an ever so slight subtlety involved in the statement that the map on colimits is an isomorphism, because the left-hand side is defined via the colimit over representations of $X$, and that on the right-hand side via representations of $X\setminus Z$. However, any representation $(H,\rho)$ of $X\setminus Z$ extends to a representation $(H,\tilde{\rho})$ of $X$ via $\tilde{\rho}(f)u:=\lim_{n\to\infty} \rho(f\psi_n)u$, where $\psi_n$ is an approximate unit of $C_0(X\setminus Z)$, see \cite[before Proposition 2.6.3]{HigsonRoe2000}.}
$$
K_*\left( \fD^{u}(X)/\fD^{u}(X,Z) \right)  \longrightarrow K_*\left( \fD^{u}(X\setminus Z)/\fD^{u}(X\setminus Z,\emptyset) \right)\; .
$$
This isomorphism may be composed with the inverse of the isomorphism $K_*(\fD^{u}(X\setminus Z)) \to K_*(\fD^{u}(X\setminus Z)/\fD^{u}(X\setminus Z,\emptyset))$ from Proposition \ref{RelEmptyset} to obtain an isomorphism
$$
K_*\left( \fD^{u}(X)/\fD^{u}(X,Z) \right) \longrightarrow K_*\left( \fD^{u}(X\setminus Z) \right) \; .
$$
Since on the level of representatives this isomorphism is little more than an inclusion map, it is not hard to see that the diagram \\
\centerline{
\xymatrix{
K_{*+1}(\fD^{u}(X)/\fD^{u}(X,Z)) \ar[r]^-{\sim} \ar[d]_{\varphi_*} & K_{*+1}(\fD^{u}(X\setminus Z))\ar[d]^{\varphi_*} \\
K_{*}^{u}(X,Z)\ar[r]_{\exc} & K_{*}^{u}(X\setminus Z)
}}
commutes. We conclude that the excision map on uniform K-homology is an isomorphism. Moreover, that it is natural is rather immediate from its definition. Thus the theorem is proved.
\end{proof}
\subsection{The long exact sequence}
We come to the long exact sequence. Again, the essential technical component is an isomorphism statement for the uniform dual algebras, this time at the level of K-theory. This isomorphism requires the assumption of transversal control. The statement is essentially proved by \v{S}pakula in Section 5 of \cite{Spakula2009}, but since some of the detail differ we give the proof here.
\begin{prop}\label{InclusionIsoDual}
Let $(X,Z)$ be a pair such that $X$ is transversally controlled near $Z$. Then, there is a natural isomorphism 
$$
K_*\left( \fD^{u}(Z)\right) \; \cong \; K_*\left( \fD^{u}(X,Z) \right) \; .
$$
\end{prop}
\begin{proof}
Let $(H_X,\rho_X)$ and $(H_Z,\rho_Z)$ be representations of $X$ and $Z$, and $V: \, (H_Z,\rho_Z) \to (H_X,\rho_X)$ an isometry uniformly covering the restriction map $\pi: \, C_0(X)\to C_0(Z)$. Let $\sigma: \, C_0^+(Z)\to C_0^*(X)$ be an approximately filtered completely positive section as in Definition \ref{ConditionGeomSubspace}. Out of these data, in the course of the proof of \cite[Lemma 5.3]{Spakula2009} \v{S}pakula constructs the following half-infinite commutative diagram \\
\begin{align}\label{LESdiagcolim}
\xymatrix{ &  K_*(\fD^{u}_{\rho_X}(X,Z)) \ar[r] ^{\Ad(SW)} \ar[d]_{\Ad(W)} & K_*(\fD^{u}_{\rho'_X}(X,Z)) \ar[r] \ar[d]& \cdots \\
K_*(\fD^{u}_{\rho_Z}(Z)) \ar[r]_{\Ad(WV)} \ar[ru]^{\Ad(V)} & K_*(\fD^{u}_{\rho'_Z}(Z))\ar[r] \ar[ru]^{\Ad(S)} & \cdots 
}
\end{align}
The representation $(H'_Z,\rho'_Z)$ of $Z$ is given by a dilation of the completely positive map $\rho_X\circ\sigma$ to a representation, meaning there is an isometry $W: \, H_X\to H_Z'$ such that $W^*\rho_Z'W = \rho_X\circ \sigma$. This dilation is provided by Stinespring's theorem. The assumption that $\sigma$ is approximately filtered ensures that $\Ad(W)$ maps $\fD^{u}_{\rho_X}(X,Z)$ to $\fD^{u}_{\rho'_Z}(Z)$, see Remark \ref{ApproxFilter}. The representation $\rho_X'$ is the direct sum of $\rho_X$ and $\rho_Z'\circ\pi$, and $S$ is the isometry including $H_X$ into $H_X\oplus H'_Z$. This isometry $S$ again covers $\pi$ uniformly, and transversal control is again needed to ensure that $\Ad(S)$ maps $\fD^{u}_{\rho'_Z}(Z)$ to $\fD^{u}_{\rho'_X}(X,Z)$. This first commutative parallelogram is extended to the right by repeating the same construction for $S$ instead of $V$ to obtain a second such parallelogram, and so on. All vertical arrows come from isometries that uniformly cover the respective identity, or are at least homotopic to conjugation by isometries uniformly covering the identity (see \cite[Section 5]{Spakula2009}). 

We extract the desired natural isomorphism from the diagram \eqref{LESdiagcolim}. Let $(H,\rho_Z)$ be a representation. By setting $\rho_X:= \rho_Z\circ\pi$ and $V:=\id_{H}$, we land in the situation of the diagram \eqref{LESdiagcolim}. We may compose the map $\Ad(V)$ with the map $K_*(\fD^{u}_{\rho_X}(X,Z))\to K_*(\fD^{u}(X,Z))$ to obtain a map $K_*(\fD^{u}_{\rho_Z}(Z))\to K_*(\fD^{u}(X,Z))$. Taking the colimit over $\rho_Z$, obtain the map $\Psi: \, K_*(\fD^{u}(Z))\to K_*(\fD^{u}(X,Z))$, which is uniquely characterized by the property of making the diagram 
\begin{align} \label{diagDefPsi}
\xymatrix{ K_*(\fD^{u}_{\rho_Z}(Z))\ar[r]^{\iota_{\rho_Z}} \ar[d]_{\Ad(V)} & K_*(\fD^{u}(Z)) \ar[d]^\Psi \\
K_*(\fD^{u}_{\rho_X}(X,Z))\ar[r]_{\iota_{\rho_X}} & K_*(\fD^{u}(X,Z))
}
\end{align}
commute for all $\rho_Z$, and $\rho_X$ and $V$ as above. Here the maps $\iota_{\rho_Z}$ and $\iota_{\rho_X}$ are the canonical maps into the respective colimit. We will show that $\Psi$ is a natural isomorphism. To that end, we note first that the diagram \eqref{diagDefPsi} still commutes for any representation $(H_X,\rho_X)$ and any isometry $V:\, (H_Z,\rho_Z)\to (H_X,\rho_X)$ uniformly covering $\pi$. This is because we have a commutative diagram \\
\centerline{
\xymatrixcolsep{4pc}\xymatrix{ K_*(\fD^{u}_{\rho_Z}(Z))\ar[r]^-{\Ad(\id_{H_Z})} \ar[d]_{\Ad(V)} & K_*(\fD^{u}_{\rho_Z\circ\pi}(X,Z)) \ar[d] \ar[ld]_{\Ad(V)} \\
K_*(\fD^{u}_{\rho_X}(X,Z))\ar[r] & K_*(\fD^{u}(X,Z))
}}
where for the vertical arrow $\Ad(V)$ is viewed as uniformly covering $\pi$, and for the diagonal arrow as uniformly covering $\id_X$. Then, the left commutative triangle comes from the factorization $V=V\circ \id_{H_Z}$ of uniform covering isometries, and the right commutative triangle comes from the definition of the colimit $K_*(\fD^{u}(X,Z))$. The combination of this latter commutative diagram with the diagram \eqref{diagDefPsi} for $\rho_X=\rho_Z\circ \pi$ and $V=\id_{H_Z}$ then gives the diagram \eqref{diagDefPsi} for arbitrary $\rho_X$ and $V$. Then, diagrams \eqref{LESdiagcolim} and \eqref{diagDefPsi} combine to yield the commutative diagram
\begin{align} \label{LESlargeDiag}
\xymatrixcolsep{4pc}\xymatrix{ K_*(\fD^{u}_{\rho_Z}(Z))\ar@/^2pc/[rr]^{\iota_{\rho_Z}}\ar[r]^{\Ad(WV)} \ar[d]_{\Ad(V)} & K_*(\fD^{u}_{\rho_Z'}(Z))\ar[r]^{\iota_{\rho_Z'}} \ar[d]^{\Ad(S)} & K_*(\fD^{u}(Z)) \ar[d]^\Psi \\
K_*(\fD^{u}_{\rho_X}(X,Z))\ar@/_2pc/[rr]_{\iota_{\rho_X}} \ar[ru]^{\Ad(W)}\ar[r]_{\Ad(SW)} & K_*(\fD^{u}_{\rho_X'}(X,Z))\ar[r]_{\iota_{\rho_X'}} & K_*(\fD^{u}(X,Z))
}
\end{align}
This diagram, considered for suitable $\rho_Z$, $\rho_X$ and $V$, will yield bijectivity of $\Psi$. Let us start with injectivity. Let $x\in K_*(\fD^{u}(Z))$ be in the kernel of $\Psi$. The element $x$ may written as $x=\iota_{\rho_Z}(u)$ for some $\rho_Z: \, C_0(Z)\to \fB(H_Z)$ and $u\in K_*(\fD^{u}_{\rho_Z}(Z))$. Let $\rho_X:= \rho_Z\circ\pi: \, C_0(X)\to\fB(H_Z)$ and $V=\id_{H_Z}$, and construct the corresponding diagram \eqref{LESlargeDiag}. Then, $\Psi(x)=0$ implies $\iota_{\rho_X}(\Ad(V)u)=0$. It may be assumed that $\Ad(V)u=0$.\footnote{An element $v\in K_*(\fD^{u}_{\rho_X}(X,Z))$ satisfies $\iota_{\rho_X}v=0$ if and only if there is a representation $\sigma_X: \, C_0(X) \to \fB(H')$ and an isometry $U: \, H_X\to H'$ uniformly covering $\id_X$, such that $\Ad(U)v=0$. This fact about colimits of abelian groups may be derived from the explicit construction of such colimits.We apply this to $v=\Ad(V)u$. The isometry $UV: \, H_Z\to H'$ uniformly covers $\pi$. Thus, me may construct the diagram \eqref{LESlargeDiag} with $UV$ instead of $V$, and $\sigma_X$ instead of $\rho_X$. Then, $\iota_{\sigma_X}(\Ad(UV)u)=0$ would hold, and $\Ad(UV)u=0$ by construction. Then, we are in the situation of the text above.} From \eqref{LESlargeDiag} it then follows that $x=\iota_{\rho_Z}(u)=(\iota_{\rho_Z'}\circ\Ad(W))(\Ad(V)u)=0$, proving injectivity. 

Let us turn to the surjectivity of $\Psi$. Any element $y\in K_*(\fD^{u}(X,Z))$ is of the form $y=\iota_{\rho_X}(v)$ for some $(H_X,\rho_X)$, and some $v\in K_*(\fD^{u}_{\rho_X}(X,Z))$. The representation can be extended to a representation of the algebra $\mathcal{L}^\infty(X)$ of bounded measurable functions on $X$, which we still denote by $\rho_X$, and this extension may then be restricted to the subalgebra $C_0(Z) \subseteq \mathcal{L}^\infty(X)$. We consider the resulting representation $\rho_Z$ of $C_0(Z)$ as acting on the Hilbert space $H_Z:=PH_X$, where $P:=\rho_X(\mathbf{1}_Z)$. Let $V$ be the inclusion $H_Z\hookrightarrow H_X$. Then,
$$
V^* \rho_X(f) V \; = \; \rho_X(\mathbf{1}_Z) \cdot \rho_X(f) \; = \; \rho_X(\mathbf{1}_Z f) \; = \; \rho_Z(\pi(f)) \; ,
$$
whence $V$ uniformly covers $\pi$. Construct the diagram \eqref{LESdiagcolim} corresponding to these choices for $\rho_Z$ and $V$. According to this diagram, we may then write 
\begin{equation} \label{PsiInv}
y=\iota_{\rho_X}(v)=\Psi((\iota_{\rho_Z'}\circ\Ad(W))(v))
\end{equation}
This proves surjectivity and hence bijectivity of $\Psi$. 

Naturality of $\Psi$, lastly, follows almost immediately from the definitions. Let $(Y,W)$ be another pair, and $\phi: \, (X,Z)\to (Y,W)$ a filtered map of pairs. By the defining properties of $\Psi$ and the maps induced by $\phi$, the diagram \\
\centerline{
\xymatrix{ K_*(\fD^{u}_{\rho_Z}(Z))\ar[ddd]_{\Ad(V)} \ar[rd]^{\iota_{\rho_Z}} \ar[rrr]^{\Ad(\id_{H_Z})}& & & K_*(\fD^{u}_{\rho_Z\circ\phi}(W)) \ar[ld]_{\iota_{\rho_Z\circ\phi}} \ar[ddd]^{\Ad(V)} \\
& K_*(\fD^{u}(Z)) \ar[d]_{\Psi} \ar[r]^{\phi_*} & K_*(\fD^{u}(W)) \ar[d]^{\Psi} & \\
& K_*(\fD^{u}(X,Z)) \ar[r]^{\phi_*} & K_*(\fD^{u}(Y,W))& \\ 
K_*(\fD^{u}_{\rho_X}(X,Z)) \ar[ru]^{\iota_{\rho_X}} \ar[rrr]_{\Ad(\id_{H_X})}& & & K_*(\fD^{u}_{\rho_X\circ\phi}(Y,W)) \ar[lu]_{\iota_{\rho_X\circ\phi}}
}}
commutes for $\rho_Z$, $\rho_X$, and $V$ uniformly covering $\pi$. Since every element of $K_*(\fD^{u}(Z))$ is in the image of some $\iota_{\rho_Z}$, it follows that the inner square commutes as well, which gives naturality of $\Psi$.
\end{proof}
Existence of the 6-term exact sequence in uniform K-homology now follows directly from the combination the 6-term exact sequence in operator K-theory and Proposition \ref{InclusionIsoDual}.
\begin{thm}
Let $(X,Z)$ be a pair such that $X$ is transversally controlled near $Z$. There exists a natural 6-term exact sequence \\
\centerline{
\xymatrix{
K_0^{u}(Z) \ar[r] & K^{u}_0(X) \ar[r] & K^{u}_0(X,Z)\ar[d]^{\partial} \\
K_1^{u}(X,Z)\ar[u]^{\partial} & K^{u}_1(X) \ar[l] &\ar[l] K^{u}_1(Z)
}}
\end{thm}
\begin{proof}
For every representation $(H,\rho)$ of $X$ there is an associated short exact sequence
$$
0 \longrightarrow \fD^{u}_\rho(X,Z) \longrightarrow \fD^{u}_\rho(X) \longrightarrow \fD^{u}_\rho(X)/\fD^{u}_\rho(X,Z) \longrightarrow 0
$$
of the uniform dual algebras. Given that the maps in this short exact sequence are simply inclusion and quotient map, the maps $\Ad(V)$ for some uniform covering isometry $V$ are clearly compatible with these exact sequences, whence we may take the colimit to obtain a short exact sequence 
$$
0 \longrightarrow \fD^{u}(X,Z) \longrightarrow \fD^{u}(X) \longrightarrow \fD^{u}(X)/\fD^{u}(X,Z) \longrightarrow 0
$$
of $C^*$-algebras. We obtain the 6-term exact sequence \\
\centerline{
\xymatrix{
K_1(\fD^{u}(X,Z)) \ar[r] & K_1(\fD^{u}(X)) \ar[r] & K_1(\fD^{u}(X)/\fD^{u}(X,Z))\ar[d]^{\partial} \\
K_0(\fD^{u}(X)/\fD^{u}(X,Z))\ar[u]^{\partial} & K_0(\fD^{u}(X)) \ar[l] &\ar[l] K_0(\fD^{u}(X,Z))
}}
in K-theory. Now, Paschke duality and the isomorphism $K_*(\fD^{u}(X,Z))\cong K_*(\fD^{u}(Z))$ provided by Proposition \ref{InclusionIsoDual} yield the 6-term exact sequence in uniform K-homology. As every part of this construction is natural in $(X,Z)$, the resulting exact sequence is also natural. 
\end{proof}
We obtain the long exact corollary by extending the six-term exact sequence via twofold periodicity.
\begin{cor}
Let $(X,Z)$ be a pair such that $X$ is transversally controlled near $Z$. There exists a natural long exact sequence 
$$
\cdots \longrightarrow K^{u}_{-p}(Z) \longrightarrow K^{u}_{-p}(X) \longrightarrow K^{u}_{-p}(X,Z) \xlongrightarrow{\partial} K^{u}_{-(p-1)}(Z)\longrightarrow \cdots \; .
$$
\end{cor}
Combining excision and the long exact sequence we may also draw the following corollary.
\begin{cor}
Let $(X,Z)$ be a pair such that $X$ is uniformly resolvable and transversally controlled near $Z$. Then, there is a natural long exact sequence 
$$
\cdots \longrightarrow K^{u}_{-p}(Z) \longrightarrow K^{u}_{-p}(X) \longrightarrow K^{u}_{-p}(X\setminus Z) \xlongrightarrow{\partial} K^{u}_{-(p-1)}(Z)\longrightarrow \cdots \; .
$$
\end{cor}
We note that contrary to usual conventions for homology theories the boundary map increases the degree. This is a consequence of the degree convention that $p$-multigraded Fredholm modules define classes in the K-homology of degree $-p$. We have adopted this convention from \cite{HigsonRoe2000}. Of course twofold periodicity in both K-homology and K-theory reduces this to little more than a formal matter. Thus, we aim to stay consistent with our degree conventions where it seems helpful, but overall prioritize readability over consistency. Also due to twofold periodicity in both (uniform) K-homology and K-theory we will not differentiate too carefully between long exact sequences and 6-term exact sequences, as one is turned into the other via periodicity.
\subsection{The relative Mayer-Vietoris sequence}
We also provide a Mayer-Vietoris sequence for relative uniform K-homology. The non-relative Mayer-Vietoris sequence was constructed by \v{S}pakula \cite{Spakula2009}, and our construction is an elaboration of the non-relative case.
\begin{definition}
Let $(X,Z)$ be a pair, and let $Y_1,Y_2\subseteq X$ be closed subspaces such that $Y_1\cup Y_2 =X$. Set $W_j:=Y_j\cap Z$ We say that $(X,Z)=(Y_1\cup Y_2,W_1\cup W_2)$ is a \emph{bounded Mayer-Vietoris decomposition} if:
\begin{itemize}
\item $d(Y_1\setminus Y_2,Y_2\setminus Y_1)>0$,
\item and the pairs $(X,Y_j)$, $(Y_j,Y_1\cap Y_2)$, $(Y_j,W_j)$, and $(Y_1\cap Y_2,W_1\cap W_2)$ are transversally controlled.
\end{itemize} 
\end{definition}
The Mayer-Vietoris sequence in relative uniform K-homology is constructed via an analogous Mayer-Vietoris sequence in the K-theory of the relative uniform dual algebras. We refer to Appendix \ref{AppendixKTheory} for a review of the Mayer-Vietoris sequence in operator K-theory, and Proposition \ref{MVQuot} will be especially relevant to us.
\begin{thm}
Let $(X,Z)=(Y_1\cup Y_2,W_1\cup W_2)$ be a bounded Mayer-Vietoris decomposition. There exists a natural Mayer-Vietoris sequence \\
\centerline{
\xymatrix{
K_0^{u}(Y_1\cap Y_2, W_1\cap W_2) \ar[r] & K^{u}_0(Y_1,W_1) \oplus K_0^{u}(Y_2,W_2) \ar[r] & K^{u}_0(X,Z)\ar[d]^{\delta} \\
K_1^{u}(X,Z)\ar[u]^{\delta} & K^{u}_1(Y_1,W_1) \oplus K_1^{u}(Y_2,W_2) \ar[l] &\ar[l] K_1^{u}(Y_1\cap Y_2, W_1\cap W_2)
}}
\end{thm}
\begin{proof}
Fix a representation $(H,\rho)$ of $X$. We consider the ideals $\fD^{u}_\rho(X,Y_j) \subseteq \fD^{u}_\rho(X)$. By our assumption that $d(Y_1\setminus Y_2,Y_2\setminus Y_1)>0$ we can find Lipschitz functions $\psi_j: \, X\to [0,1]$, $\supp(\psi_j)\subseteq W_j$, with $\psi_1+\psi_2=1$. Then, any $T\in \fD^{u}_\rho(X)$ can be written as $T=T \psi_1 + T \psi_2$. Since $T\psi_j \in \fD^{u}_\rho(X,Y_j)$ we conclude that $\fD^{u}_\rho(X)=\fD^{u}_\rho(X,Y_1)+\fD^{u}_\rho(X,Y_2)$ (compare \cite[Section 5]{Spakula2009}). Moreover, we have
\begin{align*}
\fD^{u}_\rho(X,Y_1\cap Y_2)&=  \fD^{u}_\rho(X,Y_1)\cap\fD^{u}_\rho(X,Y_2)\\
\fD^{u}_\rho(X,W_j) &= \fD^{u}_\rho(X,Y_j)\cap\fD^{u}_\rho(X,Z) \\
\fD^{u}_\rho(X,W_1\cap W_2) &= \fD^{u}_\rho(X,Y_1)\cap\fD^{u}_\rho(X,Y_2)\cap\fD^{u}_\rho(X,Z) \; .
\end{align*}
Thus Proposition \ref{MVQuot} provides a Mayer-Vietoris sequence for the decomposition
\begin{align}
\fQ^{u}_\rho(X,Z) &= \fD^{u}_\rho(X,Y_1)/\fD^{u}_\rho(X,W_1) + \fD^{u}_\rho(X,Y_2)/\fD^{u}_\rho(X,W_2) \label{MVDecomp} \\
\fD^{u}_\rho(X,Y_1\cap Y_2)/\fD^{u}_\rho(X,W_1\cap W_2) &=  \left( \fD^{u}_\rho(X,Y_1)/\fD^{u}_\rho(X,W_1)\right)\cap \left( \fD^{u}_\rho(X,Y_2)/\fD^{u}_\rho(X,W_2) \right) \; . \nonumber
\end{align}
Due to the naturality of the Mayer-Vietoris sequence we also get an exact sequence corresponding to the colimit of the decompositions \eqref{MVDecomp}.

Now, by Proposition \ref{InclusionIsoDual} -applicable due to transerval control for the various pairs involved- the inclusions $\fD^{u}(Y_j)\to \fD^{u}(X,Y_j)$, $\fD^{u}(W_j)\to \fD^{u}(Y,W_j)$ and $\fD^{u}(W_j)\to \fD^{u}(X,W_j)$ induce isomorphisms on K-theory. Combining the latter two implies that $\fD^{u}(Y_j,W_j)\to \fD^{u}(X,W_j)$ also induces an isomorphism on K-theory. The same applies to the intersections $Y_1\cap Y_2$ and $W_1\cap W_2$ as well. Naturality of the K-theory exact sequence and the 5-lemma provide us with isomorphisms
\begin{align*}
K_*(\fQ^{u}(Y_j,W_j)) & \cong K_*\left( \fD^{u}_\rho(X,Y_j)/\fD^{u}_\rho(X,W_j) \right) \\
K_*(\fQ^{u}(Y_1\cap Y_2,W_1\cap W_2)) & \cong K_*\left( \fD^{u}_\rho(X,Y_1\cap Y_2)/\fD^{u}_\rho(X,W_1\cap W_2) \right)
\end{align*}
Inserting these isomorphisms into the above Mayer-Vietoris sequence and appealing to Paschke duality yields the relative Mayer-Vietoris sequence for uniform K-homology.
\end{proof}
Using twofold periodicity we again extend the 6-term Mayer-Vietoris sequence to a Mayer-Vietoris sequence involving all degrees.
\begin{cor}
Let $(X,Z)=(Y_1\cup Y_2,W_1\cup W_2)$ be a bounded Mayer-Vietoris decomposition. There exists a natural Mayer-Vietoris sequence
$$
\cdots \to K_{-p}^{u}(Y_1,W_1)\oplus K_{-p}^{u}(Y_2,W_2) \to K^{u}_{-p}(X,Z) \xrightarrow{\delta}K^{u}_{-(p-1)}(Y_1\cap Y_2,W_1\cap W_2)\to \cdots \; .
$$
\end{cor}
\section{Suspension and the Kasparov product}
Let $X$ be a metric space. The \emph{suspension} of $X$ is defined as $SX:= (0,1)\times X$, endowed with the product metric. It arises as the complement of $X\cong \{0\}\times X$ in the cylinder $CX:= [0,1)\times X$. The algebra $C_0(CX)$ is contractible, hence $CX$ has vanishing K-homology. Consequently, the boundary map $\partial: \, K_{*-1}(SX)\to K_*(X)$ is an isomorphism, the \emph{suspension isomorphism}. We will convince ourselves that the same argument goes through in uniform K-homology.
\begin{prop}
Let $p\in\ZZ_{\geq -1}$. Let $X$ be a metric space of jointly bounded geometry, and $CX=[0,1)\times X$ the cylinder equipped with the product metric. The boundary map
$$
\partial: \, K^{u}_{-(p+1)}(SX)\cong K^{u}_{-(p+1)}\left(CX,X\right) \longrightarrow K^{u}_{-p}(X)
$$
in the long exact sequence of the pair $(CX,X)$ is a natural isomorphism.
\end{prop}
\begin{proof}
We show that $K^{u}_*(CX)=0$. The map 
$$
\phi: \, C_0(CX)\longrightarrow C_0([0,1]\times CX)=C_0([0,1]\times [0,1)\times X) \; , \; (\phi f)(t,s,x)=f(s(1-t)+t,x)
$$ 
provides a homotopy from $\mathrm{eval}_0\circ\phi=\id_{CX}$ to $\mathrm{eval}_1\circ\phi=0$, the usual contraction of $C_0(CX)$. Then $K_*^{u}(CX)=0$ follows from homotopy invariance of uniform K-homology once we prove that $\phi$ is filtered. To that end note that $\phi$ is the restriction of the map
$$
\tilde{\phi}=(m\times\id_X)*: \, C_0\left([0,1]\times X \right) \to C_0\left( [0,1]\times [0,1]\times X\right)
$$
to the kernel of $\id_{[0,1]}\times \mathrm{eval}_1\times \id_X$, where $m: \, [0,1] \times [0,1] \to [0,1]$ is given by $m(s,t)=s(1-t)+t$. The map $m$ is Lipschitz, and trivially uniformly cobounded, so that $m\times\id_X$ also has these properties. Therefore $\tilde{\phi}$ is filtered. Since $\tilde{\phi}$ maps $C_0(CX)$ into $C_0([0,1]\times CX)$ we conclude that $\phi$, the restriction of $\tilde{\phi}$ to that subalgebra, is filtered as well.
\end{proof}
There is also a natural map going in the opposite direction. Indeed the Kasparov product provides a map 
$$
K^{u}_*(X)\cong K^{u}_{-1}((0,1)) \otimes K^{u}_*(X) \longrightarrow K^{u}_{*-1}(SX) \; .
$$
The first isomorphism follows from $K^{u}_{-1}((0,1)) \cong K_{-1}((0,1))\cong \ZZ$. This map is thus given by $\alpha\mapsto d \times\alpha$ for $\alpha\in K^{u}_*(X)$, where $d$ is the generator of $K_{-1}((0,1))$ corresponding to the canonical orientation of the interval $(0,1)$. For analytic K-homology this map is (almost) inverse to the suspension isomorphism. We will prove that the same holds for uniform K-homology as well.
\begin{thm}\label{Suspension}
Let $p\in\ZZ_{\geq -1}$. Let $X$ be a proper metric space of bounded geometry. Then, the map 
$$
K^{u}_{-p}(X) \longrightarrow K^{u}_{-(p+1)}(SX) \; , \; \alpha \longmapsto d\times \alpha
$$
is inverse to the suspension isomorphism up to a sign. Indeed, it holds that
\begin{equation}
\partial \left( d\times \alpha \right) = (-1)^{p+1} \alpha \quad \forall\; \alpha\in K^{u}_{-p}(X) \; . 
\end{equation}
\end{thm}
The sign in this formula is a matter of convention. We use the boundary maps on uniform K-homology coming from Paschke duality and the K-theory boundary maps. Instead one could construct the long exact sequence using the Kasparov product and a mapping cone construction. The two resulting boundary maps $K_{-(p+1)}^{u}(SX)\to K_{-p}(X)$ differ by a sign of $(-1)^{p+1}$, so that the formula in Theorem \ref{Suspension} is sign-free if the latter boundary map is used. See \cite[Section 9.6]{HigsonRoe2000} for a discussion of this issue in analytic K-homology. However, we prefer to stay consistent with regards to the boundary map used, hence the sign. 
\begin{proof}
By the formal twofold periodicity it suffices to consider $p\geq 0$. Suppose that $x=(H^\rho,\rho,F)$ is a $(p+1)$-multigraded uniform Fredholm over $SX$, with $(H^\rho,\rho)$ a non-degenerate representation. Then, $\rho$ extends to a representation of the algebra $C_b([0,1]\times X)$ on $H^\rho$, in particular we can make sense of the the operator $G_0 := \rho(g)$ with $g(t,x):= 2t-1$. Define $G:=\kappa \varepsilon_1 G_0$, where $\kappa$ and $\varepsilon_1$ are the grading and first multigrading operator on $H$ respectively. Define
$$
s(x):= \left( H^\rho, \rho_X, V:= G+(I-G^2)^{-\frac{1}{2}}F \right) \; .
$$
Here $\rho_X$ is the representation of $X$ defined by $\rho_X(f):=\rho(\mathbf{1}\otimes f)$, and we view $H^\rho$ as $p$-multigraded by $\varepsilon_2,\cdots,\varepsilon_{p+1}$.  Then, $s(x)$ is a $p$-multigraded uniform Fredholm module over $X$, and it holds that 
$$
\partial [x] = (-1)^{p+1} [s(x)] \; \in K^{u}_{-p}(X) \; .
$$
This can be shown in the same way as for analytic K-homology.\footnote{Indeed this can be seen almost entirely on the level of Fredholm modules, where $\partial$ and $s$ are defined by the same formulas on uniform and analytic K-homology. Thus the proofs of Lemmas 9.6.9, 9.6.10, and Proposition 9.6.11 in \cite{HigsonRoe2000}, and the previous results on which they rely, which work on the Fredholm module level, go through verbatim in uniform K-homology.} Thus it suffices to prove that $[s(d\times x)]=[x]\in K_{-p}^{u}(X)$. To that end we use \cite[Lemma 9.5.7]{HigsonRoe2000}, and the last bullet point of Theorem \ref{KasparovProperties}, to write
$$
[x] = 1\times [x] = [s(d)] \times [x] \; ,
$$
where $1\in K_0^{u}(\mathrm{pt})\cong\ZZ$ is a generator. Therefore we are done once we show that $s(d\times x)$ is aligned to $s(d)$ and $x$ in the sense of Definition \ref{AlignedModules}. Indeed, by Remark \ref{PropagationKasparov} it suffices to check points (b) and (c) of that definition.

Write $d=(H:=L^2[0,1]\oplus L^2[0,1],\mu\oplus\mu,D)$ as in \cite[Definition 9.5.1]{HigsonRoe2000} with $\mu$ the representation of $C_0(0,1)$ by multiplication operators, and $D=\begin{pmatrix}
0 & -iY \\ iY& 0
\end{pmatrix}$ and $Y$ as in \cite[Definition 9.5.1]{HigsonRoe2000}. Let $N_1, N_2 \in \fB(H\hat{\otimes} H^\rho)$ as in \cite[Lemma 2.23]{Engel2019}. Then $d\times x$ is represented by 
$$
\left( H\hat{\otimes} H^\rho, (\mu\oplus\mu)\hat{\otimes} \rho, D\times F:= N_1(D\hat{\otimes} I) + N_2 (I\hat{\otimes} F) \right) \; .
$$
Consequently, 
$$
s(d\times x) = \left( H\hat{\otimes} H^\rho, I\hat{\otimes} \rho, V:= G + (I-G^2)^{-\frac{1}{2}} (D\times F) \right) \; ,
$$
with $G:= (\kappa \varepsilon_1 (\mu(g)\oplus\mu(g)))\otimes I$. Let $f\in C_0(X)$. As part of the conditions for $s(d\times x)$ to be aligned with $s(d)$ and $x$ we need to show that
\begin{equation}\label{PosModCpt}
(I\hat{\otimes} \rho(f))\left( V(I\hat{\otimes} F) + (I\hat{\otimes} F) V \right)(I\hat{\otimes} \rho(\bar{f})) 
\end{equation}
is positive modulo compacts. Inserting the definition of $V$ presents \eqref{PosModCpt} as the sum of the following terms:
\begin{itemize}
\item $(I\hat{\otimes} \rho(f))\left( G(I\hat{\otimes} F) + (I\hat{\otimes}F)G \right)(I\hat{\otimes} \rho(\bar{f}))$. Recall that $G=(\kappa \varepsilon_1 (\mu(g)\oplus\mu(g)))\otimes I$. The operator $\kappa \varepsilon_1 (\mu(g)\oplus\mu(g))$ is odd, and so is $F$. The operators $G$ and $I\hat{\otimes}F$ therefore anti-commute.\footnote{Recall that we use the graded composition rule for tensor product operators here.} Thus this summand vanishes.
\item $(I\hat{\otimes} \rho(f))\left( \sqrt{I-G^2}N_1(D\hat{\otimes} I)(I\hat{\otimes} F) + (I\hat{\otimes}F)\sqrt{I-G^2}N_1(D\hat{\otimes} I)\right)(I\hat{\otimes} \rho(f))$. Here $$
\sqrt{I-G^2}=\mu(\sqrt{1-g^2})\oplus \mu(\sqrt{1-g^2})
$$ 
is an even operator. It commutes with every other operator involved, at least up to compacts. The same holds for the even operator $N_1$. Modulo compacts this summand therefore equals to a product of operators, with one factor being $(D\hat{\otimes}I)(I\hat{\otimes} F) + (I\hat{\otimes}F)(D\hat{\otimes}I)$. Since both $D$ and $F$ are odd, this factor vanishes, and with it the entire summand.
\item $(I\hat{\otimes} \rho(f))\left( \sqrt{I-G^2}N_2(I\hat{\otimes} F)(I\hat{\otimes} F) + (I\hat{\otimes}F)\sqrt{I-G^2}N_1(I\hat{\otimes} F)\right)(I\hat{\otimes} \rho(\bar{f}))$. Again, we can factor out $\sqrt{I-G^2}$ and $N_2$ modulo compacts. The remaining term inside the brackets is given by $2(I\hat{\otimes}F)^2=2I\hat{\otimes} F^2$. Upon multiplication with $I\hat{\otimes}\rho(f)$ this becomes $2I\hat{\otimes}I=2I$ up to compacts. Modulo compacts this last summand therefore equals
$$
2(I\hat{\otimes} \rho(f))\sqrt{I-G^2}N_2(I\hat{\otimes} \rho(\bar{f})) \sim 2 \sqrt{I-G^2}N_2(I\hat{\otimes} \rho(|f|^2)) \; .
$$ 
The operators in this product commute modulo compacts, hence they are positive modulo compacts.
\end{itemize}
We conclude that \eqref{PosModCpt} is indeed positive modulo compacts. Analogously 
$$
(I\hat{\otimes} \rho(f))\left( V(D\hat{\otimes} I) + (D\hat{\otimes} I) V \right)(I\hat{\otimes} \rho(\bar{f})) 
$$
is positive modulo compact operators as well.

Lastly, $G$ and $\sqrt{I-G^2}$ are of the form $G'\hat{\otimes} I$, meaning they derive $\fK(H)\hat{\otimes}\fB(H^\rho)$, while $D\times F$ derives $\fK(H)\hat{\otimes}\fB(H^\rho)$ by definition. Thus $V$ derives $\fK(H)\hat{\otimes}\fB(H^\rho)$ as well. This concludes the proof that $s(d\times x)$ is aligned with $s(d)$ and $x$, and thus the proof of the theorem.
\end{proof}
\section{Uniform K-homology of graphs}\label{SectionSplComplex}
As a more concrete application of the formal properties of uniform K-homology we have developed in this chapter we compute the uniform K-homology of graphs of bounded geometry. We will find these uniform K-homology groups to be isomorphic to the uniformly finite graph homology groups, which we now introduce. Let $X=(V,E)$ be a directed graph of bounded geometry, meaning that only a uniformly finite number of edges are adjacent to each vertex. For an edge $e\in E$ we denote by $t(e)$ and $s(e)$ the target and source of $e$, respectively. The graph $X$ gives rise to a one-dimensional simplicial complex. The metric on $X$ is determined by the requirement that each edge be isometric to $[0,1]$, and the distance between to points be given by the length of the shortest path connecting them.
\begin{definition}
Let $X=(V,E)$ be a directed graph of bounded geometry. A \emph{uniformly finite $0$-chain} on $X$ is a formal sum $\sum_{v\in V} b_v \cdot v$, $b_v\in\ZZ$, with $|b_v|\leq C$ for some $C\geq 0$ and all $v\in V$. A \emph{uniformly finite $1$-chain} on $X$ is a formal sum $\sum_{e\in E} a_e \cdot e$ with $|a_e|\leq C$ for some $C\geq 0$ and all $e\in E$. The abelian group of all uniformly finite $k$-chains, $k=0,1$, is denoted by $C^{suf}_k(X)$. Define the differential $\partial_1: \, C^{suf}_1(X)\to C^{suf}_0(X)$ as the linear map determined by the requirement that $\partial_1(e)=t(e)-s(e)$ for every edge $e$. The \emph{simplicial uniformly finite graph homology} of $\Gamma$ is defined as
$$
H^{suf}_0(X) \; = \; C^{suf}_0(\Gamma)/\im(\partial_1) \quad , \quad H^{suf}_1(X) \; = \; \ker(\partial_1) \; . 
$$
\end{definition}
Up to canonical isomorphism the groups $H^{suf}_*(X)$ are independent of the directions of the edges. Thus one can define the simplicial uniformly finite homology of an undirected graph of bounded geometry by arbitrarily choosing directions by all of its edges.
\begin{thm}
Let $X=(V,E)$ be a graph of bounded geometry. Then, there is an isomorphism
$$
K^{u}_p(X) \; \cong \; H^{suf}_p(X)
$$
for $p=0,1$.
\end{thm}
\begin{proof}
View $V$ as a uniformly discrete subset of $X$, and consider the 6-term exact sequence of the pair $(X,V)$. By Example \ref{uniKHomDiscrete} the uniform K-homology of the uniformly discrete $V$ is given by $\ell^\infty_\ZZ(V)$ in degree $0$, and by $0$ in degree $1$. The complement $X\setminus V$ is a disjoint union of open intervals, hence its uniform K-homology is $0$ in degree $0$, and $\ell^\infty_\ZZ(E)$ in degree $1$ by the suspension isomorphism and Example \ref{uniKHomDiscrete}. Thus, the 6-term exact sequence takes the form
$$
0 \longrightarrow K^{u}_1(X) \longrightarrow \ell^\infty_\ZZ(E) \xlongrightarrow{\partial} \ell^\infty_\ZZ(V) \longrightarrow K^{u}_0(X) \longrightarrow 0 \; .
$$
It follows that $K^{u}_1(X)\cong \ker(\partial)$ and $K^{u}_0(X)\cong \ell^\infty_\ZZ(V)/\im(\partial)$. It is clear from the definitions that $\ell^\infty_\ZZ(E)$ and $C^{suf}_1(X)$ are canonically isomorphic, and so are $\ell^\infty_\ZZ(V)$ and $C^{suf}_0(X)$. Thus the proof is complete once we show that the boundary map $\partial$ on K-homology coincides with the differential $\partial_1: \, C^{suf}_1(X)\to C^{suf}_0(X)$.

To that end, recall first that the boundary map for the pair $((0,1],\{1\})$ yields an isomorphism $\partial_0: \, K_1((0,1))\to K_0(\{1\})$. Both $K_1((0,1))$ and $K_0(\{1\})$ are isomorphic to $\ZZ$, and we will identify both groups with $\ZZ$ in such a way that $\partial_0$ becomes the identity map. Next, we determine the boundary map of the pair $([0,1],\{0,1\})$. The boundary map is natural, whence there is the commutative diagram \\
\centerline{
\xymatrix{ \ZZ\cong K_1([0,1],\{0,1\}) \ar[d] \ar[r]^-\partial & K_0(\{0,1\})\ar@{=}[d] \cong \ZZ\oplus \ZZ \\
\ZZ\oplus \ZZ \cong K_1([0,\varepsilon)\cup (1-\varepsilon,1],\{0,1\})\ar[r]_-\partial & K_0(\{0,1\})\cong \ZZ\oplus \ZZ
}}
The left vertical arrow is the restriction map, it is given by the diagonal map $\ZZ\to\ZZ\oplus\ZZ$. The lower boundary map is given by $\partial = (-\partial_0) \oplus \partial_0$. Indeed, the boundary map $K_1([0,\varepsilon,),\{0\})\to K_0(\{0\})$ differs from $\partial_0$ by a reflection, and this reflection acts as the negative identity on K-homology. Overall, the boundary map $\ZZ\cong K_1([0,1],\{0,1\}) \to K_0(\{0,1\})\cong \ZZ\oplus \ZZ$ is given by $a\mapsto(-a,a)$.

Now, let $X'$ be the disjoint union $\coprod_{e\in E} [0,1]$. The space $X$ may be recovered from $X'$ by identifying the point $0$ from the interval associated to the edge $e$ with the point $1$ associated with the edge $e'$ if $t(e')=s(e)$. The quotient map $q: \, X'\to X$ is contractive, and uniformly cobounded and proper because $\Gamma$ has bounded geometry. We consider the diagram \\ 
\centerline{
\xymatrix{ K^{u}_1(X',q^{-1}(V)) \ar[d]_{q_*} \ar[r]^-\partial & K^{u}_0(q^{-1}(V))\ar[d]^{q_*}  \\
K^{u}_1(X,V)\ar[r]_-\partial & K^{u}_0(V)
}}
Note that $q$ restricts to an isomorphism $X'\setminus q^{-1}(V)\to X\setminus V$, making the left vertical arrow an isomorphism. In fact, both $K^{u}_1(X',q^{-1}(V))$ and $K^{u}_1(X,V)$ can be identified with $\ell^\infty_\ZZ(E)$ in such a way that $q_*$ becomes the identity map. The pre-image $q^{-1}(V)$ is the disjoint union $\coprod_{e\in E}\{0,1\}$. The upper boundary map is therefore simply the product of the boundary maps of each interval, hence given by $\partial (a_e)_{e\in E}=\left(-a_e,a_e\right)_{e\in E}$. The map $q_*: \, K_0^{u}(q^{-1}(V))\to K_0^{u}(V)$ maps $(a_e,b_e)_{e\in E}$ to the sequence $(c_v)_{v\in V}$ with $c_v=\sum_{t(e)=v}b_e+\sum_{s(e')=v}a_{e'}$. Composing these maps, we deduce that the boundary map $\partial: \, K_1^{u}(X,V)\to K_0^{u}(V)$ is given by
$$
\partial(a_e)_{e\in E} \; = \; (b_v)_{v\in V} \quad \mathrm{with} \quad b_v=\sum_{t(e)=v} a_e - \sum_{s(e')=v} a_{e'} \; .
$$
Under the respective identifications with $C^{suf}_1(X)$ and $C^{suf}_0(X)$, this is precisely the boundary map defining uniformly finite graph homology.
\end{proof}
The above proof may be carried out verbatim for analytic K-homology, the only difference being that the zeroth K-homology of $V$ is $\ZZ^V$, and the first K-homology of $X\setminus V$ is $\ZZ^{E}$. The groups correspond to the spaces of not necessarily uniformly finite $0$- and $1$-cycles on $X$, respectively. The homology of these cycles (using the same definition for the differential) is called the \emph{locally finite graph homology} $H^{lf}_*(X)$ (or Borel-Moore graph homology). Thus, the analogous arguments show that there are isomorphisms
$$
K_p(X) \; \cong \; H^{lf}_p(X)
$$
for $p=0,1$, where $X$ again denotes the one-dimensional simplicial complex associated to the directed graph $X$ of bounded geometry. This isomorphism is likely known already, but the author could not find reference to it in the literature.
\begin{example}[The real line]\label{UniKHomRealLine}
Consider $\RR$ with its standard metric as the graph with vertex set $\ZZ$, and edges given by the intervals $e_n:=[n,n+1]$, directed via $t(e_n)=n+1$ and $s(e_n)=n$. The boundary map $\partial: \, \ell^\infty_\ZZ(E) \to \ell^\infty_\ZZ(V)$ is then given by
$$
\partial (a_{e_n})_{n\in\ZZ} \; = \; (a_{e_n}-a_{e_{n+1}})_{n\in\ZZ} \; .
$$
Any element in the kernel of $\partial$ is constant, whence $K_1^{u}(\RR)=\ZZ$. Moreover, two elements $(b_n)_{n\in\ZZ}$ and $(b_n')_{n\in\ZZ}$ define the same class in $K^{u}_0(\RR)\cong H^{suf}_0(\RR)$ if and only if there exists $(a_{e_n})_{n\in\ZZ}$ such that $b_n+a_{e_n}=b_n'+a_{e_{n+1}}$ for all $n\in\ZZ$. In particular the chain $(b_n)_{n\in\ZZ}$ is in the image of the boundary map if and only if there exists a bounded sequence $(a_{e_n})_{n\in\ZZ}$ -identified with a sequence $(a_n)_{n\in\ZZ}$- such that $b_n=a_n-a_{n+1}$. Note that this implies that 
$$
\sum_{k=0}^N b_{n+k} = a_n - a_{n+N} \quad \; \quad \sum_{k=0}^N b_{n-k} = a_n - a_{n-N} 
$$
for all $n,k,N\in \NN$. This leads to the description $K^{u}_0(\RR)\cong \ell^\infty_\ZZ(\ZZ)/S(\ZZ)$, where $S(\ZZ)\subset \ell^\infty_\ZZ(\ZZ)$ is the subgroup of all $(b_n)_{n\in\ZZ}$ satisfying the following summability condition: For all $k\in\ZZ$ there exists $C\geq 0$ such that for all $N\in\NN$ it holds that
$$
\left| \sum_{n=0}^N b_{k+n} \right| + \left| \sum_{n=0}^N b_{k-n} \right| \; \leq \; C \; .
$$  
Note that this condition is satisfied for all $k\in\ZZ$ if and only if it is satisfied for a single $k\in\ZZ$. Thus, we have computed that
$$
K^{u}_0(\RR)=\ell^\infty_\ZZ(\ZZ)/S(\ZZ) \quad , \quad K^{u}_1(\RR) = \ZZ \; .
$$
\end{example}
\begin{example}[Closed and open half-line]
We turn to $\RR_{\geq 0}$ with the standard metric. It can be viewed as a subgraph of $\RR$, where the latter is equipped with a graph structure as in the previous example.  Analogously to the the real line, one also computes that $K^{u}_0(\RR_{\geq 0})=\ell^\infty_\ZZ(\NN_0)/S(\NN_0)$, where $S(\NN_0)$ contains those elements $(b_n)_{n\in\NN_0}$ satisfying the following summability condition: There exists $C\geq 0$ such that for all $N\in\NN_0$ it holds that
$$
\left| \sum_{n=0}^N b_{n} \right| \; \leq \; C \; .
$$ 
Moreover, any uniformly finite $1$-cycle on $\RR_{\geq 0}$ must vanish, since it must take the value $0$ on the left-most edge, and this vanishing then propagates along the half-line. Thus, we have computed that
$$
K^{u}_0(\RR_{\geq 0})=\ell^\infty_\ZZ(\NN_0)/S(\NN_0) \quad , \quad K^{u}_1(\RR_{\geq 0}) = 0 \; .
$$
Let us also compute the uniform K-homology of the open half-line $\RR_{>0}$. The 6-term exact sequence of the pair $(\RR_{\geq 0},\{0\})$ reduces to the exact sequence
$$
0 \longrightarrow K^{u}_1(\RR_{>0}) \xlongrightarrow{\partial} K^{u}_0(\{0\}) \xlongrightarrow{\iota_*} K^{u}_0(\RR_{\geq 0}) \xlongrightarrow{\pi_*} K^{u}_1(\RR_{>0}) \longrightarrow 0 \; .
$$
The map $\{0\} \hookrightarrow \RR_{\geq 0}$ factors as $\{0\} \hookrightarrow \NN_0 \hookrightarrow \RR_{\geq 0}$. The induced map $\iota_*: \, K_0^{u}(\ZZ)\to K_0^{u}(\RR_{\geq 0})$ equals the composition $\ZZ\hookrightarrow \ell^\infty_\ZZ(\NN_0) \twoheadrightarrow \ell^\infty_\ZZ(\NN_0)/S(\NN_0)$. This composition is the zero-map, as any sequence with support $\{0\}$ lies in $S(\NN_0)$. It follows that the maps $\partial$ and $\pi_*$ in the above exact sequence are isomorphisms, whence
$$
K^{u}_0(\RR_{>0}) \cong \ell^\infty_\ZZ(\NN_0)/S(\NN_0) \quad , \quad K^{u}_1(\RR_{>0}) \cong \ZZ \; .
$$
\end{example}
It turns out that $\RR_{\geq 0}$ is the essential example of a graph with vanishing first uniform K-homology. The vanishing of the zeroth uniform K-homology for graphs is in turn uniquely characterized by non-amenability.
\begin{prop}\label{AmenableGraphs}
Let $X$ be a connected graph of bounded geometry. 
\begin{itemize}
\item[(i)] $K_0^{u}(X)=0$ if and only if $X$ is not amenable.
\item[(ii)] $K_1^{u}(X)=0$ if and only $X$ is a tree with at most one end.
\end{itemize}
\end{prop}
\begin{proof}
Ad (i): Recall \v{S}pakula's result discussed in Example \ref{ExGrahKhom} that $X$ is non-amenable if and only if its fundamental class $\mathbf{S}\in K^{u}_0(X)$ vanishes. This implies that connected graphs with vanishing zeroth uniform K-homology are non-amenable. For the reverse implication we show that $\mathbf{S}=0$ if and only if $K^{u}_0(X)=0$. The class $\mathbf{S}$ is represented by a Fredholm module assigning an index $1$-operator to each vertex of $X$. Going through the proof of the isomorphism $K^{u}_0(X) \cong H^{suf}_0(X)$ it can be seen that that isomorphism maps $\mathbf{S}$ to the class $\mathbf{1}\in H^{suf}_0(X)$ of the $0$-cycle $\sum_v v$. Thus it suffices to prove that $\mathbf{1}=0$ implies $H^{suf}_0(X)=0$. We prove this in analogy to \cite[Proposition 2.3]{Block1992}. Assume that there exists $c=\sum_{e} b_e e \in C^{suf}_1(X)$ with 
$$
\partial c = \mathbf{1} = \sum_v v \quad \Leftrightarrow \quad 1 =\sum_{t(e)=v} b_e - \sum_{s(e)=v}b_v \quad \forall \, v\; .
$$
By potentially reversing the orientation of edges we may assume that $b_e\geq 0$, this does not change the validity of $\partial c=\mathbf{1}$. Consider some vertex $v_0$. There exists some edge $e_0$ with target $v_0$ such that $b_{e_0}>0$. Let $v_1$ denote the source of that edge. Again there must exist an edge $e_1$ with target $v_1$ such that $b_{e_1}>0$. Continuing this process indefinitely produces an infinite chain of edges $(e_0,e_1,\cdots)$ such that $s(e_n)=t(e_{n+1})$ and $b_{e_n}>0$. The latter property implies that there is a uniform bound on how often any particular edge can appear in this chain. In particular $t_{v_0} = \sum_{n=0}^\infty e_n$ is a well-defined element of $C^{suf}_1(X)$. Using the former property we see that $\partial t_{v_0} = t(e_0) = v_0$. Moreover since the edges appearing in the $t_v$ have been chosen from the uniformly bounded chain $c$ it is not hard to see that any given edge can only be part of a uniformly finite number of chains $t_v$, $v\in V$. Thus the element $\sum_v t_v$ is a well-defined element of $C^{suf}_1(X)$. Then, an arbitrary element $0$-chain $\sum_{v} a_v v \in C^{suf}_0(X)$ the boundary of $\sum_v a_v t_v \in C^{suf}_1(X)$, concluding the proof that $H^{suf}_0(X)=0$.

Ad (ii): Assume $K^{u}_1(X)=0$. Suppose first that $X$ contains a circuit, i.e. a closed path containing no vertex or edge more than once. Pick an edge $e$ from this circuit, and assign it the value $1$. The next edge $e'$ in the circuit satisfies either $s(e')=t(e)$ or $t(e')=t(e)$. In the former case, assign it the value $1$, in the latter case assign it $-1$. Continue around the circuit in this manner, assigning each edge either the same value as the previous edge, or the negative of the value of the previous edge, depending whether their orientation is aligned or opposite. All edges of $X$ that do not lie on the circuit are assigned $0$. The resulting element of $\ell^\infty_\ZZ(E)$ constitutes a non-zero element of the kernel of $\partial$. Thus, if $K^{u}_1(X)=0$, then $X$ contains no circuit, meaning $X$ is a tree. 

Now let $X$ be a tree. Let $\sum_e a_e e$ be an element of $\ker(\partial_1)$, and assume that $a_e\neq 0$ for some edge $e=[x,y]$. Then, $x$ must meet an edge $e'=[x,x']$ with $a_{e'}\neq 0$. Consequently, $x'$ must meet an edge $e''=[x',x'']$ with $a_{e''}\neq 0$, and this edge $e''$ cannot meet $x$ because $X$ contains no circuits. Continuing this procedure one obtains a semi-infinite path $(x,x',x'',\cdots)$ that defines an end of $X$. However, in the same way one also produces a semi-infinite path $(y,y',y'',\cdots)$, and it is clear that the two paths define different ends. This shows that $K^{u}_1(X)\neq 0$ implies $X$ to have at least two ends. Conversely, any tree with at least two ends contains a bi-infinite path, and assigning each edge on that path the value $\pm 1$ in a way compatible with orientations produces a non-zero element of $K^{u}_1(X)$.
\end{proof}
Uniformly finite graph homology exists in a wider cosmos of uniformly finite homology theories. Namely, it coincides with the simplicial uniformly finite homology of the one-dimensional simplicial complex corresponding to the graph. The simplicial uniformly finite homology of a simplicial complex of bounded geometry is defined as follows: Let $C^{suf}_n(X)$ denote the space of uniformly bounded sums $\sum_{\sigma} a_\sigma \sigma$, $a_\sigma\in\ZZ$, with $\sigma$ ranging over the $n$-simplices of $X$. Then there is the usual simplicial boundary map $\partial_n: \, C^{suf}_n(X)\to C^{suf}_{n-1}(X)$ sending a simplex $\sigma$ to the signed sum of its faces. Then $H^{suf}_n(X)=\ker(\partial_n)/\im(\partial_{n-1})$.

There is also the uniformly finite homology $H^{uf}_*(Y)$ defined on any metric space $Y$, which was introduced by Block and Weinberger \cite{Block1992}. It is a quasi-isometry invariant of the space. One way to define it is via simplicial uniformly finite homology of coarsenings of $Y$ via Rips complexes. If $Y$ is uniformly contractible, then $H^{suf}_*(Y)$ and $H^{uf}_*(Y)$ coincide. We refer to \cite{Block1992} for an account of geometric applications of uniformly finite homology that touches in particular on index theory. A detailed introduction to (simplicial) uniformly finite homology can be found in \cite{Diana2015}. 

Let us also point out a recent note by Manuilov \cite{Manuilov2024}, where maps
$$
H^{suf}_*(X) \longrightarrow K_*\left(C^*_u(X) \right) \quad \mathrm{and} \quad H^{lf}_*(X) \longrightarrow K_*\left(C^*(X) \right)
$$
were constructed for $X$ a graph of bounded geometry. Sorting through the constructions should likely reveal that these maps correspond to the compositions
$$
H^{suf}_*(X) \xlongrightarrow{\sim} K^{u}_*(X) \xlongrightarrow{\Ind} K_*(C^*_u(X))
$$
and
$$
H^{lf}_*(X) \xlongrightarrow{\sim} K_*(X) \xlongrightarrow{\mathrm{Ind}} K_*(C^*(X)) \; ,
$$
where $\Ind$ and $\mathrm{Ind}$ are the (uniform) coarse index maps.
\begin{example}[General simplicial complexes]
Let $X$ be an $n$-dimensional simplicial complex of bounded geometry, and let $X^{n-1}$ be its $(n-1)$-skeleton. The complement $X\setminus X^{n-1}$ is a disjoint union of open $n$-simplices labeled by the set $F$ of $n$-simplices in $X$. Thus $K_p^{u}(X\setminus X^{n-1})$ is isomorphic to $\ell^\infty_\ZZ(F)$ if $p=n\, \mathrm{mod} \, 2$, and $0$ otherwise (this is essentially covered by \cite[Lemma 3.34]{Engel2019}). Therefore the long exact sequence reads as
$$
0 \longrightarrow K_{-n}^{u}(X^{n-1}) \longrightarrow K_{-n}^{u}(X) \longrightarrow \ell^\infty_\ZZ(F) \xlongrightarrow{\partial} K_{-(n-1)}^{u}(X^{n-1}) \longrightarrow K^{u}_{-(n-1)}(X) \longrightarrow 0 \; .
$$ 
We conclude that $K_{-n}^{u}(X)$ contains $K^{u}_{-n}(X^{n-1})$ as a subgroup, while $K^{u}_{-(n-1)}(X)$ is a quotient of $K_{-(n-1)}^{u}(X^{n-1})$. More concretely we have short exact sequences
\begin{align*}
0 \longrightarrow K_{-n}^{u}(X^{n-1}) \longrightarrow K_{-n}^{u}(X) \longrightarrow \ker(\partial) \longrightarrow 0 \\
0 \longrightarrow \im(\partial) \xlongrightarrow{\partial} K_{-(n-1)}^{u}(X^{n-1}) \longrightarrow K^{u}_{-(n-1)}(X) \longrightarrow 0 \; .
\end{align*}
If $\partial: \, \ell^\infty_\ZZ(F)\to K^{u}_{-(n-1)}(X^{n-1})$ was known to us, we could calculate $K^{u}_{-(n-1)}(X^{n-1})$ directly. While $K^{u}_{-n}(X)$ would not be immediately accessible, we could tensor the former short exact sequence above with a field of our choice, say the rationals. We would then obtain a (non-canonical) splitting
$$
K^{u}_{-n}(X) \otimes \QQ \cong \left( K^{u}_{-n}(X^{n-1}) \otimes\QQ \right) \oplus \left( \ker(\partial) \otimes \QQ \right) \; .
$$
In this way one can get surprisingly close an isomorphism rational uniform K-homology and rational simplicial uniformly finite homology. Indeed, suppose by induction that this isomorphism holds for $X^{n-1}$. Then, if we could rationally identify the K-homological boundary map with the simplicial boundary maps, it could easily be verified that via the above sequences (splitting the latter one) $K^{u}_*(X)$ would rationally coincide with $H^{suf}_{*}(X)$. There are two issues with this: For one the algebraic tensor product with $\QQ$ is not the correct way to obtain uniform K-homology with rational coefficients, one should instead take the completed topological tensor product. This is discussed by Engel \cite{Engel2015}, \cite{Engel2019}. Then, it is not at all clear that the above sequence still splits in an appropriate sense. But even if it does, the non-canonical nature we one would obtain the isomorphism via a choice of splitting would make it practically difficult to identify the K-homological boundary map with the simplicial one in the next induction step. This second problem might be solvable by carrying this additional assumption on the boundary maps as part of the induction hypothesis. This is done in a simpler scenario in the proof of Proposition \ref{KHomExpansion} below. In any case one ultimately wants not only the abstract isomorphism statement between rational (or real or complex) uniform K-homology and simplicial uniform K-homology or a related homology theory, but instead the existence of a Chern character implementing this isomorphism whose properties one can study. For manifolds with bounded geometry Engel has constructed such a Chern character \cite{Engel2015}, but the general construction remains to be done.
\end{example}
\begin{example}[The hyperbolic plane]\label{UKHomHyperbolic}
Let us study the uniform K-homology of the hyperbolic plane $\HH^2$. We choose a triangulation of $\HH^2$ such that the metric obtained from the barycentric coordinates is quasi-isometric to the hyperbolic metric. Let $Z\subset \HH^2$ be the 1-skeleton of this triangulation. Since $\HH^2$ is non-amenable, so is $Z$, and since $Z$ is a graph of bounded geometry, this implies $K^{u}_0(Z)=0$. Moreover, $\HH^2\setminus Z$ is a disjoint union of open 2-simplices, indexed by a set $F$. Thus, $K^{u}_0(\HH^2\setminus Z)=\ell^\infty_\ZZ(F)$ and $K^{u}_1(\HH^2\setminus Z)=0$. The 6-term exact sequence of the pair $(\HH^2,Z)$ therefore reduces to
$$
0\longrightarrow K^{u}_0(\HH^2) \longrightarrow \ell^\infty_\ZZ(F)\xlongrightarrow{\partial} H^{suf}_1(Z) \longrightarrow K^{u}_1(\HH^2) \longrightarrow 0 \; .
$$
The homology $H^{suf}_1(Z)$ equals the kernel of $\partial_1: \, \ell^\infty_\ZZ(Z) \to \ell^\infty_\ZZ(Z^0)$, with $Z^0$ being the corresponding vertex set. Indeed, it can be checked that the boundary map $\partial: \, \ell^\infty_\ZZ(F)\to H^{suf}_1(Z)=\ker(\partial_1)$ coming from the 6-term exact sequence in uniform K-homology coincides with the boundary map $\partial_2: \, \ell^\infty_\ZZ(F) \to \ell^\infty_\ZZ(Z)$ used in the definition of uniformly finite simplicial homology.\footnote{It is straight-forward to check that the boundary map $\ZZ\cong K_0(\mathring{\Delta}^2)\to K_1(\partial \Delta^2)\cong \ZZ$ is an isomorphism. Identifying $K_0(\mathring{\Delta}^2)$ with the free group on one 2-simplex, and $K_1(\partial \Delta^2)$ with the subgroup generated by the unique cycle on $\partial \Delta^2$, this isomorphism can be made to correspond to the boundary map on simplicies. Gluing according to the simplicial structure as in our discussion of graphs identifies the uniform K-homological boundary map with the simplicial boundary map.} It follows that $K^{u}_0(\HH^2)=\ker(\partial_2)=H^{suf}_2(\HH^2)$, and $K^{u}_0(\HH^2)=\ker(\partial_1)/\im(\partial_2)=H^{suf}_1(\HH^2)$. Moreover, $\HH^2$ is uniformly contractible, so that its uniformly finite simplicial homology and its uniformly finite homology agree. We conclude that
$$
K_0^{u}(\HH^2)\cong H^{uf}_2(\HH^2) \; , \; K_1^{u}(\HH^2)\cong H^{uf}_1(\HH^2) \; .
$$
The uniform homology groups with real coefficients have been discussed in the literature (see for example Section 8 of \cite{Block1997} and the references given therein). In dimension $1$ one may use intersections of cycles with geodesics to produce an inclusion of $H^{uf}_0(\RR;\RR)$\footnote{The first $\RR$ denoting the space, the second the coefficient group.} into $H^{uf}_1(\HH^2;\RR)$. Since the former is infinite-dimensional, so is the second. On the other hand $H^{uf}_2(\HH^2;\RR)$ turns out to be one-dimensional. We expect the same to be true with integer coefficients as well. We also note that since $H^{uf}_*$ is a quasi-isometry invariant, and $\HH^2$ is quasi-isometric to any surface group $\Gamma_g$ of genus $\geq 2$, the groups $K^{u}_*(\HH^2)$ are isomorphic to the first and second uniform homology groups of $\Gamma_g$. The latter have a description via group homology, indeed it holds that $H^{uf}_*(\Gamma_g)\cong H_*(\Gamma_g;\ell^\infty_\ZZ(\Gamma_g))$. This is a general fact about finitely generated groups as originally observed in \cite{Wright2010}, see also \cite{Diana2015} for a detailed proof.

It should be pointed out that Engel has provided a Chern isomorphism for manifolds of bounded geometry between uniform K-homology and uniform de Rham-homology in case both carry complex coefficients \cite{Engel2015}. In the case of the hyperbolic plane the latter is isomorphic to uniformly finite homology with complex coefficients. Thus, the novelty of this example is that the isomorphism holds already with integer coefficients.
\end{example}
Much of the discussion in Example \ref{UKHomHyperbolic} can be carried out for an arbitrary non-amenable two-dimensional simplicial complex of bounded geometry. We can combine these arguments with induction to get a more general statement identifying the uniform K-homology of a simplicial complex with its simplicial uniformly finite homology, provided that the lower-dimensional $H^{suf}_*$-groups vanish.
\begin{prop}\label{KHomExpansion}
Let $X$ be an $n$-dimensional simplicial complex of bounded geometry, $n\geq 2$. If $H^{suf}_p(X)=0$ for $p=0,\cdots,n-2$, then
$$
K^{u}_{-n}(X) \cong H^{suf}_n(X) \quad , \quad K^{u}_{-(n-1)}(X) \cong H^{suf}_{n-1}(X) \; .
$$
\end{prop}
\begin{proof}
We prove this by induction over $n$. The case of $n=2$ is done exactly as in Example \ref{UKHomHyperbolic} above. Thus assume the claim is proved for dimensions $\leq n-1$, and let $X$ be an $n$-dimensional simplicial complex of bounded geometry with vanishing simplicial uniformly finite homology up to degree $n-2$. Let $Z=X^{n-1}$ denote the $(n-1)$-skeleton of $X$. Note that 
$$
H^{suf}_p(Z) = H^{suf}_p(X) = 0 
$$ 
for $p=0,\cdots, n-2$. We conclude from the induction hypothesis that $K_n^{u}(Z)=0$ and $K_{n-1}^{u}(Z)=H^{suf}_{n-1}(Z)$. Thus the long exact sequence of the pair $(X,Z)$ in uniform K-homology takes the form
$$
0 \longrightarrow K_{-n}^{u}(X) \longrightarrow \ell^\infty_\ZZ(F) \xlongrightarrow{\partial} H^{suf}_{n-1}(Z) \longrightarrow K^{u}_{-(n-1)}(X) \longrightarrow 0 \; .
$$
Suppose that we knew that the K-homological boundary map coincided with the simplicial boundary map 
$$
\partial_n: \; \ell^\infty_\ZZ(F) = C^{suf}_n(X) \longrightarrow \ker(\partial_{n-1})=H^{suf}_{n-1}(Z)\subseteq C^{suf}_{n-1}(X) \; .
$$
Then the claim would follow immediately, since by the above exact sequence we would have
$$
K^{u}_{-n}(X) \cong \ker(\partial_n)=H^{suf}_n(X) \; , \; K^{u}_{-(n-1)}(X) \cong \ker(\partial_{n-1})/\im(\partial_n) = H^{suf}_{n-1}(X) \; .
$$
That this equality of boundary maps holds can be arranged as part of the induction. Indeed, it is certainly true for $n=2$. Then, suppose it holds that the isomorphism up to dimension $n-1$ have been arranged for this to be true; we have to show that it also holds under the isomorphism in dimension $n$ constructed just now. To that end let $Y$ be an $(n+1)$-dimensional complex satisfying the assumptions of the proposition, and let $X=Y^n$ be its $n$-skeleton. The boundary maps $\partial: \, K^{u}_{-(n+1)}(Y\setminus X)\to K_{-n}^{u}(X)$ factors as \\
\centerline{\xymatrix{
 & K^{u}_{-(n+1)}(\coprod_F \mathring{\Delta}^{n+1}) \ar@{=}[d] \ar[r]^\partial & K^{u}_{-n}(\coprod_F \partial\Delta^{n+1}) \ar[d] \\
K^{u}_{-(n+1)}(Y\setminus X) \ar@{=}[r] & K^{u}_{-(n+1)}(\coprod_F \mathring{\Delta}^{n+1}) \ar[r]_\partial & K_{-n}^{u}(X)
}}
where the right-vertical arrow is the gluing map, and the set $F$ labels the $(n+1)$-simplices in $Y$. Now, $K_{-(n+1)}(\mathring{\Delta}^{n+1})\cong\ZZ$, and $K_{-n}(\partial\Delta^{n+1})$ is either isomorphic to $H_n(\partial\Delta^{n+1})\cong \ZZ$ if $n$ is odd, or to $H_n(\partial\Delta^{n+1})\oplus H_0(\partial\Delta^{n+1})\cong \ZZ\oplus\ZZ$ if $n$ is even. In the former case the boundary map $\partial: \, K_{-(n+1)}(\mathrm{\Delta}^{n+1}) \to K_{-n}(\partial\Delta^{n+1})$ is an isomorphism, which we may take to identify $1\in\ZZ$ with the fundamental class in $H_n(\partial\Delta^{n+1})$. In the latter case it is an isomorphism onto the $H_n(\partial\Delta^{n+1})$-summand, which we may again take to identify $1$ with the fundamental class. These identifications are compatible with disjoint unions, so that
$$
K_{-(n+1)}^{u}(\coprod_F \mathring{\Delta}^{n+1}) \cong \ell^\infty_\ZZ(F) \xlongrightarrow{\partial} H^{suf}_n(\coprod_F\partial \Delta^{n+1}) \subseteq K^{u}_{-n}(\coprod_F\partial \Delta^{n+1})
$$
maps the generator on each simplex to the fundamental class of its boundary. Now, the isomorphism $K^{u}_{-n}(X)$ we constructed above is implemented by the restriction map $K_{-n}^{u}(X)\to K_{-n}^{u}(X\setminus X^{n-1})\cong C^{suf}_n(X)$, and thus compatible with the gluing map. This allows us to conclude that $\partial: \, K^{u}_{-(n+1)}(Y\setminus X) \to K_{-n}^{u}(X)$ is indeed the simplicial boundary map under the respective isomorphisms to spaces of chains. This concludes the proof. 
\end{proof}
In a recent work Bottinelli and Kaiser \cite{Bottinelli2020} introduce expansion properties for simplicial complexes generalizing non-amenability that provide conditions for the vanishing of the $H^{suf}_*$-groups. Having $n$-dimensional $H^{suf}$-expansion in the sense of \cite{Bottinelli2020} implies that $H^{suf}_p=0$. Zero-dimensional $H^{suf}$-expansion is simply non-amenability. One-dimensional $H^{suf}$-expansion is equivalent to vanishing $H^{suf}_1=0$. It appears to be unknown if this is also true in higher dimensions.

Again, we anticipate a Chern isomorphism theorem for all simplicial complexes of bounded geometry, which would in particular imply the statement of the previous proposition for real or complex coefficients. Thus the last proposition is ultimately relevant only in so far that the statement already holds with integer coefficients.
\chapter{Index maps on relative uniform K-homology}\label{ChapterIndexMaps}
The coarse index map from analytic K-homology to the K-theory of the Roe algebra is fundamental to the coarse-geometric approach to index theory on non-compact spaces. There is also a uniform counterpart, mapping from uniform K-homology to the K-theory of the uniform Roe algebra. This uniform coarse index map has been been applied in connection to the Baum-Connes conjecture with coefficients \cite{Spakula2009}, \cite{Engel2019}, and as an obstruction to uniformly positive scalar curvature \cite{Engel2014}. Uniform Roe algebras and the uniform coarse index map will be reviewed in Section 3.1.

The purpose of this chapter is the discussion of uniform coarse index maps on relative uniform K-homology. We present two possibilities, though it will turn out that one is strictly stronger than the other. The first approach, which we will develop in Section 3.2, is a map to a \emph{relative uniform Roe algebra}. If $X$ is a metric space and $Z\subseteq X$ a closed subspace, then this relative uniform Roe algebra fits into a long exact sequence
$$
\cdots \longrightarrow K_*\left( C^*_u(Z) \right) \longrightarrow K_*\left( C^*_u(X) \right) \longrightarrow K_*\left( C^*_u(X,Z) \right) \xlongrightarrow{\partial} K_{*-1}\left( C^*_u(Z) \right) \longrightarrow \cdots \; .
$$
We construct a uniform coarse index map $\Ind: \, K_*^{u}(X,Z)\to K_*(C^*_u(X,Z))$. Intuitively, this index is obtained by forgetting everything close to the subspace $Z$, so that the remainder is invertible modulo $C^*_u(X,Z)$. The index in $K_*(C^*_u(X,Z))$ then measures the failure to be actually invertible. This index map fits into a commutative diagram\\
\centerline{\xymatrix{
K_*^{u}(Z) \ar[r] \ar[d]_{\Ind} & K^{u}_*(X) \ar[r] \ar[d]_{\Ind} & K^{u}_*(X,Z) \ar[r]^\partial \ar[d]_{\Ind} & K_{*-1}^{u}(Z) \ar[d]_{\Ind} \\
K_*\left(C^*_u(Z)\right) \ar[r] & K_*\left(C^*_u(X)\right) \ar[r] & K_*\left(C^*_u(X,Z)\right) \ar[r]_\partial & K_{*-1}\left(C^*_u(Z)\right)
}}
This is conceptually satisfying: Both (relative) uniform K-homology and K-theory of (relative) uniform Roe algebras form homology theories on suitable categories of metric spaces, and the uniform coarse index map provides a transformation from one theory to the other that respects long exact sequences. Unfortunately the relative group $
K_*\left(C^*_u(X,Z)\right)$ turns out to be trivial whenever $Z$ is uniformly coarsely equivalent to $X$. This happens in particular when $Z$ is the boundary of a finite-width manifold $X$ with boundary. In this case -which we would very much like to study in index theory!- we can therefore not expect any information whatsoever from the uniform coarse index in $K_*\left(C^*_u(X,Z)\right)$.

The second index map we will construct on relative uniform K-homology is in some sense closer to the non-relative uniform coarse index map. In Section 3.3 we will introduce finite-propagation counterparts $D^*_u(X)$ and $D^*_u(X,Z)$ -called the \emph{(relative) structure algebras}- to the algebras $\fD^{u}(X)$ and $\fD^{u}(X,Z)$. A finite-propagation normalization in conjunction with Paschke duality yields that
$$
K^{u}_*(X,Z) \cong K_{*+1}\left( D^*_u(X)/D^*_u(X,Z) \right) \; .
$$
Then we define the \emph{relative uniform index map} as the K-theoretical boundary map
$$
\relInd: \; K^{u}_*(X,Z) \cong K_{*+1}\left( D^*_u(X)/D^*_u(X,Z) \right) \xlongrightarrow{\partial} K_*\left(D^*_u(X,Z) \right) \; .
$$
In Section 3.4 we derive (among other things) an algebraic Mayer-Vietoris sequence expressing $K_*\left(D^*_u(X,Z) \right)$ in terms of the K-theory group of $C^*_u(Z)$, $C^*_u(X)$ and $D^*_u(Z)$. While $C^*_u(Z)$ and $C^*_u(X)$ measure only the uniform coarse geometry of $X$ and $Z$, the K-theory of the algebra $D^*_u(Z)$ depends on the metric geometry of $Z$. Thus the K-theory of the relative structure algebra is a hybrid object combining the uniform coarse geometry of $X$ with the metric geometry of $Z$. The relative uniform index inherits this hybrid nature. We will be able to fully appreciate its secondary nature in Chapter 7.

In Section 3.4 we also prove a geometric Mayer-Vietoris sequence for $K_*\left(D^*_u(X,Z) \right)$. After defining the relative uniform index map in Section 3.5 it is the main result of that section that the map $\relInd$ provides a transformation from the Mayer-Vietoris sequence in uniform K-homology to that in K-theory of the relative uniform structure algebras. This will allow us to deduce partitioned-manifold index theorems in Chapter 6. In Section 3.5 we also investigate the relation between the algebras $D^*_u(X,Z)$ and $C^*_u(X,Z)$, and between the corresponding index maps. We will see that the index map to $K_*(C^*_u(X,Z))$ factors through that to $K_*(D^*_u(X,Z))$. Indeed, the whole transformation from the long exact sequence in uniform K-homology to that in K-theory of the uniform Roe algebras factors through the algebraic Mayer-Vietoris sequence for $K_*\left(D^*_u(X,Z) \right)$, though in an unfortunately non-commuting way. Thus the valuable secondary information preserved by the index map to the relative uniform structure algebra comes at the expense of unfavorable formal properties.
\section{Uniform Roe algebras}
We begin with a review of uniform Roe algebras. Originally, these were defined for certain discrete metric spaces.\footnote{Recall our convention that any metric space is locally compact and second-countable. In particular any discrete metric space is assumed to be countable.}
\begin{definition}
A metric space $Y$ is called \emph{uniformly discrete} if $\inf \{ d(x,y) \, | \, x,y\in Y, \, x\neq y \}>0$. It is said to have \emph{bounded geometry} if $\sup_{y\in Y} |B_r(y)|<\infty$ for all $r\geq 0$. 

A subset $Y$ of a metric space $X$ is called a \emph{quasi-lattice} of $X$ if $Y$ is uniformly discrete and of bounded geometry, and there exists $c\geq 0$ such that $X=N_c(Y)$. If $X$ admits a quasi-lattice, it is said to have \emph{(coarsely) bounded geometry}.
\end{definition}
Let $Y$ be a uniformly discrete metric space of bounded geometry, and consider the Hilbert space $\ell^2(Y)$. A bounded operator on $\ell^2(Y)$ may be viewed as a potentially infinite matrix $T=(t_{xy})_{x,y\in Y}$, $t_{xy}\in \CC$. Such a matrix is said to have \emph{finite propagation} if $t_{xy}=0$ whenever $d(x,y)\geq R$ for some $R\geq 0$. Moreover, it will be said to have \emph{uniformly bounded coefficients} if $\sup_{x,y}|t_{xy}|<\infty$. If $W\subseteq Y$ is a subspace, then a bounded operator $T$ on $\ell^2(Y)$ is said to be \emph{supported near $W$} if there exists $c\geq 0$ such that $t_{xy}=0$ whenever $x$ or $y$ have distance $\geq c$ from $W$.
\begin{definition}\label{DefUniRoe}
Let $Y$ be a uniformly discrete space of bounded geometry. The \emph{uniform Roe algebra} $C^*_u(Y)$ of $Y$ is defined as the norm-closure in $\fB(\ell^2(Y))$ of finite-propagation operator with uniformly bounded coefficients.

Let $W\subseteq Y$ be a subspace. Then $C^*_u(W\subseteq Y)$ is defined as the norm-closure in $\fB(\ell^2(Y))$ of finite-propagation operator with uniformly bounded coefficients that are supported near $W$.
\end{definition} 
\begin{example}\label{ExUniRoeIntegers}
Let $\Gamma$ be a finitely generated group equipped with a word metric. Then $C^*_u(\Gamma)$ is isomorphic to the reduced crossed product $\ell^\infty(\Gamma)\rtimes_r \Gamma$ (see for example \cite[Section 5.1]{Brown2008}). Consider for example $\Gamma= \ZZ$. We wish to compute the K-theory of $C^*_u(\ZZ)\cong \ell^\infty(\ZZ)\rtimes_r\ZZ$. We appeal to the Pimsner-Vioculescu sequence (see for example \cite{Blackadar1998})\\
\centerline{\xymatrix{
K_0(\ell^\infty(\ZZ)) \ar[r]^{I-\alpha_*} & K_0(\ell^\infty(\ZZ)) \ar[r] & K_0(\ell^\infty(\ZZ)\rtimes_r \ZZ) \ar[d] \\
K_1(\ell^\infty(\ZZ)\rtimes_r \ZZ) \ar[u] & K_1(\ell^\infty(\ZZ)) \ar[l] & K_1(\ell^\infty(\ZZ)) \ar[l]^{I-\alpha_*}
}}
where $\alpha: \, \ell^\infty(\ZZ)\to \ell^\infty(\ZZ)$ is the action of $1\in\ZZ$, i.e. the right shift. It holds that $K_0(\ell^\infty(\ZZ)) \cong \ell^\infty_\ZZ(\ZZ)$ and $K_1(\ell^\infty(\ZZ))=0$ (see for example \cite{Aweygan2020}). Then, $\alpha_*: \, K_0(\ell^\infty(\ZZ)) \to K_0(\ell^\infty(\ZZ))$ is the right shift $S$ on $\ell^\infty_\ZZ(\ZZ)$. We conclude that 
$$
K_0(\ell^\infty(\ZZ)\rtimes_r \ZZ) \cong \mathrm{coker}(I-S) \quad , \quad K_1(\ell^\infty(\ZZ)\rtimes_r \ZZ) \cong \ker(I-S) \; .
$$
Comparison with the computation of $K^{u}_*(\RR)$ from Example \ref{UniKHomRealLine} reveals that $\ker(I-S)$ consists of constant sequences and is hence isomorphic to $\ZZ$, and that $\im(I-S)$ is the subgroup $S(\ZZ) \subseteq \ell^\infty_\ZZ(\ZZ)$. We conclude that
$$
K_0(C^*_u(\ZZ)) \cong \ell^\infty_\ZZ(\ZZ)/S(\ZZ) \quad , \quad K_1\left( C^*_u(\ZZ) \right) \cong \ZZ \; .
$$
In particular $K_0(C^*_u(\ZZ))$ is infinitely generated, which is in stark contrast to the Roe algebra and the group $C^*$-algebra. Observe also that $K_*(C^*_u(\ZZ))\cong K_*^{u}(\RR)$. The discussion in the remainder of this section will reveal this to be no coincidence.
\end{example}
Coarsely equivalent spaces do not generally give rise to isomorphic uniform Roe algebras. The Morita equivalence class of $C^*_u(Y)$ \emph{does} however depend only on the coarse equivalence class of $Y$, and indeed characterizes it fully.
\begin{prop}[{\cite{Brodzki2006}, \cite{Baudier2022}}]
Let $Y$, $Y'$ be uniformly discrete metric spaces of bounded geometry. Then, $Y$ and $Y'$ are coarsely equivalent if and only if $C^*_u(Y)$ and $C^*_u(Y')$ are Morita equivalent. In particular, up to isomorphism $K_*(C^*_u(Y))$ depends only on the coarse equivalence class of $Y$. 
\end{prop}
Many $C^*$-algebraic properties such as nuclearity, exactness, and so on are preserved under Morita equivalence. The uniform Roe algebra enjoying such a property tends to correspond to a geometric property of the underlying space. As one instance of this we mention the fact that a finitely generated group is exact if and only if its uniform Roe algebra is nuclear (see for example \cite{Brown2008}).

Suppose that $W$ is a subspace of $Y$. Above we introduced the algebra $C^*_u(W\subseteq Y)$ of operators near $W$. On the other hand $W$ is itself a uniformly discrete metric space of bounded geometry, and we can consider its uniform Roe algebra. The following proposition states that these are the same on the level of K-theory. This is a standard fact for non-uniform Roe algebras, and the proof is the same here.
\begin{prop}\label{IsoRoeLocalized}
Let $Y$ be a uniformly discrete metric space of bounded geometry, and $W\subseteq Y$ a subspace. The inclusion $C^*_u(W)\hookrightarrow C^*_u(Y)$ induces an isomorphism
$$
K_*\left( C^*_u(W) \right) \xlongrightarrow{\sim} K_*\left( C^*_u(W\subseteq Y) \right) \; .
$$
\end{prop}
\begin{proof}
An operator is supported near $W$ if it is supported inside some neighborhood $N_r(W)$, $r\geq 0$. Thus $C^*_u(W\subseteq Y)=\colim_{r\to \infty} C^*_u(N_r(W))$. On the other hand the inclusion $W\hookrightarrow N_{r}(W)$ is a coarse equivalence for all $r$, and thus induces an isomorphism in K-theory. Since K-theory commutes with direct limits, we conclude that
$$
K_*\left( C^*_u(W\subseteq Y)\right) = \colim_{r\to \infty} K_*\left( C^*_u(N_r(W))\right) = K_*\left( C^*_u(W)\right) \; .
$$
\end{proof}
We shall also need a different algebra. Consider the Hilbert space $\ell^2(Y\times \NN)\cong \ell^2(Y)\otimes\ell^2(\NN)$. A bounded operator on this Hilbert space may be viewed as a potentially infinite matrix $T=(t_{xy})_{x,y\in Y}$, where the entries are now bounded operators on $\ell^2(\NN)$. The properties of finite propagation, having uniformly bounded coefficients and being supported near a subspace may be defined analogously to operators on $\ell^2(Y)$. Also, we say that $(t_{xy})_{x,y\in Y}$ is \emph{locally compact} if each $t_{xy}$ is compact.
\begin{definition}[{\cite[Definition 6.6]{Spakula2009}}]
Let $Y$ be a uniformly discrete metric space of bounded geometry, and $W\subseteq Y$ a subset. Define $C^*_k(Y)$ and $C^*_k(W\subseteq Y)$ as the norm-closure in $\fB(\ell^2(Y\times\NN))$ of all locally compact finite-propagation operators (supported near $W$ for the latter) that satisfy the following uniformness condition: For all $\varepsilon>0$ there exists $N\in\NN$ such that for all $x,y\in Y$ the coefficient $t_{xy}$ is within $\varepsilon$ of a rank-$N$ operator.
\end{definition}
The algebra $C_k^*(Y)$ was introduced by \v{S}pakula in \cite{Spakula2009}. There it is shown that the inclusion $C^*_u(Y)\otimes \fK(\ell^2(\NN)) \hookrightarrow C^*_k(Y)$ is an isomorphism, so that there is no distinction between $C^*_u(Y)$ and $C^*_k(Y)$ on the level of K-theory. Clearly, this inclusion preserves the property of being supported near $W$, so we get the following:
\begin{prop}[{\cite[Lemma 6.10]{Spakula2009}}]
The inclusion $C^*_u(W\subseteq Y)\otimes \fK(\ell^2(\NN)) \hookrightarrow C^*_k(W\subseteq Y)$ is an isomorphism.
\end{prop}
We also wish to consider uniform Roe algebras on spaces that are not discrete. A first definition using quasi-lattices and corresponding partitions was given by \v{S}pakula in \cite{Spakula2009}. Engel later reformulated the definition in a more intrinsic way \cite{Engel2014}. We will present the latter of the two definitions. Also, we permit ourselves to be somewhat lax in our definitions at this point, because a more detailed treatment is contained as a special case in Section \ref{SecRelUniStructureAlg} below.

Let $X$ be a proper metric space of jointly bounded geometry. We fix an ample representation $(H,\sigma)$ of $X$.\footnote{Recall our convention that every Hilbert space is separable.} Define 
$$
C^*_u(X;\sigma) := \overline{ \left\{ T\in\fBH \, | \, T\sigma, \, \sigma T\sua 0, \, T \, \text{has finite propagation}\ \right\} } \subseteq \fBH \; .
$$
While $C^*_u(X;\sigma)$ depends on the choice of ample representation, its K-theory is independent of it (compare Proposition \ref{RelUniAmple} below). First of all we should convince ourselves that this version of the uniform Roe algebra coincides with the previously considered one -at least on the level of K-theory- if $Y$ is a uniformly discrete space of bounded geometry. To that end note that if $Y$ is such a space, then $(H=\ell^2(Y\times \NN),\sigma=\mu \otimes I)$ is an ample representation, where $\mu$ is the multiplication representation on $\ell^2(Y)$, so that $\sigma(f)$ acts as multiplication by $f(x)$ on each subspace $\ell^2(\{x\}\times\NN)$. Then it is not hard to see that the algebra $C^*_u(Y;\sigma)$ as defined here coincides with the algebra $C^*_k(Y)$ of Definition \ref{DefUniRoe}. Consequently we have $K_*(C^*_u(Y;\sigma))=K_*(C^*_u(Y))$.

The connection goes further, indeed the algebra $C^*_u(X;\sigma)$ is isomorphic to the $C^*_k$-algebra of any quasi-lattice. This is essentially shown in Section 8 of \cite{Spakula2009}, see also Section 2.7 in \cite{Engel2014}. For the sake of completeness we give a run-down of the argument nonetheless. Assume that $Y\subseteq X$ is a quasi-lattice in the proper metric space $X$ of jointly bounded geometry. We can find a Borel partition $X=\bigcup_{y\in Y} U_y$ such that each $U_y\cap Y=\{y\}$, and the $U_y$ have uniformly bounded diameter and non-empty interior. Let $P_y:=\sigma(\mathbf{1}_{U_y})$ and $H_y:=P_yH$. Each $H_y$ is a countably infinite-dimensional Hilbert space, which we can identify with $\ell^2(\NN)$. To a bounded operator $T\in\fBH$ we can thus associate the bounded operator $(t_{xy})_{x,y\in Y} \in \fB(\ell^2(Y\times\NN))$, $t_{xy}:=P_yTP_x$. We show that $T\mapsto (t_{xy})_{x,y\in Y}$ provides an isomorphism between $C^*_u(X;\sigma)$ and $C^*_k(Y)$. Clearly $T$ has finite propagation if and only if $(t_{xy})$ does, so it remains to argue for compatibility with uniform local compactness.
\begin{lemma}[Compare {\cite[Section 8]{Spakula2009}}, and {\cite[Section 2.7]{Engel2014}}] \label{LemmaLocalizingQuasilattice}
Assume that the bounded-geometry metric space $X$ is proper. Let $T$ have finite propagation. Then $T\sigma,\sigma T\sua 0$ if and only if $t_{xy}\sua 0$.\footnote{With $t_{xy}\sua 0$ we mean that each $t_{xy}$ is compact, and the uniformness property from the definition of $C^*_k(Y)$ is satisfied, i.e for all $\varepsilon>0$ there is $N$ such that each $t_{xy}$ is within $\varepsilon$ of a rank-$N$ operator.}
\end{lemma}
\begin{proof}
Assume that $T\sigma,\sigma T\sua 0$. Let $R$ be a uniform bound for the diamater of the $U_y$. For each $y\in Y$ we can find a function $\psi_y$ with $\psi_y=1$ on $U_y$, such that $\psi_y$ is $1$-Lipschitz and has support of diameter $\leq R+2$.\footnote{For example $\psi_y(x):= \min(1,d(x,N_1(U_j))$ has these properties.} Since $X$ is proper the functions $\psi_y$ are automatically compactly supported, and thus in $1\text{-}\mathrm{Lip}_{R+2}$. It follows that $t_{xy}=P_y(\sigma(\psi_y)T)P_x\sua 0$.

Conversely assume that $t_{xy}\sua 0$. For fixed $R\geq 0$ let $f$ have support of diameter $\leq R$. Because of bounded geometry there is $N\in\NN$ such that the support of $f$ intersects at most $N$ of the $U_y$. Moreover, since $T$ has finite-propagation, there is a finite number $M$ such that $\sigma(\mathbf{1}_{U_y})T\sigma(f)\neq 0$, and because of bounded geometry that number depends only on $R$. Combining these facts we obtain that $T\sigma(f)=\sum_{x,y} t_{x,y}\sigma(f)$ has at most $N\cdot M$ non-zero summands, where $N$ and $M$ depend only on $R$. Thus $T\sigma\sua 0$, and $\sigma T\sua 0$ is argued analogously. 
\end{proof}
We conclude that the mapping $T\mapsto (t_{xy})_{x,y\in Y}$ produces an isomorphism from $C^*_u(X;\sigma)$ to $C^*_k(Y)$ if $X$ is proper and has jointly bounded geometry\footnote{In fact, it would suffice to assume that $X$ has coarsely bounded geometry, i.e. admits a uniformly locally finite quasi-lattice.}. It is also clear that this mapping preserves the property of being supported near a given subset. Indeed, if $A\subseteq X$ is a (closed) subset, we can define $C^*_u(A\subseteq X;\sigma)$ to be the $C^*$-subalgebra of $C^*_u(X;\sigma)$ generated by all operators supported near $A$\footnote{The \emph{support} $\supp(T)$ of an operator $T\in\fBH$ is defined as the complement of the union of all open $U\times V \subseteq X\times X$ such that $\sigma(f)T\sigma(g)=0$ for all $f\in C_0(U)$ and $g\in C_0(V)$. Then, $T$ is said to be supported near $A$ if its support is contained in a neighborhood of $A\times A$.}. There is some subset $W\subseteq Y$ such that $A\subseteq N_c(W)$ and $W\subseteq N_c(A)$ for some $c\geq 0$. It is easily seen that an operator $T\in C^*_u(X;\sigma)$ is near $A$ if and only if $(t_{xy})_{x,y\in Y}$ is near $W$. We conclude the following:
\begin{prop}[Compare {\cite[Section 8]{Spakula2009}} and {\cite[Section 2.7]{Engel2014}}] \label{UniRoeQuasi}
Let $X$ be a proper metric space of jointly bounded geometry and $A\subseteq X$ a subset. Let $Y\subseteq X$ be a quasi-lattice, and $W\subseteq Y$ a subset such that $A\subseteq N_c(W)$ and $W\subseteq N_c(A)$ for some $c\geq 0$. Then, there is an isomorphism
$$
K_*\left( C^*_u(A\subseteq X;\sigma) \right) \xlongrightarrow{\sim} K_*\left( C^*_u(W\subseteq Y) \right) \; .
$$
\end{prop}
In light of Proposition \ref{UniRoeQuasi} we feel justified in referring to $C^*_u(X;\sigma)$ as the \emph{uniform Roe algebra} of $X$. We will also drop the representation $\sigma$ from the notation. We remark that Lemma \ref{LemmaLocalizingQuasilattice} and hence Proposition \ref{UniRoeQuasi} does not require the full strength of jointly bounded geometry, it suffices that $X$ admits a quasi-lattice (i.e. has bounded geometry in the coarse sense) and is proper. Properness of $X$ was only used for the inclusion $C^*_u(X;\sigma)\subseteq C^*_k(Y)$ in Lemma \ref{LemmaLocalizingQuasilattice}. Without assuming properness this inclusion could reasonably fail, due to a potential lack of compactly supported bump-functions. Indeed, functions in $\LLip_R(X)$ are forced to vanish when approaching regions of incompleteness, so that $Tf\sua 0$ requires less of $T$ than $t_{xy}\sua 0$ over an incomplete region. Thus $K_*(C^*_u(X;\sigma))$ and $K_*(C^*_u(Y))$ could reasonably be different if $X$ is not proper.

We come to the uniform coarse index map from uniform K-homology to the K-theory of the uniform Roe algebra. To that end we define the \emph{uniform structure algebra} as
$$
D^*_u(X) := \overline{ \left\{ T\in\fBH \, | \, [T,\sigma]\sua 0, \, T \, \text{has finite propagation}\ \right\} } \subseteq \fBH \; ,
$$
where $(H,\sigma)$ is a fixed ample representation of $X$. Again, the K-theory of $D^*_u(X)$ is independent of the choice of ample representation (compare Proposition \ref{RelUniAmple} below), so that we do not make it explicit in the notation. The uniform Roe algebra $C^*_u(X)$ is an ideal in the uniform structure algebra $D^*_u(X)$.
\begin{lemma}[{\cite[Lemma 7.2]{Engel2014}}]
Let $p\in\ZZ_{\geq -1}$. Let $X$ be a proper metric space of jointly bounded geometry. Then, there is an isomorphism
$$
K^{u}_{*}(X) \cong K_{*+1}\left( D^*_u(X)/C^*_u(X) \right) \; .
$$ 
\end{lemma}
\begin{proof}
By Paschke duality it suffices to prove that the inclusion $D^*_u(X)\to \fD^{u}(X)$ induces an isomorphism $D^*_u(X)/C^*_u(X)\cong\fD^{u}(X)/\fC^{u}(X)$. Here we use jointly bounded geometry to represent $\fD^{u}(X)$ and $\fC^{u}(X)$ via a fixed ample representation. Because $C^*_u(X)=D^*_u(X)\cap \fD^{u}(X)$ it further suffices to prove that $D^*_u(X)+\fC^{u}(X)=\fD^{u}(X)$. This can be done via truncation as in \eqref{Truncation}.
\end{proof}
With this in hand we may define the uniform coarse index map as follows.
\begin{definition}
Let $X$ be a proper metric space of jointly bounded geometry. The \emph{uniform coarse index map} is defined as the composition
$$
\Ind: \; K^{u}_*(X) \cong K_{*+1}\left( D^*_u(X)/C^*_u(X) \right) \xlongrightarrow{\partial} K_*\left( C^*_u(X)\right) \; .
$$
\end{definition}
This definition of the uniform coarse index map was given by Engel in \cite{Engel2014}. \v{S}pakula previously constructed a uniform coarse index map which goes directly into the K-theory of the uniform Roe algebra of a quasi-lattice. Of course the two versions are identified under the isomorphism of Proposition \ref{UniRoeQuasi}. 

Before we turn to uniform coarse index maps on relative uniform K-homology, we give a brief overview of the connection between uniform coarse index maps and the Baum-Connes conjecture with $\ell^\infty$-coefficients. The following theorem was proved by \v{S}pakula for torsionfree groups \cite{Spakula2009}, while the general case was provided by Engel \cite{Engel2019}.
\begin{thm}[\cite{Spakula2009}, \cite{Engel2019}]
Let $\Gamma$ be a countable discrete group. Then, there is an isomorphism
$$
\lim_{d\to \infty} K^{u}
_*(P_d\Gamma) \xlongrightarrow{\sim} K^{\Gamma}_*\left(\underline{E}\Gamma;\ell^\infty(\Gamma)\right)
$$
making the diagram \\
\centerline{\xymatrix{ 
\lim_{d\to \infty} K^{u}
_*(P_d\Gamma)\ar[r]^\sim \ar[d]_{\mu_u}  &  K^{\Gamma}_*\left(\underline{E}\Gamma;\ell^\infty(\Gamma)\right)\ar[d]^{\mu} \\
K_*\left(C^*_u(\Gamma)\right) \ar@{=}[r] & K_*\left( \ell^\infty(\Gamma)\rtimes_r \Gamma \right)
}}
commute. Here $P_d\Gamma$ are the Rips complexes of $\Gamma$, $\mu_u$ is the limit of the uniform coarse index maps, and $\mu$ is the Baum-Connes assembly map with coefficients in $\ell^\infty(\Gamma)$.
\end{thm}
Consequently the uniform coarse Baum-Connes conjecture for a group $\Gamma$, i.e. that the uniform assembly map $\mu_u: \, \lim_{d\to \infty} K^{u}
_*(P_d\Gamma) \to K_*(C^*_u(\Gamma))$ is an isomorphism, is equivalent to the Baum-Connes conjecture with coefficients in $\ell^\infty(\Gamma)$. This is a uniform analogy to a result of Yu stating that the coarse Baum-Connes conjecture for $\Gamma$ is equivalent to the Baum-Connes conjecture with coefficients in $\ell^\infty(\Gamma;\mathfrak{K})$ \cite{Yu1995}.

The Baum-Connes assembly map with arbitrary coefficients is known to an isomorphism for hyperbolic groups \cite{Kasparov2003}, \cite{Lafforgue2012}, or for groups with the Haagerup property \cite{Higson2001}. The former class of groups includes free groups, Fuchsian groups (in particular surface groups), and fundamental groups of closed manifolds with negative sectional curvature. The latter class in particular includes amenable groups. We also point out the somewhat recent survey article \cite{Aparicio2019}.

Suppose that $X$ is a manifold of bounded geometry, and $Y\subseteq X$ a quasi-lattice. If $X$ is uniformly contractible, meaning that for every $r>0$ there is $s>r$ such that the inclusion $B_r(x)\hookrightarrow B_s(x)$ is null-homotopic for every $x\in X$, then there is an isomorphism $\lim_{d\to \infty} K^{u}
_*(P_d Y)\cong K^{u}_*(M)$ \cite[Theorem 2.35]{Engel2019}. This isomorphism of course identifies the uniform coarse index map $\Ind$ with the uniform assembly map $\mu_u$. This is particularly relevant to the Baum-Connes conjecture with coefficients if $X$ carries a proper, cocompact and isometric group action by some discrete group $\Gamma$, because then $\Gamma$, identified with some orbit in $X$, is a quasi-lattice in $X$. Examples of this situation arise from universal covers of closed acyclic manifolds. Indeed, suppose that $M$ is a closed acyclic manifold, meaning that $\pi_i(M)=0$ for $i\neq 1$. Then, the universal cover $X$ of $M$ is uniformly contractible, and carries an proper, cocompact and isometric action of the fundamental group of $M$. Thus the Baum-Connes conjecture with $\ell^\infty$-coefficients translates directly to the uniform coarse Baum-Connes conjecture for $X$.
\begin{cor}[{\cite[Theorem 2.35]{Engel2019}}]
Let $X$ be the universal cover of a closed acyclic manifold $M$ with fundamental group that is hyperbolic or satisfies the Haagerup property. Then, the uniform coarse index map
$$
K_*^{u}(X) \xlongrightarrow{\Ind} K_*\left(C^*_u(X)\right)
$$
is an isomorphism.
\end{cor}
\begin{example}
Euclidean space $\RR^n$ satisfies the uniform coarse Baum-Connes conjecture because it carries a proper, cocompact and isometric action of the amenable group $\ZZ^n$, or equivalently because it is the universal cover of the acyclic manifold $\TT^n=\RR^n/\ZZ^n$ with fundamental group $\ZZ^n$. The hyperbolic plane $\HH^2$ satisfies the uniform coarse Baum-Connes conjecture because it carries a proper, cocompact and isometric action of the hyperbolic surface groups, or equivalently because it is the universal cover of the acyclic surfaces of genus $\geq 2$. More generally the universal cover of a closed manifold with negative sectional curvature satisfies the uniform coarse Baum-Connes conjecture.
\end{example}
\section{Index maps to relative uniform Roe algebras}
In this section we introduce a first construction of an index map on relative uniform K-homology, whose target is the K-theory of a certain relative uniform Roe algebra. This index map has favorable structural property, but is unfortunately trivial in relevant situations. Before we come to the definition of the relative uniform Roe algebra and the corresponding index map, we first need to transfer the result that the K-theory of the uniform Roe algebra of operators supported near a subset coincides with that of the uniform Roe algebra of that subset itself to the case of not necessarily discrete spaces. This requires an additional boundedness assumption on the subspace.
\begin{definition}
Let $X$ be a metric space of jointly bounded geometry. A subspace $Z\subseteq X$ is said to be \emph{boundedly embedded} if $Z$ also has jointly bounded geometry, and there exists a quasi-lattice $\Gamma$ in $X$ such that $\Gamma\cap Z$ is a quasi-lattice in $Z$. 
\end{definition}
\begin{example}
Any subcomplex of a simplicial complex of bounded geometry is boundedly embedded. Similarly, the boundary of a manifold with boundary and bounded geometry is boundedly embedded. More generally, if the subspace $Z\subseteq X$ admits a uniformly thick tubular neighborhood in a suitable sense, then it is boundedly embedded.
\end{example}
\begin{lemma}\label{UniRoeSubset}
Let $X$ be a proper space of jointly bounded geometry, and $Z\subseteq X$ a closed and boundedly embedded subspace. Then, 
$$
K_*\left( C^*_u(Z) \right) \cong K_*\left( C^*_u(Z\subseteq X) \right) \; .
$$
\end{lemma}
\begin{proof}
We prove this statement via reduction to the analogous statement for uniform Roe algebras (of discrete spaces). Choose an ample representation $(H,\sigma)$ of $X$ such that $(H_Z:=\sigma(\mathbf{1}_Z)H,\sigma_Z:=\sigma(\mathbf{1}_Z \cdot))$ is an ample representation of $Z$,\footnote{If $Z$ has empty interior, then this is not necessarily the case for an arbitrary ample representation of $X$. However, if $(H_X,\sigma_X)$ is an ample representation of $X$, and $(H_Z,\sigma_Z)$ an ample representation of $Z$, then $(H_X\oplus H_Z,\sigma_X\oplus \sigma_Z(\cdot|_Z))$ is again ample, and has the desired property.} and define $C^*_u(X)$ and $C^*_u(Z)$ using these representations. Then the decomposition $H=H_Z\oplus H_Z^\perp$ gives rise to an upper-corner inclusion $\fB(H_Z)\hookrightarrow \fBH$. This inclusion maps $C^*_u(Z)$ into $C^*_u(Z\subseteq X)$.

Let $\Gamma$ be a quasi-lattice in $X$ such that $\Gamma_Z:= \Gamma\cap Z$ is a quasi-lattice in $Z$. Fix Borel decompositions $X=\bigcup_{y\in \Gamma} U_y$ and $Z=\bigcup_{y\in\Gamma_Z} V_y$ such that each $U_y\cap \Gamma = \{y\}$, $V_y\cap\Gamma_Z=\{y\}$, and the $U_y$, $V_y$ have uniformly bounded diameter and non-empty interior. Let $P_y$ and $P_y^Z$ be the projections associated to the indicator function of $U_y$ and $V_y$ respectively. Then, $P_yH$ and $P_y^ZH_Z$ are both countably infinite-dimensional Hilbert spaces, which we identify with $\ell^2(\NN)$. Then, recall that $C^*_u(Z\subseteq X)$ and $C^*_u(Z)$ can be isomorphically mapped to $C^*_k(\Gamma_Z\subseteq\Gamma)$ and $C^*_k(\Gamma_Z)$ via $T\mapsto (t_{xy})$, where $t_{xy}$ is given by $P_yTP_x$ and $P^Z_yTP^Z_x$, respectively, which we view as a bounded operator on $\ell^2(\NN)$. Via this identification the inclusion $C^*_u(Z)\hookrightarrow C^*_u(Z\subseteq X)$ is turned into the inclusion $C_k^*(\Gamma_Z)\hookrightarrow C_k^*(\Gamma_Z\subseteq \Gamma)$. This latter inclusion is in turn obtained from the inclusion $C^*_u(\Gamma_Z)\hookrightarrow C^*_u(\Gamma_Z\subseteq \Gamma)$ by tensoring with the compact operators. This last inclusion induces an isomorphism on K-theory by Proposition \ref{IsoRoeLocalized}. The claim then follows from the fact that the K-theories of $C^*_u(\Gamma_Z)$ and $C^*_u(Z)$, and $C^*_u(\Gamma_Z\subseteq \Gamma)$ and $C_u^*(Z\subseteq X)$ are canonically isomorphic.
\end{proof}
Let $X$ be a proper metric space of jointly bounded geometry, and $Z\subseteq X$ a boundedly embedded subspace. We define the \emph{relative uniform Roe algebra} $C^*_u(X,Z)$ as the quotient of $C^*_u(X)$ by the ideal $C^*_u(Z\subseteq X)$. Associated to the short exact sequence
$$
0 \longrightarrow C^*_u(Z\subseteq X) \longrightarrow C^*_u(X) \longrightarrow C^*_u(X,Z) \longrightarrow 0
$$
there is a long exact sequence in K-theory. Combined with the isomorphism $K_*(C^*_u(Z))\cong K_*(C^*_u(Z\subseteq X))$ of Lemma \ref{UniRoeSubset} this yields the long exact sequence
$$
\cdots \longrightarrow K_*\left( C^*_u(Z) \right) \longrightarrow K_*\left( C^*_u(X) \right) \longrightarrow K_*\left( C^*_u(X,Z) \right) \xlongrightarrow{\partial} K_{*-1}\left( C^*_u(Z) \right) \longrightarrow \cdots \; .
$$
The uniform coarse index map should give rise to transformation of the long exact sequence in uniform K-homology to this sequence. To do this we first have to define an index map 
$$
\Ind: \; K^{u}_*(X,Z) \longrightarrow K_*\left( C^*_u(X,Z) \right) \; .
$$
Its intuitive definition goes like this: By Paschke duality we may represent $K^{u}_*(X,Z)$ by $K_{*+1}\left(\fD^{u}(X)/\fD^{u}(X,Z)\right)$, meaning that a relative uniform Fredholm module is invertible modulo $\fD^{u}(X,Z)$. Away from $Z$ operators in $\fD^{u}(X,Z)$ are uniformly locally compact. Thus, upon ignoring what happens near $Z$ a relative uniform Fredholm module is invertible modulo uniformly locally compact operators. The uniform coarse index in $K_*(C^*_u(X,Z))$ then measures the failure to be actually invertible away from $Z$. The formal construction goes as follows.
\begin{prop}\label{RelRoeIndex}
There is a natural uniform coarse index map
$$
\Ind: \; K^{u}_*(X,Z) \longrightarrow K_*\left(C^*_u(X,Z)\right) 
$$
making the diagram \\
\centerline{\xymatrix{
K^{u}_*(X) \ar[r] \ar[d]_{\Ind} & K^{u}_*(X,Z) \ar[d]^{\Ind} \\
K_*\left(C^*_u(X)\right) \ar[r] & K_*\left(C^*_u(X,Z)\right) 
}}
commute.
\end{prop}
\begin{proof}
Let $(H,\rho,F)$ be a relative uniform Fredholm module over $(X,Z)$. By Proposition \ref{RelTrunc} we may assume $F$ to have propagation $\leq r$ for some positive $r$. We may additionally take $F$ to be self-adjoint. Fix a Lipschitz function $\psi: \, X\to [0,1]$ such that $\psi(x)=0$ if $d(x,Z)\leq 2r$ and $\psi(x)=1$ if $d(x,Z)\geq 3r$. Then $F-\psi F$ is supported near $Z$. Thus, modulo such operators it holds that $I-F^2 = \psi - (\psi F)^2$. Observe that if $f\in C_0(X)$ has support inside the $r$-neighborhood of $Z$, then
$$
\rho(f) \left( \psi - (\psi F)^2 \right) = 0 = \left( \psi - (\psi F)^2 \right) \rho(f) \; .
$$
On the other hand $(I-F^2)\rho \sua 0$ over $X\setminus Z$, which implies that $(\psi - (\psi F)^2)\rho\sua 0$ over $X\setminus Z$ as well. Now, any function on $X$ can be written as the sum of a function supported in $N_r(Z)$ and a function vanishing on $Z$. Thus $(\psi - (\psi F)^2)\rho\sua 0$ actually holds over $X$. We conclude that modulo operators supported near $Z$ it holds that $I-F^2 \in C^*_u(X)$. In other words the image of $F$ in $D^*_u(X)/D^*_u(Z\subseteq X)$ is an involution modulo the image of $C^*_u(X)$ in $D^*_u(X)/D^*_u(Z\subseteq X)$. This image is precisely $C^*_u(X,Z)$. 

On more formal terms the quotient $D^*_u(X)/D^*_u(Z\subseteq X)$ contains an isomorphic image of $C^*_u(X,Z)$ as an ideal, and with the help of Paschke duality and the above argument we produce a map
$$
K^{u}_*(X,Z) \longrightarrow K_{*+1}\left(\left( D^*_u(X)/D^*_u(Z\subseteq X) \right)/C^*_u(X,Z)\right) \; .
$$
Composing this map with the K-theoretic boundary map produces the uniform coarse index map
$$
\Ind: \; K^{u}_*(X,Z) \longrightarrow K_{*+1}\left(\left( D^*_u(X)/D^*_u(Z\subseteq X) \right)/C^*_u(X,Z)\right) \xlongrightarrow{\partial} K_*\left(C^*_u(X,Z)\right) \; .
$$
Then, the claimed commutativity statement is a consequence of the naturality of the K-theoretical boundary map with respect to the quotient map.
\end{proof}
The uniform coarse index map $K^{u}_*(X,Z)\to K_*(C^*_u(X,Z))$ completes the mapping from the long exact sequence in uniform K-homology to that in K-theory of the uniform Roe algebras, and we claim that the relevant diagram commutes. 
\begin{thm}\label{CommUniIndex}
Let $X$ be a proper metric space of jointly bounded geometry, and $Z\subseteq X$ a transversally controlled and boundedly embedded subspace. Then, the diagram \\
\centerline{\xymatrix{
K_*^{u}(Z) \ar[r] \ar[d]_{\Ind} & K^{u}_*(X) \ar[r] \ar[d]_{\Ind} & K^{u}_*(X,Z) \ar[r]^\partial \ar[d]_{\Ind} & K_{*-1}^{u}(Z) \ar[d]_{\Ind} \\
K_*\left(C^*_u(Z)\right) \ar[r] & K_*\left(C^*_u(X)\right) \ar[r] & K_*\left(C^*_u(X,Z)\right) \ar[r]_\partial & K_{*-1}\left(C^*_u(Z)\right)
}}
commutes.
\end{thm}
We defer the proof to Section \ref{SectionRelIndexMap}, where it will be deduced as a corollary to a similar commutativity statement. The uniform coarse index map $K^{u}_*(X,Z)\to K_*(C^*_u(X,Z))$ and the resulting Theorem \ref{CommUniIndex} are satisfying from a structural perspective. Indeed, viewing both uniform K-homology and K-theory of the uniform Roe algebra as homology theories on suitable categories of metric spaces, the uniform coarse index map provides a transformation from one homology theory to the other. Theorem \ref{CommUniIndex} then states that this transformation is compatible with long exact sequences. Unfortunately, it is trivial as soon as $Z$ is $c$-dense in $X$ for some $c$.
\begin{lemma}
Let $X$ be a proper metric space of jointly bounded geometry, and let $X$ be a boundedly embedded subspace. If $X=N_c(Z)$ for some $c\geq 0$, then
$$
K_*\left(C^*_u(X,Z)\right) = 0 \; .
$$
\end{lemma}
\begin{proof}
In this case it holds that $C^*_u(Z\subseteq X)=C^*_u(X)$, so that the map
$$
K_*\left(C^*_u(Z)\right)\xlongrightarrow{\sim} K_*\left( C^*_u(Z\subseteq X)\right) = K_*\left( C^*_u(X) \right)
$$
is an isomorphism. It follows from the long exact sequence that $K_*(C^*_u(X,Z))=0$.
\end{proof}
This result of course makes sense from a coarse perspective; if $Z$ is $c$-dense in $X$, coarsely they are the same space. In our applications of uniform coarse index maps to bordisms we will see that one should specifically consider bordisms $\olO$ of finite-width, i.e. those for which $d(x,\partial \Omega)$ is uniformly bounded for $x\in\olO$. Since in this situation $\partial\Omega$ is $c$-dense for suitable $c$, we cannot expect the uniform coarse index map $K^{u}_*(\olO,\dO)\to K_*(C^*_u(\olO,\dO))$ to contain any information whatsoever. The index map that will be constructed and studied in the remainder of this chapter does not have this problem, it contains valuable information even in the finite-width case.
\section{Relative uniform structure algebras}\label{SecRelUniStructureAlg}
In this section we introduce the (relative) uniform structure algebra and discuss its functoriality.
\begin{definition}
 Let $(X,Z)$ be a pair, and let $(H,\rho)$ be a representation of $X$.\footnote{Recall our convention that all metric spaces are locally compact and second-countable, and all Hilbert spaces separable.} Denote by $D^*_u(X;\rho)$ the $C^*$-subalgebra of $\fD^{u}_\rho(X)$ generated by all uniformly pseudo-local operators of finite propagation. Moreover, let $D^*_u(X,Z;\rho)$ denote the $C^*$-subalgebra of $\fD^{u}_\rho(X,Z)$ generated by all uniformly pseudo-local operators of finite propagation that are uniformly locally compact over $X\setminus Z$. Just as with the non-finite propagation version $D^*_u(X,Z;\rho)$ is an ideal in $D^*_u(X;\rho)$.
\end{definition}
Note that if $(H,\sigma)$ is ample, then $D^*_u(X,\emptyset;\sigma)$ is the uniform Roe algebra $C^*_u(X;\sigma)$ discussed in Section 3.1. This section thus provides a more systematic approach to these algebras as well.

First, we will discuss how to compare the algebras coming from different representations, and then functoriality. The discussion is analogous to that for uniform dual algebras in Section \ref{SecDualAlgebra}, with suitable modifications concerning finite propagation. These are in turn analogous to those in the non-uniform setting, see for example Chapters 6 and 12 of \cite{HigsonRoe2000}. The admissible morphisms are given by uniformly coarse maps of pairs. We recall the definition of a uniformly coarse map.
\begin{definition}
Let $X$ and $Y$ be metric spaces. A map $\Phi: \, X \to Y$ will be called \emph{uniformly coarse} if for all $R\geq 0$ there exists $S\geq 0$ such that
$$
d(x,y) \leq R \; \Rightarrow \;  d(\Phi(x),\Phi(y))\leq S \quad \mathrm{and} \quad d(\Phi(x),\Phi(y))\leq R \; \Rightarrow \;  d(x,y) \leq S \; .
$$ 
\end{definition}
It is not hard to see that if $\Phi: \, X\to Y$ is a uniformly coarse map that is continuous, then $\Phi$ is uniformly continuous and uniformly cobounded. Hence, if $\Phi$ is also proper, and $X$ has jointly bounded geometry, then $\Phi$ is approximately filtered in the sense of Remark \ref{ApproxFilter} by Lemma \ref{UCapprox}. In this case, $\Phi$ therefore preserves the relation $\sua$. Thus, if $\Phi$ additionally satisfies $\Phi(Z)\subseteq W$, then $\Phi^*$ is an approximately filtered map of pairs akin to Definition \ref{DefFilteredMapsPairs}. We will refer to such $\Phi$ as \emph{uniformly coarse maps of pairs}.

Next we introduce the class of covering isometries we will use to take the colimit over all representations of $X$.
\begin{definition}
Let $(X,Z)$ and $(Y,W)$ be pairs, and $\Phi: \, (X,Z)\to (Y,W)$ a uniformly coarse map of pairs. Let $(H^\rho,\rho)$ and $(H^\sigma,\sigma)$ be representations of $X$ and $Y$ respectively. An isometry $V: \, (H^\rho,\rho) \to (H^\sigma,\sigma)$ is said to \emph{uniformly boundedly cover $\Phi$} if $V$ uniformly covers the induced map $\Phi^*: \, C_0(Y)\to C_0(X)$ in the sense of Definition \ref{DefUCI}, and it holds that $\sup_{(x,y)\in \supp(V)} d(\Phi(x),y)< \infty$.
\end{definition}
Conjugation with uniformly bounded covering isometries maps finite-propagation operators to finite-propagation operators, see for example \cite[Lemma 6.3.11]{HigsonRoe2000}.
\begin{lemma}\label{AdUBCI}
Let $(X,Z)$ and $(Y,W)$ be pairs with $X$ of jointly bounded geometry, and $\Phi: \, (X,Z)\to (Y,W)$ a uniformly coarse map of pairs. Let $(H^\rho,\rho)$ and $(H^\sigma,\sigma)$ be representations of $X$ and $Y$, and $V: \; (H^\rho,\rho) \to (H^\sigma,\sigma)$ an isometry uniformly boundedly covering $\Phi$. Then, conjugation by $V$ gives rise to well-defined maps $\Ad(V): \, D^*_u(X;\rho)\to D^*_u(Y;\sigma)$ , and $\Ad(V): \, D_u^*(X,Z;\rho)\to D_u^*(Y,W;\sigma)$. Moreover, any two uniformly bounded covering isometries $ (H^\rho,\rho) \to (H^\sigma,\sigma)$ induce the same map on K-theory.
\end{lemma}
Let $X$ be a metric space of jointly bounded geometry. Let $\fP_{ub}(X)$ denote the sub-category of $\fP(X)$ containing the same objects, i.e. representations of $X$, but as morphisms only uniformly bounded covering isometries. By virtue of Lemma \ref{AdUBCI} there are functors from $\fP_{ub}(X)$ to the category of $C^*$-algebras sending a representation $(H^\rho,\rho)$ to either $D^*_u(X;\rho)$ or $D^*_u(X,Z;\rho)$. Taking the colimit over these functors yields $C^*$-algebras
$$
D^*_u(X) := \colim\, D^*_u(X;\rho) \quad , \quad D^*_u(X,Z):= \colim\, D^*_u(X,Z;\rho) \; .
$$
Ultimately, we are mainly interested in the K-theory of these algebras, and since K-theory commutes with directed colimits, one can equally well pass to K-theory first and take the colimit over the representation after. 

Again, since $X$ has jointly bounded geometry, we can represent the K-theories of $D^*_u(X)$ and $D^*_u(X,Z)$ using a single ample representation. Indeed, as was proved by \v{S}pakula in \cite{Spakula2010} and discussed earlier, there is an isometry uniformly covering $\id_X$ from any representation to any given ample representation. Examination of the proof in \cite{Spakula2010} reveals that this covering isometry can in fact be made uniformly bounded. Appealing to the uniqueness part of Lemma \ref{AdUBCI} above we see that $K_*(D^*_u(X;\sigma))$ for $(H,\sigma)$ ample is terminal among $K_*(D^*_u(X;\rho))$, and hence represents the colimit. These remarks equally apply to $K_*(D^*_u(X,Z))$, so that we may conclude the following:
\begin{prop}\label{RelUniAmple}
Let $X$ be a metric space of jointly bounded geometry, and let $(H,\sigma)$ be an ample representation of $X$. Then
\begin{align*}
K_*\left( D^*_u(X) \right) &\cong K_*\left( D^*_u(X;\sigma) \right) \\
K_*\left( D^*_u(X,Z) \right) &\cong K_*\left( D^*_u(X,Z;\sigma) \right) \; .
\end{align*}
In particular it holds that $K_*(D^*_u(X,\emptyset))\cong K_*(C^*_u(X))$.
\end{prop}
We turn to functoriality. Let $\Phi: \, (X,Z)\to (Y,W)$ be a uniformly coarse map of pairs with $X$ having jointly bounded geometry, and $(H^\rho,\rho)$ a representation of $X$. Then $(H^\rho,\rho\circ\Phi^*)$ is a representation of $Y$, and $\id_{H^\rho}$ is an isometry $(H^\rho,\rho)\to (H^\rho,\rho\circ\Phi^*)$ that uniformly boundedly covers $\Phi$. Moreover, if $V: \, (H^\rho,\rho)\to (H^\sigma,\sigma)$ is an isometry uniformly boundedly covering $\id_X$, then $V: \, (H^\rho,\rho\circ\Phi^*)\to (H^\sigma,\sigma\circ\Phi^*)$ uniformly boundedly covers $\id_Y$, since $\Phi$ is a uniformly coarse map of pairs.

In this way the map $\Phi: \, (X,Z)\to (Y,W)$ gives rise to a functor $\Phi_*: \, \fP_{ub}(X) \to \fP_{ub}(Y)$, which induces a homomorphism from $D^*_u(X)$ to $D^*_u(Y)$, and from $D^*_u(X,Z)$ to $D^*_u(Y,W)$. We conclude that the K-theory groups
$$
K_*\left( D^*_u(X) \right) \quad , \quad K_*\left( D^*_u(X,Z) \right) 
$$
assemble into a functors from the category of metric spaces\footnote{Again, we really mean the category of locally compact second-countable metric spaces} of jointly bounded geometry, and continuous, proper, uniformly coarse maps into abelian groups in the first case, and from metric pairs with total space having jointly bounded geometry, and uniformly coarse maps of pairs into abelian groups in the second case. In the special case $Z=\emptyset$ we conclude that $K_*(C^*_u(X))$ assembles into a functor from the category of metric spaces of jointly bounded geometry and continuous, proper, uniformly coarse maps into abelian groups as well.
\section{Mayer-Vietoris sequences}\label{SecMVdual}
In this section we prove a geometric Mayer-Vietoris sequence for the relative uniform structure algebra associated to sufficiently nice decompositions into closed subspaces, as well as an algebraic Mayer-Vietoris sequence expressing the K-theory of the relative uniform structure algebra of $(X,Z)$ in terms of the uniform Roe algebras of $X$ and $Z$, and the uniform structure algebra of $Z$.

Let $(X,Z)$ be a pair, and $A\subseteq X$ a subspace. Define $D^*_u(A\subseteq X,Z)$ as the closure of the subset of $D^*_u(X,Z)$ containing those operators $T$ that are \emph{supported near} $A$, meaning that $\supp(T)\subseteq N_c(A\times A)$ for some $c\geq 0$. It is straight-forward to see that $D^*_u(A\subseteq X,Z)$ is an ideal in $D^*_u(X,Z)$.
\begin{lemma}\label{TruncDualClose}
Let $X$ be a proper metric space of jointly bounded geometry, and let $Z\subseteq A\subseteq  X$ be closed subspaces. The inclusion $D^*_u(A\subseteq X,Z) \hookrightarrow \fD^{u}(X,Z)$ induces an isomorphism
$$
D^*_u(A\subseteq X,Z)/C^*_u(A\subseteq X) \xlongrightarrow{\sim} \fD^{u}(X,Z)/\fC^{u}(X) \; .
$$
\end{lemma}
\begin{proof}
We may use a fixed ample representation to represent these algebras. The map $D^*_u(A\subseteq X,Z) \to \fD^{u}(X,Z)/\fC^{u}(X)$ induced by the inclusion is surjective. To see this first represent an element of $\fD^{u}(X,Z)/\fC^{u}(X)$ by a finite propagation operator via truncation, then cut off away from $A$ to get obtain a finite propagation operator supported near $A$. Moreover, it holds that $C^*_u(A\subseteq X)=D^*_u(A\subseteq X,Z)\cap \fC^{u}(X)$, so that the induced map between the quotients is injective as well.
\end{proof}
\begin{prop}\label{IsoStructureLocalized}
Let $X$ be a proper space of jointly bounded geometry, and let $Z\subseteq A\subseteq X$ be subspaces such that $A\subseteq X$ is boundedly embedded, and both $(X,Z)$ and $(A,Z)$ are transversally controlled. Then, the inclusion $A\hookrightarrow X$ induces an isomorphism
$$
K_*\left( D^*_u(A,Z) \right) \xlongrightarrow{\sim} K_*\left( D^*_u(A\subseteq X,Z) \right) \; .
$$
\end{prop}
\begin{proof}
We fix an ample representation $(H,\sigma)$ of $X$ such that $(H_A:=\sigma(\mathbf{1}_A)H,\sigma_A:=\sigma(\mathbf{1}_A\cdot))$ is an ample representation of $A$, and $(H_Z:=\sigma(\mathbf{1}_Z)H,\sigma_Z:=\sigma(\mathbf{1}_Z\cdot))$ an ample representation of $Z$. Then, the inclusion upper-corner inclusion $\fB(H_A)\hookrightarrow \fBH$ gives rise to the diagram \\
\centerline{\xymatrix{
0 \ar[r] & C^*_u(A;\sigma_A) \ar[r] \ar[d] & D^*_u(A,Z;\sigma_A) \ar[r] \ar[d] & D^*_u(A,Z;\sigma_A)/C^*_u(A;\sigma_A) \ar[r] \ar[d] & 0 \\
0 \ar[r] & C^*_u(A\subseteq X;\sigma) \ar[r] & D^*_u(A\subseteq X,Z;\sigma) \ar[r] & D^*_u(A\subseteq X,Z;\sigma)/C^*_u(Y\subseteq X;\sigma) \ar[r] & 0
}}
The first vertical arrow is a K-theory isomorphism because $A\subseteq X$ is boundedly embedded. The third vertical arrow is a K-theory isomorphism as well. To see this observe that it fits into the diagram\\
\centerline{\xymatrix{
D^*_u(A,Z;\sigma_A)/C^*_u(A;\sigma_A) \ar[r] \ar[d] & \fD^{u}_{\sigma_A}(A,Z)/\fC^{u}_{\sigma_A}(A) \ar[d] & \fD^{u}_{\sigma_Z}(Z)/\fC^{u}_{\sigma_Z}(Z) \ar[l] \ar[ld]  \\
D^*_u(A\subseteq X,Z;\sigma)/C^*_u(Y\subseteq X;\sigma) \ar[r] & \fD^{u}_{\sigma}(X,Z)/\fC^{u}_{\sigma}(X)  &
}}
The arrows making up the triangle are isomorphisms on K-theory by Propositions \ref{RelEmptyset} and \ref{InclusionIsoDual}, applicable due to our transversal boundedness assumption. The horizontal arrows in the square are isomorphisms even before passing to K-theory, with the lower one being covered by Lemma \ref{TruncDualClose}, and the upper one being due to a similar truncation argument. As the only remaining arrow the left vertical arrow is a K-theory isomorphism as well. Now, a combination of naturality of the 6-term exact sequence in operator K-theory and the 5-lemma lets us deduce that the map $D^*_u(A,Z)\to D^*_u(A\subseteq X,Z)$ induces an isomorphism on K-theory.
\end{proof}
We come to the geometric Mayer-Vietoris sequence. To that end recall also that a decomposition $X=A_1\cup A_2$ is called \emph{coarsely excisive} if for all $r>0$ there exists $s>0$ such that $N_r(A_1)\cap N_r(A_2)$ is contained in $N_s(A_1\cap A_2)$.
\begin{thm}\label{MVURelAlg}
Let $(X,Z)=(A_1\cup A_2,B_1\cup B_2)$ be a bounded Mayer-Vietoris decomposition, and assume that $X$ is proper and has jointly bounded geometry, and that $A_1,A_2\subseteq X$ are boundedly embedded. Assume further that $X=A_1\cup A_2$ is coarsely excisive. Then, there is a natural 6-term exact sequence\\
\centerline{\xymatrix{
K_0\left( D_u^*(A_1\cap A_2,B_1\cap B_2) \right) \ar[r] & K_0\left( D_u^*(A_1,B_1) \right) \oplus K_0\left( D_u^*(A_2,B_2) \right) \ar[r] & K_0\left( D_u^*(X,Z) \right) \ar[d]^{\delta} \\
K_1\left( D_u^*(X,Z) \right) \ar[u]^{\delta} & K_1\left( D_u^*(A_1,B_1) \right) \oplus K_1\left( D_u^*(A_2,B_2) \right) \ar[l] & K_1\left( D_u^*(A_1\cap A_2,B_1\cap B_2) \right) \ar[l]
}}
\end{thm}
\begin{proof}
Fix Lipschitz functions $\psi_j: \, X\to [0,1]$ with $\supp(\psi_j)\subseteq A_j$, $\supp(\psi_j)\cap Z\subseteq B_j$, and $\psi_1+\psi_2=1$.\footnote{For example, let $\psi_j'(x):= d(x,X\setminus (A_j \cup (Z\setminus \mathring{B_j}))$ (here $\mathring{B_j}$ denotes the interior as a subspace of $Z$) and define $\psi_j:=\psi_j'(\psi_1+\psi_2)^{-1}$. These functions are Lipschitz because $d(A_1\setminus A_2,A_2\setminus A_1)>0$ and $d(B_1\setminus B_2,B_2\setminus B_1)>0$, and the requirements on the support hold by construction.}. Let $T\in D^*_u(X,Z)$ have finite propagation, and write $T=T\psi_1+T\psi_2$. The operator $T\psi_j$ is supported near $\supp(\psi_j)\subseteq A_j$, and it is uniformly locally compact on the complement of $Z\cap \supp(\psi_j)\subseteq B_j$, so that $T\psi_j \in D^*_u(A_j\subseteq X,B_j)$. Passing to the closure conclude that $D^*_u(X,Z)=D^*_u(A_1\subseteq X,B_1)+D^*_u(A_2\subseteq X,B_2)$. Moreover, $D^*_u(A_1\subseteq X,B_1)\cap D^*_u(A_2\subseteq X,B_2) =D^*_u(A_1\cap A_2\subseteq X,B_1\cap B_2)$ because $X=A_1\cup A_2$ is coarsely excisive. The Mayer-Vietoris sequence of this decomposition combined with Proposition \ref{IsoStructureLocalized} yields the desired Mayer-Vietoris sequence.
\end{proof}
Restricting to the non-relative case $Z=\emptyset$ yields the Mayer-Vietoris sequence for the uniform Roe algebras. 
\begin{cor}
Let $X$ be a proper metric space of jointly bounded geometry, and let $X=A_1\cup A_2$ be a coarsely excisive decomposition. Then, there is a natural 6-term exact sequence\\
\centerline{\xymatrix{
K_0\left( C_u^*(A_1\cap A_2) \right) \ar[r] & K_0\left( C_u^*(A_1) \right) \oplus K_0\left( C_u^*(A_2) \right) \ar[r] & K_0\left( C_u^*(X) \right) \ar[d]^{\delta} \\
K_1\left( C_u^*(X) \right) \ar[u]^{\delta} & K_1\left( C_u^*(A_1) \right) \oplus K_1\left( C_u^*(A_2) \right) \ar[l] & K_1\left( C_u^*(A_1\cap A_2) \right) \ar[l]
}}
\end{cor}
Applying cut-offs away from $Z$ in an analogous way to the proof of Theorem \ref{MVURelAlg} one also derives an algebraic Mayer-Vietoris sequence for $D^*_u(X,Z)$. It makes precise the intuition that $D^*_u(X,Z)$ should depend on the metric geometry of $Z$ but only the uniform coarse geometry of $X$ away from $Z$.
\begin{prop}\label{AlgebraicMV}
Let $X$ be a proper metric space of jointly bounded geometry, and let $Z$ be a boundedly embedded subspace. Then, there is an algebraic Mayer-Vietoris sequence\\
\centerline{\xymatrix{
K_0\left(C^*_u(Z)\right) \ar[r] & K_0\left(C^*_u(X)\right)\oplus K_0\left(D^*_u(Z)\right) \ar[r] & K_0\left(D^*_u(X,Z)\right) \ar[d]^\delta \\
K_1\left(D^*_u(X,Z)\right) \ar[u]^\delta & K_1\left(C^*_u(X)\right)\oplus K_1\left(D^*_u(Z)\right) \ar[l] &  K_1\left(C^*_u(Z)\right) \ar[l]
}}
\end{prop}
\begin{proof}
We implicitly assume all algebras to be defined with respect to a fixed ample representation $(H,\sigma)$ of $X$ such that $(\sigma(\mathbf{1}_Z)H,\sigma|_Z)$ is an ample representation of $Z$. The algebra $D^*_u(X,Z)$ can be written as $D^*_u(X,Z)=C^*_u(X)+D^*_u(Z\subseteq X,Z)$. To see this let  $T\in D^*_u(X,Z)$ have propagation $R$, and pick a Lipschitz function $\psi: \, X\to [0,1]$ such that $\psi=1$ on $N_R(Z)$ and $\psi=0$ outside $N_{R+1}(Z)$. Write $T=\psi T+ (1-\psi)T$. Since $\psi$ is supported near $Z$ and $T$ has finite propagation the operator $\psi T$ is supported near $Z$, hence in $D^*(Z\subseteq X,Z)$. Analogously $(1-\psi)T$ is supported outside $Z$, where $T$ was locally compact to begin with, so that $(1-\psi)T\in C^*_u(X)$. Passing to the closure shows $D^*_u(X,Z)=C^*_u(X)+D^*_u(Z\subseteq X,Z)$. Both summands are ideals, and their intersection is given by $C^*_u(Z\subseteq X)$. Inserting the isomorphisms of Propositions \ref{IsoRoeLocalized} and \ref{IsoStructureLocalized} into the corresponding Mayer-Vietoris sequence concludes the proof.
\end{proof}
\begin{remark}\label{FinWidthMV}
Let us again remark on the case where $Z$ is $c$-dense, i.e. $N_c(Z)=X$. In this case the map $K_*(C^*_u(Z))\to K_*(C^*_u(X))$ is an isomorphism by Lemma \ref{UniRoeSubset}, as is the map $K_*(D^*_u(Z))\to K_*(D^*_u(X,Z))=K_*(D^*_u(Z\subseteq X,Z))$ by Proposition \ref{IsoStructureLocalized}. The latter isomorphism implies by exactness that the algebraic Mayer-Vietoris map $K_*(D^*_u(X,Z))\to K_{*-1}(C^*_u(Z))$ is zero in this case. This encodes a form of bordism invariance, see Section \ref{SectionBordismInv}.
\end{remark}
\section{The relative uniform index map}\label{SectionRelIndexMap}
We come to the definition of the relative uniform index map to the K-theory of the relative uniform structure algebra. Its definition is analogous to that of the uniform coarse index map on non-relative uniform K-homology. We just have to quickly verify that we also have a version of relative uniform Paschke duality where we replace $\fD^{u}(X)$ and $\fD^{u}(X,Z)$ by their finite-propagation counterparts.
\begin{lemma}\label{FinPropNormKHom}
Let $X$ be a proper metric space of jointly bounded geometry, and $Z\subseteq X$ a closed subspace. Then, the inclusion $D^*_u(X)\hookrightarrow \fD^{u}_\rho(X)$ induces an isomorphism
$$
D^*_u(X)/D^*_u(X,Z)\xlongrightarrow{\sim} \fD^{u}_\rho(X)/\fD^{u}_\rho(X,Z) \; .
$$
Thus, there is a natural isomorphism
$$
K^{u}_*(X,Z) \cong K_{*+1}\left( D^*_u(X)/D^*_u(X,Z) \right) \; .
$$
\end{lemma}
\begin{proof}
Since $D^*_u(X,Z)=\fD^{u}(X,Z)\cap D^*_u(X)$ it suffices to check that $D^*_u(X)+\fD^{u}(X,Z)=\fD^{u}(X)$. This can again be done via truncation, compare Lemma \ref{RelTrunc}.
\end{proof}
With this in hand we define the relative uniform index map in the usual way.
\begin{definition}
The \emph{relative uniform index map} $\relInd: \, K_*^{u}(X,Z) \to K_*(D^*_u(X,Z))$ of the pair $(X,Z)$ is defined as the composition
$$
K_*^{u}(X,Z) \xlongrightarrow{\sim} K_{*+1}\left(  D^*_u(X)/D^*_u(X,Z)\right) \xlongrightarrow{\partial} K_*\left(D^*_u(X,Z)\right) \; .
$$
\end{definition}
The relative uniform K-homology, (relative) uniform dual algebras and (relative) uniform structure algebras are functorial under uniformly coarse maps of pairs. Moreover, both Paschke duality and the boundary map of the pair $(D^*_u(X),D^*_u(X,Z))$ are natural with respect to such maps as well.
\begin{prop}\label{RelIndexNatural}
Let $X$ and $Y$ be proper metric spaces of jointly bounded geometry, and let $\Phi: \, (X,Z) \to (Y,W)$ be a uniformly coarse map of pairs. Then, the diagram\\
\centerline{\xymatrix{
K^{u}_*(X,Z) \ar[d]_{\relInd} \ar[r]^{\Phi_*} &  \ar[d]^{\relInd} K^{u}_*(Y,W) \\
K_*(D^*_u(X,Z)) \ar[r]_{\Phi_*} & K_*(D^*_u(Y,W))
}}
commutes. 
\end{prop}
We come to the main result of this chapter, which provides compatibility of the relative uniform index maps with Mayer-Vietoris boundary maps. It will be the major ingredient to our discussion of partitioned manifold indices in Chapter 6. As such the theorem can be seen as an abstract version of the relative uniform partitioned-manifold index theorem.
\begin{thm}\label{AbstractRelPM}
Let $X$ be a proper metric space of jointly bounded geometry, and $Z\subseteq X$ a closed subspace such that $X$ is transversally controlled near $Z$. Let $(X,Z)=(Y_1\cup Y_2,W_1\cup W_2)$ be a coarsely excisive and bounded Mayer-Vietoris decomposition, where $W_j=Z\cap Y_j$.  Then, the relative uniform index map takes the Mayer-Vietoris sequence in the relative uniform K-homology of the decomposition $(X,Z)=(Y_1\cup Y_2,W_1\cup W_2)$ to the Mayer-Vietoris in the K-theory of the relative uniform structure algebras of this decomposition. In particular the diagram\\
\centerline{\xymatrix{
K^{u}_{*+1}(X,Z) \ar[r]^-{\delta} \ar[d]_{\relInd} &  \ar[d]^{\relInd} K^{u}_*(Y_1\cap Y_2,W_1\cap W_2) \\
K_{*+1}(D^*_u(X,Z)) \ar[r]_-{\delta} & K_*(D^*_u(Y_1\cap Y_2,W_1\cap W_2))
}}
commutes.
\end{thm}
\begin{proof}
To derive the Mayer-Vietoris sequence for relative uniform K-homology we have used relative uniform dual algebras, but in the setting considered here one may equally well use relative uniform structure algebras. Indeed, using a partition of unity and coarse excisiveness we find that
\begin{align*}
D^*_u(X) &= D^*_u(Y_1\subseteq X,Y_1) + D^*_u(Y_2\subseteq X,Y_2) \\
D^*_u(Y_1\cap Y_2\subseteq X, Y_1\cap Y_2) &=  D^*_u(Y_1\subseteq X,Y_1) \cap D^*_u(Y_2\subseteq X,Y_2) \\
D^*_u(Y_j\subseteq X, W_j) &= D^*_u(Y_j\subseteq X,Y_j) \cap D^*_u(X,Z) \\
D^*_u(Y_1\cap Y_2\subseteq X, W_1\cap W_2) &=  D^*_u(Y_1\subseteq X,Y_1) \cap D^*_u(Y_2\subseteq X,Y_2)  \cap D^*_u(X,Z)\; .
\end{align*}
Using the facts from Section \ref{SecMVdual} and Proposition \ref{MVQuot} we find that the decomposition 
$$
D^*_u(X)/D^*_u(X,Z) = \left( D^*_u(Y_1\subseteq X,Y_1)/D^*_u(Y_1\subseteq X, W_1) \right) + \left( D^*_u(Y_1\subseteq X,Y_2)/D^*_u(Y_2\subseteq X, W_2) \right)
$$
yields the Mayer-Vietoris sequence in relative uniform K-homology. Item (ii) of Proposition \ref{MVQuot} then provides the compatibility of the relative uniform index maps with the Mayer-Vietoris boundary maps. The remaining commutativity statements required to show that the relative uniform index maps takes one Mayer-Vietoris sequence to the other are already covered by Proposition \ref{RelIndexNatural}. Thus the claim is proved.
\end{proof}
Restricting to the case $Z=\emptyset$ also gives us the compatibility of the non-relative uniform coarse index maps with Mayer-Vietoris boundary maps. Again, this is an abstract version of the (non-relative) uniform partitioned-manifold index theorem.
\begin{cor}\label{MVDiagRoe}
Let $X$ be a proper metric space of jointly bounded geometry, and $X=Y_1\cup Y_2$ a coarsely excisive and bounded Mayer-Vietoris decomposition. Then, the uniform coarse index map takes the Mayer-Vietoris sequence in the uniform K-homology of the decomposition $X=Y_1\cup Y_2$ to the Mayer-Vietoris in the K-theory of the uniform Roe algebras of this decomposition. In particular the diagram\\
\centerline{\xymatrix{
K^{u}_{*+1}(X) \ar[r]^-{\delta} \ar[d]_{\Ind} &  \ar[d]^{\Ind} K^{u}_*(Y_1\cap Y_2) \\
K_{*+1}(C^*_u(X)) \ar[r]_-{\delta} & K_*(C^*_u(Y_1\cap Y_2))
}}
commutes.
\end{cor}
There is also an analogous commutativity statement connecting the boundary maps of the long exact sequence in uniform K-homology and of the algebraic Mayer-Vietoris sequence. This statement does not follow from Theorem \ref{AbstractRelPM}, but it is proved similarly by unraveling the definitions of the objects involved, and appealing to naturality of the K-theory boundary map. For the sake of completeness we give the calculation.
\begin{thm}\label{IndexCommAlgMV}
Let $X$ be a proper metric space of jointly bounded geometry, and let $Z\subseteq X$ be a subspace such that $(X,Z)$ is transversally controlled. Then, the diagram \\
\centerline{\xymatrix{
K^{u}_{*+1}(X,Z) \ar[r]^\partial \ar[d]_{\relInd} & K^{u}_*(Z) \ar[d]^-{\Ind} \\
K_{*+1}\left( D^*_u(X,Z)\right) \ar[r]_-\delta & K_*\left( C^*_u(Z) \right)
}} 
commutes.
\end{thm}
\begin{proof}
Represent the algebras using an ample representation $(H,\sigma)$ such that $(H_Z:=P_ZH, \sigma_Z:=P_Z\sigma(\mathbf{1}_Z\cdot))$ is an ample representation of $Z$, where $P_Z:=\sigma(\mathbf{1}_Z)$. The upshot of working with such representations is this: There is an upper-corner inclusion of $\fD^{u}(Z) \subseteq \fB(H_Z)$ into $\fD^{u}(X,Z) \subseteq \fB(H=H_Z\oplus H^\perp_Z)$. Inspection of the proof of Proposition \ref{InclusionIsoDual} reveals that this upper-corner inclusion $\fD^{u}(Z)\hookrightarrow \fD^{u}(X,Z)$ realizes the isomorphism $K_*(\fD^{u}(Z))\cong K_*(\fD^{u}(X,Z))$. Restricting the upper corner inclusion gives rise to inclusions $D^*_u(Z)\hookrightarrow D^*_u(Z\subseteq X,Z)$ and $C^*_u(Z)\hookrightarrow C^*_u(Z\subseteq X)$, with the latter inducing an isomorphism on K-theory by Proposition \ref{UniRoeSubset}. Thus there is also a well-defined map 
\begin{equation}\label{TauQuot}
D^*_u(Z) / C^*_u(Z) \longrightarrow D^*_u(Z\subseteq X,Z) / C^*_u(Z\subseteq X) \; .
\end{equation}
For the sake of efficient notation we will denote all of the maps given or induced by the upper-corner inclusion by $\tau$.

Observe that by Proposition \ref{FinPropNormKHom} and naturality of the boundary map we have a commutative diagram\\
\centerline{\xymatrix{
K_{*+2}\left( D^*_u(X)/D^{*}_u(X,Z) \right) \ar@{=}[r] \ar[d]_{\relInd} & K_{*+2}\left( \fD^{u}(X)/\fD^{u}(X,Z) \right) \ar[d]^\partial \\
K_{*+1}\left( D^*_u(X,Z)\right) \ar[r]_{\iota_*} & K_{*+1}\left( \fD^{u}(X,Z) \right)
}}
where the top line represents $K^{u}_*(X,Z)$, and $\iota: \, D^*_u(X,Z) \to \fD^{u}(X,Z)$ is the inclusion map. It therefore suffices to prove that that the diagram\\
\begin{equation}\label{DiagramMV}
\xymatrix{
K_{*+1}\left( D^*_u(X,Z)\right) \ar[rr]^{\iota_*} \ar[dr]_{\delta} & & K_{*+1}\left( \fD^{u}(X,Z) \right) \ar[dl]^{\Ind} \\
& K_*(C^*_u(Z) & }
\end{equation}
commutes. Let us unravel the definitions of $\delta$ and $\Ind$. The former is given by the composition
\begin{align*}
K_{*+1}\left( D^*_u(X,Z) \right) &\xrightarrow{\omega_*} K_{*+1}\left( D^*_u(X,Z) / C^*_u(X)\right) \\
&= K_{*+1}\left( D^*_u(Z\subseteq X,Z) / C^*_u(Z\subseteq X)\right) \\
& \xrightarrow{\partial} K_*\left(C_u^*(Z\subseteq X)\right) \\
&\xrightarrow{\sim} K_*\left( C^*_u(Z) \right) \; ,
\end{align*}
with $\omega$ the quotient map. We refer to Appendix \ref{AppendixKTheory} for this description of the Mayer-Vietoris boundary map. The latter is given by 
\begin{align*}
K_{*+1}\left( \fD^{u}(X,Z) \right) &\xrightarrow{\sim} K_{*+1}\left( \fD^{u}(Z) \right) \\
&\xrightarrow{\sim} K_{*+1}\left( \fD^{u}(Z)/\fC^{u}(Z)\right) \\
&= K_{*+1}\left( D^*_u(Z)/C_u^*(Z)\right) \\
& \xrightarrow{\partial} K_*\left(C^*_u(Z)\right) \; .
\end{align*}
As discussed above the first arrow is inverse to the upper-corner inclusion. Recall also that the upper-corner inclusion induced the map \eqref{TauQuot} on quotients. Combining this with the definitions of $\delta$ and $\Ind$ we write diagram \eqref{DiagramMV} as 
\begin{equation}\label{DiagramMV2}
\xymatrix{
K_{*+1}\left( D^*_u(X,Z)\right) \ar[r]^{\iota_*} \ar[d]_{\omega_*} & K_{*+1}\left( \fD^{u}(X,Z) \right) \ar[d]^{\sim} \\
K_{*+1}\left( D_u^*(Z\subseteq X,Z)/C_u^*(Z\subseteq X)\right) \ar[d]_\partial & \ar[l]_-{\tau_*} K_{*+1}\left( D^*_u(Z)/C_u^*(Z)\right) \ar[d]^\partial \\
K_*\left( C^*_u(Z\subseteq X) \right) & \ar[l]_\sim K_*\left( C_u^*(Z) \right) }
\end{equation}
where the bottom horizontal isomorphism is also induced by the upper-corner inclusion. The lower square is commutative by naturality of the boundary map. Therefore it suffices to show commutativity of the upper square. 

To that end note first that
$$
D_u^*(Z\subseteq X,Z)/C_u^*(Z\subseteq X) = D^*_u(X,Z) / C^*_u(X) = \fD^{u}(X,Z)/\fC^{u}(X)\; .
$$
The first equality is obtained via cut-offs away from $Z$, and was in fact already used for the algebraic Mayer-Vietoris boundary map above. The second is obtained via truncation. We denote the quotient maps $\fD^{u}(X,Z)\to\fD^{u}(X,Z)/\fC^{u}(X)$ and $\fD^{u}(Z)\to \fD^{u}(Z)/\fC^{u}(Z)$ by $\omega$ as well. The diagram\\
\centerline{\xymatrix{
K_{*+1}\left( \fD^{u}(Z) \right) \ar[r]^-{\tau_*} \ar[d]_{\omega_*} & K_{*+1}\left( \fD^{u}(X,Z) \right) \ar[d]^{\omega_*} \\
K_{*+1}\left( \fD^{u}(Z)/ \fC^{u}(Z)\right)\ar[r]_-{\tau_*} & K_{*+1}\left(  \fD^{u}(X,Z) /  \fC^{u}(X)\right)
}}
commutes, where we recall that $\tau$ always denotes maps induced by upper-corner inclusion. Recall also that the upper-right isomorphism in \eqref{DiagramMV2} is given by 
$$
K_{*+1}\left( \fD^{u}(X,Z) \right) \xrightarrow{\tau_*^{-1}} K_{*+1}\left( \fD^{u}(Z) \right) \xrightarrow{\omega_*}K_{*+1}\left( \fD^{u}(Z)/ \fC^{u}(Z)\right) = K_{*+1}\left( D^*_u(Z)/ C^*_u(Z)\right) \; .
$$
The combination of these facts reduces the upper square in \eqref{DiagramMV2} to the commutativity of \\
\centerline{\xymatrix{
K_{*+1}\left( D^*_u(X,Z)\right) \ar[r]^{\iota_*} \ar[d]_{\omega_*} & K_{*+1}\left( \fD^{u}(X,Z) \right) \ar[d]^{\omega_*} \\
 K_{*+1}\left( D^*_u(X,Z) / C^*_u(X)\right) \ar@{=}[r]& K_{*+1}\left( \fD^{u}(X,Z) / \fC^{u}(X) \right)
}}
which holds even before passing to K-theory. This concludes the proof.
\end{proof}
We conclude this chapter by revisiting the uniform coarse index map $K_*^{u}(X,Z)\to K_*(C^*_u(X,Z))$. Indeed, we have not made good on our promise to provide a proof of Theorem \ref{CommUniIndex}. We will deduce it as a consequence of Theorem \ref{IndexCommAlgMV}. First, we need a lemma connecting the relative uniform structure algebra to the relative uniform Roe algebra, which was already present in the proof of the algebraic Mayer-Vietoris sequence.
\begin{lemma}
Let $X$ be a proper metric space of jointly bounded geometry, and $Z$ a closed and boundedly embedded subspace. The inclusion induces an isomorphism
$$
C^*_u(X)/C^*_u(Z\subseteq X) \xlongrightarrow{\sim} D^*_u(X,Z)/D^*_u(Z\subseteq X,Z) \; .
$$
Thus, the composition
$$
D^*_u(X,Z) \longrightarrow D^*_u(X,Z)/D^*_u(Z\subseteq X,Z) \xlongleftarrow{\sim} C^*_u(X)/C^*_u(Z\subseteq X)
$$
produces a natural forgetful map
$$
D^*_u(X,Z) \xlongrightarrow{F} C^*_u(X,Z) \; .
$$
\end{lemma}
We show that the index map to $K_*(C^*_u(X,Z))$ factors through that to $K_*(D^*_u(X,Z))$ followed by the forgetful map.
\begin{lemma}\label{ComparisonRelInd}
Let $X$ be a proper metric space of jointly bounded geometry, and $Z$ a closed and boundedly embedded subspace. The diagram \\
\centerline{\xymatrix{
 & \ar[dl]_{\relInd} K^{u}_*(X,Z) \ar[dr]^{\Ind} & \\
 K_*\left( D^*_u(X,Z) \right) \ar[rr]_F & & K_*\left(C^*_u(X,Z)\right)
}}
commutes.
\end{lemma}
\begin{proof}
This follows from the naturality of the boundary map. Indeed, there is a commutative square \\
\centerline{\xymatrix{
K_{*+1}\left( D^*_u(X)/D^*_u(X,Z) \right) \ar[r] \ar[d]_\partial & K_*\left( \left( D^*_u(X)/D^*_u(Z\subseteq X) \right)/C^*_u(X,Z) \right) \ar[d]^\partial \\
K_{*}\left( D^*_u(X,Z) \right) \ar[r] & K_*\left(C^*_u(X,Z) \right)
}}
where the horizontal arrows are induced by the quotient map $D^*_u(X)\to D^*_u(X)/D^*_u(Z\subseteq X)$, and the lower one uses the isomorphism $D^*_u(X,Z)/D^*_u(Z\subseteq X,Z)\cong C^*_u(X,Z)$. Going right then down in the square gives the uniform coarse index map $K^{u}_*(X,Z)\to K_*(C^*_u(X,Z))$, while going down then right gives the relative uniform index map $K^{u}_*(X,Z)\to K_*(D^*_u(X,Z))$ followed by the forgetful map $K_*(D^*_u(X,Z)) \to K_*(C^*_u(X,Z))$.
\end{proof}
We come to the promised proof of Theorem \ref{CommUniIndex}, which stated that the uniform coarse index maps take the long exact sequence in uniform K-homology to that in K-theory of the uniform Roe algebras.
\begin{proof}[Proof of Theorem \ref{CommUniIndex}]
We want to prove that the diagram  \\
\centerline{\xymatrix{
K_*^{u}(Z) \ar[r] \ar[d]_{\Ind} & K^{u}_*(X) \ar[r] \ar[d]_{\Ind} & K^{u}_*(X,Z) \ar[r]^\partial \ar[d]_{\Ind} & K_{*-1}^{u}(Z) \ar[d]_{\Ind} \\
K_*\left(C^*_u(Z)\right) \ar[r] & K_*\left(C^*_u(X)\right) \ar[r] & K_*\left(C^*_u(X,Z)\right) \ar[r]_\partial & K_{*-1}\left(C^*_u(Z)\right)
}}
commutes. The first square has nothing to do with relative indices, its commutativity is due to the naturality of the (non-relative) uniform coarse index map. Commutativity of the second square is part of the statement of Proposition \ref{RelRoeIndex}. Lastly, using Lemma \ref{ComparisonRelInd} write the third square as \\
\centerline{\xymatrix{
K^{u}_*(X,Z) \ar[rr]^\partial \ar[d]_{\relInd} && K_{*-1}^{u}(Z) \ar[d]^{\Ind} \\
K_*\left(D^*_u(X,Z)\right) \ar[r]_F & K_*\left(C^*_u(X,Z)\right) \ar[r]_\partial & K_{*-1}\left(C^*_u(Z)\right)
}}
The description \eqref{MVformula} of the Mayer-Vietoris boundary on operator K-theory lets us deduce that the composition $\partial\circ F: \, K_*\left(D^*_u(X,Z)\right) \to K_{*-1}\left(C^*_u(Z)\right)$ coincides with the algebraic Mayer-Vietoris boundary map $\delta: \, K_*\left(D^*_u(X,Z)\right) \to K_{*-1}\left(C^*_u(Z)\right)$. Thus, commutativity of the last square is a consequence of Theorem \ref{IndexCommAlgMV}.
\end{proof}
We summarize the situation in the following diagram: \\
\centerline{\xymatrix{
K_*^{u}(Z) \ar[r] \ar[d]_{\Ind} & K^{u}_*(X) \ar[r] \ar@{.>}[d]_{\Ind\oplus 0} & K^{u}_*(X,Z) \ar[r]^\partial \ar[d]_{\relInd} & K_{*-1}^{u}(Z) \ar[d]_{\Ind} \\ 
K_*\left(C^*_u(Z)\right) \ar[r] \ar@{=}[d] & K_*\left(C^*_u(X)\right)\oplus K_*\left(D^*_u(Z)\right) \ar@{.>}[d]_{\id\oplus 0} \ar[r] & K_*\left(D^*_u(X,Z)\right) \ar[d]_F \ar[r]^\delta & K_{*-1}\left(C^*_u(Z)\right)\ar@{=}[d] \\
K_*\left(C^*_u(Z)\right) \ar[r] & K_*\left(C^*_u(X)\right) \ar[r] & K_*\left(C^*_u(X,Z)\right) \ar[r]_\partial & K_{*-1}\left(C^*_u(Z)\right)
}}
where the first row is the long exact sequence in uniform K-homology, the second row the algebraic Mayer-Vietoris sequence for $D^*_u(X,Z)$, and the third row the long exact sequence for the uniform Roe algebras. The dotted vertical arrows indicate that the small squares containing them do \emph{not} commute due to the presence of the summand $K_*(D^*_u(Z))$. Were one to omit that summand the entire diagram would commute, but the second row would no longer be exact. The large squares from the first to the third row do however commute.

The relative uniform structure algebra together with the relative uniform index map and the algebraic Mayer-Vietoris sequence thus behave somewhat unexpectedly, at least from a homotopy-theoretical perspective. From that viewpoint one might expect $K_*(D^*_u(X,Z))$ to be the relative part of some homology theory, and the relative uniform index map to be the relative part of a transformation from uniform K-homology to this other homology theory. However, instead of some long exact sequence involving $K_*(D^*_u(X,Z))$ one has the algebraic Mayer-Vietoris sequence, and the maps from the uniform K-homology exact sequence to this sequence do not generally commute. Only if one passes further to relative uniform Roe algebra does one get the expected long exact sequence, and the uniform coarse index map acts as a transformation from one exact sequence to the other. However, as already noted above passing to relative uniform Roe algebra forgets any information located near the boundary, so that this is a trade-off between structural compatibility and informational content.  
\chapter{Sobolev spaces of domains in manifolds of bounded geometry}
Sobolev spaces are well-behaved Hilbert spaces of (weakly) differentiable functions, and as such are ubiquitous in many corners of mathematics. In our particular context they arise as domains of elliptic operators. Operator-theoretic properties of elliptic operators are therefore often established by appeal to the machinery of Sobolev spaces. Our discussion of uniform K-homology classes arising from uniformly elliptic operators will be no different. 

Sobolev spaces on manifolds of bounded geometry can be defined using local coordinate charts that are adapted to the geometry \cite{Shubin1992}. These Sobolev spaces have the same properties as the Sobolev spaces on compact manifolds\footnote{The Rellich-Kondrachov theorem is the clear exception.} or on Euclidean space. Indeed, using the definition via local coordinate charts one can always localize to a coordinate patch, appeal to facts about Sobolev spaces on Euclidean space, and afterwards glue together to get the result globally. That last step requires bounded geometry to ensure that gluing local objects over the entire space yields a bounded object. This will be no different for Sobolev spaces on domains. Given this essentially universal strategy we will not give evry proof in detail.

When it comes to manifolds with boundary and bounded geometry there are two approaches. One can consider a manifold with boundary intrinsically, and impose boundedness contitions on the geometry in the interior and on a boundary collar. This is done in \cite{Schick1998}. Alternatively, one can understand a manifold with boundary and bounded geometry externally as an open subset of a manifold of bounded geometry (without boundary), such that the boundary (in the ambient space) satisfies certain boundedness criteria. This is done in \cite{Ammann2019}, where it is also shown that both approaches are equivalent. The two approaches lead to different but equivalent definitions of Sobolev spaces.

Recall that we endeavor to establish the existence of uniform K-homology classes of uniformly elliptic operators without assumptions on the boundary of the domain. This makes it more natural for us to use the external approach to manifolds with boundary and boundary, since it is easier in that setting to relax the conditions on the boundary.

We will eventually establish the uniformness conditions of uniformly elliptic operators required to define uniform Fredholm modules by appealing a strengthened version of the Rellich-Kondrachov theorem that gives uniform summability instead of mere compactness of certain Sobolev embeddings. The fact that Sobolev embeddings on compact domains are not only compact but finitely summable is classical. The additional information we require is that the Schatten norm of the embedding can be uniformly bounded for domains in a manifold of bounded geometry. This uniformly summable Rellich-Kondrachov theorem is the content of Section \ref{SectionQuantComp}. Moreover, we will use integral representations to reduce discussion of functions of operators to resolvents, where the uniformly summable Rellich-Kondrachov theorem can be applied more directly. To make these integrals converge in the correct norm it is necessary to trade in a small bit of Sobolev degree. Thus, even if we were initially content to work only with integer-degree Sobolev space, we are forced at this point to work with Sobolev spaces of fractional degree as well. 

Our discussion of Sobolev spaces on domains in manifolds of bounded geometry is thus required to encompass domains with potentially ill-behaved boundary, and Sobolev spaces of fractional degree. These cases require an unfortunate amount of technical overhead in order to formulate and prove the desired results. The main obstacle is this: One usually treats fractional Sobolev spaces as complex interpolation spaces between integer-degree Sobolev spaces. Even in the case of domains with smooth boundary this is not without subtlety. There are multiple families of Sobolev spaces built from funcions supported in the given domain, and these families are mixed by interpolation. Even with smooth boundary these families generally do not coincide, and with rough boundary the situation is only worse. 

This chapter progresses as follows: We begin with a review of manifolds of bounded geometry in Section 4.1, where we also introduce notions of regularity for domains in such manifolds. Then we recall the basics of Sobolev spaces on manifolds of bounded geometry (without boundary) in Section 4.2. We proceed to introduce the various Sobolev spaces on domains in Section 4.3, and prove comparison statements between them. Section 4.4 is devoted to extension and trace maps. Then we discuss the interpolation properties in Section 4.5. Lastly, in Section 4.6 we discuss the uniformly summable Rellich-Kondrachov theorem.
\section{Manifolds of bounded geometry} \label{SectionBoundedGeometry}
Let us begin by recalling the definition of a manifold of bounded geometry.
\begin{definition}
A Riemannian manifold $M$ is said to have \emph{bounded geometry} if it has positive injectivity radius $\mathrm{inj}_M>0$, and its curvature and its derivatices are bounded, i.e. $\|\nabla^k R\|_\infty <\infty$ for all $k$.
\end{definition}
Note that the requirement of a positive injectivity radius implies that a manifold of bounded geometry is complete. In particular it has no boundary.
\begin{example}
Closed manifolds always have bounded geometry. A cover of a manifold with bounded geometry (endowed with the pullback metric) has bounded geometry, in particular the universal cover of a closed Riemannian manifold has bounded geometry. Homogeneous spaces with invariant metrics have bounded geometry. Manifolds with a Lie structure at infinity have bounded geometry, see for example \cite{Ammann2004}.
\end{example}
We recall the definition of vector bundles of bounded geometry.
\begin{definition}
A vector bundle $E\to M$ over a manifold of bounded geometry equipped with a bundle metric and a compatible connection is said to have \emph{bounded geometry} if the curvature of the connection and all its derivatives are bounded, i.e. $\|\nabla^k R_E \|_\infty<\infty$.
\end{definition}
This definition works equally well for real and for complex vector bundles. When it comes to Sobolev spaces we will restrict to complex vector bundles though. Note however that the complexification of a real bundle of bounded geometry again has bounded geometry.
\begin{example}
Let $M$ have bounded geometry. The tangent and cotangent bundles of $M$ have bounded geometry, as does any tensor bundle or exterior power. If $M$ is spin, then the Clifford bundle and spinor bundles have bounded geometry as well. Trivial bundles over $M$ have bounded geometry when equipped with a constant bundle metric and a flat connection. If $\pi: \, \tilde{M}\to M$ is a covering map and $E\to M$ has bounded geometry, then the pullback $\pi^*E\to \tilde{M}$ has bounded geometry. Direct sums, tensor products, exterior products, etc. of vector bundles of bounded geometry have bounded geometry.
\end{example}
Sobolev spaces on manifolds of bounded geometry are defined using well-behaved covers by normal coordinate charts. The existence of such covers is guaranteed by the following lemma.
\begin{lemma}[{\cite[Appendix 1]{Shubin1992}}\footnote{Shubin only proves the statement for $C=2$, but the argument may directly be adapted to general $C$.}]\label{CoverBdGeom}
Let $C>1$. For every $0<r<\frac{1}{C+1}\cdot \mathrm{inj}_M$ there exists a countable subset $(x_i)_{i\in I} \subseteq M$ such that the following holds:
\begin{itemize}
\item $d(x_i,x_j)\geq r$ if $i\neq j$,
\item $\{ B_r(x_i)\}_{i}$ is an open cover of $M$,
\item  $\{ B_{Cr}(x_i)\}_{i}$ is also an open cover of $M$, and there is $N\in \ZZ_{\geq 1}$ such that each point of $x$ lies in at most $N$ elements of the cover.
\end{itemize}
\end{lemma}
The geometric data of manifolds and vector bundles of bounded geometry may be converted to uniformly bounded local data using such well-behaved covers. Let us explain first how to get from normal coordinate charts to trivializations of a vector bundle. Let $E\to M$ be a smooth vector bundle equipped with a bundle metric and a metric connection. In a normal coordinate patch $B_r(x)$ one obtains a trivialization of $E$ by choosing a basis of the fiber over $x$ and moving this basis along radial geodesics via parallel transport. A trivialization of $E$ obtained this way will be called a \emph{synchronous framing}. 
\begin{lemma}\label{BGLocalData}
Let $M$ and $E\to M$ have bounded geometry\footnote{We take $E$ to be a complex vector bundle, though the statement is equally true for real vector bundles.}, and let $\{ U_i=B_{Cr}(x_i) \}_i$ be an open cover as in Lemma \ref{CoverBdGeom}. Let $\phi_i: \, U_i\times \CC^n \to E|_{U_i}$ be synchronous framings of $E$. Then:
\begin{itemize}
\item The normal coordinate maps $\kappa_i: \, B_{Cr}(0)\to U_i$ are bi-Lipschitz equivalences with Lipschitz constants uniformly bounded independently of $i$.
\item The transition functions $\kappa_j^{-1}\circ\kappa_i$ have uniformly bounded derivatives, i.e. $\|\partial^\alpha \kappa_j^{-1}\circ\kappa_i\|_\infty \leq C(\alpha)<\infty$.
\item The transition functions $\phi_j^{-1}\circ \phi_i$ have uniformly bounded derivatives in normal coordinates, i.e. $\|\partial^\alpha \phi_j^{-1}\circ \phi_i \circ (\kappa_i^{-1}\times \id_{\CC^n})\|\leq C(\alpha)<\infty$.
\item The bundle metric and the curvature tensor of $E$ as well as their derivatives are analogously uniformly bounded in synchronous framings and normal coordinates.
\item There exists a partition of unity $(\psi_i)_{i\in I}$ subordinate to $\{U_i\}_i$ such that $\supp(\psi_i)\subseteq B_{r}(x_i)$ and the derivatives of the $\psi_i$ are uniformly bounded, i.e $\|D^\alpha (\psi_i \circ \kappa_i^{-1})\|_\infty\leq C(\alpha)<\infty$.
\end{itemize}
\end{lemma}
\begin{proof}
These results are well-known, see for example \cite[Appendix A.1]{Shubin1992} and \cite[Section 2]{Engel2018}.
\end{proof}
The existence of uniformly local data is in fact equivalent to bounded geometry of the global object. Indeed, if $M$ and $E$ satisfy the conclusions of Lemmas \ref{CoverBdGeom} and \ref{BGLocalData} (and $M$ has positive injectivity radius), then they have bounded geometry. See for example \cite{Engel2018} and references therein.

We shall refer to the data $(U_i,\kappa_i,\psi_i,\phi_i)_{i\in I}$ as an \emph{admissible trivialization}\footnote{In \cite{Grosse2013} Große and Schneider introduce a more general notion of an admissible trivialization. What we call an admissible trivialization here, they call a geodesic trivialization. They show that using an admissible trivialization in their more general sense yields the same Sobolev spaces.} of $M$ and $E$, where $\{U_i=B_{Cr}(x_i)\}_{i\in I}$ is a cover by normal coordinate patches as in Lemma \ref{CoverBdGeom}, $\kappa_i$ are the corresponding normal coordinate maps, $(\psi_i)_i$ is a subordinate partition of unity, and $\phi_i$ are corresponding synchronous framings of $E$. We shall be lax in our usage of the term admissible trivialization. If there is no vector bundle involved, we will refer to the data $(U_i,\kappa_i,\psi_i)_i$ as an admissible trivialization of $M$, or if the partition of unity is not relevant, to the data $(U_i,\kappa_i)$.

Let us turn to domains in manifolds of bounded geometry, as well as different notions of boundedness and regularity of their boundary. By a \emph{domain} in the bounded-geometry manifold $M$ we mean an open subset $\Omega\subseteq M$. The first notion we will introduce is the external definition of manifold with boundary and bounded geometry as introduced in \cite{Ammann2019}.
\begin{definition}[{\cite[Definition 2.4]{Ammann2019}}]\label{BGHypersurface}
Let $M$ be a manifold of bounded geometry. A \emph{bounded-geometry hypersurface} in $M$ is a codimension-1 submanifold $N\subseteq M$ without boundary with the following properties:
\begin{itemize}
\item[(i)] $N$ is closed in $M$.
\item[(ii)] The restriction of the metric on $M$ makes $N$ into a manifold of bounded geometry.
\item[(iii)] The second fundamental form of $N$ and all its derivatives are bounded.
\item[(iv)] There exists $\delta>0$ such that the exponential map is defined and injective on $\{ \xi\in T_xM \, | \, x\in N, \, \xi\perp T_xN, \, ||\xi||<\delta \}$.
\end{itemize}
A domain $\Omega \subseteq M$ is said to have \emph{bounded geometry} if $\partial \Omega$ is a bounded geometry hypersurface in $M$. In this case $\olO$ is said to be a \emph{manifold with boundary and bounded geometry}.
\end{definition}
Let us reiterate that by a manifold of bounded geometry we always mean one without boundary. A manifold (or domain) $\Omega$ satisfying the above definition will always be referred to as a \emph{manifold with boundary and bounded geometry}. On the other hand, we will use the terms \emph{bounded-geometry domain} and \emph{manifold with boundary and bounded geometry} interchangeably, ignoring the minor technical point that the definition of the latter includes the boundary while the former does not.

As was already mentioned in the introduction to this chapter there is also an intrinsic definition of manifolds with boundary and bounded geometry given by Schick \cite{Schick2000}. The intrinsic definition takes a similar form to the definition of manifold of bounded geometry without boundary via bounds on the relevant geometric objects. There is also an equivalent local definition. It is proved in \cite{Ammann2019} that a manifold $\olO$ with boundary has bounded geometry in the intrinsic sense if and only if there exists an isometric embedding of $\olO$ into a manifold $M$ of bounded geometry (without boundary) such that the image of this embedding satisfies the definition given above.

Relaxing the smoothness requirements for a bounded-geometry domain, we introduce the notion of a (uniform) $C^{k,\alpha}$-domain. In the case of domains in Euclidean space with compact boundary this notion is already well-established in the study of Sobolev spaces. See for example \cite{McLean2000}, and also \cite{Grisvard1985} for the equivalence to so-called conic domains.
\begin{definition}\label{DefCKAlpha}
Let $k\in \ZZ_{\geq 0}$ and $\alpha\in [0,1]$. Let $M$ be an $n$-dimensional manifold of bounded geometry. A domain $\Omega\subseteq M$ is called a \emph{$C^{k,\alpha}$-domain} if the following holds: There exists a subset $\{ x_i \}_{i\in I}\subseteq M$ and radii $r_i,R_i$ with $0<\sqrt{n}r_i<R_i<\mathrm{inj}(M)$ such that $\{ B_{r_i}(x_i) \}_{i\in I}$ and $\{ B_{R_i}(x_i) \}_{i\in I}$ are locally finite collections of open sets with $\partial\Omega \subseteq \bigcup_{i\in I} B_{r_i}(x_i)$ and there are closed cubes $K_i\subset\RR^n$ centered on $0$ with side length $2c_i$ such that $B_{r_i}(0)\subset K_i \subset B_{R_i}(0)$ so that there is an orthogonal transformation taking $\kappa_i\left(\Omega \cap B_{R_i}(x_i)\right) \cap K_i$ to the set 
$$
\left\{ (x_1,\cdots,x_n)\in [-c_i,c_i]^n \, | \, x_n < \xi_i(x_1,\cdots, x_{n-1}) \right\}
$$ 
for a function $\xi_i \in C^{k,\alpha}([-c_i,c_i]^{n-1},\RR)$ taking values in $[-r_i,r_i]$.

Moreover, $\Omega$ is called a \emph{uniform $C^{k,\alpha}$-domain} if the radii $r_i,R_i$ and the side length $2c_i$ can be chosen to be independent of $i$, the cover $\{ B_{r_i}(x_i) \}_{i\in I}$ can be completed to an admissible trivialization of $M$,\footnote{More precisely, there exist $r>0$, $C>\sqrt{n}$ and a subset $\{ y_j\}_{j\in J}$ which satisfies the properties of Lemma \ref{CoverBdGeom} such that $\{x_i\}_{i\in I}\subseteq \{y_j\}_{j\in J}$, and $r_i=r$ and $R_i=R=Cr$ for all $i$.} and it holds that $\sup_{i\in I} ||\xi_i||_{C^{k,\alpha}}<\infty$.
\end{definition}
\begin{figure}
\centering
\includegraphics[scale=.6]{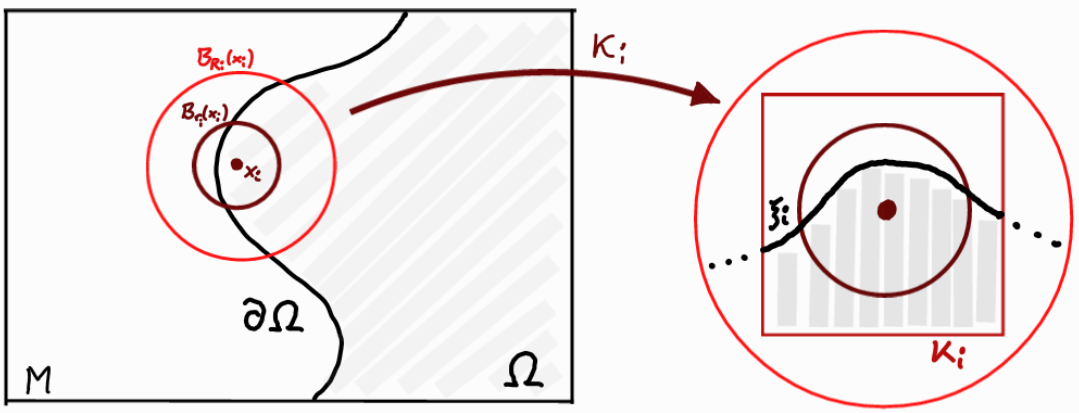}
\caption{Sketch of a uniform $C^{k,\alpha}$-domain.}\label{SketchCkAlphaDomain}
\end{figure}
If $\alpha=0$, we will simply speak of a \emph{uniform $C^k$-domain}. If $k=0$ and $\alpha=1$ we will speak of a \emph{uniform Lipschitz domain}, since $C^{0,1}$ consists precisely of Lipschitz functions, and the condition $\sup \| \xi_i \|_{C^{0,1}}<\infty$ expresses that the Lipschitz constants of the $\xi_i$ should be uniformly bounded.

While this definition might be rather technical, on an intuitive level it simply expresses that near its boundary the domain locally looks like the area below the graph of a $C^{k,\alpha}$-function. Uniformness then means that the norm of these $C^{k,\alpha}$-functions should be bounded, and that the boundary is globally well-behaved. Figure \ref{SketchCkAlphaDomain} might help parsing the definition.
\begin{example}
If $\partial \Omega$ is compact and locally given by the area under a $C^{k,\alpha}$-function, then $\Omega$ is uniformly $C^{k,\alpha}$. If $\Omega$ has bounded geometry, then it is uniformly $C^{k,\alpha}$ for all $k$ and $\alpha$. In that case if $U_\delta \cong \partial\Omega \times (-\delta,\delta)$ is a uniformly thick tubular neighborhood for sufficiently small $\delta>0$, then a uniformly $C^{k,\alpha}$-bounded perturbation of $\partial \Omega$ to a $C^{k,\alpha}$-function yields a uniformly $C^{k,\alpha}$-domain. Any reasonable definition of a manifold with corners and bounded geometry should yield a uniform Lipschitz domain, though the author has not seen such manifolds discussed in the literature. Non-examples include domains that (locally) look like the area under a sufficiently irregular function, or domains that do not look like the area under a function at all. The slit plane $\RR^2\setminus (-\infty,0]$ is an example of the latter.
\end{example}
\section{Sobolev spaces on manifolds of bounded geometry}
In this section we recall the definition and basic properties of Sobolev spaces of sections of $E$. For our presentation we draw on \cite{Grosse2013}. Let $M$ be a manifold of bounded geometry and let $E\to M$ be a vector bundle of bounded geometry\footnote{Recall that unless otherwise stated manifolds of bounded geometry do not have a boundary, and vector bundles are complex.}. We let $\Gamma_c(M;E)$ denote the space of compactly supported smooth sections of $E$. Fix an admissible trivialization $(U_i,\kappa_i,\psi_i,\phi_i)_{i\in I}$ of $M$ and $E$. For $u\in\Gamma_c(M;E)$ let $u_i := \left( \psi_i \cdot (\phi_i \circ u)\right) \circ \kappa_i^{-1}$. This is a compactly supported, smooth vector-valued function on $\RR^n$, where $n$ is the dimension of $M$. Let $s\in\RR$, and define the norm
$$
\|u\|_{H^s(M;E)} := \sum_{i\in I} \| u_i \|_{H^s(\RR^n)} \; .
$$
Here, $\| u_i \|_{H^s(\RR^n)}$ is the Sobolev norm for vector-valued functions on $\RR^n$. It is defined by $\| u_i \|_{H^s(\RR^n)}=\|(1+|\xi|^2)^\frac{s}{2} \hat{u}_i(\xi)\|^2_{L^2(\RR^n)}$.
\begin{definition}
Define $H^s(M;E)$ as the closure of $\Gamma_c(M;E)$ with respect to the norm $\|\cdot\|_{H^s(M;E)}$.
\end{definition}
The space $H^{s}(M;E)$ is independent of the choice of admissible trivializations, and any two admissible trivializations give rise to equivalent $H^{s}$-norms. The Sobolev spaces on manifolds with bounded geometry behave just as an $\RR^n$: If $s\geq 0$, then $H^{s}(\RR^n) \subseteq L^2(M;E)$ with equality for $s=0$. The space $\Gamma_c(M;E)$ is dense in each $H^{s}(M;E)$ by construction. Moreover, the norm $\|\cdot\|_{H^s(M;E)}$ is induced by the inner product
$$
\langle u,v\rangle_{H^s(M;E)} := \sum_{i\in I} \langle u_i,v_i \rangle_{H^s(\RR^n)} \; ,
$$
making $H^s(M;E)$ a Hilbert space. The $L^2$-inner product extends to a perfect pairing between $H^s(M;E)$ and $H^{-s}(M;E)$, identifying the dual space of one with the other.

Let us note that $H^s(M;E)$ can also be defined directly as a space of distributions. Indeed, if $\varphi$ is $E$-valued distribution on $M$, i.e an element of the dual space of $\Gamma_c(M;E)$, then $\varphi_i:= \left( \psi_i \cdot (\phi_i \circ \varphi_i)\right) \circ \kappa_i^{-1}$ is a compactly supported vector-valued distribution on $\RR^n$. Then, $H^s(M;E)$ coincides with the space of all distributions $\varphi$ such that $\varphi_i\in H^s(\RR^n)$, and $\sum_i \|\varphi_i\|_{H^s(\RR^n)}<\infty$. This can be seen using an approximation argument via mollifiers.

If $s=k$ is a non-negative integer, then the Sobolev space of degree $k$ can also be obtained by completing $\Gamma_c(M;E)$ with respect to a norm directly measuring the first $k$ derivatives. This is the content of the next proposition.
\begin{prop}[\cite{Grosse2013}]\label{SobolevNormConnection}
Let $k\in \ZZ_{\geq 0}$. Then, the Sobolev norm $\|u\|_{H^k(M;E)}$ is equivalent to the norm
$$
\sum_{m=0}^k \left( \int_M \|\nabla^m u\|^2 \, d\mathrm{vol} \right)^{\frac{1}{2}} \; ,
$$
where $d\mathrm{vol}$ denotes the volume density of the metric on $M$. Thus $H^k(M;E)$ coincides with the closure of $\Gamma_c(M;E)$ with respect to that norm.
\end{prop}
Next we discuss multiplication by functions.
\begin{definition}
Let $k\in \ZZ_{\geq 0}$. Define $C^{k}_b(M)$ to be the space 
$$
C^k_b(M) := \bigg\{ f\in C^k(M) \, | \, \| u \|_{C^{k}_b(M)} := \sum_{m=0}^k \| \nabla^m u \|_\infty < \infty \bigg\} \; .
$$ 
\end{definition}
\begin{lemma}\label{MultWInftyNoBdry}
Let $k\in \ZZ_{\geq 0}$, and $f\in C^k_b(M)$. Then
$$
H^s(M;E) \longrightarrow H^s(M;E) \; , \; u\longmapsto fu
$$
is a bounded map for all $-k\leq s\leq k$ with norm bounded by $\|f\|_{C^k_b(M)}$.
\end{lemma}
\begin{proof}
Let $f \in C^k_b(M)$ and $u\in \Gamma_c(M;E)$. The Leibniz rule implies that 
$$
||\nabla^m(fu)||_{L^2} \lsim \sum_{j=0}^m ||f||_{C^j_b} \| \nabla^{m-j} u\|_{L^2} \; .
$$
We conclude from Proposition \ref{SobolevNormConnection} that $||fu||_{H^k}\lsim ||f||_{C^k_b} \| u \|_{H^k}$, expressing that $f$ acts as a bounded operator on the dense subspace $\Gamma_c(M;E) \subset H^k(M;E)$. Thus, $f$ extends to a bounded operator on $H^k(M;E)$ with the same norm. Moreover, the result is clearly also true for $s=0$. Then interpolation takes care of $0<s<k$ (see Proposition \ref{InterpolationNoBdry} below), and duality of $-k\leq s <0$.
\end{proof}
\section{Sobolev spaces on domains}
In this section we introduce different Sobolev spaces on a domain, and discuss comparison results. We will focus on non-negative degrees, because those are sufficient for our purposes. The material is largely classical, at least for domains in Euclidean space or compact manifolds. See for example \cite{Lions1972} or \cite{Grisvard1985} for domains in Euclidean space under varying smoothness assumptions, or \cite{Taylor1996} for domains in compact manifolds.

Throughout let $M$ be a manifold of bounded geometry, $E\to M$ a vector bundle of bounded geometry, and $\Omega \subseteq M$ be a domain in $M$. Define $\Gamma_{cc}(\Omega;E)$ as the subspace of $\Gamma_c(M;E)$ consisting of those sections whose support is contained in $\Omega$. Moreover, let $\Gamma_c(\olO;E)$ denote the space $\Gamma_c(M;E)/\{ u\in\Gamma_c(M;E) \, | \, u|_{\olO} = 0 \}$ of restrictions of sections in $\Gamma_c(M;E)$ to $\olO$. By definition there is a restriction map $\Gamma_c(M;E) \to \Gamma_c(\olO;E)$. It is injective on $\Gamma_{cc}(\Omega;E)$, so that it makes sense to view the latter as a subspace of $\Gamma_{c}(\olO;E)$, and we shall do so implicitly.
\begin{definition}
Let $s\geq 0$. Let $M,E,\Omega$ as above, and let $B\subseteq M$ be a closed set. Define 
\begin{align*}
H^s(\olO;E)&:=H^s(M;E)/ \{ u \in H^s(M;E) \, | \, \supp(u)\cap \Omega =\emptyset \} \\
H^s_0(\Omega;E) &:= \overline{\Gamma_{cc}(\Omega;E)}^{\|\cdot\|_{H^s(\olO;E)}} \subseteq H^s(\olO;E) \\
\tilde{H}^s(\Omega,E) &:= \overline{\Gamma_{cc}(\Omega;E)}^{\|\cdot\|_{H^s(M;E)}} \subseteq H^s(M;E) \\
H^s_B(M;E) &:= \left\{ u\in H^s(M;E) \, | \, \supp(u)\subseteq B \right\}\subseteq H^s(M;E) \; .
\end{align*}
Moreover, define $H_B^s(\olO;E)\subseteq H^s(\olO;E)$ as image of 
$$
\left\{ u\in H^s(M;E) \, | \, \supp(u)\cap \olO\subseteq B \right\}
$$ 
under the restriction map $H^s(M;E)\to H^s(\olO;E)$.
\end{definition}
We endow $H^s(\olO;E)$ with the quotient norm
$$
\|u\|_{H^s(\olO;E)} := \inf \left\{ \|v\|_{H^s(M;E)} \, | \, v\in H^s(M;E), \, v|_\Omega=u \right\} \; ,
$$
and this norm is used to define the $H^s_0(\Omega;E)$. We stress that the notation $H^s(\olO;E)$ is meant to indicate the this space is obtained by restriction of functions in $H^s(M;E)$ to $\Omega$, resulting in regularity up to the boundary. However, $H^0(\olO;E)$ is given by the space $L^2(\Omega;E)$ instead of $L^2(\olO;E)$.

For a domain $\Omega \subseteq M$ we have three competing notions of what intuitively should be the Sobolev space of sections \emph{vanishing outside $\Omega$}, namely $H_0^s(\Omega;E)$, $\tilde{H}^s(\Omega;E)$ and $H^s_{\olO}(M;E)$. Their relation is fascinatingly subtle and largely depends on both the regularity of $\Omega$ and the degree $s$. 

As low-hanging fruit we may note that $H^s_{\olO}(M;E)$ is a closed subspace of $H^s(M;E)$, and it contains $\Gamma_{cc}(\Omega;E)$. Thus $\tilde{H}^s(\Omega;E)$ is always contained in $H^s_{\olO}(M;E)$. We can also consider the identification map from $\Gamma_{cc}(\Omega;E) \subseteq H^s(M;E)$ to $\Gamma_{cc}(\Omega;E)\subseteq H^s_0(\Omega;E)$. For any $u\in\Gamma_{cc}(\Omega;E)$ it holds that 
\begin{equation}\label{ComparisonSobolevSubspace}
\| u \|_{H^s(\olO;E)} := \inf \left\{ \|v\|_{H^s(M;E)} \, | \, v\in H^s(M;E), \, \supp(u-v)\subseteq M\setminus\Omega \right\} \leq \| u\|_{H^s(M;E)} \; ,
\end{equation}
as can be seen by taking $v=u$. Therefore, the identification map extends to a contractive linear map $\tilde{H}^s(\Omega;E) \to H^s_0(\Omega;E)$. Of course this is nothing but the restriction of the quotient map $H^s(M;E)\to H^s(\olO;E)$ to $\tilde{H}^s(\Omega;E)$. For $k\in\ZZ_{\geq 0}$ the estimate in \eqref{ComparisonSobolevSubspace} becomes an equality: It holds that $\|\nabla^j v\|^2_{L^2}=\|\nabla^j u\|^2_{L^2}+\|\nabla^j (v-u)\|^2_{L^2}$, because the supports of $u$ and $v-u$ don't intersect. In light of Prop. \ref{SobolevNormConnection} it follows that $\| v\|_{H^k(M;E)}$ is minimal for $v=u$, meaning that \eqref{ComparisonSobolevSubspace} is an equality. It follows that the identification map $H^k(M;E)\supseteq\Gamma_{cc}(\Omega;E) \to \Gamma_{cc}(\Omega;E)\subseteq H^k(\olO;E)$ is an isometry, and thus extends to an isometry from $\tilde{H}^k(\Omega;E)$ to $H^k_0(\Omega;E)$. We have thus proven:
\begin{lemma}\label{EqualityH0}
Let $k\in\ZZ_{\geq 0}$. Then $\tilde{H}^k(\Omega;E)$ and $H^k_0(\Omega;E)$ are isometrically isomorphic.
\end{lemma}
In the case that $\Omega$ is a smooth and bounded domain in Euclidean space there is a surprising result going back to Fujiwara \cite{Fujiwara1967} that for $s\geq 0$ the spaces $H^s_0(\Omega)$ and $\tilde{H}^s(\Omega)$ are isomorphic if and only if $s+\frac{1}{2} \notin \ZZ$. Since the transfer of such results to domains in bounded-geometry manifolds is standard fare, and the result will not be needed at any point in this thesis, we permit ourselves to be brief here. Let $\Omega$ be a uniform $C^{m,1}$-domain. Both the property of being in $H^s_0(\Omega;E)$ and of $\tilde{H}^s(\Omega;E)$ can be localized to coordinate patches, and they are preserved by local uniformly bounded $C^{m,1}$-diffeomorphisms. Moreover the original isomorphisms can be made uniform given our uniformity assumptions (see also for example \cite{Lions1972}). In this way one proves the following.
\begin{prop}
Let $\Omega$ be a uniform $C^{m,1}$-domain, and $0\leq s \leq m+1$. Then
$$
H^s_0(\Omega;E) \cong \tilde{H}^s(\Omega;E) \quad \Leftrightarrow \quad s+\frac{1}{2} \notin \ZZ \; .
$$
\end{prop}
We turn to the relation between $\tilde{H}^s(\Omega;E)$ and $H^s_{\olO}(M;E)$. There exist domains $\Omega$ that have boundary of strictly positive measure. Thus, there exist $L^2$-sections supported on $\partial \Omega$, and these cannot be limits of sections with support inside $\Omega$. One can therefore not expect the inclusion $\tilde{H}^s(\Omega;E)\subseteq H^s_{\olO}(M;E)$ to be an equality without imposing at least some regularity on $\Omega$. A too-thick boundary is not the only potential issue: Consider for example $\Omega$ to be the slit disk. Its closure is the disk, so an element of $H^s_{\olO}(M;E)$ can be non-zero on the slit. For sufficiently large $s$, say $s>\frac{1}{2}$, such elements cannot be approximated by elements of $\Gamma_{cc}(\Omega;E)$, which must vanish on the slit. Luckily, being a (uniform) $C^0-$domain suffices.
\begin{prop}[Compare {\cite[Theorem 3.29]{McLean2000}}]\label{EqTildeH}
Let $\Omega$ be a $C^0$-domain. Then $\tilde{H}^s(\Omega;E)=H^s_{\olO}(M;E)$ for all $s\geq 0$. 
\end{prop} 
We need a lemma.
\begin{lemma}\label{ApproxWithCptSupp}
For all $u\in H^s(M;E)$ and all $\varepsilon>0$ there exists $u^\varepsilon\in H^s(M;E)$ such that $\|u-u^\varepsilon\|_{H^s}<\varepsilon$ and the support of $u^\varepsilon$ is compact and contained in the support of $u$.
\end{lemma}
\begin{proof}[Proof of Lemma \ref{ApproxWithCptSupp}]
Let $(U_i,\kappa_i,\psi_i)_{i\in I}$ be an admissible trivialization. Fix $\varepsilon>0$ and $u\in H^s(M;E)$. We show that there exists a finite subset $I^\varepsilon\subseteq I$ such that $\|u-\sum_{i\in I^\varepsilon}\psi_i u\|_{H^s}<\varepsilon$. Pick $v\in \Gamma_c(M;E)$ with $\|u-v\|_{H^s}<\varepsilon$. Then there exists $I^\varepsilon\subseteq I$ finite such that $\sum_{i\in I^\varepsilon}\psi_i v=v$. Then 
$$
\|u-\sum_{i\in I^\varepsilon}\psi_iu\|_{H^s} \leq \| u-v\|_{H^s} + \left\|\sum_{i\in I^\varepsilon}\psi_i(v-u) \right\|_{H^s} \leq \left( 1+\| \sum_{i\in I^\varepsilon} \psi_i\|_{C^k_b}\right) \cdot \|u-v\|_{H^s} \; ,
$$
where $k\in \ZZ_{\geq 0}$ satisfies $k\geq |s|$. The derivatives of the $\psi_i$ are uniformly bounded and the cover is uniformly locally finite. It follows that $\| \sum_{i\in I^\varepsilon} \psi_i\|_{C^k_b}$ is bounded independently of $\varepsilon$ by some universal constant. We conclude that $\|u-\sum_{i\in I^\varepsilon}\psi_iu\|_{H^s}\lsim \varepsilon$. Indeed, this argument proves that the net of operators $\sum_{i\in I'} \psi_i$ for $I'\subset I$ finite converges in norm to the identity map on $H^s(M;E)$.
\end{proof}
\begin{proof}[Proof of Prop. \ref{EqTildeH}]
By Lemma \ref{ApproxWithCptSupp} it suffices to show that any $u\in H^s_{\olO}(M;E)$ with compact support can be approximated by elements of $\Gamma_{cc}(\Omega;E)$. Let $\{ B_{r_i}(x_i) \}$ be a collection of balls covering the boundary as in Def. \ref{DefCKAlpha}. Since $u$ has compact support, there is a finite subset $I'\subseteq I$ such that $u=0$ on $B_{R_i}(x_i)$ if $i\notin I'$. Complete $\{B_{R_i}(x_i)\}_{i\in I'}$ to a finite open cover $\{U_j\}_{j\in J}$ of $\supp(u)$ by normal coordinate charts in such a way that no $U_j$ apart from the $B_{R_i}(x_i)$ intersects the boundary of $\Omega$, and the collection obtained by replacing each $B_{R_i}(x_i)$ in $\{U_j\}_{j\in J}$with the corresponding $B_{r_i}(x_i)$ is still an open cover of $\supp(u)$. Let $\{\psi_j\}_{j\in J}$ be a subordinate partition of unity, i.e. $\psi_j \in C^\infty_{cc}(U_j)$ and $\sum_j \psi_j =1$ on $\supp(u)$, so that if $U_j=B_{R_i}(x_i)$ for some $i\in I'$, then $\psi_j$ is supported in $B_{r_i}(x_i)$. 

Fix $\varepsilon>0$ and let $u_j := \left(\phi_j \circ (\psi_j u)\right)\circ\kappa_j^{-1}$. It is an element of $H^s(\RR^n)$ with compact support in $\kappa_j\left( U_j \cap \olO \right)$. By construction it either holds that $U_j$ is fully contained in $\Omega$, or that $U_j=B_{R_i}(x_i)$ for some $i\in I$. In the former case smooth $u_j$ to a function $u_j^\varepsilon$ with $\supp(u_i^\varepsilon)\subset U_j$ and $||u_j-u_j^\varepsilon||_{H^s}<\varepsilon $.\footnote{This can be achieved via mollifiers. Indeed, if $\varphi_t$, $t>0$, denotes a family of standard mollifiers, and $u\in H^s(\RR^n)$ as well as $\varepsilon, \,\delta>0$ are given, then there exists $t_0$ such that $\| u -\varphi_t* u\|_{H^s}<\varepsilon$ and $\supp(\varphi_t*u)\subset N_\delta(\supp(u))$ for all $0<t\leq t_0$. See for example \cite[Section 5.3 and Appendix C.4]{Evans2010}.} 

In the latter case $u_j$ actually has compact support in $B_{r_i}(0)$. Moreover, by assumption there is a closed cube $K_i$ with $B_{r_i}(0)\subset K_i \subset \kappa_j(U_j)$ such that $\kappa_j\left( \Omega\cap U_j \right)\cap K_i$ is (up to an orthogonal transformation) the area below the graph of a continuous function $\xi_i$. Since $u$ is supported in $\olO$, the function $u_j$ is supported in $\{ x_n \leq \xi_i(x_1,\cdots,x_{n-1})\}$. For $\delta>0$ write $u^\delta_j(x_1,\cdots,x_n):=u_j(x_1,\cdots, x_{n-1},x_n+\delta)$. It holds that $u^\delta_j \in H^s(\RR^n)$ and $\|u_j-u_j\|_{H^s} \xrightarrow{\delta\to 0} 0$. Moreover, $u^\delta_j$ is supported within the set $\{ x_n \leq \xi_i(x_1,\cdots,x_{n-1}) - \delta\}$. Since $\xi_i$ is continuous on the compact set $[-c_i,c_i]^{n-1}$, it is even uniformly continuous. This implies that the set $\{ x_n \leq \xi_i(x_1,\cdots,x_{n-1}) - \delta\}$ has strictly positive distance from $\{ x_n \leq \xi_i(x_1,\cdots,x_{n-1})\}$. Thus the support of $u^\delta_j$ has strictly positive distance from $\kappa_j(\partial \Omega \cap U_j)\cap K_i$. Moreover, $u^\delta_j$ is still supported within $K_i$ for $\delta$ sufficiently small. Then $u^\delta_j$ may be smoothed to a smooth function $u^\varepsilon_j$ which is still compactly supported in $\{ x_n < \xi_i(x_1,\cdots,x_{n-1}) - \delta\}$. Choosing the shift $\delta$ and the amount of smoothing sufficiently small, this can be done so that $\| u-u^\varepsilon_j \|_{H^s}<\varepsilon$.

Finally, define $u^\varepsilon := \sum_{j\in J} \phi_j^{-1} \circ u^\varepsilon_j \circ \kappa_j$, which is well-defined since each $u^\varepsilon_j$ has support within $\kappa_j(U_j)$. By construction $u^\varepsilon$ is a smooth section with compact support inside $\Omega$, and it holds that 
$$
\| u -u^\varepsilon\|_{H^s} \lsim \sum_{j\in J} \left\| u_j - u_j^\varepsilon \right\| \lsim \varepsilon \; .
$$
This concludes the proof.
\end{proof}
We shall also need the fact that multiplication by $C^k_b$-functions preserves the Sobolev spaces on domains. To that end define
$$
C^k_b(\olO) := C^k_b(M)/\big\{ f\in C^k_b(M) \, | \, f=0  \text{ on }\Omega \big\} \; ,
$$
endowed with the quotient norm.
\begin{lemma}\label{MultSobolevDomain}
Let $k\in\NN$, and $f\in C^k_b(\olO)$. Then,
$$
H^{s}(\olO;E) \longrightarrow H^{s}(\olO;E) \; , \; u \longmapsto fu
$$ 
is a bounded map for all $0\leq s\leq k$ with norm bounded by $\|f\|_{C^k_b(\olO)}$.
\end{lemma}
\begin{proof}
Let $u\in H^s(\olO;E)$. Let $\tilde{u}\in H^s(M;E)$ satisfy $\tilde{u}|_\Omega=u$, and $\tilde{f} \in C^k_b(M)$ satisfy $\tilde{f}|_\Omega=f$. Then Lemma \ref{MultWInftyNoBdry} gives 
$$
\|\tilde{f}\tilde{u}\|_{H^s(M)}\leq \|\tilde{f}\|_{C^k_b}\cdot\|\tilde{u}\|_{H^s(M)} \; .
$$
Taking infima yields
$$
\| fu \|_{H^s(\olO)} \leq \inf_{\tilde{u},\tilde{f}} \|\tilde{f}\tilde{u}\|_{H^s(M)}\leq \big( \inf_{\tilde{f}}\|\tilde{f}\|_{C^k_b(M)}\big)\cdot\big( \inf_{\tilde{u}}\|\tilde{u}\|_{H^s(M)} \big) = \|f\|_{C^k_b(\olO)}\cdot \|u\|_{H^s(\olO)} \; ,
$$
proving the claim.
\end{proof}
Let us also give two lemmas that -while not immediately related to Sobolev spaces on domains- are important technical tools for our discussion of uniform K-homology on domains in Chapter 6. Again, let $M$ be a manifold of bounded geometry, and $\Omega\subseteq M$ a domain. The manifold $M$ is a metric space under the path-length metric, and we endow $\Omega \subseteq M$ with the subspace metric. Note that the subspace metric generally differs from the path-length metric on $\Omega$. Recall that $\LLip_R(M)$ denotes the set of all functions $f\in C_0(\Omega)$ with $\|f\|_\infty\leq 1$ that are $L$-Lipschitz and have compact support in $\Omega$ with diameter at most $R$. Since any $f\in\LLip_R(\Omega)$ has compact support in $\Omega$, it can be extended by zero to yield a function on $M$. Define
\begin{align*}
\LLip_{R,\delta} (\Omega) &:= \left\{ f\in \LLip_R(\Omega) \, | \, d\left( \supp(f),\dO\right)\geq\delta \right\} \\
\LLip_{R,m}(\Omega) &:= \left\{ f\in \LLip_R(\Omega)\cap C^\infty(M) \, | \, \|f\|_{C^m_b(\olO)}\leq L \right\} \\
\LLip_{R,\delta,m}(\Omega) &:= \LLip_{R,\delta}(\Omega)\cap\LLip_{R,m} (\Omega) \; .
\end{align*} 
\begin{lemma}\label{LLipStratified}
Let $\sigma$ be a uniformly continuous map from $\{ f\in C_0(M) \, | \, ||f||_\infty\leq 1 \}$ to $\fBH$ for some Hilbert space $H$. The following are equivalent:
\begin{itemize}
\item[(i)] The family $\sigma\left( \LLip_R(M) \right)$ is uniformly approximable for all $L,R\geq 0$.
\item[(ii)] The family $\sigma\left( \LLip_{R,\delta}(M) \right)$ is uniformly approximable for all $L,R\geq 0$ and $\delta>0$.
\item[(iii)] The family $\sigma\left( \LLip_{R,m}(M) \right)$ is uniformly approximable for all $L,R\geq 0$ and $m\in\ZZ_{\geq 1}$.
\item[(iv)] The family $\sigma\left( \LLip_{R,\delta,m}(M) \right)$ is uniformly approximable for all $L,R\geq 0$, $\delta>0$ and $m\in\ZZ_{\geq 1}$.
\end{itemize}
\end{lemma}
Note that the properties (iii) and (iv) are independent of the choice of admissible trivialization. We also note that the equivalence of (i) and (ii) holds for arbitrary metric spaces, and was in fact already used in the proof of the excision isomorphism.
\begin{proof}
We have inclusions $\LLip_{R,m}\subseteq \LLip_R$, $\LLip_{R,\delta}\subseteq \LLip_R$, $\LLip_{R,\delta,m}\subseteq \LLip_{R,\delta}$ and $\LLip_{R,\delta,m}\subseteq \LLip_{R,m}$, which yield the implications (i)$\Rightarrow$(ii), (i)$\Rightarrow$(iii), (ii)$\Rightarrow$(iv) and (iii)$\Rightarrow$(iv) respectively. Thus the proof is complete once we show that (ii)$\Rightarrow$(i) and that (iv)$\Rightarrow$(ii).

Let us take care of the implication (ii)$\Rightarrow$(i) first. Fix $L,R\geq 0$ and $\varepsilon>0$. Since $\sigma$ is uniformly continuous, there exists $\varepsilon'>0$ such that if $\|f-f'\|<\varepsilon'$, then $\|\sigma(f)-\sigma(f')\|<\frac{\varepsilon}{2}$. Let $f\in \LLip_R(\Omega)$. Since $f$ vanishes on $\partial\Omega$, it holds that $|f(x)|\leq L\cdot d(x,\dO)$. Set $\delta:= \frac{\varepsilon'}{L}$. Let $\psi: \, M\to [0,1]$ be a Lipschitz function such that $\psi(x)=0$ if $d(x,\dO)\leq \frac{\delta}{2}$ and $\psi(x)=1$ if $d(x,M\setminus\Omega)\geq \delta$. The Lipschitz constant $L_\delta$ of $\psi$ can be arranged to depend only on $\delta$\footnote{Take $\psi(x)$ to be the maximum of $0$ and $1-\frac{2}{\delta}\cdot d(x,\Omega\setminus N_{\delta}(M\setminus\Omega))$. This function has the desired properties and Lipschtitz constant $\frac{2}{\delta}$.}. Then, $\psi f \in (L_\delta L)\text{-}\mathrm{Lip}_R(\Omega)$ and $\| \psi f - f\|_\infty \leq \varepsilon'$. Assuming (ii) to hold, there exists $N=N(L,R,\delta)$ so that $\sigma(\psi f)$ is within $\frac{\varepsilon}{2}$ of a rank-$N$ operator. Moreover, since $\| \psi f - f\|_\infty \leq \varepsilon'$, it holds that $\|\sigma(\psi f)-\sigma(f)\|<\frac{\varepsilon}{2}$. Thus, $\sigma(f)$ is within $\varepsilon$ of a rank-$N$ operator. Since $\delta$ depends only on $L$ and $\varepsilon'$, and the latter in turn only depends on $\sigma$ and $\varepsilon$, this establishes (i).

To prove (iv)$\Rightarrow$(ii) it is sufficient to prove that for all $L,R\geq 0$ and all $\delta,\varepsilon>0$ there are $L',R'\geq 0$, $\delta'>0$ and $m\in\ZZ_{\geq 1}$ such that for every $f\in \LLip_{R,\delta}(M)$ there is $f'\in L'\text{-}\mathrm{Lip}_{R',\delta',m}(M)$ such that $\|f-f'\|_\infty<\varepsilon$, analogously to the previous implication an appeal to uniform continuity of $\sigma$ then completes the argument. This latter claim is proved by smoothing via mollifiers while controlling the norms of derivatives. Consider an admissible trivialization $(U_i=B_{2r}(x_i),\kappa_i,\psi_i)$ of $M$ such that $\supp(\psi_i)\subseteq B_{r}(x_i)$. Fix a non-negative smooth function $\varphi \in C^\infty_{c}(\RR^n)$ with $||\varphi||_{L^1}=1$ and $\supp(\varphi)\subseteq B_1(0)$. Convolution by $\varphi_t(x):= t^{-n} \varphi(t^{-1}x)$ yields a family of mollifiers. In particular, if $g\in C_c(\RR^n)$ is continuous and compactly supported, then $\varphi_t * g$ is smooth, supported in $N_{2t}(\supp(g))$, and $||g-\varphi_t * g||_\infty\to 0$ for $t\to 0$. In fact, it holds that\footnote{See for example \cite[Appendix C.4]{Evans2010}.}
\begin{equation} \label{ConvMollifier}
\left| g(x)-(\varphi_t*g)(x) \right| \leq t^{-n} \int_{B_t(x)} \left| g(x)-g(y) \right| \, dy \xlongrightarrow{t\to 0} 0 \; .
\end{equation} 
Moreover it is easily checked that 
\begin{equation} \label{DerivativeBdMollifier}
\left| D^\alpha \left( \varphi_t * g \right) \right| \leq t^{-|\alpha|} \cdot ||g||_\infty \cdot ||D^\alpha \varphi||_{L^1}
\end{equation}
for any multi-index $\alpha$. Now, set $f_t:= \sum_i \left( \varphi_t * \left( (\psi_i f)\circ\kappa_i^{-1} \right)\right) \circ \kappa_i$, where $t$ is sufficiently small so that $\supp\left(\varphi_t * \left( (\psi_i f)\circ\kappa_i^{-1} \right) \right) \subseteq B_{r+2t}(0)$ is contained in $\kappa_i(U_i)=B_{2r}(0)$. The function $f_t$ is smooth and compactly supported in $M$, and $||f-f_t||\to 0$ for $t\to 0$. By the properties of an admissible trivialization, there exists $N$, dependent only on $R$, such that at most $N$ summands of this sum are non-zero. Thus, $f_t$ is supported in the union of at most $N$ of the $U_i$, a set with diameter $\leq N\cdot r =:R'$. Moreover, there exists $t=t(\delta)$ such that $\supp(f_t) \subseteq N_{\frac{\delta}{2}}(\supp(f))$. Lastly, the estimate \eqref{DerivativeBdMollifier} yields
$$
\left| D^\alpha \left( f_t\circ\kappa_i^{-1} \right) \right| \leq \frac{||D^\alpha \varphi||}{t^{|\alpha|}} \; ,
$$
implying that $\|f_t\|_{C^m_b}$ is bounded for any $m$ by a constant dependent only on $t$ and $m$ but not $f$. It follows that for any $m$ there is $L'\geq 0$, depending on $L$, $m$ and $t$ but not the specific $f$, such that $f_t\in L'\text{-}\mathrm{Lip}_{R',\frac{\delta}{2},m}(M)$ for $t\leq t(\delta)$ sufficiently small. Moreover, inspecting the estimate \eqref{ConvMollifier} we see that $||f-f_t|| \lsim L\cdot t$. Thus, $t$ can be chosen so that $||f-f_t||\leq \varepsilon$ for all $f\in\LLip_R(M)$. This finishes the proof of the implication (iv)$\Rightarrow$(ii).
\end{proof}
The second lemma provides uniformly bounded bump functions, whose existence will be exploited for example in the proof of Theorem \ref{QuantComp}, and in Subsection \ref{SubsecFormalElliptic}.
\begin{lemma}\label{BdBump}
Let $m\in \ZZ_{\geq 0}$. Let $K\subset \Omega$ be a compact subset, and set $R:=\diam(K)$ and $\delta:=d(K,\partial\Omega)$. Then, for any $0<\delta'<\delta$ there exists a function $\psi\in C^\infty_{cc}(\Omega)$ with $\psi=1$ on $K$, $\psi=0$ outside $N_{\delta'}(K)$, and such that $\|\psi\|_{C^m_b(\olO)}$ is bounded by a constant dependent only on $R$ and $\delta'$. 
\end{lemma}
\begin{proof}
One can find a Lipschitz function $\psi' \in C_{cc}(\Omega)$ with $\psi'=1$ on $N_{\delta'/3}(K)$, $\psi'=0$ outside $N_{2\delta'/3}(K)$, such that the Lipschitz constant of $\psi'$ depends only on $\delta'$. Define $\psi \in C^\infty_{cc}(\Omega)$ by mollifying $\psi'$. Taking the mollification parameter $t$ sufficiently small it holds that $\psi=1$ on $K$ because $\psi'=1$ on $N_{\delta'/3}(K)$, and that $\psi=0$ outside $N_{\delta'}(K)$ because $\psi'=0$ outside $N_{2\delta'/3}(K)$. The estimate \eqref{DerivativeBdMollifier} yields the desired norm estimate for $\psi$.
\end{proof}
\section{Extensions and traces}
This section deals with the problems of extending Sobolev functions from a domain to the entire space, and of restricting Sobolev functions on a domain to the boundary.

There are two ways to approach the extension problem: The first is abstract, making use of the fact that surjections of Hilbert spaces split. The second is to explicitly write down an extension formula on smooth functions, and verify that it extends to Sobolev spaces. The second approach is preferable in most ways; it is concrete, the same extension formula works for many Sobolev degrees, and it works for Sobolev spaces based on $L^p$ for $p\neq 2$. The only advantage of the first approach is that it works in full generality, while the second requires uniform regularity up to some degree.

We treat the abstract approach first. Let $T: \, H\to H'$ be a continuous surjective map between Hilbert spaces. When restricted to $(\ker T)^\perp$ the map $T$ becomes bijective, and by the open mapping theorem its inverse is again continuous. Thus there exists a unique map $S: \, H'\to H$ such that $TS=\id_{H'}$ and $\im\, S =(\ker T)^\perp$. 

Now, consider the restriction map $r: \, H^s(M;E) \to H^s(\olO;E)$. An extension map $e: \, H^s(\olO;E)\to H^s(M;E)$ is the same as a right-inverse to $r$, meaning that $r(e(u))=e(u)|_\Omega=u$. The restriction map is surjective and continuous by definition. Thus, the abstract argument provides us with the following:
\begin{prop}\label{ExtensionAbstract}
Let $\Omega$ be a domain in $M$, and let $s\geq 0$. Then, there is an extension map $e: \, H^s(\olO;E)\to H^s(M;E)$ uniquely characterized by the property that $\im\, e=(\ker r)^\perp$. 
\end{prop}
To reiterate: This construction heavily relies on the fact that we work with Hilbert spaces (surjections between Banach spaces need not split), and the resulting map is highly dependent on the degree $s$. Thus, whenever possible it is preferable to use the approach described now.

Fix $k\in \ZZ_{\geq 0}$. We begin by considering the extension problem on half-space. Let $u\in C^\infty(\RR^n_+)$ be a (potentially vector-valued) smooth function on $\RR^n_+=\{x_n\geq 0\}$. The function $u$ can be extended across the plane $x_n=0$ by reflection, and via superposition of multiple such reflections the first $k$ derivatives can be made continuous. Indeed, let $e_k(u)$ be defined by $e_k(u)=u$ on $\RR^n_+$, and for $x_n<0$ by
$$
e_k(u)(x_1,\cdots,x_{n-1},x_n) := \sum_{j=1}^{k+1} a_j \cdot u(x_1,\cdots,x_{n-1},-jx_{n}) \; ,
$$
where the coefficients $a_1,\cdots,a_{k+1}$ are determined uniquely by the requirement that $e(u)$ be in $C^k(\RR^n)$. Then, the assignment $e_k$ extends to a continuous map $e_k: \, H^s(\RR^n_+)\to H^s(\RR^n)$ for all $0\leq s\leq k$. Note that the coefficients depend on the degree $k$, though a modification can be made that works for $k=\infty$ as well \cite{Seeley1964}. Note also that $e_k$ can be modified via cut-off functions so that the values of $e_k(u)$ on $\{x_n<0\}$ depend only on the values of $u$ in the strip $\{0\leq x_n \leq \delta\}$ for a fixed $\delta>0$, and such that $e_k(u)$ is supported on $\{-\delta'\leq x_n\}$ for some $\delta'>0$.

Now, consider instead the space $\olO_\xi:=\{ x_n \geq \xi(x_1,\cdots,x_{n-1}) \}\subset \RR^n$ for some bounded $C^{m,1}$-function $\xi: \, \RR^{n-1}\to \RR$ with finite $C^{m,1}$-norm. Then, the boundary of $\olO_\xi$ may be straightened out by a $C^{m,1}$-diffeomorphism $\chi_\xi: \, \RR^n \xrightarrow{\sim} \RR^n$, $\chi_\xi(\olO_\xi)=\RR^n_+$, with bounded $C^{m,1}$-norm. This $C^{m,1}$-diffeomorphism induces a bounded isomorphism
$$
H^s(\olO_\xi) \longrightarrow H^s(\RR_+^n) \; , ; u \longmapsto u \circ \chi_\xi^{-1} \; ,
$$
for all $0\leq s\leq m+1$ (see for example \cite[Theorem 3.23]{McLean2000}), and as such gives rise to extension maps
$$
e_k^\xi: \; H^s(\olO_\xi) \longrightarrow H^s(\RR^n) \; , \; e_k^\xi(u) =  e_k(u\circ \chi_\xi^{-1})\circ \chi_\xi
$$
for these $s$.

If $\Omega$ is a uniform $C^{m,1}$-domain in the bounded-geometry manifold $M$, then these consideration in Euclidean space give rise to extension maps $H^s(\olO;E)\to H^s(M)$, $0\leq s\leq m+1$. Indeed, locally the boundary of $\olO$ looks like a domain $\olO_{\xi_i}$ as considered before. Thus extension maps exist locally, and due to the various uniformness conditions the norms of the local extension maps are uniformly bounded. Using a partition of unity they may therefore be glued together to a global bounded extension map. Since the formula of the original extension map on half-space did not depend on the specific $s$, neither does the extension map for $\Omega$. Thus we obtain the following.
\begin{prop}\label{ExtensionReflection}
Let $\Omega$ be a uniform $C^{m,1}$-domain in the manifold $M$ of bounded geometry. Then, for all $0\leq s\leq m+1$ there exists a bounded extension map
$$
e_s: \; H^s(\olO;E) \longrightarrow H^s(M;E)
$$ 
such that if $s\leq s'$, then $e_{s'}(u)=e_s(u)$ for all $u\in H^{s'}(\olO;E)$. If $\Omega$ is a bounded-geometry domain, then extension maps with this property exist for all $s\geq 0$. 
\end{prop} 
The usefulness of the property that $e_{s'}(u)=e_s(u)$ for all $u\in H^{s'}(\olO;E)$ will become immediately clear when we discuss interpolation for the spaces $H^s(\olO;E)$ below. To reiterate another point, though it is less relevant for our purposes: This construction works equally well for Sobolev spaces based on $L^p$, $p\neq 2$, and is therefore more general.

Let us briefly discuss trace maps on Sobolev spaces of domains. Apart from a few examples traces will play no role in this thesis, but a treatment of Sobolev spaces on domains would be incomplete without mentioning them. Again, half-space serves as the universal model case. There one shows (see for example \cite[Proposition 4.1.6]{Taylor1996}) that for all $s>\frac{1}{2}$ there is a surjective map
$$
H^s(\RR^n_+) \longrightarrow H^{s-\frac{1}{2}}(\RR^{n-1}) \; , \; u\longmapsto u|_{\RR^{n-1}} \; .
$$
Then, if $\olO$ is a manifold with boundary and bounded geometry, one constructs a trace to the boundary by locally straightening out via uniformly bounded diffeomorphisms, taking traces to half-space locally, and gluing together via a partition of unity. The details are carried out in greater generality in \cite{Grosse2013}. Thus, one obtains the following.
\begin{thm}[{\cite[Theorem 47]{Grosse2013}}]
Let $\olO$ be a manifold with boundary and bounded geometry. Then, for all $s>\frac{1}{2}$ there is a surjective trace map
$$
H^s(\olO;E) \longrightarrow H^{s-\frac{1}{2}}(\partial\Omega;E|_{\partial\Omega}) \; , \; u\longmapsto u|_{\partial \Omega} \; .
$$
\end{thm}
In this discussion of traces we have focused on manifolds with boundary and bounded geometry, in particular to domains with smooth boundaries. This is because in this thesis we do not consider Sobolev spaces on non-smooth manifolds, so that we have no way to talk about Sobolev space on the boundary of a non-smooth domain. This is no substantial issue though. If $\Omega$ is a uniform $C^{m,1}$-domain, then its boundary is a $C^{m,1}$-manifold that already comes with a canonical set of local $C^{m,1}$-trivializations. Thus we may define $H^s(\partial\Omega)$, $0\leq s\leq m+1$, in analogy to how we defined Sobolev spaces on manifolds of bounded geometry. Using the same procedure of local straightening via uniformly bounded $C^{m,1}$-diffeomorphisms, tracing to the straightened boundary, then gluing together via a partition of unity, one obtains surjective trace maps
$$
H^s(\olO) \longrightarrow H^{s-\frac{1}{2}}(\partial \Omega)
$$
for $\frac{1}{2}<s\leq m+1$ whenever $\Omega$ is a uniform $C^{m,1}$-domain. See for example \cite{McLean2000}, particularly Theorem 3.37 therein, where this is discussed for $C^{m,1}$-domains in Euclidean space.
\section{Interpolation} \label{SecSobolevInterpolation}
At this point we will discuss interpolation of the Sobolev spaces on domains in manifolds of bounded geometry. The complex interpolation method was introduced by Calderón in \cite{Calderon1964}, and has become a standard tool in the study of Sobolev spaces and related function spaces. Knowledge about interpolation for Sobolev spaces on domains will be necessary to make exact uniform summability statements later in the context of uniform Fredholm modules. The discussion is classical in the context of Euclidean space or compact manifolds. An exhaustive treatment of interpolation for domains in Euclidean space can be found in \cite{Triebel1995}.

Let us recall the the basics. Let $E,F$ be Banach spaces such that there is a bounded injection $F\hookrightarrow E$. Denote by $\mathcal{H}(E,F)$ the space of bounded continuous functions $f$ from the complex strip $\{0\leq \mathrm{Re}(z)\leq 1\}$ to $E$ such that $f$ is analytic on the interior of $\{0< \mathrm{Re}(z)< 1\}$, and $f(1+i\RR)\subseteq F$ with $\sup_t ||f(1+it)||_F <\infty$. Equipped with the norm $||f||_{\mathcal{H}(E,F)}:= \max\{\sup_t ||f(it)||_E, \, \sup_t ||f(1+it)||_F\}$ the set $\mathcal{H}(E,F)$ becomes a Banach space. Then, for $\theta\in [0,1]$ define the \emph{(complex) interpolation space} as
$$
[E,F]_\theta:=\{ f(\theta) \, | \, f\in \mathcal{H}(E,F)\} \; ,
$$
equipped with the norm $||u||_\theta:=\inf_{f(\theta)=u} ||f||_{\mathcal{H}(E,F)}$. The interpolation space $[E,F]_\theta$ is a Banach space, and isomorphic to the quotient of $\mathcal{H}(E,F)$ by the subspace $\{ f(\theta)=0\}$. It is functorial in $E$ and $F$ in the following way. Let $E,F,E',F'$ be Banach spaces such that there are bounded injections $F\hookrightarrow E$ and $F'\hookrightarrow E'$. Let $T: \, E\to E'$ be a bounded map which restricts to a bounded map $T: \, F\to F'$. Let $\theta\in [0,1]$. Then, $T$ induces a bounded map 
$$
T_\theta : \, [E,F]_\theta \to [E',F']_\theta
$$
with norm bounded by $\|T\|_{E\to E'}^{1-\theta}\|T\|_{F\to F'}^\theta$. In particular, because composition with coordinate maps or vector bundle trivializations, and multiplication by smooth functions with bounded derivatives yield bounded maps on Sobolev spaces, they also yield maps on interpolation spaces between Sobolev spaces as well. We will make liberal use of this in the following.
\begin{prop}\label{InterpolationNoBdry}
Let $M$ be a manifold of bounded geometry, and $E\to M$ a vector bundle of bounded geometry. Then, for any $s_0,s_1\in\RR$ and $\theta\in [0,1]$ it holds that
$$
\left[ H^{s_0}(M;E), H^{s_1}(M;E) \right] = H^{(1-\theta)s_0 + \theta s_1}(M;E) \; .
$$
with equivalent norms.
\end{prop}
The proof, as usual, is a localization argument combined with an appeal to knowledge that the result holds for Euclidean space. But since this is the first appearance of interpolation spaces in this thesis we provide the details. 
\begin{proof}
For the sake of simplicity we only treat the case where $E$ is the trivial bundle, the case of general vector bundles is entirely analogous. Also, we restrict to $s_0=0$ and $s_1=k\in\NN$.

Choose an admissible trivialization $(U_i=B_R(x_i),\kappa_i,\psi_i)$ of $M$, with $\supp(\psi_i)\subset B_r(x_i)$ for $0<r<R<\mathrm{inj}_M$. Let $u\in L^2(M)$. Then $u\in H^{k\theta}(M)$ if and only if for all $i$ the function $u_i:=(\psi_i u)\circ \kappa_i^{-1}$ is in $H^{k\theta}(\RR^n)$ and it holds that $\sum_i \|u_i \|_{H^s(\RR^n)} <\infty $. Because $H^{k\theta}(\RR^n)=[L^2(\RR^n),H^k(\RR^n)]_\theta$ (see for example \cite[Section 4.2]{Taylor1996}) with equivalent norms this is equivalent to the condition that $u_i\in [L^2(\RR^n),H^k(\RR^n)]_\theta$ for all $i$ with $\sum_i \|u_i \|_{\theta} <\infty$. Thus it suffices to prove that this last condition is equivalent to $u\in [L^2(M),H^k(M)]_\theta$.

Thus suppose $u\in L^2(M)$ is such that all $u_i$ are in $[L^2(\RR^n),H^k(\RR^n)]_\theta$ with $\sum_i \|u_i \|_{\theta} <\infty$. Choose $g_i' \in \cH(L^2(\RR^n),H^k(\RR^n))$ with $g_i'(\theta)=u_i$ such that $\sum_i \|g_i'\|_{\cH}<\infty$. Then, pick a smooth function $\rho\in C^\infty_c(\RR^n)$ with $\supp(\rho)\subseteq B_R(0)$ and $\rho(x)=1$ for $|x|\leq r$. Then $\rho u_i = u_i$, and it holds that $\supp(\rho)$ is supported inside domain of the coordinate maps $\kappa_i$. Thus we may define $g_i := (\rho g_i')\circ \kappa_i$. These are elements of $\cH(L^2(M),H^k(M))$, and it still holds that
$$
\sum_i \|g_i\|_\cH \lsim \sum_i \|g_i'\|_{\cH}<\infty \; .
$$
Thus $\sum_i g_i$ is an absolutely convergent series in $\cH(L^2(M),H^k(M))$, hence converges to an element $g\in \cH(L^2(M),H^k(M))$. This element then satisfies
$$
g(\theta) = \sum_i \left( (\rho g_i')\circ \kappa_i \right)( \theta) = \sum_i (\rho u_i)\circ\kappa_i = \sum_i u_i \circ \kappa_i = u \; .
$$
Thus $u\in [L^2(M),H^k(M)]_\theta$.

Suppose conversely that $u\in [L^2(M),H^k(M)]_\theta$. Then $u_i \in [L^2(\RR^n),H^k(\RR^n)]_\theta$ as well, so that it suffices to check the finiteness condition. Due to the uniform local finiteness of the admissible trivialization there exists a decomposition $I=I_1\cup \cdots I_N$ such that the supports of $\psi_i$, $i\in I_j$, have uniformly positive distance. Since 
$$
\sum_{i\in I} \|u_i\|_\theta = \sum_{j=1}^N \left( \sum_{i\in I_j} \|u_i\|_\theta \right)
$$
it suffices to prove that each $\sum_{i\in I_j} \|u_i\|_\theta$ is finite. To that end let $g\in \cH(L^2(M),H^k(M))$ have $g(\theta)=u$. Let $r\in \{0,1\}$, and let $H_r$ be $L^2(M)$ if $r=0$, and $H^k(M)$ if $r=1$. In either case $H_r$ is a Hilbert space, and we can consider the orthogonal projection $P_i^{r}$ on $H_r$ to the subspace of functions with support in $\supp(\psi_i)$. Because the $P_i^{r}$ have mutually orthogonal ranges within each $I_j$, we have
$$
\sum_{i\in I_j} \|\psi_i g(r+it) \|_{H_r} = \sum_{i\in I_j} \|\psi_i (P_i^r g(r+it)) \|_{H_r} \lsim \sum_{i\in I_j} \|P_i^r g(r+it) \|_{H_r} \leq \|g(r+it) \|_{H_q} 
$$
for any $t\in \RR$. Taking the infimum over $g$ we conclude that
$$
\sum_{i\in I_j} \|u_i\|_\theta \lsim \sum_{i\in I_j} \|\psi_i u\|_\theta \lsim \|u\|_\theta \; ,
$$
which yields the reverse implication.
\end{proof}
An analogous argument could also be given for uniform $C^{m,1}$-domains. Instead we simply appeal to the extension map. 
\begin{prop}\label{InterpolationUpToBdry}
Let $\olO$ be a uniform $C^{m,1}$-domain. Then 
$$
[L^2(\olO;E), H^{m+1}(\olO;E)]_\theta = H^{(m+1)\theta}(\olO;E) 
$$
for all $\theta\in [0,1]$. In particular, if $\olO$ is a manifold with boundary and bounded geometry, then
$$
[L^2(\olO;E), H^{k}(\olO;E)]_\theta = H^{k\theta}(\olO;E) 
$$
for all $k\in\NN$ and $\theta\in [0,1]$.
\end{prop}
\begin{proof}
Consider the extension map of Proposition \ref{ExtensionReflection}. It provides a map $e: \, L^2(\Omega;E)\to L^2(M;E)$ that restricts to a bounded map $e: \, H^{m+1}(\olO;E) \to H^{m+1}(M;E)$. Moreover, it is right-inverse to the restriction map $r: \, L^2(M;E)\to L^2(\Omega;E)$, which also restricts to a bounded map $r: \, H^{m+1}(M;E)\to H^{m+1}(\olO;E)$. Interpolation thus provides maps 
\begin{align*}
e_\theta:& \; \left[ L^2(\Omega;E),H^{m+1}(\olO;E) \right]_\theta \longrightarrow \left[ L^2(M;E),H^{m+1}(M;E) \right]_\theta \\
r_\theta:& \; \left[ L^2(M;E),H^{m+1}(M;E) \right]_\theta \longrightarrow \left[ L^2(\Omega;E),H^{m+1}(\olO;E) \right]_\theta
\end{align*}
such that $r_\theta\circ e_\theta=\id$. In particular, $r_\theta$ is surjective. Now, we know that 
$$
\left[ L^2(M;E),H^{m+1}(M;E) \right]_\theta=H^{(m+1)\theta}(M;E) \; .
$$ 
Moreover, by definition the composition
$$
H^{(m+1)\theta}(M;E) \xlongrightarrow{r_\theta} \left[ L^2(\Omega;E),H^{m+1}(\olO;E) \right]_\theta \longrightarrow L^2(\Omega)
$$
is the restriction map. The latter map is an inclusion by definition, implying that
$$
\ker(r_\theta) = \left\{ u\in H^{(m+1)\theta}(M;E) \, | \, u|_\Omega=0 \right\} \; .
$$
We conclude that
$$
\left[ L^2(\Omega;E),H^{m+1}(\olO;E) \right]_\theta = \im(r_\theta) \cong H^{(m+1)\theta}(M;E)/\ker(r_\theta) = H^{(m+1)\theta}(\olO;E) \; ,
$$
proving the claim.
\end{proof}
In the above proof the fact that the extension map $e: \, L^2(\Omega;E)\to L^2(M;E)$ restricts to a bounded map $e: \,  H^{m+1}(\olO;E) \to H^{m+1}(M;E)$ -which is satisfied by the extension map of Proposition \ref{ExtensionReflection} but not that of Proposition \ref{ExtensionAbstract}- is crucial to obtain the map $e_\theta$. Without it we could not conclude that $r_\theta$ is surjective. Thus we would still get an inclusion of $H^{(m+1)\theta}(\olO;E)$ into $\left[ L^2(\Omega;E),H^{m+1}(\olO;E) \right]_\theta$, but could not conclude the reverse inclusion. There is no reason, therefore, to expect the conclusion of Proposition \ref{InterpolationUpToBdry} to hold without sufficient regularity assumptions.

We come to the interpolation properties of the various Sobolev spaces of sections supported inside a given domain. We begin with the general statement without any regularity assumptions. 
\begin{prop}\label{InterInclusion}
Let $\Omega$ be an arbitrary domain. Let $s_0,s_1\geq 0$ and $\theta\in [0,1]$, and set $s_\theta:=(1-\theta)s_0+\theta s_1$. Then, there is a chain of inclusions
$$
\tilde{H}^{s_\theta}(\Omega;E)\subseteq \left[ \tilde{H}^{s_0}(\Omega;E),\tilde{H}^{s_1}(\Omega;E) \right]_\theta \subseteq \left[ H^{s_0}_{\olO}(M;E),H^{s_1}_{\olO}(M;E) \right]_\theta \subseteq H^{s_\theta}_{\olO}(M;E) \; .
$$ 
\end{prop}
\begin{proof}
For the first inclusion note that $\left[ \tilde{H}^{s_0}(\Omega;E),\tilde{H}^{s_1}(\Omega;E) \right]_\theta$ is a closed subspace of $\left[ H^{s_0}(M;E),H^{s_1}(M;E) \right]_\theta=H^{s_\theta}(M;E)$, so that it suffices to show that $\Gamma_{cc}(\Omega;E)\subseteq \left[ \tilde{H}^{s_0}(\Omega;E),\tilde{H}^{s_1}(\Omega;E) \right]_\theta$. To that end let $u\in\Gamma_{cc}(\Omega;E)$ have compact support $K\subset \Omega$. Let $g\in\mathcal{H}(H^{s_0}(M;E),H^{s_1}(M;E))$ have $f(\theta)=u$. Pick a smooth function $\psi\in C^\infty_{cc}(\Omega)$ with $\psi u=u$. Then, $(\psi \cdot g)(\theta')$ has compact support in $\Omega$ for all $\theta'\in [0,1]$, and $(\psi\cdot g)(\theta)= u$. Thus $\psi g \in \mathcal{H}\left( \tilde{H}^{s_0}(\Omega;E),\tilde{H}^{s_1}(\Omega;E) \right)$, so that $u\in \left[ \tilde{H}^{s_0}(\Omega;E),\tilde{H}^{s_1}(\Omega;E) \right]_\theta$.

The second inclusion follows directly from the fact that $\tilde{H}^s(\Omega;E)$ is a subspace of $H^s_{\olO}(M;E)$ for any $s$. 

For the third inclusion note simply that 
$$
\left[ H^{s_0}_{\olO}(M;E),H^{s_1}_{\olO}(M;E) \right]_\theta\subseteq \left[ H^{s_0}(M;E),H^{s_1}(M;E) \right]_\theta=H^{s_\theta}(M;E) \; ,
$$
and that any element of $\left[ H^{s_0}_{\olO}(M;E),H^{s_1}_{\olO}(M;E) \right]_\theta$ is in particular an element of $H^{s_0}_{\olO}(M;E)$, so has support in $\olO$.
\end{proof}
In general one has no cause to expect any of these inclusions to be an equality. However, if $\Omega$ is sufficiently regular for $\tilde{H}^s(\Omega;E)$ and $H^s_{\olO}(M;E)$ to agree (see Proposition \ref{EqTildeH}), one gets a much nicer result. 
\begin{cor}
If $\Omega$ is a $C^0$-domain, then it holds that
$$
\tilde{H}^{s_\theta}(\Omega;E)= \left[ \tilde{H}^{s_0}(\Omega;E),\tilde{H}^{s_1}(\Omega;E) \right]_\theta = \left[ H^{s_0}_{\olO}(M;E),H^{s_1}_{\olO}(M;E) \right]_\theta = H^{s_\theta}_{\olO}(M;E)
$$
for all $s_0,s_1\geq 0$ and $\theta\in [0,1]$.
\end{cor}
Recall that $U_\delta(\Omega)$ and $N_\delta(\Omega)$ denote the open and closed $\delta$-neighborhood of $\Omega$ respectively. An analogous argument to Proposition \ref{InterInclusion} also provides the following statement, which may not be all that appealing but will turn out very useful nonetheless, see Lemma \ref{FormalInterpolationInterior} below.
\begin{prop}\label{SubspaceInterpolationDelta}
Let $\Omega\subseteq M$ be an arbitrary domain. Let $s_0,s_1\geq 0$ and $\theta\in [0,1]$, and set $s_\theta:=(1-\theta)s_0+\theta s_1$. Then, for any $\delta>0$ one has inclusions
$$
\tilde{H}^{s_\theta}(\Omega;E) \subseteq \left[ \tilde{H}^{s_0}(U_\delta(\Omega);E), \tilde{H}^{s_1}(U_\delta(\Omega);E) \right]_\theta
$$
and 
$$
\left[ \tilde{H}^{s_0}(\Omega;E), \tilde{H}^{s_1}(\Omega;E) \right] \subseteq \tilde{H}^{s_\theta}(U_\delta(\Omega);E) \; .
$$
The analogous statement for $H^s_{\olO}(M;E)$ and $H^s_{N_\delta(\Omega)}(M;E)$ also holds.
\end{prop}
\section{The uniformly summable Rellich-Kondrachov theorem}\label{SectionQuantComp}
We conclude our discussion of Sobolev spaces on domains with the Rellich-Kondrachov theorem. The observation that the Sobolev embedding on compact domains is not only compact but finitely summable is classical, going back to e.g. \cite{Weyl1912}. On manifolds of bounded geometry and domains in them this finite summability is even uniform. This observation will be the foundation of the uniform compactness statements required to extract uniform K-homology classes from uniformly elliptic operators.
\begin{thm}\label{QuantComp}
Let $M$ be an $n$-dimensional manifold of bounded geometry, $E\to \Omega$ a rank-$m$ vector bundle of bounded geometry, and $\Omega\subseteq M$ a domain. For all $R\geq0$, $s\geq 0$ $t>0$ and all $p>\frac{n}{t}$ there exist constants $\eta_{R,s,t,p}$ and $\eta_{R,s,t,p}^\Omega$, dependent only on $M$, $E$ and the indicated data, such that the following holds:
\begin{itemize}
\item[(i)] For all compact subsets $K\subseteq M$ with $\diam(K)\leq R$ the inclusion
$$
H^{s+t}_K(M;E) \xlongrightarrow{j} H^s_K(M;E) 
$$
is $p$-summable with $||j||_p\leq \eta_{R,s,t,p}$.
\item[(ii)] For all compact subsets $K\subseteq \olO$ with $\diam(K)\leq R$ the inclusion
$$
H^{s+t}_{K}(\olO;E) \xlongrightarrow{j} H^s_K(\olO;E)
$$
is $p$-summable with $||j||_p\leq \eta_{R,s,t,p}^\Omega$.
\end{itemize}
\end{thm}
We adapt the argument given in \cite{Rave2012}. First, we show the statement for vector-valued functions on the torus. For manifolds of bounded geometry the statement can be reduced to that case via local trivializations. The case of domains then follows via extension. 
\begin{lemma}\label{RellichTorus}
Let $\mathbb{T}^n_r$ be the flat torus obtained by identifying the opposite faces of an $n$-dimensional cube with side length $r>0$. Let $s\geq 0$ and $t>0$. Then, the embedding $j: \, H^{s+t}(\mathbb{T}^n_r;\CC^m) \to H^s(\mathbb{T}^n_r;\CC^m)$ is $p$-Schatten class for all $p>\frac{n}{t}$.
\end{lemma}
\begin{proof}
The Sobolev space $H^s(\mathbb{T}^n_r;\CC^m)$ for $s\geq 0$ consists of those $u\in L^2(\mathbb{T}^n_r;\CC^m)$ for which $\langle u,u\rangle_{H^s}<\infty$, where the $H^s$-inner product is given by
$$
\langle u,v\rangle_{H^s} = \sum_{k\in\ZZ^n} \langle \hat{u}_k, \hat{v}_k \rangle_{\CC^m} \cdot \langle k \rangle^{2s} \quad \mathrm{where} \quad \langle k\rangle := \left(1+\left( \frac{2\pi ||k||}{r} \right)^2 \right)^{1/2} \; .
$$
The elements $f^s_{k,i}(x):= r^{-1/2} \langle k\rangle^{-s}\cdot e^{2\pi i \langle k,x\rangle/r} \cdot e_i$, $k\in\ZZ^n$, $i=1,\cdots ,m$, form a maximal orthonormal system. Here $e_i$ is the $i^\mathrm{th}$ standard basis vector of $\CC^m$.

The $p$-norm of $j$ is obtained from the spectral values of the operator $|j|=\sqrt{j^*j}$ acting on $H^{s+t}(\mathbb{T}^n_r,\CC^m)$. We have $j(f^{s+t}_{k,i})=f^{s+t}_{k,i}=\langle k \rangle^{-t} f^s_{k,i}$, from which it follows that 
$$
\langle f^{s+t}_{k',i'}, (j^*j)f^{s+t}_{k,i}\rangle_{H^{s+t}} = \langle jf^{s+t}_{k',i'}, jf^{s+t}_{k,i}\rangle_{H^s} = \langle k \rangle^{-2t} \cdot \delta_{k,k'}\delta_{i,i'}  \; ,
$$
whence
$$
(j^*j)f^{s+t}_{k,i} = \langle k \rangle^{-2t} \cdot f^{s+t}_{k,i} \; .
$$ 
The $p$-norm of $j$ is therefore equal to
$$
||j||_p = \left(\sum_{k\in\mathbb{Z}^n} m\cdot \langle k \rangle^{-pt}\right)^{\frac{1}{p}} \; .
$$
To see that this sum converges for all $p>\frac{t}{n}$, rearrange the estimate
$$
\left(1+\left(\frac{2\pi k_1}{r}\right)^2\right)\cdots\left(1+\left(\frac{2\pi k_n}{r}\right)^2\right) \leq \left(1+ \left( \frac{2\pi \|k\|}{r}\right)^2\right)^n = \langle k \rangle^{2n}
$$
to obtain
$$
\langle k \rangle^{-pt} \lsim \prod_{i=1}^n \left(\frac{1}{n}+\left(\frac{2\pi k_i}{r}\right)^2\right)^{-\frac{pt}{2n}} \lsim \prod_{i=1}^n (1+|k_i|)^{-\frac{pt}{n}} \; .
$$
The series $\sum_{k_i\in\ZZ} (1+|k_i|)^{-pt/n}$ converges if and only if $p>\frac{n}{t}$. For these $p$ it then follows that
$$
||j||_p^p \lsim  \sum_{k\in\ZZ} \prod_{i=1}^n (1+|k_i|)^{-\frac{pt}{n}} = \left(\sum_{a\in\ZZ} (1+a)^{-\frac{pt}{n}} \right)^n < \infty \; .
$$
\end{proof}
We use this local computation to obtain the theorem globally on manifolds of bounded geometry.
\begin{proof}[Proof of Thm. \ref{QuantComp}]
Ad (i): Let $(U_i,\kappa_i,\psi_i,\phi_i)$ be an admissible trivialization of $M$ and $E$. There is $N\in\NN$ -dependent on $R$ but not the specific $K$- such that $K$ intersects at most $N$ of the $U_i$. Label those $U_i$ as $U_1,\cdots,U_N$. It holds that 
$$
||j||_p \leq \sum_{i=1}^N ||j(\psi_i \cdot)||_p \; .
$$
Multiplication by $\psi_i$ maps $H^{s+t}_K(M;E)$ to $H^{s+t}_{K_i}(U_i;E)$, where $K_i:= K\cap \supp(\psi_i)$ is compact in $U_i$, with norm bounded independently of $i$ by bounds on the first $k$ derivatives of $\psi_i$ for $\NN\ni k\geq s+t$. Using the coordinate maps and trivializations, transfer to vector-valued functions on Euclidean space via \\
\centerline{
\xymatrixcolsep{5pc}\xymatrix{ H^{s+t}_{K_i}(U_i;E) \ar[r]^-{\phi_i \circ (\cdot) \circ \kappa_i^{-1}} \ar[d]_j & H^{s+t}_{\kappa_i(K_i)}\left( \kappa_i(U_i) , \mathbb{C}^m\right) \ar[d]^j \\
H^s_{K_i}(U_i;E) \ar[r]_-{\phi_i \circ (\cdot) \circ \kappa_i^{-1}} & H^s_{\kappa_i(K_i)}\left( \kappa_i(U_i) , \mathbb{C}^m\right)
}}
The horizontal arrows are isomorphism with norms bounded independently of $i$ by bounds on the first $k$ derivatives of $\phi_i$ and $\kappa_i$. It therefore suffices to prove that the right vertical arrow is $p$-summable. To that end recall that the diameter of $\kappa_i(U_i)$ is strictly bounded by the injectivity radius $r:=\mathrm{inj}(M)>0$ of $M$. By putting $\kappa_i(U_i)$ inside a cube of side length $r$ and identifying opposite faces, we can view $\kappa_i(U_i)$ as an open subset of the torus $\mathbb{T}^n_r$. This allows us to view $j: \, H^{s+t}_{\kappa_i(K_i)}\left( \kappa_i(U_i) , \mathbb{C}^m\right)\to H^{s}_{\kappa_i(K_i)}\left( \kappa_i(U_i) , \mathbb{C}^m\right)$ as the restriction of the inclusion map $j: \, H^{s+t}\left( \mathbb{T}^n_r , \mathbb{C}^m\right) \to H^s(\mathbb{T}^n_r,\CC^m)$ to the closed subspace of functions with support in the compactum $\kappa_i(K_i)\subset \mathbb{T}^n_r$. That latter inclusion map is $p$-summable by Lemma \ref{RellichTorus}, hence so is its restriction. Collecting the different estimates, we conclude that 
$$
||j||_p \leq \eta_{R,s,t,p} := C(n,k)\cdot N \cdot C_k(\psi)\cdot C_k(\phi) \cdot C_k(\kappa) \cdot J(r,m,t,p) < \infty \; .
$$
Here $k\in\NN$ satisfies $k\geq t+s$, $C_k(\psi_i)$ is an $i$-independent bound on the first $k$ derivatives of $\psi_i$, $C_k(\phi)$ and $C_k(\kappa)$ are defined analogously, $C(n,k)$ is a universal constant\footnote{This universal constant arises from the specifics of how multiplication by $\psi_i$ and composition by $\phi_i$ and $\kappa_i$ can be bounded in Sobolev norm. These norms are bounded by bounds on the derivatives of these maps, times constants coming from application of the chain rule, etc.}, and $J(r,m,t,p)$ is the $p$-Schatten norm of $j: \, H^{s+t}(\mathbb{T}^n_r,\CC^m) \to H^s(\mathbb{T}^n_r,\CC^m)$.

Ad (ii): For $K\subseteq \olO$ pick a closed ball $\tilde{K}$ of radius $R+1$ such that $K$ is contained in the interior of $\tilde{K}$. We get a commutative diagram\\
\centerline{\xymatrix{
H^{s+t}(M;E) \ar[r]^{r} \ar[d]_j & H^{s+t}\left( \olO , E\right) \ar[d]^j \\
H^s(M;E) \ar[r]_{r} & H^s\left( \olO;E\right)
}}
where $r$ denotes the restriction maps. As a surjective map between Hilbert spaces $r$ admits a continuous splitting $\sigma: \, H^{s+t}(\olO;E)\to H^{s+t}(M;E)$. Use Lemma \ref{BdBump} to produce a function $C^\infty_c(M)$ with $\supp(\psi)\subseteq\tilde{K}$, $\psi|_K=1$, and such that $\|\psi\|_{C^k_b}$ and can be bounded by a constant dependent on $R$ but not the specific $K$, where $k\in\NN$, $k\geq s+t$. Then, the inclusion $j: \, H_K^{s+t}(\olO;E) \to H_K^s(\olO;E)$ can be factored as
$$
H_K^{s+t}(\olO;E) \xrightarrow{\sigma} H^{s+t}(M;E) \xrightarrow{\psi} H^{s+t}_{\tilde{K}}(M;E) \xrightarrow{\tilde{j}} H^s_{\tilde{K}}(M;E) \xrightarrow{r} H^s(\olO;E) \; . 
$$
The inclusion $\tilde{j}: \, H^{s+t}_{\tilde{K}}(M;E) \to H^s_{\tilde{K}}(M;E)$ is $p$-summable by item (i), so we get that $j$ is also $p$-summable with
$$
||j||_p \leq  ||\sigma||\cdot ||\psi||_{C^k_b} \cdot ||\tilde{j}||_p \; ,
$$
and by construction of $\psi$ and item (i) the right-hand side can be bounded by a constant $\eta^\Omega_{R,s,t,p}$ depending only on $R$ but not the specific $K$.
\end{proof}
\chapter{Uniformly elliptic operators}
It was observed by \v{S}pakula that Dirac-type operators over manifolds of bounded geometry define classes in uniform K-homology \cite{Spakula2009}. Such operators are examples of \emph{uniformly elliptic operators}, meaning that their coefficients are uniformly bounded and their principal symbol is not only invertible but invertible in a uniform way. Engel then extended \v{S}pakula's result to uniformly elliptic pseudo-differential operators \cite{Engel2018}. In this chapter we will investigate uniformly elliptic operators over domains in manifolds of bounded geometry as well as their boundary conditions. We focus on differential operators as opposed to pseudo-differential ones, though see Subsection \ref{UPsiDOsSingular}. The analysis of uniformly elliptic operators and their boundary conditions developed in this chapter combined with that of the last chapter on Sobolev spaces will enable us to prove the existence of uniform K-homology classes defined by these operators in Chapter 6.

There are two reasonable approaches for treating differential operators with Hilbert space techniques. The first is operator-theoretic: A differential operator may be viewed as an unbounded operator on $L^2$-space with domain the space of smooth compactly supported sections. The second is via Sobolev spaces, between which differential operators map as bounded operators. The two approaches turn out to be the same for uniformly elliptic operators on manifolds of bounded geometry. Indeed, any such operator -when viewed as an unbounded operator on $L^2$- has a canonical closed extension, whose domain is the Sobolev space of degree given by the order of the operator. The nice behavior of Sobolev spaces thus implies a nice $L^2$-theory for elliptic operators.

If we consider manifolds with boundary instead, or more generally domains in manifold, the situation becomes more complicated. Due to the presence of a boundary, there need not be a canonical extension, but there might be many. One should think of these extensions as arising from different boundary conditions. These extensions may or may not be contained in Sobolev space. Every boundary condition satisfies so-called \emph{interior regularity}, stating roughly that in the interior the domain of any boundary condition coincides with Sobolev space. Further, we consider \emph{regular} boundary condtition, whose domain is fully contained in Sobolev space.\footnote{Technically, we will call boundary conditions with $\dom(D_e)\subseteq H^m(\olO)$ \emph{semi-regular}, with the term \emph{regular} being reserved for boundary conditions that additionally satisfy $\dom(D_e^*)\subseteq H^m(\olO)$ as well. For the purpose of this introductory discussion we skim over this distinction.} While interior regularity suffices to prove that uniformly elliptic operators define uniform K-homology classes on the interior, regularity is needed to obtain relative uniform K-homology classes up to the boundary. This is because the uniform compactness properties of Sobolev spaces, in particular the uniformly summable Rellich-Kondrachov theorem, are necessary to obtain analogous uniform compactness statements for the domains of boundary conditions, and these are needed to prove the existence of uniform K-homology classes.

Boundary conditions for uniformly elliptic operators on manifolds with boundary and bounded geometry and their boundary conditions have been discussed in the literature already. The theory of boundary conditions and regularity for first-order elliptic operators on non-compact manifolds with boundary has been developed Bär and Bandara \cite{Baer2022}, building on earlier work for compact boundaries \cite{BaerBallmann2012}. In \cite{Grosse2014} Große and Nakad treats the special case of Dirac operators on spin$^c$ manifolds with boundary and bounded geometry. Schick formulated regularity for uniformly elliptic operators on manifolds with boundary and bounded geometry in terms of the existence of a uniform pseudo-differential solution operator, and applied this to a discussion of the Laplacian and $L^2$-Hodge theory in the presence of boundaries \cite{Schick1998}. For second order operators in this setting various notions of regularity have been discussed and compared by Große and Nistor \cite{Grosse2020}. We will discuss how these regularity notions relate to the one used in this thesis below.

In Section 5.1 we introduce uniformly elliptic operators on manifolds of bounded geometry, and on domains in these manifolds. Section 5.2 initiates the study of boundary conditions. Besides a proof of interior regularity it contains a discussion of certain grading and locality properties that will also be relevant to the existence of uniform K-homology classes. In Section 5.3 we study the regularity notion used here, and also connect it to those of \cite{Baer2022}, \cite{Schick1998}, and \cite{Grosse2020}. Section 5.4 is devoted to examples. The list of examples is by no means exhaustive, but rather supposed to illustrate the different properties of boundary conditions that have been introduced in the preceding sections. Lastly, Section 5.5 contains a brief account of the finite-propagation property of wave operators and more general functions of a local boundary condition for a first-order operator.
\section{Definition}
We begin with the definition of uniformly elliptic operators over manifolds of bounded geometry and domains in such manifolds. Let $M$ be a manifold of bounded geometry,  and $E,F\to M$ vector bundles of bounded geometry. Let $D: \, \Gamma_c(M;E)\to \Gamma_c(M;F)$ be a differential operator of order $m$ with principal symbol $\sigma$. Let $(U_i,\kappa_i,\phi_i^{E},\phi_i^F)$ be an admissible trivialization of $M$, $E$ and $F$. Consider the differential operator 
$$
D_i=\left( \phi^F_i D (\phi^{E}_i)^{-1}\right) \circ \kappa_i^{-1} =: \sum_{|\alpha|\leq m} a_{i}^\alpha \partial^\alpha
$$
over $\kappa_i(U_i)$. Its principal symbol is given by $\sigma_i:= (\phi^F_i \sigma \phi^{E}_i)^{-1}\circ \kappa_i^{-1}$. 
\begin{definition}\label{UEDOnoBdry}
The differential operator $D: \, \Gamma_c(M;E)\to \Gamma_c(M;F)$ is called \emph{bounded} if for all $k\geq 0$ there exists $C_k\geq 0$ such that
$$
\left| \partial^\beta a_i^\alpha(x) \right| \leq C_k
$$
for all $i\in I$, $|\alpha|\leq m$ and $|\beta|\leq k$.

Moreover, $D$ is called \emph{uniformly elliptic} if $D$ is bounded, the principal symbol $\sigma(\xi)$ is invertible for all $\xi\in T^*M$, $\xi\neq 0$, and there exists a constant $C>0$ such that
$$
\left| \sigma_i(x,\xi)^{-1} \right| \leq C \cdot (1+|\xi|)^{-m} 
$$
for all $i\in I$.
\end{definition}
Since the transition functions of an admissible trivialization have uniformly bounded derivatives, uniform ellipticity is independent of the choice of admissible trivialization.
\begin{definition}\label{UEDOdomain}
Let $\Omega \subseteq M$ be a domain. A differential operator $D: \, \Gamma_{cc}(\Omega;E)\to\Gamma_{cc}(\Omega;F)$ is called \emph{uniformly elliptic over $\olO$} if there exists a uniformly elliptic differential operator $\tilde{D}: \, \Gamma_{c}(M;E)\to\Gamma_{c}(M;F)$ such that $\tilde{D}u=Du$ for all $u\in \Gamma_{cc}(\Omega;E)$.
\end{definition}
This rather extrinsic definition of uniformly elliptic operators has the advantage that we can transfer knowledge of uniformly elliptic operators over manifolds of bounded geometry rather directly to such operators over domains. Thus it is technically convenient. Moreover, since this definition does not require use of any coordinates on the domain directly, it works for arbitrarily rough domains. However, for manifolds with boundary and bounded geometry there is also an intrinsic definition of uniformly elliptic operator. It is analogous to Definition \ref{UEDOnoBdry}, only that uniform boundedness and ellipticity is also required in boundary charts. This definition is essentially equivalent to that given here. More concretely, we have the following.
\begin{bonuslemma}\hspace{-5pt}\footnote{The alphabetic index means that this lemma was not present in the original version of this thesis. Alphabetic indexing was chosen so that the numerical indexing is consistent with the original version.}
Let $\olO$ be an intrinsically defined manifold with boundary and bounded geometry, and $D: \, \Gamma_{cc}(\Omega,E)\to \Gamma_{cc}(\Omega,F)$ a differential operator mapping between sections of bounded-geometry vector bundles over $\olO$.

Suppose $M$ is a manifold of bounded geometry containing an isometric copy of $\olO$ as a bounded-geometry domain. Then, if $D$ is uniformly elliptic in the in the sense of Definition \ref{UEDOdomain}, it is also uniformly elliptic in the intrinsic sense.

Conversely, if $D$ is uniformly elliptic in the intrinsic sense, there exists a manifold $M$ of bounded geometry containing an isometric copy of $\olO$ as a bounded-geometry domain, such that $D$ admits a uniformly elliptic extension to $M$, i.e. $D$ uniformly elliptic in the sense of Definition \ref{UEDOdomain}.
\end{bonuslemma}
\begin{proof}
The first statement follows directly from the fact that any uniformly elliptic extension of $D$ to $M$ must necessarily satisfy the bounds necessary for the intrinsic definition of being uniformly elliptic.

For the reverse implication, suppose for simplicity that $\olO$ has product structure on a uniformly thick collar of the boundary. This can always be arranged via a bi-Lipschitz equivalence \cite[Proposition 7.3]{Schick1998}, and the property of being uniformly elliptic is preserved under this change of metric. Form $M:=\olO\cup_{\dO}(\RR_{\geq 0}\times\dO)$ by attaching cylinders to the boundary. This is a manifold of bounded geometry containing $\olO$ as a bounded-geometry domain. Suppose $D: \, \Gamma_{cc}(\Omega,E)\to \Gamma_{cc}(\Omega,F)$ is intrinsically uniformly elliptic. Up to a uniformly bounded homotopy $E$ and $F$ are constant in normal direction on a uniformly thick tubular neighborhood of $\dO$ in $M$. By keeping them constant in that direction, $E$ and $F$ can be extended to bounded-geometry vector bundles over $M$. Moreover, near the boundary we may write $D = a_t \partial_t + A_t$ , with $A_t$ a smooth family of differential operators over $\dO$. Since $D$ is intrinsically uniformly elliptic, $a_t$ and $A_t$ can be extended to $[0,\delta)\times\dO$ for some $\delta>0$ in such a way that $a_t\partial_t+A_t$ is also uniformly elliptic for those $t$. Fix some smooth $\psi: \RR_{\geq 0}\to[0,1]$ with $\psi(t)=1$ for $t\leq\frac{\delta}{3}$, and $\psi(t)=0$ for $t\geq\frac{2\delta}{3}$. Then
$$
\tilde{D} = \psi(t)\cdot \left( a_t \partial_t + A_t \right) + (1-\psi(t))\cdot\left( a_0 \partial_t + A_0 \right)
$$
is a uniformly elliptic extension of $D$, whence $D$ is uniformly elliptic in the sense of Definition \ref{UEDOdomain}.
\end{proof}
\begin{example}
The Laplace operator, the spin and Clifford-linear Dirac operator, the Euler characteristic operator, and related geometric operators over a manifold of bounded geometry are uniformly elliptic, hence so are their restrictions to arbitrary domains. Moreover, strongly elliptic operators over domains in Euclidean space (see e.g. \cite{McLean2000}) are uniformly elliptic in our sense.
\end{example}
Associated to any differential operator $D: \, \Gamma_{cc}(\Omega;E)\to \Gamma_{cc}(\Omega,F)$ there is a differential operator $D^\dagger: \, \Gamma_{cc}(\Omega;F)\to \Gamma_{cc}(\Omega,E)$ called the \emph{formal adjoint} of $D$, which is uniquely defined by the requirement that
$$
\langle Du,v \rangle_{L^2(\Omega;F)} = \langle u, D^\dagger v\rangle_{L^2(\Omega;E)}
$$
hold for all $u\in \Gamma_{cc}(\Omega,E)$ and $v\in \Gamma_{cc}(\Omega,F)$. The operator $D$ is called \emph{symmetric} (or \emph{formally self-adjoint}) if $D=D^\dagger$. Clearly it holds that $(D^\dagger)^\dagger=D$. The coefficients of $D^\dagger$ in local charts are given by expressions involving the adjoints of those of $D$ as well as their derivatives. In particular the principal symbol $\sigma_{D^\dagger}$ of $D^\dagger$ is the point-wise adjoint $\sigma_D^*$ of that of $D$. Then, the following is easily seen by taking the formal adjoint of any uniformly elliptic extension.
\begin{lemma}
A differential operator over $\olO$ is uniformly elliptic if and only if its formal adjoint is.
\end{lemma}
Uniformly elliptic theory over manifolds of bounded geometry is well-developed. The essential property of elliptic operators over closed manifolds is that they admit pseudo-differential quasi-inverses modulo smoothing operators, and the same is true over manifolds of bounded geometry as well. The precise set-up depends on the pseudo-differential calculus one wishes to use, see for example \cite{Kordyukov1991}, \cite{Shubin1992}, \cite{Taylor2008} or \cite{Engel2018}. We will not employ pseudo-differential techniques in this thesis, but one consequence of the existence of a pseudo-differential quasi-inverse will be important to us, namely elliptic regularity. It provides a bound for the Sobolev norm of a section in terms of lower-degree Sobolev norms of that section and its image under the uniformly elliptic operator. Given our definition of uniformly elliptic operators over domains it is easily deduced that this form of regularity also holds for these operators, at least as long as the section in question is supported away from the boundary. This is regularity estimate is the core of interior regularity, which will be discussed below.
\begin{prop}\label{IntRegEst}
Let $D$ be uniformly elliptic over $\olO$. Then, for all $s\in\RR$ there exists a constant $C_s\geq 0$ such that for all $u\in \Gamma_{cc}(\Omega;E)$ it holds that
$$
\left\| u \right\|_{H^s(M;E)} \leq C_s \cdot \left( \|u\|_{H^{s-m}(M;E)} + \| Du\|_{H^{s-m}(M;F)} \right) \; .
$$
\end{prop}
\begin{proof}
This follows from applying the fundamental elliptic estimate (for example \cite[Theorem 6.6]{Engel2018}) to any uniformly elliptic extension of $D$ to $M$.
\end{proof}
We draw explicit attention to the fact that the interior regularity estimate uses the norms on $H^s(M)$ instead of $H^s(\olO)$. Of course, these are equivalent if $s$ is an integer. 
\section{Boundary conditions for uniformly elliptic operators}
In this section we study the properties of uniformly elliptic operators, viewed as unbounded operators on the Hilbert space of $L^2$-sections. Standard facts about unbounded operators are collected in Appendix \ref{AppendixUnbdOps}, and will be used without reference. In particular, the notions of closed extensions and adjoints of unbounded operators will be used heavily in our discussion of boundary conditions.

Let $M$ a manifold of bounded geometry, $E,F\to M$ vector bundles of bounded geometry, and $\Omega\subseteq M$ a domain. Let $D: \, \Gamma_{cc}(\Omega;E)\to \Gamma_{cc}(\Omega;F)$ be a differential operator. We view $D$ as a densely defined, unbounded operator from $\Gamma_{cc}(\Omega;E)\subset L^2(\Omega;E)$ to $L^2(\Omega;F)$. Similarly, the formal adjoint $D^\dagger$ is an unbounded operator from $L^2(\Omega;F)$ to $L^2(\Omega;E)$.
\begin{definition}\label{DefBdryCond}
A \emph{boundary condition} for $D$ s a closed extension $D_e$ of $D$, such that the Hilbert space adjoint $D_e^*$ is a closed extension of $D^\dagger$.
\end{definition}
It may not be clear how this definition captures what we intuitively understand as a boundary condition. We will see in Proposition \ref{MinMaxDomain} below that -at least for uniformly elliptic operators- boundary conditions do in fact only depend on the behavior on the boundary. But there is also a more sophisticated description of this phenomenon. Indeed, in favorable circumstances there is a one-to-one correspondence between boundary conditions for an elliptic operator $D$ in the above sense and closed subspaces of a certain Sobolev-like space on the boundary, so that a section is in the domain of the boundary condition if and only if its restriction to the boundary lies in the corresponding subspace. In the case of first-order operators and compact boundaries this is described in \cite{BaerBallmann2012}, and extended to non-compact boundaries in \cite{Baer2022}. Dirac operators on spin$^c$ manifolds with boundary and bounded geometry are treated in \cite{Grosse2014}. Higher-order elliptic operators on compact manifolds with boundary are treated in \cite{Bandara2023}. While this description of boundary conditions is conceptually satisfying the definition in terms of closed extensions suffices for our purposes. 
\begin{definition}
The \emph{minimal boundary condition} $D_{\min}$ of $D$ is the closure of $D$ as an unbounded operator from $L^2(\Omega;E)$ to $L^2(\Omega;F)$. The \emph{maximal boundary condition} $D_{\max}$ is the (Hilbert space) adjoint $D_{\max}=(D^\dagger)^*$ of $D^\dagger$ as an unbounded operator from $L^2(\Omega;F)$ to $L^2(\Omega;E)$.
\end{definition}
A section $u\in L^2(\Omega;E)$ is in the domain of $D_{\min}$ if and only if there exists a sequence $(u_n)_{n\in\NN}$, $u_n \in \Gamma_{cc}(\Omega;E)$, such that $||u-u_n||_{L^2} \to 0$ and $\sup_n ||Du_n||_{L^2}< \infty$. Equivalently, $\dom(D_{\min})$ is the closure of $\Gamma_{cc}(\Omega;E)$ in the graph norm $||u||_D=||u||_{L^2}+ ||Du||_{L^2}$, and $D_{\min}$ is the unique extension of $D$ to this closure. A section $u\in L^2(\Omega;E)$ is in the domain of $D_{\max}$ if and only if there exists $v\in L^2(M;F)$ such that for all $w\in \Gamma_{cc}(\Omega;F)$ it holds that
$$
\langle u , D^\dagger w\rangle_E = \langle v, w \rangle_F \; .
$$
Then $D_{\max}u= v$. 

Both $D_{\min}$ and $D_{\max}$ are closed extensions of $D$, and it holds that $D_{\min} \subseteq D_{\max}$. Moreover, it holds that $D_{\max}=(D^\dagger)^*=(D_{\min}^\dagger)^*$, which implies that both the minimal and the maximal boundary conditions are indeed boundary conditions in the sense of Definition \ref{DefBdryCond}. Moreover, any $D_e$ boundary condition for $D$ satisfies $D_{\min}\subseteq D_e\subseteq D_{\max}$. In fact, a closed extension satisfies Definition \ref{DefBdryCond} if and only if this condition holds. This explains the nomenclature of minimal and maximal boundary condition.

The following proposition states that on the interior of the domain $\Omega$ there is no difference between minimal and maximal boundary condition. Consequently, all boundary conditions do indeed behave the same away from the boundary.
\begin{prop}\label{MinMaxDomain}
Let $D: \, \Gamma_{cc}(\Omega;E)\to\Gamma_{cc}(\Omega;F)$ be uniformly elliptic over $\olO$. Let $u \in L^2(M;E)$ have $\supp(u)\subset \Omega$. Then
$$
u \in \dom(D_{\min}) \quad \Leftrightarrow \quad u\in \dom(D_{\max})  \; .
$$
\end{prop}
\begin{proof}
Let $\tilde{D}: \, \Gamma_c(M;E)\to \Gamma_c(M;F)$ be a uniformly elliptic extension of $D$. It holds that $\dom(\tilde{D}_{\min})=\dom(\tilde{D}_{\max})$ \cite[Proposition 4.1]{Shubin1992}. If $u\in L^2(\Omega;E)$ has support inside $\Omega$, then it is straight-forward to see that $u\in\dom(D_{\max})$ if and only if $u \in \dom(\tilde{D}_{\max})$. Thus, in that case there exists a sequence $(v_n)_{n\in\NN}$ in $\Gamma_c(M;E)$ converging to $u$ in the graph norm of $\tilde{D}$. The supports of the $v_n$ can be arranged to lie inside $\Omega$ via cut-offs. Since the graph norm of $\tilde{D}$ on $\Gamma_{cc}(\Omega;E)$ coincides with the one of $D$, it follows that $u\in \dom(D_{\min})$.
\end{proof}
\begin{remark}
The statement of Proposition \ref{MinMaxDomain} holds for first-order differential operators without the assumption of ellipticity, see for example \cite[Lemma 10.2.5]{HigsonRoe2000}. However, it can fail for non-elliptic higher-order operators, see for example \cite{Tausk2010}. 
\end{remark}
The next proposition identifies the domain of the minimal boundary condition with the Sobolev space $\tilde{H}^m(\Omega;E)$. This essential fact allows for the transfer of knowledge about Sobolev spaces to knowledge about boundary conditions of uniformly elliptic operators. Note that by Lemma \ref{EqualityH0} we could also use the Sobolev space $H_0^m(\Omega;E)$ instead. However, since we use the $H^m(M;E)$-norm instead of the $H^m(\olO;E)$-norm here, and the difference between the two becomes substantial when we consider fractional degrees later, we prefer to use $\tilde{H}^m(\Omega;E)$.
\begin{prop} \label{SobolevDomain}
Let $D: \, \Gamma_{cc}(\Omega;E)\to\Gamma_{cc}(\Omega;F)$ be uniformly elliptic over $\olO$. Then, it holds that
$$
\dom(D_{\min}) = \tilde{H}^m(\Omega;E) 
$$
with equivalent norms.
\end{prop} 
\begin{proof}
The estimate of Proposition \ref{IntRegEst} for the special case $s=m$ expresses that 
$$
||u||_{H^m(M;E)} \lsim \|u\|_{L^2(M;E)} + \|Du\|_{L^2(M;E)} =  ||u||_D
$$
for $u\in\Gamma_{cc}(\Omega;E)$. Conversely, $D$ is bounded over $\olO$, we have $||Du||_{L^2(M;E)}\lsim ||u||_{H^m(M;E)}$ and hence $||u||_D\lsim ||u||_{H^m(M;E)}$ for $u\in\Gamma_{cc}(\Omega;E)$. We conclude that the $H^m(M;E)$-norm and the graph norm of $D$ are equivalent norms on $\Gamma_{cc}(\Omega;E)$, meaning that closures of $\Gamma_{cc}(\Omega;E)$ with respect to the two norms agree. This proves the claim.
\end{proof}
We refer to the combination of Propositions \ref{MinMaxDomain} and \ref{SobolevDomain} as \emph{interior regularity}. It expresses that for a uniformly elliptic operator a section supported in the interior is in the domain of any boundary condition if and only if it is in the Sobolev space $\tilde{H}^m(\Omega;E)$, and that the graph norm and the Sobolev norm are equivalent for these sections. Interior regularity allows the transfer of the nice properties of Sobolev spaces to the domains of arbitrary boundary conditions, at least over the interior. This will be an important ingredient in our discussion of the uniform K-homology classes of boundary conditions for uniformly elliptic operators.

Let us already note a first consequence of interior regularity concerning multiplication by functions supported in the interior.
\begin{prop}
Let $D$ be uniformly elliptic over $\olO$, and let $f\in C^\infty_{cc}(\Omega)$. Multiplication by $f$ yields a bounded map
$$
\dom(D_{\max}) \longrightarrow \dom(D_{\min}) \; , \; u\longmapsto f\cdot u
$$
with norm proportional $\|f\|_{C^m_b(\olO)}$.
\end{prop}
\begin{proof}
Consider the subspace of $\dom(D_{\max})$ containing all sections supported at least distance $\delta$ from $\partial\Omega$ for some $0<\delta<d(\supp(f),\partial\Omega)$, and let $P: \, \dom(D_{\max})\to \dom(D_{\max})$ be the orthogonal projection onto that subspace.\footnote{Here we consider $\dom(D_{\max})$ as a Hilbert space under the inner product $\langle \cdot,\cdot\rangle_D=\langle\cdot,\cdot\rangle_{L^2} + \langle \cdot , |D_{\max}|\cdot \rangle_{L^2}$, which induces the norm $\|\cdot \|_D$.} Note that if $u\in \dom(D_{\max})$ lies in the kernel of $P$, then its support is contained in the $\delta$-neighborhood of $\partial\Omega$. By Proposition \ref{MinMaxDomain} $\im(P)$ is contained in $\dom(D_{\min})$, so that $u\mapsto fu, \, \im(P) \to \dom(D_{\min})$ is bounded with norm bound proportional to $\|f\|_{C^m_b(\olO)}$ by Proposition \ref{SobolevDomain} and Lemma \ref{MultWInftyNoBdry}. On the other hand, if $u\in \dom(D_{\max})$ lies in the kernel of $P$, then its support is contained in the $\delta$-neighborhood of $\partial\Omega$, whence $fu=0$. We conclude that $u\mapsto fu$ factors as 
$$
\dom(D_{\max})\xlongrightarrow{P} \im(P) \xlongrightarrow{f} \dom(D_{\min})
$$ 
with claimed norm bound. The claim follows.
\end{proof}
From the previous proposition we can draw conclusions about arbitrary boundary conditions. Indeed, we have seen that multiplication by a smooth function with compact support in the interior maps $\dom(D_{\max})$ to $\dom(D_{\min})$. Since the domain of any boundary condition is contained in that of the maximal boundary condition, we conclude the following.
\begin{cor}\label{MultInt}
Suppose $D$ is uniformly elliptic, and $D_e$ a boundary condition for $D$. Let $f\in C^\infty_{cc}(M)$. Then, multiplication by $f$ is a bounded map
$$
\dom(D_e) \longrightarrow \dom(D_{\min})
$$
with norm proportional to $\|f\|_{C^k_b(\olO)}$.
\end{cor}
In order for a boundary condition to define a class in certain uniform K-homology groups it needs to satisfy certain properties. Besides regularity, which will be discussed in the next section, this concerns compatibility with (multi)gradings, but also with multiplication by functions. We conclude this section with a discussion of these compatibility conditions. 
\begin{definition}
Let $D: \, \Gamma_{cc}(\Omega;E)\to \Gamma_{cc}(\Omega;E)$ be a differential operator acting on sections of a $p$-multigraded vector bundle, $p\geq 0$. Suppose that $D$ is odd and $p$-multigraded. A boundary condition $D_e$ is called \emph{odd} if the grading operator $\kappa$ of $E$ preserves $\dom(D_e)$, and $D_e \kappa = - \kappa D_e$. The boundary condition is called \emph{$p$-multigraded} if every multigrading operator $\varepsilon_j$ of $E$ preserves $\dom(D_e)$ and $D_e\varepsilon_j=\varepsilon_j D_e$. 
\end{definition}
The following lemma is deduced rather directly from the definition of the adjoint boundary condition.
\begin{lemma}\label{AdjointMultigraded}
A boundary condition is odd resp. mulitgraded if and only if its adjoint is odd resp. multigraded.
\end{lemma}
Compatibility with multiplication by smooth functions is expressed through the following notion of locality.
\begin{definition}
A boundary condition $D_e$ for $D$ is called \emph{local} if $f\cdot\dom(D_e)\subseteq \dom(D_e)$ for all $f\in C^\infty_c(\olO)$. 
\end{definition} 
The nomenclature of calling such boundary conditions local is adopted from \cite{BDT1989}. However, we advise caution that for higher-order operators there can be differential boundary conditions that are not local in this sense (see Example \ref{ExNeumann} below). This situation does not appear for first-order operators (which is the case considered in \cite{BDT1989}). There any differential boundary condition must have order $0$, i.e. be given by the kernel of a bundle map on the boundary, and these boundary conditions are easily seen to be local. The first-order case is also convenient in the following way.
\begin{prop}\label{LocalFirstOrder}
Let $D_e$ be a boundary condition for the first-order differential operator $D$. Then, $D_e$ is local if and only if $D_e^*$ is local.
\end{prop}
\begin{proof}
The roles of $D_e$ and $D_e^*$ are symmetric, so it suffices to prove that if $D_e$ is local, so is $D_e^*$. It further suffices to prove that any real-valued function $f\in C^\infty_c(\olO)$ preserves $\dom(D^*_e)$. Let $u\in \dom(D_e^*)$. We need to produce a $v\in L^2(\Omega)$ such for every $w\in\dom(D_e)$ it holds that $\langle v,w\rangle = \langle fu,Dw\rangle$. Write 
$$
\langle fu,Dw\rangle = \langle u,fDw\rangle = \langle u, D(fw)\rangle + \langle u,[f,D]w\rangle =\langle fD_e^*u, w\rangle + \langle u,[f,D]w\rangle  \; .
$$
Since $D$ is first-order $[f,D]=-i\sigma_D(df)$ extends to a bounded operator on $L^2(\Omega)$, so that $[f,D_e]$ is bounded as well. Thus, $v:=fD_e^*u + [f,D_e]^*u$ is a well-defined element of $L^2(\Omega)$, for which it holds that
$$
\langle fu,Dw\rangle =\langle fD_e^*u, w\rangle + \langle u,[f,D]w\rangle = \langle v,w\rangle \; .
$$
This concludes the proof.
\end{proof}
For higher-order operators the proof fails because $[f,D]$ is no longer a bounded operator. Indeed Proposition \ref{LocalFirstOrder} is not generally true for higher-order operators, see Example \ref{ExMinMax} for an example. 
\section{Regularity of boundary conditions}
We come to regularity of boundary conditions. Regularity is crucial for our treatment of uniform K-homology classes, as it allows for the transfer of the favorable properties of Sobolev spaces to the entire domain of the boundary conditions, as opposed to interior regularity, which provides these properties only in the interior.

Throughout we consider a domain $\Omega$ in the manifold $M$ of bounded geometry, and the uniformly elliptic operator $D: \, \Gamma_{cc}(\Omega;E)\to \Gamma_{cc}(\Omega;F)$ over $\olO$ mapping between sections of the bounded-geometry vector bundles $E,F\to M$.
\begin{definition}\label{DefReg}
The boundary condition $D_e$ is called \emph{semi-regular} if there is a bounded inclusion $\dom(D_e)\hookrightarrow H^m(\olO;E)$ with equivalent norms. It is called \emph{regular} if both $D_e$ and $D_e^*$ are semi-regular.
\end{definition}
We phrase the definition this way to draw explicit attention to the fact that the inclusion $\dom(D_e)\hookrightarrow H^m(\olO;E)$ comes with a norm equivalence attached. It is a consequence of the open mapping theorem that an inclusion $\dom(D_e)\subseteq H^m(\olO;E)$ as a closed subspace already gives rise to such a norm equivalence.
\begin{lemma}\label{RegClosed}
Let $D_e$ be a boundary condition. Then $D_e$ is semi-regular if and only if $\dom(D_e)$ is a closed subspace of $H^m(\olO;E)$.
\end{lemma}
\begin{proof}
If $D_e$ is semi-regular, then $\dom(D_e) \subseteq H^m(\olO;E)$ with equivalent norms on the image. As the domain of a closed extension $\dom(D_e)$ is complete with respect to the graph norm, and hence also with respect to the Sobolev norm. Thus $\dom(D_e) \subseteq H^m(\olO;E)$ is closed.

Suppose conversely that $\dom(D_e) \subseteq H^m(\olO;E)$ is closed. Then, $\dom(D_e)$ is complete with respect to both the graph norm and the Sobolev norm. Moreover, because $D$ induces bounded maps between Sobolev spaces it holds that
$$
\|u\|_D = \|u\|_{L^2(\Omega;E)} + \|Du\|_{L^2(\Omega;F)} \lsim \|u\|_{H^m(\olO;E)} 
$$ 
for $u\in \dom(D_e)$. Consequently the identity map is bounded as a map from $(\dom(D_e),\|\cdot\|_{H^m(\olO;E)})$ to $\dom(D_e),\|\cdot \|_D)$. As a bounded bijection between Banach spaces it is already a Banach space isomorphism as well, meaning there is the converse estimate
\begin{equation}\label{RegularityDomain}
\|u\|_{H^m(\olO;E)} \lsim \|u\|_{L^2(\Omega;E)} + \|Du\|_{L^2(\Omega;F)} = \|u\|_D 
\end{equation}
for all $u\in\dom(D_e)$.
\end{proof}
To reiterate: A semi-regular boundary condition automatically comes with a regularity estimate as in \eqref{RegularityDomain}. Conversely, regularity estimates are sufficient to provide semi-regularity of a boundary condition. Let us state this in the particular case of pseudo-differential boundary conditions.
\begin{prop}\label{EstimateGivesReg}
Let $\olO$ be a manifold with boundary and bounded geometry. Let $D: \, \Gamma_{cc}(\Omega;E)\to \Gamma_{cc}(\Omega;F)$ be a uniformly elliptic operator of order $m$ over $\olO$. Moreover, let $P: \, \Gamma_c(\partial\Omega;E|_{\partial\Omega})\to \Gamma_c(\partial\Omega;V)$ be a uniform pseudo-differential operator\footnote{The precise pseudo-differential calculus is irrelevant here, any of the calculi of \cite{Kordyukov1991}, \cite{Taylor2008} or \cite{Engel2018} work. All that is needed is that $P$ gives rise to a bounded map $H^{m-\frac{1}{2}}(\dO;E|_{\dO})\to H^{m-\ell-\frac{1}{2}}(\dO,V)$. } of order $0\leq \ell \leq m-1$. Suppose that there exists $C>0$ such that
$$
\|u\|_{H^m(\olO;E)} \leq C\cdot\left(\|u\|_{L^2(\Omega;E)} + \|Du\|_{L^2(\Omega;F)} + \|Pu|_{\partial\Omega}\|_{H^{m-\ell-\frac{1}{2}}(\partial\Omega;V)}   \right)
$$
for all $u\in H^k(\olO;E)$. Then
$$
\dom(D_P) := \left\{ u\in H^k(\olO;E) \, | \, Pu|_{\partial\Omega}=0 \right\} \subseteq L^2(\Omega;E)
$$
defines a semi-regular boundary condition for $D$. 
\end{prop}
\begin{proof}
The space $\dom(D_P)\subseteq H^m(\olO;E)$ is closed in the Sobolev norm. We need to prove that $\dom(D_P)$ is closed in the graph-norm $\|\cdot\|_D$, i.e. that it indeed is a boundary condition. If this is done, semi-regularity follows from Lemma \ref{RegClosed}. To that end note first that the Sobolev norm dominates the graph norm. Conversely, for $u\in\dom(D_P)$ it holds that $Pu|_{\dO}=0$ so that the assumed regularity estimate for $u$ reads
$$
\|u\|_{H^m(\olO;E)} \lsim \|u\|_{L^2(\Omega;E)} + \|Du\|_{L^2(\Omega;E)} = \|u\|_D \; .
$$
We conclude that $H^m$-norm and graph norm are equivalent on $\dom(D_P)$. Since $\dom(D_P)$ is closed in the former, it is also closed in the latter. This proves the claim.
\end{proof}
Norm estimates as in Proposition \ref{EstimateGivesReg} have been studied extensively. There are different approaches to obtain them. One is through the existence of a pseudo-differential quasi-inverse (see for example \cite{Hörmander1963} for the classical case). This is analogous to regularity for elliptic operators on manifolds without boundary, where the existence of a pseudo-differential quasi-inverse also gave rise to a regularity estimate. In the case of manifolds with boundary and bounded geometry this approach to regularity of boundary conditions has been studied by Schick in \cite{Schick1998}. There a boundary condition admitting a uniform pseudo-differential quasi-inverse is called \emph{uniformly elliptic}, and it is shown that such boundary conditions admit a regularity estimate. 

Another approach is to require the existence of regularity estimate of associated local problems or problems with frozen coefficients that are suitable uniform. From this one then produces a global regularity estimate (see for example \cite[Chapter 5]{Taylor1996} for the compact case). For second-order operators on manifolds with boundary and bounded geometry this approach (among others) was investigated by Große and Nistor in \cite{Grosse2020}. For example they present a uniform version of the classical Shapiro-Lopatinski condition, and show that any boundary condition satisfying this condition admits a regularity estimate.

Let us also mention the work of Bär and Bandara on boundary conditions for first-order elliptic operators on manifolds with non-compact operators \cite{Baer2022}. This extends earlier work by Bär and Ballmann treating the case of compact boundaries \cite{BaerBallmann2012}, and by Große and Nakad treating Dirac operators on spin$^c$ manifolds with boundary and bounded geometry \cite{Grosse2014}. Bär and Bandara parametrize boundary conditions via trace maps in terms of closed subspaces of a certain Sobolev-like space of Cauchy-data on the boundary. There, a boundary condition is called semi-regular if its associated subspace is (locally) contained in $H^{\frac{1}{2}}(\partial\Omega)$.

The following proposition illuminates the connection between these different regularity notions present in the literature to the one used here.
\begin{prop}
Let $\olO$ be a manifold with boundary and bounded geometry. Let $D_e$ be a boundary condition for the uniformly elliptic operator $D$ over $\olO$.
\begin{itemize}
\item[(i)] Suppose that $D$ is first-order. Then, the following are equivalent:
\begin{itemize}
\item $D_e$ is (semi-)regular in the sense of Definition \ref{DefReg} above.
\item $D_e$ is $A$-(semi-)regular in the sense of \cite[Definition 6.8]{Baer2022}.
\end{itemize}
These equivalent conditions imply that $D_e$ is (semi-)regular in the sense of \cite[Definition 6.5]{Baer2022}.
\item[(ii)] Suppose that $D_e=D_P$ is defined by a differential operator $P$ on the boundary. If $(D,P)$ is a uniformly elliptic boundary condition in the sense of \cite[Definition 4.11]{Schick1998}, then $D_P$ is semi-regular.
\item[(iii)] Suppose that $D$ is a second-order differential operator in divergence form, and that $D_e=D_P$ is defined by a differential operator $P$ on the boundary. If $(D,P)$ has uniform regularity of the local problems in the sense of \cite[Theorem 5.6]{Grosse2020}, satisfies the uniform Shapiro-Lopatinski condition of \cite[Definition 6.8]{Grosse2020}, or is well-posed in energy space in the sense of \cite[Theorem 7.2]{Grosse2020}, then $D_P$ is semi-regular.
\end{itemize}
\end{prop}
\begin{proof}
Ad (i): We freely use the notation from \cite{Baer2022}. Under our assumption $A$ is a uniformly elliptic first-order operator on the bounded-geometry manifold $\dO$. Hence $\dom(A)=H^1(\dO;E|_{\dO})$. Then $\dom(|A|^{\frac{1}{2}})=H^{\frac{1}{2}}(\dO;E|_{\dO})$ as well.\footnote{It holds that $\dom(|A|^\mu)=[L^2(\dO;E),\dom(A)]_\mu=[L^2(\dO;E|_{\dO}),H^1(\dO;E|_{\dO})]_\mu=H^\mu(\dO;E|_{\dO})$ for $0\leq\mu\leq 1$, where we used Propositions \ref{InterpolationNoBdry} and \ref{FormalSobolevStructure} below, as well as the fact that $\dom(|A|^\mu)=\dom((A^*A+I)^{\mu/2})$.} Suppose that $D_e$ is $A$-semi-regular, meaning that the trace of $\dom(D_e)$ is contained in $\dom(|A|^{\frac{1}{2}})=H^{\frac{1}{2}}(\dO,E|_{\dO})$. By the trace theorem for Sobolev spaces this implies $\dom(D_e)\subseteq H^1(\olO;E)$. We have to argue that it is closed in the Sobolev norm. To that end note that $\|\cdot\|_D \lsim \|\cdot \|_{H^1}$. Then, if $(u_n)_{n\in\NN}$ is Cauchy with respect to the Sobolev norm, then it is also Cauchy with respect to the graph norm. Since $\dom(D_e)$ is closed with respect to the latter by definition, $(u_n)_n$ has a graph-norm limit $u\in \dom(D_e)$. On the other hand $H^1(\olO;E)$ is a Banach space, so that $(u_n)_{n}$ also has a Sobolev norm limit in $H^1(\olO;E)$. This second limit must coincide with the graph-norm limit $u\in\dom(D_e)$. We conclude that $\dom(D_e)$ is closed in the $H^1$-norm as well, hence semi-regular in our sense. Conversely, if $D_e$ is semi-regular in our sense, then $\dom(D_e)\subseteq H^1(\olO;E)$ has trace contained in $H^{\frac{1}{2}}(\dO;E|_{\dO})=\dom(|A|^{\frac{1}{2}})$, hence it is also $A$-semi-regular. Any $A$-semi-regular boundary condition is automatically semi-regular in the sense of \cite{Baer2022}. This is Corollary 6.10 in \cite{Baer2022}, but also follows directly from the inclusion $\dom(|A|^{\frac{1}{2}})=H^{\frac{1}{2}}(\dO;E|_{\dO})\subseteq H^{\frac{1}{2}}_{loc}(\dO;E|_{\dO})$ in our setting.

Ad (ii), (iii): In either case the stated condition is sufficient for a regularity estimate as in Proposition \ref{EstimateGivesReg}. In the case of (ii) this is \cite[Theorem 4.14]{Schick1998}. In the case of (iii) these are Theorems 5.6, 6.9, and 7.2 in \cite{Grosse2020}, respectively.
\end{proof}
If $D_e$ is semi-regular in the sense of \cite{Baer2022}, then $\dom(D_e)$ is contained in $\dom(D_{\max})\cap H^1_{loc}(\olO;E)$. There is no immediate reason to expect such a boundary condition to be contained in $H^1(\olO;E)$. Generally we can therefore not expect a boundary condition that is (semi-)regular in the sense of \cite{Baer2022} to be (semi-)regular in our sense as well.

Before discussing another more technical regularity condition let us briefly note that the norm of multiplication operators on the domains of semi-regular boundary conditions can be bounded by the $C^m_b$-norm.
\begin{lemma}\label{MultSemiReg}
Let $D_e$ be a semi-regular boundary condition, and let $f\in C^\infty_c(\olO)$ preserve the domain of $D_e$, i.e. $f\cdot\dom(D_e)\subseteq \dom(D_e)$. Then,
$$
\dom(D_e) \longrightarrow\dom(D_e) \; , \; u\longmapsto fu
$$ 
is a bounded map with norm bounded by a multiple of $\|f\|_{C^m_b(\olO)}$.
\end{lemma}
\begin{proof}
Since $D_e$ is semi-regular, its domain is a closed subspace of $H^m(\olO;E)$, and Sobolev and graph norm are equivalent. The statement thus follows from the fact that the corresponding norm estimate holds on Sobolev spaces, see Lemma \ref{MultSobolevDomain}.
\end{proof}
In our application to uniform K-homology we intend to use interpolation techniques. To that end we shall need another stronger notion of regularity that also provides uniform finite summability statements in intermediate degrees. Over rough domains we cannot expect Sobolev spaces to interpolate favorably, so that we are forced to impose this regularity in intermediate degrees manually.
\begin{definition}\label{DefSemiInterReg}
A boundary condition $D_e$ for $D$ is said to be \emph{semi-interpolation-regular} if for all $0\leq \mu\leq 1$ there is an extension map
$$
\dom(|D_e|^\mu) \xlongrightarrow{e_\mu} H^{m\mu}(M;E) \; , \; e_\mu(u)|_\Omega = u \; ,
$$
which is a Banach space isomorphism onto its image, such that whenever $\mu\geq\mu'$ it holds that $e_\mu(u)=e_{\mu'}(u) \in H^{m\mu'}(M;E)$ for all $H^{m\mu}(M;E)$.

Moreover, $D_e$ will be said to be \emph{interpolation-regular} if both $D_e$ and $D_e^*$ are semi-interpolation-regular.
\end{definition}
If the boundary of $\Omega$ is sufficiently smooth, then any semi-regular boundary condition is automatically semi-interpolation-regular, owing to the fact that Sobolev spaces interpolate favorably on sufficiently smooth domains.
\begin{prop}
Let $D$ be a uniformly elliptic operator of order $m$ over $\olO$. Any semi-interpolation-regular boundary condition for $D_e$ is also semi-regular. If $\Omega$ is a uniform $C^{m-1,1}$-domain, then the converse holds as well.
\end{prop}
\begin{proof}
Suppose that $D_e$ is semi-interpolation-regular. In particular, there is an extension map $e_1: \, \dom(D_e) \to H^m(M;E)$. Composing it with the restriction $H^m(M;E)\to H^m(\olO;E)$ produces a bounded inclusion $\dom(D_e)\subseteq H^m(\olO;E)$. Then, it holds that $\|u\|_D\lsim \|u\|_{H^m(\olO;E)}$ for all $u\in\dom(D_e)$, because $D$ is a bounded operator $H^m(\olO;E)\to L^2(\Omega;F)$. Put together we deduce that $\dom(D_e)\subseteq H^m(\olO;E)$ with equivalent norms, i.e. that $D_e$ is semi-regular. 

Suppose conversely that $D_e$ is semi-regular, and also that $\Omega$ is a uniform $C^{m-1,1}$-domain. The latter assumption implies that $[L^2(\Omega;E),H^m(\olO;E)]_\mu =H^{m\mu}(\olO;E)$. It also holds that $\dom(|D_e|^\mu)=[L^2(\olO;E),\dom(D_e)]_\mu$, see Proposition \ref{FormalSobolevStructure} below. Thus, if there is a continuous inclusion $\dom(D_e)\hookrightarrow H^m(\olO;E)$ with equivalent norms, then there is also a continuous inclusion
$$
\dom(|D_e|^\mu)=[L^2(\olO;E),\dom(D_e)]_\mu \hookrightarrow [L^2(\Omega;E),H^m(\olO;E)]_\mu =H^{m\mu}(\olO;E)
$$
with equivalent norms. Composing this inclusion map with the extension map $H^{m\mu}(\olO;E)\to H^{m\mu}(M;E)$ from Proposition \ref{ExtensionReflection} produces the desired extension map $\dom(|D_e|^\mu)\to H^{m\mu}(M;E)$. That it has the necessary properties follows from the fact that the extension map $H^{m\mu}(\olO;E)\to H^{m\mu}(M;E)$ has those properties. 
\end{proof}
Our notion of semi-interpolation-regularity therefore only becomes relevant on rough domains. However, we better show that even on rough domains there always exist such boundary conditions. 
\begin{prop}\label{MinSemiInterpol}
The minimal boundary condition $D_{\min}$ is semi-interpolation-regular on any domain $\Omega\subseteq M$.
\end{prop} 
\begin{proof}
It holds that $\dom(|D_{\min}|^\mu)=[L^2(\olO),\dom(D_{\min})]_\mu = [L^2(\olO),\tilde{H}^m(\Omega;E)]_\mu$ with equivalent norms. Fix some $\delta>0$. Then Proposition \ref{SubspaceInterpolationDelta} provides us with an inclusion $\dom(|D_{\min}|^\mu) \subseteq \tilde{H}^{m\mu}(U_\delta(\Omega);E)$ with equivalent norms. Since the latter space is in turn a subspace of $H^{m\mu}(M;E)$ this yields an extension map $\dom(|D_{\min}|^\mu)\to H^{m\mu}(M;E)$. Since it is the identity map on the set-level, the required compatibility with inclusions for smaller $\mu'$ holds trivially.
\end{proof}
\section{Examples}\label{SectionExamples}
The purpose of this section is to give examples of boundary conditions for different geometrically relevant operators, and discuss which of the above-discussed properties they possess, concretely if they are self-adjoint, local, (semi-)regular, odd or multigraded. The goal is by no means an exhaustive account but rather concrete illustrations. Thus, for our examples concerning Dirac operators we focus on spin Dirac operators, though the informed reader can no doubt extrapolate to spin$^c$ Dirac operators or more generally Dirac-type operators wherever appropriate (we refer also to \cite{Grosse2014} for a detailed treatment of boundary conditions for spin$^c$ Dirac operators). A summary of the examples can be found in Table \ref{TableExamples}.
\begin{remark}[Self-adjoint boundary conditions]
Before properly diving into the examples we wish to remark on self-adjointness. Coming from boundaryless manifolds, where every symmetric differential operator automatically has a unique self-adjoint extension (see e.g. \cite[Corollary 10.2.6]{HigsonRoe2000}), it is easy to take self-adjointness for granted. This attitude is misleading in the presence of boundaries. Indeed, there are symmetric differential operators that do not admit any self-adjoint extension! The most basic example of this is the operator $i\frac{d}{dx}$ on $[0,\infty)$, see \cite[Example 3.2 and Theorem 13.10]{Schmudgen2012}. Even if self-adjoint boundary conditions exist, self-adjointness may be incompatible with properties like locality and regularity. Corollary \ref{ObstructionRegular} below provides an obstruction for a uniformly elliptic operator to admit a self-adjoint, local and regular boundary condition. Since this obstruction is always non-zero for the Dirac operator, it can never admit such a boundary condition.
\end{remark}
\begin{example}[Minimal boundary condition]
Every (uniformly elliptic) differential operator admits a minimal boundary condition. If the operator is symmetric, then the adjoint of the minimal boundary condition is the maximal boundary condition. Minimal and maximal boundary condition are never equal unless the boundary is empty, hence neither is self-adjoint. Proposition \ref{SobolevDomain} above states that $\dom(D_{\min})=\tilde{H}^m(\Omega;E)$. It follows that $D_{\min}$ is semi-regular, and we have seen in Proposition \ref{MinSemiInterpol} above that it is even semi-interpolation-regular. Moreover, it is not hard to see that $\tilde{H}^m(\Omega;E)$ is closed under multiplication by functions from $C^\infty_c(\olO)$, meaning the minimal boundary condition is local. Lastly, if $D$ is odd resp. multigraded, these properties extend to the closure, i.e. the minimal boundary condition, so that $D_{\min}$ is also odd resp. multigraded in this case.
\end{example}
\begin{example}[Maximal boundary condition]\label{ExMinMax}
As discussed in the previous example the maximal boundary condition is generally not self-adjoint. It is also not regular unless $\olO$ is one-dimensional. For an explicit example consider the cylinder $\olO=\RR_{\geq 0}\times S^1$ with the uniformly elliptic operator $D=i \partial_\vartheta + \partial_y$. Here $\vartheta$ is the circle-variable and $y$ the height-variable. Fix some $u_0 \in L^2(S^1)$, $u\notin H^{\frac{1}{2}}(S^1)$, and define $\tilde{u}(\vartheta,y):= \left( e^{-y\partial_\vartheta^2}u_0 \right)(\vartheta)$. The function $\tilde{u}$ arises via heat-flow from $u_0$, thus is smooth for $y>0$. Smoothly cutting off $\tilde{u}$ at $y\geq 1$ results in an $L^2$-function $u$ on $\olO$ that is smooth for all $y>0$. Thus $u$ is in $\dom(D_{\max})$. If $u$ were also in $H^{1}(\olO)$, it would have a trace $u(\cdot,0)$ in $H^{\frac{1}{2}}(S^1)$. Since the heat kernel converges to the identity for $y\to 0$, this trace would coincide with $u_0$. But we have explicitly chosen $u_0$ not to lie in $H^{\frac{1}{2}}(S^1)$, making this a contradiction.

If $D$ is first-order, then $D_{\max}$ is local, because $D_{\min}$ is local, see Proposition \ref{LocalFirstOrder} above. For higher-order operators the maximal boundary condition need no longer be local. We give an example that is adapted from \cite[Remark 3.10]{Fries2025}. Consider the Laplace operator $\Delta$ on functions on the square $\Omega =\{0< x,y< 1\} \subseteq \RR^2$. The function $u(x,y)=\log(x^2+y^2)$ is smooth away from $0$ and square-integrable on $\Omega$. Moreover, it holds that $\Delta u=0$ away from $0$. Thus $u$ is in the kernel - hence in particular in the domain- of $\Delta_{\max}$. Consider also the smooth function $f(x,y):=\frac{1}{4}\left((1-x)^2-y^2 \right)$. An explicit calculation shows that
$$
\Delta(fu) = -\frac{2x}{x^2+y^2} + 2+ u \; .
$$
The function $x(x^2+y^2)^{-1}$ is not integrable over $\Omega$, hence neither is $\Delta(fu)$. Thus we have found $u\in \dom(\Delta_{\max})$ and $f\in C^\infty(\olO)$ such that $fu \notin \dom(\Delta_{\max})$.

Lastly, if $D$ is odd resp. multigraded, then so is $D_{\max}$. This follows from Lemma \ref{AdjointMultigraded} and the corresponding property of the minimal boundary condition.
\end{example}
\begin{example}[APS-boundary condition for the Dirac operator]
Assume that $\olO$ is a spin manifold with boundary and bounded geometry, and assume that $\olO$ has product form near the boundary. Let $\slashed{S}(\olO)\to \olO$ denote the complex spinor bundle, and $\Dirac: \, \Gamma_{cc}(\Omega;\slashed{S}(\olO))\to \Gamma_{cc}(\Omega;\slashed{S}(\olO))$ the spin Dirac operator. If the dimension of $\Omega$ is odd, then $\slashed{S}(\olO)|_{\partial \Omega}$ is naturally isomorphic to $\slashed{S}(\dO)$, and if the dimension is even, then it has two graded components, both of which are isomorphic to $\slashed{S}(\dO)$. Fix an inward-pointing unit normal $\nu$. On the boundary $\Dirac$ takes the form
$$
\Dirac = c(\nu) \left( \partial_\nu + A \right)  \; ,
$$
where $c(\nu)$ denotes Clifford multiplication by $\nu$, and $c(\nu)A=\Dirac_{\partial\Omega}$ if $\dim(\Omega)$ is odd, and $c(\nu)A=\begin{pmatrix}
0 & \Dirac_{\partial\Omega} \\ \Dirac_{\partial \Omega} & 0
\end{pmatrix}$ if it is even. In either case $A$ is an essentially self-adjoint operator, and thus admits a spectral projection $\pi_{(-\infty,0)}=\mathbf{1}_{(-\infty,0)}(A)$. The \emph{APS-boundary condition} is defined by
$$
\dom(\Dirac_\mathrm{APS}) = \left\{ u\in\dom(\Dirac_{\max})\, | \, \pi_{(-\infty,0)}u|_{\partial\Omega} = 0 \right\}
$$
The spectral projection $\pi_{(-\infty,0)}$ is a \emph{non-local} operator, therefore $\Dirac_\mathrm{APS}$ is not local. It is however regular \cite[Proposition 6.23]{Baer2022}. The domain of $\dom(\Dirac_\mathrm{APS}^*)$ equals $\dom(\Dirac_\mathrm{APS}) \oplus \ker(A)$, so that the APS-boundary condition is self-adjoint if and only if $A$ has trivial kernel, which in turn is the case if and only if $\Dirac_{\dO}$ has trivial kernel. Moreover, if $\olO$ is even-dimensional, then the spinor bundle is graded by the complex volume element, which locally can be written as $\kappa=\pm i c(e_1)\cdots c(e_{n-1}) c(\nu)$, where $e_1,\cdots,e_{n-1}$ form a local orthonormal frame of $\nu^\perp=T(\dO) \subseteq T\olO$. The operator $\Dirac$ is odd with respect to this grading. Since both $c(\nu)$ and the operator $\begin{pmatrix}
0 & \Dirac_{\partial\Omega} \\ \Dirac_{\partial \Omega} & 0
\end{pmatrix}$ are odd with respect to this grading, the operator $A$ is commutes with the grading operator $\kappa$. Hence the same holds for $\pi_{(-\infty,0)}$, and we conclude that $\kappa$ preserves $\dom(\Dirac_\mathrm{APS})$. Hence the APS-condition is odd in even dimensions. 
\end{example}
\begin{example}[Chiral boundary conditions for the Dirac operator]
We still consider the spin Dirac operator $\Dirac$ on a spin manifold with boundary and bounded geometry with inward-pointing unit normal $\nu$. Then, the Clifford multiplication by $i\nu$ is a self-adjoint involution on $\slashed{S}(\olO)|_{\partial \Omega}$. Let $P_\pm$ denote projection onto its $(\pm 1)$-eigenspace, and define the \emph{chirality boundary conditions} $\Dirac_\pm$ by
$$
\dom(\Dirac_\pm) = \left\{ u\in \dom(\Dirac_{\max}) \, | \, P_\pm u|_{\partial\Omega} = 0 \right\} \; .
$$
Since they are defined by bundle morphisms the chirality boundary conditions are local. It holds that $\Dirac_\pm^*= \Dirac_\mp$ \cite[Lemma 5.1]{Grosse2014}. The chirality boundary conditions are regular \cite[Proposition 6.32 and Corollary 6.10]{Baer2022}. However, they are not graded. Indeed, suppose that $\dim(\Omega)$ is even. Then the spinor bundle is graded by the complex volume element $\kappa$, which we recall is locally given by $\kappa=\pm i c(e_1)\cdots c(e_{n-1})c(\nu)$. Then, $\kappa c(\nu) = - c(\nu)\kappa$. Thus the grading operator $\kappa$ permutes the kernels of $P_\pm$, and we conclude that $\kappa \Dirac_\pm = -\Dirac_\mp \kappa$.
\end{example}
\begin{example}[Absolute and relative boundary condition for the Euler characteristic operator]
Assume that $\olO$ is a manifold with boundary and bounded geometry. We consider $d+d^*$ acting on the bundle $\bigwedge\olO=\bigoplus_{p=0}^{\dim(\Omega)} \bigwedge^p T^*\olO$, graded by even/odd degree. The resulting graded operator will be denoted by $D^\chi$. It is odd if $\olO$ is even-dimensional. Fix an inward-pointing unit normal $\nu$. Then, there is an orthogonal splitting
$$
\left(\bigwedge  T^*\olO\right)\bigg|_{\partial \Omega} = \left( \bigwedge T^*(\partial \Omega)  \right) \oplus \left( \bigwedge T^*(\partial \Omega)  \right) \wedge (\CC\cdot \nu^\flat)
$$
Denote by $P_1$ the projection onto the first summand, and by $P_2$ that onto the second. The principal symbol $\sigma_{D^\chi}(\nu^\flat)$ is the sum of the wedge product and the interior product by $\nu$. It therefore permutes the kernels of $P_1$ and $P_2$. On the other hand $\sigma_{D^\chi}(\xi)$ preserves the kernels of $P_1$, $P_2$ if $\xi \in T^*(\partial \Omega)$. Define the \emph{absolute} and \emph{relative} boundary conditions $D^\chi_\mathrm{abs}$ and $D^\chi_\mathrm{rel}$ by
\begin{align*}
\dom(D^\chi_\mathrm{abs}) &= \left\{ \alpha\in\dom(D^\chi_{\max}) \, | \, P_2\alpha = 0 \right\} \\
\dom(D^\chi_\mathrm{rel}) &= \left\{ \alpha\in\dom(D^\chi_{\max}) \, | \, P_1\alpha = 0 \right\} \; .
\end{align*}
They are defined by bundle morphism, and hence local. Moreover, they are both self-adjoint and regular.\footnote{Both self-adjointness and regularity follow from results in \cite{Baer2022}. All references in this footnote refer to the corresponding part of that paper. The adapted boundary operator is given by $A=\sigma_{D^\chi}(\nu)^{-1}D^\chi_{\partial \Omega}$. Then, the operator $\Xi=P_1 \oplus (-P_2)$ anti-commutes with $A$, and it squares to the identity. Thus $\Xi$ is a chirality operator for $A$ in the sense of Definition 6.31. Its $(\pm 1)$-eigenspaces are the kernels of $P_2$ and $P_1$, respectively. Proposition 6.32 and Corollary 6.10 combine to show that the absolute and relative boundary conditions are regular. Using that $\Xi$ is self-adjoint and anti-commutes with $\sigma_0= \sigma_{D^\chi}(\nu)$ we conclude self-adjointness from Remark 6.34(b).} The grading by degree preserves the domains of $D^\chi_\mathrm{abs/rel}$, hence $D^\chi_\mathrm{abs/rel}$ is odd in even dimensions.
\end{example}
\begin{table}[t]
\centering
\begin{tabular}{|c|c|c|c|c|c|c|}
\hline
 & self-adjoint & local & semi-reg. & reg. & odd & multigraded \\ \hline
 $D_{\min}$ & \xmark & \cmark & \cmark & \xmark & (\cmark) & (\cmark) \\  \hline
 $D_{\max}$ & \xmark & If $m=1$ & \xmark & \xmark & (\cmark) & (\cmark) \\  \hline
 $\Dirac_\mathrm{APS}$ & Iff $\ker(\Dirac_{\partial\Omega})=0$ & \xmark & \cmark & \cmark & (\cmark) &  \\  \hline
  $\Dirac_{\pm}$ & $\Dirac_\pm^*=\Dirac_\mp$ & \cmark & \cmark & \cmark & ($\kappa\Dirac_\pm=-\Dirac_\mp\kappa$) &  \\  \hline
  $D^\chi_\mathrm{abs}$ & \cmark & \cmark & \cmark & \cmark & If $\dim(\Omega)$ even &  \\  \hline
  $D^\chi_\mathrm{rel}$ & \cmark & \cmark & \cmark & \cmark & If $\dim(\Omega)$ even  &  \\  \hline
  $\Delta_\mathrm{Dir}$ & \cmark & \cmark & \cmark & \cmark &  &  \\  \hline
  $\Delta_\mathrm{Neu}$ & \cmark & \xmark & \cmark & \cmark &  &  \\  \hline
\end{tabular}
\caption{Summary of the properties of the boundary conditions discussed in this section. A check-mark in brackets indicates that the property is possessed if the prerequisites are met. For example the mark (\cmark) for $D_{\min}$ in the \glqq odd\grqq -column means that if $D$ is odd, then so is $D_{\min}$. Similarly, no mark is made if the grading property cannot meaningfully applied to the operator to begin with. For example the box for $\Delta_\mathrm{Dir}$ in the \glqq odd\grqq -column is left empty because the Laplacian is not a graded operator.}\label{TableExamples}
\end{table}
\begin{example}[Dirichlet and Neumann Laplacian]\label{ExNeumann}
Assume that $\olO$ is a manifold with boundary and bounded geometry. Consider the Laplace operator $\Delta=d^*d$ on functions. Again, we fix an inward-pointing unit normal $\nu$. Define the \emph{Dirichlet and Neumann boundary conditions} by
\begin{align*}
\dom(\Delta_\mathrm{Dir}) &= \left\{ u\in\dom(\Delta_{\max}) \, | \, u|_{\partial\Omega} = 0 \right\} \\
\dom(\Delta_\mathrm{Neu}) &= \left\{ u\in\dom(\Delta_{\max}) \, | \, (\partial_\nu u)|_{\partial\Omega} = 0 \right\} \; .
\end{align*}
Both are self-adjoint \cite[Theorem 1.2]{Ammann2019} and regular \cite[Proposition 5.9]{Schick1998}. The Dirichlet boundary condition is local in our sense, but the Neumann condition is not. Indeed, a function $f\in C^\infty_c(\olO)$ preserves $\dom(\Delta_\mathrm{Neu})$ if and only if $(\partial_\nu f)|_{\partial \Omega}=0$. This illustrates that differential boundary conditions need not be local in our sense. Indeed, if the differential operator defining the boundary condition has order at least $1$, there is no reason to expect the boundary condition to be local in this sense.
\end{example}
\section{Finite propagation}\label{SectionFinProp}
In this section we establish the finite-propagation property for wave operators and certain functions of local self-adjoint boundary conditions for uniformly elliptic differential operators of first-order. Using these results we will obtain finite-propagation representatives of the uniform K-homology classes of such operators in Chapter 6. While we know from Lemma \ref{RelTrunc} that any uniform K-homology class admits a finite-propagation representative, the truncation procedure that gives rise to that lemma is not exactly convenient to use in practice. In contrast the propagation results of this section will allow us to work with functions of boundary conditions directly.

We assume throughout that $D$ is a uniformly elliptic first-order differential operator. Since $D$ is in particular bounded the quantity
$$
c_D \; := \; \sup_{x\in M} \;  \sup_{\substack{\xi\in T^*_xM \\ ||\xi||=1}} \left|\left| \sigma_D(x,\xi) \right|\right|
$$
is finite. It is called the \emph{propagation speed} of $D$. The following proposition establishes that solutions of the wave equation propagate at most with speed $c_D$.
\begin{prop}\label{FinPropwaveOp}
Let $D_e$ be a local self-adjoint boundary condition for $D$. Then, for any compactly supported $u\in \dom(D_e)$ it holds that
$$
\supp(e^{itD_e}u)\subseteq N_{c_D|t|}\left( \supp(u) \right) \; .
$$
\end{prop}
The proof is exactly analogous to the situation without boundary. Since it might not be entirely standard that this result carries over to local boundary conditions, we repeat it here. We follow the proof provided in \cite{Willett2019}, Proposition 8.2.1.  
\begin{proof}
We may take $t$ to be positive, otherwise replace $D$ by $-D$. Set $K:=\supp(u)$ and $u_t:=e^{itD_e}u$. Let $\varepsilon>0$ and $\delta>0$, and choose $g\in C^\infty_c(\olO)$ such that
\begin{itemize}
\item[(i)] $g|_K \equiv 1$ and $g(x)=0$ if $d(x,K)\geq \delta$,
\item[(ii)] $\sup_{d(x,K)\leq \delta} \left| g(x)-\frac{1}{\delta}\left(1-d(x,K)\right)\right|<\varepsilon$,
\item[(iii)] $||dg||<\varepsilon+\frac{1}{\delta}$.
\end{itemize}
Choose $f:\, \RR \to [0,1]$ smooth and non-decreasing such that $f(t)=1$ if and only if $t\geq 1$. Furthermore, we set
\begin{align*}
\gamma \; &:= \; c_D(1+\delta\varepsilon) \\
h_t \; &:= \; f\left( g(x) +\frac{\gamma}{\delta}t\right) \; .
\end{align*}
Observe that $h_t(x)$ is smooth in $t$ and $x$, and that 
$$
h^{-1}_t\left(\{1\}\right) \; \subseteq \; N_{\gamma t + \varepsilon} (K) \; .
$$
Moreover, we may calculate the derivatives of $h_t(x)$:
\begin{align*}
\dot{h}_t(x) \; &= \; \frac{\gamma}{\delta} \cdot  f'\left( g(x)+\frac{\gamma}{\delta} t \right) \\
dh_t(x) \; &= \; f'\left( g(x)+\frac{\gamma}{\delta}t \right)\cdot dg(x) \; .
\end{align*}
Since $D_e$ is local, multiplication by $h_t$ preserves $\dom(D_e)$, so that we may form the commutator $[D_e,h_t]$. Then
$$
[D_e,h_t] = -i \sigma_D(dh_t) =  -\frac{i\delta}{\gamma} \dot{h}_t(x) \cdot \sigma_D(dg) \; .
$$
As per our assumptions on $g$ we have
$$
||\sigma_D(dg)|| \; \leq \; \frac{\gamma}{1+\delta \varepsilon} ||dg|| \; \leq \;  \frac{\gamma}{1+\delta \varepsilon} (\frac{1}{\delta}+\varepsilon) \; = \; \frac{\gamma}{\delta} \; ,
$$
from which we may deduce that the operator
$$
\dot{h}_t-i[D_e,h_t] \; = \; \dot{h}_t\left( I - \frac{\delta}{\gamma}\sigma_D(dg) \right)
$$
is positive, because $\dot{h}_t \geq 0$ as well. \\
By definition $u_t$ satisfies $\dot{u_t}=iD_eu_t$. Conclude that that
\begin{align*}
\frac{\partial}{\partial t} \langle h_tu_t,u_t \rangle \; &= \; \langle \dot{h}_t u_t,u_t \rangle + \langle h_t (iD_e) u_t,u_t \rangle + \langle h_tu_t, (iD_e)u_t\rangle \\
&= \; \langle \dot{h}_t - i[D_e,h_t]u_t,u_t \rangle \\
&\geq \; 0 \; .
\end{align*}
Moreover, since $e^{itD_e}$ is unitary, we have
$$
\langle u_t,u_t \rangle \; = \; \langle u_0,u_0\rangle \; = \; \langle h_0u_0,u_0 \rangle \; \leq \; \langle h_t u_t ,u_t \rangle \; . 
$$
On the other hand, Cauchy-Schwarz gives
$$
\langle h_t u_t ,u_t \rangle \; \leq \; ||h_tu_t||\cdot ||u_t|| \; \leq \; ||u_t||^2 \; ,
$$
where we used that $h_t\leq 1$. It follows that all inequalities in the above two displays are actually equalities. Because equality in Cauchy-Schwarz occurs if and only if the two vectors are linearly dependent, it follows that $h_tu_t = u_t$. We deduce that
$$
\supp(u_t) \; \subseteq \; h_t^{-1}\left( \{1\}\right) \; \subseteq \; N_{\gamma t + \varepsilon}(K) \; .
$$
Letting $\varepsilon$ tend to $0$ completes the proof.
\end{proof}
Proposition \ref{FinPropwaveOp} also holds for non-local self-adjoint boundary conditions if we assume that $u$ is supported in the interior, and that $c_D|t|\leq d(\supp(u),\dO)$. In this case $e^{itD_e}u$ is actually independent of the boundary condition $D_e$. Using this one may actually establish the existence of the K-homology class of a first-order elliptic operator over the interior of a manifold with boundary, and show that it is independent of the choice of boundary condition. This is the strategy used by Baum, Douglas and Taylor in \cite{BDT1989}. We will however employ a different strategy that also works for higher-order operators. 

For Proposition \ref{FinPropwaveOp} to hold up to the boundary it is however necessary that $D$ is first-order, and that $D_e$ is local. That the order of the operator matters can be seen like this: The wave operator $e^{itD}$ is the solution operator of the equation $\frac{\partial}{\partial t}u=iDu$. If $D$ is first-order, this is a hyperbolic equation, which have finite propagation speed. If $D$ is of higher order, then the equation is diffusive instead, and such equations generally have infinite propagation speed.

To see that some form of locality is necessary, consider $D=i\frac{d}{dx}$ on $M=[0,1]$, and $D_e$ the self-adjoint extension of $D$ determined by
$$
\dom(D_e)=\{ \, u\in H^1[0,1] \, | \, u(0)=u(1) \, \} \; .
$$
Then $(e^{itD_e}u)(x) = \tilde{u}(t-x)$, where $\tilde{u}$ is the $1$-periodic continuation of $u$ to a function on $\RR$. Thus, a function propagating into the left side of the interval immediately appears at the right side of the interval. Considering functions with supports getting closer to the left end of $[0,1]$, it takes less time for $e^{itD_e}$ to propagate them to the right side of the interval. Thus $e^{itD_e}$ has infinite propagation speed.

Using Fourier inversion we extend the finite-propagation result from wave operators to more general functions of $D_e$. Indeed, suppose that $f$ is a bounded Borel function with Fourier transform supported inside $[-R,R]$. Then we have the Fourier inversion formula
$$
f(D_e)u= \frac{1}{\sqrt{2\pi}} \int_{-R}^R \hat{f}(t) e^{itD_e}u \, dt
$$
for $u\in L^2(\Omega)$. Then Proposition \ref{FinPropwaveOp} implies the following.
\begin{cor}[Compare {\cite[Lemma 10.5.5]{HigsonRoe2000}}] \label{finPropFourier}
Let $D_e$ be a local self-adjoint boundary condition for $D$. Let $f$ be a bounded Borel function with $\supp(\hat{f})\subseteq [-R,R]$. Then, for any compactly supported $u\in \dom(D_e)$ it holds that
$$
\supp\left(f(D_e)u\right)\subseteq N_{c_DR}\left( \supp(u) \right) \; .
$$
\end{cor}
\chapter{Elliptic operators and uniform K-homology}
At this point all the technical machinery is in place, and we are ready to prove the existence and properties of the uniform K-homology classes of uniformly elliptic operators over arbitrary domains in manifolds of bounded geometry. 

To prove the existence of the uniform K-homology classes we use a formal operator calculus. This calculus is motivated and described further in the beginning of Section 6.1. That section contains the prerequisite definitions and facts concerning this calculus, as well as most calculations necessary to deduce existence of the uniform K-homology classes. The results of this section rely on the analysis done in the previous two chapters. As such Section 6.1 is the technical core of this chapter, and indeed the entire thesis. Existence of uniform K-homology classes is proved in Section 6.2, where we also discuss finite-propagation representatives that will be useful in our discussion of index maps later. In Section 6.3 we prove compatibility of the classes with restrictions and products, and we prove the boundary-of-Dirac-is-Dirac formula. Lastly, in Section 6.4 we use boundary-of-Dirac-is-Dirac to deduce bordism invariance of the relative uniform index, and combine it with the abstract results of Chapter 3 to deduce the uniform partitioned-manifold index theorems.
\section{The formal operator calculus}\label{SectionOpCalc}
Let $D_e$ be a boundary condition for a uniformly elliptic symmetric operator over a domain $\olO$, and let $F=D_e(D_e^*D_e+I)^{-\frac{1}{2}}$. We wish to show that $F$ defines a uniform Fredholm module, i.e. that 
$$
[f,F] \quad , \quad \left(I-F^2\right) f \quad , \quad \left(F-F^*\right) f
$$ 
are uniformly approximable for all relevant functions $f$. Recall from the introduction how this can be done in the case where $D_e=D$ is an elliptic operator over a closed manifold. We noted that in this case $(D^*D+I)^{-\frac{1}{2}}$ is an elliptic pseudo-differential operator of order $-m$ (with $m$ the order of $D$ as usual). Then, the commutator 
$$
[f,F]= [f,D](D^*D+I)^{-\frac{1}{2}}  + D\left[f,(D^*D+I)^{-\frac{1}{2}}\right]
$$
extends to a bounded map $L^2(\Omega)\to H^1(\Omega)$, because the commutator of a (pseudo-)differential operator with a smooth function is an operator of one order lower. Compactness then follows from the Rellich-Kondrachov theorem. Similarly, $I-F^2=(D^*D+I)^{-1}$ is a pseudo-differential operator of order $-2m$, hence extends to a bounded map $L^2(\Omega)\to H^{2m}(\Omega)$. Compactness again follows from the Rellich-Kondrachov theorem. Lastly, $F=F^*$ yields the compactness of the thrid operator trivially.

The used result that $(D^*D+I)^{-\frac{1}{2}}$ is an elliptic pseudo-differential operator of order $-m$ is already non-trivial for closed manifolds (see \cite{Seeley1967}). While interesting a discussion of functions of elliptic operators on possibly non-compact manifolds with boundary is miles beyond the scope of this thesis (though see \cite{Ammann2004a} for the case of boundaryless manifolds of bounded geometry). Therefore we adopt the formal argument, but instead of mapping properties of pseudo-differential operators on the Sobolev scale we use a formal scale $\cH^\alpha(D_e)$, $\alpha\in\RR$, of spaces associated to the fixed boundary condition $D_e$. In its properties this formal scale is very similar to the Sobolev scale, but its intrinsically operator-theoretic nature makes it behave much nicer under functional calculus. For example the statement that $(D_e^*D_e+I)^{-\frac{1}{2}}$ has negative order is true almost by definition on the formal scale of $D_e$. The formal scale of a closed operator and its basic mapping properties are introduced in Subsection 6.1.1.

To deduce the properties necessary for the existence of uniform K-homology classes we need uniform compactness statements for operators that have negative order on the formal scale. This will be achieved by connecting the formal scale of a boundary condition back to the Sobolev scale via regularity. This is the content of Subsection 6.1.2.

The uniform pseudo-locality of $F$ will require us to handle commutators of functions of $D_e$ with functions on $\olO$. We will do this by means of holomorphic functional calculus, as this provides a convenient tool to reduce the case of more general functions to resolvents. This is done in Subsection 6.1.3.

The uniform compactness of $\left( F-F^* \right)f$ was trivial in the absence of boundaries, but becomes a subtle matter in their presence: It holds that $F^*=D_e^*(D_eD_e^*+I)^{-\frac{1}{2}}$, so that $F^*$ is better interpreted as arising from the adjoint boundary condition $D_e^*$. Viewed in this light the uniform local compactness of $F-F^*$ becomes a comparison statement for different boundary conditions. This introduces technical complications, as the spaces $\cH^\alpha(D_e)$ will be different for different boundary conditions. The common ground is given by interior regularity; on the interior all boundary conditions are the same. This observation enables us make a more general statement: If $D_e$ and $D_e'$ are boundary conditions for \emph{different} uniformly elliptic operators, but we assume that $D$ and $D'$ differ only by a lower order operator -meaning they have the same principal symbol- then on the interior the domains of $D_e$ and $D_e'$ are both given by Sobolev space, and there $D_e-D_e'$ extends to an operator of lower order. From this observation we can use holomorphic functional calculus piece together the statement that $\left( D_e(D_e^*D_e+I)^{-\frac{1}{2}}- (D_e'((D_e')^*D_e' +I)^{-1/2} \right)f$ has negative order and thus is uniformly locally compact. This will imply that $D_e$ and $D_e'$ define the same class in uniform K-homology. We carry out this calculation in Subsection 6.1.4. Lastly, Subsection 6.1.5 contains a few remarks on extensions of these results to uniformly elliptic pseudo-differential operators, and elliptic operators with certain types of singularities on the boundary.

The use of the formal operator calculus as a more convenient substitute for pseudo-differential calculus and the Sobolev scale is well-established. Its use in index theory goes back at least to Connes and Moscovici \cite{Connes1995}, and is also used in the recent article by Fries \cite{Fries2025} extending the results of Baum, Douglas and Taylor \cite{BDT1989} to higher-order operators. For an example outside index theory, see for example the investigations of the Calderón projection in \cite{BoossBavnbek2009}.
\subsection{The abstract framework}
We begin with a general discussion of the scale associated to a fixed closed operator $T$ on a Hilbert space $H$. For $\alpha\geq0$ let $\cH^\alpha(T)$ denote the domain of $(T^*T+I)^{\frac{\alpha}{2}}$. Equip $\cH^\alpha(T)$ with the norm $||u||_{\cH^\alpha(T)} := ||(T^*T+I)^{\alpha/2} u||$, it is equivalent to the graph norm of $(T^*T+I)^{\frac{\alpha}{2}}$. Moreover, we set $\cH^{-\alpha}(T):= (\cH^\alpha(T))^*$. The spaces $\cH^\alpha(T)$ are Hilbert spaces, and if $\alpha\leq\alpha'$ there is a bounded inclusion $\cH^{\alpha'}(D_e)\hookrightarrow \cH^\alpha(D_e)$. Lastly, set $\cH^\infty(T) := \bigcap_{\alpha\in\RR}\cH^\alpha(T)$ and $\cH^{-\infty}(T):= \bigcup_{\alpha\in\RR}\cH^\alpha(T)$. When the operator $T$ is clear from context we shall take the liberty of writing $\cH^\alpha$, $\cH^\infty$, etc. instead of $\cH^\alpha(T)$, $\cH^\infty(T)$, etc.
\begin{prop}\label{FormalSobolevStructure}
\begin{itemize}
\item[(i)] For all $\alpha\geq 0$ the inner product on $H$ extends to a perfect pairing between $\cH^\alpha(T)$ and $ \cH^{-\alpha}(T)$.
\item[(ii)] For all $\alpha,\beta\in\RR$, $\alpha\leq \beta$, the embedding $\cH^\beta(T)\hookrightarrow\cH^\alpha(T)$ has dense image. In fact, $\cH^\infty(T)$ is dense in each $\cH^\alpha(T)$. 
\item[(iii)] For $\alpha \geq 0$, the norm of $u\in H$ as an element of $\cH^{-\alpha}(T)$ is given by 
$$
||u||_{\cH^{-\alpha}(T)} =  \left|\left| \left(T^*T+I\right)^{-\frac{\alpha}{2}}u \right|\right| \; .
$$
\item[(iv)] For all $\alpha < \beta$ the interpolation space $[\cH^\alpha(T),\cH^\beta(T)]_\theta$ is given by $\cH^{(1-\theta)\alpha+\theta\beta}(T)$.
\end{itemize}
\end{prop}
\begin{proof}
Item (i), (ii) and (iii) are proved in Appendix A.1 of Bär-Bandara. Item (iv) is proved in Proposition 4.2.2 of Taylor\footnote{Strictly speaking, in that proposition it is only proved that $[H,\cH^k]_\theta=\cH^{k\theta}$. This takes care of the case $\alpha\geq 0$.For $\beta <0$ note that via the dense embedding $H \hookrightarrow \cH^{\alpha}$, $\dom(T)=\dom((T^*T+I)^{1/2})$ is also the domain of the operator $T^{(\alpha)}:=(T^*T+I)^{\frac{1-\alpha}{2}}$ viewed as an unbounded operator on $\cH^\alpha$. In fact, since $(T^*T+I)^{-\alpha/2}$ is an isometric isomorphism from $H$ to $\cH^\alpha$, the operator $T^{(\alpha)}$ is closed, and it holds that $\cH^{\gamma}(T^{\alpha})=\cH^{(1-\alpha)\gamma+\alpha}(T)$. The case $\alpha <0$ may thus be reduced to the case $\alpha=0$.}. 
\end{proof}
Let $A: \, \cH^\infty(T)\to \cH^{-\infty}(T)$ be a linear map. The majority of our investigations will revolve around the question whether there is $\beta\in\RR$ such that $A$ extends to a bounded map $\cH^\alpha(T)\to \cH^{\alpha-\beta}(T)$ for some or even all $\alpha$. In the case where $T=D_e$ is a suitable boundary condition, such operators have similar mapping properties as pseudo-differential operators on the Sobolev scale, so we should think of them as operators of order $\beta$ on the scale associated to $T$. Our first lemma is an application of duality and interpolation. 
\begin{lemma}\label{DualityInterpolation}
Let $A,B: \, \cH^\infty \to H$ be linear maps such that $\langle u,Av \rangle_H= \langle Bu,v\rangle_H$. Let $\alpha,\beta\in \RR$. 
\begin{itemize}
\item[(i)] If $A$ extends to a bounded map from $\cH^\alpha$ to $\cH^{\alpha-\beta}$, then $B$ extends to a bounded map from $\cH^{-\alpha+\beta}$ to $\cH^{-\alpha}$.
\item[(ii)] If $A$ and $B$ extend to bounded maps $\cH^{\alpha} \to \cH^{\alpha-\beta}$, then $A$ and $B$ extend to bounded maps $\cH^{\gamma}\to \cH^{\gamma-\beta}$ for all $\gamma$ between $-\alpha+\beta$ and $\alpha$.
\end{itemize}
\end{lemma}
\begin{proof}
Ad (i): $\cH^{-\alpha}$ is the dual of $\cH^\alpha$, by definition if $\alpha\geq 0$ and by reflexivity if $\alpha<0$. Consider the dual operator $A'$ to $A: \, \cH^\alpha \to \cH^{\alpha-\beta}$, it maps $(\cH^{\alpha-\beta})^*=\cH^{-\alpha+\beta}$ to $(\cH^{\alpha})^*=\cH^{-\alpha}$. By Proposition \ref{FormalSobolevStructure} it is fully determined by the images of $\varphi_u \in (\cH^{\alpha})^*$, $\varphi_u(v)=\langle u,v\rangle_H$, for $u\in\cH^\infty\subseteq \cH^{\alpha-\beta}$. For $v\in \cH^\infty\subseteq \cH^{-\alpha+\beta}$ this image is given by
$$
\left(A'\varphi_u \right) v = \varphi_u(Av) = \langle u,Av \rangle_H = \langle Bu,v \rangle_H = \varphi_{Bu}(v) \; .
$$
Note that since $B$ extends to a bounded map $\cH^\alpha \to \cH^{\alpha-\beta}$, $\varphi_{Bu}$ is an element of the dual of $\cH^{-\alpha+\beta}$. Since $A'$ is bounded as the dual operator to a bounded operator, it follows that $B$ extends to a bounded operator $\cH^{-\alpha+\beta}\to \cH^{-\alpha}$ via $\varphi_u \mapsto \varphi_{Bu}$. 

Ad (ii): If $A,B$ both extend to bounded maps from $\cH^\alpha$ to $\cH^{\alpha-\beta}$, then both maps also extend to bounded maps from $\cH^{-\alpha+\beta}$ to $\cH^{-\alpha}$ by item (i). Therefore $A$ and $B$ also extend to bounded maps from $[\cH^{-\alpha+\beta},\cH^\alpha]_\theta$ to $[\cH^{-\alpha},\cH^{\alpha-\beta}]_\theta$ by interpolation, and by Proposition \ref{FormalSobolevStructure} these exhaust all intermediate exponents.
\end{proof}
Next we come to functions of the operators. Applying a function to the operator $T$ produces an operator of order $\beta$ on the scale of $T$, where the order $\beta$ depends on the growth of the function. Concretely, for $\beta \in\RR$ let $\cF^\beta$ denote the space of measurable functions $f$ on $\RR$ for which $(1+x^2)^{-\frac{\beta}{2}} f \in \cL^\infty(\RR)$. We equip $\cF^\beta$ with the norm $||f||_{\cF^\beta}:= ||(1+x^2)^{-\frac{\beta}{2}} f||_\infty$. Then, we have the following.
\begin{prop}\label{MappingPropertiesT}
Let $\alpha,\beta\in\RR$ and $f\in\cF^\beta$. 
\begin{itemize}
\item[(i)] $T$ is a bounded operator from $\cH^\alpha(T)$ to $\cH^{\alpha-1}(T^*)$, and $T^*$ a bounded operator from $\cH^\alpha(T^*)$ to $\cH^{\alpha-1}(T)$.
\item[(ii)] $f(|T|)$ is a bounded operator from $\cH^\alpha(T)$ to $\cH^{\alpha-\beta}(T)$, and $||f(|T|)||\leq||f||_{\cF^\beta}$ independently of $\alpha$.
\item[(iii)] Suppose $T$ is self-adjoint. Then, $f(T)$ is a bounded operator from $\cH^\alpha(T)$ to $\cH^{\alpha-\beta}(T)$, and $||f(T)||\leq ||f||_{\cF^\beta}$ independently of $\alpha$.
\end{itemize}
\end{prop}
\begin{proof}
Ad (i): Let $u\in \cH^\alpha(T)\cap \cH^1(T)$. Then,
$$
||Tu||_{\cH^{\alpha-1}(T^*)} = ||(TT^*+I)^{\frac{\alpha-1}{2}} Tu|| = || T(T^*T+I)^{\frac{\alpha-1}{2}} u || \leq || (T^*T+I)^{\frac{\alpha}{2}} u || \leq ||u||_{\cH^\alpha(T)} \; ,
$$
for any $\alpha\in\RR$. Using that $\cH^1$ is dense in $\cH^\alpha(T)$ for $\alpha\leq 1$ the claim for $T$ follows. For $T^*$ it is proved analogously.

Ad (ii): By definition of $\cF^\beta$ and the properties of the functional calculus the operator $f(|T|)$ can be written as the product of $(T^*T+I)^{\beta/2}$ with a bounded operator. For $u\in \cH^\infty(T)$ we have
$$
||f(|T|)u||_{\cH^{\alpha-\beta}} = \left|\left| \left((T^*T+I)^{-\beta/2} f(|T|) \right)\cdot (T^*T+I)^{\alpha/2} u \right|\right| \leq ||f||_{\cF^\beta}\cdot ||u||_{\cH^\alpha(T)} \; ,
$$
where we used commutativity and isometry of the functional calculus. Using that $\cH^\infty(T)$ is dense in $\cH^\alpha(T)$, the claim for $f(|T|)$ follows. The proof of (iii) is completely analogous.
\end{proof}
\subsection{The case of uniformly elliptic operators}\label{SubsecFormalElliptic}
In our application to uniform K-homology we shall need uniform compactness statements related to boundary conditions for uniformly elliptic operators. These must come from the corresponding statements for Sobolev spaces. We will handle calculations with boundary conditions with the formal operator calculus. Therefore we need to connect the scale associated to suitable boundary conditions to the Sobolev scale, so that we can then deduce uniform compactness statements based on behavior in the operator calculus. More concretely we will see how statements that operators have negative order in the abstract operator calculus of a boundary condition imply uniform compactness statements. 

Throughout this section we consider two situations simultaneously: The behavior of arbitrary boundary conditions in the interior, and that of semi-interpolation-regular boundary conditions up to the boundary. In either case we can connect the domain of the boundary condition to Sobolev space, in the former case by interior regularity and in the latter by assumed semi-regularity. In either case we will deduce the desired statements for the operator calculus from the corresponding statement for Sobolev spaces. The behavior of arbitrary boundary conditions on the interior will eventually give rise to the fact that any boundary condition defines a uniform K-homology class on the interior, while that of (local) semi-regular boundary conditions will yield the relative uniform K-homology classes associated to such boundary conditions.

Let $M$ be a manifold of bounded geometry, $\Omega \subseteq M$ a domain, and $E\to M$ vector bundles of bounded geometry. Let $D: \, \Gamma_{cc}(\Omega;E)\to\Gamma_{cc}(\Omega;E)$ be a differential operator of order $m$ that is uniformly elliptic over $\olO$. Let $D_e$ be a boundary condition for $D$. Then, $D_e$ is a closed operator on the Hilbert space $L^2(\Omega)$, and we will investigate its formal scale $\cH^\alpha(D_e)$. We begin with a lemma concerning multiplication by functions.
\begin{lemma}
Let $D_e$ be a boundary condition for $D$, and $f\in C^\infty_c(\olO)$. Suppose that either $D_e$ is arbitrary and $\supp(f)\subset\Omega$, or $D_e$ is semi-regular and $f\cdot\dom(D_e)\subseteq \dom(D_e)$. Then multiplication by $f$ is a bounded map 
$$
\cH^{\alpha}(D_e) \longrightarrow \cH^{\alpha}(D_e)
$$
for all $0\leq\alpha\leq 1$, with norm proportional to $\|f\|_{C^m_b(\olO)}$.
\end{lemma}
\begin{proof}
For $\alpha=1$ the statement holds by Corollary \ref{MultInt} or Lemma \ref{MultSemiReg}, for $\alpha=0$ it is true because $\cH^0(D_e)=L^2(\Omega;E)$. For intermediate values it follows by interpolation.
\end{proof}
The following lemma may appear to be an unassuming consequence of interior regularity, and while this is not incorrect, we point out that the proof uses the interpolation behavior of the Sobolev spaces $\tilde{H}^s(\Omega)$. Indeed, being able to make inferences for fractional degrees on the formal scale of boundary conditions was the motivation for our investigations into interpolation o Sobolev spaces in the first place.
\begin{lemma}\label{FormalInterpolationInterior}
Let $D_e$ be any boundary condition for $D$. Let $K\subset \Omega$ be a compact subspace in the interior. Then,
$$
\cH^\alpha(D_e) \cap L^2_K(\Omega;E) = H^{m\alpha}_K(M;E)
$$
with equivalent norms. The constants in the norm equivalence depend only on the distance of $K$ to $\partial\Omega$ but not the specific $K$. 
\end{lemma}
\begin{proof}
For $\alpha\in \ZZ_{\geq 0}$ a non-negative integer this is a consequence of interior regularity, because in these cases $\cH^\alpha(D_e)$ is the domain of a uniformly elliptic operator. Indeed, if $\alpha=2k$, $k\in\NN$, then $\cH^{\alpha}(D_e)=\dom((D^*_eD_e)^k)$. If $\alpha=2k+1$, $k\in \NN$, then $\cH^{\alpha}(D_e)=\dom(D_e(D_e^*D_e)^k)$. The claim then follows from a combination of Propositions \ref{MinMaxDomain} and \ref{SobolevDomain}.

For intermediate values of $\alpha$ we interpolate. Using a bump function around $K$ and interior regularity we find that 
\begin{align*}
\left[L^2(\Omega;E),\cH^{2k}(D_e)\right]_\theta \cap L^2_K(\Omega;E) &= \left[L^2(\Omega;E),\dom((D^\dagger D)_{\min})\right]_\theta \cap L^2_K(\Omega;E) \\
&= \left[L^2(\Omega;E),\tilde{H}^{2km}(\Omega;E)\right]_\theta \cap L^2_K(\Omega;E)
\end{align*}
with equivalent norms, where the constants in the norm equivalence depend on $K$ only through $d(K,\partial\Omega)$, which enters through the norm of the bump function. Then, fix $\delta>0$. By means of Proposition \ref{SubspaceInterpolationDelta} we conclude for $0<\theta<1$ that
\begin{align*}
\cH^{2k\theta}(D_e)\cap L^2_K(\olO;E) &= \left[L^2(\Omega;E),\cH^{2k}(D_e)\right]_\theta \cap L^2_K(\Omega;E) \\
&= \left[L^2(\Omega;E),\tilde{H}^{2mk}(\Omega;E)\right]_\theta \cap L^2_K(\Omega;E) \\
&\subseteq \tilde{H}^{2mk\theta}(U_\delta(\Omega);E) \cap L^2_K(\Omega;E) \\
&= H^{2mk\theta}_K(M;E)
\end{align*}
with equivalent norms, where the constants of the equivalence still depend on $K$ only through $d(K,\partial\Omega)$. This concludes the proof.
\end{proof}
In the following proposition we transfer the uniformly summable Rellich-Kondrachov theorem to the spaces $\cH^\alpha(D_e)$. This forms the technical backbone of our subsequent investigations into the uniform K-homology classes of uniformly elliptic operators. Again, we feel obligated to point out that this result is the derived from technical overhead developed in Chapters 4 and 5, thereby necessitating analysis developed in these chapters.
\begin{prop}
Let $K\subseteq \olO$ be a compact subset. Assume that either:
\begin{itemize}
\item[(a)] $K\subset \Omega$, and $D_e$ is arbitrary,
\item[(b)] or $K$ is arbitrary, and $D_e$ is semi-interpolation-regular.
\end{itemize}
Let $0<\mu\leq 1$. Set $R:=\diam(K)$, and $\delta:=d(K,\dO)$. The maps 
$$
\cH^\mu(D_e)\cap L^2_K(\Omega;E) \xlongrightarrow{j_\mu} L^2_K(\Omega;E) \quad \mathrm{and} \quad L^2_K(\Omega;E) \xlongrightarrow{j_{-\mu}} \cH^{-\mu}(D_e)
$$
are $p$-summable for all $p>\frac{n}{m\mu}$. In case (a) it holds that
$$
\|j_\mu\|_p, \, \| j_{-\mu}\|_p \leq C(\mu,p,R,\delta) \; .
$$
In case (b) it holds that
$$
\|j_\mu\|_p, \, \| j_{-\mu}\|_p \leq C'(\mu,p,R) \; .
$$
Here, $C$ and $C'$ depend only on $\olO$, $D_e$, and the indicated data, but not the specific $K$.
\end{prop}
\begin{proof}
We begin with analysis of the map $j_\mu$. In case (a) the claim follows from a combination of Lemma \ref{FormalInterpolationInterior} and the uniformly summable Rellich-Kondrachov theorem \ref{QuantComp}. In case (b) one instead uses the extension map $e_\mu: \, \cH^\mu(D_e)\to H^{m\mu}(M;E)$ of Definition \ref{DefSemiInterReg} to again reduce to the uniformly summable Rellich-Kondrachov theorem. The argument is completely analogous to the reduction to the boundaryless case in the proof of Theorem \ref{QuantComp}. 

We turn to the map $j_{-\mu}$. Summability and norm estimate for this map follow algebraically from that of $j_\mu$, so that we may treat both cases for $K$ and $D_e$ simultaneously. For $u,v\in L^2(\Omega)$ we have
$$
\langle u, (j_{-\mu})^* j_{-\mu} v\rangle_{L^2} = \langle j_{-\mu} u, j_{-\mu} v\rangle_{\cH^{-\mu}} = \langle u, (D_e^*D_e+I)^{-\mu} v \rangle_{L^2} \; ,
$$
so that $(j_{-\mu})^* j_{-\mu} = (D_e^*D_e+I)^{-\mu}$. The polar decomposition of $j_{-\mu}$ thus takes the form $j_{-\mu}=U (D_e^*D_e+I)^{-\frac{\mu}{2}}$ with $U: \, L^2(\Omega;E) \to \cH^{-\mu}(D_e)$ a partial isometry. Let $P_K=\mathbf{1}_K$ denote the projection $L^2(\Omega;E)\to L^2_K(\Omega;E)$. For $u\in L^2_K(\Omega;E)$ we may then write
$$
j_{-\mu} u = j_\mu P_K u = U (D_e^*D_e+I)^{-\frac{\mu}{2}} P_K u = U \left( (P_K j_{\mu})\cdot (D_e^*D_e+I)^{-\frac{\mu}{2}} \right)^* u \; .
$$
The map $P_K j_\mu: \, \cH^\mu(D_e) \to L_K^2(\Omega;E)$ coincides with the map
$$
\cH^\mu(D_e) \xlongrightarrow{P_K^\mu} \cH^\mu(D_e)\cap L^2_K(\Omega;E) \xlongrightarrow{j_\mu} L^2_K(\Omega;E) \; ,
$$
where $P_K^\mu$ denotes the orthogonal projection onto the closed subspace $\cH^\mu(D_e)\cap L^2_K(\Omega;E)$. Thus $P_Kj_\mu$ is $p$-summable for $p>\frac{n}{m\mu}$, and satisfies the corresponding norm estimate. It follows that the same is true for $j_{-\mu}: \, L^2_K(\Omega;E)\to \cH^{-\mu}(D_e)$.
\end{proof}
Using this adaptation of the uniformly summable Rellich-Kondrachov theorem we deduce that if an operator has negative order on the formal scale of a boundary condition, then it is uniformly locally compact.  
\begin{prop}\label{NegOrderUA}
Let $D_e$ be a boundary condition for $D$. Let $\mu>0$ and $k\in\NN$. Let $T\in \fB(L^2(\Omega;E))$, and $\mu>0$. 
\begin{itemize}
\item[(i)] If $Tf$ extends to a bounded map $\cH^{-\mu}(D_e)\to L^2(\Omega;E)$ for all $f\in C^\infty_{cc}(\Omega)$ with norm satisfying 
$$
\|Tf\| \leq \|f\|_{C^k_b(\olO)} \cdot C(\mu,R,\delta) 
$$
for some function $C$ of $\mu$, the diameter $R$ of $\supp(f)$, and the distance $\delta$ of $\supp(f)$ to $\partial \Omega$, then $Tf\sua 0$ for $f\in C_0(\Omega)$.
\item[(ii)] Suppose that $D_e$ is semi-interpolation-regular. If $Tf$ extends to a bounded map $\cH^{-\mu}(D_e)\to L^2(\Omega;E)$ for all $f\in C^\infty_{c}(\olO)$ with norm satisfying 
$$
\|Tf\| \leq \|f\|_{C^k_b(\olO)} \cdot C'(\mu,R)
$$
for some function $C'$ of $\mu$ and the diameter $R$ of $\supp(f)$, then $Tf\sua 0$ for $f\in C_0(\olO)$.
\end{itemize}
\end{prop}
\begin{proof}
We prove (i), the proof of (ii) is analogous. Fix $L,R\geq 0$ and $\delta>0$, and suppose $f\in \LLip_{R,\delta}(\Omega)$. Lemma \ref{BdBump} provides a function $\psi\in C^\infty_{cc}(\Omega)$ with $\psi=1$ on $K:=\supp(f)$ with $\supp(\psi)\subseteq N_{\delta/2}(K)$ and $\|\psi\|_{C^k_b(\olO)}$ bounded by some constant dependent only on $L$, $R$ and $\delta$. It follows that $T\psi$ extends to a bounded operator $\cH^{-\mu}(D_e)\to L^2(\Omega;E)$ with norm bounded by some constant also depending only on $L$, $R$ and $\delta$. Then $Tf=(T\psi)f$ factors as
$$
L^2(\Omega;E) \xlongrightarrow{f} L^2_K(\Omega;E) \xlongrightarrow{j_{-\mu}} \cH^{-\mu}(D_e) \xlongrightarrow{T\psi} L^2(\Omega;E) \; .
$$
The composition of the first two maps is $p$-summable for all sufficiently large $p$, hence so is the composition of all maps. Combining this with the relevant norm estimates we conclude that $Tf$ is $p$-summable for all $f\in \LLip_{R,\delta}(\Omega)$ with $p$-norm bounded by a uniform constant depending only on $L$, $R$ and $\delta$. By Lemmas \ref{LLipStratified} and \ref{UniformlySchattenClassUniformlyApprox} this suffices to prove the claim.
\end{proof}
We note that there is an analogous version of Proposition \ref{NegOrderUA} for operators $T$ with the property that $fT$ extends to a bounded map $L^2(\Omega;E)\to \cH^{\mu}(D_e)$ with an analogous norm estimate. Such operators then satisfy $fT\sua 0$. This remark also applies to the following corollary. Observe that the conditions of Proposition \ref{NegOrderUA} are in particular satisfied if $T$ itself extends to a bounded map $\cH^{-\mu}(D_e)\to L^2(\Omega;E)$ for some $\mu>0$. Thus we may draw the following corollary.
\begin{cor}\label{CorNegOrder}
Suppose $T\in \fB(L^2(\Omega;E))$ extends to a bounded map $\cH^{-\mu}(D_e)\to L^2(\Omega;E)$ for some $\mu >0$. Then $Tf\sua 0$ for $f\in C_0(\Omega)$. If $D_e$ is semi-interpolation-regular, then $Tf\sua 0$ for $f\in C_0(\olO)$.
\end{cor}
An analogous statement also holds for commutators: If the commutator of an operator with all functions has negative order, then that operator is uniformly pseudo-local.
\begin{prop}\label{CommLowerOrder}
Let $D_e$ be a boundary condition $D$. Let $\mu>0$ and $k\in\NN$. Let $T\in \fB(L^2(\Omega;E))$. 
\begin{itemize}
\item[(i)] If the commutators $[f,T]$ and $[f,T^*]$ extend to bounded maps $\cH^{-\mu}(D_e)\to L^2(\Omega;E)$ for all $f\in C^\infty_{cc}(\Omega)$ with norms satisfying
$$
\left\|[f,T]\right\|, \, \left\|[f,T^*]\right\| \leq \|f\|_{C^k_b(\olO)}\cdot C(\mu,R,\delta) \; ,
$$
for some function $C$ of $\mu$, the diameter $R$ of $\supp(f)$, and the distance $\delta$ of $\supp(f)$ to $\partial \Omega$, then $[f,T]\sua 0$ for $f\in C_0(\Omega)$.
\item[(ii)] Suppose that $D_e$ is semi-interpolation-regular. If the commutators $[f,T]$ and $[f,T^*]$ extend to bounded maps $\cH^{-\mu}(D_e)\to L^2(\Omega;E)$ for all $f\in C^\infty_{c}(\olO)$ with norms satisfying
$$
\left\|[f,T]\right\|, \, \left\|[f,T^*]\right\| \leq \|f\|_{C^k_b(\olO)}\cdot C'(\mu,R) \; ,
$$
for some function $C'$ of $\mu$ and the diameter $R$ of $\supp(f)$, then $[f,T]\sua 0$ for $f\in C_0(\olO)$.
\end{itemize}
\end{prop}
\begin{proof}
We prove (i), the proof of (ii) is again analogous. Fix $L,R\geq 0$ and $\delta>0$, and let $f\in \LLip_{R,k,\delta}$. Use Lemma \ref{BdBump} to produce a function $\psi \in C^\infty_{cc}(\Omega)$ such that $\psi=1$ on the support of $f$, $\supp(\psi)\subseteq N_{\delta/2}(\supp(\psi))$, and such that $\|\psi\|_{C^k_b(\olO)}$ is bounded by a constant depending only on $L$, $R$ and $\delta$. Write
$$
[f,T] = [\psi f, T] = \psi [f,T] + [\psi,T]f = -\left( [f,T^*] \bar{\psi} \right)^*+ [\psi,T]f \; .
$$
Both $[f,T^*]$ and $[\psi,T]$ extend to bounded operators $\cH^{-\mu}(D_e) \to L^2(\Omega;E)$ by assumption, and their norms are bounded by constants dependent only on $L$, $R$ and $\delta$. The claim then follows from Corollary \ref{CorNegOrder}. To be a bit more precise: The operator $[\psi,T]f$ factors as
$$
L^2(\Omega;E) \xlongrightarrow{f} L^2_K(\Omega;E) \xlongrightarrow{j_{-\mu}} \cH^{-\mu}(D_e) \xlongrightarrow{[\psi,T]} L^2(\Omega;E) \; ,
$$
with $K$ denoting the support of $f$. The middle map is $p$-summable for sufficiently large $p$ with $p$-norm dependent only on $R$, hence so is the entire composition, and the norm estimate for $[f,T]$ implies the $p$-norm to be bounded by a constant depending only on $L$, $R$ and $\delta$. The term $[f,T^*]\bar{\psi}$ admits an analogous conclusion. Proposition \ref{UniformlySchattenClassUniformlyApprox} then yields the claim. 
\end{proof}
\subsection{Commutators with functions}\label{SectionFormalPseudolocal}
In this section we investigate commutators of functions of a boundary condition $D_e$ with smooth functions on $\Omega$. The result will be analogous to the fact that the commutator of an operator of order $\beta$ with a smooth function has order $\beta -1$. In particular we will be able to conclude the uniform pseudo-locality of the bounded transform $D_e(D_e^*D_e+I)^{-\frac{1}{2}}$.

Holomorphic functional calculus is our tool of choice. Suppose that $T$ is a closed operator, $U\subseteq \CC$ an open neighborhood of the spectrum $\sigma(T)$, $\psi$ a holomorphic function on $U$, and $\Gamma$ a curve in $U\setminus \sigma(T)$ circumscribing $\sigma(T)$ in positive direction. Then the holomorphic functional calculus is defined as 
$$
\psi(T) = \frac{i}{2\pi} \int_\Gamma \psi(z) (T-z)^{-1} dz \; .
$$ 
If $T$ is self-adjoint, then it coincides with the Borel functional calculus. In this case the holomorphic functional calculus is a convenient way to obtain an integral representation in terms of resolvents, and this is the way in which we use it. Indeed, it essentially suffices to prove the desired uniform pseudo-compactness for resolvents, the statement for holomorphic functions then follows via holomorphic functional calculus. The use of integral representations to reduce pseudo-locality to the case of resolvents is a standard technique (see for example \cite{Kasparov1988} or \cite{BDT1989}). We use the holomorphic functional calculus here because it covers a reasonably large class of functions\footnote{Note that the holomorphic functional calculus can also be extended to functions with almost-analytic extensions (see for example \cite[Chapter 8]{Dimassi1999}), including for example compactly supported smooth functions, so that our results can be extended to this class of functions. We will not have need for this however.}, but the choice is ultimately arbitrary. 

As usual let $D_e$ be  boundary condition for the uniformly elliptic operator $D$ of order $m$. Let $f\in C^\infty_c(\olO)$ be a real-valued function. Throughout this section we shall assume that either:
\begin{itemize}
\item $D_e$ is arbitrary, and $f\in C^\infty_{cc}(\Omega)$, or
\item $D_e$ is local and semi-interpolation-regular, and $f$ is arbitrary.
\end{itemize}
Throughout, we also let $\mu:= \frac{1}{m}$.
\begin{lemma}\label{CommutatorBasicLemma}
The commutator $[f,D_e]$ extends to a bounded operator from $\cH^{1-\mu}(D_e)$ to $L^2(\Omega;E)$, and $[f,D_e^*D_e]$ extends to a bounded operator $\cH^{1-\mu}(D_e)\to \cH^{-1}(D_e)$. In either case the operator norm of the commutator has a bound proportional to $\|f\|_{C^m_b(\olO)}$, where the proportionality constant may depend on the distance of $\supp(f)$ to $\dO$ if $D_e$ is arbitrary and $f$ supported in the interior.
\end{lemma}
\begin{proof}
First consider the case that $D_e$ is arbitrary and $f\in C^\infty_{cc}(\Omega)$. It suffices to prove that $\|[f,D_e]u\|_{L^2} \lsim \|u\|_{\cH^{1-\mu}(D_e)}$ for all $u\in \cH^1(D_e)$. To that end note first that $[f,D_e]u=0$ if $u$ has support disjoint from that of $f$. Using a smooth cut-off we see that it suffices to prove the claim for $u$ supported in a $\delta$-neighborhood of the support of $f$, where $\delta$ can be chosen sufficiently small so that this neighborhood does not intersect $\partial\Omega$. By Lemma \ref{FormalInterpolationInterior} we have $\cH^\alpha(D_e)\cap L^2_{N_\delta(K)}(\Omega;E) \subseteq H^{m\alpha}_{N_\delta(K)}(M;E)$ with equivalent norms for all $\alpha\geq 0$, where the constants in the norm equivalence may depend on $\delta$, but not the specific $f$. Therefore it holds that
$$
\|[f,D_e]u\|_{L^2(\Omega;E)} \lsim \|u\|_{\cH^{1-\mu}(D_e)} \quad \Leftrightarrow \quad \|[f,D_e]u\|_{L^2(\Omega;E)}=\|[f,D]u\|_{L^2} \lsim \|u\|_{H^{m-1}(\olO;E)}
$$
for $u\in \cH^{1}(D_e)$. The latter estimate holds because $[f,D]$ is an operator of order $m-1$. This proves that $[f,D_e]$ extends to a bounded map $\cH^{1-\mu}(D_e)\to L^2(\Omega;E)$ in this case.

In the case where $D_e$ is local and semi-interpolation-regular the proof is similar. Pick an extension $\tilde{f}\in C^\infty_c(M)$ of $f$ such that $\|\tilde{f}\|_{C^m_b(M)}\leq 2\|f\|_{C^m_b(\olO)}$. Also, pick an extension $\tilde{D}$ of $D$ to a uniformly elliptic operator over $M$. Recall that by semi-interpolation-regularity there are bounded extension maps $e_\alpha: \, \cH^\alpha(D_e)\to H^{m\alpha}(M;E)$, $0\leq\alpha\leq 1$. Then, it holds that 
\begin{align*}
\|[f,D_e] u \|_{L^2(\Omega;E)} &\lsim \left\| [\tilde{f},\tilde{D}] e_1(u) \right\|_{L^2(M;E)} \\
&\lsim \|f\|_{C^m_b(\olO)}\cdot\|e_1(u)\|_{H^{m-1}(M;E)} \\
&\lsim \|f\|_{C^m_b(\olO)}\cdot \|u\|_{\cH^{1-\mu}(D_e)} \; ,
\end{align*}
where we have used that $e_1(u)=e_{1-\mu}(u)\in \cH^{1-\mu}(D_e)$, and that the map $e_{1-\mu}$ is a Banach space isomorphism onto its image. The shows that $[f,D_e]$ extends to a bounded map $\cH^{1-\mu}(D_e) \to L^2(\Omega;E)$ in this case as well. 

We deduce the mapping properties of $[f,D^*_eD_e]$ from the mapping properties on Sobolev spaces again. We prove it in the case where $D_e$ is local and semi-interpolation-regular. The other case is analogous. Note first that it suffices by Proposition \ref{FormalSobolevStructure} that
$$
\left| \left\langle [f,D_e^*D_e] u, v \right\rangle_{L^2(\Omega)} \right| \lsim \|f\|_{C^m_b(\olO)} \cdot \|u\|_{\cH^{1-\mu}(D_e)} \cdot \|v\|_{\cH^1(D_e)}
$$
for all $u,v \in \cH^1(D_e)$. One quickly computes that
$$
\left\langle [f,D_e^*D_e] u, v \right\rangle_{L^2(\Omega)} = \left\langle [f,D_e]u,D_ev \right\rangle_{L^2(\Omega)} - \left\langle D_eu,[f,D_e]v\right\rangle_{L^2(\Omega)} \; .
$$
The first summand satisfies the desired estimate because of the first part of the claim. For the second summand we again make use of the extension maps. With $\tilde{f}$ and $\tilde{D}$ as above one has
\begin{align*}
\left| \left\langle D_eu,[f,D_e]v\right\rangle_{L^2(\Omega)} \right|^2 &\leq \int_{\Omega}  \left|\langle D_eu,[f,D_e]v\rangle_{E}\right| d\mathrm{vol} \\
&\leq \int_{\Omega}  \left|\langle \tilde{D}e_{1}(u),[\tilde{f},\tilde{D}]e_1(v)\rangle_{E}\right| d\mathrm{vol} \\
&\leq \int_{M}  \left|\langle \tilde{D}e_{1}(u),[\tilde{f},\tilde{D}]e_1(v)\rangle_{E}\right| d\mathrm{vol} \\
&\lsim \|f\|_{C^m_b(\olO)}\cdot  \left\|  e_1(u) \right\|_{H^{m-1}(M)} \cdot \|e_1(v)\|_{H^m(M)} \; .
\end{align*}
To get to the last line we used that $[\tilde{f},\tilde{D}]$ is a bounded map $H^m(M;E)\to H^1(M;E)$ with norm bounded by a multiple of $\|\tilde{f}\|_{C^m_b(M)}\leq 2\|f\|_{C^m_b(\olO)}$, that $\tilde{D}$ is a bounded map $H^{m-1}(M;E)\to H^{-1}(M;E)$, and that the $L^2$-inner product is a perfect pairing between $H^{-1}(M;E)$ and $H^{1}(M;E)$. We conclude that
\begin{align*}
\left| \left\langle D_eu,[f,D_e]v\right\rangle_{L^2(\Omega)} \right|^2 &\lsim \|f\|_{C^m_b(\olO)}\cdot  \left\|  e_1(u) \right\|_{H^{m-1}(M)} \cdot \|e_1(v)\|_{H^m(M)} \\
&\lsim \|f\|_{C^m_b(\olO)}\cdot \|u\|_{\cH^{1-\mu}(D_e)}\cdot \|v\|_{\cH^1(D_e)} \; .
\end{align*}
Combination with the corresponding estimate for the first summand we deduce the claim for the commutator $[f,D_e^*D_e]$.
\end{proof}
From here on out the calculations are purely formal, and dependent only on the validity of the above lemma. We first deduce that the commutator of the resolvent with $f$ has lower order.
\begin{lemma}\label{CommutatorResolvent}
Let $z\in \CC\setminus [0,\infty)$ and let $\gamma\geq 0$. The commutator $\left[ f, \left(D_e^*D_e-z\right)^{-1} \right]$ extends to a bounded operator $\cH^{-1-\mu}(D_e)\to \cH^{1-\gamma}(D_e)$. Its norm satisfies
\begin{equation}\label{NormCommutatorResolvent}
\left\| \left[ f, \left(D_e^*D_e-z\right)^{-1} \right] \right\| \lsim \left\| f \right\|_{C^m_b(\olO)} \cdot |z|^{-\frac{\gamma}{2}} \; ,
\end{equation}
where the implicit constant in \eqref{NormCommutatorResolvent} may depend on the distance of $\supp(f)$ to $\dO$ if $D_e$ is arbitrary and $f$ supported in the interior.
\end{lemma}
\begin{proof}
The formula
$$
\left(D_e^*D_e-z\right)\left[ f, \left(D_e^*D_e-z\right)^{-1} \right] \left(D_e^*D_e-z\right) = \left[ D_e^*D_e,f \right]
$$
is straightforwardly verified to hold on $\cH^{1}(D_e)$. It follows that the left-hand side extends to a map $\cH^{1-\mu}(D_e) \to \cH^{-1}(D_e)$, or equivalently that $\left[ f, \left(D_e^*D_e+I\right)^{-1} \right]$ extends to a bounded operator $\cH^{-1-\mu}(D_e)\to \cH^{1}(D_e)$. Thus it also extends to a bounded operator $\cH^{-1-\mu}(D_e)\to \cH^{1-\gamma}(D_e)$ for all $\gamma\geq 0$. Moreover, the norm of $(D_e^*D_e-z)^{-1}$ as a bounded map $\cH^{\alpha}(D_e)\to \cH^{\alpha+2-\gamma}(D_e)$ is bounded by $\|(x^2-z)^{-1}\|_{\cF^{2-\gamma}}$ by Proposition \ref{MappingPropertiesT}. An elementary calculation shows this to be proportional to $|z|^{-\frac{\gamma}{2}}$. Together with the norm estimate for $[f,D^*_eD_e]$ this implies the estimate \eqref{NormCommutatorResolvent}.
\end{proof}
From resolvents we then extend to holomorphic functions via the holomorphic functional calculus. Let $\beta\in\RR$, and define $\cO^{-\beta}$ as the space of holomorphic functions $\psi: \, U\to \CC$ with $U\subseteq \CC$ an open neighborhhod of the positive real axis, such that $\sup_{z\in U} |z|^{\beta}|\psi(z)| <\infty$.
\begin{prop}\label{FormalPseudoLocality}
Let $\psi \in \cO^{-\beta}$, $0<\beta\leq 1$. Then $\left[  f, \psi(D_e^*D_e+I)\right]$ extends to a bounded operator $\cH^{-1-\mu}(D_e)\to \cH^{-1+2\beta-\gamma}(D_e)$ for all $\gamma>0$, with norm bound proportional to $\|f\|_{C^m_b(\olO)}$, where the proportionality constant may depend on the distance of $\supp(f)$ to $\dO$ if $D_e$ is arbitrary and $f$ supported in the interior.
\end{prop}
\begin{proof}
According to the holomorphic functional calculus we may write
$$
\psi(D_e^*D_e+I) = \frac{i}{2\pi} \int_\Gamma \psi(z)\cdot (D_e^*D_e+1-z)^{-1} \, dz \; ,
$$
so that 
$$
\left[f,\psi(D_e^*D_e+I)\right] = \frac{i}{2\pi} \int_\Gamma \psi(z)\cdot \left[f,(D_e^*D_e+1-z)^{-1}\right] \, dz$$
By Lemma \ref{CommutatorResolvent} $\left[ f, \left(D_e^*D_e-z\right)^{-1} \right]$ extends to a bounded operator $\cH^{-1-\mu}(D_e)\to \cH^{1-\gamma'}(D_e)$ for all $\gamma'\geq 0$, and its norm is bounded by $\left\| f\right\|_{C^m_b(\olO)} \cdot |z|^{-\frac{\gamma'}{2}}$. It follows that $\left[  f, \psi(D_e^*D_e+I)\right]$ extends to a bounded operator $\cH^{-1-\mu}(D_e)\to \cH^{1-\gamma'}(D_e)$ whenever the integral $\int_\Gamma \left|\psi(z)\right|\cdot |z|^{-\frac{\gamma'}{2}}\, dz$ converges, and that in this case the norm is bounded by $\|f\|_{C^m_b(\olO)}$ times the value of that integral. Since $\psi \in \cO^{-\beta}$, meaning that $|\psi(z)|\lsim |z|^{-\beta}$, this is the case if $-\beta-\frac{\gamma'}{2}<-1$, i.e. if $\gamma'>2(1-\beta)$. Defining $\gamma:=2(\beta-1)+\gamma'$ the claim follows.
\end{proof}
We cannot extend the result of Proposition \ref{FormalPseudoLocality} to cover $\gamma=0$. This is because an infinitessimal amount of $\gamma$ is needed to ensure that the integral representation of $[f,\psi(D_e^*D_e+I)]$ to converge in the operator norm of maps from $\cH^{-1-\mu}\to \cH^{-1+2\beta-\gamma}$. This is not a fundamental limitation, but rather an artifact of our choice to treat $\psi(D_e^*D_e+I)$ via an integral representation. However, using duality and interpolation we can extend the statement to more domains. Indeed, suppose that $\psi\in\cO^{-\beta}$. It holds that 
$$
\left\langle \left[f, \psi(D_e^*D_e+I)\right]u,v \right\rangle_H = \left\langle u,\left[\bar{\psi}(D_e^*D_e+I),f\right]v \right\rangle_H \; ,
$$
for all $u,v \in \cH^\infty(D_e)$. It (almost) never holds that $\bar{\psi}$ is in $\cO^{-\beta}$, but if $\psi$ is real-valued on $\RR_{\geq 0}$, then $\bar{\psi}(D_e^*D_e+I)=\psi(D_e^*D_e+I)$. In this case the conditions for Lemma \ref{DualityInterpolation} are satisfied, and Proposition \ref{FormalPseudoLocality} can be strengthened as follows.
\begin{cor}\label{FormalPseudoLocStrengthened}
Let $\psi \in \cO^{-\beta}$, $\beta\geq 0$, and suppose that $\psi$ is real-valued on $\RR_{\geq 0}$. Then $\left[  f, \psi(D_e^*D_e+I)\right]$ extends to a bounded operator $\cH^{\alpha}(D_e)\to \cH^{\alpha+2\beta+\mu-\gamma}(D_e)$ for all $\alpha\in[-1-\mu,1-2\beta]$ and all $\gamma>0$, with norm bound proportional to $\|f\|_{C^m_b(\olO)}$, where the proportionality constant may depend on the distance of $\supp(f)$ to $\dO$ if $D_e$ is arbitrary and $f$ supported in the interior.
\end{cor}
The additional range provided by the previous corollary allows us to prove that the the commutator of the bounded transform with $f$ has negative order.
\begin{cor}\label{BoundedTransLowerOrder}
The commutator $\left[f, D_e(D_e^*D_e+I)^{-\frac{1}{2}}\right]$ extends to a bounded operator $\cH^{-\gamma}(D_e)\to L^2(\Omega;E)$ for all $0<\gamma<\mu$, with norm bound proportional to $\|f\|_{C^m_b(\olO)}$, where the proportionality constant may depend on the distance of $\supp(f)$ to $\dO$ if $D_e$ is arbitrary and $f$ supported in the interior.
\end{cor}
\begin{proof}
Write
$$
\left[f, D_e(D_e^*D_e+I)^{-\frac{1}{2}}\right] = \left[ f,D_e\right](D_e^*D_e+I)^{-\frac{1}{2}} + D_e \left[ f, (D_e^*D_e+I)^{-\frac{1}{2}} \right] \; .
$$
The first summand extends to a bounded operator $\cH^{-\mu}(D_e)\to L^2(\Omega;E)$, with norm bounded by a multiple of $\|f\|_{C^m_b(\olO)}$. For the second summand note that $(z^2+1)^{-\frac{1}{2}} \in \cO^{-\frac{1}{2}}$ is real-valued for $z\in\RR_{\geq 0}$. The commutator $\left[ f, (D_e^*D_e+I)^{-\frac{1}{2}} \right]$ thus extends to a bounded operator $\cH^{-\gamma}(D_e)\to \cH^{-\gamma+1+\mu-\gamma'}(D_e)$ for all $0\leq\gamma\leq 1+\mu$ and all $\gamma'>0$ by the previous corollary, again with norm bounded by a multiple of $\|f\|_{C^m_b(\olO)}$. In either case the norm bound may depend on the distance of $\supp(f)$ to $\dO$ if $D_e$ is arbitrary and $f$ supported in the interior. Taking $0<\gamma'<\gamma<\mu$ we get that $\left[ f, (D_e^*D_e+I)^{-\frac{1}{2}} \right]$ extends to a bounded operator $\cH^{-\gamma}(D_e) \to \cH^1(D_e)$. Since $D_e$ maps $\cH^1(D_e)$ to $L^2(\Omega;E)$ the claim follows.
\end{proof}
Then, the appeal to Proposition \ref{CommLowerOrder} lets us conclude the desired uniform pseudo-locality.
\begin{cor}\label{BdTComm}
Let $D_e$ be a boundary condition for the uniformly elliptic operator $D$ over $\olO$.
\begin{itemize}
\item[(i)] $\left[ f,D_e(D_e^*D_e+I)^{-\frac{1}{2}} \right] \sua 0$ for $f\in C_0(\Omega)$.
\item[(ii)] If $D_e$ is local and semi-interpolation-regular, then $\left[ f,D_e(D_e^*D_e+I)^{-\frac{1}{2}} \right] \sua 0$ for $f\in C_0(\olO)$.
\end{itemize}
\end{cor}
\subsection{Lower-order perturbations}\label{SectionLowerOrderPerturbations}
We come to the comparison of boundary conditions for two uniformly elliptic operators that differ only by lower-order terms. In particular this includes the case of two boundary conditions for the same operator. Our strategy is analogous to that of the previous subsection: Use comparison to Sobolev spaces to argue that the difference between the operators has lower order on the formal scale, deduce that the analogous statement holds for resolvents as well, then extend to holomorphic functions via holomorphic functional calculus, and finally derive the desired result for the bounded transforms. This process is slightly complicated by the fact that the formal scales of the two boundary conditions are generally different, so that there is no common scale on which to work. This issue is solved by interior regularity. Thus the results of this subsection work only on the interior $\Omega$. While an extension to $\olO$ might be possible in favorable circumstances, we will not investigate this further. 

Throughout this section we consider the following situation: 
Let $D$ and $\tilde{D}$ be uniformly elliptic operators over $\olO$ acting on sections of the same bundles. Assume that $D$ and $\tilde{D}$ both have order $m$, and that $D-\tilde{D}$ is a differential operator of order $\leq m-1$. Let $D_e$ and $\tilde{D}_e$ be boundary conditions for $D$ and $\tilde{D}$ respectively. Again, we set $\mu:=\frac{1}{m}$. Lastly, let $f\in C^\infty_{cc}(\Omega)$ be a smooth function with compact support in the interior.

Let us start with the following observation. Interior regularity for the uniformly elliptic operators $D$ and $\tilde{D}$ implies that $f$ maps the domains of $D_e$ and $\tilde{D}_e$ to $\tilde{H}^m(\Omega;E) \subseteq \dom(D_e)\cap \dom(\tilde{D}_e)$, and the norm of $f$ is bounded by a multiple of $\|f\|_{C^m_b(\olO)}$ with regards to any of the equivalent norms to put on $\tilde{H}^m(\Omega;E) \subseteq \dom(D_e)\cap \dom(\tilde{D}_e)$. Similarly, $f$ maps the domains of $D^*_e D_e$ and $\tilde{D}^*_e \tilde{D}_e$ to $\tilde{H}^{2m}(\Omega;E) \subseteq \dom(D^*_e D_e)\cap\dom(\tilde{D}^*_e \tilde{D}_e)$ with norm bounded by a multiple of $\|f\|_{C^{2m}_b(\olO)}$. In either case the norm bound may depend on the distance of $\supp(f)$ to $\dO$. As a consequence of the latter, it makes sense to consider $(D_e^*D_e - \tilde{D}_e^*\tilde{D}_e)f$ as a bounded operator from $\cH^2(D_e)=\dom(D_e^*D_e)\to L^2(\Omega;E)$. Then, by arguing via Sobolev spaces in an analogous manner to the last subsection one proves the following.
\begin{lemma}
The operators $[f,D_e^*D_e]$ and $\left( D_e^*D_e - \tilde{D}_e^*\tilde{D}_e \right)f$ extend to bounded operators $\cH^{2-\mu}(D_e)\to L^2(\Omega;E)$, with norms bounded by $\|f\|_{C^{2m}_b(\olO)}$ and a constant depending on the distance of $\supp(f)$ to $\dO$.
\end{lemma} 
From this we again derive the corresponding statement for resolvents. Now we multiply the difference of the resolvents from both sides by functions with compact support in the interior. This produces an operator that acts on a single scale. 
\begin{lemma}
Let $f,g\in C^\infty_{cc}(\Omega)$, and let $\gamma\geq 0$, $z\in \CC\setminus[0,\infty)$. The operator 
$$
g\left( (D_e^*D_e-z)^{-1} - (\tilde{D}^*_e\tilde{D}_e-z)^{-1} \right)f
$$ 
extends to bounded operators $\cH^{-1-\mu}(D_e)\to \cH^{1-\gamma}(D_e)$, with norm bounded by $\|g\|_{C^{2m}_b(\olO)}\cdot \|f\|_{C^{2m}_b(\olO)}\cdot |z|^{-\frac{\gamma}{2}}$, and a constant depending on the distances of the supports of $f$ and $g$ to $\dO$.
\end{lemma}
\begin{proof}
Write $R:= D_e^*D_e-z$ and $\tilde{R}:=\tilde{D}_e^*\tilde{D}_e-z$. On the intersection $\dom(R)\cap \dom(\tilde{R})$ the formula
$$
\left( R^{-1}- \tilde{R}^{-1} \right) R = \tilde{R}^{-1} \left( \tilde{R}- R \right)
$$
holds. Since $f$ maps $\dom(R)$ to $\dom(R)\cap \dom(\tilde{R})$ we can thus calculate as follows:
\begin{align*}
\left( R^{-1} - \tilde{R}^{-1} \right)f &= \left( R^{-1} - \tilde{R}^{-1} \right) \left( RfR^{-1} + \left[f,R\right]R^{-1} \right) \\
&= \tilde{R}^{-1}\left(\left( \tilde{R}-R \right) f - \left[f,R\right]\right)R^{-1} + R^{-1}\left[f,R\right]R^{-1} \; .
\end{align*}
By the previous lemma the operator $\left(\left( \tilde{R}-R \right) f - \left[f,R\right]\right)R^{-1}$ extends to a bounded operator $\cH^{-\mu}(D_e)\to L^2(\Omega;E)$ with  norm bounded by $\|f\|_{C^{2m}_b(\olO)}$ and a constant dependent on the distance of $\supp(f)$ to $\dO$, so that the first summand extends to a bounded operator $\cH^{-\mu}(D_e)\to \cH^{2-\gamma}(\tilde{D}_e)$ with norm bounded by $\|f\|_{C^{2m}_b(\olO)} \cdot |z|^{-\frac{\gamma}{2}}$. The second summand extends to a bounded operator $\cH^{-\mu}(D_e)\to \cH^{2-\gamma}(D_e)$, again with an analogous norm bound. Now, multiplication by $g$ is a bounded map on $\cH^{2-\gamma}(D_e)$ with norm bounded by $\|g\|_{C^{2m}_b(\olO)}$ and a constant depending on the distance of $\supp(g)$ to $\dO$. We conclude that $g(R^{-1}-\tilde{R}^{-1})f$ extends to a bounded operator $\cH^{-\mu}(D_e)\to \cH^{2-\gamma}(D_e)$ with norm bounded by $\|g\|_{C^{2m}_b(\olO)} \|f\|_{C^{2m}_b(\olO)} |z|^{-\frac{\gamma}{2}}$ and a constant depending on the distance of the respective supports to $\dO$. Note that the same is also true for 
$$
\left( g\left(R^{-1}-\tilde{R}^{-1}\right)f \right)^* = \bar{f}\left(R^{-1}-\tilde{R}^{-1}\right)\bar{g} \; .
$$
Thus Lemma \ref{DualityInterpolation} implies that $g(R^{-1}-\tilde{R}^{-1})f$ actually extends to a bounded map $\cH^{\alpha}(D_e)\to \cH^{\alpha+\mu+2-\gamma}(D_e)$ for all $\alpha$ between $-2$ and $-\mu$ with an analogous norm bound. The claim is the special case of $\alpha=-1-\mu$.
\end{proof}
Analogously to Proposition \ref{FormalPseudoLocality} and Corollary \ref{FormalPseudoLocStrengthened} the result may again be extended to functions of $D_e^*D_e$ and $\tilde{D}_e^*\tilde{D}_e$ via holomorphic functional calculus. 
\begin{prop}\label{LowerOrderPsi}
Let $f,g\in C^\infty_{cc}(\Omega)$, and let $\psi \in \cO^{-\beta}$, $0\leq\beta\leq 1$, be real-valued on $\RR_{\geq 0}$. Then, the operator 
$$
g\left( \psi(D_e^*D_e+I) - \psi(\tilde{D}_e^*\tilde{D}_e+I) \right)f
$$
extends to bounded maps $\cH^{\alpha}(D_e) \to \cH^{\alpha+\mu+2\beta-\gamma}(D_e)$ and $\cH^{\alpha}(\tilde{D}_e) \to \cH^{\alpha+\mu+2\beta-\gamma}(\tilde{D}_e)$ for all $\alpha\in[-1-\mu,1-2\beta]$ and all $\gamma>0$, with norms bounded by $\|g\|_{C^{2m}_b(\olO)}\cdot \|f\|_{C^{2m}_b(\olO)}$, and a constant depending on the distances of the supports of $f$ and $g$ to $\dO$.
\end{prop}
Again, we use this to derive the uniform local compactness statement for the bounded-transforms that we are ultimately after.
\begin{cor}\label{BdTLowerOrderPerturbation}
$$
\left( D_e\left(D_e^*D_e+I\right)^{-\frac{1}{2}} - \tilde{D}_e\left(\tilde{D}_e^*\tilde{D}_e + I\right)^{-\frac{1}{2}} \right)f \sua 0 \quad \mathrm{for} \; f\in C_0(\Omega) \; .
$$
\end{cor}
\begin{proof}
Fix $L,R\geq 0$ and $\delta>0$, and let $f\in \LLip_{R,2m,\delta}(\Omega)$. Also, fix a bump function $\psi\in C^\infty_{cc}(\Omega)$ with $\psi=1$ on the support of $f$, $\supp(\psi)\subseteq N_{\delta/2}(\supp(f))$, and $\|\psi\|_{C^{2m}_b(\olO)}$ bounded by a constant dependent only on $L$, $R$ and $\delta$. Then, write
\begin{align*}
&\left( D_e\left(D_e^*D_e+I\right)^{-\frac{1}{2}} - \tilde{D}_e\left(\tilde{D}_e^*\tilde{D}_e + I\right)^{-\frac{1}{2}} \right)f \\
&= D_e f \left(D_e^*D_e+I\right)^{-\frac{1}{2}} - \tilde{D}_e f \left(\tilde{D}_e^*\tilde{D}_e + I\right)^{-\frac{1}{2}}\psi \\
&\quad - D_e\left[f,\left(D_e^*D_e+I\right)^{-\frac{1}{2}}\right] + \tilde{D}_e \left[ f,\left(\tilde{D}_e^*\tilde{D}_e + I\right)^{-\frac{1}{2}} \right]\psi \\
&= \left( D_e - \tilde{D}_e\right) f \left(D_e^*D_e+I\right)^{-\frac{1}{2}} + \tilde{D}_e f \left( \left(D_e^*D_e+I\right)^{-\frac{1}{2}} - \left(\tilde{D}_e^*\tilde{D}_e + I\right)^{-\frac{1}{2}} \right)\psi \\
&\quad - D_e\left[f,\left(D_e^*D_e+I\right)^{-\frac{1}{2}}\right] + \tilde{D}_e \left[ f,\left(\tilde{D}_e^*\tilde{D}_e + I\right)^{-\frac{1}{2}} \right]\psi
\end{align*}
The first summand extends to a bounded map $\cH^{-\mu}(D_e)\to L^2(\Omega;E)$, because $(D_e-\tilde{D}_e)f$ extends to a bounded map $\cH^{1-\mu}(D_e)\to L^2(\Omega;E)$. The third summand extends to a bounded map $\cH^{-\mu}(D_e)\to L^2(\Omega;E)$ by the results of the previous subsection. The fourth summand analogously extends to a bounded operator $\cH^{-\mu}(\tilde{D}_e)\to L^2(\Omega;E)$. The norm of each summand is bounded by $\|\psi\|_{C^{2m}_b(\olO)}\cdot \|f\|_{C^{2m}_b(\olO)}$, and a constant depending on $\delta$. It remains to treat the second summand. Proposition \ref{LowerOrderPsi} implies via an analogous argument to that in the proof of Corollary \ref{BoundedTransLowerOrder} that 
$$
f \left( \left(D_e^*D_e+I\right)^{-\frac{1}{2}} - \left(\tilde{D}_e^*\tilde{D}_e + I\right)^{-\frac{1}{2}} \right)\psi
$$
extends to a bounded operator $\cH^{-\gamma}(\tilde{D}_e)\to \cH^1(\tilde{D}_e)$ for $0<\gamma<\mu$, with norm again bounded by $\|\psi\|_{C^{2m}_b(\olO)}\cdot \|f\|_{C^{2m}_b(\olO)}$ and a $\delta$-dependent constant. Thus the second summand extends to a bounded map $\cH^{-\gamma}(\tilde{D}_e)\to L^2(\Omega;E)$ with an analogous norm bound. Now, $\|\psi\|_{C^{2m}_b(\olO)}\cdot \|f\|_{C^{2m}_b(\olO)}$ is bounded by a constant dependent only on $L$, $R$ and $\delta$. We conclude that each of the summands is $p$-summable for all $p>\frac{n}{m\mu}$ with $p$-norm bounded by a constant dependent only on $L$, $R$ and $\delta$. Proposition \ref{UniformlySchattenClassUniformlyApprox} thus yields the claim.
\end{proof}
\subsection{A note on pseudo-differential operators}\label{UPsiDOsSingular}
Pseudo-differential operators are not treated in this thesis, but let us remark upon them here briefly. Consider a domain $\Omega$ in a manifold of bounded geometry $M$. In line with Definition \ref{UEDOdomain} we could call a linear operator $P: \, \Gamma_{cc}(\Omega) \to \Gamma(\olO)$ a \emph{uniformly elliptic pseudo-differential operator of order $m$ over $\olO$} if there exists a uniform pseudo-differential operator\footnote{There are slightly different notions of pseudo-differential operators over manifolds of bounded geometry in use \cite{Kordyukov1991}, \cite{Shubin1992}, \cite{Taylor2008}, \cite{Engel2018}, and hence also of uniformly elliptic pseudo-differential operators. The author had the definitions used in \cite{Engel2018} in mind, but all that we need here is that any uniform pseudo-differential operator can be made to have arbitrarily small propagation by a uniformly smoothing perturbation (for example \cite[Remark 4.3]{Engel2018}), and that uniformly elliptic pseudo-differential operators satisfy an elliptic estimate (for example \cite[Theorem 6.7]{Engel2018}).} $\tilde{P}: \, \Gamma_c(M)\to \Gamma(M)$ of order $m$ such that $Pu=(\tilde{P}u)|_{\Omega}$ for all $u\in\Gamma_{cc}(\Omega)$. We do not claim that this is necessarily the correct notion for uniformly elliptic pseudo-differential operators on domains, but it is the one that is most immediately accessible by our techniques. Note that at least any uniformly elliptic differential operator over $\olO$ has a uniformly elliptic parametrix (i.e. uniformly pseudo-differential quasi-inverse) of this type, obtained by restricting a uniformly elliptic parametrix (existent for example by \cite[Theorem 6.5]{Engel2018}) for any uniformly elliptic extension $\tilde{P}: \, \Gamma_c(M)\to \Gamma_c(M)$ to $\Omega$. In any case, such an operator $P$ satisfies an interior regularity estimate of the form
$$
\|u\|_{H^s(M)} \leq C_{\delta,s} \left( \|Pu\|_{H^{s-m}(M)} + \|u\|_{H^{s-m}(M)} \right)
$$
for all $u\in\Gamma_{cc}(\Omega)$ with $d(\supp(u),\partial\Omega)\geq\delta$. This estimate can be derived from the elliptic estimate for any uniformly elliptic extension $\tilde{P}$ and the fact that $\tilde{P}$ and hence $P$ can be made to have propagation $<\delta$ by a uniformly smoothing perturbation. The constant $C_{\delta,s}$ generally depends on $\delta$, and there is no a priori reason why it should not blow up as $\delta\to 0$. In particular we cannot conclude that the minimal domain of $P$, again defined as the graph-norm closure of $\Gamma_{cc}(\Omega)$ in $L^2(\Omega)$, coincides with $\tilde{H}^m(\Omega)$. This $\delta$-stratified interior regularity estimate suffices to show that $\dom(P_{\min})\cap L^2_{\Omega\setminus U_\delta(\partial\Omega)}(\Omega)$ and $\tilde{H}^m(\Omega\setminus N_\delta(\partial\Omega))$ are equal with equivalent norms for every $\delta$. Using Lemma \ref{FormalInterpolationInterior} this is enough to extend the conclusions of this section regarding uniform compactness over the interior $\Omega$ to uniformly elliptic pseudo-differential operators over $\olO$. In principle the uniform compactness of certain commutators over $\olO$ could also be concluded for these operators, if we were willing to assume the existence of semi-interpolation-regular boundary condition. However, given that we cannot even say that the minimal boundary condition has this property the author sees no reason to expect such extension to exist.

In a similar way one could also embed elliptic differential operators exhibiting certain kinds of singularity on the boundary into our setting. Indeed, let $D: \, \Gamma_{cc}(\Omega) \to \Gamma_{cc}(\Omega)$ such that for every $\delta>0$ the operator $D$ is uniformly elliptic over $\olO\setminus U_\delta(\partial \Omega)$. Such an operator admits a $\delta$-stratified interior regularity estimate just as that for uniformly elliptic pseudo-differential operators above. Analogously the conclusion of this section regarding uniform compactness on the interior also hold for these operators. Examples include operators with singular coefficients, such as $x\partial_x + i\partial_y$ and $x^{-1}\partial_x + i\partial_y$ on the half-plane $\{x>0\}$, or operators that fail to admit extensions to $M$ for topological reasons, such as $e^{i\varphi/2}(\partial_x + i \partial_y)$, $\varphi$ being the argument of $(x,y)$, on the slit plane $\RR^2\setminus \RR_{\leq 0}$, viewed as a domain in $\RR^2$.
\section{The uniform K-homology classes of an elliptic operator}
With the tools developed in the previous section we are ready to prove the existence of the uniform K-homology classes associated to uniformly elliptic operators over domains. There is one more technicality to take care of first: the metric in use. Let $M$ be a manifold of bounded geometry, and $\Omega\subseteq M$ a domain. So far we have exclusively worked with the subspace metric on $\olO$. However, if $\olO$ is a manifold with boundary and bounded geometry, then the Riemannian metric induces a path-length metric on $\olO$. Especially when taking the intrinsic viewpoint on manifolds with boundary and bounded geometry this is the natural metric to put on $\olO$. The path-length metric is different from the subspace metric unless $\olO$ is geodesically convex in $M$. It is always true that the subspace metric is dominated by the path-length metric, but generally the two do not even have to be equivalent.\footnote{As an example consider $\Omega\subset \RR^2$ to be the complement of $\{|x|\leq 1\}\cup \{x_1\geq 0, \, |x_2|\leq 1\}$. Technically speaking the boundary $\dO$ is differentiable but not smooth in this example, but this does not effect the spirit of the example. Then, we consider the points $x_n,y_n\in \Omega$, $x_n=(n,2)$, $y_n=(n,-2)$. Then, the difference between $x_n$ and $y_n$ in the subspace metric induced by the standard metric on $\RR^2$ is constant in $n$. On the other hand any path in $\Omega$ connecting $x_n$ and $y_n$ must run around the origin, so that the distance between the two in the path-length metric of $\Omega$ diverges as $n\to \infty$. Thus the two metrics cannot be equivalent.} Since we would also like to state that uniformly elliptic operators define uniform Fredholm modules over $\olO$ equipped with the path-length metric, this requires further discussion.

There is an ad-hoc solution. Any Riemannian metric making $\olO$ into a manifold with boundary and bounded geometry is boundedly homotopic to one that has product structure near the boundary \cite[Proposition 7.3]{Schick1998}. Then, $\olO$ can be isometrically embedded as a geodesically convex submanifold by attaching cylinders with product metric to the boundary. In this case the subspace and path-length metric agree, and the path-length metric of the homotoped Riemannian metric is bi-Lipschitz equivalent to the original one.

Instead we can also note that due to bounded geometry the path-length and the subspace metric are equivalent on sufficiently small scales. Together with the knowledge that the subspace metric is dominated by the path-length metric this suffices to argue that the identity map provides a filterd map from $\olO$ with the former metric to $\olO$ with the latter. Thus, if we show that uniformly elliptic operators give rise to uniform Fredholm modules over $\olO$ with the subspace metric, then they automatically do so with the path-length metric as well. 

This and related statements are collected in the following lemma. Let us introduce some notation first. Denote by $d_{\olO}$ the path-length metric on $\olO$ induced by the restriction of the Riemannian metric of $M$ to $\olO$. Denote by $d_M$ the subspace metric on $\olO$ induced by the path length metric on $M$. Moreover, denote by $d_{\dO}$ the path-length metric on $\dO$ induced by the restriction of the Riemannian metric on $M$ to $\dO$. Further denote by $d_M$ and $d_{\olO}$ denote the subspace metrics on $\dO$ coming from the respective inclusions $\dO \subseteq M$ and $\dO \subseteq \olO$. Lastly, let $d_\Omega$ denote the path-length metric on $\Omega$, and $d_{\olO}$ the subspace metric coming from the inclusion $\Omega \subseteq \olO$.
\begin{lemma}\label{DifferentMetrics}
Let $M$ be a manifold of bounded geometry, and $\Omega\subseteq M$ a bounded-geometry domain.
\begin{itemize}
\item[(i)] The identity map $(\olO,d_M) \to (\olO,d_{\olO})$ is filtered. 
\item[(ii)] The identity maps $(\dO,d_M) \to (\dO,d_{\dO})$ and $(\dO,d_{\olO})\to (\dO,d_{\dO})$ are filtered.
\item[(iii)] The metrics $d_\Omega$ and $d_{\olO}$ are equal. In particular the identity map $(\Omega,d_M)\to (\Omega,d_\Omega)$ is also filtered.
\item[(iv)] Let $M'$ be another manifold of bounded geometry, and let $d_1=d_M+ d_{M'}$ denote the product metric, and $d_g$ the path-length metric induced by the product Riemannian metric $g_M\oplus g_{M'}$. Then the identity map $(M\times M',d_g)\to (M\times M',d_1)$ is a bi-Lipschitz equivalence.
\end{itemize}
\end{lemma}
\begin{proof}
We prove (i), the proof of (ii) is similar. Fix $L,R\geq 0$. It suffices to prove that there exist $L',R'\geq 0$ such that $\LLip_R(\olO,d_{\olO}) \subseteq L'\text{-}\mathrm{Lip}_{R'}(\olO,d_M)$. Let $f\in\LLip_R(\olO,d_{\olO})$. Because $d_M\leq d_{\olO}$ it holds that 
$$
\diam_{d_M}\left(\supp(f)\right) \leq \diam_{d_{\olO}}\left( \supp(f) \right) \leq R \; .
$$
For the Lipschitz property fix some small $\delta>0$, say smaller than the injectivity radii of $M$ and $\dO$ as well as half the thickness of a uniformly thick tubular neighborhood of $\dO$. If $d_M(x,y)\geq \delta$, then  
$$
\left| f(x) -f(y) \right| \leq 2\|f\|_\infty \leq \frac{2}{\delta} d_M(x,y) \; ,
$$
where we used that $f$ as the element of some $\LLip_R$ has $\|f\|_\infty\leq 1$. Now suppose that $d_M(x,y)<\delta$. Due to the bounded geometry of $M$ and $\olO$ there is a uniform bi-Lipschitz equivalence from $B_\delta^{d_M}(x)\cap \olO$ to the intersection of a Euclidean ball with half-space. The distance between $x$ and $y$ in this (potentially capped) ball is realized by the segment between them. Mapping this back to $\olO$ produces a curve $\gamma$ in $\olO$ connecting $x$ and $y$. Thus
$$
d_{\olO}(x,y) \leq \mathrm{length}(\gamma) \leq C\cdot d_M(x,y) \; ,
$$
where $C$ is a uniform constant coming from the uniform local bi-Lipschitz equivalence to half-space. We conclude that
$$
\left| f(x)-f(y) \right| \leq L\cdot d_{\olO}(x,y) \leq LC\cdot d_M(x,y) \; .
$$
Putting $L':=\max\{2\delta^{-1},LC\}$ we have found that $f\in L'\text{-}\mathrm{Lip}_R(\olO,d_M)$.

For (iii) simply note that every curve in $\olO$ between points in $\Omega$ can be approximated arbitrarily well by curves through $\Omega$. For (iv) note that $d_1$ is equivalent to $d_2=\sqrt{d_M^2 + d_{M'}^2}$. The metrics $d_2$ and $d_g$ can be shown to be equivalent by explicit inspection of curve lengths (see for example \cite{Nguyen2022}).
\end{proof}
It is a consequence of Lemma \ref{DifferentMetrics} (i) that if $(H,\rho,F)$ is a uniform Fredholm module over $(\olO,d_M)$, then it is also a uniform Fredholm module over $(\olO,d_{\olO})$. Its class in $K^{u}_*(\olO,d_{\olO})$ is the image of the class in $K^{u}_*(\olO,d_M)$ under the push-forward along the map $\id: \, (\olO,d_M)\to (\olO,d_{\olO})$. Analogous remarks apply to the remaining items of that lemma. Consequently, we shall always work with the subspace metric on any domain, and leave implicit that on bounded-geometry domains any resulting (relative) uniform K-homology class can be pushed forward to the path-length metric as well.
\begin{remark}
Item (i) of the Lemma \ref{DifferentMetrics} straight-forwardly generalizes to the case where $\olO$ is a uniform Lipschitz domain. In this case $\dO$ could be understood as a Lipschitz hypersurface of bounded geometry. Since Lipschitz maps are smooth almost everywhere, so is $\dO$. Thus we could make sense of an a.e. defined Riemannian metric on $\dO$, which is still induces a path-length metric. Then item (ii) carries over to this situation as well. We will not pursue the geometry on Lipschitz hypersurfaces in this thesis, though see Remark \ref{LipschitzBdryMap} below. However, the second part of item (iii) cannot hold with out at least some boundedness and regularity assumption on the boundary, consider $\Omega = \RR^2\setminus \RR_{\leq 0}$ in $M=\RR^2$.
\end{remark}
With the metric-related technicalities sorted out we are ready to state and prove our cornerstone theorem on the existence of uniform K-homology classes of uniformly elliptic operators over domains.
\begin{thm}\label{ExistencKHomClasses}
Let $M$ be a manifold of bounded geometry, and $\Omega\subseteq M$ a domain. Let $E\to M$ be a $p$-multigraded vector bundle of bounded geometry, $p\geq -1$. Let $D: \, \Gamma_{cc}(\Omega;E)\to \Gamma_{cc}(\Omega;E)$ be a symmetric and uniformly elliptic differential operator over $\olO$. Let $D_e$ be a boundary condition for $D$. If $p\geq 0$, assume that $D$ and $D_e$ are odd and multigraded. Set $H:= L^2(\Omega;E)$, and let $\rho$ the representation of either $C_0(\olO)$ or $C_0(\Omega)$ on $H$ by multiplication operators. Lastly, set 
$$
F_e:= D_e(D_e^*D_e+I)^{-\frac{1}{2}} \; .
$$
Then:
\begin{itemize}
\item[(i)] $(H,\rho,F_e)$ is a $p$-multigraded uniform Fredholm module over $\Omega$.
\item[(ii)] Assume that $D_e$ is local and semi-interpolation-regular. Then, $(H,\rho,F_e)$ is a $p$-multigraded relative uniform Fredholm module over $(\olO,\partial \Omega)$.
\item[(iii)] Assume that $D_e$ is local, self-adjoint and interpolation-regular. Then, $(H,\rho,F_e)$ is a $p$-multigraded uniform Fredholm module over $\olO$.
\end{itemize}
In each case $\Omega$ and $\olO$ is endowed with the subspace metric from the inclusions $\Omega,\olO\subseteq M$. The resulting uniform K-homology classes enjoy the following properties:
\begin{itemize}
\item[(a)] Assume that $D_e$ is local and semi-interpolation-regular. The excision map $\exc: \, K_{-p}^{u}(\olO;\partial \Omega) \to K^{u}_{-p}(\Omega)$ maps the class defined by (ii) to the class defined by (i).
\item[(b)] Assume that $D_e$ is local, self-adjoint and interpolation-regular. The restriction map $K^{u}_{-p}(\olO) \to K^{u}_{-p}(\olO;\partial \Omega)$ maps the class defined by (iii) to the class defined by (ii).
\item[(c)] The classes in $K^{u}_{-p}(\Omega)$ and $K^{u}_{-p}(\olO;\dO)$ defined by (i) and (ii), respectively, are independent of the boundary condition, and depends only on the principal symbol of $D$.
\item[(d)] Assume that the boundary condition $D_e$ is self-adjoint. Then, the conclusions of (i), (ii), and (iii) also hold with $D_e(D_e^*D_e+I)^{-\frac{1}{2}}$ replaced by $\chi(D_e)$ for any continuous odd function $\chi: \, \RR\to\RR$ with $\lim_{\lambda\to\infty}\chi(\lambda)=1$. The resulting uniform K-homology classes are independent of $\chi$.
\end{itemize}
\end{thm}
\begin{proof}
The properties necessary to establish (i), (ii) and (iii) were essentially proved over the course of Section \ref{SectionOpCalc}. Let us begin with (i). Then $[F_e,f]\sua 0$ for $f\in C_0(\Omega)$ by Corollary \ref{BdTComm}. Furthermore, $I- F_e^*F_e  = (D_e^*D_e+I)^{-1}$ extends to a bounded map $\cH^{-2}(D_e)\to L^2(\Omega;E)$, so that Corollary \ref{CorNegOrder} yields that $(I-F_e^*F_e)f\sua 0$ for $f\in C_0(\Omega)$. Now, suppose that $D_d$ is another boundary condition for $D$, and let $F_d:=D_d(D_d^*D_d+I)^{-1/2}$. It follows from Corollary \ref{BdTLowerOrderPerturbation} that $(F_e-F_d)f \sua 0$ for $f\in C_0(\Omega)$. Since $D$ is assumed symmetric, $D_d:=D_e^*$ is also a boundary condition for $D$, and this particular case yields $(F-F^*)f\sua 0$ for $f\in C_0(\Omega)$. Lastly, if $p\geq 0$ and $D_e$ is odd and multigraded, then $D^*_eD_e$ is even and multigraded, hence so is $(D_e^*D_e+I)^{-\frac{1}{2}}$. It follows that $F_e$ is odd and multigraded. This proves (i). Moreover, we have already seen that if $F_d$ is the bounded transform of a different extension $D_d$, then $(F_e-F_d)f\sua 0$ for $f\in C_0(\Omega)$. Hence $F_e$ and $F_d$ differ by a uniformly compact perturbation, and thus define the same class in $K^{u}_*(\olO,\dO)$. More generally, Corollary \ref{BdTLowerOrderPerturbation} yields that if $D$ and $\tilde{D}$ are two uniformly elliptic operators of order $m$ over $\olO$ that differ by an operator of order $\leq m-1$, then the Fredholm modules defined by any boundary conditions for $D$ and $\tilde{D}$ differ by a uniformly compact perturbation over $\Omega$, proving (c) in the case of absolute classes over the interior.

For (ii) it remains to show that $[f,F_e]\sua 0$ for $f\in C_0(\olO)$ if $D_e$ is local and semi-interpolation-regular. This was shown in Corollary \ref{BdTComm}. Note that (a) then holds because the excision map is the identity at the level of Fredholm modules. If $\tilde{D}_d$ is a local and semi-interpolation-regular boundary condition of another uniformly elliptic operator with the same principal symbol (or in particular another such boundary condition for $D$), Corollary \ref{BdTLowerOrderPerturbation} again gives that $F_e$ and $\tilde{F}_d$ differ by a uniformly compact perturbation. Hence the relative case of (c) follows.

For (iii) we have already noted that $[f,F_e]\sua 0$ for $f\in C_0(\olO)$, and $(F_e-F_e^*)f\sua 0$ holds trivially due to the assumed self-adjointness. It remains to show $(I-F_e^*F_e)f\sua 0$ for $f\in C_0(\olO)$. Since $I-F_e^*F_e=(D_e^*D_e+I)^{-1}$ extends to a bounded map $\cH^{-2}(D_e)\to L^2(\Omega;E)$, this is again implied by Corollary \ref{CorNegOrder} via (semi-)interpolation-regularity. This concludes the proof of (iii). Then (b) again holds because the restriction map $K^{u}_*(\olO)\to K_*^{u}(\olO,\dO)$ is the identity at the level of Fredholm modules.

It remains to prove (d). Note that the operator $F_e$ considered thus far corresponds to the normalizing function $\chi_0(x)=x(x^2+1)^{-\frac{1}{2}}$. By a density argument it suffices to consider normalizing functions $\chi$ that differ from $\chi_0$ by a function in $\cF^{-\beta}$ for some $\beta>0$. Then, $(\chi-\chi_0)(D_e)$ extends to a map $\cH^{-\beta}(D_e)\to L^2(\Omega;E)$ by Proposition \ref{MappingPropertiesT}. If $D_e$ is arbitrary, we conclude from Corollary \ref{CorNegOrder} that $\chi(D_e)=F_e+(\chi-\chi_0)(D_e)$ is uniformly pseudo-local over $\Omega$, and differs from $F_e$ by a uniformly compact perturbation. If $D_e$ is local and semi-interpolation-regular, we conclude from the same corollary that $\chi(D_e)$ is uniformly pseudo-local over $\olO$, and that it differs from $F_e$ by a uniformly compact perturbation over $\olO$. This proves (d), completing the proof of the theorem.
\end{proof}
Items (i) and (ii) of Theorem \ref{ExistencKHomClasses} in conjunction with property (c) provide well-defined classes
$$
[D] \in K^{u}_{-p}(\Omega) \quad \text{and} \quad [D] \in K^{u}_{-p}(\olO,\dO)
$$
defined as the uniform K-homology classes of any boundary condition with the required properties. To make sure that the latter class always exists, i.e. that any uniformly elliptic operator admits a local, semi-interpolation-regular and suitably graded boundary condition, it suffices to recall from Chapter 5 that the minimal boundary condition always has these properties.
\begin{example}[Minimal boundary condition]
The minimal boundary condition for any uniformly elliptic operator is local and semi-interpolation-regular, and possesses the same grading properties as $D$, and thus defines a relative uniform Fredholm module over $(\olO,\partial \Omega)$. Thus any uniformly elliptic operator $D$ over $\olO$ admits a class in $[D]\in K^{u}_*(\olO,\dO)$.
\end{example}
Let us consider some of the examples of Section \ref{SectionExamples}.
\begin{example}[APS- and chirality boundary conditions]\label{ExChirality1}
The APS-boundary conditions are self-adjoint (if $\ker(\Dirac_{\dO})=0$) and regular, but not local. Thus they define a uniform Fredholm module over $\Omega$, but not one over $\olO$. The chirality boundary conditions are local and regular, but not self-adjoint. In even dimensions they further lack the correct grading properties, but in odd dimensions this is no problem. In this case they define an ungraded relative uniform Fredholm module over $(\olO,\dO)$, but due to lack of self-adjointness we cannot conclude from Theorem \ref{ExistencKHomClasses} (iii) that they define a uniform Fredholm module over $\olO$ in this case. In fact, we will see in Example \ref{ExChirality2} below that they cannot define an ungraded uniform Fredholm module over $\olO$. Thus, while we do not expect self-adjointness to be generally required to get a uniform Fredholm module over $\olO$, some condition on the adjoint will be necessary. 
\end{example}
\begin{example}[Absolute and relative boundary conditions]
Both the absolute and the relative boundary conditions for the Euler characteristic operator are local, self-adjoint and regular, and thus define uniform Fredholm modules over $\olO$. On a closed connected even-dimensional manifold $M$ the index of the Euler characteristic operator $D^\chi$ equals the Euler characteristic $\chi(M)$. In fact it holds that
$$
[D^\chi] = \chi(M) \cdot i_*(\mathbf{1}) \; \in \; K_0(M) \; ,
$$ 
where $i_*(\mathbf{1})$ is the image of the generator of $K_0(\mathrm{pt})$ under the inclusion of a point \cite{Rosenberg1999}. Thus the K-homology class of the Euler-characteristic operator contains no information beyond its index. It would be interesting to know to what extent this is true in the bounded-geometry setting, and for manifolds with boundary. One would conjecture that the uniform coarse index of $D^\chi_{abs/rel}$ computes some sort of absolute/relative $L^2$-Euler characteristic, but of course both of these objects require correct interpretation before we can even ask whether they are equal. Let us mention the related result of Schick that the kernel of the absolute/relative boundary conditions compute the absolute/relative $L^2$-cohomology via bounded Hodge-theory \cite{Schick1998}.

In any case absolute and relative boundary condition in particular define K-homology classes over a compact manifold with boundary. Their indices are the absolute and relative Euler characteristic respectively, hence generally not equal. We conclude that even in the compact case $[D^\chi_{abs}]\neq [D^\chi_{rel}]$, so that Theorem \ref{ExistencKHomClasses} (c) is not generally true for the classes over $\olO$ provided by (iii).
\end{example}
\begin{example}[Dirichlet Laplacian]\label{DirichletTrivial}
The Dirichlet Laplacian is local, regular, and self-adjoint. Thus $\left( L^2(\olO),\rho, \Delta_\mathrm{Dir}(\Delta_\mathrm{Dir}^2+I)^{-1/2} \right)$ is an ungraded uniform Fredholm module over $\olO$. Unfortunately, its class in uniform K-homology is trivial. Indeed, let $F=\Delta_\mathrm{Dir}(\Delta_\mathrm{Dir}^2+I)^{-1/2}$ and note that $F$ is positive. Then, we see that
$$
F_t = \left( \sin(t) + \cos(t)F \right)\left( I + \cos(t)\sin(t) F \right)^{-\frac{1}{2}}
$$
is an operator-homotopy from $F$ to the degenerate Fredholm module $(L^2(\olO),\rho,I)$. More generally suppose that $D$ is any positive symmetric uniformly elliptic operator, and $D_e \geq 0$ a positive self-adjoint boundary condition. Then, the class $[D_e] \in K^{u}_1(\Omega)$ vanishes. If $D_e$ is also local and regular, then the class $[D_e]\in K^{u}_1(\olO)$ vanishes as well.
\end{example}
\begin{remark}
Following the discussion of Subsection \ref{UPsiDOsSingular} one can extend Conclusion (i) of Theorem \ref{ExistencKHomClasses} to uniformly pseudo-differential operators over $\olO$, and to elliptic differential operators over $\Omega$ exhibiting uniformly controlled singularities on the boundary. In principle, conclusions (ii) and (iii) could also hold for these operators provided one assumes the existence of (semi-)regular boundary conditions, but as discussed in that subsection this seems unlikely.
\end{remark}
In their work \cite{BDT1989} on K-homology classes for elliptic operators on manifolds with boundary Baum, Douglas and Taylor approach the case of relative classes slightly differently. They consider graded odd operators $D=\begin{pmatrix}
0 & D^- \\ D^+ & 0
\end{pmatrix}$, and instead of looking at semi-regular boundary condition for $D$ directly -which is what we did above- they consider semi-regular boundary conditions $D_e^+$ for $D^+$, and form the boundary condition
$$
D_e = \begin{pmatrix}
 0 & (D^+_e)^* \\ D_e^+ & 0
\end{pmatrix} \; .
$$
It is not semi-regular unless $D^+_e$ is already regular. They then show that the boundary condition $D_e$ defines a relative Fredholm module. We can do the same in our setting.
\begin{prop}\label{KHomGradedComp}
Let $D$ be a graded uniformly elliptic operator, and let $D^\pm$ denote the graded components of $D$. Let $D_e^+$ be a local and semi-interpolation-regular boundary condition for $D^+$, and set 
$$
D_e = \begin{pmatrix}
 0 & (D^+_e)^* \\ D_e^+ & 0
\end{pmatrix} \; .
$$
Then, $\left( L^2(\Omega;E),\rho,\chi(D_e)\right)$ is a graded relative uniform Fredholm module over $(\olO,\partial \Omega)$. Its class in $K_0^{u}(\olO,\partial\Omega)$ coincides with that defined in Theorem \ref{ExistencKHomClasses}. In particular it is independent of the choice of normalizing function.
\end{prop}
We need a lemma.
\begin{lemma}\label{SemiRegGradedComp}
Let $D$ and $D_e$ be as above. Then, for any odd function $C_0(\RR)$ with $\varphi(0)=0$ it holds that $\varphi(D_e)\in \fD^{u}(\olO,\partial\Omega)$.
\end{lemma}
\begin{proof}
By way of a density argument it suffices to consider the case where $\varphi(x)=x(x^2+\lambda)^{-\beta}$ for $\lambda >0$ and $\beta>\frac{1}{2}$. In this case
$$
\varphi(D_e) =  \begin{pmatrix}
 0 & \left( D_e^+((D_e^+)^*D^+_e+\lambda)^{-\beta}\right)^* \\  D_e^+((D_e^+)^*D^+_e+\lambda)^{-\beta}
\end{pmatrix} \; ,
$$
so that it suffices to prove that $D_e^+((D_e^+)^*D^+_e+\lambda)^{-\beta} \in \fD^{u}(\olO,\partial\Omega)$. To see this it suffices to note that $(x^2+\lambda)^{-\beta} \in \cF^{-2\beta}$, so that $D_e^+((D_e^+)^*D^+_e+\lambda)^{-\beta}$ extends to a bounded operator $\cH^{-2\beta+1}(D_e^+)\to \cH^0((D_e^+)^*)=L^2(\Omega,E^-)\subseteq L^2(\Omega;E)$. Then Proposition \ref{NegOrderUA} implies that $D_e(D_e^*D_e+\lambda)^{-\beta} \in \fD^{u}(\olO,\partial\Omega)$, completing the proof.
\end{proof}
\begin{proof}[Proof of Prop. \ref{KHomGradedComp}]
Since $D_e$ is an odd boundary condition $D$ the triple $(L^2(\Omega;E),\rho,\chi(D_e))$ is a graded uniform Fredholm module over $\Omega$ by Theorem \ref{ExistencKHomClasses}, (i). Thus it remains to show that $\chi(D_e)$ is also in uniformly pseudo-local over $\olO$. First, consider $\chi_0(x)=x(x^2+1)^{-1/2}$. Then, $\chi_0(D_e) = \begin{pmatrix}
0 & F^*\\ F& 0
\end{pmatrix}$ with $F=D_e^+((D_e^+)^*D^+_e+I)^{-1/2}$. Since $D^+_e$ is local and semi-interpolation-regular, $F$ and hence also $F^*$ are in $\fD^{u}(\olO)$. Thus $\chi_0(D_e)\in\fD^{u}(\olO)$ as well. Now, an arbitrary normalizing function $\chi$ differs from $\chi_0$ by an odd function vanishing both at infinity and at $0$. Thus Lemma \ref{SemiRegGradedComp} implies that $\chi(D_e) \in\fD^{u}(\olO,\partial\Omega)$ for arbitrary $\chi$. This proves that $(L^2(\Omega;E),\rho,\chi(D_e))$ is a graded relative uniform Fredholm module over $(\olO,\partial\Omega)$. That it defines the same relative uniform K-homology class as that constructed in Theorem \ref{ExistencKHomClasses} (ii), follows from excision.
\end{proof}
Property (d) of Theorem \ref{ExistencKHomClasses} allows for the representation of the uniform K-homology classes of uniformly elliptic operators by means of functional calculus if the boundary condition is self-adjoint. However, self-adjoint extensions are potentially sparse, and sometimes it is technically convenient to have a representation via functional calculus nonetheless. This is what the following lemma provides.
\begin{lemma}\label{KHomDilation}
Let $D_e$ be a semi-interpolation-regular, odd and $p$-multigraded boundary condition for the odd and $p$-multigraded uniformly elliptic operator $D$. Let $\hat{D_e}$ be the self-adjoint operator
$$
\hat{D_e} = \begin{pmatrix}
0 & D_e^* \\ D_e & 0
\end{pmatrix}
$$ 
on the Hilbert space $\hat{H}:=L^2(\Omega;E)\oplus L^2(\Omega;E)$. Let $P: \, \hat{H}\to L^2(\Omega;E)$ be the projection onto the first summand. Then, for any normalizing function $\chi$ the triple 
$$
\left(L^2(\Omega;E),\rho,(I-P)\chi(\hat{D}_e)P\right)
$$ 
is a $p$-multigraded relative uniform Fredholm module that differs from 
$$
\left(L^2(\Omega;E),\rho,D_e(D_e^*D_e+I)^{-\frac{1}{2}}\right)
$$ 
by a uniformly compact perturbation. In particular
$$
\left[ L^2(\Omega;E),\rho, D_e(D_e^*D_e+I)^{-1/2} \right] = \left[ L^2(\Omega;E),\rho, (I-P) \chi(\hat{D}_e) P  \right] \; \in \; K^{u}_{-p}(\olO,\partial\Omega) \; .
$$
\end{lemma}
\begin{proof}
Set $H=L^2(\Omega;E)$ so that $\hat{H}:=H\oplus H$. The space $\hat{H}$ is graded by $\hat{\kappa}=\kappa\oplus\kappa$ and admits a $p$-multigrading with multigrading operators $\hat{\varepsilon}_j:=\varepsilon_j\oplus \varepsilon_j$, $j=1,\cdots,p$.
The operator $\hat{D}_e$ is odd and $p$-multigraded. If $\chi$ is any normalizing function, then $\chi(\hat{D_e})$ is odd and $p$-multigraded as well. In particular it takes the form $\chi(\hat{D_e})=\begin{pmatrix}
0 & F^* \\ F & 0
\end{pmatrix}$. Consider the projection $P$ onto the first summand $H\subseteq \hat{H}$. It is even and $p$-multigraded as a map $\hat{H}\to H$. It follows that $F=(I-P)\chi(\hat{D}_e)P \in\fBH$ is odd and $p$-multigraded. By Proposition \ref{KHomGradedComp} $(\hat{H},\rho\oplus\rho,\chi(\hat{D}_e))$ is a relative uniform Fredholm module over $(\olO,\dO)$, and thus the same holds for $(H,\rho,F)$. This shows that $(H,\rho,F)$ is a $p$-multigraded relative uniform Fredholm module.

By Lemma \ref{SemiRegGradedComp} the operator $\chi(\hat{D}_e)$ is independent of the choice of $\chi$ up to a uniformly compact perturbation, thus the same is true for $F$. Choosing $\chi(x)=x(x^2+1)^{-1/2}$ leads to $F= D_e(D_e^*D_e+I)^{-1/2}$. This concludes the proof.
\end{proof}
For applications to relative uniform index maps it will be useful to have know explicit representatives of the uniform K-homology classes of a uniformly elliptic operator that have finite propagation, or are limits of operators with this property. For first-order operators this can be achieved by means of the propagation estimates discussed in Section \ref{SectionFinProp}.
\begin{prop}
Let $D$ be a symmetric first-order uniformly elliptic operator, and $D_e$ a local semi-interpolation-regular boundary condition for $D$. Then, $D_e(D_e^*D_e+I)^{-\frac{1}{2}}$ is an element of $D^*_u(\olO)$, and is an involution modulo $D^*_u(\olO,\partial\Omega)$.
\end{prop}
\begin{proof}
Let $\hat{D}_e$ and $P$ be as in Lemma \ref{KHomDilation}. Then $\hat{D}_e$ is a local and self-adjoint extension of a uniformly elliptic first-order differential operator. Thus, if $\chi$ is a normalizing function with compactly supported Fourier transform, the operator $\chi(\hat{D}_e)$ has finite propagation by Corollary \ref{finPropFourier}. In combination with Proposition \ref{KHomGradedComp} we conclude that $\chi(\hat{D}_e) \in D^*_u(\olO)$. The projection $P$ has finite propagation, in fact zero propagation, so that $(I-P)\chi(\hat{D}_e)P$ is a uniformly pseudo-local finite-propagation operator as well. The same holds for $P\chi(\hat{D}_e)(I-P)$. Since moreover $(I-P)\chi(\hat{D}_e)P-P\chi(\hat{D}_e)(I-P) \in \fD^{u}(\olO,\partial\Omega)$ it follows that $(I-P)\chi(\hat{D}_e)P$ is self-adjoint modulo $D^*_u(\olO,\partial\Omega)$. Now, normalizing functions with compactly supported Fourier transform are dense in the space of all normalizing functions. Thus $D_e(D_e^*D_e+I)^{-1/2}=(I-P)\chi_0(\hat{D}_e)P$ (with $\chi_0(x)=x(x^2+1)^{-1/2}$) is the limit of uniformly pseudolocal finite-propagation operators that are self-adjoint modulo $D^*_u(\olO,\partial\Omega)$. Hence it is also in $D^*_u(\olO)$, and self-adjoint modulo $D^*_u(\olO,\partial\Omega)$. Lastly, an analogous approximation argument shows that the operator
$$
\left(D_e^*D_e+I\right)^{-1} = I - \left( D_e(D_e^*D_e+I)^{-\frac{1}{2}} \right) ^* \left( D_e(D_e^*D_e+I)^{-\frac{1}{2}} \right)
$$
is in $D_u^*(\olO,\partial\Omega)$, proving that $D_e(D_e^*D_e+I)^{-1/2}$ is an involution modulo $D_u^*(\olO,\partial\Omega)$.
\end{proof}
We remark that in case one only cares about first-order differential operators, one can use propagation arguments to show that such operators give rise to classes in K-homology. In the boundary-less case that approach is taken in \cite{HigsonRoe2000}. Baum, Douglas and Taylor also use these techniques in \cite{BDT1989} to treat classes in the interior, and obtain independence of the choice of boundary condition. It seems likely that these arguments can be made uniform to prove the analogous statements for uniform K-homology. In fact Engel used wave operators to treat uniformly elliptic pseudo-differential operators on manifolds of bounded geometry \cite{Engel2018}. As indicated above these wave operators no longer have finite propagation (even if one restricts to differential operators), but they are quasi-local, meaning there is a uniform bound on the amount of mass they propagate over a certain distance.  

Let us conclude this section with a brief discussion of Dirac operators. We recall some basics here, but refer the reader to Appendix \ref{MultigradingAppendix} for a more detailed account. So far we have encountered the spinor Dirac operator $\Dirac$ over (domains in) a spin manifold. Suppose that $\Omega$ is a domain in an $n$-dimensional spin manifold of bounded geometry. If $\Omega$ is even-dimensional, then $\Dirac$ is graded and odd, so that it defines a class
$$
[\Dirac] \in K_0^{u}(\olO,\dO) \; .
$$
If $\dim(\Omega)=n$ is instead odd, $\Dirac$ is ungraded, hence defines a class
$$
[\Dirac] \in K_1^{u}(\olO,\dO) \; .
$$
There is also the Clifford-linear Dirac operator $\slashed{\mathfrak{D}}$ acting on sections of the Clifford-module bundle $\slashed{\mathfrak{S}}(\olO)$. It is $n$-multigraded and odd, hence defines a class
$$
[\slashed{\mathfrak{D}}] \in K_{-n}^{u}(\olO,\dO) \; .
$$
The operators $\Dirac$ and $\slashed{\mathfrak{D}}$ are connected via two-fold periodicity. Indeed, applying two-fold periodicity $\lceil\frac{n}{2}\rceil$ times to the $n$-multigraded Clifford-module bundle $\slashed{\mathfrak{S}}(\olO)$ produces a graded bundle if $n$ is even, and an ungraded one if $n$ is odd, and in either case this bundle equals the spinor bundle $\slashed{S}(\olO)$. Under this repeated application of two-fold periodicity $\slashed{\mathfrak{D}}$ and $\Dirac$ correspond to each other. Thus we conclude the following.
\begin{prop}[Compare {\cite[Section 11.3]{HigsonRoe2000}}]
Let $\Omega$ be a domain in an $n$-dimensional spin manifold of bounded geometry. Then,
$$
[\slashed{\mathfrak{D}}] \in K_{-n}^{u}(\olO,\dO) \quad \text{and} \quad [\Dirac] \in K_{n \, \mathrm{mod} \, 2}^{u}(\olO,\dO)
$$
are mapped to each other by the two-fold periodicity isomorphisms.
\end{prop}
Thus it is essentially a matter of convenience whether we work with $\Dirac$ or $\slashed{\mathfrak{D}}$ (in the complex case, see Remark \ref{RealKHom} below). For the most part we will work with the latter, as it behaves favorably with respect to products. To that end we first recall multigraded tensor products of vector bundles and differential operators. Again, a more detailed account can be found in Appendix \ref{MultigradingAppendix}. Suppose that $M$ and $M'$ are manifolds (with or without boundary), and $E\to M$ and $E'\to M$ are $p$- and $p'$-multigraded vector bundles. Suppose further that $D$ and $D'$ are $p$- and $p'$-multigraded differential operators acting on sections of $E$ and $E'$ respectively. Then, one may form the exterior multigraded tensor product
$$
S \hat{\boxtimes} S' \longrightarrow M \times M' \; ,
$$
which is a $(p+p')$-multigraded vector bundle. Using the multigraded tensor product one then forms the product operator
$$
D \times D' := D \hat{\otimes} I + I \hat{\otimes} D' \; ,
$$
which is a $(p+p')$-multigraded differential operator acting on sections of $S \hat{\boxtimes} S'$. Now, if $M$ and $M'$ are spin manifold, then $M\times M'$ carries a canonical spin structure, and it holds that
$$
\slashed{\mathfrak{S}}(M\times M') = \slashed{\mathfrak{S}}(M) \hat{\boxtimes} \slashed{\mathfrak{S}}(M') \quad , \quad \slashed{\mathfrak{D}}_{M\times M'} = \slashed{\mathfrak{D}}_{M} \times \slashed{\mathfrak{D}}_{M'} \; .
$$
In particular, if $\olO$ is a manifold with boundary and bounded geometry that has product structure near the boundary, then on a uniformly thick tubular neighborhood $U\cong [0,\delta)\times \dO$ it holds that
$$
\slashed{\mathfrak{D}}_\Omega = \slashed{\mathfrak{D}}_\RR \times \slashed{\mathfrak{D}}_{\dO} \; .
$$
We will prove below that on uniform K-homology the multigraded tensor product of operators is implemented by the Kasparov product. This will be a crucial component of our proof of boundary-of-Dirac-is-Dirac. Because the Clifford-linear Dirac operators behave conveniently in regards to products, we will use them for our considerations below. Indeed, we will use Clifford-linear Dirac operators as a model case to phrase the condition that an operator has product structure near the boundary or more generally near a hypersurface. This multigraded framework is more convenient to work with, but one can always use two-fold periodicity to reduce to graded and ungraded operators only.
\begin{remark}[Real uniform K-homology and periodicity]\label{RealKHom}
The difference between $\Dirac$ and $\slashed{\mathfrak{D}}$ is only not felt because we are working with complex instead of Real K-homology. There is a Real\footnote{We use the convention of referring to objects like Hilbert space, $C^*$-algebras, vector bundles, etc., as \emph{Real} if they are defined over the field of complex numbers, and additionally equipped with an anti-linear involution. An operator on a Real space is \emph{Real} if it is complex-linear and commutes with the anti-linear involution. In contrast a \emph{real} Hilbert space, $C^*$-algebra, vector bundle, etc., is simply one defined over the field of real numbers.} version of K-homology, whose cycles are given by Real multigraded Hilbert spaces together with a Real representation and a Real operator satisfying the usual conditions. Thus we can also define \emph{uniform Real K-homology} in the obvious way. The contents of this chapter up to this point then hold verbatim for uniform Real K-homology if we additionally require the uniformly elliptic operators to be Real. In particular  both the spinor Dirac operator and the Clifford-linear Dirac operator define classes in (relative) uniform Real K-homology. However instead of a two-fold periodicity Real K-homology satisfies an eight-fold periodicity. Thus, the uniform Real K-homology class of $\slashed{\mathfrak{D}}$ contains information that can generally not be recovered from $\Dirac$.

Let us remark that transferring the results of Chapter 2 and 3 to uniform Real K-homology is not quite as straight-forward as adding the word \glqq Real\grqq \, in the appropriate places. This is because Paschke duality, which laid the foundation of our discussion in Chapter 2 and hence Chapter 3 as well, uses two-fold periodicity in an essential way. One first would have to prove a Real version of Paschke duality, see for example \cite{Roe2004} for the corresponding discussion for analytic K-homology. If this was done, then the author anticipates the Real counterparts of our results to follow rather straight-forwardly.
\end{remark} 
\section{Properties of the uniform K-homology classes}
This section is devoted to establishing the compatibility of the uniform K-homology classes of uniformly elliptic operators over domains with some of the basic uniform K-homological operations, namely restriction to subsets, the Kasparov product, and the boundary map. These compatibility statements provide crucial computational tools that will allow us to deduce index theorems in the following section. 
\subsection{Restrictions}
Let $\Omega$ and $\Omega'$ be domains in the manifold $M$ of bounded geometry. As usual we endow them with the subspace metric coming from $M$. If $\Omega'\subseteq \Omega$, then the inclusion map $C_0(\Omega') \hookrightarrow C_0(\Omega)$ is filtered, thus it induces a restriction map $i^!: \, K^{u}_*(\Omega)\to K^{u}_*(\Omega')$. In case domain and target of the restriction map need to be made explicit we also write $i^!_{\Omega\to\Omega'}$ for this map.
\begin{prop}\label{Restriction}
Let $\Omega' \subseteq \Omega$ be domains in the bounded-geometry manifold $M$. Let $D$ be a symmetric uniformly elliptic differential operator over $\olO$. Then
$$
i^![D_\Omega] = [D_{\Omega'}] \; \in K_*^{u}(\Omega') \; .
$$
\end{prop}
\begin{proof}
It suffices to prove the claim for the restriction from $M$ to a domain $\Omega$. Indeed, assume the claim holds for the restriction from $M$ to any domain. Then,
$$
i^!_{\Omega\to \Omega'} [D_\Omega] = i^!_{\Omega\to \Omega'} i^!_{M\to \Omega}[D_M] = i^!_{M\to \Omega'}[D_M] = [D_{\Omega'}] \; .
$$
Thus, let $D$ be symmetric uniformly elliptic over $M$, and let $\Omega$ be a domain. We consider the minimal boundary condition $D_0:=D_{\min}$ for $D$ as an operator over $\Omega$. Then, $[D_\Omega]$ is represented by $(L^2(\Omega;E),\rho_\Omega,D_0(D_0^*D_0+I)^{-1/2})$, where $\rho_\Omega$ is the representation by multiplication operators. On the other hand $i^![D_M]$ is represented by $(L^2(M;E),\rho_M|_\Omega,D(D^*D+I)^{-1/2})$. This Fredholm module is a uniformly compact perturbation of the sum of $(L^2(\Omega;E),\rho_\Omega,qD(D^*D+I)^{-1/2}q)$ and a degenerate, where $q$ is the projection given by the indicator function of $\Omega$. We will show that $D_0(D_0^*D_0+I)^{-1/2}$ and $qD(D^*D+I)^{-1/2}q$ differ by a uniformly compact operator.

To that end let $D'$ denote the restriction of $D$ to $M\setminus \olO$, and let $D'_0$ denote its minimal boundary condition. Consider the closed operator $T=D_0 \oplus D'_0$ on $L^2(M;E)=L^2(\Omega;E)\oplus L^2(M\setminus\olO;E)$.\footnote{For this decomposition to be valid we have to assume that $\partial\Omega$ has measure zero. If this is not the case, fix $\delta>0$. By Lemma \ref{LLipStratified} it suffices to consider functions in $C^\infty_{cc}(\Omega)$ that are supported distance $\geq\delta$ from $\partial\Omega$. Then $\partial\Omega$ can be modified within a $\delta/2$-neighborhood to have measure zero, and we can proceed with the same proof.} Observe that $(D-T)f=0=(D^*D-T^*T)f$ for $f\in C^\infty_{cc}(\Omega)$. Using the arguments of Section \ref{SectionLowerOrderPerturbations}\footnote{While the results of that Section were formulated for the difference of boundary conditions for two uniformly elliptic operators, everything beyond Lemma \ref{CommutatorBasicLemma} is an abstract argument and hence also valid for an operator like $T$ that still satisfies the relevant analytical properties, but is strictly speaking not a boundary condition.} one may then show that 
$$
\left( D\left(D^*D+I\right)^{-\frac{1}{2}} - T\left(T^*T+I\right)^{-\frac{1}{2}} \right) f \sua 0
$$
for $f\in C_0(\Omega)$. Hence the same holds true when multiplying both operators with $q$ on both sides. Since $qT(T^*T+I)^{-1/2}q=D_0(D_0^*D_0+I)^{-1/2}$, this concludes the proof.
\end{proof}
\subsection{Products}
We come to the relation between the multigraded tensor product of operators and the Kasparov product. Recall that if $D$ and $D'$ are multigraded operators acting on sections of the multigraded bundles $E$ and $E'$, then $D\times D'$ denotes the operator $D\hat{\otimes}I+I\hat{\otimes}D'$ acting on sections of the multigraded exterior tensor product $E\hat{\boxtimes}E'\to M\times M'$. Recall also from Lemma \ref{DifferentMetrics} above that if $M$ and $M'$ are Riemannian manifolds, then the product metric $d_{M}+d_{M'}$ and the path length metric $d_{M\times M'}$ induced by the product Riemannian metric $g_{M\times M'}=g_M\oplus g_{M'}$ are equivalent. Hence either may be used for the uniform K-homology of the product $M\times M'$.
\begin{prop}\label{ProductComplete}
Let $M$, $M'$ be manifolds of bounded geometry, and $D$, $D'$ symmetric uniformly elliptic operators over $M$ and $M'$ respectively. Then,
$$
[D] \times [D'] = [D\times D'] \; \in K_*^{u}(M\times M') \; .
$$
\end{prop}
\begin{proof}
The proof is analogous to that in analytic K-homology. Indeed, the discussion in Section 10.7 in \cite{HigsonRoe2000} covers our case as well.\footnote{Technically, Higson and Roe treat the case where $D$ and $D'$ are first-order differential operators, but it is easily verified that their arguments hold verbatim for higher-order operators as well. Indeed, their discussion is mostly operator-theoretic in nature, apart from using at one point that elements in the domain of the unique closure of $D$, $D'$ can be approximated by smooth sections. This is true in our setting as well.} We conclude that $\chi(D\times D')$ is aligned to $\chi(D)$ and $\chi(D')$ in the sense of \cite[Definition 9.2]{HigsonRoe2000} for any normalizing function $\chi$. Apart from finite-propagation requirements that definition is the same as Definition \ref{AlignedModules}, and these requirements are taken care of by Remark \ref{PropagationKasparov}. Thus $\chi(D\times D')$ represents the Kasparov product of $\chi(D)$ and $\chi(D')$ on uniform K-homology as well.
\end{proof}
We deduce the corresponding product formula on domains by use of restriction maps. Note that in order to consider the Kasparov product on domains in the first place we need to assume them to have jointly bounded geometry. This is satisfied in particular for bounded-geometry domains. Also, we again endow the product of two domains with the product metric of the respective subspace metrics. This coincides with the subspace metric induced by the product metric on the product of the ambient manifolds. 
\begin{cor}\label{ProductDomain}
Let $\Omega \subseteq M$, $\Omega' \subseteq M'$ be domains in manifolds of bounded geometry. Assume that $\Omega$, $\Omega'$ have jointly bounded geometry. Let $D$, $D'$ be symmetric uniformly elliptic operators over $\olO$, $\olO'$ respectively. Then, 
$$
[D] \times [D'] = [D\times D'] \; \in K_*^{u}(\Omega\times \Omega')
$$
\end{cor}
\begin{proof}
By assumption there are symmetric uniformly elliptic pseudo-differential operators $\tilde{D}$, $\tilde{D}'$ over $M$, $M'$ that restrict to $D$ and $D'$ respectively. The Kasparov product is natural with respect to the restriction maps by Theorem \ref{KasparovProperties}. Combining this with Propositions \ref{Restriction} and \ref{ProductComplete} we compute that
$$
[D] \times [D'] = (i^! \times (i')^!)\left( [\tilde{D}] \times [\tilde{D}'] \right) = (i^! \times (i')^!) [\tilde{D}\times \tilde{D}'] = [D \times D'] \; .
$$
\end{proof}
In the special of Dirac operators we deduce the following from their combatibility with product.
\begin{cor}
Let $\Omega$ and $\Omega'$ be domains of jointly bounded geometry in spin manifolds of bounded geometry and dimension $n$ and $n'$ respectively. Then,
$$
[\slashed{\mathfrak{D}}_{\Omega\times\Omega'}] = [\slashed{\mathfrak{D}}_{\Omega}] \times [\slashed{\mathfrak{D}}_{\Omega'}] \; \in \; K_{-(n+n')}^{u}(\Omega\times\Omega') \; .
$$
Thus, under two-fold periodicity it also holds that
$$
[\Dirac_{\Omega\times\Omega'}] = [\Dirac_{\Omega}] \times [\Dirac_{\Omega'}] \; \in \; K_{(n+n')\, \mathrm{mod}\, 2}^{u}(\Omega\times\Omega')
$$
\end{cor}
\subsection{Boundary maps}
The goal of this section is to determine the image of the relative class defined by a uniformly elliptic operator under the boundary map of the long exact sequence in uniform K-homology. We will restrict ourselves to first-order differential operators, though see Remark \ref{ReductionToFirstOrder} below. 

In our discussion we will require the operator to have product structure near the boundary. This is defined as follows.
\begin{definition}
Let $\olO$ be a manifold with boundary and bounded geometry, $E\to \olO$ a vector bundle of bounded geometry, and $D: \, \Gamma_{cc}(\Omega;E)\to\Gamma_{cc}(\Omega;E)$ a differential operator. Let $p\in\ZZ_{\geq -1}$, and suppose that $E$ and $D$ are $(p+1)$-multigraded, and that $D$ is odd. We say that $D$ has \emph{product structure} near $\partial\Omega$ if there exists a uniformly thick collar neighborhood of $\partial \Omega$, and a multigraded unitary isomorphism between $E$ and $\slashed{\mathfrak{S}}(\RR)\hat{\boxtimes} E_{\partial\Omega}$ for some $p$-multigraded vector bundle $E_{\partial\Omega}\to \partial\Omega$ over that collar neighborhood, so that under this isomorphism $D$ is given by 
$$
D=\slashed{\mathfrak{D}}_\RR \times D_{\partial\Omega} 
$$
for some odd $p$-multigraded differential operator $D_{\partial \Omega}: \, \Gamma_c(\partial\Omega;E_{\partial\Omega})\to \Gamma_c(\partial\Omega;E_{\partial\Omega})$.
\end{definition}
We can now formulate the boundary-of-Dirac-is-Dirac formula in uniform K-homology. Its proof is rather straight-forward application of Theorem \ref{Suspension}, together with the compatibility with restrictions and products we have derived in the previous subsections. We consider $\dO$ to carry the subspace metric coming from the inclusion $\dO\subseteq M$, though see Lemma \ref{DifferentMetrics}. 
\begin{thm}\label{BdryOfDirac}
Let $\olO$ be a manifold with boundary and bounded geometry. Let $D$ be a uniformly elliptic $(p+1)$-multigraded first order differential operator over $\olO$, and assume that $D$ has product structure $D=\slashed{\mathfrak{D}}_\RR \times D_{\dO}$ near $\dO$. Then
$$
\partial [D] = (-1)^{p+1} [D_{\dO}] \; \in \; K^{u}_{-p}(\dO) \; .
$$
\end{thm}
\begin{proof}
By naturality of the boundary map we may restrict to a uniformly thick collar neighborhood. Indeed, $\dO$ admits a uniformly thick collar neighborhood $U$ that bi-Lipschitz equivalent to $[0,\delta)\times \dO$ with the product metric (see \cite[Proposition 7.3]{Schick1998}). Then, the boundary map factors as \\
\centerline{\xymatrix{
K^{u}_{-(p+1)}(\olO,\dO) \ar[r]^\partial \ar[d]_{i^!} & K^{u}_{-p}(\dO) \ar@{=}[d] \\
K^{u}_{-(p+1)}(U,\dO)\ar[r]^\partial \ar@{=}[d] & K^{u}_{-p }(\dO) \ar@{=}[d] \\
K^{u}_{-(p+1)}([0,\delta)\times \dO,\dO) \ar[r]^-\partial & K^{u}_{-p}(\dO) 
}}
By Proposition \ref{Restriction} the restriction map $i^!$ maps $[D]$ to $[D_U] \in K^{u}_{-(p+1)}(U,\dO)$. Since we assume $D$ to have product structure near the boundary $[D_U]$ is identified with 
$$
[\slashed{\mathfrak{D}}_\RR\times D_{\dO}] = [\slashed{\mathfrak{D}}_\RR] \times [D_{\dO}] \; \in \; K^{u}_{-(p+1)}([0,\delta)\times \dO,\dO) \; ,
$$
where we used Corollary \ref{ProductDomain}. Now, $[\slashed{\mathfrak{D}}_\RR]$ coincides with the generator $d\in K^{u}_{-1}([0,\delta),\{0\})$. Thus 
$$
\partial \left( [\slashed{\mathfrak{D}}_\RR] \times [D_{\dO}] \right) = (-1)^{-(p+1)}[D_{\dO}] \; \in \; K^{u}_{-p}(\dO)
$$
by Theorem \ref{Suspension}. It follows that $\partial[D]=(-1)^{p+1}[D_{\dO}]$ as well, proving the claim.
\end{proof}
\begin{remark}[Regarding the sign]\label{SignBdryDirac}
The sign in the formula $\partial [D] = (-1)^{p+1} [D_{\dO}]$ is an artifact of our chosen construction of the boundary map. As opposed to using Paschke duality and the K-theoretical boundary map (which is what we did) it is also possible to use the Kasparov product and a mapping cone construction to construct the long exact sequence in K-homology, see \cite[Section 9.5]{HigsonRoe2000}. The resulting boundary maps differ by a sign of $(-1)^{p+1}$. For the latter the formula $\partial [D] = [D_{\dO}]$ thus holds without the sign, which is not surprising given that its construction is organized around the Kasparov product. We chose to stay consistent in the use of the boundary map constructed in Chapter 2, hence the sign.
\end{remark}
We obtain the following corollary for Dirac operators.
\begin{cor}
Let $\olO$ be an $(n+1)$-dimensional spin manifold with boundary and bounded geometry that has product structure near the boundary. Then,
$$
\partial [\slashed{\mathfrak{D}}_\Omega] = (-1)^{n+1} [\slashed{\mathfrak{D}}_{\partial\Omega}] \; \in \; K_{-n}^{u}(\dO) \; .
$$
Thus, under two-fold periodicity it also holds that
$$
\partial [\Dirac_\Omega] = (-1)^{n+1} [\Dirac_{\partial\Omega}] \; \in \; K_{n\, \mathrm{mod}\, 2}^{u}(\dO) \; .
$$
\end{cor}
By exactness of the long exact sequence we know that $[D]\in K^{u}_*(\olO,\dO)$ can only have non-zero image under the boundary map if $[D]$ does not lift to a class in $K_*^{u}(\olO)$. If $D$ has product structure, the class $[D_{\partial\Omega}]=\partial[D]$ therefore obstructs the existence of boundary conditions for $D$ that would allow for the construction of such a lift.
\begin{cor}\label{ObstructionRegular}
Let $D$ be a uniformly elliptic first-order differential operator with product structure $D=\slashed{\mathfrak{D}}_\RR\times D_{\partial\Omega}$ near the boundary. If $D$ admits a boundary condition that defines a uniform Fredholm module over $\olO$, then 
$$
[D_{\partial\Omega}] = 0 \; \in \; K_*(\dO) \; .
$$
In particular, if $[D_{\partial\Omega}]\neq 0$, then $D$ admits no local, regular and self-adjoint boundary condition.
\end{cor}
\begin{example}[Chirality boundary conditions revisited]\label{ExChirality2}
In Example \ref{ExChirality1} above we noted that even though they are local and regular the chirality boundary conditions for the spin Dirac operator on an odd-dimensional manifold are not covered by Theorem \ref{ExistencKHomClasses} due to lack of self-adjoitness. In principle they could of course define a uniform Fredholm module nonetheless. It is a consequence of Corollary \ref{ObstructionRegular} that they cannot. This implies that some condition relating a local regular boundary condition to its adjoint is required to get a uniform Fredholm module over $\olO$. Recalling how we proved that $D_e(D^*_eD_e+I)^{-1/2}$ is self-adjoint modulo $\mathfrak{C}^{u}(\Omega)$ for any boundary condition $D_e$ this makes intuitive sense. Indeed, this uniformly essential self-adjointness was a consequence of the fact that different boundary conditions (in particular $D_e$ and $D_e^*$) agree on the interior. This no longer holds over the entire manifold, as is witnessed by the fact that different boundary conditions can define different K-homology classes. In view of Section \ref{SectionLowerOrderPerturbations} it seems reasonable to conjecture that a regular differential boundary condition defines a uniform Fredholm module over $\olO$ if the principal symbol of the differential operator defining the boundary condition is self-adjoint.
\end{example}
\begin{remark}[Higher-order operators]\label{ReductionToFirstOrder}
We discussed boundary maps only first-order elliptic operators with product structure near the boundary, meaning essentially that $D$ is the sum of a tangential and a normal differential operator. If $D$ is a higher-order differential operator, then splitting $D$ into a tangential and a normal is not immediately possible. However, at least on Euclidean space an $m$th-order differential equation can be transformed into a system of $m$ first-order equations (see for example Section VII.12 in \cite{Taylor2011}), and it makes sense to require this system of first-order equations to have product form. Generally, these first-order equations will be pseudo-differential instead of differential. If $D$ is a suitable $m$th-order operator on $\Omega$, we can locally transform $D$ into a system of $m$ first-order equations in product form encoded by some operator $\tilde{D}=\slashed{\mathfrak{D}}_\RR\times P$, where $P$ is a first-order pseudo-differential operator, and we ignore details about (multi)gradings. It seems plausible that under the usual uniformness conditions, these local transformations can be patched together to globally represent $D$ by a system $\tilde{D}=\slashed{\mathfrak{D}}_\RR\times P$ of $m$ first-order pseudo-differential equations in product form. Then, one could perform the same analysis as in the first-order case. Thus, it seems plausible that for a higher-order differential operator $D$ the boundary $\partial[D]$ is given by the class $[P]$. We point to \cite{Fries2025} for a recent treatment of the boundary of the relative class of higher-order operators in the compact case, which uses the Calderón projection instead.
\end{remark}
\begin{remark}[Less regular boundaries]\label{LipschitzBdryMap}
Even if $\partial\Omega$ is not smooth, there is still a boundary map $K^{u}_{*+1}(\olO,\partial\Omega)\to K^{u}_*(\partial\Omega)$, and we can ask where uniformly elliptic operator over $\olO$ get mapped to. Let us suppose that $\olO$ is at least a uniform Lipschitz domain, so that $\partial\Omega$ is smooth almost everywhere. Then, derivatives and normal vector fields can still be made sense of a.e. on $\partial\Omega$, and Sobolev spaces $H^s(\partial\Sigma)$ can be defined at least for $|s|\leq 1$. Then we may still speak of a first-order operator having product form near the boundary, where tangential operator and normal derivative are now a.e. defined. See for example \cite{McLean2000} for a discussion in the case of domains in Euclidean space. The author expects that upon suitable interpretation Theorem \ref{BdryOfDirac} and the corollary results still hold in the Lipschitz setting. 
\end{remark}
\section{Bordism invariance and partitioned manifolds}
In the last section of this chapter we combine the results on uniform K-homology classes of uniformly elliptic operators over manifolds with boundary with our previous results on (relative) uniform index maps. First we deduce the bordism invariance of the uniform coarse index, then we prove a non-relative and a relative uniform partitioned-manifold index theorem.
\subsection{Bordism invariance of the uniform coarse index}\label{SectionBordismInv}
We begin with bordism invariance. As in the compact case it is straight-forwardly deduced from boundary-of-Dirac-is-Dirac. We will say that a manifold $\olO$ with boundary and bounded geometry has \emph{finite width} if the function $d(x,\dO)$ remains uniformly bounded for $x\in \olO$. If $\olO$ has finite width, the inclusion $\dO\hookrightarrow \olO$ is a uniform coarse equivalence, which is necessary for the proof of bordism invariance below. Furthermore it is also necessary to have an interesting bordism theory, because otherwise any manifold $\Sigma$ would be null-bordant via the infinite-width null-bordism $\olO=\RR_{\geq0}\times\Sigma$. Note also that the following proposition generally fails in this case.
\begin{prop}\label{BordismInvariance}
Let $\olO$ be a finite-width manifold with boundary and bounded geometry. Let $D$ be a $(p+1)$-multigraded uniformly elliptic operator over $\olO$ that has product structure $D=\slashed{\mathfrak{D}}_\RR \times D_{\dO}$ near the boundary. Then,
$$
\Ind(D_{\dO}) = 0 \; \in \; K_{-p}\left( C^*_u(\dO) \right) \; .
$$
\end{prop}
\begin{proof}
Due to the assumption of finite width the inclusion $\iota: \, \dO\hookrightarrow \olO$ is a uniform coarse equivalence, hence induces an isomorphism $C^*_u(\dO) \cong C^*_u(\olO)$. By naturality of the index map (or Theorem \ref{CommUniIndex}) we have a diagram \\
\centerline{\xymatrix{
K^{u}_{-(p+1)}(\olO,\dO) \ar[r]^-\partial & K^{u}_{-p}(\dO) \ar[r]^{\iota_*} \ar[d]_{\Ind} & K^{u}_{-p}(\olO) \ar[d]^{\Ind} \\
& K_{-p}\left( C^*_u(\dO) \right) \ar[r]^\sim_{\iota_*} & K_{-p}\left(C^*_u(\olO)\right)
}}
Using Theorem \ref{BdryOfDirac} and exactness of the top row we compute that
$$
\iota_*\left( \Ind(D_{\dO})\right) = \Ind\left(\iota_* [D_{\dO}]\right) = (-1)^{p+1} \Ind\left( (\iota_* \circ \partial)[D] \right) = 0 \; .
$$
Since $\iota_*: \, K_{-p}\left( C^*_u(\dO) \right) \to K_{-p}\left(C^*_u(\olO)\right)$ is an isomorphism, it follows that $\Ind(D_{\dO})=0$ as well.
\end{proof}
Suppose that $\olO$ is a manifold with boundary and bounded geometry so that $\dO$ is the disjoint union $\Sigma \sqcup \Sigma'$. If the functions $x\mapsto d(x,\Sigma)$ and $x\mapsto d(x,\Sigma')$ are both uniformly bounded\footnote{Note that this requirement is stronger than requiring that $\olO$ have finite width.} for $x\in\olO$, then $\Sigma\hookrightarrow \olO$ and $\Sigma'\hookrightarrow\olO$ are both uniform coarse equivalences. Thus there are canonical identifications $K_{*}(C^*_u(\Sigma))\cong K_{*}(C^*_u(\olO))\cong K_{*}(C^*_u(\Sigma'))$, by means of which we may compare indices over $\Sigma$ and $\Sigma'$. The previous proposition then has the the following corollary for Dirac operators.
\begin{cor}
Let $\Sigma$ and $\Sigma'$ be $n$-dimensional spin manifolds of bounded geometry. Suppose there exists a spin manifold $\olO$ with boundary and bounded geometry that has product structure near the boundary, so that $\partial \Omega = \Sigma \sqcup (-\Sigma')$, where $-\Sigma'$ indicates $\Sigma'$ with reversed orientation. Suppose that the functions $x\mapsto d(x,\Sigma)$ and $x\mapsto d(x,\Sigma')$ are both uniformly bounded for $x\in\olO$. Then
$$
\Ind(\slashed{\mathfrak{D}}_\Sigma) = \Ind(\slashed{\mathfrak{D}}_{\Sigma'})
$$
upon identifying $K_{-n}(C^*_u(\Sigma))\cong K_{-n}(C^*_u(\olO))\cong K_{-n}(C^*_u(\Sigma'))$. Thus, under two-fold periodicity it also holds that
$$
\Ind(\Dirac_\Sigma) = \Ind(\Dirac_{\Sigma'}) \; .
$$
\end{cor}
\begin{proof}
The Clifford-linear Dirac operator on $\olO$ has product structure $\slashed{\mathfrak{D}}_\Omega= \slashed{\mathfrak{D}}_\RR \times \slashed{\mathfrak{D}}_{\partial\Omega}$ near the boundary. Hence $\Ind(\slashed{\mathfrak{D}}_{\partial \Omega})=0$ in $K_{-n}(C^*_u(\partial\Omega))$. The uniform coarse equivalences $\Sigma,\Sigma' \hookrightarrow \partial\Omega$ induce canonical isomorphisms on the corresponding uniform Roe algebras. Under this identification it holds that $\Ind(\slashed{\mathfrak{D}}_\Sigma)-\Ind(\slashed{\mathfrak{D}}_{\Sigma'})=\Ind(\slashed{\mathfrak{D}}_{\partial \Omega})=0$. The claim follows.
\end{proof}
The uniform coarse index $\Ind(\Dirac_\Sigma)$ on a spin manifold $\Sigma$ of bounded geometry is an obstruction to the existence of a upsc-metric of bounded geometry in the uniform coarse equivalence class of $\Sigma$ \cite[Proposition 7.10]{Engel2014}. Via bordism invariance of the uniform coarse index we conclude that $\Sigma$ can already not admit such a metric if it is bordant to a spin manifold with non-vanishing uniform coarse index.
\begin{cor}
Let $\Sigma$ be a spin manifold of bounded geometry. Suppose $\Sigma$ is spin-bordant via a spin bordism of bounded geometry as in the previous corollary to a spin manifold $\Sigma'$ of bounded geometry such that
$$
\Ind(\Dirac_{\Sigma'})\neq 0 \; \in \; K_*\left(C^*_u(\Sigma')\right) \; .
$$
Then $\Sigma$ does not admit a upsc-metric of bounded geometry in its uniform coarse equivalence class.
\end{cor}
Before we turn to partitioned-manifold indices let us note one more consequence of the boundary-of-Dirac-is-Dirac formula pertaining to relative uniform indices. It follows directly from the combination of Theorems \ref{IndexCommAlgMV} and \ref{BdryOfDirac}.
\begin{prop}
Let $\olO$ be a manifold with boundary and bounded geometry. Let $D$ be a $(p+1)$-multigraded uniformly elliptic operator that has product structure $D =\slashed{\mathfrak{D}}_\RR\times D_{\partial\Omega}$ near the boundary. Then, the algebraic Mayer-Vietoris boundary map
$$
\delta: \; K_{-(p+1)}\left( D^*(\olO,\partial\Omega)\right) \longrightarrow K_{-p}\left( C^*_u(\partial\Omega) \right)
$$
satisfies
$$
\delta\big( \relInd(D_\Omega) \big)=(-1)^{p+1}\Ind(D_{\partial\Omega}) \; \in \; K_{-p}\left( C^*_u(\partial\Omega) \right) \; .
$$
\end{prop}
We note that this also implies bordism invariance of the uniform coarse index. Indeed, if $\olO$ has finite width, then the algebraic Mayer-Vietoris map is zero, as we have noted in Remark \ref{FinWidthMV}. Hence $\Ind(D_{\dO})=(-1)^{p+1}\delta(\relInd(D))=0$ as well.
\subsection{Uniform partitioned-manifold index theorems}
To conclude this chapter we prove partitioned-manifold index theorems that relate the uniform coarse index on a manifold to that on a certain hypersurface. We begin by introducing the necessary terminology.
\begin{definition}
Let $M$ be a manifold of bounded geometry. A \emph{bounded-geometry partitioning hypersurface} in $M$ is a bounded-geometry hypersurface $N\subseteq M$ in the sense of Definition \ref{BGHypersurface} such that $M\setminus N$ has two-connected components. Note that each of these connected components is a bounded-geometry domain in $M$. 

Let $D_M$ be a differential operator over $M$ acting on the sections of a bounded-geometry vector bundle $E$. Let $p\in\ZZ_{\geq -1}$, and suppose that $E$ and $D$ are $(p+1)$-multigraded, and that $D$ is odd. We say that $D$ has \emph{product structure} near the bounded-geometry partitioning hypersurface $N\subseteq M$ if there exists a uniformly thick tubular neighborhood of $N$, and a multigraded unitary isomorphism between $E$ and $\slashed{\mathfrak{S}}(\RR)\hat{\boxtimes} E_{N}$ for some $p$-multigraded vector bundle $E_{N}\to N$ over that tubular neighborhood, so that under this isomorphism $D_M$ is given by 
$$
D_M=\slashed{\mathfrak{D}}_\RR \times D_{N} 
$$
for some odd $p$-multigraded differential operator $D_{N}$ acing on sections of $E_N$.
\end{definition}
The canonical example is again given by the Clifford-linear Dirac operator.
\begin{example}
If $M$ is a spin manifold of bounded geometry and $N\subseteq M$ a bounded-geometry partitioning hypersurface such that $M$ has product structure near $N$, then the Clifford-linear Dirac operator $\slashed{\mathfrak{D}}_M$ has product structure $\slashed{\mathfrak{D}}_M=\slashed{\mathfrak{D}}_\RR\times \slashed{\mathfrak{D}}_N$ near $N$
\end{example}
Given such a partitioning hypersurface $N\subseteq M$ and a uniformly elliptic operator $D_M$ over $M$ that has product structure $D_M=\slashed{\mathfrak{D}}_\RR\times D_N$ near $N$, we wish to relate the uniform coarse index of $D_M$ to that of $D_N$. Associated to the partition of $M$ into the two sides of $N$ there are Mayer-Vietoris boundary maps both in uniform K-homology and in K-theory of the uniform Roe algebras, and we know from Corollary \ref{MVDiagRoe} that the uniform coarse index map intertwines the two, i.e. there is a diagram\\
\centerline{\xymatrix{
K^{u}_{*+1}(M) \ar[r]^\delta \ar[d]_{\Ind} & K^{u}_*(N) \ar[d]^{\Ind} \\
K_{*+1}\left(C^*_u(M)\right) \ar[r]_\delta & K_{*}\left(C^*_u(N)\right)
}}  
Define the \emph{uniform partitioned manifold index map} as the composition
$$
\Ind^\mathrm{pm} := \delta\circ \Ind : \; K_{*+1}^{u}(M) \longrightarrow K_*\left( C^*_u(N) \right) \; .
$$
To compute it we may equivalently compute the index of the element $\delta[D_M]\in K^{u}_*(N)$. It is a consequence of the boundary-of-Dirac-is-Dirac formula that this element is given by $[D_N]$ (up to a sign).
\begin{prop}
Let $M$ be a manifold of bounded geometry, and $N$ a bounded-geometry partitioning hypersurface. Let $D$ be a $(p+1)$-multigraded differential operator, and assume that $D$ has product structure $D=\slashed{\mathfrak{D}}_\RR\times D_N$ near $N$. Then the Mayer-Vietoris boundary map $\delta: \, K^{u}_{-(p+1)}(M)\to K^{u}_{-p}(N)$ associated to the decomposition of $M$ into the two sides of $N$ satisfies
$$
\delta[D] = (-1)^{p+1} [D_N] \; \in \; K^{u}_{-p}(N) \; .
$$ 
\end{prop}
\begin{proof}
The Mayer-Vietoris boundary map on uniform K-homology equals the composition
\begin{equation} \label{MVComp}
K^{u}_{-(p+1)}(M) \xlongrightarrow{i^!_*} K^{u}_{-(p+1)}(M_+)\cong K^{u}_{-(p+1)}(\olM_+,N) \xlongrightarrow{\partial} K^{u}_{-p}(N)) \; ,
\end{equation}
where $M_+$ is the connected component of $M\setminus N$ with the property that if $\nu$ is a normal vector field on $N$ pointing into $M_+$, and $(e_1,\cdots,e_n)$ is an oriented frame on $N$, then $(\nu,e_1,\cdots,e_n)$ is an oriented frame on $M$. Then claim then follows from Proposition \ref{Restriction} and Theorem \ref{BdryOfDirac}.
\end{proof}
This computation together with the commutativity of the above diagram lets us conclude that the uniform partitioned-manifold index of $D_M$ agrees with the uniform coarse index of $D_N$ up to a sign.
\begin{thm}[Uniform partitioned-manifold index theorem]
Let $M$ be a manifold of bounded geometry, and $N\subseteq M$ a bounded-geometry partitioning hypersurface. Let $D_M$ be a $(p+1)$-multigraded elliptic operator over $M$, and assume that $D_M$ has product structure $D_M=\slashed{\mathfrak{D}}_\RR \times D_N$ near $N$. Then,
$$
\Ind^\mathrm{pm}(D_M) = (-1)^{p+1} \Ind(D_N) \quad \in \; K_{-p}\left(C^*_u(N)\right) \; .
$$
\end{thm}
\begin{remark}
Let us give a more explicit description of the partitioned manifold index in the case of the spinor Dirac operator on an odd-dimensional spin manifold $M$. Then $\Ind(\Dirac_M) \in K_{1}(C^*_u(M))$ can be represented by the Cailey transform $U=(\Dirac_M - i)(\Dirac_M+i)^{-1}$.\footnote{Recall that the uniform coarse index map is the boundary map $K^{u}_{1}(M)\cong K_0(D^*_u(M)/C^*_u(M)) \to K_1(C^*_u(M))$. The class of $\Dirac_M$ in $K_0(D^*_u(M)/C^*_u(M))$ is represented by the projection $(I+\chi(\Dirac_M))/2$ in $D^*_u(M)/C^*_u(M)$. Its boundary is given by the class of the unitary $\exp(\pi i (I+\chi(\Dirac_M))$. For the normalizing function $\chi(\lambda)=\frac{2}{\pi} \arctan(\lambda)$ this is precisely equal to the Cailey transform.} Let $\psi_\pm$ be the characteristic function of $M_\pm$. Then $\psi_- + \psi_+ U$ defines an element in $K_1(C^*_u(M_+\subseteq M)/C^*_u(N \subseteq M))$, and its image under the boundary map to $K_0(C^*_u(N\subseteq M))\cong K_0(C^*_u(N))$ represents that of $[U]$ under the Mayer-Vietoris boundary map. Thus we arrive at the formula
$$
\Ind^\mathrm{pm}(\Dirac_M) = \partial \left[ \psi_- + \psi_+ U \right] \; \in \; K_0\left( C^*_u(N) \right) \; .
$$
This formula is structurally identical with definitions of the (non-uniform) partitioned-manifold index that do not appeal explicitly to the Mayer-Vietoris boundary map, see for example \cite{Higson1991} in the case of compact hypersurfaces, and the recent pre-print \cite{Hochs2025} for noncompact ones.
\end{remark}
Restricting to Dirac operators we obtain the following.
\begin{cor}
Let $M$ be a spin manifold of bounded geometry, and $N\subseteq M$ a bounded-geometry partitioning hypersurface such that $M$ has product structure near $N$. Then,
$$
\Ind^\mathrm{pm}(\slashed{\mathfrak{D}}_M) = (-1)^{n+1} \Ind(\slashed{\mathfrak{D}}_{N}) \quad \in \; K_{-n}\left(C^*_u(N)\right) \; .
$$
Thus, under two-fold periodicity it also holds that
$$
\Ind^\mathrm{pm}(\Dirac_M) = (-1)^{n+1} \Ind(\Dirac_{N}) \quad \in \; K_{n \, \mathrm{mod}\, 2}\left(C^*_u(N)\right) \; .
$$
\end{cor}
Recall that if $M$ is spin, then the uniform coarse index of its (spinor) Dirac operator obstructs the existence of bounded-geometry metric with uniformly positive scalar curvature in its uniform coarse equivalence class. In particular, if $M$ admits such a metric, then the uniform partitioned manifold index $\Ind^{pm}(\Dirac_M)$ vanishes. The uniform partitioned-manifold index theorem therefore provides us an obstruction to the existence of such a metric that is based solely on the hypersurface $N$.
\begin{cor}
Let $M$ be an $(n+1)$-dimensional spin manifold of bounded geometry, and let $N \subseteq M$ be a bounded-geometry partitioning hypersurface such that $M$ has product structure near $N$. If
$$
\Ind(\Dirac_N)\neq 0 \; \in \; K_{n \, \mathrm{mod}\, 2}\left(C^*_u(N)\right) \; ,
$$
then $M$ admits no bounded-geometry metric with uniformly positive scalar curvature in its uniform coarse equivalence class.
\end{cor}
We can also carry out this discussion in the presence of a boundary. First, we need the definition of a partitioning hypersurface in this case. We find it easier to give the definition from the extrinsic viewpoint. Though it might be somewhat artificial, it is more convenient on a technical level.
\begin{definition}
Let $\Omega$ be a domain in the manifold $M$ of bounded geometry. Suppose that $N'\subseteq M$ is a bounded-geometry partitioning hypersurface in $M$. Then, $N:= N'\cap \olO$ will be called a \emph{bounded-geometry partitioning hypersurface} in $\olO$.

Similarly, if $D$ is a uniformly differential operator over $\olO$ that has an extension to a uniformly elliptic differential operator over $M$ that has product structure near $N'$, we will say that $D$ has \emph{product structure} near $N$. 
\end{definition}
Given such a bounded-geometry hypersurface $N\subseteq \olO$ the decomposition of $\olO$ into the two sides of $N$ gives rise to a Mayer-Vietoris sequence in relative uniform K-homology and one in K-theory of the relative uniform structure algebras. By Theorem \ref{AbstractRelPM} the relative uniform index maps intertwine these two sequences. In particular there is a diagram \\
\centerline{\xymatrix{
K^{u}_{*+1}(\olO,\partial\Omega) \ar[r]^\delta \ar[d]_{\relInd} & K^{u}_*(N,\partial N) \ar[d]^{\relInd} \\
K_{*+1}\left(D^*_u(\olO,\partial\Omega)\right) \ar[r]_\delta & K_{*}\left(D^*_u(N,\partial N)\right)
}}
Analogously to the boundaryless case we define the \emph{relative uniform partitioned-manifold index map} as the composition
$$
\Ind^\mathrm{rel,pm}:= \delta\circ \relInd: \; K^{u}_{*+1}(\olO,\partial\Omega) \longrightarrow K_{*}\left(D^*_u(N,\partial N)\right) \; .
$$
For its computation it is again sufficient to consider the Mayer-Vietoris boundary map on relative uniform K-homology. The latter map is computed by another boundary-of-Dirac-is-Dirac formula.
\begin{prop}
Let $\olO$ be a domain in a bounded-geometry manifold, and let $N\subseteq \olO$ be a bounded-geometry partitioning hypersurface. Let $D_\Omega$ be a $(p+1)$-multigraded uniformly elliptic operator over $\olO$ that has product structure $D_\Omega=\slashed{\mathfrak{D}}_\RR \times D_N$ near $N$. Then, 
$$
\delta [D_\Omega] = (-1)^{p+1} [D_N] \; \in \; K^{u}_{-p}(N,\partial N) \; .
$$
\end{prop}
\begin{proof}
The proof is a slightly more involved version of its boundaryless counterpart. Again, let $\Omega_+\subset \Omega$ denote the that connected component of $\Omega\setminus N$ with the property that if $\nu$ is a normal vector field on $N$ pointing into $\Omega_+$, and $(e_1,\cdots,e_n)$ is an oriented frame on $N$, then $(\nu,e_1,\cdots,e_n)$ is an oriented frame on $\Omega$. Let $M_+ \subseteq M\setminus N'$ denote the connected component containing $\Omega_+$. Consider the following diagram \\
\centerline{\xymatrix{
K^{u}_{-(p+1)}(\olO,\dO) \ar[d]_{\sim} \ar[rrr]^\delta & & & K^{u}_{-p}(N,\partial N) \ar[d]^\sim \\
K^{u}_{-(p+1)}(\Omega) \ar[r]^-{i^!} & K^{u}_{-(p+1)}(\Omega_+) & \ar[l]_-\sim K^{u}_{-(p+1)}(\olO_+\setminus\dO,N\setminus\partial N) \ar[r]^-\partial & K^{u}_{-p}(N\setminus\partial N)  \\
& K^{u}_{-(p+1)}(M_+) \ar[u]_{i^!} & \ar[l]^-\sim K^{u}_{-(p+1)}(\olM_+,N') \ar[r]_-\partial & K^{u}_{-p}(N') \ar[u]_{i^!}
}}
The large top square commutes by definition of the Mayer-Vietoris boundary map. The bottom square commutes due to naturality of the excision and boundary maps. It suffices to calculate that the image of $[D_\Omega]\in K^{u}_{-(p+1)}(\Omega)$ under the middle row coincides with $(-1)^{p+1}[D_N]\in K^{u}_{-p}(N,\partial N)$. 

To that end extend $D_\Omega$ to a uniformly elliptic operator $D_M$ over $M$ that has product structure $D_M=\slashed{\mathfrak{D}}_\RR\times D_{N'}$ near $N'$. Note that $D_{N'}=D_N$ on $N$. Then, $i^![D_\Omega]=i^![D_M] \in K^{u}_{-(p+1)}(\Omega_+)$ by Proposition \ref{Restriction}. Moreover, the combination of Theorem \ref{BdryOfDirac} and Proposition \ref{Restriction} computes the image of $[D_M]$ under the composition
$$
K^{u}_{-(p+1)}(M_+) \xlongrightarrow{\sim} K^{u}_{-(p+1)}(\olM_+,N') \xlongrightarrow{\partial} K^{u}_{-p}(N') \xlongrightarrow{i^!} K^{u}_{-p}(N\setminus \partial N)
$$
to be the class $(-1)^{p+1}[D_N] \in K^{u}_{-p}(N\setminus \partial N)$. Commutativity of the diagram lets us conclude that the image of $[D_\Omega]\in K^{u}_{-(p+1)}(\Omega)$ under the middle row is given by $(-1)^{p+1}[D_N]\in K^{u}_{-p}(N,\partial N)$ as well, concluding the proof.
\end{proof}
The appeal to Theorem \ref{AbstractRelPM} then lets us deduce the relative uniform partitioned-manifold index theorem.
\begin{thm}[Relative uniform partitioned-manifold index theorem]
Let $\olO$ be a domain in a manifold of bounded geometry, and let $N\subseteq \olO$ be a bounded-geometry partitioning hypersurface. Let $D_\Omega$ be a $(p+1)$-multigraded uniformly elliptic operator over $\olO$ that has product structure $D_\Omega=\slashed{\mathfrak{D}}_\RR \times D_N$ near $N$. Then, 
$$
\Ind^\mathrm{rel,pm}(D_\Omega) = (-1)^{p+1} \relInd(D_N) \; \in \; K_{-p}\left(D^*_u(N,\partial N)\right) \; .
$$
\end{thm}
The relative uniform indices are more subtle objects; besides a uniform coarse contribution from the interior they also have a metric contribution from the boundary. This complicates their interpretation. Say for example that $\olO$ is a spin manifold with boundary and bounded geometry, whose metric has uniformly positive scalar curvature. Does $\relInd(\Dirac_\Omega)$ vanish? This question and related ones will be investigated in the next chapter. The answer turns out to be negative, but also more interesting. Nonetheless, the subtler nature of the relative uniform index forbids any simple corollary of the relative uniform partitioned-manifold index theorem as a upsc-obstructions.
\chapter{Relative uniform indices and scalar curvature}\label{ChapterScalCurv}
In this last chapter we return to the geometric application that provided a guiding theme for much of the introduction, namely metrics with (uniformly) positive scalar curvature. On spin manifolds the existence of such metrics is obstructed by the index of the spin Dirac operator. 

Let $M$ be a spin manifold that has uniformly positive scalar curvature outside a closed subset $Z\subseteq M$. There is a well-known localization principle (see for example \cite{Roe2016}) for coarse indices stating that in this case the Dirac operator admits a localized index $\mathrm{Ind}^Z(\Dirac) \in K_*(C^*(Z\subseteq M))$ that gets mapped to $\mathrm{Ind}(\Dirac) \in K_*(C^*(M))$ under the inclusion map $K_*(C^*(Z\subseteq M)) \to K_*(C^*(M))$. In Section \ref{SectionCoercive} we prove a uniform analogue of this localization principle that also works in the presence of boundaries. It states that if a spin manifold $\olO$ with boundary and bounded geometry (or more generally a domain inside a spin manifold of bounded geometry) has uniformly positive scalar curvature outside a closed subset $Z\subseteq \olO$, then its relative uniform index refines to an index $\Ind^Z(\Dirac) \in K_*(D^*_u(\dO\cup Z\subseteq \olO,\dO))$. Note that even though the scalar curvature might be uniformly positive near $\dO$, the localized index is still also supported near $\dO$. The reason is this: The localization principle is derived from the coercivity condition 
$$
\|\Dirac u\|_{L^2} \gsim \|u\|_{L^2} \quad \forall \, u\in \Gamma_{cc}(\Omega\setminus Z;\slashed{S}(\olO)) \; ,
$$ 
which is provided by the Lichnerowicz formula. This condition is not necessarily satisfied if the support of $u$ intersects the boundary, due to the presence of boundary terms from Green's formula. Therefore $\dO$ also has to be treated as an area of non-coercivity. Since the localization principle for Dirac operators is the consequence of a coercivity statement, we prove it more generally for uniformly elliptic first-order operators whose minimal boundary condition is coercive with respect to a closed subspace. The focus on the minimal boundary condition is dictated by our wish to apply it to Dirac operators. On the other hand, the focus on first-order operators is a consequence of their favorable propagation properties that are not generally enjoyed by higher-order operators, see Section 5.5.

In Section \ref{SectionRhoInv} we compute the relative uniform index of the Dirac operators on spin manifolds under certain scalar curvature assumption. These indices turn out to be closely connected to uniform analogues of the higher $\rho$-invariants of the boundary. In particular, they are no obstructions to uniformly positive scalar curvature, but rather secondary measures of the geometry. Concretely, we will prove that if $\olO$ is an even-dimensional spin manifold with boundary and bounded geometry, then the localized relative uniform index of $\Dirac_\Omega$ coincides with the uniform higher $\rho$-invariant of the boundary, i.e.
$$
\Ind^{\dO}(\Dirac_\Omega) = \rho^{u}(\Dirac_{\dO}) \; \in \; K_0\left( D^*_u(\dO) \right) \; .
$$
If $\olO$ is only assumed to have upsc on the boundary, then there is a correction term coming from the lack of upsc in the interior. This correction term is a uniform analogue of the APS-index of Piazza and Schick \cite{Piazza2014}. The formula for the relative uniform index then reads as
$$
\relInd(\Dirac_\Omega) = \Ind^{APS}(\Dirac_\Omega) + \rho^{u}(\Dirac_{\dO}) \; \in \; K_0\left(D^*_u(\olO,\dO)\right) \; .
$$
Modulo a sign convention this formula yields a uniform analogue of the delocalized APS-index theorem, that was proved by Piazza and Schick \cite{Piazza2014} in the even-dimensional case, and by Xie and Yu \cite{Xie2014} in the general case. Indeed, by exactness the relative uniform index lies in the kernel of the inclusion map $K_*(D^*_u(\olO,\dO)) \to K_*\left( D^*_u(\olO) \right)$. Thus, our index formula implies that
$$
\Ind^{APS}(\Dirac_\Omega) = - \rho^{u}(\Dirac_{\dO}) \; \in \; K_0\left(D^*_u(\olO)\right) \; .
$$ 
The non-uniform version of this formula is the delocalized APS-index theorem. We only prove our version in the even-dimensional case. The difficulties in odd dimensions are discussed below.
\section{Coercivity and localized indices}\label{SectionCoercive}
As usual we let $M$ be a manifold of bounded geometry, $\Omega \subseteq M$ a domain, and $D$ as symmetric and uniformly elliptic operator over $\olO$. While $D$ generally acts on sections of a vector bundle we do not make this explicit in the notation.
\begin{definition}
Let $Z\subseteq \olO$ be a closed subset. A boundary condition $D_e$ for $D$ is said to be \emph{coercive with respect to} $Z$ if there exists a constant $c>0$ such that
$$
\|D_eu\| \geq c \| u \| \quad \forall \; u\in\dom(D_e), \, \supp(u)\cap Z = \emptyset \; .
$$
If $Z$ is empty, we will say that $D_e$ is \emph{coercive}. If there exists a compact $Z$ such that $D_e$ is coercive with respect to $Z$, we will say that $D_e$ is \emph{coercive at infinity}.
\end{definition}
\begin{example}[Dirac operators on manifolds with uniformly positive scalar curvature]\label{CoerciveUPSC}
Suppose that $\olO$ is a spin manifold with boundary and bounded geometry. The Lichnerowicz formula asserts that
\begin{equation}\label{LichnerowiczBdry}
\Dirac^2 u = \nabla^*\nabla u + \frac{\mathrm{scal}}{4} u \quad \forall \, u \in \Gamma_{cc}(\olO,\slashed{S}(\olO)) \; .
\end{equation}
If $\olO$ has uniformly positive scalar curvature, i.e. $\mathrm{scal}\geq 4c>0$, then $\Dirac^2 u \geq c\|u\|$ for all $u\in \Gamma_{cc}(\olO,\slashed{S}(\olO))$. It follows that the minimal boundary condition $\Dirac_{\min}$ is coercive. This generalizes to the case where upsc holds only away from a closed subset, i.e. if $\olO$ has upsc outside a closed $Z\subseteq \olO$, then $\Dirac_{\min}$ is coercive with respect to $Z$.

For sections that do not vanish on the boundary \eqref{LichnerowiczBdry} is no longer valid due to the appearance of boundary terms from Green's formula. Thus coercivity is generally not to be expected for any other boundary condition. It turns out that under the additional assumption that the mean curvature of the boundary is non-negative both the APS- and the chirality boundary conditions are coercive (see \cite[Section 6]{Grosse2014} and \cite[Section 8]{Baer2022}).
\end{example}
\begin{example}[Dirichlet Laplacian on finite-width manifolds]
The Dirichlet Laplacian on a manifold of finite width is coercive as a consequence of the Poincaré equality. In the case of compact domains this is classical, for manifolds with boundary and bounded geometry this was proved in much greater generality by Ammann, Große and Nistor in \cite{Ammann2019}.
\end{example}
The main result of this section is a localization principle for first-order operators. It states that if the minimal boundary condition is coercive with respect to some $Z\subseteq \olO$, then there is a relative uniform index supported near $Z\cup \dO$.
\begin{thm}\label{LocInd}
Let $D$ be a first-order uniformly elliptic differential operator. Suppose that $D_{\min}$ is coercive with respect to a closed subset $Z\subseteq \olO$. Then, the class $[D]\in K_*^{u}(\olO,\partial\Omega)$ has a representative that is invertible modulo $D^*_u(\partial\Omega\cup Z\subseteq \olO,\partial \Omega)$. Consequently, it admits an index class
$$
\Ind^Z (D) \; \in \; K_*\left( D^*_u(\partial\Omega\cup Z\subseteq \olO,\partial \Omega) \right) 
$$
that is mapped to $\relInd(D)$ under the inclusion map 
$$
K_*(D^*_u(\partial\Omega\cup Z\subseteq\olO,\partial \Omega))\to K_*(D^*_u(\olO,\partial\Omega)) \; .
$$
\end{thm}
We first prove a lemma.
\begin{lemma}\label{CoerciveSupport}
Let $D$ be a first-order uniformly elliptic differential operator. Suppose that $D_{\min}$ is coercive with respect to a closed subset $Z\subseteq \olO$, i.e. that 
$$
\|Du\|\geq c\|u\| \quad \forall \, u\in \tilde{H}^1(\Omega), \, \supp(u)\cap Z = \emptyset
$$
for some $c>0$. Suppose $\varphi\in C_0(\RR)$ has support inside $(-c,c)$. Then 
$$
\varphi(D_e)\in C^*_u(Z\cup \partial \Omega \subseteq \Omega)
$$ 
for any self-adjoint boundary condition $D_e$.
\end{lemma}
We draw attention to the fact that $\varphi(D_e)$ is only claimed to be uniformly locally compact over the interior, not up to the boundary.
\begin{proof}
The proof is analogous to that of the corresponding non-uniform statement in the absence of boundaries, see for example \cite[Lemma 2.3]{Roe2016}. Fix $\varepsilon>0$, and let $\psi\in C_0(\RR)$ be such that $\|\varphi-\psi\|<\varepsilon
$ and $\supp(\hat{\psi})\subseteq [-r,r]$ for some $r\geq c> 0$. Pick a Lipschitz function $g: \, \olO \to [0,1]$ such that $g=1$ on $N_{4c_Dr}(\partial \Omega \cup Z)$ and $g=0$ outside $N_{5c_Dr}(\partial \Omega\cup Z)$. Here $c_D$ denotes the propagation speed of $D$. Then write
$$
\varphi(D_e) = g \psi(D_e) g + (1-g) \psi(D_e) g + \psi(D_e) (1-g) + \left( \varphi(D_e) - \psi(D_e) \right) \; .
$$
The first summand has finite propagation by Corollary \ref{finPropFourier}, is uniformly locally compact over $\Omega$,\footnote{Approximate $\psi$ by functions $\psi'$ in $\cF^{-\beta}$ for some $\beta>0$. Then $\psi'(D_e)$ is uniformly locally compact over $\Omega$ by Proposition \ref{MappingPropertiesT} and Corollary \ref{CorNegOrder}. As the limit of uniformly locally compact operators, $\psi(D_e)$ is also uniformly locally compact.} and it is supported near $\partial \Omega \cup Z$. We will show that the remaining terms are $\lsim \varepsilon$. This is clear for the last summand.

To treat the remaining summands we will show that if $\|\psi(D_e) f\| \lsim \|f\|\cdot \sup_{|x|\geq c} |\psi(x)|$ whenever $f\in C_{cc}(\Omega)$ is supported in distance at least $4c_Dr$ from $\partial\Omega \cup Z$. Since $1-g$ has support of distance $\geq 4c_D r$ from $\partial \Omega \cup Z$ this will yield the claim, since 
$$
\sup_{|x|\geq r} |\psi(x)| \leq \varepsilon + \sup_{|x|\geq r} |\varphi(x)| = \varepsilon \; .
$$
According to \cite[Lemma 2.6]{Roe2016} it suffices to prove this for even functions $\psi$, and $d(\supp(f),\partial \Omega \cup Z)\geq 2c_Dr$. Then, the proof of \cite[Lemma 2.5]{Roe2016} applies verbatim. Indeed, letting $\Omega_j$ denote the complement of $N_{jc_Dr}(\partial \Omega \cup Z)$ in $\Omega$, $j=1,2$, we consider the operator $D^2$ on $L^2(\Omega_1)$ with domain $\Gamma_{cc}(\Omega_1)$. Due to our coercivity assumption it is strictly positive, and thus has a self-adjoint Friedrichs extension $T=(D^2)_F$ satisfying $T\geq c^2 I$. For $u\in \Gamma_{cc}(\Omega_1)$ it holds that $\sqrt{T}u=\pm Du$. In particular $\sqrt{T}$ has the same propagation speed as $D$. Thus we can conclude that
$$
\cos(t\sqrt{T})u = \cos(tD_e)u \quad \forall \; u\in L^2(\Omega_2), \, |t|\leq r \; .
$$
Using Fourier inversion and the assumption that $\psi$ is even, we conclude that
$$
\psi(D_e) f = \int_{-r}^r \hat{\psi}(t) \cos(tD_e) f \, dt = \int_{-r}^r \hat{\psi}(t) \cos(t\sqrt{T}) f \, dt = \psi(\sqrt{T}) f \; .
$$
Appealing to the fact that $\sqrt{T}\geq c$ concludes the proof.
\end{proof}
With this fact in hand we can proceed to prove Theorem \ref{LocInd}.
\begin{proof}[Proof of Thm. \ref{LocInd}]
We wish to apply Lemma \ref{CoerciveSupport}, but since $D$ generally does not admit a self-adjoint extension that also defines an invertible element in $D^*_u(\olO)/D^*_u(\olO;\partial\Omega)$ we cannot apply it directly. Instead consider the operator $\hat{D}=\begin{pmatrix}
0 & D \\ D & 0
\end{pmatrix}$ acting on sections of $E\oplus E$. It still satisfies $\|\hat{D}u\|\geq c\|u\|$ for all $u \in \Gamma_{cc}(\Omega\setminus Z)$. Then, consider the self-adjoint extension $\hat{D}_e=\begin{pmatrix}
0 & D_{\max} \\ D_{\min} & 0
\end{pmatrix}$. Fix a normalizing function $\chi$ so that $1-\chi^2$ is supported inside $(-c,c)$. Then Lemma \ref{CoerciveSupport} yields that $\chi(\hat{D}_e)$ is invertible modulo $C^*_u(Z\cup\partial\Omega\subseteq \Omega;\hat{\rho})$. Here we chose to make explicit the representation $\hat{\rho}=\rho\oplus\rho$ on $L^2(\Omega,E\oplus E)$. Let $P$ denote the projection on $L^2(\Omega,E\oplus E)$ onto $L^2(\Omega,E)$. Then, compression by $P$ constitutes a map\footnote{This is because the inclusion $V: \, (L^2(\Omega,E),\rho)\hookrightarrow (L^2(\Omega,E\oplus E),\hat{\rho})$ uniformly boundedly covers the identity on $\olO$, and $P$ is given by $VV^*$.}
$$
T\mapsto PTP , \, C^*_u(Z\cup\partial\Omega\subseteq \Omega;\hat{\rho}) \to C^*_u(Z\cup\partial\Omega \subseteq \Omega;\rho) \; .
$$
We conclude that $I-P\chi(\hat{D}_e)^2 P \in C^*_u(Z\cup\partial\Omega \subseteq \Omega;\rho)$. On the other hand $(I-P)\chi(\hat{D}_e)P$ is a representative of $[D]\in K^{u}_*(\olO,\partial\Omega)$, so that
$$
I -\left( (I-P)\chi(\hat{D_e})P\right)^*\left( (I-P)\chi(\hat{D_e})P \right) = I- P\chi(\hat{D_e}) P \; \in \fD^{u}_\rho(\olO,\partial\Omega) \; .
$$
We conclude that $I-P\chi(\hat{D}_e)^2 P$ lies in 
$$
C^*_u(Z\cup\partial\Omega\subseteq \Omega) \cap \fD^{u}(\olO;\partial\Omega)= D^*_u(Z\cup\partial\Omega\subseteq \olO,\partial\Omega) \; ,
$$ 
meaning that $[D]$ has a representative that is invertible modulo $D^*_u(Z\cup\partial\Omega\subseteq \olO,\partial\Omega)$. Thus, $[D]\in K^{u}_*(\olO,\partial\Omega)$ lies in the image of the map
$$
K_{*+1}\left( D^*_u(\olO)/D^*_u(Z\cup\partial\Omega\subseteq \olO,\partial\Omega) \right) \longrightarrow K_{*+1}\left( D^*_u(\olO)/D^*_u(\olO,\partial\Omega) \right) \cong K_*^{u}(\olO;\partial\Omega)
$$
induced by inclusion. It fits into a commutative diagram \\
\centerline{ \xymatrix{
K_{*+1}\left( D^*_u(\olO)/D^*_u(Z\cup\partial\Omega\subseteq\olO,\partial\Omega) \right) \ar[r] \ar[d]_\partial & K_*^{u}(\olO;\partial\Omega) \ar[d]^{\relInd} \\
K_*\left( D^*_u(Z\cup\partial\Omega\subseteq \olO,\partial\Omega) \right) \ar[r] & K_*\left( D^*_u(\olO,\partial\Omega)\right)
}}
where the bottom horizontal arrow is also induced by inclusion. Defining 
$$
\Ind^Z(D) \in K_*\left( D^*_u(Z\cup\partial\Omega\subseteq\olO,\partial\Omega) \right)
$$ 
as the image of 
$$
[(I-P)\chi(\hat{D}_e)P] \in K_{*+1}\left( D^*_u(\olO)/D^*_u(Z\cup\partial\Omega\subseteq\olO,\partial\Omega) \right)
$$ 
under the boundary map concludes the proof.
\end{proof}
If the subspace $Z$ is compact, then being near $\partial \Omega\cup Z$ is equivalent to being near $\partial \Omega$. Thus any boundary condition that is coercive at infinity has an index class in $K_*(D^*_u(\partial\Omega\subseteq \olO,\partial \Omega))$. Recalling from Proposition \ref{IsoStructureLocalized} that the inclusion induces an isomorphism $K_*(D^*_u(\partial\Omega))\cong K_*(D^*_u(\partial\Omega\subseteq \olO,\partial \Omega))$ if $\dO$ is sufficiently regular, we obtain the following.
\begin{cor}
Suppose that $D_{\min}$ is coercive at infinity, and that $\Omega$ is a bounded-geometry domain. Then, $D$ admits an index class
$$
\Ind^{\partial\Omega} (D) \; \in \; K_*\left( D^*_u(\partial \Omega) \right) 
$$
that is mapped to $\relInd(D)$ under the inclusion map $K_*(D^*_u(\partial \Omega))\to K_*(D^*_u(\olO,\partial\Omega))$. In this case
$$
\Ind(D) = 0 \; \in \; K_*\left( C^*_u(\olO,\dO) \right) \; .
$$
\end{cor}
\begin{proof}
The only part left to show is that $\Ind(D) = 0$ in $K_*\left( C^*_u(\olO,\dO) \right)$. Recall from Section 3.5 that this uniform coarse index map factors as the composition
$$
K^{u}_*(\olO,\dO) \xlongrightarrow{\relInd} K_*\left(D^*_u(\olO,\dO)\right) \longrightarrow K_*\left(  D^*_u(\olO,\dO)/D^*_u(\dO\subseteq \olO,\dO\right)
= K_*\left(C^*_u(\olO,\dO)\right) \; .
$$
Since $\relInd(D)$ is in the image of $D^*_u(\dO)\cong D^*_u(\dO\subseteq \olO,\dO)$, its image under the quotient map $K_*\left(D^*_u(\olO,\dO)\right) \to K_*\left(  D^*_u(\olO,\dO)/D^*_u(\dO\subseteq \olO,\dO\right)$ vanishes, and the claim follows.
\end{proof}
Recall from Example \ref{CoerciveUPSC} above that the minimal boundary condition for the spinor Dirac operator is coercive with respect to $Z$ if the scalar curvature is uniformly positive outside $Z$. Thus the uniform coarse index $\Ind(\Dirac) \in K_*(C^*_u(\olO,\dO))$ obstructs the existence of a upsc-metric.
\begin{cor}
Let $\Omega$ be a domain in a spin manifold of bounded geometry,, and suppose the scalar curvature on $\olO$ is uniformly positive outside a compact subset of $\olO$. Then,
$$
\Ind(\Dirac) = 0 \; \in \; K_*\left( C^*_u(\olO,\dO) \right) \; .
$$
\end{cor}
Thus the uniform coarse index $\Ind(\Dirac) \in K_*\left( C^*_u(\olO,\dO) \right)$ functions as an obstruction to the existence of a upsc-metric of bounded geometry in the uniform coarse equivalence class of $\olO$. However, we recall that $K_*\left( C^*_u(\olO,\dO) \right)=0$ if $\olO$ has finite width, so that this obstruction is only ever non-trivial if $\olO$ has infinite width. Even in this case it seems to be not especially interesting, since all information near the boundary is forgotten. One might therefore hope that the relative uniform index of the Dirac operator provides a more interesting obstruction. Surprisingly, it turns out to be no obstruction at all, though it is interesting nonetheless. 
\section{Uniform higher $\rho$-invariants and the delocalized APS index theorem}\label{SectionRhoInv}
Let $\Sigma$ be an $n$-dimensional spin manifold of bounded geometry. If $\Sigma$ carries a metric of uniformly positive scalar curvature, then $\Dirac_\Sigma$ is invertible. Not only does its uniform coarse index $\Ind(\Dirac_\Sigma)\in K_*(C^*_u(\Sigma))$ vanish, but $\chi(\Dirac_\Sigma)$ is an involution in $D^*_u(\Sigma)$ whenever the normalizing function $\chi$ is such that $1-\chi^2$ is supported in $(-c,c)$ for $c>0$ such that $\Dirac^2_\Sigma\geq c^2$. Indeed, in this case $\chi(D_e)$ is self-adjoint and squares to the identity. Thus it defines a class 
$$
\rho^{u}(\Dirac_\Sigma):= \left[ \chi(\Dirac_\Sigma) \right] \; \in \; K_{n+1}(D^*_u(\Sigma)) \; .
$$
It is called the \emph{uniform higher $\rho$-invariant} of the upsc-metric on $\Sigma$. The (non-uniform) higher $\rho$-invariant is an important tool in the study of the space of upsc-metrics on $\Sigma$, and it is to be expected that its uniform counterpart plays an analogous role. However, the primary use of the uniform higher $\rho$-invariants here is that they appear in index formulas for relative uniform indices.
\begin{thm}\label{RelIndRho}
Let $\olO$ be an even-dimensional spin manifold with boundary and bounded geometry that carries a upsc metric having product form near the boundary. Then
$$
\Ind^{\dO}(\Dirac_\Omega) = \rho^{u}(\Dirac_{\dO}) \; \in \; K_{0}\left( D^*_u(\dO)\right) \; .
$$ 
\end{thm}
\begin{proof}
By use of a collar neighborhood it suffices to treat the case $\olO=\RR_{\geq 0}\times \Sigma$, where $\Sigma$ is an odd-dimensional spin manifold of bounded geometry. Then, the Dirac operator $\Dirac$ on $\RR_{\geq 0}\times\Sigma$ takes the form
$$
\Dirac = \begin{pmatrix}
0 & \Dirac^- \\ \Dirac_+ & 0
\end{pmatrix} = \begin{pmatrix}
0 & -\partial_t + \Dirac_\Sigma \\ \partial_t + \Dirac_\Sigma  & 0
\end{pmatrix} \; .
$$
We consider the boundary condition
$$
\Dirac_e = \begin{pmatrix}
0 & \Dirac^-_{\min} \\ \Dirac_{\max}^+ & 0
\end{pmatrix} \; .
$$ 
By assumption of upsc there is $c>0$ such that $\Dirac^+_{\max}\Dirac_{\min}^-\geq c^2$ and $\Dirac_\Sigma^2\geq c^2$. Pick a normalizing function $\chi$ such that $1-\chi^2$ is supported inside $(-c,c)$. Write $\chi(\Dirac_e)=\begin{pmatrix}
0 & U^* \\ U & 0
\end{pmatrix}$. Then, 
$$
\begin{pmatrix}
 U^*U & 0 \\ 0 & UU^*
\end{pmatrix} = \chi^2(\Dirac_e) = I - P_{\ker(\Dirac_e)} = \begin{pmatrix}
 I-P_{\ker(\Dirac^+_{\max})} & 0 \\ 0 & I
\end{pmatrix}  \; .
$$
We conclude that $U$ is a partial isometry with kernel $\ker(\Dirac_{\max}^+)^\perp$ and trivial cokernel. Moreover, $I-\chi(\Dirac_e)^2 \in D^*_u(\Sigma\subseteq\RR_{\geq 0}\times\Sigma,\Sigma)$, meaning that the projection $P_{\ker(\Dirac_{\max}^+)}$ is in $D^*_u(\Sigma\subseteq\RR_{\geq 0}\times\Sigma,\Sigma)$. Now, recall that the map $\Ind^\Sigma$ is defined as 
$$
K_{1}\left( D^*_u(\RR_{\geq 0}\times\Sigma)/D^*_u(\Sigma\subseteq\RR_{\geq 0}\times\Sigma,\Sigma) \right) \xlongrightarrow{\partial} K_{0}\left( D^*_u(\Sigma\subseteq\RR_{\geq 0}\times\Sigma,\Sigma) \right) \cong K_{0}\left(D^*_u(\Sigma) \right) \; .
$$
The operator $\chi(\Dirac_e)$ is first identified with the K-theory class 
$$
[U] \in K_{1}\left( D^*_u(\RR_{\geq 0}\times\Sigma)/D^*_u(\Sigma\subseteq\RR_{\geq 0}\times\Sigma,\Sigma) \right) \; .
$$
Its image under the boundary map is the class $[P_{\ker(\Dirac_{\max}^+)}] \in K_0( D^*_u(\Sigma\subseteq\RR_{\geq 0}\times\Sigma,\Sigma) )$. The kernel of $\Dirac_{\max}^+=\left(\partial_t + \Dirac_\Sigma\right)_{\max}$ consists of all elements of the form $u(t,\cdot) = e^{-t\Dirac_\Sigma} u_0$, where $u_0\in L^2(\Sigma;\slashed{S}(\Sigma))$ is in the image of the projection $P_0 = \mathbf{1}_{\RR_{\geq 0}}(\Dirac_\Sigma)$. Observe that $P_0 = \frac{1}{2}(I+\chi(\Dirac_\Sigma))$ represents the uniform higher $\rho$-invariant $\rho^{u}(\Dirac_\Sigma)$. In spirit this concludes the proof, as we have found $\Ind^\Sigma(\Dirac)$ and $\rho^{u}(\Dirac_\Sigma)$ to be projections onto subspaces that can be canonically identified. In practice we have to put in some legwork to show that this identification can be accomplished in a way that gives rise to an equality of K-theory classes.

Fix some small $\varepsilon>0$. Because $P_{\ker(\Dirac_{\max}^+)} \in D^*_u(\Sigma\subseteq \RR_{\geq 0}\times\Sigma,\Sigma)$ is the limit of operators supported near $\Sigma$, there exists some $R\geq 0$ such that $\|P_{\ker(\Dirac_{\max}^+)}u\|<\varepsilon\|u\|$ whenever $d(\supp(u),\Sigma) \geq R$. Then, let $P_R$ be the projection onto the subspace containing all functions $u(t,\cdot)=\psi_R\cdot e^{-t\Dirac_\Sigma}u_0$ for $u_0\in\im(P_0)$, where $\psi_R=\mathbf{1}_{[0,R]}(t)$. Then 
$$
\|(P_{\ker(\Dirac_{\max}^+)} - P^R)u\|=\|P_{\ker(\Dirac_{\max}^+)}u\|<\varepsilon\|u\|
$$
if $d(\supp(u),\Sigma)\geq R$. Suppose on the other hand that $\supp(u)\subseteq [0,R]\times\Sigma$. Write $u=\psi_R e^{-t\Dirac_\Sigma}u_0 + u^\perp$, where $u^\perp$ is uniquely determined by the requirement that $\langle \psi_R e^{-t\Dirac_\Sigma} v_0,u^\perp\rangle=0$ for all $v_0\in\im(P_0)$. Note that in particular $u^\perp=\psi_R u^\perp$, so that it follows that $\langle e^{-t\Dirac_\Sigma} v_0,u^\perp\rangle=0$ for all $v_0\in\im(P_0)$ as well. Thus
\begin{align*}
\left\| (P_{\ker(\Dirac_{\max}^+)} - P^R)u \right\| &= \left\| (P_{\ker(\Dirac_{\max}^+)} - I)\psi_R e^{-t\Dirac_\Sigma}u_0 \right\| \\
&= \left\| P_{\ker(\Dirac_{\max}^+)} e^{-t\Dirac_\Sigma}u_0 - P_{\ker(\Dirac_{\max}^+)}\left( (1-\psi_R)e^{-t\Dirac_\Sigma}u_0 \right) - \psi_R e^{-t\Dirac_\Sigma}u_0 \right\| \\
&= \left\| e^{-t\Dirac_\Sigma}u_0 - P_{\ker(\Dirac_{\max}^+)}\left( (1-\psi_R)e^{-t\Dirac_\Sigma}u_0 \right) - \psi_R e^{-t\Dirac_\Sigma}u_0 \right\| \\
&\leq 2\cdot \left\| (1-\psi_R)e^{-t\Dirac_\Sigma}u_0 \right\| \; .
\end{align*}
It holds that 
$$
\frac{\varepsilon^2}{4} \int_{t=0}^R e^{-2t\lambda} \, dt \geq \int_{t=R}^\infty e^{-2t\lambda} \, dt
$$
for all $\lambda\geq c$ whenever $R$ is sufficiently large, say $R\gsim \log(\varepsilon^{-2}) c^{-1}$. Let $E_\lambda$ denote the spectral measure of $\Dirac_\Sigma$ (see Theorem \ref{funcCalc}). Then, 
\begin{align*}
\left\| (1-\psi_R)e^{-t\Dirac_\Sigma}u_0 \right\|^2 &=  \int_{t=R}^\infty \left\| e^{-t\Dirac_\Sigma}u_0 \right\|^2 \, dt \\
&= \int_{t=R}^\infty \left( \int_{\lambda=c}^\infty e^{-t\lambda} d\langle E_\lambda u_0,u_0\rangle \right) \, dt \\
&= \int_{\lambda=c}^\infty \left( \int_{t=R}^\infty e^{-t\lambda}   \, dt \right) d\langle E_\lambda u_0,u_0\rangle \\
&\leq \frac{\varepsilon^2}{4} \int_{\lambda=c}^\infty \left( \int_{t=0}^R e^{-t\lambda}   \, dt \right) d\langle E_\lambda u_0,u_0\rangle \\
&= \frac{\varepsilon^2}{4} \left\| \psi_R e^{-t\Dirac_\Sigma}u_0 \right\|^2 \\
&\leq \frac{\varepsilon^2}{4} \|u\|^2 \; .
\end{align*}
The interchange of order of integration is justified by Tonelli's theorem. Thus
$$
\left\| (P_{\ker(\Dirac_{\max}^+)} - P^R)u \right\| \leq \varepsilon \|u\| \; .
$$
Overall, we conclude that $P_{\ker(\Dirac_{\max}^+)}$ and $P_R$ are projection whose difference has norm strictly less than $1$, so that they are homotopic (see for example \cite[Proposition 4.1.7]{HigsonRoe2000}).

Let $H_R=L^2([0,R]\times \Sigma)$, and denote by $V_R$ inclusion map $H_R\hookrightarrow L^2(\RR_{\geq 0} \times \Sigma)$. Then $V_R$ uniformly boundedly covers the inclusion $[0,R]\times\Sigma \hookrightarrow \RR_{\geq 0}\times\Sigma$. Also consider $H_\Sigma = L^2(\Sigma)$. The isometry $V_\Sigma: \, H_\Sigma \to H_R$, $(V_\Sigma v)(t,x)=R^{-1/2}v(x)$ uniformly covers the inclusion $\Sigma \hookrightarrow [0,R]\times\Sigma$. The conjugation maps corresponding to these isometries induce the isomorphisms
$$
K_*\left(D^*_u(\Sigma) \right) \xlongrightarrow{\sim} K_*\left( D^*_u([0,R]\times\Sigma,\Sigma) \right) \xlongrightarrow{\sim} K_*\left( D^*_u(\Sigma \subseteq \RR_{\geq 0}\times\Sigma,\Sigma \right) \; .
$$
Now, observe that $P_R$ is in the image of $\Ad(V_R)$, namely it is the image of the projection onto the subspace of $H_R$ containing all functions $u(t,\cdot)=e^{-t\Dirac_\Sigma}u_0$, $u_0 \in\im(P_0)$, $t\in [0,R]$. Denote that projection by $P_R$ as well. Now, there is a homotopy of projections $P_R^s$, $s\in[0,1]$, where $P_R^s$ projects onto the subspace of $H_R$ containing all functions $u(t,\cdot)=e^{-(1-s)t\Dirac_\Sigma}u_0$, $u_0 \in\im(P_0)$. Then, $P_R^0=P_R$, while $P_R^1=V_\Sigma P_0 V^*_\Sigma$. Because $V_\Sigma$ and $V_R$ induce isomorphisms, we conclude that
$$
[P_{\ker(\Dirac_{\max}^+)}]=[P_R]= [P_0] \; \in \; K_*\left(D^*_u(\Sigma) \right) \cong K_*\left( D^*_u(\Sigma \subseteq \RR_{\geq 0}\times\Sigma,\Sigma \right) \; .
$$
Since we have already noted that the former projection represents $\Ind^\Sigma(\Dirac)$ and the latter $\rho^{u}(\Dirac_\Sigma)$, this concludes the proof.
\end{proof}
We can derive a form of bordism invariance for the uniform higher $\rho$-invariant. Suppose that $\olO$ is a spin manifold with boundary and bounded geometry that has uniformly positive scalar curvature and product structure near the boundary. Then Theorem \ref{RelIndRho} implies in particular that $\rho^{u}(\Dirac_{\dO})$ in in the image of the relative index map. Exactness of 
$$
K^{u}_*(\olO,\dO) \xlongrightarrow{\relInd} K_*\left(D^*_u(\olO,\dO)\right) \longrightarrow K_*\left(D^*_u(\olO)\right)
$$
then implies the following.
\begin{cor}\label{BordismInvRho}
Let $\olO$ be an even-dimensional spin manifold with boundary and bounded geometry that has uniformly positive scalar curvature and product structure near the boundary. Then, 
$$
\rho^{u}(\Dirac_{\partial\Omega}) \in \ker\left( K_0\left(D^*_u(\dO)\right) \to K_0\left(D^*_u(\olO)\right) \right) \; .
$$
\end{cor}
Again, we interpret this statement as a form of bordism invariance: The uniform $\rho$-invariant of a upsc-null-bordant manifold vanishes upon being pushed forward to the bordism. It is admittedly unnatural that one first has to push forward to the bordism, though note that the uniform $\rho$-invariant itself can never vanish. The statement becomes more natural if one includes actions by some group $\Gamma$, and considers the push-forwards of the (uniform) higher $\rho$-invariants to $E\Gamma$. If one requires the upsc-null-bordism to extend the reference map of the boundary, then one obtains a genuine vanishing statement for the push-forward of the higher $\rho$-invariant of a upsc-null-bordant manifold. We refer to \cite{Piazza2014} for the non-uniform setting.

We move on to the case where we require $\olO$ to have uniformly positive scalar curvature only on the boundary. In this case we can still consider the uniform higher $\rho$-invariant of the boundary. Due to the potential lack of upsc in the interior, we expect $\relInd(\Dirac)$ to have contributions from the interior. To get an idea of what to expect consider the following diagram \\
\centerline{\xymatrix{
 & K^{u}_{*+1}(\olO,\dO) \ar[d]_{\relInd}\ar[r]^-\partial & K^{u}_{*}(\dO)\ar[d]^{\Ind} \\
 K_{*+1}\left(C^*_u(\olO)\right)\oplus K_{*+1}\left( D^*_u(\dO) \right) \ar[r] & K_{*+1}\left( D^*_u(\olO,\dO) \right) \ar[r]^-\delta \ar[d] & K_{*}\left( C^*_u(\dO)\right) \\
 & K_{*+1}\left( D^*_u(\olO) \right) &
}}
The upper-right square is the content of Theorem \ref{IndexCommAlgMV}, the horizontal row comes from the algebraic Mayer-Vietoris sequence of Theorem \ref{AlgebraicMV}, and the horizontal row from the exact sequence defining the relative uniform index map. The image of the class $[\Dirac_{\Omega}] \in K_{*+1}^{u}(\olO,\dO)$ under the boundary map on uniform K-homology is the class $[\Dirac_{\dO}]\in K_{*}^{u}(\dO)$. Because $\dO$ has upsc the uniform coarse index $\Ind([\Dirac_{\dO}])$ vanishes. Thus
$$
\delta\left(\relInd( [\Dirac_{\olO}])\right) = \Ind\left(  [\Dirac_{\dO}]\right) = 0\; .
$$
We conclude by exactness that the relative uniform index of $[\Dirac_{\olO}]$ is the sum of an element of $K_{*+1}(D^*_u(\dO))$ and an element of $K_{*+1}(C^*_u(\olO))$. We have a good guess of what the former should be, namely the uniform $\rho$-invariant $\rho^{u}(\Dirac_{\dO})$. The latter turns out to be a uniform analogue of the APS-index considered by Piazza and Schick \cite{Piazza2014}. It is defined as follows: Since we assume $\olO$ to have product structure near the boundary, we may attach a cylinder (endowed with the product metric) to $\dO$ to form the manifold $M=\olO \cup_{\dO}(\RR_{\leq 0}\times\dO)$ of bounded geometry. If we assume upsc on $\dO$, then $M$ has upsc outside the closed subspace $\olO\subseteq M$. By the localization principle the Dirac operator on $M$ admits a localized index class
$$
\Ind^{APS}(\Dirac_\Omega) := \Ind^{\olO}(\Dirac_M) \; \in \; K_{n+1}\left( C^*_u(\olO\subseteq M)\right) \cong K_{n+1}\left( C^*_u(\olO)\right) \; .
$$  
We call it the \emph{uniform APS-index}. It provides the other component to the relative uniform index of $\Dirac_\Omega$.
\begin{thm}\label{DelocIndexFormula}
Let $\olO$ be an even-dimensional spin manifold with boundary and bounded geometry. Assume that the metric on $\olO$ has product structure near the boundary, and that the scalar curvature of $\dO$ is uniformly positive. Then,
$$
\relInd(\Dirac_\Omega) = \Ind^{APS}(\Dirac_\Omega) + \rho(\Dirac_{\dO}) \; \in \; K_{0}\left( D^*_u(\olO,\dO)\right) \; .
$$
\end{thm}
Strictly speaking the right-hand side of the the index formula of Theorem \ref{DelocIndexFormula} is given by the push-forwards of $\Ind^{APS}(\Dirac_\Omega)$ and $\rho^{u}(\Dirac_{\dO})$ under the respective inclusion maps
$$
K_*\left(C^*_u(\olO)\right) \to K_*\left( D^*_u(\olO,\dO) \right) \quad , \quad K_*\left(D^*_u(\dO)\right) \to K_*\left( D^*_u(\olO,\dO) \right) \; .
$$
To keep notation uncluttered we have kept these push-forward maps implicit, and we will do so in the following. 
\begin{proof}
Write $\Sigma:= \dO$, and $\Sigma_\infty := \RR_{\leq 0} \times \Sigma \subseteq M$. Let $\chi$ be a normalizing function such that $1-\chi^2$ is supported in $(-c,c)$, where $c>0$ is such that $\Dirac_\Sigma^2\geq c^2$. Then, $\chi(\Dirac_M) \in D^*_u(M)$ is an involution modulo $C^*_u(\olO\subseteq M)$. Let $\psi_{\olO}$ and $\psi_{\Sigma_\infty}$ denote the characteristic functions of $\olO\subseteq M$ and $\Sigma_\infty $ respectively. Then, $\psi_{\olO}$ and $\psi_{\Sigma_\infty}$ are mutually orthogonal projections summing to the identity. 

We claim that $\psi_{\olO} \chi(\Dirac_M) \psi_{\Sigma_\infty}$ and $\psi_{\Sigma_\infty} \chi(\Dirac_M) \psi_{\olO}$ are in $D^*_u(\Sigma\subseteq M,\Sigma)$. We show this for the former, the latter is treated analogously. Let $\chi_n$ be a sequence of normalizing functions with compactly supported Fourier transforms that converge to $\chi$. Then, each $\chi_n(\Dirac_M)$ has some finite propagation $R_n$, so that if $f\in C_0(\Sigma_\infty)$ and $g\in C_0(\olO)$ are each supported in distance $\geq R_n$ from $\Sigma$, then $g\left( \psi_{\olO} \chi_n(\Dirac_M) \psi_{\Sigma_\infty} \right)f=0$. Thus $\psi_{\olO} \chi_n(\Dirac_M) \psi_{\Sigma_\infty}$ is supported near $\Sigma$. Since $\psi_{\olO}$ and $\psi_{\Sigma_\infty}$ are in $D^{*}_u(M)$, so is $\psi_{\olO} \chi_n(\Dirac_M) \psi_{\Sigma_\infty}$. Thus it remains to prove that $\psi_{\olO} \chi(\Dirac_M) \psi_{\Sigma_\infty} \sua 0$ over $M\setminus \Sigma$. To that end fix $L,R\geq 0$ and $\delta>0$, and let $f\in \LLip_R(M)$ have distance $\geq \delta$ from $\Sigma$. Let $\chi'$ be a normalizing function so that $\chi'(\Dirac_M)$ has propagation $<\delta$. Then, $\left(\psi_{\olO} \chi'(\Dirac_M) \psi_{\Sigma_\infty}\right) f=0$. On the other hand $\chi_n(\Dirac_M)-\chi'(\Dirac_M) \in C^*_u(M)$.\footnote{By an approximation argument it suffices to consider the case $\chi_n-\chi' \in \cF^{-\beta}$ for some $\beta>0$. Then $(\chi_n(\Dirac_M)-\chi'(\Dirac_M))f \sua 0$ for $f\in C_0(M)$ by Corollary \ref{CorNegOrder}. Alternatively appeal to \cite[Section 7]{Engel2018}.} We conclude that
$$
\left(\psi_{\olO} \chi_n(\Dirac_M) \psi_{\Sigma_\infty}\right) f = \psi_{\olO}\left( \chi_n(\Dirac_M)-\chi'(\Dirac_M) \right)\left( \psi_{\Sigma_\infty} f \right) \sua 0 \; ,
$$
and $f\left( \psi_{\olO} \chi_n(\Dirac_M) \psi_{\Sigma_\infty} \right)\sua 0$ is shown analogously. With the help of Lemma \ref{LLipStratified} we conclude that $\psi_{\olO} \chi(\Dirac_M) \psi_{\Sigma_\infty} \in D^*_u(\Sigma\subseteq M,\Sigma)$.

Thus, with respect to the decomposition into the images of $\psi_{\olO}$ and $\psi_{\Sigma_\infty}$ we may write
$$
\chi(\Dirac_M) = \begin{pmatrix}
\psi_{\olO} \chi(\Dirac_M) \psi_{\olO} & 0 \\ 0 & \psi_{\Sigma_\infty} \chi(\Dirac_M) \psi_{\Sigma_\infty}
\end{pmatrix} \; \mathrm{mod} \; D^*_u(\Sigma\subseteq M , \Sigma) \; .
$$
Since $\chi(\Dirac_M)$ is an involution modulo $C^*_u(\olO\subseteq M)$, we deduce that $\psi_{\olO} \chi(\Dirac_M) \psi_{\olO} \in D^*_u(\olO)$ is an involution modulo $D^*_u(\olO,\Sigma)$, and $\psi_{\Sigma_\infty} \chi(\Dirac_M) \psi_{\Sigma_\infty} \in D^*_u(\Sigma_\infty)$ one modulo $D^*_u(\Sigma \subseteq \Sigma_\infty,\Sigma)$. Consider the commutative diagram\\
 \centerline{\xymatrix{
K_{1}\left(D^*_u(\olO)/D^*_u(\olO,\Sigma)\right)\ar[d]_{\partial}\ar@{}[r]|{\oplus} & K_{1}\left(D^*_u(\Sigma_\infty)/D^*_u(\Sigma\subseteq\Sigma_\infty,\Sigma)\right)\ar[d]^{\partial} \ar[r] & K_{1}\left(D^*_u(M)/D^*_u(\olO\subseteq M,\Sigma)\right)\ar[d]^{\partial} \\
K_{0}\left(D^*_u(\olO,\Sigma)\right)\ar@{}[r]|{\oplus} & K_{0}\left(D^*_u(\Sigma\subseteq \Sigma_\infty,\Sigma)\right) \ar[r] & K_{0}\left(D^*_u(\olO\subseteq M,\Sigma)\right)
}} 
where the horizontal arrows are the sums of the respective inclusion maps. Note that the first vertical arrow is precisely the relative uniform index map. Inserting the above decomposition for $\chi(\Dirac_M)$ into this diagram lets us conclude that
$$
\partial (\chi(\Dirac_M)) = \relInd \left( \psi_{\olO} \chi(\Dirac_M) \psi_{\olO} \right) + \partial \left( \psi_{\Sigma_\infty} \chi(\Dirac_M) \psi_{\Sigma_\infty} \right) \; ,
$$
Abusing notation we view $\chi(\Dirac_M)$ etc. as elements of the respective K-theory groups directly, whereas they should strictly speaking be identified with the corresponding unitary via Paschke duality first. 

Now, $\partial (\chi(\Dirac_M))$ coincides with the image of $\Ind^{APS}(\Dirac_\Omega)$ under the inclusion of $C^*_u(\olO)$ into $D^*_u(\olO\subseteq M,\Sigma)$ by definition of the uniform APS-index. Moreover, $\psi_{\olO} \chi(\Dirac_M) \psi_{\olO}$ gives rise to a uniform Fredholm module over $\Omega$ that coincides up to a degenerate Fredholm module with the image of $\chi(\Dirac_M)$ under the restriction map $K^{u}_{0}(M) \to K^{u}_{0}(\Omega)$. Thus Proposition \ref{Restriction} yields $[\psi_{\olO} \chi(\Dirac_M) \psi_{\olO}]=[\Dirac_\Omega]$ in $K^{u}_{0}(\Omega)$. By excision this equality also holds in $K^{u}_{0}(\olO,\Sigma)$, so that 
$$
\relInd \left( \psi_{\olO} \chi(\Dirac_M) \psi_{\olO} \right) = \relInd(\Dirac_\Omega) \; .
$$
Analogously, it holds that $[\psi_{\Sigma_\infty} \chi(\Dirac_M) \psi_{\Sigma_\infty}]= [\Dirac_{\Sigma_\infty}]\in K^{u}_{0}(\Sigma_\infty,\Sigma)$. These two operators differ only by an operator in $D^*_u(\Sigma\subseteq \Sigma_\infty,\Sigma)$, so that this equality actually holds in $K_{1}(D^*_u(\Sigma_\infty)/D^*_u(\Sigma\subseteq \Sigma_\infty,\Sigma))$. Thus we conclude from Theorem \ref{RelIndRho} that
$$
\partial \left( \psi_{\Sigma_\infty} \chi(\Dirac_M) \psi_{\Sigma_\infty} \right) = \Ind^{\Sigma}(\Dirac_{\Sigma_\infty}) = -\rho^{u}(\Dirac_\Sigma) \; .
$$
The sign comes from the fact that the attached cylinder $\Sigma_\infty=\RR_{\leq 0}\times \Sigma$ carries the opposite orientation of that considered in Theorem \ref{RelIndRho}. In combination we conclude that
$$
\Ind^{APS}(\Dirac_\Omega) = \relInd(\Dirac_\Omega) - \rho^{u}(\Dirac_\Sigma)
$$
in $K_0(D^*_u(\olO\subseteq M,\Sigma))$. Since this K-theory group is isomorphic to $K_0(D^*_u(\olO,\Sigma))$ by Proposition \ref{IsoStructureLocalized}, this concludes the proof.
\end{proof}
Again, we can appeal to exactness of
$$
K^{u}_*(\olO,\dO)\xlongrightarrow{\relInd} K_*\left(D^*_u(\olO,\dO)\right) \longrightarrow K_*\left(D^*_u(\olO)\right)
$$
to conclude that the right-hand side of the index formula in Theorem \ref{DelocIndexFormula} vanishes upon being pushed forward to $K_*(D^*_u(\olO))$. This yields the uniform version of the delocalized APS-index theorem.
\begin{cor}[Uniform delocalized APS-index theorem]
Let $\olO$ be an even-dimensional spin manifold with boundary and bounded geometry. Assume that the metric on $\olO$ has product structure near the boundary, and that the scalar curvature of $\dO$ is uniformly positive. Then,
$$
\Ind^{APS}(\Dirac_\Omega) = - \rho(\Dirac_{\dO}) \; \in \; K_{0}\left( D^*_u(\olO)\right) \; .
$$
\end{cor} 
The non-uniform version of this index theorem was proved by Piazza and Schick \cite{Piazza2014} in even dimensions, and by Xie and Yu in the general case \cite{Xie2014}.\footnote{We have already remarked in the introduction that the non-uniform delocalized APS-index theorem does not have the sign that ours does. This is a consequence of a difference in choice of orientation. Indeed, in the other treatments of this theorem the orientation on $M$ is chosen such that the cylinder is positively oriented, meaning that the orientation on $\olO$ is reversed. While this results in a sign-free delocalized APS-index theorem, it is the less natural choice from our point of view. Since our primary interest was the relative index $\relInd(\Dirac_\Omega)$, it made sense for us to choose that orientation on $M$ which coincides with the given one on $\olO$. This resulted in the sign.} The author expects this theorem to hold in odd dimensions as well. Indeed, the proof of Theorem \ref{DelocIndexFormula} relied solely on the validity of Theorem \ref{RelIndRho}, so it would suffice to prove that theorem in odd dimensions. Let us sketch the difficulty involved. If $\olO$ is odd-dimensional, the Dirac operator $\Dirac=\Dirac_\Omega$ is ungraded. To get a finite-propagation representative proceed as in Lemma \ref{KHomDilation} and consider the operator $\hat{D}_e=\begin{pmatrix}
 0 & \Dirac_{\max} \\ \Dirac_{\min} & 0
\end{pmatrix}$. Since $\Dirac_{\min}$ is coercive $\hat{D}_e$ has $0$ as an isolated spectral value due to the potential presence of $\ker(\Dirac_{\max})$. Choose a normalizing function that is equal to $\pm 1$ on $\sigma(\hat{D}_e)\setminus \{0\}$ , then 
$$
\chi(\hat{D}_e) = \begin{pmatrix}
 0 & V^* \\ V & 0
\end{pmatrix} \; ,
$$
where $V$ is a generally non-surjective isometry in $D^*_u(\olO)$ that is self-adjoint and invertible modulo $D^*_u(\Sigma\subseteq \olO,\Sigma)$. Then $V$ represents $[\Dirac]\in K^{u}_1(\olO,\Sigma)$. Its relative uniform index is computed as follows. Form the projection $\frac{1}{2}(I+V)$ in $D^*_u(\olO)/D^*_u(\olO,\Sigma)$ and pick a self-adjoint lift $F\in D^*_u(\olO)$. Then $\exp(2\pi i F)$ represents $\Ind^\Sigma(\Dirac)$. However no natural self-adjoint lift presents itself, beyond the standard $F=\frac{1}{2}(I+\frac{V+V^*}{2})$. The author didn't manage to extract any useful information from this, since he couldn't find a workable description of how $V$ and $V^*$ differ concretely. A reasonable guess would be that $\partial[\frac{1}{2}(I+V)]$ is related to the cokernel projection $I-VV^*$, because it provides a measure of how far $\frac{1}{2}(I+V)$ is from an actual projection. Since $I-VV^*$ is coincides with the projection onto $\ker(\Dirac_{\max})$, this would likely allow for an analogous proof to that given in the even case. But the author could not get the computations to work.

It speaks to a somewhat systemic issue that Piazza and Schick could also treat only the even-dimensional case in their work on the delocalized APS-theorem \cite{Piazza2014}. The general case of that theorem was proved by Xie and Yu using localization algebras instead of Fredholm modules \cite{Xie2014}. The localization algebra picture of K-homology has the advantage that the description of the coarse index map is easier from a practical perspective, and that product operations are more transparently incorporated. In this context we also point to Zeidler`s work on product formulas and secondary invariants \cite{Zeidler2016}. Investigating a localization algebra picture in the uniform setting would be very worthwhile, since localization algebra techniques regularly succeed in situations that are practically inaccessible (or at least much harder) in the Fredholm module picture. We discuss this further in the outlook below.
\chapter{Outlook}
We conclude this thesis with an outlook. Since uniform K-homology and uniform index theory has been the subject of a comparatively small amount of research there is much that can be done, both in the direction of application to geometric questions and the development of further theory.
\subsection*{Field testing}
The aim of this thesis was the development of theory concerning the uniform K-homology classes of uniformly elliptic operators on manifolds with boundary. Even as far as applications to scalar curvature and Dirac operators were concerned, they were of a theoretical nature, equating relative indices with secondary invariants and the like. At this point concrete and practical applications are sorely needed. Given that scalar curvature considerations have been a driving motivation throughout this thesis it is natural to pose the following problem.
\begin{problem}
Use uniform coarse index theory to rule out the existence of upsc-metrics on a manifold (or class of manifolds) for which it has not previously been known whether it admits a upsc-metric. For instance, find a spin manifold $M$ of bounded geometry for which 
$$
\mathrm{Ind}(\Dirac_M) = 0 \; \in \; K_*\left(C^*(M)\right) \; ,
$$
but 
$$
\Ind(\Dirac_M) \neq 0 \; \in \; K_*\left(C^*_u(M)\right) \; .
$$
\end{problem}
Note that the situation that $\mathrm{Ind}(\Dirac_M)=0$ but $\Ind(\Dirac_M)\neq 0$ could also reflect that $M$ admits a upsc-metric, though not a upsc-metric of bounded geometry. It would also be interesting to have an example where the (uniform) coarse index detects this. It seems to be more immediately practical to search for examples in an intrinsically non-compact setting. Indeed, if $M$ is the universal cover of a closed manifold with fundamental group $\Gamma$, then we know that the uniform coarse index of the Dirac operator on $M$ is the image of the $\Gamma$-equivariant index under the map $K_*(C^*_r\Gamma)=K_*(C^*(M)^\Gamma)\to K_*(C^*_u(M))$. Thus, if the uniform coarse index is non-zero, the $\Gamma$-equivariant index (or indeed the Rosenberg index) must already have been non-zero. The Gromov-Lawson-Rosenberg posits that a vanishing Rosenberg index is a sufficient condition for the existence of a psc-metric on a closed manifold of dimension at least $5$. While proved for a number of groups, it is known to be false in general \cite{SchickGLR}. Thus, there is a chance that uniform coarse index theory provides insights to the psc-question on closed manifolds that are overlooked by the Rosenberg index. Again, since the Rosenberg index refines the uniform coarse index of the equivariant Dirac operator, such an application must necessarily arise through breaking equivariance, and over a space not known to satisfy the Gromov-Lawson-Rosenberg conjecture.

Let us sketch in some more detail how this might work. Let $M$ be the universal cover of a closed spin manifold with fundamental group $\Gamma$, and let $\tilde{\Dirac}$ be the $\Gamma$-equivariant Dirac operator on $M$. If $E\to M$ is a flat vector bundle, then the Lichnerowicz formula shows that $\Ind(\tilde{\Dirac}_E)$ still obstructs uniformly positive scalar curvature on $M$. If $E$ is not $\Gamma$-invariant, then $\tilde{\Dirac}_E$ is not $\Gamma$-equivariant, meaning that $\Ind(\tilde{\Dirac}_E)$ is an obstruction to upsc that is intrinsic to the uniform setting, and not accessible by equivariant methods. The idea of twisting by flat vector bundles to investigate geometric questions was already suggested by Gromov \cite{Gromov1991}. A potential source of interesting flat vector bundles is group-cohomological data. In \cite{Mathai2006} Mathai constructed an index for Dirac operators twisted by flat line bundles coming from projective representations constructed via group cocycles. These twisted indices live in the K-theory of a twisted group $C^*$-algebra, and they obstruct the existence of a psc-metric on the base manifold. Projectively equivariant operators are still uniform, so that the K-theory uniform Roe algebra could provide a universal receptacle for the twisted indices arising from different twists. In this way uniform coarse index theory could contribute fruitfully to the psc-question for closed manifolds.
\begin{problem}
Find a closed spin manifold of dimension at least $5$ (denote its fundamental group by $\Gamma$ as usual), and a flat vector bundle $E\to M$ of bounded geometry, such that\footnote{We gloss over the detail here that one really should consider the Real version of these indices.} 
$$
\mathrm{Ind}^\Gamma(\tilde{\Dirac}) = 0 \; \in \; K_*\left(C^*_r\Gamma\right) \; ,
$$
but 
$$
\Ind(\tilde{\Dirac}_E) \neq 0 \; \in \; K_*\big( C^*_u(\tilde{M}) \big) \; .
$$
One should look for $M$ among the known counterexamples to the Gromov-Lawson-Rosenberg conjecture (e.g. \cite{SchickGLR}), and candidates for $E$ are flat line bundles associated to group cocycles.
\end{problem}
There is also the secondary question of whether two metrics of (uniformly) positive scalar curvature on a manifold agree up to a given equivalence (concordance, psc-bordism, etc.). The higher $\rho$-invariant provides a potent tool to study this problem, and more generally to describe the moduli space of metrics of positive scalar curvature. We refer to \cite{Xie2021} for a discussion that is particularly close to the content of this thesis. Here too could one seek for situations in which the uniform machinery detects structure the coarse machinery misses. Again, it would be more immediately natural to search for such instances in an intrinsically non-compact context, but potential applications to closed manifolds should not be counted out. These could again arise from twisting with interesting non-equivariant structures. We point to a discussion of twisted secondary invariants and APS-indices by Azzali and Wahl \cite{Azzali2019}.

There are of course also other interesting applications of index theory besides scalar curvature. For one, higher indices of the signature operator are of immediate relevance to higher signatures and thus the Novikov conjecture, which is closely connected to the Baum-Connes conjecture and the assembly map. See for example the survey \cite{Rosenberg2016}. Moreover, the vanishing of the index of the Euler characteristic operator -i.e the Euler characteristic- carries geometric significance in the form of the potential existence of an affine structure, and a vast amount of index theory has been applied to this problem. We refer to the diagrammatic overview in \cite{Löh2023} and the multitude of references therein. These too are situations to which uniform coarse index theory might be fruitfully applied. Irrespective of the particular geometric problem on which one decides to field-test uniform coarse index theory it would be highly beneficial to gain an understanding of practical situations in which uniform coarse index theory provides insight that cannot be gleamed with either coarse or equivariant methods. Phrased succinctly: We should understand the niche occupied by uniform coarse index theory.
\subsection*{Understanding the relative index class}
In this thesis we have introduced a relative index map to a relative version of the structure algebra that contains secondary information on the boundary. In the case of the Dirac operator on a manifold (of even dimension) with uniformly positive scalar curvature on its boundary we understand this index satisfactorily through the index formula of Theorem \ref{DelocIndexFormula}. Here the secondary boundary information is the (uniform) $\rho$-invariant of the boundary. Even in this situation there remains much to be done; first of all the extension of Theorem \ref{DelocIndexFormula} to odd-dimensional, but also a secondary partitioned manifold theorem for uniform $\rho$-invariants (see \cite{Siegel2012}, \cite{Piazza2014}, \cite{Zeidler2016}). These problems may benefit from an adaptation of localization algebra techniques, see our discussion below. But the computation of these relative indices should be interesting beyond scalar curvature. What secondary information could the relative index of the Dirac operator contain in absence of uniformly positive scalar curvature? What about other operators like the signature operator? 

Another question is the relation of our relative index to other indices already investigated in the literature, primarily the relative index of Chang, Weinberger and Yu. The author is lacking any concrete idea of how to relate the two, but a vague conjecture could go like this.
\begin{conjecture}[Relation to the Chang-Weinberger-Yu index]
Let $\olO$ be a compact manifold with boundary and fundamental group $\Gamma$, and $\dO$ its boundary with fundamental group $\Gamma'$. Let $\olO^\sim$ be the universal cover of $\olO$. There exists an equivariant refinement $D^*(\olO^\sim,\partial(\Omega^\sim))^{(\Gamma,\Gamma')}$ of the relative structure algebra, together with a map
$$
K_*\left( D^*(\olO^\sim,\partial(\Omega^\sim))^{(\Gamma,\Gamma')} \right) \longrightarrow K_*\left( C^*_{\max}(\Gamma,\Gamma') \right)
$$
that maps the relative index class in $K_*( D^*(\olO^\sim,\partial(\Omega^\sim))^{(\Gamma,\Gamma')} )$ to the relative index of Chang, Weinberger and Yu.
\end{conjecture}
A potential map $K_*( D^*(\olO^\sim,\partial(\Omega^\sim))^{(\Gamma,\Gamma')}) \longrightarrow K_*( C^*_{\max}(\Gamma,\Gamma'))$ would necessarily have to forget any secondary information on the boundary, because the Chang-Weinberger-Yu index vanishes on a upsc-manifold, whereas our relative index produces the (uniform) $\rho$-invariant, which never vanishes. The Chang-Weinberger,Yu index is therefore closer to the index in $K_*(C^*(\olO,\dO))$, which at least vanishes on upsc-manifolds. Thus we might expect the map $K_*( D^*(\olO^\sim,\partial(\Omega^\sim))^{(\Gamma,\Gamma')}) \longrightarrow K_*( C^*_{\max}(\Gamma,\Gamma'))$ to factor through an equivariant version of the forgetful map
$$
K_*\left( D^*(\olO^\sim,\partial(\Omega^\sim)\right)^{(\Gamma,\Gamma')}) \longrightarrow K_*\left( C^*(\olO^\sim,\partial(\Omega^\sim))^{(\Gamma,\Gamma')}\right) \; ,
$$
where $C^*(\olO^\sim,\partial(\Omega^\sim))^{(\Gamma,\Gamma')}$ is a suitable equivariant version of $C^*(\olO^\sim,\partial(\Omega^\sim))$. Given the background that the K-theory of the equivariant Roe algebra on a universal cover reproduces the K-thery of the group $C^*$-algebra, it would make conceptual sense if the K-theory of the relative group $C^*$-algebra could be recovered by the K-theory of an equivariant relative Roe algebra of the universal cover of a compact manifold with boundary. However, we have remarked repeatedly on the forgetful nature of $K_*(C^*(\olO,\dO))$, and since the boundary of the universal cover of a compact manifold with boundary is $c$-dense, the K-theory of the relative Roe algebra would vanish in this case. Of course, this need no longer be true if equivariance is added to the picture, but it is completely unclear to the author if this could realistically recover the relative group $C^*$-algebra.

Moreover, using the maximal relative group $C^*$-algebra is likely to be wrong, and one should rather use the reduced relative $C^*$-algebra under the additional assumption that $\Gamma'$ is a subgroup of $\Gamma$. In this situation too Chang, Weinberger and Yu provide a relative index, and use of the reduced version would be more in line with the fact that the equivariant Roe algebra on a universal cover is Morita-equivalent to the reduced group $C^*$-algebra of the fundamental group.
\subsection*{Mapping bounded geometry to analysis with $\ell^\infty$-coefficients}
In the introduction to this thesis we discussed the framework of mapping positive scalar curvature to analysis. It consists of Stolz' positive-scalar-curvature sequence on  the geometric side, the analytic surgery exact sequence on the analytic side, and maps between the two that are given by K-homology classes, higher indices and higher $\rho$-invariants. This framework makes explicit the various connections between scalar curvature, higher indices and secondary invariants, and assembly.

In this thesis we have introduced uniform analogues to the higher $\rho$-invariants and APS-indices that are used to map geometry to analysis. Moreover various bordism invariance statements were proved in the uniform setting, whose non-uniform counterparts give well-definedness of the maps from geometry to analysis. Indeed, only the bordism invariance of the higher APS-index has not been proved in the uniform setting yet. Thus we are rather close to mapping geometry to uniform analysis, where we keep Stolz' positive-scalar-curvature sequence on the the geometric side, but take the uniform counterpart of the analytic surgery exact sequence as the analytic side, and realize the mappings from the former to the latter via uniform K-homology classes, uniform indices and uniform higher secondary invariants. Under the simplifying assumption that $E\Gamma$ can be realized as a locally compact second-countable metric space the uniform analogue of the analytic surgery exact sequence takes the form
\begin{equation}\label{AnaSESuniform}
\cdots \longrightarrow K_{*+1}\left(D^*_u(E\Gamma)\right) \longrightarrow K_*^{u}(E\Gamma) \xlongrightarrow{\Ind} K_*\left( C^*_u(E\Gamma) \right) \longrightarrow \; .
\end{equation}
The K-theory group $K_*\left( C^*_u(E\Gamma) \right)$ is isomorphic to $K_*\left( C^*_u(\Gamma) \right)=K_*\left(\ell^\infty(\Gamma)\rtimes_r\Gamma\right)$, and the uniform coarse index map $K_*^{u}\left( E\Gamma\right) \to K_*\left(\ell^\infty(\Gamma)\rtimes_r\Gamma\right)$ is closely connected to the Baum-Connes assembly map with $\ell^\infty(\Gamma)$-coefficients. Therefore we feel justified in calling the exact sequence \eqref{AnaSESuniform} the \emph{analytic surgery exact sequence with $\ell^\infty$-coefficients}.

Admittedly, mapping the positive-scalar-curvature sequence to the analytic surgery exact sequence with $\ell^\infty$-coefficients would in itself not be all that interesting. Since all analytic objects arising from the former sequence are built from $\Gamma$-equivariant Dirac operators, mapping to analysis with $\ell^\infty$-coefficients instead of analysis simply means forgetting equivariance. Thus the maps from geometry to analysis with $\ell^\infty$-coefficients factor through those mapping geometry to analysis. 

The issue is that in the above set-up we no longer insist on equivariance on the analytic side, but require it on the geometric side. The path forward then seems clear: No longer insist on equivariance on the geometric side either. Instead of compact manifolds together with maps to $B\Gamma$, which correspond to $\Gamma$-coverings of compact manifolds, we can consider general manifolds of bounded geometry, together with (sufficiently nice) maps to $E\Gamma$. This yields bounded-geometry analogues to the terms of the positive-scalar-curvature exact sequence. We conjecture that this results in a framework of mapping bounded geometry to analysis with $\ell^\infty$-coefficients. 
\begin{conjecture}[Mapping bounded geometry to analysis with $\ell^\infty$-coefficients]
Let $\Gamma$ be a discrete group, and assume for simplicity that $E\Gamma$ can be realized as a locally compact second-countable metric space. Define the following objects:
\begin{itemize}
\item Let $\mathrm{Pos}^{bg,spin}_n(E\Gamma)$ consist of $n$-dimensional spin manifolds equipped with a upsc-metric of bounded geometry and a uniformly cobounded Lipschitz map to $E\Gamma$, modulo finite-width spin-bordisms equipped with a upsc-metric of bounded geometry and a uniformly cobounded Lipschitz map to $E\Gamma$.
\item Let $\Omega^\mathrm{bg,spin}_n(E\Gamma)$ consist of $n$-dimensional spin manifolds of bounded geometry together with a uniformly cobounded Lipschitz map to $E\Gamma$, modulo finite-width spin-bordisms of bounded geometry and a uniformly cobounded Lipschitz map to $E\Gamma$.
\item Let $R^{\mathrm{bg,spin}}_n(E\Gamma)$ consist of spin-manifolds with possibly non-empty boundary, equipped with a bounded-geometry metric having upsc on the boundary (if existent) and a uniformly cobounded Lipschitz map to $E\Gamma$, modulo the bounded-geometry version of the relation for $R^\mathrm{spin}_n$.
\end{itemize}
Then, there is a commutative diagram \\
\centerline{\xymatrix{
\cdots \ar[r] & \mathrm{Pos}^\mathrm{bg,spin}_{*}(E\Gamma) \ar[d] \ar[r] & \Omega_*^\mathrm{bg,spin}(E\Gamma) \ar[d] \ar[r] & R^\mathrm{bg,spin}_*(E\Gamma) \ar[r] \ar[d] & \cdots \\
\cdots \ar[r] & K_{*+1}\left(D^*_u(E\Gamma) \right) \ar[r] & K_*^{u}(E\Gamma) \ar[r] & K_*\left(\ell^\infty\Gamma \rtimes_r \Gamma \right) \ar[r] & \cdots
}}
with exact rows. Here the uniform index map $K_*^{u}(E\Gamma) \to K_*(C^*_u(E\Gamma)))= K_*\left(\ell^\infty(\Gamma) \rtimes_r \Gamma \right)$ factors as the composition
$$
K^{u}_*(E\Gamma) \longrightarrow K^{u}_*(\underline{E}\Gamma)=K^{\Gamma}_*\left(\underline{E}\Gamma;\ell^\infty(\Gamma)\right) \xlongrightarrow{\mu_{\ell^\infty(\Gamma)}} K_*\left(\ell^\infty(\Gamma) \rtimes_r \Gamma \right) \; ,
$$
where the first map is induced by the canonical map $E\Gamma\to \underline{E}\Gamma$, and the second map is the Baum-Connes assembly map with $\ell^\infty(\Gamma)$-coefficients. The vertical arrows are given by the push-forwards of the uniform $\rho$-invariant, the uniform K-homology class of the Dirac operator, and the uniform APS-index along the maps uniformly cobounded Lipschitz maps to $E\Gamma$, respectively. 
\end{conjecture}
Much concerning the vertical arrows in the above diagram was already done over the coarse of this thesis. Well-definedness of the uniform $\rho$-invariant mapping $\mathrm{Pos}^\mathrm{bg,spin}_{*}(E\Gamma) \to K_{*+1}(D^*_u(E\Gamma))$ is a consequence of Corollary \ref{BordismInvRho} (compare \cite[Corollary 1.26]{Piazza2014}). Well-definedness of the map $\Omega_*^\mathrm{bg,spin}(E\Gamma) \to K_*^{u}(E\Gamma)$ is a consequence of boundary-if-Dirac-is-Dirac (Theorem \ref{BdryOfDirac}). The commutativity of the above diagram should largely be by definition, or a consequence of the uniform delocalized APS-theorem. The well-definedness of the uniform APS-index map $R^\mathrm{bg,spin}_*(E\Gamma) \to K_*\left(\ell^\infty(\Gamma)\rtimes_r \Gamma\right)$ was not investigated in this thesis, though it seems likely that it is within reach of the tools developed here.

Part of this conjecture is also the exactness of the bounded-geometry positive-scalar-curvature sequence. This is a purely geometric matter, and therefore somewhat disjoint from the analytic focus of this thesis. Proving exactness requires understanding of how the existence of upsc-metrics on a manifold of bounded geometry is preserved under surgeries and bordisms. The author cannot claim to be at all familiar with the state of the literature on these matters.
\subsection*{Chern character, cyclic homology, numerical invariants}
In the introduction to this thesis we have discussed how the usual $\rho$-invariant can be recovered from the higher $\rho$-invariants by means of a relative trace on the K-theory of the structure algebra. This indicated a general principle by which K-theoretical invariants can be turned into numerical ones by means of traces. Indeed, traces are a particular instance of cyclic cocycles, and there is a non-commutative Chern character from the K-theory of any $C^*$-algebra to its cyclic homology, which can be paired with cyclic cocycles. Examples of cyclic cocycles can come from algebraic data associated with, say, the fundamental group of the space in question, but for example also from averaging procedures over a compact exhaustion of an amenable space. Index-theoretical instances of the former are discussed at length in \cite{Piazza2025}, whereas (among other things) an interpretation of Roe's index theorem on amenable manifolds in light of the latter can be found in \cite{Engel2015}, a little more on this below.

An investigation into these matters in the context of manifolds with boundary could enable a deeper understanding of the objects discussed in this thesis. For example we could ask to understand the Chern character of the relative uniform index, an element of the cyclic homology of $D^*_u(\olO,\partial\Omega)$, or that of the uniform $\rho$-invariants. Are there natural cyclic cocycles with which they can be paired to obtain interesting numerical invariants? Investigations in this direction could promise a more concrete understanding of these objects, which at present is missing. 

A related question is that of a geometric Chern character from uniform K-homology to a suitable version of singular homology (say uniformly locally finite homology). Engel has constructed such a Chern character on manifolds with bounded geometry \cite{Engel2015}. It seems likely that a similar construction can also be carried out for manifolds with boundary and bounded geometry. Moreover, it seems natural to ask for a relative Chern character from relative uniform K-homology to relative uniformly locally finite homology, as well as a Chern character on more general spaces than  manifolds. Such Chern characters could be used to turn K-homological invariants into homological ones, which are generally more amenable to direct computations via combinatorial methods. 

On a related note Engel also constructed index classes in uniform de Rham-cohomology for a uniformly elliptic operator on an oriented manifold of bounded geometry \cite{Engel2015}. For the classical operators these are extensions of the usual characteristic classes like Todd- and $\hat{A}$-class to the bounded-geometry setting. Again, it seems natural to ask for a relative counterpart in the presence of a boundary. If one could associate an index class in relative uniform de Rham-cohomology to, say, the Dirac operator on a spin manifold with boundary and bounded geometry, it would be reasonable to expect it to be some relative bounded version of the $\hat{A}$-class. We could then ask how this relative index class relates to the relative uniform index in $K_*\left(D^*_u(\olO,\partial\Omega)\right)$, providing another potential avenue to understand the relative uniform index.
\subsection*{New perspectives on uniform K-homology}
The definition of uniform K-homology is an adaptation of analytic K-homology. However, analytic K-homology is not the only model of K-homology. There is the homotopy-theoretic definition via the K-theory spectrum, there is geometric K-homology whose description is closer to bordism homology, there are Yu's localization algebras, and there are more. Different models of K-homology are better suited for different applications. We have discussed an instance of this in the context of index maps and products, where the technically more straight-forward description of the objects provided by the localization algebra picture enabled treatment of situations that were practically inaccessible in the Fredholm module picture. As different models provide different avenues to a given problem it is also desirable to have different models for uniform K-homology.

Let us start with localization algebras, which were introduced by Yu in the context of the coarse Baum-Connes conjecture \cite{Yu1997}. Let $X$ be a proper metric space, and fix some ample representation of $X$ on the Hilbert space $H$. The \emph{localization algebra} $C^*_L(X)$ of $X$ is defined as the $C^*$-algebra generated by all functions $[0,\infty)\to \fBH, \, t \mapsto T_t$, such that $T_t$ is locally compact, has finite propagation, and $\mathrm{prop}(T_t)\to 0$ for $t\to \infty$. It is true but not obvious that there is a natural isomorphism $K_*(C^*_L(X))\cong K_*(X)$. Evaluation at $0$ provides a map $K_*(C^*_L(X))\to K_*(C^*(X))$, this is the coarse index map in the localization algebra picture. An immediate uniformization of the localization algebra presents itself: Instead of requiring the map $t\mapsto T_t$ to have locally compact operators as images, use uniformly locally compact operators instead. Denoting the resulting algebra by $C^*_{L,u}(X)$ evaluation at zero yields a map $K_*(C^*_{L,u}(X))\to K_*(C^*_u(X))$, which presents an analogue to the uniform coarse index maps. The isomorphism $K_*(C^*_L(X))\cong K_*(X)$ is not obvious, and it seems reasonable though not immediately clear that it generalizes to the respective uniform versions. We refer to our discussion in Chapter 7 for potential applications to index theory and secondary invariants that this isomorphism would have. 

We come to geometric K-homology. A cycle for the geometric K-homology of $X$ is a triple $(M,E,f)$ consisting of a closed spin$^c$ manifold $M$, a vector bundle $E\to M$, and a continuous map $f: \, M\to X$. The sum of geometric cycles is given by disjoint union, and the relations imposed to obtain geometric K-homology are bordism and so-called vector bundle modification. The map from geometric to analytic K-homology is easily described on the level of cycles: The spin$^c$-Dirac operator on $M$ twisted by $E$ defines a class $[\Dirac_E]\in K_*(M)$, which then gets pushed forward along the map $f$ to the class $f_*[\Dirac_E]\in K_*(M)$. Well-definedness of this assignment under bordism follows from the boundary-of-Dirac-is-Dirac formula, well-definedness under vector bundle modification is a consequence of compatibility with the Kasparov product. That this map gives an isomorphism between geometric and analytic K-homology is far less clear, and the proof actually goes through the topological model for K-homology \cite{Baum2007}.

The question of a uniform version of geometric K-homology has already been raised by Engel \cite{Engel2014}. Writing down a candidate theory is not too difficult. If $X$ is a metric space, let us define a \emph{uniform geometric cycle} for $X$ to be a triple $(M,E,f)$ consisting of a spin$^c$ manifold $M$ of bounded geometry, a vector bundle $E\to M$ of bounded geometry, and a suitably uniform map $f: \, M\to X$, say uniformly cobounded and Lipschitz. Impose the relations of finite-width bordism and vector bundle modification to obtain uniform geometric K-homology. On the level of cycles we can get an assignment to uniform K-homology just as in the non-uniform case: The twisted Dirac operator over $M$ defines a class $[\Dirac_E]\in K^{u}_*(M)$, which can be pushed forward to a class $f_*[\Dirac_E] \in K^{u}_*(X)$. That this assignment respects the bordism relation is a consequence of the uniform boundary-of-Dirac-is-Dirac formula, i.e. our Theorem \ref{BdryOfDirac}. That it respects vector bundle modifications should again follow from compatibility with the Kasparov product (our Proposition \ref{ProductComplete}, compare \cite[Proposition 3.6]{Baum2007}). 
\begin{conjecture}
The natural map 
$$
K^{u,geo}_*(X) \longrightarrow K^{u}_*(X) \; , \; [M,E,f] \longmapsto f_*[\Dirac_E]
$$
is an isomorphism for reasonable spaces $X$.
\end{conjecture}
Notice that the geometric model of uniform K-homology would enable an extension of the Chern character defined in \cite{Engel2015} for manifolds of bounded geometry to more general spaces via $\mathrm{ch}_*([M,E,f]) := f_* \mathrm{ch}_*([\Dirac_E]) \in H^{ulf}_*(X)$. Of course, it would need to be checked that this definition is well-defined and coincides with the one already given for manifolds of bounded geometry. Nonetheless this seems like a promising approach to the extension of Chern characters to more general spaces, the usefulness of which was indicated in the previous section.

As noted above the proof of this isomorphism in the non-uniform setting uses the topological model of K-homology as a mediator between the geometric and the analytic model. It is much less clear how a uniform version of topological K-homology would look like. In \cite{Engel2015} Engel proposed the development of a more general uniform homotopy theory, into which uniform K-homology could enter as the uniform homology theory corresponding to a uniform K-theory spectrum. See \cite{Bunke2020} for corresponding work in the coarse category. This ties into considerations regarding the category on which uniform K-homology should be considered. Originally, this was the category of metric spaces and uniformly cobounded Lipschitz maps \cite{Spakula2009}, \cite{Engel2019}. However, in this thesis we have noted that functoriality under a more general class of maps (which we called (approxiamtely) filtered maps) holds, and that the additional flexibility provided by this more general class of morphisms can be useful. However, the definition of these morphisms is very ad-hoc, and likely not very natural outside the specific context of uniform K-homology. Thus it would be desirable to find a geometrically relevant class of morphisms that is more general than uniformly cobounded Lipschitz maps but still falls inside the framework of (approximately) filtered morphisms.
\subsection*{Local regular boundary conditions}
In our applications of the uniform K-homology classes of uniformly elliptic operators we have focused on relative classes. This is because Dirac operators, which facilitate the connection to scalar curvature, do not admit local, self-adjoint and regular boundary conditions. There are however interesting operators that do admits such boundary conditions. A prime example is the Euler characteristic operator with its absolute and relative boundary conditions. In the compact case its index is the absolute and relative Euler characteristic respectively. On manifolds with boundary and bounded geometry Schick has identified the solution space of absolute/relative boundary condition with the absolute/relative $L^2$-de Rham cohomology of the underlying manifold \cite{Schick1998}. Therein curvature obstructions in terms of $L^2$-Betti numbers where given as well. It seems likely that cohomological information can be gleamed from the uniform K-homology class of the Euler characteristic operator with absolute/relative boundary conditions, perhaps via its image under the Chern character.

The Atiyah-Bott theorem calculates the index of a regular differential boundary conditions in terms of topological data built from the symbols of both the operator on the total space and the operator defining the boundary condition. In particular the index depends only on these symbols. We have shown that that the class of a uniformly elliptic operator in relative uniform K-homology depends only on its principal symbol. It seems feasible that analogous techniques can be used to show that the absolute class of a regular differential boundary condition depends only on the symbols of the elliptic operator and the boundary operator. The analysis will no doubt be made more difficult by the general dependence of the class on the choice of boundary condition.

Roe proved an index theorem for uniformly elliptic Dirac-type operators on amenable manifolds of bounded geometry \cite{Roe1988}. The numerical indices were defined via an averaging procedure of the usual analytic and topological indices over a F\o lner sequence, and Roe's theorem asserts that the averages of the two indices coincide. The proof is an adaption of the usual heat kernel proof. Engel put these index maps on a uniform K-homological basis, and strengthened Roe's index theorem to treat general uniformly elliptic pseudo-differential operators \cite{Engel2014}. The foundations of these index maps on uniform K-homology were established for general amenable spaces of jointly bounded geometry, and are therefore also available on manifolds with boundary and bounded geometry. Moreover, there is a heat kernel proof for the Atiyah-Bott theorem as well (see for example \cite{Gilkey1996}). Thus, mimicking the analysis of Roe and Engel could produce an analogue of the Atiyah-Bott theorem on amenable manifolds with boundary and bounded geometry obtained via index maps on uniform K-homology and averaging over a F\o lner sequence.
\appendix
\chapter{Gradings and multigradings} \label{MultigradingAppendix}
Multigradings are used to define K-homology groups beyond degrees zero and one, and we have cause occasionally to use them in the main part of this thesis. Therefore we collect the basic definitions and facts in this appendix. All conventions and definitions are taken from \cite{HigsonRoe2000}.
\subsection*{Multigradings on vector spaces}
Let $V$ be an inner-product space. A \emph{grading} on $H$ is a unitary operator $\kappa$ on $H$ such that $\kappa^2=1$. The graded space $V$ can be decomposed into an orthogonal direct sum $V=V_+\oplus V_-$ of the $(\pm 1)$-eigenspaces $V_\pm$ of $\kappa$. The orthogonal projection onto $V_\pm$ is given by $P_\pm=\frac{1}{2}(I\pm \kappa)$. Elements of $V$ that lie in either $V_+$ or $V_-$ are called \emph{homogeneous}. The \emph{degree} $\deg(v)$ of a homogeneous element $v$ is defined to be $0$ if $v\in V_+$ and $1$ if $v\in V_-$.

If $\kappa$ is a grading on $V$, then so is $-\kappa$. The space $V$, graded by $-\kappa$, will be denoted by $V^\mathrm{op}$. If $V$ is graded by $\kappa$, and $V'$ by $\kappa'$, then $V\oplus V'$ is graded by $\kappa\oplus \kappa'$.  
\begin{definition}
Let $p\in\ZZ_{\geq 1}$. A \emph{$p$-multigrading} on the inner product space $V$ consists of a grading on $V$, together with $p$ odd unitary operators $\varepsilon_1,\cdots,\varepsilon_p$ on $V$ such that $\varepsilon_i\varepsilon_j+ \varepsilon_j\varepsilon_i=-2\delta_{ij}$. 

We will also call a grading on $H$ a \emph{$0$-multigrading}. If $H$ is ungraded, we will also say that $H$ is \emph{$(-1)$-multigraded}.
\end{definition}
We come to tensor products of (multi)graded spaces. Let $V$ be $p$-multigraded, and $V'$ $p'$-multigraded. Then $V\otimes V'$ is graded by $\kappa\otimes \kappa'$. The space of even elements is given by $(V_0\otimes V_0')\oplus (V_1\otimes V_1')$, that of odd elements by $(V_0\otimes V_1')\oplus (V_1\otimes V_0')$ Moreover, $V\otimes V'$ carries a $(p+p')$-multigrading $\hat{\varepsilon}_1,\cdots,\hat{\varepsilon}_{p+p'}$. For homogeneous elements $v,v'$ they are defined as
$$
\hat{\varepsilon}_j (v\otimes v') = (-1)^{\deg(v')} (\varepsilon_j v)\otimes v' \; , \; \hat{\varepsilon}_{p+j'} (v\otimes v') = v\otimes (\varepsilon'_{j'}v') \; .
$$
We will write $V\hat{\otimes}V'$ for the tensor product equipped with this multigrading structure. If $V$ and $V'$ are Hilbert spaces, then the multigrading structure extends to the Hilbert space tensor product of $V$ and $V'$, and we will also denote this by $V\hat{\otimes}V'$.

Suppose $V$ is graded by $\kappa$. An operator $T$ on $V$ is called \emph{even} if $T\kappa = \kappa T$, and odd if $T\kappa = -\kappa T$. Indeed, conjugation by $\kappa$ defines a grading on the space $\Hom(V)$ of linear operators on $V$, so that even operators are those of degree $0$ and odd operators those of degree $1$.

If $V'$ is another graded space, we can form the graded tensor product $\Hom(V)\hat{\otimes}\Hom(V')$, and endow it with the product defined by 
$$
(T\hat{\otimes} T')(S\hat{\otimes}S') = (-1)^{\deg(T')\deg(S)} (TS)\hat{\otimes}(T'S')
$$
on homogeneous elements. Then, $\Hom(V)\hat{\otimes}\Hom(V')$ becomes a graded algebra. The symbol $\hat{\otimes}$ will be used whenever we consider graded tensor products, and it should be understood that in this context composition of tensor product operators is given according to the \emph{graded} product given above.

If $V$ and $V'$ are Hilbert spaces, the above also applies to the algebras $\fB(V)$, $\fB(V')$ and their algebraic tensor product.

Lastly, if $V$ is $p$-multigraded, then $T\in\Hom(V)$ is said to be \emph{$p$-multigraded} if $T\varepsilon_i = \varepsilon_i T$ for $i=1,\cdots, p$. Any operator is defined to be $0$- and $(-1)$-multigraded.
\subsection*{Multigradings on vector bundles}
Gradings and multigradings can equally well be defined on vector bundles. Let $X$ be a space and $E\to X$ a vector bundle, equipped with a bundle metric. A \emph{grading} on $E$ is a bundle endomorphism $\kappa \in \Hom(E)$ such that the restriction to each fiber $\kappa_x \in \Hom(E_x)$, $x\in X$, is a grading on the inner-product space $E_x$. Analogously, a \emph{$p$-multigrading} on $E$ is a collection of bundle morphisms $\varepsilon_1,\cdots,\varepsilon_p \in \Hom(E)$ that restrict to a $p$-multigrading of each fiber. Direct sums and tensor products of (multi)graded vector bundles can be formed just as for inner-product spaces. Additionally, if $E\to X$ and $E'\to X'$ are $p$- and $p'$-multigraded vector bundles, respectively, then one may form thee $(p+p')$-multigraded exterior tensor product $E\hat{\boxtimes} E' \to X\times X'$, whose fiber at a point $(x,x')$ is the multigraded tensor product $E_x \hat{\otimes} E'_{x'}$.

Suppose that $M$ is a Riemannian manifold, and $E\to M$ a smooth vector bundle equipped with a smooth bundle metric. Then the space $\Gamma_c(M;E)$ of smooth compactly supported sections becomes an inner product space. A (multi)grading on $E$ induces a (multi)grading on $\Gamma_c(M;E)$. A smooth bundle endomorphism of $E$ is even/odd/multigraded if and only if the corresponding operator on $\Gamma_c(M;E)$ is even/odd/multigraded.

Let $D: \, \Gamma_c(M;E)\to \Gamma_c(M;E)$ be a differential operator. We call $D$ even/odd/multigraded if it is even/odd/multigraded as an operator on the inner product space $\Gamma_c(M;E)$. If $D$ is even/odd/multigraded, then its principal symbol is even/odd/multigraded as a bundle morphism on $E$.

Suppose that $D: \, \Gamma_c(M;E)\to \Gamma_c(M;E)$ and $D': \, \Gamma_c(M';E')\to \Gamma_c(M';E')$ are $p$- and $p'$-multigraded differential operators, respectively. Then,
$$
D\times D' := D \hat{\otimes} I + I \hat{\otimes} D'
$$
is a $(p+p')$-multigraded differential operator acting on sections of $E\hat{\boxtimes} E'$. If $D$ and $D'$ are odd, then so is $D\times D'$.
\begin{example}
Consider the $1$-multigraded vector bundle $\slashed{\mathfrak{S}}(\RR):=\CC^2 \to \RR$, equipped with the grading coming from the decomposition $\CC^2=\CC\oplus\CC$ and the multigrading operator $\varepsilon_1= \begin{pmatrix}
0 & -1 \\ 1 & 0
\end{pmatrix}$. Then, 
$$
\slashed{\mathfrak{D}}_\RR = \begin{pmatrix}
0 & -\frac{d}{dt} \\ \frac{d}{dt} & 0
\end{pmatrix}
$$
is an odd and $1$-multigraded differential operator. Let $D: \, \Gamma_c(M;E)\to \Gamma_c(M;E)$ be a $(-1)$-multigraded (that is to say ungraded) differential operator.  Then, $\slashed{\mathfrak{S}}(\RR)\hat{\boxtimes} E \to \RR\times M$ is a $0$-multigraded (i.e. graded) vector bundle, whose fiber at a point $(t,x)$ is $E_x\oplus E_x$. With respect to this decomposition $\slashed{\mathfrak{D}}_\RR\times D$ takes the form
$$
\slashed{\mathfrak{D}}_\RR\times D = \begin{pmatrix}
0 & -\partial_t + D \\ \partial_t + D & 0
\end{pmatrix}
$$
\end{example}
\subsection*{Spinor bundles and Dirac operators}
The essential source of multigraded differential operators is via Clifford-linear Dirac operators on spin-manifolds. We use this section to recall basic facts about spinor bundles and Dirac operators. The standard reference for the matters is \cite{LawsonMichelson1989}.

Let $M$ be an oriented $n$-dimensional Riemannian manifold, with or without boundary. Assume that $M$ is spin, and that a spin structure in form of a lift of the oriented orthonormal frame bundle to a $\Spin(n)$-bundle $P_{\Spin}(M)\to M$ has been chosen. Any representation of $\Spin(n)$ gives rise to an associated vector bundle. Two instances of this are particularly important for our purposes. 

First, consider the complexified Clifford algebra $\CC\ell(n)$ associated to the Euclidean inner product on $\RR^n$. If $n=2m$ is even, then there is a canonical isomorphism from $\CC\ell(n)$ to $\Hom(\CC^{2^m})$, if $n=2m+1$ is odd, a canonical isomorphism to $\Hom(\CC^{2^m})\oplus \Hom(\CC^{2^m})$. In either case one obtains a representation of $\Spin(n)\subseteq \CC\ell(n)$ on the vector space $\Delta_n:=\CC^{2^m}$, in the even case via the isomorphism $\CC\ell(n)\cong \Hom(\Delta_n)$ directly, in the odd case via $\CC\ell(n)\cong \Hom(\Delta_n) \oplus \Hom(\Delta_n)$ followed by the projection onto the first component. This representation is called the \emph{(complex) spinor representation}. It carries a canonical inner product, with respect to which the action of $\CC\ell(n)$ is graded-unitary. The spinor representation gives rise to the \emph{(complex) spinor bundle}
$$
\slashed{S}(M) := P_{\Spin}(M) \times_{\Spin(n)} \Delta_n \; .
$$
In the even case $\slashed{S}(M)$ is graded by the \emph{complex volume element} $\kappa$. Fiberwise it is given by the action of $\kappa=i^m e_1\cdots e_{2m}$, where $(e_1,\cdots,e_{2m})$ is any oriented orthonormal frame of the tangent space of $M$. This is independent of the choice of oriented orthonormal frame, and assembles into a well-defined self-adjoint involution on $\slashed{S}(M)$. Its $(\pm 1)$-eigenbundles are denoted by $\slashed{S}^\pm(M)$. In the odd case $\slashed{S}(M)$ is ungraded.

Second, $\Spin(n)$ acts on $\CC\ell(n)$ by left-multiplication. This action is unitary with respect to the canonical inner product on $\CC\ell(n)$. It gives rise to the vector bundle
$$
\slashed{\mathfrak{S}}(M) := P_{\Spin}(M) \times_{\Spin(n)} \CC\ell(n) \; .
$$ 
Recall that $\CC\ell(n)$ is graded by degree. The group $\Spin(n)$ is contained in the even-degree part of $\CC\ell(n)$, and thus preserves the grading. In this way $\slashed{\mathfrak{S}}(M)$ becomes graded. Moreover, $\CC\ell(n)$ acts on it self by right-multiplication. This right-action is graded, and commutes with the left-action of $\Spin(n)$. In this way $\slashed{\mathfrak{S}}(M)$ becomes an $n$-multigraded vector bundle, where the multigrading operators are given by the right-action of the $n$ generators of $\CC\ell(n)$.

The bundles $\slashed{S}(M)$ and $\slashed{\mathfrak{S}}(M)$ are connected via two-fold periodicity. This can already be seen at the level of representations. In the even case the $2m$-multigraded inner-product space $\CC\ell(2m)$ corresponds up to isomorphism under $m$-fold application of the two-fold periodicity to the $0$-multigraded (i.e. graded) inner-product space $\Delta_{2m}$. In the odd case the $(2m+1)$-multigraded space $\CC\ell(2m+1)$ is reduced up to isomorphism via $m$-fold periodicity to the $1$-multigraded space $\Delta_{2m+1}\oplus \Delta_{2m+1}$ (with $\varepsilon_1=\begin{pmatrix}
0 & i \\ i & 0
\end{pmatrix}$). A final appeal to the equivalence between $1$-multigraded and $(-1)$-multigraded (i.e. ungraded) spaces yields the equivalence of $\CC\ell(2m+1)$ and $\Delta_{2m+1}$ under $(m+1)$-fold periodicity.

Both $\slashed{S}(M)$ and $\slashed{\mathfrak{S}}(M)$ admit a canonical metric connection that is compatible with with the Clifford multiplication (from the left). Denote these by $\nabla^{\slashed{S}}$ and $\nabla^{\slashed{\mathfrak{S}}}$ respectively. The various (multi)grading operators are parallel with respect to these connections. Define the \emph{spin Dirac operator} $\Dirac$ on $\slashed{S}(M)$ and the \emph{$\CC\ell(n)$-linear Dirac operator} $\slashed{\mathfrak{D}}$ on $\slashed{\mathfrak{S}}(M)$ by the local formulas
$$
\Dirac = \sum_{j=1}^n c(e_j) \nabla^{\slashed{S}}_{e_j} \quad , \quad \slashed{\mathfrak{D}} = \sum_{j=1}^n c(e_j) \nabla^{\slashed{\mathfrak{S}}}_{e_j} \; ,
$$
where $(e_1,\cdots,e_n)$ is a local oriented orthonormal frame, and $c(\cdot)$ denotes Clifford multiplication. Both $\Dirac$ and $\slashed{\mathfrak{D}}$ are formally self-adjoint, and elliptic with principal symbol $ic(\xi)$. If $n$ is even, so that $\slashed{S}(M)$ is graded, $\Dirac$ is odd. The operator $\slashed{\mathfrak{D}}$ is odd and $n$-multigraded.

Let $M'$ be another spin manifold of dimension $n'$ with chosen spin structure. The spin structures on $M$ and $M'$ determine a unique spin structure on $M\times M'$. There is an isomorphism $\CC\ell(n+n')\cong \CC\ell(n) \hat{\otimes} \CC\ell(n')$, where $\hat{\otimes}$ denotes the multigraded tensor product as above. It induces a well-defined isomorphism $\slashed{\mathfrak{S}}(M\times M')\cong \slashed{\mathfrak{S}}(M) \hat{\boxtimes} \slashed{\mathfrak{S}}(M')$. Under that isomorphism the Clifford-linear Dirac operators are related by the formula
$$
\slashed{\mathfrak{D}}_{M\times M'} = \slashed{\mathfrak{D}}_M \times \slashed{\mathfrak{D}}_{M'} = \slashed{\mathfrak{D}}_M \hat{\otimes} I + I \hat{\otimes} \slashed{\mathfrak{D}}_{M'} \; .
$$
Let us also draw attention to a special case of this. Suppose that $M$ has a boundary $\partial M$. The spin structure on $M$ induces a canonical spin structure on $\partial M$. In a collar neighborhood $[0,1)\times\partial M \cong U\subseteq M$ of $\partial M$ the spin structure on $M$ is then given by the product of that on $\partial M$ and the canonical spin structure on $[0,1)\subseteq \RR$. We conclude that
$$
\slashed{\mathfrak{D}}_M = \slashed{\mathfrak{D}}_\RR \hat{\otimes} I + I \hat{\otimes} \slashed{\mathfrak{D}}_{\partial M}
$$
on that collar neighborhood.
\chapter{Uniform approximability and finitely summable operators}\label{AppendixUA}
In this appendix we collect some basic facts about uniformly approxible families. Moreover we discuss finitely summable (or Schatten class) operators, because Schatten norm estimates present the main method of proving uniform approximability in the main body of this thesis.
\begin{definition}
Let $H$ be a Hilbert space. A subset $\mathfrak{A}\subseteq \fKH$ is called \emph{uniformly approximable} if for all $\varepsilon>0$, there exists an $N\in\mathbb{N}$, such that for all $T\in \mathfrak{A}$ there exists a rank-$N$ operator $K$ such that 
$$
||T-K|| < \varepsilon \; .
$$
\end{definition}
The following two lemmas are easily proven via $\frac{\varepsilon}{2}$-arguments.
\begin{lemma}\label{UAdense}
Let $\fA \subset \fKH$ be a family of compact operators. If there exists a subset $\fA_0 \subseteq \fA$, which is dense in $\fA$ and uniformly approximable, then $\fA$ is uniformly approximable. 
\end{lemma}
\begin{lemma}\label{UAapprox}
Let $\fA\subseteq \fKH$ be a family of compact operators. Suppose that for every $\varepsilon >0$ there exists a uniformly approximable family $\fA_\varepsilon$ such that for every $K\in \fA$ there exists $K_\varepsilon \in \fA_\varepsilon$ such that $||K-K_\varepsilon||\leq \varepsilon$. Then $\fA$ is also uniformly approximable.
\end{lemma}
The next lemma yields the compatibility of the notion of uniform approximability with the algebraic operations. 
\begin{lemma}\label{UAIdeal}
Let $(T_\alpha)_{\alpha\in I}$ and $(S_\alpha)_{\alpha\in I}$ be uniformly approximable families indexed by the same index set, and let $\lambda\in\CC$. 
\begin{itemize}
\item[(i)] The family $(T^*_\alpha)_{\alpha \in I}$ is uniformly approximable.
\item[(ii)] The family $(\lambda\cdot T_\alpha)_{\alpha \in I}$ is uniformly approximable.
\item[(iii)] The family $(T_\alpha+S_\alpha)_{\alpha\in I}$ is uniformly approximable.
\item[(iv)] Let $(R_\alpha)_{\alpha\in I}$ be any family of bounded operators with $\sup_{\alpha\in I} ||R_\alpha||=:C<\infty$. Then, the families $(R_\alpha T_\alpha)_{\alpha\in I}$ and $(T_\alpha R_\alpha)_{\alpha\in I}$ are uniformly approximable.
\end{itemize}
\end{lemma}
\begin{proof}
(i) follows from the fact that the adjoint of a rank-$N$ operator has rank $N$, and $||T_\alpha - K||=||T^*_\alpha-K^* ||$. (ii) is a special case of (iv). (iii) is proved by an easy $\varepsilon/2$-argument. \\
To prove (iv), fix $\varepsilon>0$. By assumption there exists $N\in\NN$ such that for any $\alpha$ there is a rank-$N$ $K_\alpha$ with $||T_\alpha-K_\alpha||<\varepsilon$. The operators $R_\alpha K_\alpha$, $K_\alpha R_\alpha$ have rank at most $N$, and $||R_\alpha T_\alpha-R_\alpha K_\alpha||, \, |||T_\alpha R_\alpha-K_\alpha R_\alpha||<C\varepsilon$. The claim follows.
\end{proof}
Consider the $C^*$-algebra $\ell^\infty(I,\fBH)$ of bounded maps from an index set $I$ to the bounded operators on some Hilbert space $H$. We can consider the subset of those $(T_\alpha)_{\alpha\in I} \in \ell^\infty(I,\fBH)$ which are uniformly approximable. The previous lemmas yield that this is an ideal\footnote{By an ideal in a $C^*$-algebra we always mean a closed two-sided $^*$-ideal} in $\ell^\infty(I,\fBH)$. In particular it is itself a $C^*$-algebra.

In the context of uniform Fredholm modules we have a slightly different situation though. There we consider maps $C_0(X)\to \fBH$ which are uniformly approximable when restricted to $\LLip_R(X)$ for any fixed $L,R\geq 0$ (see equation \eqref{LLip} for the definition of this set). These maps are usually linear, but occasionally quadratic instead. To subsume these cases, consider the $C^*$-algebra $\ell^\infty(B_X,\fBH)$, where $B_X$ denotes the unit ball in $C_0(X)$. Note that it contains $\fB(C_0(X),\fBH)$ isometrically as a subspace. Note further that $\LLip_R(X)$ is a subset of $B_X$ for any $L,R\geq 0$. The previous paragraph implies that the subset of those maps in $\ell^\infty(B_X,\fBH)$ whose restriction to $\LLip_R(X)$ is uniformly approximable forms an ideal. Let $\mathfrak{U}(X,H)$ denote the subset of those maps in $\ell^\infty(B_X,\fBH)$ whose restriction to $\LLip_R(X)$ is uniformly approximable for all $L,R\geq 0$. As an intersection of ideals it is an ideal as well. We have thus proven the following:
\begin{prop}
The subspace $\mathfrak{U}(X,H)$ is an ideal in $\ell^\infty(B_X,\fBH)$. In particular it is a $C^*$-algebra.
\end{prop}
A consequence of this proposition is the closure of $\mathfrak{U}(X,H)$ under functional calculus. For example if $\sigma: \, B_X \to \fBH$ is a bounded map taking values in the positive operators such that $\sigma^*\sigma \in \mathfrak{U}(X,H)$, then $\sigma \in \mathfrak{U}(X,H)$ as well.
\begin{definition}
Let $\sigma,\sigma' \in \ell^\infty(B_X,\fBH)$. We write $\sigma \sua \sigma'$ if $\sigma-\sigma' \in \mathfrak{U}(X,H)$. If the underlying space is not clear from context, we also write that $\sigma \sua \sigma'$ over $X$.
\end{definition}
Next, we recall the definition of finitely summable (or Schatten class) operators.
\begin{definition}
Let $s\in[1,\infty)$. A compact operator $K \in \fKH$ is called \emph{$s$-summable} if 
$$
||K||_s := \left( \sum_{\lambda\in\sigma(|K|)} \mathrm{mult}(\lambda)\cdot\lambda^s \right)^{\frac{1}{s}} \; <\infty \; .
$$
Here $\mathrm{mult}(\lambda)$ denotes the multiplicity of the eigenvalue $\lambda$. The collection of all $s$-summable operators will be denoted by $\fL^s(H)$. An operator will be called \emph{finitely summable} if it is $s$-summable for some $s\in[1,\infty)$.
\end{definition}
Note that being $1$-summable is the same as being trace-class, and being $2$-summable is the same as being Hilbert-Schmidt.
\begin{prop}[{\cite[Section 1.1]{Rave2012}}]
\begin{itemize}
\item [(i)] $\fL^s(H)$ is a vector space, and $||\cdot||_s$ defines a complete norm on it.
\item[(ii)] The set $\fL^s(H)$ is a $^*$-ideal in $\fBH$. For $K\in \fL^s(H)$ and $T,S\in\fBH$, the estimate
$$
||TKS||_s \leq ||T||\cdot ||K||_s \cdot ||S||
$$
holds.
\item[(iii)] For $s,s' \in [1,\infty)$ and $s<s'$, it holds that $||K||_s \leq ||K||_{s'}$ for all $K\in\fKH$, and thus $\fL^s(H) \subseteq \fL^{s'}(H)$.
\end{itemize}
\end{prop}
Finite summability provides a useful criterion to prove uniform approximability. The essence is \glqq Uniformly finitely summable implies uniformly approximable \grqq.
\begin{lemma}\label{UniformlySchattenClassUniformlyApprox}
Let $\fA \subseteq \fKH$ be a family of compact operators. If there exists $s\in[1,\infty)$ such that
$$
\sup_{K \in \fA} \, ||K||_s \;  = \; C\; < \infty \; ,
$$
then the family $\fA$ is uniformly approximable.
\end{lemma}
\begin{proof}
We may assume the elements of $\fA$ to be self-adjoint, otherwise consider the family $\{ \mathrm{Re}(K), \, \mathrm{Im}(K) \, | \, K\in\fA\}$ instead. Uniform approximability of this latter family implies uniform approximability of $\fA$ via Lemma \ref{UAIdeal}. 

Fix $\varepsilon>0$, and let $K\in\fA$. Because $K$ is self-adjoint, we have $||K||_s= \left(\sum_i |\lambda_i(K)|^s \right)^{1/s}$, where the $\lambda_i(K)$ are the eigenvalues of $K$, counted with multiplicity. There are only finitely many $\lambda_i(K)$ with absolute value $\geq \varepsilon$, the number of these $\lambda_i(K)$ is in fact bounded uniformly on $\fA$. This can be deduced from the estimates
$$
\big| \left\{ \, i \, | \, |\lambda_i(K)|\geq\varepsilon\, \big\} \right| \cdot \varepsilon^s \leq \sum_{|\lambda_i(K)|\geq\varepsilon} |\lambda_i(K)|^s \leq C^s \; ,
$$
implying that at most $(\frac{C}{\varepsilon})^s$ many $\lambda_i(K)$ have absolute value $\geq\varepsilon$. Pick a natural number $N\geq (\frac{C}{\varepsilon})^s$. Choose a maximal orthonormal system $\{ e_i\}_{i}$ of eigenvectors for $K$, $Ke_i= \lambda_i(K) \cdot e_i$, and define $F\in\fBH$ via
$$
Fe_i := \begin{cases} \lambda_i(K)\cdot e_i \quad & |\lambda_i(K)|\geq \varepsilon \\ 0& \mathrm{else}  \end{cases} \; .
$$
The rank of $F$ is (at most) $N$, and $||K-F||<\varepsilon$. Because the rank of $F$ does not depend on $K$, this completes the proof.
\end{proof}
\chapter{Unbounded operators} \label{AppendixUnbdOps}
Unbounded operators feature heavily in our treatment of boundary conditions. We recall basic facts about these operators here, and use them mostly implicitly throughout the main part of the thesis.

Let $H$ be a Hilbert space. An \emph{unbounded operator} on $H$ is a linear map $T: \, \dom(T) \to H$ defined on a linear subspace $\dom(T)\subseteq H$. If $S$ is another unbounded operator on $H$ such that $\dom(T)\subseteq \dom(S)$ and $Su=Tu$ for $u\in\dom(T)$, we say that $S$ is an\emph{ extension} of $T$, and write $T\subseteq S$. The operators are equal if $T\subseteq S$ and $S\subseteq T$.

The operator $T$ is called \emph{closed} if the graph of $T$ is closed in $H\times H$. This is equivalent to the condition that the graph norm $||\cdot||_T=||\cdot||+||T\cdot||$ define a complete norm on $\dom(T)$. Moreover, $T$ is called \emph{closable} if $T$ has a closed extension. In that case the minimal closed extension $\bar{T}$ of $T$ is called the \emph{closure} of $T$.
\begin{lemma}
Let $T$ be closable, and $u\in H$. Then, the following are equivalent:
\begin{itemize}
\item[(i)] $u\in \dom(\bar{T})$.
\item[(ii)] There exists a sequence $(u_n)_{n\in\NN}$ in $\dom(T)$ such that $||u-u_n||_H\to 0$ and $(u_n)_n$ is a Cauchy sequence with respect to the graph norm of $T$. Equivalently, there exists $v\in H$ such that $(u,v)$ lies in the closure of $\mathrm{Graph}(T)\subseteq H\oplus H$. In this case $\bar{T}u = v = \lim_{n\to \infty} Tu_n$.
\item[(iii)] There exists a sequence $(u_n)_{n\in\NN}$ in $\dom(T)$ such that $||u-u_n||_H\to 0$ and $\sup_{n\in\NN} ||u_n||_H<\infty$.
\end{itemize}
\end{lemma}
\begin{proof}
The equivalence of (i) and (ii) is Proposition 1.4 in \cite{Schmudgen2012}, the equivalence of (i) and (iii) is Lemma 1.8.1 in \cite{HigsonRoe2000}.
\end{proof}
Assume that $T$ is densely defined, i.e. $\dom(T)\subseteq H$ is dense. The adjoint $T^*$ of $T$ is defined to have domain
$$
\dom(T^*) := \left\{ \, u\in H \, | \, \exists v \in H \, \forall \, w\in \dom(T) : \, \langle u,Tw\rangle = \langle v,w\rangle\, \right\} \; .
$$
Since $T$ is densely defined, $v$ is uniquely determined by this condition, and we set $T^*u:=v$. The operator $T$ is called \emph{symmetric} if $T\subseteq T^*$, \emph{self-adjoint} if $T^*=T$, and \emph{essentially self-adjoint} if $T$ is closable and $\bar{T}$ is self-adjoint.
\begin{lemma}[{\cite[Prop. 1.6, Theorem 1.8]{Schmudgen2012}}]
\begin{itemize}
\item $T^*$ is a closed operator.
\item $T$ is closable if and only if $T^*$ is densely defined. In that case, $T^*=\bar{T}^*$ and $\bar{T}=T^{**}$.
\item If $T\subseteq S$, then $S^*\subseteq T^*$.
\end{itemize}
\end{lemma}
Let $T$ be an unbounded operator on $H$. The \emph{resolvent set} of $T$ is the set $\rho(T)$ containing those $\lambda\in\CC$ for which $T-\lambda: \, \dom(T)\to H$ has a bounded inverse $(T-\lambda)^{-1}: \, H\to\dom(T)$. The \emph{spectrum} of $T$ is the complement $\sigma(T):=\CC\setminus\rho(T)$. The spectrum of an operator is a closed set. It is not generally true that the spectrum of a symmetric operator is real. In fact it holds for a closed symmetric operator $T$ that $T$ is self-adjoint if and only if $\sigma(T)\subseteq \RR$.
\begin{lemma}\label{ResolventGraphNorm}
Let $T$ be a closed unbounded operator on $H$. If $\lambda \notin \sigma(T)$, then $(T-\lambda)^{-1}$ is continuous as a map $(H,||\cdot||_H) \to (\dom(T),||\cdot||_T)$.
\end{lemma}
\begin{proof}
The map $T-\lambda: \, (\dom(T),||\cdot||_T)\to (H,||\cdot||_H)$ is continuous. Since $T$ is closed $(\dom(T),||\cdot||_T)$ is a Banach space. Thus the open mapping theorem implies that $(T-\lambda)^{-1}: \, (H,||\cdot||_H) \to (\dom(T),||\cdot||_T)$ is continuous as well.
\end{proof}
\begin{lemma}[{\cite[Prop. 3.18, Lemma 5.8]{Schmudgen2012}}]\label{BdTrans}
Let $T$ be a closed operator. 
\begin{itemize}
\item[(i)] $T^*T$ and $TT^*$ are positive self-adjoint operators.
\item[(ii)] The operator $(T^*T+I)^{-1}$ is bounded and positive. Its square root $(T^*T+I)^{-1/2}$ has image contained in $\dom(T)$.
\item[(iii)] The operator $Z_T:=T(T^*T+I)^{-1/2}$ is bounded with $||Z_T||\leq 1$. It is called the bounded transform of $T$.
\item[(iv)] It holds that $T(T^*T+I)^{-1/2}=(TT^*+I)^{-1/2}T^*$, i.e. $Z_T=(Z_{T^*})^*$. In particular, if $T$ is self-adjoint, so is $Z_T$.
\end{itemize}
\end{lemma}
\begin{thm}[{\cite[Theorem 5.7, Theorem 5.9]{Schmudgen2012}}]\label{funcCalc}
Let $T$ be a self-adjoint operator. There exists a unique spectral measure $E_T$ on the Borel $\sigma$-algebra such that
$$
T=\int_\RR \lambda \, dE_T(\lambda) \; .
$$
Then, if $f: \, \RR\to \CC$ is a measurable function, the operator 
$$
f(T) := \int_\RR f(\lambda) \, dE_T(\lambda) \quad \mathrm{with} \quad \dom(f(T)) = \left\{ u\in H \, | \, \int_\RR |f(\lambda)|^2 \, d\langle E_T(\lambda)u,u\rangle <\infty  \right\}
$$ 
is a normal operator with dense domain. For $u\in\dom(f(T))$ it then holds that
$$ 
\left\| f(T)u \right\|^2 = \int_{\RR} |f(\lambda)|^2 \, d\langle E_T(\lambda)u,u\rangle \; .
$$
If $f=\sum a_i x^{i}$ is a polynomial, then $f(T)=\sum_ia_i  T^{i}$. If $g$ is another measurable function, then $(fg)(T)=f(T)g(T)$. In particular, $f(T)$ and $g(T)$ commute. If $f$ is bounded, so is $f(T)$, and the assignment
$$
\cL^\infty(\RR) \longrightarrow \fBH \; , \; f\longmapsto f(T)
$$
is a $^*$-homomorphism.
\end{thm}
\chapter{Operator K-theory}\label{AppendixKTheory}
In this appendix we review some relevant basics on operator K-theory. We use the definitions of the K-theory groups from \cite{HigsonRoe2000}. For additional details on unitization and induced morphisms, see for example \cite{WeggeOlsen1993}.
\begin{definition}
Let $A$ be a unital $C^*$-algebra. Define $K_0(A)$ as the free abelian group with one generator $[p]$ for each projection $p$ in some matrix algebra $M_n(A)$, modulo the following relations:
\begin{itemize}
\item If $p_t$ is a continuous path of projections from $p_0$ to $p_1$ in $M_n(A)$, then $[p_0]=[p_1]$.
\item $[0]=0$ for the $0$-matrix of any size.
\item $[p\oplus q]=[p]+[q]$ for any projections $p,q$.
\end{itemize}
\end{definition}
Let $\alpha: \, A\to B$ be a unital $^*$-homomorphism between unital $C^*$-algebras. Via entry-wise application, $\alpha$ gives rise to $^*$-homomorphisms $M_n(A)\to M_n(B), \, (a_{ij})\mapsto (\alpha(a_{ij}))$, which map projections to projections. These maps descend to give a well-defined homomorphism $\alpha_*: \, K_0(A) \to K_0(B), \, [p] \mapsto [\alpha(p)]$. \\
Recall that the unitization of a $C^*$-algebra $J$ is defined as the unital $C^*$-algebra $J^+:= J \times \CC$ with sum and adjoint defined component-wise, multiplication given by $(a,\lambda)(b,\mu):=(ab+\lambda b + \mu a,\lambda\mu)$, and equipped with a suitable norm. The map $\pi: \, J^+\to\CC, \, (a,\lambda)\mapsto \lambda$ is a unital $^*$-homomorphism. If $\beta: \, J\to I$ is a $^*$-homomorphism, then $\beta^+: \, J^+\to I^+$, $\beta^+(a,\lambda):= (\beta(a),\lambda)$ is a unital $^*$-homomorphism. If $J$ is already unital, then $J^+$ is isomorphic to the direct sum $J\oplus\CC$. An isomorphism is given by $J^+\to J\oplus \CC,\, (a,\lambda)\mapsto (a-\lambda \mathbf{1}_J,\lambda)$. If $\beta: \, J\to I$ is a not necessarily unital $^*$-homomorphism between unital $C^*$-algebras, then the induced map on unitizations becomes $\beta^+: \, J\oplus \CC \to I\oplus \CC, \, (a,\lambda)\mapsto (\phi(a)+\lambda(\mathbf{1}_I-\phi(\mathbf{1}_J)),\lambda)$. 
\begin{definition}
Let $J$ be a $C^*$-algebra. Define $K_0(J)$ as the kernel of the map $\pi_*: \, K_0(J^+)\to K_0(\CC)=\ZZ$. If $\beta: \, J\to I$ is a $^*$-homomorphism, then define $\beta_*: \, K_0(J)\to K_0(I)$ as the restriction of $\beta^+_*: \, K_0(J^+)\to K_0(I^+)$ to $K_0(J)$. Then, the zeroth K-theory assembles into a functor $K_0: \, C^*\mathrm{Alg}\to \mathrm{Ab}$.
\end{definition}
If $J$ is already unital, then the second definition of $K_0(J)$ coincides with the first. If $\beta: \, J\to I$ is a non-unital $^*$-homomorphism between unital $C^*$-algebras, then $\beta_*: \, K_0(J) \to K_0(I)$ still maps the class of a projection $p$ over $J$ to the class of the projection $\beta(p)$ over $I$.
\begin{definition}
Let $A$ be a unital $C^*$-algebra. Define $K_1(A)$ as the free abelian group with one generator $[u]$ for each projection $u$ in some matrix algebra $M_n(A)$, modulo the following relations:
\begin{itemize}
\item If $u_t$ is a continuous path of unitaries from $u_0$ to $u_1$ in $M_n(A)$, then $[u_0]=[u_1]$.
\item $[\mathbf{1}]=0$ for the identity matrix of any size.
\item $[u\oplus v]=[u]+[v]$ for any unitaries $u,v$.
\end{itemize}
\end{definition}
Let $\alpha: \, A \to B$ be a unital $^*$-homomorphism between unital $C^*$-algebras. The induced unital $^*$-homomorphism between the matrix algebras then takes unitaries to unitaries, thus descends to a homomorphism $\alpha_*: \, K_1(A) \to K_1(B), \, [u]\mapsto [\alpha(u)]$. 

Let $J$ be a not necessarily unique $C^*$-algebra. Define $U^+_n(J)$ to be the set of unitaries in $u\in M_n(J^+)$ with $\pi(u)=\mathrm{1}_n$. Note that if $\beta: \, J \to I$ is any $^*$-homomorphism, then the unital $^*$-homomorphism $\beta^+: \, J^+\to I^+$ maps $U^+_n(J)$ to $U^+_n(I)$.
\begin{definition}
Let $J$ be a $C^*$-algebra. Define $K_1(J)$ as the free abelian group with one generator for each element in $U^+_n(J)$, modulo the following relations:
\begin{itemize}
\item If $u_t$ is a continuous path of unitaries from $u_0$ to $u_1$ in $U^+_n(J)$, then $[u_0]=[u_1]$.
\item $[\mathbf{1}_n]=0$ for $\mathbf{1}_n \in U^+_n(J)$ with entries $(0,\delta_{ij})$. 
\item $[u\oplus v]=[u]+[v]$ for any unitaries $u,v$.
\end{itemize}
If $\beta: \, J \to I$ is any $^*$-homomorphism, then define $\beta_*: \, K_1(J) \to K_1(I), \, [u] \mapsto [\beta^+(u)]$. Then, the first K-theory assembles into a functor $K_1: \, C^*\mathrm{Alg} \to \mathrm{Ab}$.
\end{definition}
If $J$ is already unitary, then the elements of $U^+_n(J)$ are given by the sum of a unitary in $M_n(J)$ and the identity matrix in $\CC$. Moreover, under the isomorphism $M_n(J^+)\cong M_n(J)\oplus M_n(\CC)$ the matrix $\mathbf{1}_n \in M_n(J^+)$ with entries $(0,\delta_{ij})$ corresponds to the sum of the $J$-valued identity matrix with the $\CC$-valued identity matrix. Thus, the second definition of $K_1(J)$ coincides with the first. Let $\beta: \, J\to I$ be a not necessarily unital $^*$-homomorphism between unital $C^*$-algebras. Then, if $(u,\mathbf{1}_n) \in U^+_n(J)$, the element $\beta^+(u,I) \in U^+_n(I)$ corresponds to the unitary $\beta(u)+(\mathbf{1}_I - \beta(\mathbf{1}_J)) \in M_n(I)$. Thus, the induced map on K-theory is $\beta_*: \, K_1(J)\to K_1(I), \, [u]\mapsto [\beta(u)+(\mathbf{1}_I - \beta(\mathbf{1}_J))]$.
\begin{thm}
Let $A$ be a $C^*$-algebra, and $J\subseteq A$ an ideal\footnote{Again, by an ideal we mean a closed two-sided $^*$-ideal.}. There is a natural 6-term exact sequence \\
\centerline{
\xymatrix{
K_0(J) \ar[r]^-{\iota_*} & K_0(A)\ar[r]^-{\pi_*} & K_0(A/J)\ar[d]^\partial \\
K_1(A/J) \ar[u]^{\partial} & \ar[l]^-{\pi_*} K_1(A) & \ar[l]^-{\iota_*} K_1(J)
}
} 
\end{thm}
Let us describe the boundary maps. For $\partial: \, K_0(A/J)\to K_1(J)$ let $p$ be a projection in $A/J$, and let $x$ a self-adjoint lift to $A$. Note that $\exp(2\pi i x)$ is a unitary in $J^+$ which equals $1$ modulo $J$. Then
$$
\partial [p] = [\exp(2\pi i x)] \; .
$$
For $\partial: \, K_1(A/J)\to K_0(J)$ let $u$ be a unitary in $A/J$, and $a\in A$ a lift of $u$ with $||a||=1$. Let 
$$
P := \begin{pmatrix}
aa^* & a(1-a^*a)^{1/2} \\ a^*(1-aa^*)^{1/2} & 1-a^*a
\end{pmatrix} \; \in M_2(J^+) \; .
$$
Modulo $J$, $P$ equals $Q=\begin{pmatrix}
1&0 \\ 0&0
\end{pmatrix}$. Then,
$$
\partial[u] = [P]-[Q] \; .
$$
In particular, if $a$ can be chosen to be a partial isometry, so that $1-a^*a$ and $1-aa^*$ are projections in $J$, then
$$
\partial[u] = [1-a^*a]-[1-aa^*] \; .
$$
In our discussion around relative uniform index maps we also need the Mayer-Vietoris sequence in operator K-theory. It can be constructed straight-forwardly from the 6-term exact sequence.
\begin{prop}
Let $A$ be a $C^*$-algebra, and $I,J \subseteq A$ ideals such that $I+J=A$. Denote by $i: \, I\hookrightarrow A$ and $j: \, J\hookrightarrow A$ the respective inclusions. Then, there is a Mayer-Vietoris sequence \\
\centerline{
\xymatrix{
K_0(I\cap J) \ar[r]^-{i_* + j_*} & K_0(I)\oplus K_0(J)\ar[r]^-{i_*-j_*} & K_0(A)\ar[d]^{\delta_{MV}} \\
K_1(A) \ar[u]^{\delta_{MV}} & \ar[l]^-{i_*-j_*}K_1(I)\oplus K_1(J) & \ar[l]^-{i_*+j_*} K_1(I\cap J)
}}
\end{prop}
\begin{proof}
Consider the $C^*$-algebra $I\oplus J$. It contains $I\oplus J$ as an ideal via $I\cap J \hookrightarrow I\oplus J, \, a\mapsto (a,a)$. The condition that $A=I+J$ translates to the fact that $(I\oplus J)/(I\cap J) \to A, \, [a,b]\mapsto a-b$ is a well-defined isomorphism. Then, the 6-term exact sequence of the pair $(I\oplus J,I\cap J)$ produces the Mayer-Vietoris sequence.
\end{proof}
We also need a formula connecting the Mayer-Vietoris boundary map to the boundary map of the 6-term exact sequence. By definition $\delta_{MV}$ is the boundary map of the 6-term exact sequence of the pair $(I\oplus J,I\cap J)$. Consider the map $\pi: \, (I\oplus J)/(I\cap J) \to I/(I\cap J)$ induced by projection onto the first summand. Note also that the inclusion $I\hookrightarrow A$ induces an isomorphism $I/(I\cap J)\cong A/J$. Naturality of the boundary map of the 6-term exact sequence provides us with the diagram\\
\centerline{\xymatrix{
K_{*+1}(A) \ar@{=}[r] \ar[dr]_{\pi_*} & K_*\left( (I\oplus J)/(I\cap J) \right) \ar[d]_{\pi_*} \ar[r]^-{\delta_{MV}} & K_*(I\cap J) \ar@{=}[d] \\
& K_{*+1}\left(I/(I\cap J)\right) \ar[r]_-\partial & K_*(I\cap J) 
}}
Phrased differently, the Mayer-Vietoris boundary map is given by the composition
\begin{equation}\label{MVformula}
K_{*+1}(A) \xlongrightarrow{\pi_*} K_{*+1}(A/J)=K_{*+1}(I/(I\cap J)) \xlongrightarrow{\partial} K_*(I\cap J) \; .
\end{equation}
There is also a different construction of the Mayer-Vietoris boundary map, which uses suspensions and a sort of mapping-cone construction. This picture is generally less amenable to direct computation, but has the advantage that here the Mayer-Vietoris boundary map is induced by a $^*$-homomoprhism between auxiliary algebras. A detailed description of this approach to the Mayer-Vietoris sequence as well as the equivalence to the boundary map described here can be found in Sections 5.1 and 5.3 of \cite{Siegel2012}. 

We conclude this brief review of operator K-theory with a proposition on the compatibility of the Mayer-Vietoris sequence with quotients and the 6-term exact sequence. The statement well-known (it is certainly implicit in \cite{Siegel2012}), but the author could not find a reference.
\begin{prop}\label{MVQuot}
Let $A$ be a $C^*$-algebra, and let $B,I_1,I_2 \subseteq A$ be ideals\footnote{Recall our convention that by an ideal in a $C^*$-algebra we always mean a two-sided $^*$-ideal.} such that $A=I_1+I_2$. Set $J_k = B\cap I_k$, $k=1,2$.
\begin{itemize}
\item[(i)] Let $Q_k$, $k=1,2$, denote the images of the maps $I_k/J_k \to A/B$. Then $Q_1$, $Q_2$ are ideals in $A/B$ such that $A/B=Q_1+Q_2$, and $Q_1\cap Q_2$ coincides with the image of $I_1\cap I_2/J_1\cap J_2 \to A/B$.
\item[(ii)] The diagram \\
\centerline{\xymatrix{
K_{*+1}(A/B) \ar[r]^-\partial \ar[d]_\delta & K_*(B) \ar[d]^\delta \\ 
K_*(I_1\cap I_2/J_1\cap J_2) \ar[r]_-\partial & K_{*-1}(J_1\cap J_2)
}}
commutes, where $\partial$ are the boundary maps from the respective 6-term exact sequences, and $\delta$ are the respective Mayer-Vietoris boundary maps.
\end{itemize}
\end{prop}
\begin{proof}
Ad (i): Note first that $J_k$ are ideals in $B$, and that $B=J_1+J_2$. Using that $I_k$ are ideals in $A$ and $J_k$ ones in $B$ it is straight-forward to check that $Q_k$ are ideals in $A/B$. Note also that by definition the quotient map $I_k/J_k \to Q_k$ is an isomorphism. Moreover, the image of $I_1\cap I_2$ in $A/B$ coincides with $Q_1\cap Q_2$. That the former is contained in the latter is clear. The reverse inclusion follows from the fact that an element of $Q_1\cap Q_2$ is represented by an element $a\in A$ such that $a+b_k \in I_k$ for suitable $b_k\in B$. Then, we may pick a partition of unity $i_1^2+i_2^2=I$, $i_k\in I_k$, and write
$$
a = i_1^2 a i_1^2 + i_1^2 a i_2^2 + i_2^2 a i_1^2 + i_2^2 a i_2^2 \; .
$$
The middle summands are automatically in $I_1\cap I_2$, and using that $a+ b_k\in I_k$ it follows that the first and last summand are in $I_1\cap I_2$ modulo $B$. We conclude that $Q_1\cap Q_2 \subseteq I_1\cap I_2/B$.

Ad (ii): We use the aforementioned description of the Mayer-Vietoris boundary map via suspensions and mapping cones. The reader is referred to \cite{Siegel2012} as well as \cite{HigsonRoe2000} for some standard facts about these constructions, which we will simply take for granted here. Let $\Omega(A,I_1,I_2)$ denote the $C^*$-algebra of continuous maps $f: \, [0,1]\to A$ with $f(0)\in I_1$, $f(1)\in I_2$. There is a natural inclusion map from $SA=C_0(0,1)\otimes A$ into $\Omega(A,I_1,I_2)$. Moreover, there are natural isomorphisms $K_{*}(SA)\cong K_{*+1}(A)$ and $K_{*}(\Omega(A,I_1,I_2))\cong K_*(I_1\cap I_2)$. Then, the composition
$$
K_{*+1}(A) \cong K_{*}(SA) \longrightarrow K_{*}(\Omega(A,I_1,I_2))\cong K_*(I_1\cap I_2)
$$
realizes the Mayer-Vietoris boundary map. There is a diagram \\
\centerline{\xymatrix{
0 \ar[r] & SB \ar[r] \ar[d]& SA \ar[r]\ar[d] & S(A/B) \ar[d]\ar[r] & 0 \\
0 \ar[r] & \Omega(B,J_1,J_2) \ar[r] & \Omega(A,I_1,I_2) \ar[r] & \Omega(A/B,Q_1,Q_2) \ar[r] & 0
}}
whose rows are exact. Passing to K-theory and exploiting naturality of the boundary map provides us with the commutative diagram \\
\centerline{\xymatrix{
K_{*+1} (S(A/B)) \ar[r]^\partial \ar[d] & K_*(SB) \ar[d] \\
K_{*+1} (\Omega(A/B,Q_1,Q_2)) \ar[r]_\partial & K_*(\Omega(B,J_1,J_2))
}}
The identifications $K_{*+1} (S(A/B))\cong K_{*+2}(A/B)$ and $ K_*(J_1\cap J_2)\cong K_*(\Omega(B,J_1,J_2))$ are each compatible with boundary maps, in the first case by definition, in the second because it is induced by the inclusion $J_1\cap J_2\hookrightarrow\Omega(B,J_1,J_2)$. The claim follows.
\end{proof}
\chapter{Geometry of metric pairs}\label{AppendixUnifPairs}
In our discussion of excision and the long exact sequence in uniform K-homology we had to assume that the geometry of the pair $(X,Z)$ is well-behaved in a certain technical way. It is our goal in this appendix to provide sufficient conditions for this well-behavedness criterion in terms of more standard notions. Let us first re-state the relevant definitions. 
\begin{definition}\label{DefUnifPair}
Let $(X,Z)$ be a pair. We say that $X$ is \emph{uniformly resolvable} near $Z$ if there is a partition of unity $\{\psi_i\}_i$ on $X\setminus Z$ with the following properties:
\begin{itemize}
\item[(1)] For all $\delta>0$ there exists $L(\delta)\geq 0$ such that $\sqrt{\psi_i}$ is $L(\delta)$-Lipschitz if $d(\supp(\psi_i),Z)\geq\delta$.
\item[(2)] For all $R\geq 0$ and $\delta>0$, there exists $M(R,\delta)\in\NN$ such that for all compact $K\subseteq X\setminus Z$ with $\diam(K)\leq R$ and $d(K,Z)\geq\delta$ it holds that
$$
 \left| \left\{ i \, | \, \psi_i|_K \neq 0 \right\}\right| \leq M(R,\delta) \; .
$$
\item[(3)] For all $\varepsilon>0$ there exists $\delta>0$ such that for all $i$ it holds that
$$
\diam(\supp(\psi_i)) \geq \varepsilon \; \Rightarrow \; d(\supp(\psi_i),Z) \geq \delta \; .
$$
\end{itemize}
Furthermore we say that $X$ is \emph{transversally controlled} near $Z$ if the restriction map $C_0^+(X)\to C_0^+(Z)$ admits an approximately filtered positive section, meaning that there exists a positive map $\sigma: \, C_0^+(Z)\to C_0^+(X)$ such that $\sigma(f)|_Z=f$ for all $f\in C_0^+(Z)$, and for all $L,R\geq 0$ and $\varepsilon>0$ there exist $L',R'\geq 0$ such that $\sigma\big(\LLip_R(Z)\big)$ is contained in the $\varepsilon$-neighborhood of $L'\text{-}\mathrm{Lip}_{R'}(X)$.
\end{definition}
We first give a sufficient condition for uniform resolvability in terms of the existence of certain open covers.
\begin{lemma}\label{LemmaLipPU}
Suppose there exists an open cover $\cU$ of $X\setminus Z$ with the following properties:
\begin{itemize}
\item[(C1)] For all $R\geq 0$ and $\delta>0$, there exists $M(R,\delta)\in\NN$ such that for all compact $K\subseteq X\setminus Z$ with $\diam(K)\leq R$ and $d(K,Z)\geq\delta$ it holds that
$$
 \left| \left\{ U\in\cU \, | \, U\cap K\neq\emptyset \right\}\right| \leq M(R,\delta) \; .
$$
\item[(C2)] For all $\varepsilon>0$ there exists $\delta>0$ such that for all $U\in\cU$ it holds that
$$
\diam(U) \geq \varepsilon \; \Rightarrow \; d(U,Z) \geq \delta(\varepsilon) \; .
$$
\item[(C3)] For all $\delta>0$ there exists $\lambda>0$ such that for all $x\in X$, $d(x,Z)\geq\delta$ there exists $U\in \cU$ with $B_\lambda(x)\subseteq U$.
\end{itemize} 
Then $X$ is uniformly resolvable near $Z$.
\end{lemma}
\begin{proof}
Let $\cU=\{U_i\}_i$ be an open cover of $X\setminus Z$ satisfying (C1), (C2) and (C3). Define 
$$
\psi_i(x):= d(x,X\setminus U_i)^2 \bigg( \sum_j d(x,X\setminus U_j)^2 \bigg)^{-1} \; .
$$
By construction $\psi_i^{-1}(0,1]\subseteq U_i$, so that points (2) and (3) of the definition of uniform resolvability follow from properties (C1) and (C2) respectively. Thus it remains to prove (1). To that end suppose that $d(U_i,Z)\geq\delta$ calculate that
$$
|\sqrt{\psi_i}(x)-\sqrt{\psi_i}(y)| \leq L_i \cdot \left|d(x,X\setminus U_i)-d(y,X\setminus U_i)\right| \leq L_i\cdot d(x,y)
$$
by the reverse triangle inequality, where $L_i:=\sup_{x\in U_i}\left( \sum_j d(x,X\setminus U_j)^2 \right)^{-1/2}$. By property (C3) there exists $\lambda>0$ such that the $\lambda$-ball around any $x\in U_i$ is contained in some $U_j$. It follows that $\sum_j d(x,X\setminus U_j)^2 \geq \lambda^2$, so that $L_i\leq \lambda^{-1}$.
\end{proof}
If $X$ has jointly bounded geometry (Definition \ref{DefJBB}) and $Z$ is an arbitrary closed subspace, one can always produce such a cover by means of an iterative refinement process, starting with a cover from a quasi-lattice and successively refining ever closer to $Z$. This is done in the following proposition.
\begin{prop}
Let $X$ be a metric space of jointly bounded geometry, and $Z\subseteq X$ an arbitrary closed subspace. Then there exists an open cover $\cU$ of $X\setminus Z$ with properties (C1), (C2), (C3). 

Consequently, a metric space of jointly bounded geometry is uniformly resolvable near any closed subspace.
\end{prop}
\begin{proof}
Fix a Borel decomposition $X=\bigcup_{i\in I} X_i$ and a compatible quasilattice $\Gamma\subseteq X$ as in the definition of jointly bounded geometry. Let $c>0$ be sufficiently large that the open $c$-balls around $\Gamma$ cover $X$, and let $\cU_0$ be the open cover of $X\setminus Z$ obtained by intersecting each $B_c(y)$, $y\in\Gamma$, with $X\setminus Z$. Because of jointly bounded geometry the cover $\cU_0$ has the property
\begin{itemize}
\item[(U)] For every $R\geq 0$ there exists $N\in\NN$ such that every compactum $K\subseteq X\setminus Z$, $\diam(K)\leq R$, intersects at most $N$ elements of the open cover.
\end{itemize}
Fix a strictly decreasing null-sequence $(\delta_n)_{n\in\NN}$. We will construct a sequence of open covers $\cU_n$, $n\in\NN$, of $X\setminus Z$ such that
\begin{itemize}
\item  $\cU_{n+1}$ refines $\cU_n$.
\item  $\cU_n$ satisfies (U).
\item  If $U\in\cU_{n+1}$ satisfies $d(U,Z)\geq\delta_n$, then $U\in\cU_{n}$.
\item  If $d(U,Z)<\delta_n$ for $U\in\cU_{n+1}$, then $\diam(U)\leq \delta_n$.
\end{itemize}

We have already constructed $\cU_0$. Assume that the cover $\cU_n$ has already been constructed. Suppose $U\in\cU_n$ satisfies $d(U,Z)<\delta_n$. Because $U_n$ refines $\cU_0$ the set $U$ is contained in a $c$-ball around an element of $\Gamma$. Thus there is a uniform bound on the number of $X_i$ that $U$ intersects. By assumption on $X_i$ there exists a $\frac{\delta_{n+1}}{4}$-net $\Gamma_{n+1}^{i}$ of uniformly bounded cardinality in each $X_i$. Thus there is a uniform bound on how many sets $B_{\delta_{n+1}}(y)\cap U$, $y\in\Gamma_{n+1}^{i}$, $i\in I$, are non-empty. Construct $\cU_{n+1}$ by replacing each $U\in\cU_{n}$ with $d(U,Z)<\delta_n$ by the collection of those $B_{\delta_{n+1}}(y)\cap U$, $y\in\Gamma_{n+1}^{i}$, $i\in I$, that are non-empty. Then $\cU_{n+1}$ is a refinement of $\cU_n$, and the third and fourth bullet point above hold by construction. Moreover, since $\cU_n$ satisfies (U) and every element of $\cU_n$ is replaced by a uniformly finite number of elements in $\cU_{n+1}$ the cover $\cU_{n+1}$ again satisfies (U).  

Now define $\cU$ as the limit of $\cU_n$: An open set $U\subseteq X\setminus Z$ is in $\cU$ if and only if $d(U,Z)>0$ and $U\in\cU_n$ for one and hence any $n\in\NN$ such that $d(U,Z)\geq\delta_n$. To see that this is indeed an open cover of $X\setminus Z$, let $x$ be a point with $d(x,Z)=\delta>0$ and let $\delta_n$ be less than $\delta$. There is an element $U\in\cU_{n+1}$ containing the point $x$. If $d(U,Z)\geq\delta_n$, then $U\in\cU$ by definition. If $d(U,Z)<\delta_n$, then $\diam(U)\leq\delta_n$. Since on the other hand $d(x,Z)>\delta_n$ we conclude that $d(U,Z)\geq\delta_m>0$ for some $m>n$. Since $\cU_m$ refines $\cU_n$, $x$ is contained in an element of $\cU_m$, and that element has distance $\geq\delta_m$ to $Z$. Thus it lies in $\cU$. This concludes the argument that $\cU$ is an open cover of $X\setminus Z$. 

We will prove that $\cU$ satisfies (C1), (C2) and (C3). For (C1) let $K$ be a compactum in $X\setminus Z$ with $\diam(K)\leq R$ and $d(K,Z)\geq\delta$. Let $\delta_n$ be strictly less than $\delta$. Because $\cU_{n+1}$ satisfies (U) the set $K$ intersects at most $N_{n+1}(R)$ elements of $\cU_{n+1}$. Analogously to the argument above, there must be $m\geq n+1$ such that each of these open sets has distance $\geq \delta_m$ to $Z$. Additionally choose $m$ such that $2\delta_m<\delta_n$. The compactum $K$ intersects at most $N_m(R)$ many elements of $\cU_m$. Now, assume $U\in\cU$ intersects $K$. Either $d(U,Z)\geq\delta_m$, whence $U\in\cU_m$, or $d(U,Z)<\delta_m$. In this case it must hold that $\diam(U)\leq \delta_m$. This would however yield
$$
\delta_n \leq d(K,Z) \leq \diam(U)+d(U,Z) < 2\delta_m \; ,
$$
a contradiction. Any element of $\cU$ intersecting $K$ must therefore be contained in $\cU_m$, and there can be at most $N_m(R)$ such elements - a number depending only on $\cU$, $R$ and $\delta$, but not the specific $K$. This establishes (C1).

For (C2) suppose that $U\in\cU$ has diameter $\geq \varepsilon>0$. The set $U$ lies in some $\cU_{n+1}$, where we may take $n$ such that $\delta_n\leq \varepsilon$. The properties of $\cU_{n+1}$ then yield that $d(U,Z)\geq\delta_n$.

Lastly, we establish (C3). Let $x\in X$ have $d(x,Z)\geq\delta>0$. As a consequence of our argument for (C1) it holds that every $U\in \cU$ that contains $x$ is an element of $\cU_m$ for a suitable $m$. Thus it suffices to prove that the covers $\cU_m$ has positive Lebesgue number. We will do this via induction. We defined $\cU_0$ as the collection of $c$-balls around the quasi-lattice $\Gamma$. After potentially increasing $c$ we may assume that there is $0<c'<c$ such that the $c'$-balls around $\Gamma$ still cover $X$. Then, every $x$ is in the $c'$-ball of some $y\in\Gamma$, and the $\frac{c-c'}{2}$-ball around $y$ is contained in $B_c(y)\in \cU_0$.

Now assume that $\cU_m$ has positive Lebesgue number. Recall that $\cU_m$ was constructed by replacing certain $U\in \cU_m$ with sets $B_{\delta_{m+1}}(y)\cap U$, where $y\in\Gamma^{i}_{m+1}$ is a $\frac{\delta_{m+1}}{4}$-net in $X_i$. Take $x\in X_i$, and let $U\in\cU_m$ contain a $\lambda$-ball around $x$. If $U$ is not changed upon passing to $\cU_{m+1}$, there is nothing to show, so we suppose that it is. By construction $x$ is within $\frac{\delta_{m+1}}{4}$ of an element $y\in\Gamma_{m+1}^{i}$, so that the $\frac{\delta_{m+1}}{4}$-ball around $x$ is contained in $B_{\delta_{m+1}/2}(y)$. Choosing $\lambda'$ smaller than both $\lambda$ and $\frac{\delta_{m+1}}{4}$ we conclude that $B_{\lambda'}(x)\subseteq B_{\delta_{m+1}/2}(y)\cap U$. This proves (C3) and hence the proposition.
\end{proof}
Next we give a sufficient condition for transversal control in proper metric spaces in terms of the existence of $\varepsilon$-dense subsets of $Z$ that are uniformly locally finite in $X$. The proof is an adaptation of the proof that surjections of unital function algebras admit (completely) positive splits found in \cite[Section 3.2]{HigsonRoe2000}.
\begin{prop}\label{FilteredCP}
Let $X$ be a proper metric space, and $Z\subseteq X$ a closed subspace such that
\begin{itemize}
\item[(F)] For all $\varepsilon>0$ there exists an $\varepsilon$-dense subset\footnote{Which is to say that $N_\varepsilon(\Gamma)=Z$.} $\Gamma\subseteq Z$ such that 
$$
\sup_{x\in X} \left| \left\{ y\in\Gamma \, | \, d(x,y)<2\varepsilon \right\} \right| < \infty \; .
$$
\end{itemize} 
Then $X$ is transversally controlled near $Z$.
\end{prop}
\begin{proof}
The cases $Z=\emptyset$ and $Z=X$ are trivial, we will therefore take $Z$ to be a non-empty proper subset of $X$. We fix the following data: Let $\{\varepsilon_n\}_{n\in\NN}$ be a strictly decreasing null-sequence such that $2(\varepsilon_n+\varepsilon_{n+1})<2^{-n}$. Let $\Gamma_n=\{x^{(n)}_i\}_i$ be an $\varepsilon_n$-dense subset of $Z$ provided by (F). Consider the open cover $\mathcal{U}_n$ of $X$ given by the open $2\varepsilon_n$-ball in $X$ centered on the elements of $\Gamma_n$, together with the complement of the closed $\frac{\varepsilon_n}{2}$-neighborhood of $\Gamma_n$ in $X$. The assumption (F) guarantees that every $\cU_n$ is a uniformly locally finite cover. Choose a subordinate partition of unity $\{\psi_i^{(n)}\}_i\cup \{\psi_\infty^{(n)}\}$ such that each $\psi_i^{(n)}$ is $l_n$-Lipschitz with uniform Lipschitz constant $l_n<\infty$. Here $\psi^{(n)}_i$ is supported in $U_{2\varepsilon_n}(x^{(n)}_i)$ and $\psi_\infty^{(n)}$ is supported in the complement of the $\frac{\varepsilon_n}{2}$-neighborhood of $\Gamma_n$. Note that the restrictions of the $\{\psi^{(n)}_i\}$ to $Z$ form a partition of unity of $Z$. To see that the $\psi^{(n)}_i$ can be chosen to be uniformly Lipschitz, it suffices to check that the cover $\cU_n$ has non-zero Lebesgue number, compare the proof of Lemma \ref{LemmaLipPU} above. To that end, let $x\in X$. Assume first that there is $y\in Z$ with $d(x,y)<\frac{3\varepsilon_n}{4}$. Because $\Gamma_n$ is $\varepsilon_n$-dense, $y$ is contained in the $\varepsilon$-ball of some $y'\in \Gamma_n$. In this case, $d(x,y')<\frac{7\varepsilon_n}{4}$, so that the $\frac{\varepsilon_n}{4}$-ball around $x$ is contained in $U_{2\varepsilon_n}(y')$. If on the other hand $d(x,Z)\geq\frac{3\varepsilon_n}{4}$, then the $\frac{\varepsilon_n}{4}$-ball around $x$ is contained in the complement of the closed $\frac{\varepsilon_n}{2}$-neighborhood of $Z$. Hence $\cU_n$ has Lebesgue number at least $\frac{\varepsilon_n}{4}$.

Define $\sigma_n':\, C_0^+(Z)\to C_0^+(X)$ via $\sigma_n'(\mathbf{1}_Z):=\mathbf{1}_X$ and 
\begin{equation}
\sigma_n'(g) := \sum_{i} g(x_i^{(n)}) \cdot \psi^{(n)}_i 
\end{equation} 
if $g\in C_0(Z)$. Note that $\sigma_n'$ is a positive and contractive map. If $g$ has support of radius $R$, then $\sigma_n'(g)$ has radius at most $R+4\cdot\varepsilon_n$. Moreover, $\sigma_n'(g)$ is Lipschitz for any $g$:
\begin{align*}
\left|\sigma_n'(g)(x)-\sigma_n'(g)(y)\right| &= \left| \sum_i g(x_n^{i})\cdot\left( \psi^{(n)}_i(x) - \psi^{(n)}_i(y) \right) \right| \\
&\leq ||g||_\infty \cdot \sum_i \left| \psi^{(n)}_i(x) - \psi^{(n)}_i(y) \right| \\
&\leq 2l_n \cdot \mathrm{mult}(\cU_n)\cdot ||g||_\infty \cdot d(x,y) \; .
\end{align*}
In particular there exist $L_n',R_n'\geq 0$ such that $\sigma'_n$ maps $\LLip_R(Z)$ to $L'_n\text{-}\mathrm{Lip}_{R_n'}(X)$, using also the $\psi^{(n)}_i$ are compactly supported. We stress that this statement thus uses both the property (F) and the properness of $X$. Furthermore, $\sigma_n'(g)$ converges to $g$ on $Z$, and this convergence is uniform on $\LLip_R(Z)$ for all $L,R$:
\begin{align*}
|| \sigma_n'(g)|_Z - g || &= \sup_{x\in Z} \left| \sum_i \left(g(x^{(n)}_i) - g(x) \right)\cdot\psi^{(n)}_i(x) \right| \\
&\leq \sup_{x\in Z} \sum_i L \cdot d(x,x^{(n)}_i) \cdot \psi^{(n)}_i(x) \\
&\leq 2L\cdot \varepsilon_n \cdot \sup_{x\in Z} \sum_i  \psi^{(n)}_i(x) \\
&= 2L\cdot \varepsilon_n \xlongrightarrow{n\to \infty} 0 \; .
\end{align*}
By construction $\sigma_n'(g)$ is supported in the $2\varepsilon_n$-neighborhood of $Z$, so we cannot hope for convergence of $\sigma_n'(g)$ to a continuous function outside of $Z$. To remedy that, pick Lipschitz functions $\varphi_n: \, X\to [0,1]$ such that $\varphi_n= 1$ on $X\setminus N_{\frac{\varepsilon_n}{2}}(Z)$ and and $\varphi_n=0$ on $N_{\frac{\varepsilon_{n+1}}{2}}(Z)$. Then, set $\sigma_1:=\sigma_1'$ and inductively define 
\begin{equation}
\sigma_{n+1} := (1-\varphi_{n+1})\cdot \sigma_{n+1}' + \varphi_{n+1}\cdot \sigma_n \; .
\end{equation}
Each $\sigma_n$ is again positive and contractive, and $\sigma_n(\mathbf{1}_Z)=\mathbf{1}_X$. Moreover, since $\varphi_{n+1}$ is bounded and Lipschitz, $\sigma_n$ maps $\LLip_R(Z)$ to $\hat{L}_n\text{-}\mathrm{Lip}_{\hat{R}_n}(X)$ for appropriate $\hat{L}_n,\hat{R}_n\geq 0$. Since $\varphi_{n+1} = 0$ on $Z$ it holds that $\sigma_n(g)|_Z=\sigma_{n+1}'(g)|_Z$, whence $\sigma_n(g)|_Z\to g$ and this convergence is uniform on each $\LLip_R(Z)$. We will now compute that $\sigma_n(g)$ forms a Cauchy sequence whose rate of convergence is again uniform on each $\LLip_R(Z)$. Note first that $(1-\varphi_{n+1})\sigma_n=(1-\varphi_{n+1})\sigma_n'$, following from the fact that $(1-\varphi_{n+1})\varphi_{n}=0$ by construction. Then proceed as follows:
\begin{align*}
\left|\left|\sigma_{n+1}(g)-\sigma_n(g)\right|\right| &= \left|\left| (1-\varphi_{n+1})\sigma_{n+1}'(g) + \varphi_{n+1}\cdot \sigma_n(g) - \sigma_n(g)\right|\right| \\
&= \left|\left|(1-\varphi_{n+1})(\sigma_{n+1}'(g)-\sigma_n(g))\right|\right| \\
&= \left|\left|(1-\varphi_{n+1})(\sigma_{n+1}'(g)-\sigma_n'(g))\right|\right| \\
&\leq \sup_{x\in N_{\frac{\varepsilon_{n+1}}{2}}(Z)} \left|\sigma_{n+1}'(g)(x)-\sigma_n'(g)(x)\right| \\
&= \sup_{x\in N_{\frac{\varepsilon_{n+1}}{2}}(Z)} \left|\sum_j g(x^{(n+1)}_j)\cdot\psi^{(n+1)}_j(x) - \sum_i g(x^{(n)}_i)\cdot\psi^{(n)}_i(x)\right| 
\end{align*}
Now, $\{\psi^{(n)}_i\}_i$ and $\{\psi^{n+1}_j\}_j$ form a partition of unity on $N_{\frac{\varepsilon_{n+1}}{2}}(Z)$, as $\psi^{(n)}_\infty=0=\psi^{(n+1)}_\infty$ there. Thus we compute further:
\begin{align*}
\left|\left|\sigma_{n+1}(g)-\sigma_n(g)\right|\right| &= \sup_{x\in N_{\frac{\varepsilon_{n+1}}{2}}(Z)} \left|\sum_{i,j} \left(g(x^{(n+1)}_j) - g(x^{(n)}_i)\right)\cdot\psi^{(n+1)}_j(x) \cdot\psi^{(n)}_i(x)\right| \\
&\leq \sup_{x\in N_{\frac{\varepsilon_{n+1}}{2}}(Z)} \sum_{i,j} L\cdot \left(d(x^{(n+1)}_j,x) + d(x,x^{(n)}_i)\right)\cdot\psi^{(n+1)}_j(x) \cdot\psi^{(n)}_i(x) \\
&\leq 2L(\varepsilon_{n+1}+\varepsilon_n) \cdot \sup_{x\in N_{\varepsilon_{n+1}}(Z)} \sum_{i,j} \psi^{(n+1)}_j(x) \cdot\psi^{(n)}_i(x) \\
&= 2L(\varepsilon_{n+1}+\varepsilon_n) < \frac{L}{2^n} \; .
\end{align*}
It follows that $||\sigma_{n+m}(g)-\sigma_n(g)||\leq L\sum_{k=n}^\infty 2^{-n} \xrightarrow{n\to\infty} 0$. Since the union $\bigcup_{L,R}\LLip_R(Z)$ is dense in the unit ball of $C_0(Z)$, and each $\sigma_n$ is linear and contractive, it follows that $\sigma(g):=\lim_{n\to\infty}\sigma_n(g)$ exists for all $g\in C_0(Z)$. The resulting map $\sigma: \, C_0^+(Z)\to C_0^+(X)$ is again positive and contractive. Since $\sigma_n(g)|_Z\to g$ it holds that $\sigma(g)|_Z=g$. Lastly, for all $L,R\geq 0$ and all $\varepsilon>0$ there exists $n\in\NN$ such that $||\sigma(g)-\sigma_n(g)||<\varepsilon$ for all $g\in \LLip_R(Z)$, and by construction $\sigma_n$ maps $\LLip_R(Z)$ to $L'\text{-}\mathrm{Lip}_{R'}(X)$ for appropriate $L',R'\geq 0$. This completes the construction.  
\end{proof}
The next next lemma confirms that uniform resolvability and transversal control are indeed notions that depend only on the geometry near $Z$. 
\begin{lemma}
Let $(X,Z)$ be a pair, and $U\subseteq X$ a neighborhood of $Z$ such that $d(Z,X\setminus U)>\delta$. If $U$ is uniformly resolvable resp. transversally controlled near $Z$, then so is $X$.
\end{lemma}
\begin{proof}
Let $\delta:=d(Z,X\setminus U)$, and fix a Lipschitz map $\varphi: \, X\to [0,1]$ such that $\sqrt{\varphi}$ is also Lipschitz, $\varphi=1$ on $N_{\delta/2}(Z)$ and $\varphi=0$ on $X\setminus U_\delta(Z)$. Let $\{\psi_i\}_i$ be a partition of unity over $U\setminus Z$ as in Definition \ref{DefUnifPair}. Then $\{\psi_i\varphi\}_i \cup\{1-\varphi\}$ is a partition of unity over $X\setminus Z$ that also satisfies the requirements of Definition \ref{DefUnifPair}. Moreover, suppose $\sigma: \, C_0^+(Z)\to C_0^+(U)$ is an approximately filtered positive section. Then, the map $\sigma': \, C_0^+(Z)\to C_0^+(X)$ defined by
$$
\sigma'(\lambda \mathbf{1}_Z + g ) := \varphi\cdot \sigma(\lambda \mathbf{1}_Z + g) + \lambda\cdot (1-\varphi)
$$
for $\lambda\in\CC$, $g\in C_0(Z)$, is also an approximately filtered positive section. This concludes the proof.
\end{proof}
We come to more concrete examples.
\begin{lemma}\label{TCcone}
Let $Z$ be arbitrary, and consider $CZ=[0,1)\times Z$, equipped with the product metric. Then, $CZ$ is transversally controlled near $Z\cong Z\times \{0\}$.
\end{lemma}
\begin{proof}
Fix a Lipschitz function $\varphi: \, [0,1]\to [0,1]$ with $\varphi(0)=1$ and $\varphi(x)=0$ for $x\geq \frac{1}{2}$. Then, $\sigma: \, C_0^+(Z)\to C_0^+(CZ)$ defined by 
$$
\sigma(\mathbf{1}_Z) = \mathbf{1}_{CZ} \; , \; \sigma(g)(t,x) = \varphi(t)\cdot g(x) 
$$
for $g\in C_0(Z)$ has the desired properties.
\end{proof}
The argument of this lemma can be generalized. Instead of crossing with $[0,1)$ one can instead cross with an arbitrary (non-empty) metric space $Y$. Embedding $Z$ as $Z\cong \{y\}\times Z$ for some $y\in Y$ and taking $\varphi: \, Y\to [0,1]$ a compactly supported Lipschitz function with $\varphi(y)=1$ allows the same argument to prove that $Y\times Z$ is transversally controlled near $Y\cong\{y\}\times Z$. We conclude that any metric space is transversally controlled near any subspace with uniformly thick tubular neighborhoods (interpreted suitably).
\begin{cor}\label{TubNbhd}
Let $X$ be any metric space, and $Z\subseteq X$ a closed subspace. Suppose there exists an open neighborhood $U\subseteq X$ of $Z$ such that $d(Z,X\setminus U)>0$, and $U$ is bi-Lipschitz equivalent to $Y\times Z$ with product metric for some $Y$ via a map that takes $Z\subseteq U$ to $\{y\}\times Z \subseteq Y\times Z$ for some $y\in Y$. Then $X$ is transversally controlled near $Z$.
\end{cor}
Let us conclude this appendix with the observation that every pair of interest in the body of this thesis has the properties of Definition \ref{DefUnifPair}.
\begin{cor}
Let $(X,Z)$ be one of the following:
\begin{itemize}
\item $X$ is a manifold of bounded geometry, and $Z$ a bounded-geometry hypersurface.
\item $X$ is a simplicial complex of bounded geometry, and $Z$ is a subcomplex.
\item $Z$ has jointly bounded geometry, $X=CZ=[0,1)\times Z$ with the product metric, and $Z\cong Z\times\{0\} \subset CZ$.
\end{itemize}
Then $X$ is uniformly resolvable and transversally controlled near $Z$.
\end{cor}
\begin{proof}
Manifolds of bounded geometry have jointly bounded geometry \cite[Example 2.21]{Engel2019}, hence are uniformly resolvable near any closed subspace. Moreover, the metric on $X$ is bi-Lipschitz equivalent to one that has product structure on a uniformly thick tubular neighborhood of $Z$ (\cite[Prop. 7.3]{Schick1998} and \cite[Theorem 2.10]{Ammann2019}). By Lemma \ref{DifferentMetrics}(iv), this metric is bi-Lipschitz equivalent to the product metric on that tubular neighborhood. Hence, Corollary \ref{TubNbhd} implies that $X$ is transversally controlled near $Z$.
 
Simplicial complexes of bounded geometry have jointly bounded geometry \cite[Example 2.21]{Engel2019}, hence are uniformly resolvable near any closed subspace. Moreover, any subcomplex $Z$ satisfies the property (F) of Proposition \ref{FilteredCP}. To see this, fix $\varepsilon>0$. Let $n$ be the dimension of $X$. Inductively choose finite $\varepsilon$-dense subsets on the $k$-simplex, $k\leq n$, such that the intersection of this subset of the $k$-simplex with any face coincides with the chosen subset on the $(k-1)$-simplex. Then, define $\Gamma$ to be the subset of $Z$ that coincides with the chosen subset on any simplex in $Z$. Because $X$ has bounded geometry, there is a uniform bound on how many elements of $\Gamma$ can be within $2\varepsilon$ of a given $x\in X$. Hence $Z\subseteq X$ has the property (F). Since $X$ is also proper, Proposition \ref{FilteredCP} gives the desired result. 

Lastly, if $Z$ has jointly bounded geometry, so does $CZ$. Indeed, if $\Gamma\subseteq Z$ is a quasi-lattice in $Z$, it is also one in $CZ$. Moreover, if $\{X_i\}$ is a countable Borel decomposition of $Z$ as in Definition \ref{DefJBB}, then $\{[0,1)\times X_i\}$ is a decomposition of $CZ$ with the analogous properties. Hence, $CZ$ is uniformly resolvable near any closed subspace. It is transversally controlled near $Z$ by Lemma \ref{TCcone}. 
\end{proof}
\addcontentsline{toc}{chapter}{Bibliography}
\printbibliography
\end{document}